\documentclass[a4paper,11pt]{book}

\usepackage[T1]{fontenc}
\usepackage[utf8]{inputenc}
\usepackage[french,german,english]{babel}

\usepackage{fourier} 
\usepackage{setspace} 

\usepackage{graphicx,xcolor}
\graphicspath{{images/}}

\usepackage{subfig}
\usepackage{booktabs}
\usepackage{lipsum}
\usepackage{microtype}
\usepackage{url}
\usepackage[final]{pdfpages}

\usepackage{fancyhdr}

\fancypagestyle{plain}{
	\fancyhf{}

	\fancyfoot[R]{\thepage}
	}
\fancypagestyle{addpagenumbersforpdfimports}{
	\fancyhead{}
	
	\fancyfoot{}
    \fancyfoot[R]{\thepage}
}

\usepackage{listings}
\usepackage{hyperref}
\hypersetup{pdfborder={0 0 0},
	colorlinks=true,
	linkcolor=black,
	citecolor=black,
	urlcolor=black}
\makeatletter
\def\cleardoublepage{\clearpage\if@twoside \ifodd\c@page\else
    \hbox{}
    \thispagestyle{empty}
    \newpage
    \if@twocolumn\hbox{}\newpage\fi\fi\fi}
\makeatother \clearpage{\pagestyle{plain}\cleardoublepage}

\usepackage{color}
\usepackage{tikz}
\usepackage[explicit]{titlesec}
\newcommand*\chapterlabel{}
\titleformat{\chapter}[display]  
	{\normalfont\bfseries\Huge} 
	{\gdef\chapterlabel{\thechapter\ }}     
 	{0pt} 
 	  {\begin{tikzpicture}[remember picture,overlay]
    \node[yshift=-8cm] at (current page.north west)
      {\begin{tikzpicture}[remember picture, overlay]
        \draw[fill=black] (0,0) rectangle(35.5mm,15mm);
        \node[anchor=north east,yshift=-7.2cm,xshift=34mm,minimum height=30mm,inner sep=0mm] at (current page.north west)
        {\parbox[top][30mm][t]{15mm}{\raggedleft $\phantom{\textrm{l}}$\color{white}\chapterlabel}};  
        \node[anchor=north west,yshift=-7.2cm,xshift=37mm,text width=\textwidth,minimum height=30mm,inner sep=0mm] at (current page.north west)
              {\parbox[top][30mm][t]{\textwidth}{\color{black}#1}};
       \end{tikzpicture}
      };
   \end{tikzpicture}
   \gdef\chapterlabel{}
  } 

\titlespacing*{\chapter}{0pt}{50pt}{30pt}
\titlespacing*{\section}{0pt}{13.2pt}{*0}  
\titlespacing*{\subsection}{0pt}{13.2pt}{*0}
\titlespacing*{\subsubsection}{0pt}{13.2pt}{*0}

\newcounter{myparts}
\newcommand*\partlabel{}
\titleformat{\part}[display]  
	{\normalfont\bfseries\Huge} 
	{\gdef\partlabel{\thepart\ }}     
 	{0pt} 
 	  {\setlength{\unitlength}{20mm}
	  \addtocounter{myparts}{1}
	  \begin{tikzpicture}[remember picture,overlay]
    \node[anchor=north west,xshift=-65mm,yshift=-6.9cm-\value{myparts}*20mm] at (current page.north east) 
      {\begin{tikzpicture}[remember picture, overlay]
        \draw[fill=black] (0,0) rectangle(62mm,20mm);   
        \node[anchor=north west,yshift=-6.1cm-\value{myparts}*20mm,xshift=-60.5mm,minimum height=30mm,inner sep=0mm] at (current page.north east)
        {\parbox[top][30mm][t]{55mm}{\raggedright \color{white}Part \partlabel $\phantom{\textrm{l}}$}};  
        \node[anchor=north east,yshift=-6.1cm-\value{myparts}*20mm,xshift=-63.5mm,text width=\textwidth,minimum height=30mm,inner sep=0mm] at (current page.north east)
              {\parbox[top][30mm][t]{\textwidth}{\raggedleft \color{black}#1}};
       \end{tikzpicture}
      };
   \end{tikzpicture}
   \gdef\partlabel{}
  } 

\usepackage{amsmath,amssymb,amsfonts}%
\usepackage{amsthm}%
\makeatletter
\def\resetMathstrut@{%
  \setbox\z@\hbox{%
    \mathchardef\@tempa\mathcode`\(\relax
      \def\@tempb##1"##2##3{\the\textfont"##3\char"}%
      \expandafter\@tempb\meaning\@tempa \relax
  }%
  \ht\Mathstrutbox@1.2\ht\z@ \dp\Mathstrutbox@1.2\dp\z@
}
\makeatother
\usepackage{mathtools}
\usepackage{imakeidx}
\usepackage{needspace}

\newtheoremstyle{break}
  {\topsep}{\topsep}%
  {\itshape}{}%
  {\bfseries}{}%
  {\newline}{}%
\theoremstyle{break}

\newtheorem{theorem}{Theorem}[chapter]
\newtheorem{lemma}[theorem]{Lemma}
\newtheorem{proposition}[theorem]{Proposition}
\newtheorem{corollary}[theorem]{Corollary}
\newtheorem{definition}[theorem]{Definition}
\newtheorem{notation}[theorem]{Notation}
\newtheorem{remark}[theorem]{Remark}

\DeclareMathOperator{\spec}{spec}
\DeclarePairedDelimiter\ceil{\lceil}{\rceil}
\DeclarePairedDelimiter\floor{\lfloor}{\rfloor}
\DeclareMathOperator{\dom}{Dom}
\let\Im\relax
\DeclareMathOperator{\Im}{Im}
\let\Re\relax
\DeclareMathOperator{\Re}{Re}
\DeclareMathOperator{\dist}{dist}
\DeclareMathOperator{\Tr}{Tr}
\DeclareMathOperator{\supp}{supp}
\usepackage{xparse, etoolbox}
\DeclarePairedDelimiterX\innerp[2]{\langle}{\rangle}{#1,#2}
\DeclareMathOperator*{\fsum}{\sum^{\mathrm{finite}}}
\DeclareMathOperator{\kadm}{\mathit{k}-adm}

\numberwithin{equation}{chapter}

\makeatletter
\newcommand{\subalign}[1]{%
  \vcenter{%
    \Let@ \restore@math@cr \default@tag
    \baselineskip\fontdimen10 \scriptfont\tw@
    \advance\baselineskip\fontdimen12 \scriptfont\tw@
    \lineskip\thr@@\fontdimen8 \scriptfont\thr@@
    \lineskiplimit\lineskip
    \ialign{\hfil$\m@th\scriptstyle##$&$\m@th\scriptstyle{}##$\hfil\crcr
      #1\crcr
    }%
  }%
}
\makeatother

\makeindex[columns=2, title=Index of Notations]
\usepackage[version=4]{mhchem}
\usepackage{chemist}
\usepackage{chemfig}

\begin{document}
\frontmatter
\begin{titlepage}
\begin{center}
\sffamily

\null\vspace{2cm}
{\huge On blow-up dynamics for the Nonlinear Wave Equation and the Harmonic Map Heat Flow} 
    
\vfill

\begin{tabular} {cc}
\parbox{0.3\textwidth}{\includegraphics[width=4cm]{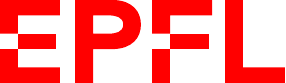}}
&
\parbox{0.7\textwidth}{%
	Thèse présentée le 30 Avril 2026\\
	à la Faculté des Sciences de Base\\
	chaire des Équations Différentielles Partielles\\
	programme doctoral en Mathématiques\\
%
	École Polytechnique Fédérale de Lausanne\\[6pt]
	pour l'obtention du grade de Docteur ès Sciences\\
	par\\ [4pt]
	\null \hspace{3em} Dylan Samuelian\\[9pt]
\small
acceptée sur proposition du jury:\\[4pt]
    Prof Eisenbrand Friedrich, président du jury\\
    Prof Krieger Joachim, directeur de thèse\\
    Prof Moschidis Georgios, rapporteur\\
    Prof Schoerkhuber Birgit, rapporteur\\
    Prof Lawrie Andrew, rapporteur\\[12pt]
Lausanne, EPFL, 2026}
\end{tabular}
\end{center}
\vspace{2cm}
\end{titlepage}

\cleardoublepage
\thispagestyle{empty}

\vspace*{3cm}

\begin{raggedleft}
\raggedleft{
\begin{tabular}{l@{}}
    	Que fait une onde pour se réchauffer ? \\
        Elle se concentre. \quad --- Antoine Schaub\\
\end{tabular}\par}
\end{raggedleft}

\vspace{4cm}

\begin{center}
    A ma famille et, en particulier, à mon défunt père qui serait certainement \\
    le plus fier des papas en me voyant devenir Docteur.
\end{center}

\setcounter{page}{0}
\chapter*{Acknowledgements}
~\newline~\newline~
\markboth{Acknowledgements}{Acknowledgements}
\addcontentsline{toc}{chapter}{Acknowledgements}

I'm deeply thankful to my advisor Joachim Krieger, without whom this thesis would have been impossible to complete. His deep knowledge of mathematics, even outside of his specialization, is impressive and I was always able to count on his support when faced with a roadblock. He is kind-hearted and I was able to work without feeling any sort of pressure. I also thank him for giving me the opportunity to do a PhD after doing my master project in his chair. 

I also want to thank my advisor and my colleagues (Yang, Katie, Jaime, Tim, Tobias, Shengquan, David and Jose) for making the work environment a friendly place, where I don't struggle to wake up and go in the morning. We joke a lot about my schedule and me not being in the office, but this is definitely because of my activities and not the work environment! 

Talking about activities, I want to thank the people who introduced me to music, as well as dancing, during my thesis: my voice teacher (Hélène), my music theory teacher (Maude), my piano teachers (Corentin, Florian, Frédéric), my oboe teacher (Moana) and oboe duet partner (Léa), the university choir (CHUL) and the Rock Affinity association. Music and dancing have been an incredible discovery in my life and I would not have been able to relieve the pressure from the PhD without them. I've also met incredible friends thanks to these social activities. It's just so much fun and I cannot stop now!

Finally, I thank my family, especially my mom (Silva), whose logistical and financial support at home was of great help. It saved me countless hours and I was able to focus on my family, my PhD and my social activities without worrying about having a roof or a warm meal. She is much more than just some help. I love her and I hope she is proud of having a Doctor in the family! 

\bigskip
 
\noindent\textit{Montreux, 30 avril 2026}
\hfill Dylan Samuelian


\cleardoublepage
\chapter*{Abstract}
~\newline~\newline~
\addcontentsline{toc}{chapter}{Abstract (English/Français)} 

This thesis concerns the blow-up analysis, and more precisely, the asymptotic decomposition into profiles of solutions to nonlinear partial differential equations. In Chapter 1, we present a unifying framework for finding profile decompositions in Hilbert spaces. The next chapters are then devoted to my publications.

In the second chapter, based on \cite{samuelian2026construction}, we construct solutions $u(x,t)$ to the focusing, energy-critical, nonlinear wave equation 
\begin{equation} \partial_{tt}u - \Delta u - |u|^{p-1}u = 0, \quad t \geq 0, \ x \in \mathbb{R}^d, \ d \geq 3, \ p = (d+2)/(d-2) 
\end{equation} in dimension $d \in \{4,5\}$, exhibiting finite-time Type II blow-up precisely at $x = t = 0$ with a single profile and a prescribed polynomial blow-up rate of $t^{-1-\nu}$, where $\nu > 1$ for $d = 4$ and $\nu > 3$ for $d = 5$. Such solutions have been constructed by Krieger-Schlag-Tataru for $d = 3$ and by Jendrej for $d = 5$. The major difference between dimension $4$ and $5$ consists in the renormalization procedure. In $d = 4$, we essentially follow the Krieger-Schlag-Tataru scheme developed for the 3-dimensional equation. This scheme has been applied with success for other equations such as the 3D-critical NLS, Schrödinger maps or wave maps. In all of these cases, the polynomial structure of the nonlinearity permits the use of simple algebraic manipulations to control error terms. By contrast, the case $d = 5$ requires a modified setup due to the lower regularity of the nonlinearity, which complicates the treatment of nonlinear error terms.

The remaining chapters are based on \cite{samuelian2025blowuptrees}. We first provide, in Chapter 3, a modern presentation of classical $L^p$ and Schauder estimates for parabolic equations. Then, in Chapter 4, we consider finite-time and $k$-equivariant solutions to the harmonic map heat flow from $B^2$ to $S^2$: 
\begin{align}
    v_t &= \Delta v + |\nabla v|^2 v, \quad (x,t) \in B^2 \times  (0,+\infty)  \\
    v(x,t) &= v_0(x,t), \quad (x,t) \in \partial B^2 \times [0,+\infty) \cup B^2 \times \{0\}. \notag
\end{align}
Under general time-dependent, $k$-equivariant and smooth boundary data, we prove that the profile decomposition at the singular time contains only one bubble. The method relies on the maximum and comparison principle. We also prove a global existence criterion and exhibit solutions blowing up in infinite time for any degree of equivariance. Finally, the appendix presents the profile decomposition for the harmonic map heat flow proved by Struwe and Qing using energy arguments and the Palais-Smale condition.

\vskip0.5cm
Key words: partial differential equations, dispersive, nonlinear wave equation, energy-critical, focusing, parabolic, harmonic map heat flow, profile decomposition, bubble, soliton, blow-up

\begin{otherlanguage}{french}
\cleardoublepage
\chapter*{Résumé}
~\newline~\newline~

Cette thèse porte sur l’analyse des explosions et, plus précisément, sur la décomposition asymptotique en profils de solutions d’équations aux dérivées partielles non linéaires. Dans le Chapitre 1, nous présentons un cadre général pour les décompositions en profils dans les espaces de Hilbert. Les chapitres suivants sont consacrés à mes publications.

Dans le Chapitre 2, basé sur \cite{samuelian2026construction}, nous construisons des solutions $u(t,x)$ à l’équation d'onde non linéaire focalisante, énergie-critique,
\begin{equation}
    \partial_{tt}u - \Delta u - |u|^{p-1}u = 0,
    \quad t \geq 0,\quad x \in \mathbb{R}^d,\quad d \geq 3,\quad p = \frac{d+2}{d-2},
\end{equation}
en dimensions $d \in \{4,5\}$, qui présentent une explosion de Type II en temps fini précisément en $(t,x)=(0,0)$, avec un seul profil et un taux polynomial prescrit de $t^{-1-\nu}$, où $\nu > 1$ pour $d = 4$ et $\nu > 3$ pour $d = 5$. De telles solutions avaient précédemment été construites par Krieger-Schlag-Tataru pour $d = 3$ et par Jendrej pour $d = 5$. La différence principale entre les dimensions $4$ et $5$ réside dans la procédure de renormalisation. En dimension $d = 4$, nous suivons essentiellement la méthode de Krieger-Schlag-Tataru développée pour l’équation en 3 dimensions. Cette méthode a été appliquée avec succès à d’autres équations, telles que l'équation de Schrödinger non linéaire (3D critique), les applications de Schrödinger et les applications d’ondes. Dans tous ces cas, la structure polynomiale de la non-linéarité permet d’utiliser des manipulations algébriques simples pour contrôler les termes d’erreur. En revanche, le cas $d = 5$ nécessite un cadre modifié en raison de la plus faible régularité de la non-linéarité, ce qui complique le traitement des termes d’erreur non linéaires.

Les chapitres restants sont basés sur \cite{samuelian2025blowuptrees}. Nous présentons d’abord, dans le Chapitre 3, une exposition moderne des estimations de Schauder classiques pour les équations paraboliques. Ensuite, dans le Chapitre 4, nous considérons des solutions $k$-équivariantes avec singularité en temps fini du flot de chaleur des applications harmoniques de $B^2$ vers $S^2$,
\begin{align}
    v_t &= \Delta v + |\nabla v|^2 v,
    \quad (x,t) \in B^2 \times (0,+\infty), \\
    v(x,t) &= v_0(x,t),
    \quad (x,t) \in \partial B^2 \times [0,+\infty) \cup B^2 \times \{0\}. \notag
\end{align}
Nous travaillons avec des données de bord générales, dépendant du temps, $k$-équivariantes, lisses et nous démontrons que la décomposition en profils au temps singulier ne contient qu’une seule bulle d'énergie. La méthode repose sur les principes du maximum et de comparaison. Nous prouvons également un critère d'existence global et exhibons des solutions qui explosent en temps infini pour tout degré d’équivariance. Enfin, l’annexe présente la décomposition en profils pour le flot de chaleur des applications harmoniques démontrée par Struwe et Qing à l’aide d’arguments d’énergie et de la condition de Palais-Smale.

\vskip0.5cm
Mots-clés : équations aux dérivées partielles, équations dispersives, équation d'onde non linéaire, énergie-critique, focalisante, équations paraboliques, flot de chaleur des applications harmoniques, décomposition en profils, bulle, soliton, explosion.
\end{otherlanguage}


\tableofcontents

\setlength{\parskip}{1em}

\mainmatter
\chapter{Introduction : Lack of compactness and profile decompositions}

This thesis concerns the analysis of solutions $u(t,x)$ to evolution partial differential equations which exhibit a singular behaviour either in finite or global time, while keeping a finite "energy". If $T \in \mathbb R \cup \{\pm \infty\}$ is a singular time, the fundamental issue in understanding the asymptotic behaviour of $u(t,x)$ as $t \to T$ stems from the lack of compactness in infinite-dimensional Banach spaces, such as Sobolev spaces. On contrary to the finite-dimensional case, the classical Bolzano-Weierstrass theorem fails. Focusing on the case of Hilbert spaces, a bounded sequence in a Hilbert space $H$ might not have a convergent subsequence, even though Banach-Alaoglu's theorem ensures it has a weakly-convergent subsequence. This is called \textit{lack of compactness}. One way to fix this is to look at a compact embedding $H \hookrightarrow B$. If such a compact embedding exists, then the subsequence does converge strongly in $B$. However, the issue remains when there is no such embedding, which is typically the case for Sobolev embeddings into $L^p(\mathbb R^n)$ spaces. For those spaces, the lack of compactness is easy to exhibit: take a compactly supported function $u(x) \in W^{k,p}(\mathbb R^n)$, $1 < p < +\infty$, and a sequence of points $x_n \in \mathbb R^n$ for which $|x_n| \to +\infty$. Then the sequence of functions $u_n(x) = u(x + x_n)$ converges weakly to zero, while $||u_n||_{W^{k,p}} = ||u||_{W^{k,p}} \neq 0$ stays constant.

Both the lack of compactness of this embedding and the lack of compactness for bounded sequences in a Sobolev space is due, at least partly, to the action of a group of isometries on the space. In such a setting, the lack of compactness can be described by the following profile decomposition theorem (Theorem \ref{thm:abstract profile decomposition}).

\begin{definition}[Group of dislocations]
    Let $H$ be a Hilbert space. Let $G$ be a locally sequentially compact group acting unitarily on $H$ and assume that the group action map $G \times H \to H$, $(g,u) \mapsto gu$ is jointly continuous. We say that $G$ is a group of \textit{dislocations} if for any sequence $(g_n) \in G$ which escapes any compact set, the sequence of induced isometries converges operator-weakly to zero, i.e., $g_nu \to 0$ weakly in $H$ for any $u \in H$.
\end{definition}

\begin{remark}
    The trivial group is a group of dislocations on any Hilbert space $H$. The translation action $u(x) \mapsto u(x - x_0)$, $x_0 \in \mathbb R^n$, and the scaling action $u(x) \mapsto \lambda^{-\frac{n-2s}{2}}u(x/\lambda)$, $\lambda > 0$, (or the composition of both) create a group of dislocations for the homogeneous Sobolev space $\dot{H}^s(\mathbb R^n)$, as well as $L^{\frac{2n}{n-2s}}(\mathbb R^n)$ (this is technically not a Hilbert space, but the definition can be generalized, see Remark \ref{rem:generalization dislocation}).
\end{remark}

\begin{theorem}[Profile Decomposition]\label{thm:abstract profile decomposition}
    Let $H$ be a Hilbert space and $G$ be a group of dislocations acting on $H$. Let $(x_n)_{n \in \mathbb N} \subset H$ be a bounded sequence. After passing to a subsequence (still denoted with index $n$) of $(x_n)_{n \in \mathbb N}$, there exists profiles $(\phi_{m})_{m \in \mathbb N} \subset H$, group elements $((g_{m,n})_{m \in \mathbb N})_{n \geq m} \subset G$, remainders  $((w_{m,n})_{m \in \mathbb N})_{n \geq m} \subset H$ for which
$$
x_{n} = \sum_{m=1}^M g_{m,n} \phi_m + w_{M,n}, \quad M \in \mathbb N_{\geq 1}, n \geq M.
$$
Moreover, one has:
\begin{enumerate}
    \item \textbf{Pythagorean decomposition of energy:}
    $$
    \limsup_{n \to +\infty} ||x_{n}||^2  =  \sum_{m=1}^M||\phi_m||^2 + \limsup_{n \to +\infty} ||w_{M,n}||^2, \quad M \in \mathbb N_{\geq 1}.
    $$
    \item \textbf{Finiteness of decomposition and non-zero profiles:} Let $w_{0,n} = x_n$ and
    \begin{equation}
    D[(x_n)] := \sup\{||\phi|| : \exists \phi \in H, \exists (n_k)_{k \in \mathbb N} \subset \mathbb N, (g_k)_{k \in \mathbb N} \subset G, g^{-1}_{k}x_{n_k} \to \phi \text{ weakly } \}. \label{eq:D(x_n) def}
\end{equation}
    For $M \in \mathbb N_{\geq 0}$, the extracted profile $\phi_{M+1}$ is non-zero if and only if $D[(w_{M,n})] > 0$. Moreover, if $D[(w_{M,n})] = 0$, then all the subsequent profiles $\phi_{M+1} = \phi_{M+2} = ... = 0$ are zero and the decomposition is finite. 
    \item \textbf{Separation of scales:} Fix $i \neq j$ and assume that $\phi_i, \phi_j$ are non-zero profiles. Then $(g^{-1}_{i,n} g_{j,n})_{n \in \mathbb N}$  escapes any compact set.
    \item \textbf{Weak convergence to profile:} For any $M \in \mathbb N_{\geq 1}$, if either $M = 1$ or $\phi_M$ is a non-zero profile, then $g^{-1}_{M,n}x_n \to \phi_M$ weakly.
    \item \textbf{Smallness of remainder:} For $1 \leq m \leq M$, $g^{-1}_{m,n}w_{M,n} \to 0$ weakly. Moreover, the remainders also converge weakly to zero in the following sense:
    $$
    \lim_{M\to +\infty}\sup_{||y|| = 1}  \limsup_{n \to +\infty} \sup_{g \in G} |\langle g^{-1}w_{M,n}, y\rangle| = 0.
    $$
\end{enumerate}
\end{theorem}

\begin{remark}\label{rem:generalization dislocation}
    The abstract decomposition in a Hilbert space and the notion of dislocation group were developped by Fieseler-Schindler-Tintarev (\cite{schindlertintarev2002abstract}, \cite{TintarevFieseler2007}). The theorem generalizes to Banach spaces and to dislocation sets (i.e. without the group structure), see the book of Tintarev (\cite{tintarev2020concentration}), as well as the papers \cite{okumura2022profile} and \cite{schindlertintarev2002abstract}, for an exposition on the topic.
\end{remark}

\begin{remark}
    Using a trivial group of dislocations, one recovers Banach-Alaoglu's theorem on $H$. In particular, passing first to a subsequence $x_{n_k}$ which converges weakly to $x_0$, one can always apply the theorem to $x_{n_k} - x_0$ to a get a profile decomposition of $x_{n_k}$ in terms of the weak limit $x_0$, a sum of "rescaled" profiles, and a remainder (see Corollary \ref{thm:abstract profile decomposition with weak limit first}).
\end{remark}

The work on understanding the lack of compactness in $L^p$ and Sobolev spaces goes back to the eighties with the pioneer work of Lions \cite{Lions1984}. Given a sequence of nonnegative probability densities $\rho_n$, Lions' first concentration-compactness lemma establishes a strict trichotomy \cite[Lemma I.1]{Lions1984}: either (i) the sequence is uniformly vanishing on balls of fixed radius (i.e., $\sup_{y}\int_{B_R(y)}\rho_n(x)dx \to 0$ as $n \to +\infty$), either (ii) the mass (up to some $\varepsilon$) is concentrated on a sequence of translated balls of fixed radius $R(\varepsilon)$ (compactness) or (iii) the sequence splits into two densities with almost-disjoint support (dichotomy). In the language of profile decompositions, this would roughly translate to either no non-trivial profile (vanishing), one dominant profile (compactness) or at least two asymptotically orthogonal profiles (dichotomy). At the same time, Brezis-Lieb Lemma (\cite{brezislieb1983}) provided a prototype for energy-decoupling, showing how the norm of a bounded sequence in $L^p$ splits when it converges pointwise almost-everywhere to a weak limit. 

Around the same period, Struwe (\cite{struwe1984global}) provided the first application to critical elliptic variational equations, using the so-called Palais-Smale condition. For a critical variational problem defined by an energy functional $E$, he showed that any sequence $u_n$ with bounded energy $\sup \limits_n E(u_n) < +\infty$ and $dE(u_n) \to 0$ admits a profile decomposition. The variational structure allows to identify the weak-limit of $u_n$ as a solution to the problem associated to $E$ and the profiles as solutions of a limiting autonomous equation. Moreover, these profiles carry a lower-bounded mass determined by the optimal Sobolev constant, which implies that the decomposition is finite. The approach of Struwe was later generalized to parabolic problems, an approach described in Appendix \ref{section:struwe energy hmhf}. Fifteen years later, Bahouri and Gérard (\cite{bahourigerard1999}) adapted this framework to critical dispersive equations. They proved a general profile decomposition for bounded sequences in $L^2(\mathbb R^n)$ using Fourier techniques and applied it to describe the asymptotic behaviour of sequences of solutions to the defocusing, energy-critical, 3D-nonlinear wave equation. 

From this abstract machinery, several natural questions arise:
\begin{enumerate}
    \item How many profiles can occur ?
    \item What are the profiles ? Are they solutions of a limiting problem ?
    \item Is the decomposition unique (up to reordering and symmetries), i.e., do we get the same decomposition if we change the subsequence ? In case of a solution $u(t,x)$ to an evolution problem for which we can get a profile decomposition along a sequence of times $t_n \to T$, can the decomposition be taken continuously in time ?
    \item In PDE applications, is the decomposition stable under perturbations of the initial data or of the equation ?
\end{enumerate}

The two papers presented in this thesis study complementary sides of this profile-decomposition problem. The wave equation paper \cite{samuelian2026construction} constructs a finite-time blow-up solution which, near the blow-up time, admits a continuous one-profile decomposition with a prescribed blow-up rate. The harmonic map heat flow paper \cite{samuelian2025blowuptrees} starts instead from the existence of a profile decomposition for a finite-time singular solution and proves, using the maximum principle, that such a decomposition can only contain one profile. Thus, the first paper constructs a prescribed continuous-in-time profile dynamics, while the second proves a rigidity theorem for the possible profile configurations.

We now complete the proof sketch for Theorem \ref{thm:abstract profile decomposition} given by \cite[Theorem 4.5.3]{tao2013compactness}, starting with the following lemmas. Section 2 in \cite{okumura2022profile} was also useful for filling in some gaps.

\begin{lemma}\label{lemma:profile decomposition, first profile}
    Let $H$ be a Hilbert space and $G$ be a group of dislocations acting on $H$. Let $(x_n)_{n \in \mathbb N} \subset H$ be a bounded sequence. For every $\varepsilon > 0$, there exists a subsequence $(n_k)_{k \in \mathbb N} \subset \mathbb N$, a sequence $(g_k)_{k \in \mathbb N} \subset G$, a profile $\phi \in H$ for which
    $$
    x_{n_k} = g_{k} \phi + w_k, \quad k \in \mathbb N
    $$
    and $g^{-1}_k w_k \to 0$ weakly. Moreover,
   \begin{align*}
       \limsup_{k \to +\infty} ||w_{k}||^2 &\leq  \left( \limsup_{n \to +\infty} ||x_n||^2 \right) - ||\phi||^2 \\
       \limsup_{n \to +\infty} \sup_{g \in G} |\langle g^{-1}x_n, y\rangle| &\leq (1-\varepsilon)^{-1} ||\phi||
   \end{align*}
\end{lemma}

\begin{proof}
Let
\begin{equation*}
    D(x_n) = \sup\{||\phi|| : \exists \phi \in H, \exists (n_k)_{k \in \mathbb N} \subset \mathbb N, (g_k)_{k \in \mathbb N} \subset G, g^{-1}_{k}x_{n_k} \to \phi \text{ weakly } \}
\end{equation*}
If $D(x_n) = 0$, take $\phi_1 = 0$ and $x_n = w_n$. If $D(x_n) > 0$, let $\phi_1$ be chosen near-maximal in the sense that
$$
(1-\varepsilon) D(x_n) < ||\phi_1|| \leq D(x_n).
$$
In particular, for any other weak limit $g^{-1}_{k}x_{n_k} \to \phi$, one has
$$
 ||\phi|| \leq D(x_n) < (1-\varepsilon)^{-1} ||\phi_1||.
$$
No matter the value of $D(x_n)$, it also holds that
$$
\limsup_{n \to +\infty} \sup_{g \in G} |\langle g^{-1}x_n, y\rangle| \leq (1-\varepsilon)^{-1} ||\phi_1||
$$
for any $y \in H, ||y|| = 1$. If that were not the case, there would be $y \in H$, $||y|| = 1$, a subsequence $(n_k)_{k \in \mathbb N} \subset \mathbb N$ and a sequence $(g_k)_{k \in \mathbb N} \subset G$ for which
$$
|\langle g_k^{-1}x_{n_k}, y\rangle| > (1-\varepsilon)^{-1} ||\phi_1||
$$
 Passing to a further subsequence (still denoted with the same indices) which converges weakly to a limit $\phi$ (thanks to Banach-Alaoglu), one would get
 $$
||\phi|| = \sup_{||z|| = 1}|\langle \phi,z \rangle| \geq |\langle \phi, y \rangle| = \lim_{k \to +\infty} |\langle g_k^{-1}x_{n_k}, y\rangle| \geq (1-\varepsilon)^{-1} ||\phi_1|| > D(x_n),
 $$
 which contradicts the definition of $D(x_n)$. 
 
 Consider $(g_{n})_{n \in \mathbb N} \subset G$ and a subsequence of $(x_n)_{n \in \mathbb N}$ (still denoted with the same index $n$) such that $g_{n}^{-1}x_n$ converges weakly to $\phi_1$ by definition. Decompose
$$
x_n = g_{n} \phi_1 + w_{n}
$$
and observe that $g^{-1}_{n}w_{n}$ converges weakly to zero. Moreover,
\begin{align*}
    ||x_n||^2 = ||g_{n} \phi_1 + w_{n}||^2 &= ||g_{n}\phi_1||^2 + ||w_{n}||^2 + 2\Re \langle g_{n} \phi_1,  g_{n}  g_{n}^{-1} w_{n} \rangle  \\
    &=||\phi_1||^2 + ||w_{n}||^2 + 2\Re \langle \phi_1, g_{n}^{-1} w_{n} \rangle
\end{align*}
as unitary operators preserve the inner product (thanks to the Polarization Identity). Hence,
$$
\limsup_{n \to +\infty} ||w_{n}||^2 \leq  \left( \limsup_{n \to +\infty} ||x_n||^2 \right) - ||\phi_1||^2
$$
as $g_{n}^{-1} w_{n}$ converges weakly to zero. 

\end{proof}

\begin{lemma}\label{lemma:g_ig_j unboundedness}
     Let $H$ be a Hilbert space and $G$ be a group of dislocations acting on $H$. Assume that there are $\phi_1, \phi_2 \in H$, $\phi_2 \neq 0$, $(x_n)_{n \in \mathbb N} \subset H$ and $(g_{1,n})_{n \in \mathbb N},  (g_{2,n})_{n \in \mathbb N} \subset G$ for which
     $$
    g^{-1}_{1,n} x_n \to \phi_1, \quad g^{-1}_{2,n} (x_n-g_{1,n}\phi_1) \to \phi_2 \quad \text{weakly}.
     $$
     Then $(g^{-1}_{1,n}g_{2,n})_{n \in \mathbb N}$ escapes any compact set. In particular, $g^{-1}_{2,n} x_n \to \phi_2$ weakly. 
\end{lemma}

\begin{remark}
    The sequence $(g^{-1}_{2,n}g_{1,n})_{n \in \mathbb N}$ escapes any compact set as well, by continuity of the inverse map $g \mapsto g^{-1}$.
\end{remark}

\begin{proof}
    Assume not. Up to passing to a subsequence (still denoted with index $n$), $(g^{-1}_{2,n}g_{1,n})_{n \in \mathbb N}$ converges to some $g_* \in G$. Then
   \begin{align*}
           g^{-1}_{2,n}x_n &= g^{-1}_{2,n} (x_n-g_{1,n}\phi_1)  + g^{-1}_{2,n}g_{1,n}\phi_1 \to \phi_2 + g_*\phi_1 \quad \text{weakly}, \\
    g^{-1}_{1,n}x_n &= (g^{-1}_{1,n}g_{2,n})(g^{-1}_{2,n}x_n) \to g_*^{-1}(\phi_2 + g_*\phi_1) = g_*^{-1}\phi_2 + \phi_1 \quad \text{weakly}, \\
    g^{-1}_{1,n}x_n &\to \phi_1 \quad \text{weakly}.
   \end{align*}
    As $G$ acts unitarily, one must have $\phi_2 = 0$, which is a contradiction. 
    
    Finally, since $G$ is a dislocation group and $(g^{-1}_{2,n}g_{1,n})_{n \in \mathbb N}$ escapes any compact set, $g^{-1}_{2,n}g_{1,n} \phi_1 \to 0$ weakly by definition of $G$.
\end{proof}

\begin{corollary}\label{cor:g_ig_j unboundedness}
     Let $H$ be a Hilbert space and $G$ be a group of dislocations acting on $H$. Assume that there are $\phi_1, ..., \phi_M \in H$, $\phi_j \neq 0$ when $j \neq 1$, $(x_n)_{n \in \mathbb N} \subset H$ and $(g_{1,n})_{n \in \mathbb N}, ..., (g_{M,n})_{n \in \mathbb N} \subset G$, $M \geq 2$, such that for all $L = 1, ..., M$:
     $$
     g^{-1}_{L,n} \left(x_n- \sum_{m=1}^{L-1} g_{m,n}\phi_m \right) \to \phi_L \quad \text{weakly}.
     $$
     Then for any $1 \leq i \neq j \leq M$, $(g^{-1}_{i,n}g_{j,n})_{n \in \mathbb N}$ escapes any compact set. In particular, $g^{-1}_{L,n} x_n \to \phi_L$ weakly for any $L = 1, ..., M$.
\end{corollary}

\begin{proof}
    This is a proof by induction. The base case $M = 2$ comes from Lemma \ref{lemma:g_ig_j unboundedness}. Let us prove the $M+1 > 2$ case, while assuming that the Lemma holds with a fewer number of elements and sequences. Define
    $$
    y_n = x_n- \sum_{m=1}^{M-1}g_{m,n}\phi_m, \quad g^{-1}_{M,n}y_n \to \phi_M \quad \text{weakly}.
    $$
    Observe that
    $$
g^{-1}_{M+1,n}(y_n  - g_{M}\phi_M)= g^{-1}_{M+1,n}\left( x_n- \sum_{m=1}^{M}g_{m,n}\phi_m\right) \to \phi_{M+1} \quad \text{weakly}
    $$
    and apply Lemma \ref{lemma:g_ig_j unboundedness}. We obtain that $(g^{-1}_{M,n}g_{M+1,n})_{n \in \mathbb N}$ and $(g^{-1}_{M+1,n}g_{M,n})_{n \in \mathbb N}$ also escape any compact set and 
    $$
g^{-1}_{M+1,n}y_n = g^{-1}_{M+1,n} \left( x_n- \sum_{m=1}^{M-1}g_{m,n}\phi_m \right) \to \phi_{M+1} \quad \text{weakly}.
    $$
    Hence, we can apply Lemma \ref{lemma:g_ig_j unboundedness} once again to $\phi_1,...,\phi_{M-1},\phi_{M+1}$, $(x_n)_{n \in \mathbb N}$ and $(g_{1,n})_{n \in \mathbb N}$, ..., $(g_{M-1,n})_{n \in \mathbb N}$, $(g_{M+1,n})_{n \in \mathbb N} \subset G$. Here, $\phi_{M+1}$ and $g_{M+1,n}$ replaces $\phi_{M}$ and $g_{M,n}$ when applying Lemma \ref{lemma:g_ig_j unboundedness}. This finishes the proof.
\end{proof}
Now, we are ready to prove the profile decomposition.

\begin{proof}[Profile Decomposition]
We iterate Lemma \ref{lemma:profile decomposition, first profile}. First, we find a subset $I_1 \subset \mathbb N$ and a decomposition
$$
x_{n} = g_{1,n} \phi_1 + w_{1,n}, \quad n \in I_1
$$
where $\lim \limits_{n \in I_1, n \to +\infty} g^{-1}_{1,n} w_{1,n} = 0$ weakly and
   \begin{align*}
       \limsup_{n \to +\infty} ||w_{1,n}||^2 &\leq  \left( \limsup_{n \to +\infty} ||x_n||^2 \right) - ||\phi_1||^2 \\
       \limsup_{n \to +\infty} \sup_{g \in G} |\langle g^{-1}x_{n}, y\rangle| &\leq 2||\phi_1||, \quad ||y|| = 1.
   \end{align*}
Applying the Lemma to $w_{1,n}$ and then to each new remainder, one finds a sequence of decreasing unbounded subsets $\mathbb N \supset I_1 \supset I_2 \supset ...$, a sequence of profiles $(\phi_m)_{m \in \mathbb N} \subset H$, group elements $g_{m,n} \subset G$, $m \in \mathbb N$, $n \in I_m$, and remainders $w_{m,n} \subset H$, $m \in \mathbb N$, $n \in I_m$, for which
$$
x_{n} = \sum_{m=1}^M g_{m,n} \phi_m + w_{M,n}, \quad M \in \mathbb N, n \in I_M.
$$
Moreover, for any $M \in \mathbb N_{\geq 1}$, $\lim \limits_{n \in I_m, n \to +\infty} g^{-1}_{M,n}w_{M,n} = 0$ weakly,
$$
       0 \leq \limsup_{n \to +\infty} ||w_{M,n}||^2 \leq  \left( \limsup_{k \to +\infty} ||w_{M-1,n}||^2 \right) - ||\phi_M||^2 \leq ... \leq \left( \limsup_{n \to +\infty} ||x_{n}||^2 \right) - \sum_{m=1}^M||\phi_m||^2
$$
and
\begin{equation}
    \sup_{||y|| = 1} \limsup_{n \to +\infty} \sup_{g \in G} |\langle g^{-1}w_{M-1,n}, y\rangle| \leq 2 ||\phi_M||  \label{eq:smallness of remainder}
\end{equation}
using the convention that $w_{0,n} = x_n$, $w_{m,n} = 0$ if $n \notin I_m$. Observe that the decomposition is finite (i.e., $\phi_{M+1} = \phi_{M+2} = ... = 0$) if one reaches $M$ for which $D[(w_{M,n})_{n \in \mathbb N}] = 0$, where $D[\cdot]$ is defined as in (\ref{eq:D(x_n) def}). Otherwise, the decomposition is infinite and the series
$$
\sum_{m=1}^{\infty}||\phi_m||^2 < +\infty
$$
is convergent. In both cases, we deduce from Lemma \ref{lemma:profile decomposition, first profile} (if the decomposition is finite) and (\ref{eq:smallness of remainder}) (if the decomposition is infinite) that
$$
\lim_{M\to +\infty}\sup_{||y|| = 1}  \limsup_{n \to +\infty} \sup_{g \in G} |\langle g^{-1}w_{M-1,n}, y\rangle| = 0.
$$
To finish the proof, we pass to a subsequence $(n_k)_{k \in \mathbb N}$ for which $n_j \in I_k$ for all $j \geq k$. The decomposition rewrites as follows
\begin{equation}
    x_{n_k} = \sum_{m=1}^M g_{m,n_k} \phi_m + w_{M,n_k}, \quad M \in \mathbb N, k \geq M. \label{eq:final decomposition}
\end{equation}
Since $\lim \limits_{k \to +\infty} g^{-1}_{M,n_k}w_{M,n_k} = 0$ weakly for any $M \in \mathbb N_{\geq 1}$, Corollary \ref{cor:g_ig_j unboundedness} implies the separation of scales and $g^{-1}_{M,n_k} x_{n_k} \to \phi_{M}$ weakly as $k \to +\infty$ (as long as the profiles $\phi_2, ..., \phi_M$ are non-zero). It follows directly from those weak convergences and the decomposition (\ref{eq:final decomposition}) that for $1 \leq m < M$, $g^{-1}_{m,n_k}w_{M,n_k} \to 0$ weakly as well. Finally, for $M \in \mathbb N$ fixed, one has the Pythagorean decomposition of the energy:
\begin{align*}
    ||x_{n_k}||^2 &= \langle \sum_{m=1}^M g_{m,n_k} \phi_m + w_{M,n_k}, \sum_{m=1}^M g_{m,n_k} \phi_m + w_{M,n_k} \rangle, \quad k \geq M \\
    &= \sum_{m=1}^M\sum_{l=1}^M \langle g_{l,n_k}^{-1}g_{m,n_k}\phi_m,\phi_l \rangle + ||w_{M,n_k}||^2 + 2\sum_{m=1}^M\Re \langle \phi_m,  g_{m,n_k}^{-1}w_{M,n_k} \rangle.
\end{align*}
Taking $k \to +\infty$, one obtains
$$
\limsup_{k \to +\infty} ||x_{n_k}||^2 = \sum_{m=1}^{M} ||\phi_m||^2 + \limsup_{k \to +\infty} ||w_{M,n_k}||^2  \quad M \in \mathbb N,
$$
which finishes the proof.
\end{proof}

As a corollary, we show that one can always choose a decomposition for which the first extracted profile is the weak limit of the subsequence. 

\begin{corollary}[Profile Decomposition with explicit weak limit]\label{thm:abstract profile decomposition with weak limit first}
    Let $H$ be a Hilbert space and $G$ be a group of dislocations acting on $H$. Let $(x_n)_{n \in \mathbb N} \subset H$ be a bounded sequence. After passing to a subsequence (still denoted with index $n$) of $(x_n)_{n \in \mathbb N}$ converging weakly to $x_* \in H$, there exists profiles $(\phi_{m})_{m \in \mathbb N} \subset H$, group elements $((g_{m,n})_{m \in \mathbb N})_{n \geq m} \subset G$, remainders  $((w_{m,n})_{m \in \mathbb N})_{n \geq m} \subset H$ for which
$$
x_{n} = x_* + \sum_{m=2}^M g_{m,n} \phi_m + w_{M,n}, \quad M \in \mathbb N_{\geq 1}, n \geq M.
$$
Moreover, one has:
\begin{enumerate}
    \item \textbf{Pythagorean decomposition of energy:}
    $$
    \limsup_{n \to +\infty} ||x_{n}||^2  =  ||x_*||^2 + \sum_{m=2}^M||\phi_m||^2 + \limsup_{n \to +\infty} ||w_{M,n}||^2, \quad M \in \mathbb N_{\geq 1}.
    $$
    \item \textbf{Finiteness of decomposition and non-zero profiles:} Let $w_{1,n} = x_n-x_*$ and
    \begin{equation}
    D[(x_n)] := \sup\{||\phi|| : \exists \phi \in H, \exists (n_k)_{k \in \mathbb N} \subset \mathbb N, (g_k)_{k \in \mathbb N} \subset G, g^{-1}_{k}x_{n_k} \to \phi \text{ weakly } \}.
\end{equation}
    For $M \in \mathbb N_{\geq 1}$, the extracted profile $\phi_{M+1}$ is non-zero if and only if $D[(w_{M,n})] > 0$. Moreover, if $D[(w_{M,n})] = 0$, then all the subsequent profiles $\phi_{M+1} = \phi_{M+2} = ... = 0$ are zero and the decomposition is finite. 
    \item \textbf{Separation of scales:} Let $g_{1,n} = \mathrm{Id}_G$. Fix $i \neq j$ and assume either $\phi_i, \phi_j$ are non-zero profiles or $i = 1$ and $\phi_j$ is a non-zero profile. Then, $(g^{-1}_{i,n} g_{j,n})_{n \in \mathbb N}$ escapes any compact set.
    \item \textbf{Weak convergence to profile:} For any $M \in \mathbb N_{\geq 1}$, if either $M = 1$ or $\phi_M$ is a non-zero profile, then $g^{-1}_{M,n}x_n \to \phi_M$ weakly. In particular, $x_n \to \phi_1$ weakly.
    \item \textbf{Smallness of remainder:} For $1 \leq m \leq M$, $g^{-1}_{m,n}w_{M,n} \to 0$ weakly. Moreover, the remainders also converge weakly to zero in the following sense:
    $$
    \lim_{M\to +\infty}\sup_{||y|| = 1}  \limsup_{n \to +\infty} \sup_{g \in G} |\langle g^{-1}w_{M,n}, y\rangle| = 0.
    $$
\end{enumerate}
\end{corollary}

\begin{proof}
    Let $x_*$ be a weak limit for $x_n$ after passing to a first subsequence. Apply Theorem \ref{thm:abstract profile decomposition}  using the trivial dislocation group. Hence,
    $$
    x_n = x_* + (x_n - x_*)
    $$
    with energy decoupling. Apply Theorem \ref{thm:abstract profile decomposition} to $x_n -x_*$:
  \begin{equation}
      x_{n} -x_* = \sum_{m=2}^M g_{m,n} \phi_m + w_{M,n}, \quad M \in \mathbb N_{\geq 1}, n \geq M, \label{eq: profile decomposition of x_n-x_*}
  \end{equation}
where we order the profiles starting at $m = 2$ and define $g_{1,n} = \mathrm{Id}_G$.  The energy decoupling follows easily.

Next, we know from the theorem that $(g_{i,n}^{-1}g_{j,n})_{n \in \mathbb N}$ is not contained in any compact set for $i \neq j$, $i,j \geq 2$, if $\phi_i, \phi_j$ are non-zero profiles. Moreover, $g_{j,n}^{-1}(x_n-x_*) \to \phi_j$ weakly if $j = 2$ or $\phi_j$ is a non-zero profile. We first show that $(g_{j,n})_{n \in \mathbb N} = (g_{1,n}^{-1}g_{j,n})_{n \in \mathbb N}$ escapes any compact set for $j \neq 1$ if $\phi_j$ is a non-zero profile. Assume not. After passing to a further subsequence, one can assume that $g_{j,n}$ converges to a group element $g_{j,*} \in G$ in the group topology. As the group action is jointly continuous and $G$ acts unitarily,
$$
\langle g_{j,n}^{-1}(x_n-x_*), \phi \rangle_{H} -\langle g_{j,*}^{-1}(x_n-x_*), \phi \rangle_{H} = \langle (x_n-x_*), (g_{j,n}-g_{j,*})\phi \rangle_{H} \to 0 \quad \forall \phi \in H 
$$
by Cauchy-Schwarz, while 
$$
\langle g_{j,*}^{-1}(x_n-x_*), \phi \rangle_{H} = \langle (x_n-x_*), g_{j,*} \phi \rangle_{H} \to 0 \quad \forall \phi \in H
$$
by weak convergence of $x_n \to x_*$. In other words, $g_{j,n}^{-1}(x_n-x_*) \to 0$ weakly. Yet, by construction of the profiles, $g_{j,n}^{-1}(x_n-x_*) \to \phi_j$ weakly, contradicting the fact that $\phi_j$ is  a non-zero profile.

Furthermore, we show that $g_{j,n}^{-1} x_n \to \phi_j$ weakly if $\phi_j$ is a non-zero profile. We know that $g_{j,n}^{-1} (x_n-x_*) \to \phi_j$ weakly by construction and $g_{j,n}^{-1}$ escapes any compact set. By definition of a dislocation group, $g_{j,n}^{-1}x_* \to 0$ weakly, which proves the claim.

Finally, we need to show that $g_{1,n}^{-1}w_{M,n} = w_{M,n} \to 0$ weakly. Looking at the decomposition (\ref{eq: profile decomposition of x_n-x_*}), it is sufficient to prove that $g_{j,n}\phi_j \to 0$ if $\phi_j$ is a non-zero profile. This follows from the definition of a dislocation group, since $(g_{j,n})_{n \in \mathbb N}$ escapes any compact set. This finishes the proof. 
\end{proof}

If one has a sufficiently nice embedding of $H$ into some other normed space $Y$, the weak convergence of the remainder to zero can be strengthened to a convergence statement in $Y$, which is more useful for applications.

\begin{definition}[$G$-weak convergence]
    Let $H$ be a Hilbert space and $G$ be a group acting unitarily on $H$. A sequence $(x_n)_{n \in \mathbb N} \subset H$ is said to converge $G$-weakly to a limit $x \in H$ if
    $$
\lim_{n \to +\infty} \sup_{g \in G}|\langle y, g^{-1}(x_n-x) \rangle| = 0, \quad \forall y \in H.
    $$
\end{definition}

\begin{remark}
    In particular, such a sequence converges weakly to $x$ in the usual sense.
\end{remark}

\begin{definition}[$G$-complete continuity]
    Let $H$ be a Hilbert space, $G$ be a group of dislocations acting on $H$, $Y$ be a normed space and $T : H \rightarrow Y$ be a continuous linear operator. The operator $T$ is said to be $G$-completely continuous if for any sequence $(x_n)_{n \in \mathbb N} \subset H$ which converges $G$-weakly to a limit $x \in H$, the sequence $(Tx_n)_{n \in \mathbb N} \subset Y$ converges strongly to $Tx \in Y$.
\end{definition}

\begin{proposition}
Under the same assumptions as Theorem \ref{thm:abstract profile decomposition}, if $H$ is separable, $Y$ is a normed space and $H \hookrightarrow Y$ is a $G$-completely continuous embedding, then
$$
\lim_{M \to +\infty} \limsup_{n \to \infty}||w_{M,n}||_{Y} = 0.
$$
\end{proposition}

\begin{remark}
    This is Theorem 2.6 in \cite{okumura2022profile}. However, the original proof contains a technical gap: the constructed subsequences depend on a fixed element $\phi \in H^* = H$. We provide a corrected argument under the additional assumption that $H$ is separable.
\end{remark}

\begin{proof}
  Assume not. Then
    $$
\limsup_{M \to +\infty} \limsup_{n \to \infty}||w_{M,n}||_{Y} \geq \varepsilon > 0.
    $$
    Let $(x_k)_{k \in \mathbb N} \subset B_1^{H}(0)$ be a countable and dense subset of the unit ball of $H$. Consider a subsequence $M_k \to +\infty$ along which
\begin{align*}
   \limsup_{n \to \infty}||w_{M_k,n}||_{Y} &\geq \varepsilon/2, \quad k \in \mathbb N \\
\sup_{||y|| = 1}  \limsup_{n \to +\infty} \sup_{g \in G} |\langle g^{-1}w_{M_k,n}, y\rangle| &\leq 2^{-k}, \quad k \in \mathbb N,
\end{align*}
    using the smallness of the remainder obtained in Theorem \ref{thm:abstract profile decomposition}. In particular, for each $y \in B_1^{H}(0)$, one has
    $$
 \limsup_{n \to +\infty} \sup_{g \in G} |\langle g^{-1}w_{M_k,n}, y\rangle| \leq 2^{-k}||y|| \leq 2^{-k}, \quad k \in \mathbb N.
    $$
    Therefore, there exists a subsequence $n_k \to +\infty$ along which
\begin{align*}
    ||w_{M_k,n_k}||_{Y} &\geq \varepsilon/4, \quad k \in \mathbb N \\
 \sup_{g \in G} |\langle g^{-1}w_{M_k,n_k}, x_l\rangle| &\leq 2^{2-k}, \quad k \in \mathbb N, 1 \leq l \leq k.
\end{align*}
In particular, 
$$
\lim \limits_{k \to +\infty}  \sup_{g \in G} |\langle g^{-1}w_{M_k,n_k}, x_l\rangle| = 0 \quad \forall l \in \mathbb N.
$$
If $x \in B_1^{H}(0)$ is not an element from the countable dense subset and $|x-x_l| \leq \delta$ for some $l \in \mathbb N$, then
\begin{align*}
    \sup_{g \in G} |\langle g^{-1}w_{M_k,n_k}, x\rangle|  &\leq  \sup_{g \in G} |\langle g^{-1}w_{M_k,n_k}, x-x_l\rangle| +  \sup_{g \in G} |\langle g^{-1}w_{M_k,n_k}, x_l\rangle| \\
    &\leq \delta ||w_{M_k,n_k}||_{H} +  \sup_{g \in G} |\langle g^{-1}w_{M_k,n_k}, x_l\rangle| \\
    &\leq \delta \left(\limsup_{n \to +\infty} ||x_n||_{H} \right) +  \sup_{g \in G} |\langle g^{-1}w_{M_k,n_k}, x_l\rangle|.
\end{align*}
Taking $k \to +\infty$ and then $\delta \to 0$ implies that 
$$
\lim \limits_{k \to +\infty}  \sup_{g \in G} |\langle g^{-1}w_{M_k,n_k}, x\rangle| = 0 \quad \forall x \in B_1^{H}(0)
$$
and this actually holds no matter the norm of $x$ by rescaling. Hence, we have proved that $(w_{M_k,n_k})_{k \in \mathbb N} \subset H$ converges $G$-weakly to zero. By $G$-complete continuity of the embedding $H \hookrightarrow Y$, $(w_{M_k,n_k})_{k \in \mathbb N} \subset Y$  converges strongly to zero in $Y$. In other words,
$$
\lim_{k \to +\infty} ||w_{M_k,n_k}||_{Y} = 0,
$$
which is a contradiction.
\end{proof}

\chapter{Construction of Type II blow-up solutions for the Nonlinear Wave Equation}
~\newline~

The content of this chapter is taken from my article \cite{samuelian2026construction}, published in \textit{Annales Henri Poincaré}, with no significant change.

\section{Introduction}

We study the focusing, energy-critical, nonlinear wave equation
\begin{equation}
    \partial_{tt}u - \Delta u - |u|^{p-1}u = 0, \quad t \geq 0, \ x \in \mathbb R^d, \ d \geq 3, \ p = (d+2)/(d-2) \quad \text{(NLW)}.
\label{wave eq}
\index{p@$p = (d+2)/(d-2)$, the critical exponent}
\end{equation}
The equation is called energy-critical because it is invariant under the scaling $u_{\lambda}(x,t) = \lambda^{(d-2)/2}u(\lambda x, \lambda t)$. Moreover, the $\dot{H}^1 \times L^2$ norm of the rescaled data $(u_{\lambda}, \partial_t u_{\lambda})$ at $t = 0$ is independent of $\lambda$. The equation is locally well-posed (in the sense that a solution given by the Duhamel formula exists) for initial data in $\dot{H}^1 \times L^2$. If a solution fails to be global, then its $L_{t,x}^{\frac{2(d+1)}{(d-2)}}$-norm blows up (\cite{Kenig_Merle_local_wellposed_low_dim}, \cite{Critical_local}). Changing the sign of the nonlinearity leads to the defocusing case, which is known to be globally well-posed (\cite{Critical_defocusing_global}).

A classical example of blow-up is due to Levine (\cite[Section 12.5.1]{evans10}). Utilizing the conserved energy functional
\begin{equation}
E(u(t),\partial_tu(t)) = \int_{\mathbb R^d} \underbrace{\frac{1}{2} |\nabla_{t,x}u(t)|^2}_{kinetic \ en.} \underbrace{- \frac{1}{p+1}|u(t)|^{p+1}}_{potential \ en.}dx
\index{E(u(t),\partial_tu(t))@$E(u(t),\partial_tu(t))$, the conserved energy}
\end{equation}
of the solution, Levine showed that if the initial data is smooth and compactly supported with negative energy $E(u(0),\partial_tu(0)) < 0$, then the solution cannot exist globally in time. Denoting by $T_{\max} \in [0,+\infty]$ the maximal forward time of existence of the solution, $$(u(t),\partial_tu(t)) \in C([0,T_{\max}), \dot{H}^1 \times L^2),$$ 
 two distinct types of solutions are distinguished:
\begin{enumerate}
    \item Type I: $||(u(t),\partial_tu(t))||_{L^{\infty}([0,T_{\max}), \dot{H}^1 \times L^2)} = +\infty,$
    \item Type II: $||(u(t),\partial_tu(t))||_{L^{\infty}([0,T_{\max}), \dot{H}^1 \times L^2)} < +\infty.$
\end{enumerate}
Examples of Type I blow-up solutions can be obtained by considering the solution
$$
u_T(x,t) = c_p(T-t)^{-\frac{2}{p-1}}, \quad c_p = \left(\frac{2(p+1)}{(p-1)^2}\right)^{\frac{1}{p-1}}
$$
to the wave equation. Selecting initial data $(u_T(x,0), \partial_t u_T(x,0))$ with an appropriate spatial cutoff produces such a blow-up at a finite time $T_{\max} \in (0,T]$ (\cite[Section 6]{Subcritical_local}).

In the following, we are interested in constructing radially symmetric Type II blow-ups. The study of Type II radial solutions is closely related to the stationary radial solution
\begin{equation}
    W(x) = \left( 1 + \frac{|x|^2}{d(d-2)} \right)^{-\frac{d-2}{2}}
    \index{W(x)@$W(x)$, the ground state}
\end{equation}
called the ground state (or ``bubble'' of energy), which is also an extremizer for the Sobolev embedding $\dot{H}^1 \hookrightarrow L^{p + 1}$. Roughly speaking, the soliton resolution conjecture (recently solved in \cite{soliton_odd}, \cite{soliton_jendrej}) asserts that any radial Type II solution behaves asymptotically as a superposition of dynamically scaled bubbles plus a radiation term, which is given by a free wave if the solution is global and a stationary element in $\dot{H}^1$ otherwise. 

The goal of this paper is to construct explicit solutions $u(x,t)$ of (\ref{wave eq}) in dimension $d \in \{4,5\}$, exhibiting finite-time Type II blow-up precisely at the origin $x = t = 0$, with a prescribed polynomial blow-up rate $t^{-1-\nu}$, where $\nu > 1$ for $d = 4$ and $\nu > 3$ for $d = 5$. Such solutions were previously constructed by Krieger-Schlag-Tataru for $d = 3$ (\cite{Krieger_2012}) with $\nu > 0$ and by Jendrej for $d = 5$ (\cite{Jendrej_blowup5d}) with $\nu \in \{3\} \cup (8,+\infty)$.

As in \cite{Krieger_2012}, where they prove that the previously known range of exponents $\nu > 1/2$ for $d = 3$ can be relaxed to $\nu > 0$, it is likely that our assumption $\nu > 1$ in dimension $4$ can be relaxed to $\nu > 0$. This restriction is a technicality arising from a nonlinear Sobolev estimate. In Proposition \ref{prop:local lipschitz}, we exploit the embedding of the algebra $H^{1 + \frac{(6-d)\nu}{2}}(\mathbb R^d)$ into $C^0_b(\mathbb R^d)$ to establish a local Lipschitz property on a nonlinear operator, which is essential to apply the Banach Fixed Point theorem. A similar limitation arises in the work of Jendrej \cite[Lemma 4.6, Proposition 4.7]{Jendrej_blowup5d}, where it results in the condition $\nu > 8$ for $d = 5$. Our analysis improves the permissible range to all $\nu > 3$ in the fifth dimension. However, in contrast to the fourth dimension, the restriction $\nu > 3$ is more than technical: while it also occurs in the nonlinear estimates from Proposition \ref{prop:local lipschitz}, the restriction is crucial in ensuring the positivity of the approximation $u_2$ (see (\ref{H(z) definition}) and (\ref{v_2 dominant component})), which is needed to handle the absolute value in the nonlinearity $F(x) = |x|^{p-1}x$ without losing regularity. 
    
To our knowledge, the $d = 4$ case with prescribed polynomial blow-up rate has not been previously addressed, but existence is already known (\cite{HillairetRaphael2012}). We also remark that in dimension $d = 6$, infinite-time superposition of two bubbles have been constructed (\cite{soliton_jendrej}), but finite-time blow-up has not been exhibited yet. Although the soliton resolution conjecture is now proven, explicit constructions of finite- or infinite-time solutions exhibiting a dynamically scaled bubble profile have only been achieved in low dimensions $3 \leq d \leq 6$. Whether such constructions can be extended to higher dimensions, where the nonlinearity lacks twice differentiability, remains an open question. 

The main result of this paper is the following theorem:

\begin{notation}
    We write $u(x) \in H^{s-}(\mathbb R^d)$ or $|u(x)| \leq |1-x|^{s-}$ if the property holds with exponent $s-\delta$ for all sufficiently small $\delta > 0$ instead of $s$. A similar meaning applies to expressions such as $u(x) \in H^{s+}(\mathbb R^d)$ or $|u(x)| \leq |1-x|^{s+}$. 
    \index{Hs-@$H^{s\pm}(\mathbb R^d)$}
    \index{u(x) < 1-x@$\lvert u(x) \rvert \leq \lvert 1-x \rvert^{s\pm}$}
\end{notation}

\begin{theorem}\label{thm:blow-up, main thm}
    Let $d \in \{4,5\}$, $\nu > \nu_0(d)$ where $\nu_0(4) = 1$ and $\nu_0(5) = 3$, $\delta > 0$, $N_0 \gg 1+\nu$ be fixed. There exists a radial solution $u(x,t)$ of (\ref{wave eq}) on $\mathbb R^d \times [0,t_0] $, $t_0 \ll 1$, which has the form:
\begin{enumerate}
    \item $u(x,t) = \lambda(t)^{\frac{d-2}{2}}W(\lambda(t)x) + \eta(x,t)$ inside the backward light cone $$C = \{(x,t): 0 \leq |x| \leq t, 0 < t \leq t_0\}.$$ 
    Moreover, if $d = 5$, $u(x,t) > 0$ on $C$.
    \item $\lambda(t) = t^{-1-\nu}$ and the solution blows up at $r = t = 0$.
    \item $\eta(x,t)$ can be decomposed as
$\eta(x,t) = u^e(x,t) + \varepsilon(x,t)$, where $u^e \in C^{\frac{1}{2}+\frac{6-d}{2}\nu-}(C)$ and
\begin{align*}
\sup_{0 < t < t_0} t^{-\frac{6-d}{2}\nu-1}||u^e||_{H^{1+\frac{6-d}{2}\nu-}(\mathbb R^d)} +  t^{-\frac{6-d}{2}\nu} ||\partial_t u^e||_{H^{\frac{6-d}{2}\nu-}(\mathbb R^d)} &< +\infty \\
\sup_{0 < t < t_0} t^{-N_0}||\varepsilon||_{H^{1+\frac{6-d}{2}\nu-}(\mathbb R^d)} + t^{-N_0+1}||\partial_t \varepsilon||_{H^{\frac{6-d}{2}\nu-}(\mathbb R^d)} &< +\infty
\end{align*}
    \item The local energy 
$$E_{loc} = \int_{|x| < t} \eta_t^2 + |\nabla \eta|^2  + |\eta|^{p+1} dx$$
of $\eta(x,t)$ in the light cone $|x| \leq t$ vanishes as $t \to 0$
    \item Outside the light cone, the energy of the solution can be controlled
$$\int_{|x| \geq t} u_t^2 + |\nabla u|^2  + |u|^{p+1} dx \leq \delta$$
for all sufficiently small $t  > 0$.
\end{enumerate}
\index{N0@$N_0$, the smallness of the approximation error $e_k$ for large $k$, as well as the $\varepsilon$ part from Theorem \ref{thm:blow-up, main thm}}
\end{theorem}

Our proof proceeds in two main steps. In sections \ref{section:renormalization step, basic idea} to \ref{section:renormalization step, next iterates}, starting from the bubble $u_0 = \lambda^{\frac{d-2}{2}}W(\lambda x)$, we linearize and simplify (NLW) to iteratively construct a sequence of approximate solutions $u_k = u_0 + v_1 + ... + v_k$ to (NLW) on a cone $0 < r < t < t_0 \ll 1$. The correction terms $v_k$ are smooth on $0 \leq r < t < t_0$, except for a logarithmic-power singularity of the form $(1-a)^{\frac{1}{2}+\frac{1}{2}\nu}\log(1-a)^j$, $a = r/t$ at the boundary $r = t$. Moreover, the pointwise error $|F(u_k)-\square u_k|$ decreases at each iteration. The term $u^e(x,t)$ from Theorem \ref{thm:blow-up, main thm} is precisely the difference $u_k(x,t)-u_0(x,t)$ for sufficiently large $k$, extended from the cone to all of $\mathbb R^d$ while keeping the same size and regularity.

 In sections \ref{section:spectral theory of the linearized operator} to \ref{section:exact solution fourier method}, we find an exact solution $\lambda^{\frac{d-2}{2}}W + u^e + \varepsilon$ within cone by solving a fixed-point problem in a generalized Fourier space $L^2(\mathbb R, d\rho(\xi))$. This space arises naturally when analyzing the spectral properties of the perturbed Schrödinger operator $$
 \mathcal{L} = -\partial_{RR} - pW(R)^{p-1} + \frac{1}{R^2} \cdot \left( \frac{(d-3)(d-1)}{4} \right), \quad R = \lambda |x|, 
 $$ which emerges from the linearization of (NLW) around the ground state $u_0$. Finally, we show in section \ref{section:end of proof} using finite propagation of speed and well-posedness theory that this solution extends outside the cone.

The principal difference between the cases $d = 4$ and $d = 5$ in this paper lies precisely in this renormalization procedure. For $d = 4$, discussed separately in Section \ref{section:appendix, dimension 4}, our method closely follows the approach of Krieger-Schlag-Tataru developed originally the 3-dimensional case. Their scheme has been successfully applied to various other equations, including the 3D-critical NLS, Schrödinger maps, wave maps (see, e.g.,  \cite{ortoleva2012nondispersive}).  In all of these cases, the polynomial structure of the nonlinearity permits the use of simple algebraic manipulations to control error terms. By contrast, the case $d = 5$ requires a modified setup due to the lower regularity of the nonlinearity, which complicates the treatment of nonlinear error terms. We address this difficulty by carefully distinguishing three distinct spatial regions:
$$
R \lesssim (t\lambda)^{\frac{2}{3}}, \quad (t\lambda)^{\frac{2}{3}} \lesssim R \lesssim (t\lambda)^{\frac{2}{3}+\varepsilon}, \quad R \gtrsim (t\lambda)^{\frac{2}{3}+\varepsilon}
$$
using cutoff functions. This allows for the use of convergent multinomial expansions on each region to treat the nonlinear errors. Yet, the renormalization step encounters serious challenges in higher dimensions $d \geq 6$, the main difficulty being that solving equation (\ref{even v}) introduces singularities at the tip of the cone $r = t$ unless $\nu > 0$ is very small (Remark \ref{rmk:loss of regularity high dimension}). In that scenario, the resulting solutions exhibit low regularity, which would necessitate modifying the Banach spaces employed in the fixed-point argument. Additionally, positivity of the approximations $u_k$ must be carefully verified to remove the absolute value in the nonlinearity during the approximation step, which might further restrict the range of admissible values for $\nu$. Indeed, this explicitly occurs in dimension $d = 5$ where one must impose the condition $\nu > 3$ to ensure positivity. Finally, for $d \geq 6$, the generalized eigenfunction $\phi(R,\xi)$, $R,\xi \geq 0$, of the perturbed Schrödinger operator has a singularity at $\xi \sim 0$, $R^2\xi \sim  1$, which gets worse as the dimension increases (Proposition \ref{prop:spectral density properties}, Corollary \ref{phi upper bound R^2z > 1}). Such singularities could yield less favorable estimates in the transference identity (Section \ref{section:transferrance identity}), when passing back and forth from physical to generalized Fourier space.







\section{Renormalization Step: Basic idea} \label{section:renormalization step, basic idea}
This section outlines the strategy for constructing an approximate solution to the nonlinear wave equation (\ref{wave eq}) within the cone $0 \leq |x| \leq t$ for small times $0 < t \leq t_0$. The core of our method is an iterative process starting from the scaled ground state $u_0$, and successively adding correction terms $v_k$. At each step, we produce an improved approximation $u_{k+1} = u_k + v_k$ where the corresponding approximation error $e_k = F(u_k) - \square u_k$, $\square = \partial_{tt} - \Delta$, $F(x) = |x|^{p-1}x$,  has decreased. \index{square@$\square = \partial_{tt} - \Delta$, d'Alembert operator} \index{F(x)@$F(x) = \lvert x \rvert^{p-1}x$, the nonlinearity} These corrections are determined by solving a pair of linearized equations: an elliptic-type equation to cancel the dominant error near the origin, and a wave-type equation to cancel the dominant error near the boundary of the backward light cone. 

Let $R = \lambda(t)r$, $u_0(R) = \lambda^{(d-2)/2}(t)W(R)$, $d \in \{4,5\}$.
\index{R@$R = \lambda(t)\lvert x \rvert$, a radial variable}
\index{lambda t@$\lambda(t) = t^{-1-\nu}$, polynomial rate of blow-up}
\index{u0@$u_0(R) = \lambda^{(d-2)/2}(t)W(R)$} From the current approximation $u_k$, we set $u_{k+1} = u_k + v_k$ where $v_k$ is a correction term and at each step, we compute the error $e_k = F(u_k) - \square u_k$. If $u$ were an exact solution to $\square u = F(u)$, then the difference $\varepsilon = u - u_{k-1}$ would satisfy
\begin{align*}
    -\square \varepsilon = - \square u + \square u_{k-1} = -F(u) - e_{k-1} + F(u_{k-1}),
\end{align*}
i.e.,
\begin{align} \label{space-variable equation epsilon}
    -\square \varepsilon = -F(u_{k-1}+\varepsilon) - e_{k-1} + F(u_{k-1}).
\end{align}
Linearizing around $\varepsilon = 0$, we obtain the approximation
\begin{equation*}
    -\square \varepsilon + F'(u_{k-1}) \varepsilon + e_{k-1} \approx 0.
\end{equation*}
Further approximating $u_{k-1}$ by $u_0$, we simplify this to
\begin{equation}
    -\square \varepsilon + pu_0^{p-1} \varepsilon + e_{k-1} \approx 0.
\end{equation}
Thus, the correction terms $v_k$ are constructed, roughly, as follows:
\begin{align}
    \Delta v_{1} + pu_0^{p-1} v_{1} + e_{0} &= 0, \quad e_0 = u_0^p - \square u_0
    \label{first v} \\
    - \square v_{2k} + e^0_{2k-1} &= 0, \quad k \geq 1
    \label{even v}
\end{align}
and
\begin{equation}
    \Delta v_{2k+1} + pu_0^{p-1} v_{2k+1} + e^0_{2k} = 0, \quad k \geq 1
    \label{odd v}
\end{equation}
in radial coordinates with zero Cauchy data at the origin. From this point onwards, we focus exclusively on the case $d = 5$ and readers are referred to Section \ref{section:appendix, dimension 4} for $d = 4$. We note that equation (\ref{first v}), which can already be solved explicitly using $e_0$ (see Section \ref{section:renormalization step, first iterate}), is treated separately from equations (\ref{even v}) and (\ref{odd v}), for which a careful analysis of the nonlinearity is necessary to isolate a suitable forcing term $e^0_k$. The overall strategy is the same in both dimensions $d \in \{4,5\}$, but the definition of the forcing terms $e_k^0$ differs slightly since no cutoff is used in dimension $d = 4$.

Equation (\ref{even v}) will be solved in the self-similar variable $a = r/t$, $a \in (0,1)$. This allows us to improve the approximation error near the tip of the cone. The forcing term $e^0_{2k-1}$ for (\ref{even v}) is extracted from $e_{2k-1}$ by keeping only the non-negligible component near the tip of the cone. The remainder 
$$
t^2e^1_{2k-1}:= t^2[e_{2k-1}-e^0_{2k-1}]
$$
is then negligible near the cone tip, and we subsequently improve upon it near the origin in the next iteration. Thus, the updated error is given by
\begin{align*}
    t^2e_{2k} &= t^2[F(u_{2k}) - \square u_{2k}] = t^2[F(v_{2k} + u_{2k-1}) - \square (v_{2k} + u_{2k-1})] \\
    &= t^2[e_{2k-1} - \square v_{2k}] + t^2[F(v_{2k} + u_{2k-1}) - F(u_{2k-1})] \\
    &= t^2[e_{2k-1}-e^0_{2k-1}] + t^2[F(v_{2k} + u_{2k-1}) - F(u_{2k-1})] \\
    &= t^2e^1_{2k-1} + t^2[F(v_{2k} + u_{2k-1}) - F(u_{2k-1})],
\end{align*}
where we used $F(u_{2k-1}) - \square u_{2k-1} = e_{2k-1}$.
\index{a@$a = \lvert x\rvert/t$, a self-similar variable}
The nonlinear part $F(v_{2k} + u_{2k-1}) - F(u_{2k-1})$, which is supported near the tip of the cone, is smaller in magnitude compared to $e_{2k-1}$. It is included within $e^1_{2k}$ and will be further improved upon when constructing the subsequent correction $v_{2k+2}$.

Equation (\ref{odd v}) is solved in the variables $(R,t) = (r\lambda(t),t)$, treating $t$ as a parameter. It allows improving our current error near the origin. The forcing term $e^0_{2k}$ for (\ref{odd v}) is precisely given by the non-negligible component of $e^1_{2k-1}$. The resulting remainder 
$$
t^2e^1_{2k}:= t^2[e_{2k}-e^0_{2k}] = t^2[e^1_{2k-1}-e^0_{2k}] + t^2[F(v_{2k} + u_{2k-1}) - F(u_{2k-1})]
$$
has better smallness properties throughout the entire cone, thus requiring no immediate further improvement. The new error becomes
\begin{align*}
    t^2e_{2k+1} &= t^2[F(u_{2k+1}) - \square u_{2k+1}] \\
    &= t^2[F(v_{2k+1} + u_{2k}) - F(u_{2k}) + F(u_{2k}) - \square (v_{2k+1} + u_{2k})] \\
    &= t^2e^1_{2k} - t^2\partial_t^2 v_{2k+1} + t^2[F(v_{2k+1} + u_{2k}) - F(u_{2k}) - F'(u_0)v_{2k+1}],
\end{align*}
where we have used $F(u_{2k}) - \square u_{2k} = e_{2k}$ and $-\square v_{2k+1} = -\partial_t^2v_{2k+1} - e^0_{2k} - pu_0^{p-1}v_{2k+1}$.At this stage, the error $e_{2k+1}$ has better smallness on the whole cone compared to $e_{2k}$. 

This iterative process is then repeated finitely many times. Specifically, if $1/3 > \varepsilon > 0$, $N_0  \gg \nu > 3$ are fixed, then performing $K_0 = K_0(N_0,\varepsilon) \in \mathbb N$ iterations, where $$2 + \left( \frac{2}{3}-2\varepsilon \right)(K_0-1) \geq N_0$$ leads to an approximate solution of (\ref{wave eq}) with an error of order $\lambda^{\frac{3}{2}}(t\lambda)^{-N_0}$.
\index{N0@$N_0$, the smallness of the approximation error $e_k$ for large $k$, as well as the $\varepsilon$ part from Theorem \ref{thm:blow-up, main thm}}

\section{Renormalization Step: The First Iterate} \label{section:renormalization step, first iterate}

As previously mentioned, (\ref{first v}) can be explicitly solved using a power-series Ansatz (also known as Frobenius method). In this section, we explicitly compute the first correction $v_1$ by solving the elliptic-type equation (\ref{first v}). We then carefully analyze its analytic properties and asymptotic behavior of this first correction, since this first correction forms the foundation for all subsequent steps in the renormalization procedure.

Let us set $u_0(R,t) = \lambda^{\frac{3}{2}}W(R)$, where $W$ is the ground state solution. We define the constants
\begin{equation}
    C_1(\nu) = \frac{105}{128}\pi \nu (1+\nu), \quad C_2(\nu) = \frac{1}{4}(\nu-3)(\nu-5)C_1(\nu),
\index{C1@$C_1(\nu)$, $C_2(\nu)$, some constants coming from $v_1$} \label{eq:constant c_1,c_2}
\end{equation}
which will appear later. First, observe that both $u_0$, $t^2 e_0 \in \lambda^{\frac{3}{2}}C^{\omega}([0,+\infty])$, meaning that they are real-analytic, with an even expansion at $R = 0$ and a regular expansion at $R = +\infty$ with dominant term of order $R^{-3}$. Explicitly,
\begin{align*}
    t^2 e_0(R,t) &= -\lambda^{\frac{3}{2}} \cdot \frac{45 \sqrt{15} (\nu +1) \left(225 (3 \nu +5)+(3 \nu +1) R^4-210 (\nu +1) R^2\right)}{4 \left(R^2+15\right)^{7/2}} \\
    &= \lambda^{\frac{3}{2}} \cdot E_0(R).
\end{align*}
In radial coordinates, (\ref{odd v}) is expressed as
$$t^2 \mathcal{L}_r v_1(r,t) = t^2 e_0(r,t), \quad r \geq 0, \quad
\mathcal{L}_r = -\partial_r^2 - \frac{4}{r} \partial_r - pW(r)^{p-1},$$ where $t$ is treated as a parameter. We seek a solution in the variables $(R,t) = (r\lambda,t)$ variables and rewrite the equation as
\begin{equation}
    (t \lambda)^2 \mathcal{L} v_1(R,t) = t^2 e_0(R,t), \quad R \geq 0, \quad \mathcal{L} = -\partial_R^2 - \frac{4}{R} \partial_R - pW(R)^{p-1}. \label{odd v, change of variables}
\end{equation}
Since $t$ is a parameter and the variables of the forcing term $t^2e_0(R,t)$ are separated, we can expect to find a solution $v_1(R,t)$ in the form
$$v_1(R,t) = \frac{\lambda^{\frac{3}{2}}}{(t\lambda)^2} V_1(R),$$
where $V_1(R)$ is analytic on $[0,+\infty)$ with an even expansion at $R = 0$ starting from order $2$. At $R = +\infty$, an expansion involving powers and logarithms is also expected, consistent with the Frobenius method applied to find solutions to ordinary differential equations with regular singular points.

To facilitate our analysis, instead of restricting ourselves to the positive real line, we shall adopt a complex-analytic framework. Setting $R = z$, we note that the operator $\mathcal{L}$ has a regular singularity at $z = 0$. Consequently, $V_1(z)$ can be found near the origin by looking for a power series solution as in Theorem \ref{thm:inhomogeneous fuchs ode} ($r_1 = 0, r_2 = -3, \beta = 2$). Since 
 $u_0$ and $e_0$ both extend holomorphically to the half-plane
\begin{align*}
    \{z \in \mathbb C: |z| < \sqrt{15} \} \cup \{z \in \mathbb C: \Re(z) > 0\},
\end{align*}
\index{z@$z$, a complex variable which replaces $R$ or $a$ depending on the context}
standard ODE theory ensures that $V_1$ admits a holomorphic extension to the same domain. To study the behaviour of $V_1$ at infinity on the right half-plane, we consider the change of variables $V(z) = V_1(z^{-1})$. Then $V(z)$ solves
\begin{equation}\label{step 1 at infinity}
    V''(z) - \frac{2}{z}V'(z) + z^{-4}pW(z^{-1})^{p-1}V(z) = -z^{-4}E_0(z^{-1}), \quad |z| < \frac{1}{\sqrt{15}}, \ \Re(z) > 0.
\end{equation}
This is again a regular singular equation at $z = 0$, because
$$(z^{-2}+15)^{-1/2} = \frac{z}{(1+15z^2)^{1/2}}$$
on this side of the half-plane. Therefore, $V(z)$ coincides with the particular solution
$$z \sum_{n=0}^{+\infty}v_nz^n + v_{-1}u_1(z)\log(z)$$
given by Theorem \ref{thm:inhomogeneous fuchs ode} ($r_1=3, r_2 = 0, \beta = 1$), modulo some linear combination of the fundamental system $\{u_1(z),u_2(z)\}$ (see (\ref{fuchs fundamental system}) from Section \ref{section:appendix, ode}) which introduces a dominant term of order $z^0$. On the real-line, one explicitly has
\begin{equation} \label{v_1 formula}
    V_1 (R) = \frac{\begin{aligned}
                   &\sqrt{15} (\nu +1) \Big[-360 \left(R^2-15\right) R^5-100800 \nu  \left(R^2-15\right) R^3 \text{arccoth} \left(\frac{30}{R^2}+1\right) \\
                   &-75 \nu  \left(13 R^6-1397 R^4+6195 R^2+4725\right) R \\
                   &+7 \sqrt{15} \nu  \left(R^8+300 R^6-20250 R^4+67500 R^2+50625\right) \arctan \left(\frac{R}{\sqrt{15}}\right)\Big]
                \end{aligned}}{64 R^3 \left(R^2+15\right)^{5/2}}.
\end{equation}
It follows that $V_1(R)$ is real-analytic on $[0,+\infty)$ with an even Taylor expansion at zero starting from order $2$ with a positive coefficient. At infinity, it has an asymptotic expansion of the form:
$$\underbrace{\frac{105}{128} \pi \nu (1+\nu)}_{=:C_1(\nu)} R^0 - \frac{45}{8}\sqrt{15}(1+\nu)(1+3\nu)R^{-1} + \frac{55125}{256}\pi \nu(1+\nu)R^{-2} + \mathcal{O}(\log(R)R^{-3}),$$
which is positive as well. However, $V_1(R)$ may take negative values on a fixed compact interval in $(0,+\infty)$; this does not pose a problem because if $t_0$ is small enough, then $u_0(R,t) + v_1(R,t) > 0$ on $[0,+\infty) \times (0,t_0]$.

Functions such as $u_0$ and $v_1$, which have an even power series expansion at $R = 0$ and a series of power and logarithms (with bounded logarithmic exponent) at infinity, are regrouped into the following function space $S^{2n}(R^I,\log(R)^J)$. This function space is used to describe correction and error terms near the origin of the light cone.

\begin{definition}[Space $S^{2n}(R^I,\log(R)^J)$] \label{space s^m(r^k)}
Let $R_0 = 4\sqrt{15}$ be fixed. \index{R0@$R_0 = 4\sqrt{15}$, the threshold distinguishing small and large $R$} \index{S-space@$S^{2n}(R^I,\log(R)^J)$, a vector space of smooth functions with respect to the radial variable $R$} For $I,J,n \in \mathbb N_{\geq 0}$, define $S^{2n}(R^I,\log(R)^J)$ as the vector space of smooth functions $w(R)$ on $(0,+\infty)$ which satisfy:
\begin{enumerate}
    \item $w$ admits a holomorphic extension to an open neighbourhood of the set
$$\{z \in \mathbb C: |z| \leq \sqrt{15}/2 \} \cup \{z \in \mathbb C: \Re(z) > 0, |z| \geq R_0\}.$$
    \item $w$ has a zero of order $2n$ at $z = 0$ with an even Taylor expansion.
    \item For $\Re(z) > 0$ and $|z| \geq R_0$, $w(z)$ has the following expansion:
    $$w(z) =  z^I \sum_{j=0}^{J}w_{j}(z^{-1})
 \log(z)^j,$$
    where the functions $w_{j}(y)$ are holomorphic on a neighbourhood of $\{y \in \mathbb C: |y| \leq R_0^{-1}\}$.
\end{enumerate}
For $I < 0$, the space $S^{2n}(R^I,\log(R)^J)$ is defined as the subspace of $S^{2n}(R^0,\log(R)^J)$ for which each function $w_{j}$ in the expansion at infinity has a zero of order at least $|I|$. 
\end{definition}

\begin{remark} \label{rmk: s^2n properties}
 We start with a few remarks regarding the spaces $S^{2n}(R^I, \log(R)^J)$.
\begin{enumerate} 
\item For any $l \in \mathbb N$, $(R\partial_R)^l S^{2n}(R^I,\log(R)^J) \subset S^{2n}(R^I,\log(R)^J)$. 
    \item An element $w \in S^{2n}(R^m,\log(R)^J)$ is not necessarily holomorphic on the whole half-plane $\Re(z) > 0$, allowing room for localized cutoff.
    \item The subspaces $S^{2n}(R^I,\log(R)^J)$ with $I < 0$ are only relevant for the functions $u_0$, $t^2e_0$ and $u_0^p$. In these cases, one has $u_0, t^2e_0 \in \lambda^{\frac{3}{2}} S^0(R^{-3}, \log(R)^0)$, $u_0^p \in \lambda^{\frac{3p}{2}} S^0(R^{-7},\log(R)^0)$ and $v_1 \in \lambda^{\frac{3}{2}} (t\lambda)^{-2} S^2(R^{0},\log(R))$, where no logarithm occurs in the dominant term at infinity. More precisely, $w_0(y) = \mathcal{O}(1)$ is the dominant component, while $w_{j}(y) = \mathcal{O}(y^3)$, $j \neq 0$.
    \item One has
    \begin{align*}
        t^2 u_0^p \left( \frac{v_1}{u_0} \right)^{n} &\in \dfrac{\lambda^{\frac{3}{2}}}{(t\lambda)^{2(n-1)}} S^{2n}(R^{-7+3n},\log(R)^{n} ).
    \end{align*}
\end{enumerate} 
\end{remark}

\section{Renormalization Step: The Hypergeometric Function and the second iterate}

Having established the appropriate function spaces required for solving (\ref{odd v}), it is instructive, before proceeding with the full iterative scheme, to motivate and introduce the function spaces $\mathcal{Q}$ that will occur when solving equation (\ref{even v}) near the tip of the cone. The simplest case for equation (\ref{even v}), when rewritten in radial coordinates, corresponds to a forcing term that is a pure power of $t$:
$$
t^2\left(-\partial_t^2 + \partial_r^2 + \frac{d-1}{r}\partial_r\right) v =t^s.
$$
This scenario will appear as the dominant part, near the tip of the cone, of the error $e_1$ that remains after the first iteration. Equivalently, if we look for a solution $v(r,t) = t^s w(r,t)$, then: 
$$
t^2 \left( - \left( \partial_t + \frac{s}{t} \right)^2 + \partial_r^2 + \frac{d-1}{r}\partial_r \right)w  = 1.
$$
Introducing the self-similar variable $a = r/t$, this equation becomes $L_{s} w(a) = 1$, $0 < a < 1$, where
$$L_{s} = (1-a^2)\partial_{aa} + ((d-1)a^{-1} + 2as - 2a)\partial_a + (s - s^2).$$
Seeking a solution of the form $w(a) = W(a^2)$, we naturally reduce to an hypergeometric equation for $W(z)$:
$$z(1-z)W''(z) + \left( \frac{d}{2} + z \left(s - \frac{3}{2} \right) \right) W'(z) + \frac{s - s^2}{4} = 1, \quad 0 < z < 1$$
whose parameters are
$$
\tilde{\alpha} = -\frac{s}{2}, \tilde{\beta} = -\frac{s}{2} + \frac{1}{2}, \tilde{\gamma} = \frac{d}{2},
$$
and for which a particular solution is explicitly given by
$$
\frac{1}{\alpha \beta}[F(\tilde{\alpha},\tilde{\beta};\tilde{\gamma},z) - 1],
$$
where $F$ is the Gauss hypergeometric function defined as follows.
\begin{definition}\label{gauss hypergeometric definition}
Let $\alpha, \beta, \gamma \in \mathbb R$, $\gamma \notin \mathbb Z_{\leq 0}$. The Gauss hypergeometric function $F$ is defined by the series
\begin{equation*} 
F(\alpha,\beta;\gamma;z) = \sum_{n=0}^{+\infty} \frac{(\alpha)_n(\beta)_n}{(\gamma)_nn!}z^n, \quad |z| < 1, \quad (x)_n = x(x+1)...(x+n-1), \quad (x)_0 = 1
\end{equation*}
\index{F(alpha)@$F(\alpha,\beta;\gamma;z)$, the Gauss Hypergeometric function}
which converges absolutely for $|z| < 1$, as well as for $|z| \leq 1$ if $\gamma-\alpha-\beta > 0$.
\end{definition}

The hypergeometric function satisfies the hypergeometric equation
\begin{equation}\label{hypergeometric ode}
    z(1-z)w''(z) + (\gamma - (\alpha + \beta + 1)z)w'(z) - \alpha \beta w(z) = 0, \quad 0 < z <1.
\end{equation}
This equation has regular singular points at $z = 0$ and $z = 1$. The roots of the indicial equation at these points are respectively $\{0,1-\gamma\}$ and $\{0,\gamma-\alpha-\beta\}$. Fundamental systems can be found as described in Section \ref{section:appendix, ode} (see (\ref{fuchs fundamental system})).

In our application, we focus on the case where $\gamma > 0$, $\alpha \beta > 0$, which ensures monotonicity and positivity of $F(\alpha,\beta;\gamma,z) - 1$ on $[0,1)$, as shown in Corollary \ref{monotonicity hypergeom function}. This property is essential to ensure positivity of the second correction term $v_2$. Additionally, we will also impose $\gamma - \alpha - \beta > 0$, which allows us to determine which indicial root is the largest at $z = 1$ and guarantees continuity of $F(\alpha,\beta;\gamma,z)$ up to the boundary $|z| = 1$.

\begin{lemma}[A monotonicity and positivity result]\label{monotonicity hypergeom eq}
Let $\alpha \beta > 0$ and $g(z) \in C^0([0,1))$ be non-negative. Suppose $w(z)$ is a $C^1([0,1)) \cap C^2((0,1))$-solution to the inhomogeneous hypergeometric equation
\begin{equation*}
    z(1-z)w''(z) + (\gamma - (\alpha + \beta + 1)z)w'(z) - \alpha \beta w(z) = g(z), \quad 0 \leq z < 1.
\end{equation*}
with initial conditions $w(0) \geq  0, w'(0) > 0$. Then $w(z)$ is strictly increasing on $[0,1)$.

\end{lemma}

\begin{proof}
At $z = 0$, we have that $w'(0) > 0$, which implies that $w(z)$ is strictly increasing near $z = 0$. Let $I = [0,z_0)$ be the maximal interval on which $w'(z) > 0$. Suppose for a contradiction that $z_0 < 1$. Then $w'(z_0) = 0$ and $w(z_0) > w(0) \geq 0$ by continuity and strict monotonicity. Evaluating the differential equation at $z_0$, we obtain:
$$w''(z_0) = \frac{\alpha \beta}{z_0(1-z_0)}w(z_0) + \frac{g(z_0)}{z_0(1-z_0)} > 0,$$
which implies that $z_0$ is a local strict minimum of $w(z)$ and $w(z)$ is strictly decreasing on the left of $z_0$: a contradiction.
\end{proof}

\begin{corollary}\label{monotonicity hypergeom function}
If $\gamma > 0$ and $\alpha \beta > 0$, then $z \mapsto F(\alpha,\beta;\gamma,z) - 1$ is strictly increasing on $[0,1)$.
\end{corollary}

\begin{proof}
    The proof follows from Lemma \ref{monotonicity hypergeom eq} since $w(z) = F(\alpha,\beta;\gamma,z) - 1$ solves
    $$z(1-z)w''(z) + (\gamma - (\alpha + \beta + 1)z)w'(z) - \alpha \beta w(z) = \alpha \beta, \quad 0 \leq z < 1$$
    with $w(0) = 0$ and $w'(0) = \alpha \beta \gamma^{-1} > 0$.
\end{proof}


We now define a specific hypergeometric-type function
    \begin{align} \label{H(z) definition}
    H(z) &= 4 C_1(\nu) (F(\tilde{\alpha},\tilde{\beta};\tilde{\gamma},z)-1 ), \quad 0 \leq z < 1 \\
    s &=  \frac{3}{2}(-1-\nu) + 2\nu, \tilde{\alpha} = -\frac{s}{2}, \tilde{\beta} = -\frac{s}{2} + \frac{1}{2}, \tilde{\gamma} = \frac{5}{2} \notag
    \end{align} \index{H(z)@$H(z)$, an hypergeometric-like function}
     with $C_1$ defined as in (\ref{eq:constant c_1,c_2}). It turns out that the term
     $$
      \frac{\lambda^{\frac{2}{3}}}{(t\lambda)^2}H(a^2) 
     $$
    will be the dominant component of the second correction term $v_2$. The parameters satisfy $\tilde{\alpha} \tilde{\beta} > 0$ if $\nu > 3$, $\tilde{\gamma} > 0$ and $\tilde{\gamma}-\tilde{\alpha}-\tilde{\beta} = \frac{1}{2} + \frac{1}{2}\nu > 0$. Therefore, $H(z)$ is strictly increasing on $[0,1)$ (Corollary \ref{monotonicity hypergeom function}), holomorphic on $|z| < 1$ and continuous on $|z| \leq 1$. Using a fundamental system for the hypergeometric equation (see (\ref{fuchs fundamental system})), we obtain the singular expansion at $z = 1$:
    \begin{align} \label{expansion of H(z) near z = 1}
    F(\tilde{\alpha},\tilde{\beta};\tilde{\gamma},z)-1 = q_0(1-z)+q_1(1-z)(1-z)^{\frac{1}{2}+\frac{1}{2}\nu}+q_2(1-z)(1-z)^{\frac{1}{2}+\frac{1}{2}\nu}\log(1-z),
    \end{align}
    where $q_i(y)$ is holomorphic on $|y| < 1$. Moreover, since $H(1) > 0$, one has $q_0(0) > 0$ as well. 
    
    As with the space $S^{2n}(R^I, \log(R)^J)$, we now define a family of functions $\mathcal{Q}$ in the self-similar variable $a = r/t$ having a series of power and logarithms at $a = 0$ and $a = 1$. These functions will be used to describe correction and error terms near the tip of the cone. Their structure and regularity are motivated by the Frobenius method applied to the hypergeometric equation. However, extra care must be taken with regard to the domain on which our functions are holomorphic. Although $H(z)$ is holomorphic on $|z| < 1$, the nonlinear term
    $$
    \left( u_0 + v_1 + \frac{\lambda^{\frac{2}{3}}}{(t\lambda)^2}H(a^2) \right)^p
    $$
    has, near $r \sim t$, $a = 1$, a dominant component of the form
    $$
    \frac{\lambda^{\frac{2}{3}}}{(t\lambda)^2}\left( H(a^2) + C_1(\nu) \right)^p
    $$
    with the constant term coming from $v_1$. This expression is holomorphic only on the subset of the unit disk where $H(z) + C_1(\nu) \neq 0$, which may be smaller than $|z| < 1$. Thus, to include this function in the family $\mathcal{Q}$, we cannot demand holomorphy on the entire unit disk. Fix $U$ to be an open, simply-connected neighbourhood $U \subset B(0,1) \subset \mathbb C$ such that
    \begin{enumerate}
        \item $(0,1) \subset U$,
        \item $U$ contains $\{a \in B(0,1): |a-0| \leq a_0 \text{ or } |a-1| \leq a_0\}$ for some $0 < a_0 \ll 1$,
        \item $H(a)$ and $H(a) + C_1(\nu)$ (defined in (\ref{H(z) definition})) are non-zero everywhere on $U \setminus \{0\}$,
        \item If $q_0(y)$ is the analytic function appearing in the expansion of $H(z)$ near $z = 1$ (see (\ref{expansion of H(z) near z = 1})), then $q_0(y)$ and $q_0(y) + C_1(\nu)$ are non-zero everywhere on $|y| \leq a_0$.
    \end{enumerate}
    
\begin{definition}[Space $\tilde{\mathcal{Q}}$ and $\mathcal{Q}$]\label{definition of tilde Q}
\index{Q@$\tilde{\mathcal{Q}}$, $\mathcal{Q}$, families of holomorphic functions with respect to the self-similar variable $a$}
Let $U$ be as described above. Let $\tilde{\mathcal{Q}}_{\beta}$, $\beta \in \mathbb R$, be the vector space of real-analytic functions $q(a): (0,1) \rightarrow \mathbb R$ satisfying:
\begin{enumerate}
\item $q(a)$ extends holomorphically to $U$.
    \item Near $a = 0$, one has the finite expansion:
    $$q(a) = q_0(a) + \sum_{j=1}^{L}q_{j}(a) \log(a)^j,$$
    where $0 \leq L < +\infty$ and each $q_{j}$ is holomorphic on a neighbourhood of $|a-0| \leq a_0$.
    \item Near $a = 1$, one has the expansion:
    $$q(a) = q_{0,0}(1-a) + \sum_{i=1}^{+\infty}(1-a)^{\beta(i)} \sum_{j=0}^{L_i}q_{i,j}(1-a)\log(1-a)^j,$$
    where $L_i \geq 0$, $\sup_{i\geq 1}L_i/i < +\infty$ , the functions $q_{i,j}(y)$ are holomorphic on a neighbourhood of $|y| \leq a_0$, $\beta(i) \geq \beta$ and either the expansion is finite, i.e.,
    $$|\{(i,j) \in \mathbb N^2: q_{i,j}(y) \not \equiv  0\}| < +\infty, $$
    or the growth of $i \mapsto \beta(i)$ is at least linear, i.e., 
    $$\beta(i+1) > \beta(i), \quad \lim \limits_{i \to +\infty}\beta(i) = +\infty, \quad \inf_{i \geq 2}\left(\frac{\beta(i)-\beta}{i}\right) > 0.$$
    Moreover, there should exist $\varepsilon > 0$ and $C > 0$ for which
    $$|\beta(i)| + ||q_{i,j}||_{A(|y| \leq a_0+\varepsilon)} \leq C^i \quad \forall 0 \leq j \leq L_i, \forall i \in \mathbb N_{\geq 1}.$$
    where $A\left(|y| \leq R \right)$ is the Wiener algebra (see Definition \ref{wiener space}). In other words, the growth of the sequences $\beta(i)$ and $q_{i,j}$ must be at most exponential exponential.
\end{enumerate}
Finally, we define $\mathcal{Q}_{\beta}$ as the vector space of real-analytic functions on $(0,1)$ of the form $a \mapsto \tilde{q}(a^2)$ for some $\tilde{q} \in \tilde{\mathcal{Q}}_{\beta}$.
\end{definition}

\begin{remark}
We start with some remarks concerning the spaces $\mathcal{Q}$ and $\tilde{\mathcal{Q}}$.
\begin{enumerate}
    \item We will primarily work with the spaces $\tilde{\mathcal{Q}}_{\frac{1}{2} + \frac{1}{2}\nu}$ and $\mathcal{Q}_{\frac{1}{2} + \frac{1}{2}\nu}$, but it is convenient, for theoretical results, to allow any $\beta \in \mathbb R$.
    \item As established in (\ref{H(z) definition}) and (\ref{expansion of H(z) near z = 1}), $H(a) \in \tilde{\mathcal{Q}}_{\frac{1}{2} + \frac{1}{2}\nu}$ with a finite expansion at $a = 1$. Moreover, for any exponent $e \in \mathbb R$, $(H(a) + C_1(\nu))^e \in \tilde{\mathcal{Q}}_{\frac{1}{2} + \frac{1}{2}\nu}$ as well. The expansion near $a = 1$ is obtained via a binomial expansion:
    \begin{align*}
        (H(z) + C_1(\nu))^e = \sum_{n=0}^{\infty} \binom{e}{n} &\left( q_0(1-z) + C_1(\nu) \right)^{e-n} \cdot \\
        &(1-z)^{\left( \frac{1}{2} + \frac{1}{2}\nu\right) n} \left( q_1(1-z) + q_2(1-z)\log(1-z) \right)^n.
    \end{align*}
    \item We do not require that the expansion of a $\tilde{Q}$-element near $a = 1$ converges everywhere on $|a-1| < a_0$. However, the coefficient functions $q_{i,j}$ must all be defined and holomorphic around the disk $|a-1| \leq a_0$. The growth condition on $q_{i,j}$ ensures uniform and absolute convergence of the expansion on a smaller ball (depending only on $a_0$) around $a = 1$.
    \item The linear growth estimate on $\beta(i)$ is essential. If we solve the hypergeometric equation near $z = 1$ with a forcing term of the form $(1-z)^{\beta(i)}q_{i,j}(1-z)\log(1-z)^j$, then the estimate (\ref{exponential upper bound on w_k}) on the solution from Theorem \ref{thm:inhomogeneous fuchs ode} takes the form:
    \begin{equation*}
          ||w(1-z)||_{A(|1-z| < a_0 + \varepsilon)}  \leq C^i \cdot ||q_{i,j}||_{L^{\infty}(|1-z| < a_0 + \varepsilon}.
   \end{equation*}
   Thus, if the growth of $q_{i,j}$ is at most exponential in $i$, then so is the growth of the solution. Accordingly, if we solve (\ref{even v}) with a forcing term from a family $\mathcal{Q}_{\beta}$, we expect the solution to remain within $\mathcal{Q}_{\tilde{\beta}}$.
\end{enumerate}
\end{remark}

\begin{proposition}[Product rules] \label{q algebra calculation rules}
    The following algebraic rules hold: 
    \begin{enumerate}
        \item[A.] Differentiation: $\partial_a \left( a^{\delta} \mathcal{Q}_{\beta} \right)  \subset \delta a^{\delta-1} \mathcal{Q}_{\beta} + a^{\delta-1} \mathcal{Q}_{\beta-1}$ for any $\beta, \delta \in \mathbb R$. 
        \item[B.] Summation: $\mathcal{Q}_{\beta_1} + \mathcal{Q}_{\beta_2} \subset \mathcal{Q}_{\min\{\beta_1,\beta_2\}}$ for any $\beta_1,\beta_2 \in \mathbb R$. 
        \item[C.] Product: $\mathcal{Q}_{\beta_1} \cdot \mathcal{Q}_{\beta_2} \subset  \mathcal{Q}_{\min\{\beta_1,\beta_2,\beta_1+\beta_2\}}$ for any $\beta_1,\beta_2 \in \mathbb R$. In particular, if $\beta \geq 0$, then $\mathcal{Q}_{\beta}$ is an algebra.
    \end{enumerate}
\end{proposition}

\begin{proof}
    We prove only the differentiation rule. Suppose $q(a) = a^{\delta}Q(a^2)$, where $Q \in \tilde{\mathcal{Q}}_{\beta}$. Then, it holds that
    $$
    \partial_a q(a) = \delta a^{\delta-1}Q(a^2) + 2a^{\delta+1}Q'(a^2) = a^{\delta-1} \left( \delta Q(a^2) + 2a^2Q'(a^2) \right).
    $$
    It suffices to show that $a\partial_a Q(a) \in \tilde{\mathcal{Q}}_{\beta-1}$. Since $Q(a)$ is holomorphic on $U$, so is $a\partial_aQ$. Near $a = 0$, one has
    \begin{align*}
        Q(a) &= q_0(a) + \sum_{j=1}^{L}q_{j}(a) \log(a)^j, \\
        aQ'(a) &=aq_0'(a) + \sum_{j=1}^{L}aq_{j}'(a) \log(a)^j + \sum_{j=1}^{L} j q_{j}(a) \log(a)^{j-1},
    \end{align*}
    where $0 \leq L < +\infty$ and each $q_{j}, a\partial_a q_j$ is holomorphic on a neighbourhood of $|a-0| \leq a_0$. Near $a = 1$,  
   \begin{align*}
       Q(a) &= q_{0,0}(1-a) + \sum_{i=1}^{+\infty}(1-a)^{\beta(i)} \sum_{j=0}^{L_i}q_{i,j}(1-a)\log(1-a)^j,
   \end{align*}
    where $0 \leq L_i \leq L(i+1)$ for some $L > 0$, each $q_{i,j}(y)$ is holomorphic on $|y| < a_0 + \varepsilon$, $\beta(i) \geq \beta$ and we assume without loss of generality that we have an infinite series with 
    $$\beta(i+1) > \beta(i), \quad \lim \limits_{i \to +\infty}\beta(i) = +\infty, \quad \inf_{i \geq 2}\left(\frac{\beta(i)-\beta}{i}\right) > 0,$$
    as well as
    $$|\beta(i)| + ||q_{i,j}||_{A(|y| \leq a_0+\varepsilon)} \leq C^i \quad \forall 0 \leq j \leq L_i, \forall i \in \mathbb N_{\geq 1},$$
    for some constants $C, \varepsilon > 0$. Formally consider the sum of derivatives:
    \begin{align*}
           S(a) = q_{0,0}'(1-a) &- \sum_{i=1}^{+\infty}(1-a)^{\beta(i)-1} \sum_{j=0}^{L_i}\beta(i)q_{i,j}(1-a)\log(1-a)^j \\
           &-\sum_{i=1}^{+\infty}(1-a)^{\beta(i)} \sum_{j=0}^{L_i}q_{i,j}'(1-a)\log(1-a)^j \\
           &-\sum_{i=1}^{+\infty}(1-a)^{\beta(i)-1} \sum_{j=0}^{L_i}jq_{i,j}(1-a)\log(1-a)^{j-1}.
    \end{align*}
    For all $0 \leq j \leq L_i, \forall i \in \mathbb N_{\geq 1}$, one has
    $$
    ||\beta(i)q_{i,j}||_{A(|y| \leq a_0+\varepsilon)} +  ||j q_{i,j}||_{A(|y| \leq a_0+\varepsilon)}  \leq C^iC^i + L(i+1)C^i \leq \tilde{C}^i,
    $$
    as well as
    $$
||q_{i,j}'||_{A(|y| \leq a_0+\varepsilon/2)} \lesssim_{a_0,\varepsilon} ||q_{i,j}||_{A(|y| \leq a_0+\varepsilon)} \leq \tilde{C}^i.
    $$
    Hence, the sum of derivatives converges normally on some ball $|a-1| \lesssim \min\{\tilde{C}^{-1},a_0\}$. Integrating, we recover $Q(a)$ modulo some additive constant. Thus, $S(a) = Q'(a)$ by the Identity Theorem, and $Q'(a)$ has the desired expansion near $a = 1$.
\end{proof}

\section{Renormalization Step: Preliminaries for the next iterates} \label{section:renormalization step, preliminaries next iterates}

This section establishes the technical framework required for all subsequent renormalization steps. We first partition the light cone into three distinct regions to analyze nonlinear error terms such as:
$$
t^2e_1(R,t) = t^2[F(u_0 + v_1) - F(u_0) - F'(u_0)v_1] - t^2\partial_{tt} (v_1(r\lambda,t))
$$
through a multinomial expansion within each region. The idea is that it will be easier to solve (\ref{even v}) and (\ref{odd v}) for each term in the expansion and sum everything back. The multinomial expansion depends on whether $u_0 \geq v_1$, $u_0 \sim v_1$ or $u_0 \leq v_1$, which corresponds to the three regions $R \lesssim (t\lambda)^{\frac{2}{3}}, R \sim (t\lambda)^{\frac{2}{3}}, R \gtrsim (t\lambda)^{\frac{2}{3}}$ of the light cone. 

The first step is to fix an appropriate constant $m \ll 1$ and define the region $R \leq m(t\lambda)^{\frac{2}{3}}$. This constant is chosen so that $(u_0 + v_1)^p$ can be reliably expanded around $u_0$ in that region $R \leq m(t\lambda)^{\frac{2}{3}}$ using a binomial expansion. Moreover, $m$ should be sufficiently small so that one can keep only the terms from the expansion with a bounded logarithmic exponent. This restriction is important, as solving equation (\ref{odd v, change of variables}) at $R = +\infty$ with a logarithmic forcing term $R^{I}\log(R)^J$, $I \leq 0$, introduces factors of order $J!$ in the solution. Even if $|I| \gg J$, the bound (\ref{exponential upper bound on w_k}) from Theorem \ref{thm:inhomogeneous fuchs ode} is not available. This makes controlling convergence of solutions more difficult if one allows for a sequence of logarithmic terms $\log(R)^J$, $J \to +\infty$, in the binomial expansion.

Once the constant $m$ is defined, we will divide the cone into three main regions, define the correction and error function spaces and prove multiple computation rules that will facilitate the iteration scheme.

\begin{definition}[Constant $m$]
    Let $t_0 \ll 1$ be sufficiently small so that
    \begin{enumerate}
        \item $|u_0(z,t)| \geq 2|v_1(z,t)|$ on $$\left( \{z \in \mathbb C: |z| \leq \sqrt{15}/2\} \cup \{z \in \mathbb R: z \in [0,2\sqrt{15}]\} \right) \times (0,t_0].$$
        \item $u_0(R,t) + v_1(R,t) > 0$ on $[0,+\infty) \times (0,t_0]$.
    \end{enumerate}
    For $\Re(z) > 0, |z| > \sqrt{15}$, consider the expansion:
     $$\frac{v_1(z,t)}{u_0(z,t)} = \frac{z^{3}}{(t\lambda)^{2}} \left(w_0(z^{-1}) + z^{-3}w_1(z^{-1})\log(z) \right),$$
 where $w_j(y)$ is holomorphic on $|y| < (\sqrt{15})^{-1}$, $w_j(0) \neq 0$. Fix $m \ll 1$ any constant for which
    \begin{equation}
        (2m)^3 \max_{j \in \{0,1\}}\left[1+||w_j||_{A\left(|y| \leq 2/(3\sqrt{15}) \right)}\right] \leq \frac{1}{4}. \label{m definition}
    \end{equation}
    \index{m@$m$, the constant in $R \leq m(t\lambda)^{\frac{2}{3}}$} \sloppy In particular, this choice ensures:
    $|u_0(z,t)| \geq 2|v_1(z,t)|$
 on $$\{(z,t) \in \mathbb C \times (0,1): \Re(z) > 0, 2\sqrt{15} \leq |z| \leq m(t\lambda)^{\frac{2}{3}}, 0 < t \leq t_0\},$$
    which is a non-empty set for $t_0$ small enough. 
\end{definition}

Given this constant $m$, we can decompose the cone into four regions.

\begin{definition}[$k$-admissible pairs]
Define the overlapping regions:
\begin{align*}
    C_{\mathrm{ori}} &= \{(R,t): 0 < t \leq t_0, 0 \leq R \leq m (t\lambda)^{\frac{2}{3}}\}, \\
    C_{\mathrm{mid}} &= \{(R,t): 0 < t \leq t_0, \frac{m}{2}(t\lambda)^{\frac{2}{3}} \leq R \leq 2(t\lambda)^{\frac{2}{3}+\varepsilon}\}, \\
    C_{\frac{2}{3}+\varepsilon} &= \{(R,t): 0 < t \leq t_0,(t\lambda)^{\frac{2}{3}+\varepsilon} \leq R \leq 2(t\lambda)^{\frac{2}{3}+\varepsilon}\}, \\
    C_{\mathrm{tip}} &= \{(R,t): 0 < t \leq t_0, (t\lambda)^{\frac{2}{3}+\varepsilon} \leq R \leq (t\lambda)\}. 
\end{align*}
\index{C1@$C_{\mathrm{ori}}$, the subset $0 \leq R \leq m (t\lambda)^{\frac{2}{3}}$ of the cone}
\index{C2@$C_{\mathrm{mid}}$, the subset $\frac{m}{2}(t\lambda)^{\frac{2}{3}} \leq R \leq 2(t\lambda)^{\frac{2}{3}+\varepsilon}$ of the cone}
\index{C3@$C_{\frac{2}{3}+\varepsilon}$, the subset $(t\lambda)^{\frac{2}{3}+\varepsilon} \leq R \leq 2(t\lambda)^{\frac{2}{3}+\varepsilon}$ of the cone}
\index{C4@$C_{\mathrm{tip}}$, the subset $(t\lambda)^{\frac{2}{3}+\varepsilon} \leq R \leq (t\lambda)$ of the cone}
A pair of indices $(\alpha,i) \in \mathbb R \times \mathbb Z$ will be called $k$-admissible for $k > 1$ on a region $C_x \neq C_{\mathrm{ori}}$ if 
$$\forall (R,t) \in C_{x}: \frac{|R|^i}{(t\lambda)^{\alpha}} \lesssim \frac{1}{(t\lambda)^{2+\left( \frac{2}{3}-2\varepsilon \right) \cdot (k-1)}}.$$
\index{k-admi@$k$-admissible pairs}
When $k = 1$, the pair $(\alpha,i)$ is called 1-admissible if $(\alpha,i) = (2,0)$ or if it is 2-admissible.

On $C_{\mathrm{ori}}$, the pair will be called $1$-admissible if $(\alpha,i) = (2,0)$ or 
$$\forall (R,t) \in C_{\mathrm{ori}}: \frac{|R|^i}{(t\lambda)^{\alpha}} \lesssim \frac{1}{(t\lambda)^{2+\frac{2}{3}}}$$
and $k$-admissible for $k > 1$ if $(\alpha-2,i-2)$ is $(k-1)$-admissible. In particular, this implies
$$\forall (R,t) \in C_{\mathrm{ori}}: \frac{|R|^i}{(t\lambda)^{\alpha}} \lesssim \frac{1}{(t\lambda)^{2+ \frac{2}{3} \cdot (k-1)}}.$$
We will usually omit the region $C_x$ as it will be clear from the context.
\end{definition}

The notion of admissible pairs allows to easily quantify the smallness of our correction and error terms. We also observe that if $(\alpha,i)$ is $k$-admissible on $C_{\mathrm{ori}}$, then $i \geq 0$ and $(\alpha,i)$ is also $k$-admissible on $C_{\mathrm{mid}}$. Similarly, if $(\alpha,i)$ is $k$-admissible on $C_{\mathrm{mid}}$ or $C_{\mathrm{tip}}$, then $(\alpha,i)$ is also $k$-admissible on $C_{\frac{2}{3}+\varepsilon}$.

Next, we endow the space $w \in S^{2n}(R^I, \log(R)^J)$ with a norm giving control on $w$ on $R \leq m(t\lambda)^{\frac{2}{3}}$. This norm is useful to consider series of $S^{2n}(R^I, \log(R)^J)$ elements and ensure their convergence.

\begin{definition}[Norm on the space $S^{2n}(R^I,\log(R)^J)$] \label{norm of space s^m(r^k)}
 For $I,J,k,n \in \mathbb N_{\geq 0}$, consider an element $w \in S^{2n}(R^I,\log(R)^J)$ and its expansion:
 $$w(z) =  z^I \sum_{j=0}^{J}w_{j}(z^{-1})\log(z)^j$$
 on $R \geq R_0$. Define the following semi-norms:
\begin{align*}
||w||_{S, \mathrm{ori}} &= ||w(z)||_{A(|z| \leq \sqrt{15}/2)} + ||w(z)||_{L^{\infty}(z \in [\sqrt{15}/2, R_0] \subset \mathbb R)}, \\
||w||_{S,I,J,\infty} &=  m^I \max_{0 \leq j \leq J} ||w_{j}(y)||_{A(|y| \leq R_0^{-1})},
\end{align*}
\index{S-space norms@$\lvert \lvert w \rvert \rvert_{S, \mathrm{ori}}, \lvert \lvert w \rvert \rvert_{S,I,J,\infty}$, semi-norms on $S^{2n}(R^I,\log(R)^J)$}
the sum of which creates a norm.
\end{definition}

\begin{remark}
 We start with a few observations regarding the semi-norms on $S^{2n}(R^I, \log(R)^J)$.
\begin{enumerate} 
    \item One has a product structure: if $v \in S^{2n_1}(R^{I_1},\log(R)^{J_1})$ and $w \in S^{2n_2}(R^{I_2},\log(R)^{J_2})$ for $n_1,n_2,I_1,I_2,J_1,J_2 \geq 0$, then $v \cdot w \in S^{2(n_1+n_2)}(R^{I_1+I_2},\log(R)^{J_1+J_2})$ and
    \begin{align*}
||v \cdot w||_{S, \mathrm{ori}} &\leq ||v||_{S, \mathrm{ori}} \cdot ||w||_{S, \mathrm{ori}}, \\
||v \cdot w||_{S, I_1 + I_2, J_1 + J_2, \infty} &\leq  ||v||_{S,I_1,J_1, \infty} \cdot ||w||_{S,I_2,J_2, \infty}.
\end{align*}
    \item We will not define, nor use, a norm $|| \cdot ||_{S,I,J,\infty}$ with $I < 0$.
    \item If $v \in S^{0}(R^{-I_1},\log(R)^{J_1})$ is fixed and $w \in S^{2n}(R^{I_2},\log(R)^{J_2})$ for $n,I_1,I_2,J_1,J_2 \geq 0$, then $v \cdot w \in S^{2n}(R^{-I_1+I_2},\log(R)^{J_1+J_2})$ and
    \begin{align*}
||v \cdot w||_{S, \mathrm{ori}} &\leq C(I_1, m, R_0, v) \cdot ||w||_{S, \mathrm{ori}}, \\
||v \cdot w||_{S, -i + I_2, J_1 + J_2, \infty} &\leq C(I_1, m, R_0, v) \cdot ||w||_{S,I_2,J_2, \infty} \quad \forall 0 \leq i \leq \min\{I_1,I_2\},
\end{align*}
i.e., the product map is continuous.
\end{enumerate} 
\end{remark}

\begin{notation}
Throughout the whole paper, $\chi_{[a,+\infty)}: \mathbb R \rightarrow [0,1]$ denotes a fixed smooth transition function which satisfies $\chi(x) = 0$ in a neighbourhood of $x = a$, $\chi(x) = 1$ in an open neighbourhood of $[2a,+\infty)$ and $\chi(x) > 0$ on $\mathrm{int}(\mathrm{supp}(\chi_{[a,+\infty)}))$. In particular, such a transition function is supported on $[a,+\infty)$. Explicitly, up to an affine transformation (depending only on $a \in \mathbb R$), one can choose
$$
\chi(x) = \frac{f(x)}{f(x) + f(1-x)}, \quad f(x) = \begin{cases}
    e^{-\frac{1}{x}}, \quad &x > 0 \\
    0, \quad &x \leq 0
\end{cases}
$$
\index{chi@$\chi_{[a,+\infty)}$, a fixed transition function}
\end{notation}

On our three main regions of the light cone, observe that:
 \begin{enumerate}
     \item $C_{\mathrm{ori}}$: $F(u_0+v_1) = (u_0 + v_1)^p$ can be expanded around $u_0$.
     \item $C_{\mathrm{mid}}$:  $(u_0+v_1)^p$ can be expanded around the dominant component $\lambda^{\frac{3}{2}}(R^{-3} (15)^{\frac{3}{2}} + (t\lambda)^{-2} C_1(\nu))$ of $u_0 + v_1$.
     \item $C_{\mathrm{tip}}$:  $(u_0+v_1)^p$ can be expanded around the dominant component $\lambda^{\frac{3}{2}}(t\lambda)^{-2}C_1(\nu)$ of $v_1$.
 \end{enumerate}
 On $C_{\mathrm{mid}}$ and $C_{\mathrm{tip}}$, we are able to terminate the expansions of $(u_0 + v_1)^p$ as each term $(v_1/u_0)^n$ exhibits improved decay in either $R$ or $(t\lambda)$. Near the origin, thanks to $m$ being chosen small enough, we will be able to discard any logarithmic power above some threshold. Indeed, consider the binomial expansion:
 $$
 t^2(u_0 + v_1)^p = t^2u_0^p\left( 1 + \frac{v_1}{u_0} \right) = t^2 u_0^p \sum_{n=0}^{+\infty} \binom{p}{n} \left( \frac{v_1}{u_0} \right)^n,
 $$
 where we recall that $|v_1/u_0| \leq 1/2$ on $C_{\mathrm{ori}}$. We aim to replace $(v_1/u_0)^n$ by a suitable approximation where the logarithmic exponents are bounded.  When $\Re(z) > 0$ and $|z| \geq 2\sqrt{15}$, we have:
 \begin{align*}
     \left( \frac{v_1}{u_0} \right)^n &= \frac{z^{3n}}{(t\lambda)^{2n}} \left(w_0(z^{-1}) + z^{-3}w_1(z^{-1})\log(z)  \right)^n \\
     &= \frac{z^{3n}}{(t\lambda)^{2n}} \sum_{\substack{i+j = n \\ i,j \geq 0}} \binom{n}{i,j} w_0(z^{-1})^i \left[ z^{-3}w_1(z^{-1})\log(z) \right]^j.
 \end{align*}
 Thanks to $m$ being chosen small enough in (\ref{m definition}), the error part 
 $$E_{n} = \frac{z^{3n}}{(t\lambda)^{2n}} \sum_{\substack{i+j = n \\ i,j \geq 0 \\ j \geq \frac{3}{2}N_0}} \binom{n}{i,j} w_0(z^{-1})^i \left[ z^{-3}w_1(z^{-1})\log(z) \right]^j  $$
 can be estimated by
 $$
 \sum_{\substack{i+j = n \\ i,j \geq 0 \\ j \geq \frac{3}{2}N_0}} \binom{n}{i,j} \frac{1}{4^n} \left[ z^{-3} \log(z) \right]^j = \mathcal{O}(2^{-n} z^{-3N_0}) = \mathcal{O}(2^{-n}(t\lambda)^{-N_0})
 $$
 on the region $2m(t\lambda)^{\frac{1}{3}} \leq |z| \leq m(t\lambda)^{\frac{2}{3}}, 0 < t \leq t_0$ when $n \geq N_0$. Similarly, on the region $3\sqrt{15}/2 \leq |z| \leq 2m(t\lambda)^{\frac{1}{3}}, 0 < t \leq t_0$,  we have the estimate
  $$
(t\lambda)^{-n} \sum_{\substack{i+j = n \\ i,j \geq 0 \\ j \geq \frac{3}{2}N_0}} \binom{n}{i,j} \frac{1}{4^n} \left[ z^{-3} \log(z) \right]^j = \mathcal{O}(2^{-n}(t\lambda)^{-N_0}), \quad n \geq N_0,
 $$
Using Cauchy's Integral Formula, for any multi-index $(l_1,l_2) \in \mathbb N^2$, we obtain
  \begin{equation}
      (t\partial_t)^{l_1}\partial_z^{l_2} E_n = \mathcal{O}(2^{-n}(t\lambda)^{-N_0}), \quad n \geq N_0, \label{eq: E_n uniform bound}
  \end{equation}
 on the smaller region $2\sqrt{15} \leq |z| \leq m(t\lambda)^{\frac{1}{3}}, 0 < t \leq t_0$.  Hence, we define the ``truncation'' operator:
 \begin{align}
     T\left[ \left( \frac{v_1}{u_0}\right)^n \right] &= 
 \begin{cases}  \left( \frac{v_1}{u_0}\right)^n &\text{ if } 0 \leq n < N_0 \\
     \left( \frac{v_1}{u_0}\right)^n - \chi_{[2\sqrt{15},+\infty)}(|z|) E_{n} &\text{ otherwise}
 \end{cases} \label{truncation}
 \index{truncation@$T$, the truncation operator}
 \end{align}
 On $4\sqrt{15} = R_0 \leq |z| \leq m(t\lambda)^{\frac{2}{3}}$, the error part $E_n$ is completely removed and the loagrithmic powers are capped to $\log(R)^{3N_0}$. In particular, we have:
    \begin{align*}
        t^2 u_0^p \left( \frac{v_1}{u_0} \right)^{n} &\in \dfrac{\lambda^{\frac{3}{2}}}{(t\lambda)^{2(n-1)}} S^{2n}(R^{-7+3n},\log(R)^{n} ), \\
        t^2 u_0^p  T \left[ \left( \frac{v_1}{u_0} \right)^{n} \right] &\in \dfrac{\lambda^{\frac{3}{2}}}{(t\lambda)^{2(n-1)}} S^{2n}(R^{-7+3n},\log(R)^{3N_0} ). 
    \end{align*}
    As the function $t^2u_0^p (v_1/u_0)^n$ is holomorphic around $|z| \leq \sqrt{15}/2$ and $z \in [\sqrt{15}/2,R_0]$, there exists $C(l_1,l_2) > 1$ for which
    \begin{align*}
        \left[ \dfrac{\lambda^{\frac{3}{2}}}{(t\lambda)^{2(n-1)}} \right]^{-1} \cdot \left| \left| (t\partial_t)^{l_1}( z  \partial_z)^{l_2} t^2 u_0^p \left( \frac{v_1}{u_0} \right)^{n} \right| \right|_{S, \mathrm{ori}} \leq C^n 
    \end{align*}
    using the product rule when $l_2 = 0$ and the complex-differentiability when $l_2 > 0$. Combining this together with the estimate (\ref{eq: E_n uniform bound}), we conclude
    \begin{align*}
        \left[ \dfrac{\lambda^{\frac{3}{2}}}{(t\lambda)^{2(n-1)}} \right]^{-1} \cdot \left| \left| (t\partial_t)^{l_1}( z  \partial_z)^{l_2} t^2 u_0^p  T \left[ \left( \frac{v_1}{u_0} \right)^{n} \right] \right| \right|_{S, \mathrm{ori}} \leq C^n 
    \end{align*}
    since differentiating the cutoff $\chi_{[2\sqrt{15},+\infty)}$ introduces no issue. Near the origin, we allow an exponential growth in $n$, because the prefactor $(t\lambda)^{-2(n-1)}$ cancels this growth on a sufficiently small cone. When $3n - 7  > 0$, we also obtain
\begin{align*}
    \left[ \dfrac{\lambda^{\frac{3}{2}}}{(t\lambda)^{2(n-1)}} \right]^{-1} \cdot \left| \left| (t\partial_t)^{l_1}( z \partial_z)^{l_2}   t^2 u_0^p T \left[ \left( \frac{v_1}{u_0} \right)^{n} \right] \right| \right|_{S,-7+3n,3N_0,\infty} \\
    = \left[ \dfrac{\lambda^{\frac{3}{2}}}{(t\lambda)^{2(n-1)}} \right]^{-1} \cdot \left| \left| (t\partial_t)^{l_1}( z  \partial_z)^{l_2}  t^2 u_0^p \left( \frac{v_1}{u_0} \right)^{n} \right| \right|_{S,-7+3n,3N_0,\infty}  \leq c^n
\end{align*}
    for some $0 < c(l_1,l_2) < 1$. For $l_2 = 0$, this is a consequence of how small $m$ was chosen in (\ref{m definition}) and the product rule. For $l_2 > 0$, this follows from the complex-differentiability. These bounds ensure the convergence of the series
 $$
 t^2(u_0 + v_1)^p \approx t^2 u_0^p \sum_{n=0}^{+\infty} \binom{p}{n} T \left[ \left( \frac{v_1}{u_0} \right)^n \right]
 $$
 and its derivatives. The difference between $t^2(u_0+v_1)^p$ and its truncation, i.e.,
    $$
     t^2(u_0 + v_1)^p - t^2 u_0^p \sum_{n=0}^{+\infty} \binom{p}{n} T \left[ \left( \frac{v_1}{u_0} \right)^n \right]
    $$
    is negligible in the following sense, which makes the truncation a suitable approximation for use in the iteration scheme.
    
\begin{definition}[Negligible terms]\label{negligible simeq definition}
    We say that $f(R,a,t) \in \mathcal{E}_{N_0,\nu}$ if the function $f\left(R,\frac{R}{(t\lambda)},t\right)$ is smooth on $C = \{(R,t): 0 < t < t_0, 0 < R < (t\lambda)\}$ and if for any indices $i,j$,
    $$
    \left|(t \partial_t)^i (\langle R \rangle \partial_R)^j f\left(R,\frac{R}{(t\lambda)},t\right) \right| \lesssim \lambda^{\frac{3}{2}}(t\lambda)^{-N_0} \left[1 + \left(1-\frac{R}{(t\lambda)} \right)^{\frac{1}{2}+\frac{1}{2}\nu-i-j-} \right]
    $$
    on the entire cone. In other words, $f$ has the desired smallness.
    \index{E(N0,nu)@$\mathcal{E}_{N_0,\nu}$, a space of negligible terms}
\end{definition}

As noted earlier,
$$
     t^2(u_0 + v_1)^p - t^2 u_0^p \sum_{n=0}^{+\infty} \binom{p}{n} T \left[ \left( \frac{v_1}{u_0} \right)^n \right] = t^2 u_0^p [1-\chi_{[2\sqrt{15},+\infty)}(|z|)] \sum_{n=0}^{+\infty} \cdot E_{n} \in \mathcal{E}_{N_0,\nu}
$$
due to (\ref{eq: E_n uniform bound}), with no singularity at $R = (t\lambda)$. The use of three different binomial expansions for $t^2(u_0+v_1)^p$  across different regions of the cone motivates the following definitions of correction and error spaces.

\begin{definition}[Correction space $V_k$]
\index{Vk@$V_{2k-1}$, $V_{2k}$, the vector spaces of correction terms}
    Let $V_{2k-1}$, $k \geq 1$, be the vector space spanned by the set of smooth functions $v(R,t) = v_{\mathrm{ori}}(R,t) + v_{\mathrm{mid}}(R,t) + v_{\mathrm{tip}}(R,t) + \eta(R,t)$ on the cone $0 \leq R \leq (t\lambda)$, $0 < t \leq t_0$, where each component admits a decomposition as specified below:
    \begin{enumerate}
    \item $0 \leq R \leq m(t\lambda)^{\frac{2}{3}}$: The component $v_{\mathrm{ori}}$ is given by a convergent sum of the form: $$v_{\mathrm{ori}}(R,t) = \left( \sum_{n = 0}^{+\infty} \ \frac{\lambda^{\frac{3}{2}}}{(t\lambda)^{\alpha+2n}} w_{n}(R) \right) \chi_{[1/m,+\infty)} 
 \left(\frac{(t\lambda)^{\frac{2}{3}}}{R} \right)^{1+j},$$
     where $j \in \mathbb N_{\geq 0}$,  $w_{n} \in S^{2(k-1)}(R^{I+3n}, \log(R)^{J})$ for some common $I,J \in \mathbb N_{\geq 0}$ and $(\alpha,I)$ is $k$-admissible on $C_{\mathrm{ori}}$. The convergence is understood as follows: for any fixed derivative, there exists constants $0 < c(l) < 1 < C(l)$ and $n_0(l) > 0$ for which
    \begin{align*}
    \left| \left| (R\partial_R)^l w_{n} \right| \right|_{S, \mathrm{ori}}  \leq C^n, \quad \left| \left| (R\partial_R)^l w_{n} \right| \right|_{S, I+3n, J, \infty }  \leq c^n
    \end{align*}
    for all $n \geq n_0(l)$.
    \item $\frac{m}{2}(t\lambda)^{\frac{2}{3}} \leq R \leq 2(t\lambda)^{\frac{2}{3}+\varepsilon}$: The component $v_{\mathrm{mid}}$ is a single term of the form:
    $$v_{\mathrm{mid}}(R,t) = \frac{\lambda^{\frac{3}{2}}}{(t\lambda)^{\alpha}}R^i \log(R)^{j_1} \log(t\lambda)^{j_2} g \left( \frac{(t\lambda)^{\frac{2}{3}}}{R} \right) (1-\chi_{[1,+\infty)})\left( \frac{R}{(t\lambda)^{\frac{2}{3}+\varepsilon}} \right)^{1+j_3},$$
where $(\alpha,i)$ is $k$-admissible on $C_{\mathrm{mid}}$, $j_l \in \mathbb N_{\geq 0}$, $g(y)$ is any smooth function on $(0,+\infty)$ which is zero when $y \geq 2/m$ and expands as a finite sum of holomorphic functions and logarithms near $y = 0$.\footnote{More precisely, there should exist some $y_0 > 0$ and finitely many holomorphic functions $g_0,g_1,...$ around $|y| \leq y_0$ for which we can write $g(y) = \sum_{j=0}^{J}g_j(y) \log(y)^j$ whenever $|y| \leq y_0$}  

We note that we can always build a basis of our vector space by omitting the powers of $\log(t\lambda)$ in $v_{\mathrm{mid}}$ as they can always be rewritten as $\frac{3}{2}\left( \log(y) - \log(R) \right)$, $y = (t\lambda)^{\frac{2}{3}}/R$, and the $\log(y)$ part can be included in the $g(y)$-type functions.
\item $(t\lambda)^{\frac{2}{3}+\varepsilon} \leq R \leq (t\lambda)$: The component $v_{\mathrm{tip}}$ is a single term of the form:
    $$v_{\mathrm{tip}}(R,t) = \frac{\lambda^{\frac{3}{2}}}{(t\lambda)^{\alpha}}R^i \log(R)^{j_1} \log(t\lambda)^{j_2} h \left( \frac{R}{(t\lambda)^{\frac{2}{3}+\varepsilon}} \right),$$
where $(\alpha,i)$ is $k$-admissible on $C_{\mathrm{tip}}$, $j_l \in \mathbb N_{\geq 0}$ and $h(y)$ is any smooth function which is constant outside $[1,2]$ and zero when $y$ is in a neighbourhood of $1$.

\item $\eta \in \mathcal{E}_{N_0,\nu}$ and $\eta$ has no singularity at the boundary $R = (t\lambda)$. 
\end{enumerate}

Similarly, let $V_{2k}$, $k \geq 1$, be the vector space generated by smooth functions $v(R,t)$ on the cone $0 \leq R < (t\lambda)$, $0 < t \leq t_0$, whose derivatives up to order $\floor{\left(\frac{1}{2}+\frac{1}{2}\nu\right)-}$ are continuous at the boundary $R = (t\lambda)$, and which are given by a finite sum of the form:
    $$v(R,t) = \frac{\lambda^{\frac{3}{2}}}{(t\lambda)^{\alpha}}R^i \log(R)^{j_1} \log(t\lambda)^{j_2} q \left( \frac{R}{(t\lambda)} \right)  \chi_{[1,+\infty)}\left(  \frac{R}{(t\lambda)^{\frac{2}{3}+\varepsilon}} \right)^{1+j_3},$$
where $(\alpha,i)$ is $k$-admissible on $C_{\mathrm{tip}}$, $j_l \in \mathbb N_{\geq 0}$ and $q(a) \in a^2\mathcal{Q}_{k,\frac{1}{2}+\frac{1}{2}\nu}$. We note that for these functions, the $\log(t\lambda)$ powers can always be rewritten so that they do not appear in the finite sum.
\end{definition}

As a consistency check, we verify that the first correction term $v_1(R,t)$ belongs to $V_1$. If $V_1(R)$ is defined as in (\ref{v_1 formula}), then
\begin{equation}
     v_1(R,t):= \frac{\lambda^{\frac{3}{2}}}{(t\lambda)^{2}} V_1(R) \in V_1  \label{eq: v_1 belongs to V_1}
\end{equation}
as a consequence of the following proposition applied for $k = 1$, $(\alpha,I) = (2,0)$.

\begin{proposition}[Examples of $V_{2k-1}$ functions]\label{inclusion S^2n into V_2k-1}
    Let $(\alpha,I)$ be $k$-admissible on $C_{\mathrm{mid}}$ and $C_{\mathrm{tip}}$, $k \geq 1$. For $J \in \mathbb N_{\geq 0}$, 
    $$
    \frac{\lambda^{\frac{3}{2}}}{(t\lambda)^{\alpha}}S^0(R^I, \log(R)^J) \cdot (1- \chi_{[1/m,+\infty)}) 
 \left(\frac{(t\lambda)^{\frac{2}{3}}}{R} \right) \subset V_{2k-1}.
    $$
    Furthermore, if $(\alpha,I)$ is also $k$-admissible on $C_{\mathrm{ori}}$, then
    $$
\frac{\lambda^{\frac{3}{2}}}{(t\lambda)^{\alpha}}S^0(R^I, \log(R)^J) \subset  V_{2k-1}.
    $$
\end{proposition}

\begin{proof}
    The second statement follows from the first one by writing
    \begin{align*}
            \frac{\lambda^{\frac{3}{2}}}{(t\lambda)^{\alpha}}S^0(R^I, \log(R)^J) &= \frac{\lambda^{\frac{3}{2}}}{(t\lambda)^{\alpha}}S^0(R^I, \log(R)^J)  \cdot \chi_{[1/m,+\infty)} \\
            &\phantom{=}+\frac{\lambda^{\frac{3}{2}}}{(t\lambda)^{\alpha}}S^0(R^I, \log(R)^J)  \cdot (1- \chi_{[1/m,+\infty)}). 
    \end{align*}
    Let $V(R) \in \frac{\lambda^{\frac{3}{2}}}{(t\lambda)^{\alpha}} S^0(R^I, \log(R)^J)$. On $R \geq R_0$, we can write
    $$
    V(R) = R^I \sum_{j=0}^{J}w_{j}(R^{-1})\log(R)^j,
    $$
    where $w_{j}(y)$ is holomorphic on $|y| < R_0^{-1} + \delta$ for some $0 < \delta \ll 1$. Now, rewrite $V(R)$ as
    $$
    V(R) = \sum_{i=0}^{|I|+N_0} \sum_{j=0}^J c_{k,j}R^{I-k} \log(R)^j +  \sum_{j=0}^{J}\tilde{w}_{j}(R^{-1})\log(R)^j,
    $$
    where each $\tilde{w}_j(y)$ is holomorphic with a zero of order at least $y^{N_0+1}$ at $y = 0$. Consider a smooth cutoff function $\chi(y)$ such that $\chi(y) = 1$ for $|y| < R_0^{-1} + \delta/2$ and $\chi(y) = 0$ for $|y| \geq R_0^{-1} + \delta$. Then, we have
    $$
    V(R) = \sum_{i=0}^{|I|+N_0} \sum_{j=0}^J c_{k,j} R^{I-k} \log(R)^j +  \chi(R^{-1}) \sum_{j=0}^{J}\tilde{w}_{j}(R^{-1})\log(R)^j, \quad R \gtrsim (t\lambda)^{\frac{2}{3}} \gg R_0.
    $$
    Therefore, on $R \gtrsim (t\lambda)^{\frac{2}{3}}$, we obtain
    $$
    V(R) = \sum_{i=0}^{|I|+N_0} \sum_{j=0}^J c_{k,j} R^{I-k} \log(R)^j +  S^0(R^{-N_0-1}, \log(R)^J) =: \tilde{V}(R) +  S^0(R^{-N_0-1}, \log(R)^J).
    $$
    Consider the expression:
    $$
\frac{\lambda^{\frac{3}{2}}}{(t\lambda)^{\alpha}} S^0(R^{-N_0-1}, \log(R)^J)  \cdot (1- \chi_{[1/m,+\infty)}) 
 \left(\frac{(t\lambda)^{\frac{2}{3}}}{R} \right).
    $$
    This expression is in $\mathcal{E}_{N_0,\nu}$ with no singularity at $a = 1$ since no $\mathcal{Q}$ functions is involved. Indeed, this remainder is sufficiently small in a pointwise sense. Moreover, for all indices $l_1,l_2 \geq 0$, we have:
\begin{align}
    (t\partial_t)^{l_1} (R\partial_R)^{l_2} \frac{\lambda^{\frac{3}{2}}}{(t\lambda)^{\alpha}}R^i \log(R)^{j_1} \log(t\lambda)^{j_2} &= c_{l_1,l_2,\alpha,i,j_2,j_2} \frac{\lambda^{\frac{3}{2}}}{(t\lambda)^{\alpha}}R^i \log(R)^{j_1} \log(t\lambda)^{j_2} \notag\\
    (t\partial_t)^{l_1} (R\partial_R)^{l_2} S^{2n}(R^I,\log(R)^J)&\subset S^{2n}(R^I,\log(R)^J) \notag \\
    (t\partial_t)^{l_1} (R\partial_R)^{l_2} g \left( \frac{(t\lambda)^{\frac{2}{3}}}{R} \right) &= (-1)^{l_1+l_2}\left[\frac{2 \nu}{3} \right]^{l_1} \left[ (y\partial_y)^{l_1+l_2}g \right] \left( \frac{(t\lambda)^{\frac{2}{3}}}{R} \right) \notag \\
     (t\partial_t)^{l_1} (R\partial_R)^{l_2} h \left(  \frac{R}{(t\lambda)^{\frac{2}{3}+\varepsilon}} \right) &= \left[\nu \left(\frac{2}{3}+\varepsilon \right) \right]^{l_1} \left[(y\partial_y)^{l_1+l_2}h \right] \left(  \frac{R}{(t\lambda)^{\frac{2}{3}+\varepsilon}} \right), \label{eq: derivative t R in terms of a, y}
\end{align}
so that the smallness of the expression is preserved under differentiation with $(t\partial_t)^{l_1} (R\partial_R)^{l_2}$. 

Finally, we consider the finite sum $\tilde{V}(R)$. We extract a $v_{\mathrm{tip}}$ and $v_{\mathrm{mid}}$ component by decomposing:
\begin{align*}
   \frac{\lambda^{\frac{3}{2}}}{(t\lambda)^{\alpha}}  \tilde{V}(R)  \cdot (1- \chi_{[1/m,+\infty)}) 
 \left(\frac{(t\lambda)^{\frac{2}{3}}}{R} \right)  =\frac{\lambda^{\frac{3}{2}}}{(t\lambda)^{\alpha}}  \tilde{V}(R)  \cdot \chi_{[1,+\infty)}\left( \frac{R}{(t\lambda)^{\frac{2}{3}+\varepsilon}} \right) \\
+\frac{\lambda^{\frac{3}{2}}}{(t\lambda)^{\alpha}}  \tilde{V}(R)  \cdot (1- \chi_{[1/m,+\infty)}) 
 \left(\frac{(t\lambda)^{\frac{2}{3}}}{R} \right) \cdot (1-\chi_{[1,+\infty)})\left( \frac{R}{(t\lambda)^{\frac{2}{3}+\varepsilon}} \right).
\end{align*}
\end{proof}

We now define the error spaces, whose structure is similar to that of the correction spaces.
\begin{definition}[Error space $E_k$]
\index{Ek@$E_{\mathrm{ori},k}$, $E_{\mathrm{tip},k}$, the vector spaces of error terms}
    Let $E_{\mathrm{ori},k}$, $k \geq 1$, be the vector space spanned by the set of smooth functions $e(R,t) = e_{\mathrm{ori}}(R,t) + e_{\mathrm{mid}}(R,t) + e_{\mathrm{tip}}(R,t) + \eta(R,t)$ on the cone $0 \leq R \leq (t\lambda)$, $0 < t \leq t_0$, where each component admits a decomposition as specified below:
    \begin{enumerate}
    \item $0 \leq R \leq m(t\lambda)^{\frac{2}{3}}$: The component $e_{\mathrm{ori}}$ is given by a convergent sum of the form: $$e_{\mathrm{ori}}(R,t) = \left( \sum_{n = 0}^{+\infty} \ \frac{\lambda^{\frac{3}{2}}}{(t\lambda)^{\alpha+2n}} w_{n}(R) \right) \chi_{[1/m,+\infty)} 
 \left(\frac{(t\lambda)^{\frac{2}{3}}}{R} \right)^{1+j},$$
     where $j \in \mathbb N_{\geq 0}$,  $w_{n} \in S^{2(k-1)}(R^{I+3n}, \log(R)^{J})$ for some common $I,J \in \mathbb N_{\geq 0}$ and $(\alpha,I)$ is $k$-admissible on $C_{\mathrm{ori}}$. Moreover, for any fixed derivative, there exists constants $0 < c(l) < 1 < C(l)$ and $n_0(l) > 0$ for which
    \begin{align*}
    \left| \left| (R\partial_R)^l w_{n} \right| \right|_{S, \mathrm{ori}}  \leq C^n, \quad \left| \left| (R\partial_R)^l w_{n} \right| \right|_{S, I+3n, J, \infty }  \leq c^n
    \end{align*}
    for all $n \geq n_0(l)$. 
    \item $\frac{m}{2}(t\lambda)^{\frac{2}{3}} \leq R \leq 2(t\lambda)^{\frac{2}{3}+\varepsilon}$: The component $e_{\mathrm{mid}}$ is given by a single term of the form: 
    $$e_{\mathrm{mid}}(R,t) = \frac{\lambda^{\frac{3}{2}}}{(t\lambda)^{\alpha}}R^i \log(R)^{j_1} \log(t\lambda)^{j_2} g \left( \frac{(t\lambda)^{\frac{2}{3}}}{R} \right) (1-\chi_{[1,+\infty)})\left( \frac{R}{(t\lambda)^{\frac{2}{3}+\varepsilon}} \right)^{1+j_3},$$
where $(\alpha,i)$ is $k$-admissible on $C_{\mathrm{mid}}$, $j_l \in \mathbb N_{\geq 0}$, $g(y)$ is any smooth function on $(0,+\infty)$ which is zero when $y \geq 2/m$ and expands as a finite sum of holomorphic functions and logarithms near $y = 0$.\footnote{More precisely, there should exist some $y_0 > 0$ and a finite number of holomorphic functions $g_1,g_2,...$ around $|y| \leq y_0$ for which we can write $g(y) = \sum_{j=0}^{J}g_j(y) \log(y)^j$ whenever $|y| \leq y_0$}  

We note that we can always build a basis of our vector space by omitting the powers of $\log(t\lambda)$ in $v_{\mathrm{mid}}$ as they can always be rewritten as $\frac{3}{2}\left( \log(y) - \log(R) \right)$, $y = (t\lambda)^{\frac{2}{3}}/R$, and the $\log(y)$ part can be included in the $g(y)$-type functions.
\item $(t\lambda)^{\frac{2}{3}+\varepsilon} \leq R \leq (t\lambda)$: The component $e_{\mathrm{tip}}$ is given by a single term of the form: 
    $$e_{\mathrm{tip}}(R,t) = \frac{\lambda^{\frac{3}{2}}}{(t\lambda)^{\alpha}}R^i \log(R)^{j_1} \log(t\lambda)^{j_2} h \left( \frac{R}{(t\lambda)^{\frac{2}{3}+\varepsilon}} \right),$$
where $(\alpha,i)$ is $k$-admissible on $C_{\frac{2}{3}+\varepsilon}$, $j_l \in \mathbb N_{\geq 0}$ and $h(y)$ is any smooth function with compact support in $(1,2)$. We note that the $\log(t\lambda)$ powers can always be rewritten so that they do not appear in the finite sum for $e_{\mathrm{tip}}$.

\item $\eta \in \mathcal{E}_{N_0,\nu}$ and $\eta$ has no singularity at the boundary $R = (t\lambda)$. 
\end{enumerate}

The error space $E_{\mathrm{tip},k}$ is defined as the vector space of functions $E_{\mathrm{tip},k} = \frac{1}{a^2} V_{2k} + \mathcal{E}_{N_0,\nu}$. 
\end{definition}

 \begin{notation}
The symbol
$$
\fsum_{\alpha \: \ ...}
$$
indicates summation over a finite set of indices $\alpha$ satisfying a certain relation. If the exact set of indices does not matter, such a notation will be used. It will mostly be used when considering a sum over a finite set of $k$-admissible pairs $(\alpha,i)$, which are not explicitly defined.
\index{finite sum@$\fsum$, any kind of finite sum}
\end{notation}

\begin{remark} We make a few remarks regarding the definitions of the correction and error spaces.
    \begin{enumerate}
        \item To simplify the definitions of $V_{2k-1}$, $V_{2k}$, $E_{\mathrm{ori},k}$ and $E_{\mathrm{tip},k}$, we gave the definition in terms of a basis. That means that the correction terms and error terms on each region are always given by a finite sum of such elements. For example, a term $v_{\mathrm{mid}}$ has the form:
        $$\fsum_{\substack{(\alpha,i) \kadm \\ j_1,j_2,j_3 \geq 0}}\frac{\lambda^{\frac{3}{2}}}{(t\lambda)^{\alpha}}R^i \log(R)^{j_1} \log(t\lambda)^{j_2} g_{\alpha,i,j_1,j_2,j_3} \left( \frac{(t\lambda)^{\frac{2}{3}}}{R} \right) (1-\chi_{[1,+\infty)})\left( \frac{R}{(t\lambda)^{\frac{2}{3}+\varepsilon}} \right)^{1+j_3}.$$
        The cutoff functions are fixed, but the functions $w_n(R) = w_{n}^{\alpha,I,J,j}(R)$, \sloppy $g(y) = g_{\alpha,i,j_1,j_2,j_3}(y)$, $h(y) = h_{\alpha,i,j_1,j_2}(y)$ and $q(a) = q_{\alpha,i,j_1,j_2}(a)$ can be anything with the desired properties.
        \item The only singularities that can happen are at $a = 1$, when considering functions from $V_{2k}$ and $E_{\mathrm{tip},k}$. In particular, restricting such a function near the origin or the middle region of the cone removes the singularity.
        \item The main difference between $V_{2k-1}$ and $E_{\mathrm{ori},k}$ lies near the tip of the cone $(t\lambda)^{\frac{2}{3}+\varepsilon} \leq R \leq (t\lambda)$. Up to a negligible part, error terms are supported on $(t\lambda)^{\frac{2}{3}+\varepsilon} \leq R \leq 2(t\lambda)^{\frac{2}{3}+\varepsilon}$, while correction terms are supported on the whole cone. 
        \item For the error terms, we write $e_{k} \simeq \tilde{e}_k$ if $e_k - \tilde{e}_k \in \mathcal{E}_{N_0,\nu}$. In other words, the equality $e_k(R,a,t) = \tilde{e}_k(R,a,t)$ holds up to negligible terms, which do not affect the subsequent analysis. This negligible difference is carried over to all subsequent error terms $(e_{k+1},e_{k+2})$, but it does not affect the algorithm, so we will omit them when describing the next error terms. It is important to note that in the renormalization step, we never differentiate an error term. Hence, the singularity of an $E_{\mathrm{tip},k}$ element at $a = 1$ cannot get worse. We do differentiate even correction terms $v_{2k}$, but in that case, the support is always restricted to the middle or the origin part of the cone, which removes the singularity.
    \end{enumerate}
\end{remark}

We conclude this section by proving some useful computation rules concerning these correction and error spaces. 

\begin{proposition}[Stability under differentiation]\label{prop:stability under differentiation}
For any $l_1,l_2,k \in \mathbb N$, 
$$(t \partial_t)^{l_1}(R\partial_t)^{l_2} V_{2k-1} \subset V_{2k-1}, \quad (t \partial_t)^{l_1}(R\partial_t)^{l_2} E_{\mathrm{ori},k} \subset E_{\mathrm{ori},k}.$$
\end{proposition}

\begin{proof}
    We prove the inclusion $(t \partial_t)^{l_1}(R\partial_t)^{l_2} V_{2k-1} \subset V_{2k-1}$. Recall first the identities (\ref{eq: derivative t R in terms of a, y}) which hold for any exponents and any smooth functions $g$ and $h$. 
\item[\textbf{Tip part:}] Suppose 
$$v_{\mathrm{tip}}(R,t) = \frac{\lambda^{\frac{3}{2}}}{(t\lambda)^{\alpha}}R^i \log(R)^{j_1} \log(t\lambda)^{j_2} h \left( \frac{R}{(t\lambda)^{\frac{2}{3}+\varepsilon}} \right).$$
Then, applying the derivatives yields:
$$
(t \partial_t)^{l_1}(R\partial_t)^{l_2} v_{\mathrm{tip}} = \fsum_{k_1,k_2 \geq 0} c_{k_1,k_2}\frac{\lambda^{\frac{3}{2}}}{(t\lambda)^{\alpha}}R^i \log(R)^{j_1} \log(t\lambda)^{j_2} \left[ (y\partial y)^{k_1+k_2} h \right] \left( \frac{R}{(t\lambda)^{\frac{2}{3}+\varepsilon}} \right),
$$
which is a finite sum of elements $\tilde{v}_{\mathrm{tip},k_1,k_2}$ having a suitable form. 

\item[\textbf{Middle part:}] Similarly, suppose 
$$v_{\mathrm{mid}}(R,t) = \frac{\lambda^{\frac{3}{2}}}{(t\lambda)^{\alpha}}R^i \log(R)^{j_1} \log(t\lambda)^{j_2} g \left( \frac{(t\lambda)^{\frac{2}{3}}}{R} \right) (1-\chi_{[1,+\infty)})\left( \frac{R}{(t\lambda)^{\frac{2}{3}+\varepsilon}} \right)^{1+j_3}.$$
Then: 
\begin{align*}
    (t \partial_t)^{l_1}(R\partial_t)^{l_2}v_{\mathrm{mid}}(R,t) = (1-\chi_{[1,+\infty)})^{1+j_3} (t \partial_t)^{l_1}(R\partial_t)^{l_2} \frac{\lambda^{\frac{3}{2}}}{(t\lambda)^{\alpha}}R^i \log(R)^{j_1} \log(t\lambda)^{j_2} g + \\ 
\phantom{=} \sum_{l = 1}^{l_1+l_2} \left[ (y \partial_y)^l (1-\chi_{[1,+\infty)})^{1+j_3} \right] \cdot \fsum_{k_1,k_2 \geq 0} c_{l,k_1,k_2} (t \partial_t)^{k_1}(R\partial_t)^{k_2} \frac{\lambda^{\frac{3}{2}}}{(t\lambda)^{\alpha}}R^i \log(R)^{j_1} \log(t\lambda)^{j_2} g.
\end{align*}
The first part is treated as for $v_{\mathrm{tip}}$. It produces a finite sum of $\tilde{v}_{\mathrm{mid}}$ elements. As for the second part, it yields a finite sum of elements of the form:
$$
\tilde{v}_{\mathrm{mid}} \cdot h \left( \frac{R}{(t\lambda)^{\frac{2}{3}+\varepsilon}} \right) := \frac{\lambda^{\frac{3}{2}}}{(t\lambda)^{\alpha}}R^i \log(R)^{j_1} \log(t\lambda)^{j_2} g \left( \frac{(t\lambda)^{\frac{2}{3}}}{R} \right) h \left( \frac{R}{(t\lambda)^{\frac{2}{3}+\varepsilon}} \right), 
$$
where $h$ is smooth with compact support in $(1,2)$. Then $R \sim (t\lambda)^{\frac{2}{3}+\varepsilon}$ and $(t\lambda)^{\frac{2}{3}}R^{-1} \sim (t\lambda)^{-\varepsilon} \sim 0$. Near $y = 0$, one can expand
$$
g(y) = g_0(y) + \sum_{j=1}^{L}g_{j}(y) \log(y)^j
$$
for some functions $g_j$ holomorphic around zero. Consider the $M$-th order Taylor polynomial $P_j(y)$ of $g_j(y)$ near $y = 0$ with $M = \ceil{N_0\varepsilon^{-1}}$, meaning that the remainder $\eta_j(y) = g_j(y) - P_j(y)$ is a holomorphic function around zero with $(y\partial_y)^l [g_j(y) - P_j(y)] = \mathcal{O}(y^{M})$ for any fixed $l \geq 0$. Then, we obtain
$$
h \cdot \tilde{v}_{\mathrm{mid}}(R,t) = 
   \frac{\lambda^{\frac{3}{2}}}{(t\lambda)^{\alpha}}R^i \log(R)^{j_1} \log(t\lambda)^{j_2} \left(P_0(y) + \sum_{j=1}^{L}P_{j}(y) \log(y)^j \right) h \left(  \frac{R}{(t\lambda)^{\frac{2}{3}+\varepsilon}} \right)+ \eta 
$$
for $y = (t\lambda)^{\frac{2}{3}}R^{-1}$. The remainder
\begin{align}
    \eta =  \frac{\lambda^{\frac{3}{2}}}{(t\lambda)^{\alpha}}R^i \log(R)^{j_1} \log(t\lambda)^{j_2} \left(\eta_0(y) + \sum_{j=1}^{L}\eta_{j}(y) \log(y)^j \right) h \left(  \frac{R}{(t\lambda)^{\frac{2}{3}+\varepsilon}} \right) \label{eq:remainder on mid/tip is negligible}
\end{align}
is supported on $R \sim (t\lambda)^{\frac{2}{3}+\varepsilon}$, so it has no singularity at the tip of the cone.  We verify that $\eta \in \mathcal{E}_{N_0,\nu}$. The $\log(y)$ rewrites as $-\log(R) + \frac{2}{3}\log(t\lambda)$, so assume without loss of generality that we have only one term:
$$
\eta =  \frac{\lambda^{\frac{3}{2}}}{(t\lambda)^{\alpha}}R^i \log(R)^{j_1} \log(t\lambda)^{j_2} \eta_0(y)h \left(  \frac{R}{(t\lambda)^{\frac{2}{3}+\varepsilon}} \right).
$$
Given that  $(\alpha,i)$ is $k$-admissible and $y^M \sim (t\lambda)^{-N_0}$ on $R \sim (t\lambda)^{\frac{2}{3}+\varepsilon}$, the remainder have the desired smallness. Furthermore, the smallness is left unchanged under $(t \partial_t)^{l_1}(R\partial_t)^{l_2}$, as a consequence of the equalities (\ref{eq: derivative t R in terms of a, y}).

\item[\textbf{Origin part:}] Finally, consider a term of the form:
$$
v_{\mathrm{ori}}(R,t) = \left( \sum_{n = 0}^{+\infty} \ \frac{\lambda^{\frac{3}{2}}}{(t\lambda)^{\alpha+2n}} w_{n}(R) \right) f
 \left(\frac{(t\lambda)^{\frac{2}{3}}}{R} \right),
 $$
     where $f = \chi_{[1/m,+\infty)}^{1+j}$, $w_{n} \in S^{2(k-1)}(R^{I+3n}, \log(R)^{J})$ for some common $I,J \in \mathbb N_{\geq 0}$ with $(\alpha,I)$ being $k$-admissible on $C_{\mathrm{ori}}$. One has
         \begin{align*}
     (R \partial_R) v_{\mathrm{ori}}
 &=  f
 \left(\frac{(t\lambda)^{\frac{2}{3}}}{R} \right) \cdot    \sum_{n = 0}^{+\infty}  \frac{\lambda^{\frac{3}{2}}}{(t\lambda)^{\alpha+2n}} (R\partial_R)w_n(R) \\
 &\phantom{=}+  \sum_{n = 0}^{+\infty} \frac{\lambda^{\frac{3}{2}}}{(t\lambda)^{\alpha+2n}} w_{n}(R) \cdot f'
 \left(\frac{(t\lambda)^{\frac{2}{3}}}{R} \right) \cdot \left(\frac{(t\lambda)^{\frac{2}{3}}}{R} \right) \\
 &= f
 \left(\frac{(t\lambda)^{\frac{2}{3}}}{R} \right) \cdot    \sum_{n = 0}^{+\infty}  \frac{\lambda^{\frac{3}{2}}}{(t\lambda)^{\alpha+2n}} (R\partial_R)w_n(R) \\
 &\phantom{=}+  \tilde{f}  \left(\frac{(t\lambda)^{\frac{2}{3}}}{R} \right) \sum_{n = 0}^{+\infty} \frac{\lambda^{\frac{3}{2}}}{(t\lambda)^{\alpha+2n}} w_{n}(R),
    \end{align*}
    where $\tilde{w}_{n} \in S^{2(k-1)}(R^{I+3n}, \log(R)^{J})$  and $(\alpha,I)$ remains $k$-admissible. A similar equality holds for $(t\partial_t) v_{\mathrm{ori}}$ and, by induction, 
   \begin{align*}
       (t \partial_t)^{l_1}(R \partial_R)^{l_2} v_{\mathrm{ori}} &= 
      f  \left(\frac{(t\lambda)^{\frac{2}{3}}}{R} \right) \sum_{n = 0}^{+\infty} (t \partial_t)^{l_1}\frac{\lambda^{\frac{3}{2}}}{(t\lambda)^{\alpha+2n}} (R \partial_R)^{l_2} w_{n}(R) 
       \\ 
&\phantom{=}+\sum_{j=1}^M \tilde{f}_j  \left(\frac{(t\lambda)^{\frac{2}{3}}}{R} \right) \sum_{n = 0}^{+\infty} \frac{\lambda^{\frac{3}{2}}}{(t\lambda)^{\alpha+2n}} w_{n,j}(R),
   \end{align*}
    where $M = M(l_1,l_2) < +\infty$ and each $\tilde{f}_j$ is compactly supported on $(1/m,2/m)$. To conclude, it suffices to show that a term of the form: 
    \begin{align}
\tilde{f} \left(\frac{(t\lambda)^{\frac{2}{3}}}{R} \right) \sum_{n = 0}^{+\infty} \frac{\lambda^{\frac{3}{2}}}{(t\lambda)^{\alpha+2n}} w_{n}(R) \label{eq:approximation of restricted ori term}
    \end{align}
can be approximated by a finite sum if $\tilde{f}$ is compactly supported, so that they yield a finite number of $\tilde{v}_{\mathrm{mid}}$ components. Such a term is supported on $\frac{m}{2}(t\lambda)^{\frac{2}{3}} \leq R \leq m(t\lambda)^{\frac{2}{3}}$. On this region, 
     $$
     w_n(z) =  z^{I+3n} \sum_{j=0}^{J}w_{n,j}(z^{-1})
 \log(z)^j
     $$
     and $w_{n}(z)$ is well-approximated by a Taylor polynomial of degree $3N_0 + 1$,
$$
\tilde{w}_{n,j}(z) =   \sum_{l=0}^{3N_0 + 1} w_{n,j,l}z^{l}.$$
We check that summing all the errors that we do with this approximation is negligible in the sense that it belongs to $\mathcal{E}_{N_0,\nu}$ (and it has no singularity at the tip of the cone as the error is supported on $R \sim (t\lambda)^{\frac{2}{3}}$). Recall that:
$||(z\partial_z)^{l}w_{n,j}||_{A(|y| \leq R_0^{-1})} \leq m^{-3n}c_l^n$
for some $0 < c_l < 1$ and all $n \geq n_l$. As $(z\partial_z)^l$ is a linear combination of $z\partial_z$, $z^2\partial_z^2$, ..., $z^l \partial_z^l$, it follows by induction that $$
R_0^{-l} ||\partial_z^l w_{n,j}||_{A(|y| \leq R_0^{-1})}  = ||z^l\partial_z^{l}w_{n,j}||_{A(|y| \leq R_0^{-1})}  \leq m^{-3n}\tilde{c}_l^n
$$
for some $0 < \tilde{c}_l < 1$ and all $n \geq \tilde{n}_l$. Working with $W_{n,j,k} = (z\partial_z)^k w_{n,j}$, one similarly deduces
$$
||\partial_z^l W_{n,j,k}||_{A(|y| \leq R_0^{-1})} = ||\partial_z^l (z\partial_z)^k w_{n,j}||_{A(|y| \leq R_0^{-1})} \leq R_0^{l}m^{-3n}\tilde{c}_{l,k}^n
$$ 
for some $0 < \tilde{c}_{l,k} < 1$ and for all $n \geq \tilde{n}_{l,k}$. Then $(z\partial_z)^{l}\tilde{w}_{n,j}$ is the Taylor polynomial of degree $3N_0 + 1$ of $(z\partial_z)^{l}w_{n,j}$ and
$$
(z\partial_z)^{l}(w_{n,j} - \tilde{w}_{n,j}) = \int_{[0,z]} \frac{\partial_y^{(3N_0+2)} (y\partial_y)^l w_{n,j}(y)}{(3N_0+2)!}(z-y)^{3N_0+2}dy = \mathcal{O}(m^{-3n}c_l^{n}z^{3N_0+2})
$$
for some $0 < c_l < 1$ and all $n$ large enough. On $\frac{m}{2}(t\lambda)^{\frac{2}{3}} \leq R = z^{-1} \leq m(t\lambda)^{\frac{2}{3}}$,
$$
(z\partial_z)^{l}(w_{n,j} - \tilde{w}_{n,j}) = \mathcal{O}\left( c_l^n (t\lambda)^{-N_0-2/3} \right)
$$
The final error in approximating the expression
    $$
\tilde{f} \left(\frac{(t\lambda)^{\frac{2}{3}}}{R} \right) \sum_{n = 0}^{+\infty} \frac{\lambda^{\frac{3}{2}}}{(t\lambda)^{\alpha+2n}} w_{n}(R) 
$$ 
is given by
    $$
\tilde{f} \left(\frac{(t\lambda)^{\frac{2}{3}}}{R} \right) \sum_{n = 0}^{+\infty} \frac{\lambda^{\frac{3}{2}}}{(t\lambda)^{\alpha+2n}} z^{i+3n} \sum_{j=0}^{J}\left( w_{n,j}(z^{-1}) - \tilde{w}_{n,j}(z^{-1}) \right) \log(z)^j.
$$ 
Given the compact support of $\tilde{f}$, the estimate on $(z\partial_z)^{l}(w_{n,j} - \tilde{w}_{n,j})$ and the admissibility of the pair $(\alpha,I)$, the desired smallness is achieved.
\end{proof}

For $V_{2k}$ and $E_{\mathrm{ori},k}$, such stability under differentiation does not hold in general, because if one applies the operator $(z\partial_z)^l$, then one gets an element of similar form but the $q(a)$ coefficient belongs to $\mathcal{Q}_{\frac{1}{2}+\frac{1}{2}\nu-l}$ instead of $\mathcal{Q}_{\frac{1}{2}+\frac{1}{2}\nu}$. However, smallness is preserved in the following sense:

\begin{proposition}[Smallness is preserved under differentiation]\label{prop:stability of smallness under differenatiation}
Let $v(R,t) \in V_{2k} \cup E_{\mathrm{tip},k}$. For any $l_1,l_2,k \in \mathbb N$, it holds that
$$\left|(t\partial_t)^{l_1}(R\partial_R)^{l_2}v \right| \lesssim \frac{\lambda^{\frac{3}{2}}}{(t\lambda)^{2+\left( \frac{2}{3}-2\varepsilon \right) \cdot (k-1)}}\left[1 + \left(1-\frac{R}{(t\lambda)} \right)^{\frac{1}{2}+\frac{1}{2}\nu-l_1-l_2-} \right].
    $$
\end{proposition}

\begin{proof}
    The proof is similar to that of Proposition \ref{prop:stability under differentiation}. The singularity comes from differentiating the $\mathcal{Q}_{\frac{1}{2}+\frac{1}{2}\nu}$ elements.
\end{proof}

\begin{proposition}[Product Rules] \label{product rules for V, E}
   Let $k,k_1,k_2 \geq 1$, $g(y)$ be any smooth function on $(0,+\infty)$ which is zero when $y \geq 2/m$ and expands as a finite sum of holomorphic functions and logarithms near $y = 0$, and let $h(y)$ be any smooth function with support in $(-\infty,2)$. The following product rules hold:
    \begin{enumerate}
        \item 
        $$V_{2k} \subset E_{\mathrm{tip},k}, \quad V_{2k-1} \subset E_{\mathrm{ori},k} + E_{\mathrm{tip},k}.$$ \label{prop:product rule 1}
        \item 
        $$g \left( \frac{(t\lambda)^{\frac{2}{3}}}{R}  \right) (V_{2k-1} \cup E_{\mathrm{ori},k}) \subset E_{\mathrm{ori},k}, \quad h \left(  \frac{R}{(t\lambda)^{\frac{2}{3}+\varepsilon}} \right) (V_{2k} \cup E_{\mathrm{tip},k}) \subset E_{\mathrm{ori},k}.$$ \label{prop:product rule 2}
        \item
        \begin{align*}
         \left[\frac{\lambda^{\frac{3}{2}}}{(t\lambda)^{2}} \right]^{-1} ( V_{2k_1-1} \cup E_{\mathrm{ori},k_1} ) \cdot  \left[\frac{\lambda^{\frac{3}{2}}}{(t\lambda)^{2}} \right]^{-1} E_{\mathrm{ori},k_2} &\subset \left[\frac{\lambda^{\frac{3}{2}}}{(t\lambda)^{2}} \right]^{-1} E_{\mathrm{ori},k_1+k_2},\\
            \left[\frac{\lambda^{\frac{3}{2}}}{(t\lambda)^{2}} \right]^{-1} ( V_{2k_1} \cup E_{\mathrm{tip},k_1} ) \cdot \left[\frac{\lambda^{\frac{3}{2}}}{(t\lambda)^{2}} \right]^{-1} E_{\mathrm{tip},k_2} &\subset \left[\frac{\lambda^{\frac{3}{2}}}{(t\lambda)^{2}} \right]^{-1} E_{\mathrm{tip},k_1+k_2}, \\
             \left[\frac{\lambda^{\frac{3}{2}}}{(t\lambda)^{2}} \right]^{-1} ( V_{2k_1-1} \cup E_{\mathrm{ori},k_1} ) \cdot \left[\frac{\lambda^{\frac{3}{2}}}{(t\lambda)^{2}} \right]^{-1} E_{\mathrm{tip},k_2} &\subset \left[\frac{\lambda^{\frac{3}{2}}}{(t\lambda)^{2}} \right]^{-1} E_{\mathrm{ori},k_1+k_2}.
        \end{align*} \label{prop:product rule 3}
    \end{enumerate}
\end{proposition}

\begin{remark}
    From (\ref{prop:product rule 1}) and (\ref{prop:product rule 3}), one also deduces the inclusion
    $$
    \left[\frac{\lambda^{\frac{3}{2}}}{(t\lambda)^{2}} \right]^{-1}  V_{2k_1-1}  \cdot  \left[\frac{\lambda^{\frac{3}{2}}}{(t\lambda)^{2}} \right]^{-1} V_{2k_2} \subset \left[\frac{\lambda^{\frac{3}{2}}}{(t\lambda)^{2}} \right]^{-1} (E_{\mathrm{ori},k_1+k_2} +  E_{\mathrm{tip},k_1+k_2})
    $$
    and the same conclusion holds for the other pairs $(V_{2k_1-1}, V_{2k_2-1})$ or $(V_{2k_1}, V_{2k_2})$.
\end{remark}

\begin{proof}
\item[\ref{prop:product rule 1}.] 
    The inclusion $V_{2k} \subset E_{\mathrm{tip},k}$ is straightforward. If $v \in V_{2k-1}$ has parts $v_{\mathrm{ori}}, v_{\mathrm{mid}}, v_{\mathrm{tip}}$ where
     $$v_{\mathrm{tip}}(R,t) = \fsum_{\substack{(\alpha,i) \kadm \\ j_1,j_2 \geq 0}} \frac{\lambda^{\frac{3}{2}}}{(t\lambda)^{\alpha}}R^i \log(R)^{j_1} \log(t\lambda)^{j_2} h_{\alpha,i,j_1,j_2} \left( \frac{R}{(t\lambda)^{\frac{2}{3}+\varepsilon}} \right),$$
     the other inclusion $V_{2k-1} \subset E_{\mathrm{ori},k} + E_{\mathrm{tip},k}$ is obtained by writing
     \begin{align*}
         e_{\mathrm{ori}} = v_{\mathrm{ori}} &+ v_{\mathrm{mid}}   \\
         &+ \chi_{[1,+\infty)}(x)  \fsum_{\substack{(\alpha,i) \kadm \\ j_1,j_2 \geq 0}} \frac{\lambda^{\frac{3}{2}}}{(t\lambda)^{\alpha}}R^i \log(R)^{j_1} \log(t\lambda)^{j_2} \left[ h_{\alpha,i,j_1,j_2} (x) - h_{\alpha,i,j_1,j_2}(2) \right]
     \end{align*}
     for the variable
     $$
x = \frac{R}{(t\lambda)^{\frac{2}{3}+\varepsilon}}
     $$
     and then $e_{\mathrm{tip}} = v - e_{\mathrm{ori}}$. This proves the first assertion.

   \item[\ref{prop:product rule 2}.] To prove the inclusion $$h \left(  \frac{R}{(t\lambda)^{\frac{2}{3}+\varepsilon}} \right) (V_{2k} \cup E_{\mathrm{tip},k}) \subset E_{\mathrm{ori},k},$$
   consider an element
   $$
   v(R,t) = 
   \frac{\lambda^{\frac{3}{2}}}{(t\lambda)^{\alpha}}R^i \log(R)^{j_1} \log(t\lambda)^{j_2} q \left( \frac{R}{(t\lambda)} \right)  \chi_{[1,+\infty)}\left(  \frac{R}{(t\lambda)^{\frac{2}{3}+\varepsilon}} \right)^{1+j_3},
   $$
where $(\alpha,i)$ is $k$-admissible on $C_{\mathrm{tip}}$, $j_l \in \mathbb N_{\geq 0}$ and $\tilde{h}(y)$ is any smooth function which is constant outside $[1,2]$ and zero when $y$ is in a neighbourhood of $y = 1$. The proof is identical for an error term. The product $\tilde{h}(y) = h \cdot  \chi_{[1,+\infty)}^{1+j_3}$ is compactly supported in $(1,2)$. Hence, $R \sim (t\lambda)^{\frac{2}{3}+\varepsilon}$ and $R(t\lambda)^{-1} \sim (t\lambda)^{-\frac{1}{3}+\varepsilon} \sim 0$. One can expand
$$
q(a) = q_0(a) + \sum_{j=1}^{L}q_{j}(a) \log(a)^j
$$
for some functions $q_j$ holomorphic around $a = 0$. Consider the $M$-th order Taylor polynomial $P_j(a)$ of $q_j(a)$ near $a = 0$ with $M = \ceil{N_0(\frac{1}{3}-\varepsilon)^{-1}}$, meaning that the remainder $\eta_j(a) = q_j(a) - P_j(a)$ is a holomorphic function around zero with $(a\partial_a)^l [q_j(a) - P_j(a)] = \mathcal{O}(a^{M})$ for any fixed $l \geq 0$. Therefore,
$$
h \cdot v(R,t) = 
   \frac{\lambda^{\frac{3}{2}}}{(t\lambda)^{\alpha}}R^i \log(R)^{j_1} \log(t\lambda)^{j_2} \left(P_0(a) + \sum_{j=1}^{L}P_{j}(a) \log(a)^j \right) \tilde{h} \left(  \frac{R}{(t\lambda)^{\frac{2}{3}+\varepsilon}} \right)+ \eta 
$$
for $a = R(t\lambda)^{-1}$ and $\eta \in \mathcal{E}_{N_0,\nu}$ element (by a reasoning analogous to the proof for (\ref{eq:remainder on mid/tip is negligible})).

For the other inclusion
   $$g \left(  \frac{(t\lambda)^{\frac{2}{3}}}{R}  \right) (V_{2k-1} \cup E_{\mathrm{ori},k}) \subset E_{\mathrm{ori},k},$$
    the only issue is to show that, for the origin component $v_{\mathrm{ori}}$ of an element in $V_{2k-1}$ (resp. $E_{\mathrm{ori},k}$), then
    $$g \left( \frac{(t\lambda)^{\frac{2}{3}}}{R} \right) v_{\mathrm{ori}}$$
    can be approximated by a finite sum with the desired smallness on $C_{\mathrm{mid}}$. This follows from the proof for (\ref{eq:approximation of restricted ori term}).
   \item[\ref{prop:product rule 3}.]  For the first identity, write
    $$v_1 = v_{1,\mathrm{ori}} + v_{1,\mathrm{mid}} + v_{1,\mathrm{tip}} + \eta_1 \in V_{2k_1-1},$$
    and similarly with $e_2 \in E_{\mathrm{ori},k_2}$. Define the decomposition:
    \begin{align*}
        e_{\mathrm{ori}} &= v_{1,\mathrm{ori}} \cdot e_{2,\mathrm{ori}}, \\
        e_{\mathrm{mid}} &= v_{1,\mathrm{mid}} \cdot e_{2,\mathrm{mid}} + \tilde{v}_{1,\mathrm{ori}} \cdot e_{2,\mathrm{mid}} + v_{1,\mathrm{mid}} \cdot \tilde{e}_{2,\mathrm{ori}}, \\
        e_{\mathrm{tip}} &= v_{1,\mathrm{tip}} \cdot e_{2,\mathrm{tip}} + v_{1,\mathrm{mid}} \cdot e_{2,\mathrm{tip}} + v_{1,\mathrm{tip}} \cdot e_{2,\mathrm{mid}}, \\
        \eta &= v_1 \cdot \eta_2 + \eta_1 \cdot e_2 + (v_{1,\mathrm{ori}} - \tilde{v}_{1,\mathrm{ori}}) \cdot e_{2,\mathrm{mid}} + v_{1,\mathrm{mid}} \cdot (e_{2,\mathrm{ori}} - \tilde{e}_{2,\mathrm{ori}}).  
    \end{align*}
     Then, we claim that:
     \begin{align*}
         \left[\frac{\lambda^{\frac{3}{2}}}{(t\lambda)^{2}} \right]^{-1} ( e_{\mathrm{ori}} + e_{\mathrm{mid}} + e_{\mathrm{tip}} + \eta)  \in E_{\mathrm{ori},k_1+k_2}. 
     \end{align*}
     We only treat the mixed term $v_{1,\mathrm{ori}} \cdot e_{2,\mathrm{mid}}$. The other terms are handled similarly and do not introduce additional difficulties because for each component of $v$ and $e$ on $C_{\mathrm{ori}}, C_{\mathrm{mid}}, C_{\mathrm{tip}}$, there is a natural product structure coming from the definition. 
     
     The mixed term $v_{1,\mathrm{ori}} \cdot e_{2,\mathrm{mid}}$ is supported on $\frac{m}{2}(t\lambda)^{\frac{2}{3}} \leq R \leq m(t\lambda)^{\frac{2}{3}}$. If $g(y)$ is a smooth function coming from $e_{2,\mathrm{mid}}$, then
     $$v_{1,\mathrm{ori}} \cdot g\left( \frac{(t\lambda)^{\frac{2}{3}}}{R} \right)$$
     can be approximated by a finite sum as described in (\ref{eq:approximation of restricted ori term}). Then, it is only a matter of multiplying a finite number of $k_1$-admissible terms on $C_{\mathrm{ori}}$ (hence $k_1$-admissible on $C_{\mathrm{mid}}$) coming from $v_{1,\mathrm{ori}}$ with a finite number of $k_2$-admissible terms on $C_{\mathrm{mid}}$ coming from $e_{2,\mathrm{mid}}$.

     The second identity follows by multiplying both finite sums and using the fact that $\mathcal{Q}_{\frac{1}{2}+\frac{1}{2}\nu}$ is an algebra. The third identity follows by multiplying both finite sums and approximating the $q(a) \in \mathcal{Q}_{\frac{1}{2}+\frac{1}{2}\nu}$ functions by finite sums as in (\ref{prop:product rule 2}).
\end{proof}

\begin{corollary}[Nonlinear rules] \label{nonlinear rules for V,E}
    Let $u_0 = \lambda^{\frac{3}{2}}W(R)$,
    \begin{align*}
        w_1 &= \frac{\lambda^{\frac{3}{2}}}{(t\lambda)^2} C_1(\nu)  + V_3\\
        w_2 &= \frac{\lambda^{\frac{3}{2}}}{(t\lambda)^2} H(a^2) \cdot \chi_{[1,+\infty)} \left( \frac{R}{(t\lambda)^{\frac{2}{3}+\varepsilon}} \right) + V_4, \quad a = R/(t\lambda),
    \end{align*}
    where $H(z)$ is defined in (\ref{H(z) definition}), and $w_{k} \in V_{2k-1} \cup V_{2k}$, $k \geq 2$. Then, it holds that
    $$
    t^2[F(u_0 + w_1 + w_2 + w_k) - F(u_0 + w_1 + w_2)]  \cdot (1-\chi_{[1/m,+\infty)}) \left( \frac{(t\lambda)^{\frac{2}{3}}}{R} \right) \in E_{\mathrm{ori},k+1} + E_{\mathrm{tip},k+1},
    $$
    as well as
\begin{align*}
    t^2F(u_0 + w_1 + w_2) \cdot (1-\chi_{[1/m,+\infty)}) &\in E_{\mathrm{ori},2} + E_{\mathrm{tip},2}, \\ 
    t^2F(u_0 + w_1) \cdot (1-\chi_{[1/m,+\infty)}) &\in E_{\mathrm{ori},2} + E_{\mathrm{tip},2}.
\end{align*}
\end{corollary}

\begin{proof}
    We distinguish the two regions $\frac{m}{2}(t\lambda)^{\frac{2}{3}} \leq R \leq 2(t\lambda)^{\frac{2}{3}+\varepsilon}$ and $R \gtrsim 2(t\lambda)^{\frac{2}{3}+\varepsilon}$ using a cutoff $\chi_{[1,+\infty)}$. On the first region, the nonlinearity contributes to the middle part of an $E_{\mathrm{ori},k+1}$ element, and on the second region, we get the tip part of an $E_{\mathrm{ori},k+1}$ element plus some $E_{\mathrm{tip},k+1}$ element. More precisely, we proceed as follows:
     \item[\textbf{1.} $\frac{m}{2}(t\lambda)^{\frac{2}{3}} \leq R \leq 2(t\lambda)^{\frac{2}{3}+\varepsilon}$:] We perform a multinomial expansion around the dominant component 
        $$ \frac{\lambda^{\frac{3}{2}}}{(t\lambda)^2} g\left( \frac{(t\lambda)^{\frac{2}{3}}}{R} \right), \quad g(y) = y^3  (15)^{\frac{3}{2}} + C_1(\nu),$$
        of $u_0 + w_1 + w_2 + w_k$ and $u_0 + w_1 + w_2$. Define the cutoff:
        \begin{align*}
         \chi(R,t) &:=(1-\chi_{[1,+\infty)})\left( \frac{R}{(t\lambda)^{\frac{2}{3}+\varepsilon}} \right) \cdot (1-\chi_{[1/m,+\infty)} ) \left( \frac{(t\lambda)^{\frac{2}{3}}}{R} \right),\\
         \chi(R,t) &= \left[ \chi(R,t)^{\frac{1}{n}} \right]^{n}, \quad n \geq 0.
       \end{align*}
       As $\chi_{[a,+\infty)}(R,t)$ was chosen so that $\chi > 0$ on $\mathrm{int}(\mathrm{supp}(\chi_{[a,+\infty)}))$, the $n$-th root of $\chi$ remains a product of two smooth transition functions.  Observe that:
        \begin{align}
            \left| \left( u_0 + w_1 + w_2 - \frac{\lambda^{\frac{3}{2}}}{(t\lambda)^2} g \right) \right| &\lesssim \left( \frac{\lambda^{\frac{3}{2}}}{R^4} + \frac{\lambda^{\frac{3}{2}}}{(t\lambda)^{2 + \left( \frac{2}{3} - 2\varepsilon \right)}} + a^2 \frac{\lambda^{\frac{3}{2}}}{(t\lambda)^2} \right) \notag \\
            &\lesssim \frac{\lambda^{\frac{3}{2}}}{(t\lambda)^{2 + \left( \frac{2}{3} - 2\varepsilon \right)}}, \notag \\
            |w_k | &\lesssim \frac{\lambda^{\frac{3}{2}}}{(t\lambda)^{2 + \left( \frac{2}{3} - 2\varepsilon \right) \cdot (k-1) }}, \quad \frac{m}{2}(t\lambda)^{\frac{2}{3}} \leq R \leq 2(t\lambda)^{\frac{2}{3}+\varepsilon}, \label{eq:smallness of E_1, E_2}
        \end{align}
        and the smallness is preserved under $(t \partial_t)^{l_1}(R \partial_R)^{l_2}$ as in Proposition \ref{prop:stability under differentiation} and \ref{prop:stability of smallness under differenatiation} (note that $a = 1$ is not included in this region, hence there is no singularity of type $(1-a)^{\beta}$). Define: 
        $$
        E_1:=\left(\frac{ u_0 + w_1 + w_2 - \frac{\lambda^{\frac{3}{2}}}{(t\lambda)^2} g  }{\frac{\lambda^{\frac{3}{2}}}{(t\lambda)^2} g} \right), \quad E_2:=  \left( \frac{w_{k}}{\frac{\lambda^{\frac{3}{2}}}{(t\lambda)^2} g} \right).
        $$
        Using the product rules (\ref{prop:product rule 2}) and (\ref{prop:product rule 3}) from Proposition \ref{product rules for V, E},
        $$
       E_1^i \cdot \chi^{\frac{i}{i+j}} =  \left[ \frac{\lambda^{\frac{3}{2}}}{(t\lambda)^2}  \right]^{-i} \left( u_0 + w_1 + w_2 - \frac{\lambda^{\frac{3}{2}}}{(t\lambda)^2} g \right)^i \cdot \chi^{\frac{i}{i+j}} \in \left[ \frac{\lambda^{\frac{3}{2}}}{(t\lambda)^2}  \right]^{-1} E_{\mathrm{ori},i}, \ i \geq 1,
        $$
        as well as
        $$
        E_2^j \cdot \chi^{\frac{j}{i+j}} = \left[ \frac{\lambda^{\frac{3}{2}}}{(t\lambda)^2}  \right]^{-j} w_k^j \cdot \chi^{\frac{j}{i+j}} \in \left[ \frac{\lambda^{\frac{3}{2}}}{(t\lambda)^2}  \right]^{-1} E_{\mathrm{ori},jk}, \quad j \geq 1,
         $$
         since $w_k \in V_{2k-1} \cup V_{2k}$.
       Write
       \begin{align*}
        N(R,t) &:= t^2[F(u_0 + w_1 + w_2 + w_k) - F(u_0 + w_1 + w_2)].
       \end{align*}
       We find that
       \begin{align*}
           N \cdot \chi &\simeq  \frac{\lambda^{\frac{3}{2}}}{(t\lambda)^\frac{8}{3}} \sum_{\substack{N_0 \gtrsim i \geq 0 \\ N_0 \gtrsim j \geq 1} } \binom{p}{i,j} g^p \left(\frac{ u_0 + w_1 + w_2 - \frac{\lambda^{\frac{3}{2}}}{(t\lambda)^2} g  }{\frac{\lambda^{\frac{3}{2}}}{(t\lambda)^2} g} \right)^i \chi^{\frac{i}{i+j}} \left( \frac{w_{k}}{\frac{\lambda^{\frac{3}{2}}}{(t\lambda)^2} g} \right)^j \chi^{\frac{j}{i+j}}  \\
           &\simeq   \frac{\lambda^{\frac{3}{2}}}{(t\lambda)^\frac{8}{3}}  \sum_{\substack{N_0 \gtrsim i \geq 0 \\ N_0 \gtrsim j \geq 1} } \left[ \frac{\lambda^{\frac{3}{2}}}{(t\lambda)^2}  \right]^{-1} E_{\mathrm{ori},i}  \left[ \frac{\lambda^{\frac{3}{2}}}{(t\lambda)^2}  \right]^{-1} E_{\mathrm{ori},jk}  \in E_{\mathrm{ori},k+1}
       \end{align*}
       using the product rule (\ref{prop:product rule 3}) from Proposition \ref{product rules for V, E}. The remainder term
       $$
        \eta \cdot \chi = \frac{\lambda^{\frac{3}{2}}}{(t\lambda)^\frac{8}{3}} \sum_{\substack{i \gtrsim N_0 \\ j \gtrsim N_0} } \binom{p}{i,j} g^p E_1^i E_2^j \chi
        $$
       belongs to $\mathcal{E}_{N_0,\nu}$ with no singularity since the support is restricted away from the tip of the cone. Indeed, 
      \begin{align*}
          (t\partial_t)^{l_1}(R\partial_R)^{l_2} \left[ E_1^i \right] \lesssim_{l_1,l_2} \frac{i^{l_1+l_2}}{(t\lambda)^{\left( \frac{2}{3} - 2\varepsilon \right) \cdot i}}, \quad \frac{m}{2}(t\lambda)^{\frac{2}{3}} \leq R \leq 2(t\lambda)^{\frac{2}{3}+\varepsilon} 
      \end{align*}
       using Faa di Bruno's formula, as smallness is preserved under differentiation. A similar estimate holds for $E_2^j$. Derivatives falling on $g^p$ or $\chi$ cause no loss of smallness, as was shown in (\ref{eq: derivative t R in terms of a, y}). Hence, 
       \begin{align*}
           (t\partial_t)^{l_1}(R\partial_R)^{l_2}[\eta  \cdot \chi] &\lesssim_{l_1,l_2,p} \frac{\lambda^{\frac{3}{2}}}{(t\lambda)^\frac{8}{3}} \sum_{\substack{i \gtrsim N_0 \\ j \gtrsim N_0} } \binom{p}{i,j} i^{l_1+j_2} j^{l_1+j_2}  \left( \frac{1}{(t\lambda)^{\left( \frac{2}{3} - 2\varepsilon \right) }} \right)^{(i+j)} \\
           &\lesssim_{l_1,l_2,p} \frac{\lambda^{\frac{3}{2}}}{(t\lambda)^\frac{8}{3}} \sum_{\substack{i \gtrsim N_0 \\ j \gtrsim N_0} }\binom{p}{i,j}  \left( \frac{2}{(t\lambda)^{\left( \frac{2}{3} - 2\varepsilon \right) }} \right)^{(i+j)}  \\
           &\lesssim_{l_1,l_2,p} \frac{\lambda^{\frac{3}{2}}}{(t\lambda)^{\frac{8}{3}+N_0}} \sum_{\substack{i \gtrsim N_0 \\ j \gtrsim N_0} }\binom{p}{i,j}  \left( \frac{2}{(t\lambda)^{\left( \frac{2}{3} - 2\varepsilon -\delta \right) }} \right)^{(i+j)} \\
           &\lesssim_{l_1,l_2,p} \frac{\lambda^{\frac{3}{2}}}{(t\lambda)^{\frac{8}{3}+N_0}} \left( 1 + \frac{4}{(t\lambda)^{\left( \frac{2}{3} - 2\varepsilon -\delta \right) }} \right)^{p} 
       \end{align*}
       given $\delta > 0$ small enough so that $\frac{2}{3}-2\varepsilon-\delta > 0$ and $i,j \geq \delta^{-1}N_0$. Similarly, the elements
      \begin{align*}
       t^2 F(u_0 + w_1 + w_2) \cdot \chi &\simeq  \frac{\lambda^{\frac{3}{2}}}{(t\lambda)^\frac{8}{3}} \sum_{N_0 \gtrsim i \geq 0} \binom{p}{i} g^p \left(\frac{ u_0 + w_1 + w_2 - \frac{\lambda^{\frac{3}{2}}}{(t\lambda)^2} g  }{\frac{\lambda^{\frac{3}{2}}}{(t\lambda)^2} g} \right)^i  \cdot \chi, \\
       t^2 F(u_0 + w_1) \cdot \chi &\simeq  \frac{\lambda^{\frac{3}{2}}}{(t\lambda)^\frac{8}{3}} \sum_{N_0 \gtrsim i \geq 0} \binom{p}{i} g^p \left(\frac{ u_0 + w_1 - \frac{\lambda^{\frac{3}{2}}}{(t\lambda)^2} g  }{\frac{\lambda^{\frac{3}{2}}}{(t\lambda)^2} g} \right)^i  \cdot \chi 
      \end{align*}
       are in $E_{\mathrm{ori},2}$ due to the presence of the $\dfrac{\lambda^{\frac{3}{2}}}{(t\lambda)^\frac{8}{3}}$ factor.
       
       \item[\textbf{2.} $2(t\lambda)^{\frac{2}{3}+\varepsilon} \leq R \leq (t\lambda)$:] We perform a multinomial expansion around the dominant component
        $$ \frac{\lambda^{\frac{3}{2}}}{(t\lambda)^2} q \left( \frac{R}{(t\lambda)} \right), \quad q(a) = C_1(\nu) + H(a^2)$$
        of $u_0 + w_1 + w_2 + w_k$ and $u_0 + w_1 + w_2$. Observe that
        \begin{align*}
            \left| \left( u_0 + w_1 + w_2 - \frac{\lambda^{\frac{3}{2}}}{(t\lambda)^2} q \right) \right| &\lesssim \left( \frac{\lambda^{\frac{3}{2}}}{R^3} + \frac{\lambda^{\frac{3}{2}}}{(t\lambda)^{2 + \left( \frac{2}{3} - 2\varepsilon \right)}} + (1-\chi_{[1,+\infty)}) a^2 \frac{\lambda^{\frac{3}{2}}}{(t\lambda)^2} \right) \\
            &\lesssim \frac{\lambda^{\frac{3}{2}}}{(t\lambda)^{2 + \left( \frac{2}{3} - 2\varepsilon \right)}}, \quad 2(t\lambda)^{\frac{2}{3}+\varepsilon} \leq R \leq (t\lambda)
        \end{align*}
        and the smallness is preserved under $(t \partial_t)^{l_1}(R \partial_R)^{l_2}$ as in Proposition \ref{prop:stability under differentiation} and \ref{prop:stability of smallness under differenatiation} (up to increasing the singularity at $a = 1$ as in Proposition \ref{prop:stability of smallness under differenatiation}). We remark that
      \begin{align*}
        \left( w_2 - \dfrac{\lambda^{\frac{3}{2}}}{(t\lambda)^2} q \right) \chi_{[1,+\infty)}^{\frac{1}{i+j}} &= \frac{\lambda^{\frac{3}{2}}}{(t\lambda)^2} H(a^2) \cdot (1- \chi_{[1,+\infty)}) \left( \frac{R}{(t\lambda)^{\frac{2}{3}+\varepsilon}} \right)  \chi_{[1,+\infty)} \left( \frac{R}{(t\lambda)^{\frac{2}{3}+\varepsilon}} \right)^{\frac{1}{i+j}}  \\
        &\phantom{=}+ V_4 \in E_{\mathrm{ori},1} +  E_{\mathrm{tip},1}. 
      \end{align*}
      Similarly, it holds that
      $$
(u_0+w_1) \chi_{[1,+\infty)}^{\frac{1}{i+j}} \in V_1 \subset E_{\mathrm{ori},1} +  E_{\mathrm{tip},1}.
      $$
        Then, applying the product rules from Proposition \ref{product rules for V, E}:
        $$
            \left[ \frac{\lambda^{\frac{3}{2}}}{(t\lambda)^2}  \right]^{-i} \left( u_0 + w_1 + w_2 - \frac{\lambda^{\frac{3}{2}}}{(t\lambda)^2} q \right)^i \cdot \chi_{[1,+\infty)}^{\frac{i}{i+j}} \in \left[ \frac{\lambda^{\frac{3}{2}}}{(t\lambda)^2}  \right]^{-1} (E_{\mathrm{ori},i} + E_{\mathrm{tip},i}), \quad i \geq 1.
        $$
 We conclude as in the first part that
        $$  t^2[F(u_0 + w_1 + w_2 + w_k) - F(u_0 + w_1 + w_2)] \cdot \chi_{[1,+\infty)} \in E_{\mathrm{ori},k+1} + E_{\mathrm{tip},k+1}.$$
        Similarly, it holds that
      \begin{align*}
       t^2 F(u_0 + w_1 + w_2) \cdot \chi_{[1,+\infty)} &\simeq  \frac{\lambda^{\frac{3}{2}}}{(t\lambda)^\frac{8}{3}} \sum_{N_0 \gtrsim i \geq 0} \binom{p}{i} t^2q^p \left(\frac{ u_0 + w_1 + w_2 - \frac{\lambda^{\frac{3}{2}}}{(t\lambda)^2} q  }{\frac{\lambda^{\frac{3}{2}}}{(t\lambda)^2} q} \right)^i  \cdot \chi_{[1,+\infty)} \\
       &\in E_{\mathrm{ori},2} + E_{\mathrm{tip},2}, 
      \end{align*}
      as well as
      \begin{align*}
             t^2 F(u_0 + w_1) \cdot \chi_{[1,+\infty)}  &\simeq  \frac{\lambda^{\frac{3}{2}}}{(t\lambda)^\frac{8}{3}} \sum_{N_0 \gtrsim i \geq 0} \binom{p}{i} t^2C_1(\nu)^p \left(\frac{ u_0 + w_1 - \frac{\lambda^{\frac{3}{2}}}{(t\lambda)^2} C_1(\nu)  }{\frac{\lambda^{\frac{3}{2}}}{(t\lambda)^2} C_1(\nu)} \right)^i  \cdot \chi_{[1,+\infty)} \\
       &\in E_{\mathrm{tip},2}.
       \end{align*}
\end{proof}

In the next section, we prove the following theorem:

\begin{theorem}[Construction of an approximate solution] \label{Thm:Approximate Solution}
        Assume $d = 5$. The successive errors and correction terms satisfy the following properties when $k \geq 1$:
    \begin{enumerate}
    \item  $t^2e_{2k-1} = t^2e_{2k-1}^0 + t^2e_{2k-1}^1\in E_{\mathrm{tip},k} + E_{\mathrm{ori},k}$.
    
    \item $v_{2k} \in V_{2k}$.
    
    Moreover, the function $v_2$ is non-negative everywhere on $0 \leq R  \leq (t\lambda)$, $0 < t \leq t_0 \ll 1$ and takes the form:
    \begin{align*}
        v_2(R,t) &= \frac{\lambda^{\frac{3}{2}}}{(t\lambda)^{2}}H(a^2)\chi_{[1,+\infty)}\left(  \frac{R}{(t\lambda)^{\frac{2}{3}+\varepsilon}} \right) \\
        &\phantom{=}+ \fsum_{\substack{(\alpha,i) \mathrm{2-adm} \\ j \geq 0}}\frac{\lambda^{\frac{3}{2}}}{(t\lambda)^{\alpha}}R^i \log(R)^{j} Q_{\alpha,i,j} \left(a \right)  \chi_{[1,+\infty)}\left(  \frac{R}{(t\lambda)^{\frac{2}{3}+\varepsilon}} \right),
    \end{align*}
    where $a = R/(t\lambda)$, $H(z)$ is defined as in (\ref{H(z) definition}) and is a positive function on $(0,1)$.
        
    \item $t^2e_{2k} = t^2e_{2k}^0 + t^2e_{2k}^1 \in E_{\mathrm{ori},k} + E_{\mathrm{tip},k+1}$.
    
    \item $v_{2k+1} \in V_{2k+1}$.
    \end{enumerate}
    
    In particular, the approximate solution
    $$u_{k} = u_0 + v_1 + v_2 + \sum_{i=3}^k v_i, \quad k \geq 3$$
    is positive everywhere on $0 \leq R < (t\lambda)$, $0 < t \leq t_0$, and has asymptotics
    \begin{align*}
    | (\langle R \rangle^i \partial_R^i)(t^j\partial_t^j) u_k| &\lesssim \begin{dcases}
        \lambda^{\frac{3}{2}} \quad &0 \leq R \lesssim 1 \\
        \frac{\lambda^{\frac{3}{2}}}{1+R^3} \quad &1 \lesssim R \lesssim (t\lambda)^{\frac{2}{3}} \\
        \frac{\lambda^{\frac{3}{2}}}{(t\lambda)^2}\left[1 + \left(1-\frac{R}{(t\lambda)} \right)^{\frac{1}{2}+\frac{1}{2}\nu-i-j-} \right] \quad &(t\lambda)^{\frac{2}{3}} \lesssim R < (t\lambda)
    \end{dcases} \\
    |(R^i \partial_R^i) (t^j\partial_t^j)(u_k-u_0)| &\lesssim \begin{dcases}
        \frac{\lambda^{\frac{3}{2}}}{(t\lambda)^2} R^{i+\max\{2-i,0\}} \quad &0 \leq R \lesssim 1 \\
        \frac{\lambda^{\frac{3}{2}}}{(t\lambda)^2}\left[1 + \left(1-\frac{R}{(t\lambda)} \right)^{\frac{1}{2}+\frac{1}{2}\nu-i-j-} \right] \quad &1 \lesssim R < (t\lambda)
    \end{dcases}
    \end{align*}
    where $\lesssim$ can be replaced by $\asymp$ when $i = j = 0$.
\end{theorem}

The analogous theorem in dimension 4 is stated in Theorem \ref{Thm:Approximate Solution d = 4}.

\begin{remark}
    we observe that 
     $$
     \left(  \int_{|x| < t} (u_k-u_0)^2 dx \right)^{\frac{1}{2}} \lesssim \frac{\lambda^{\frac{3}{2}}}{(t\lambda)^2} \left( \int_{0}^t r^4 dr \right)^{\frac{1}{2}} \lesssim t^{\frac{1}{2}\nu + 1}
     $$
    with a similar estimate for the $R$ derivatives and $t$ derivatives (we lose one power of $t$ when differentiating with respect to $t$). 

    Extending $u^e = u_k-u_0$ to all of $\mathbb R^d \times [0,t_0]$ as a function of the same size and regularity, supported on $0 < |x| < 2t$ (as described in Remark \ref{extension of v_k, e_k}), one obtains the energy decay for $u^e$ claimed in Theorem \ref{thm:blow-up, main thm}.
\end{remark}

\section{Renormalization Step: Next Iterates} \label{section:renormalization step, next iterates}

We now perform the main inductive argument of the renormalization procedure in dimension $d = 5$, explaining how to construct the even correction terms $v_{2k}$ from the error $e_{2k-1}$ by solving a wave-like equation in self-similar coordinates, and the odd correction terms $v_{2k+1}$ from the error $e_{2k}$ using an elliptic-like equation. We prove that at each step, there is a systematic decrease in the error, thereby completing the proof of Theorem \ref{Thm:Approximate Solution}.

\subsection{Proof of Theorem \ref{Thm:Approximate Solution}}
The proof is done by induction. Assuming the claimed decomposition of $t^2e_{2k-1}$ (which will be proven to be true for $k = 1$), we show that $v_{2k}$, $t^2e_{2k}$, $v_{2k+1}$, $t^2e_{2k+1}$ all have the desired form.

\subsubsection{Construction of $e_1$ from $v_1$}

The error $t^2e_1(R,t)$ is given exactly by 
\begin{align*}
    t^2e_1(R,t) &= \underbrace{t^2[F(u_0 + v_1) - F(u_0) - F'(u_0)v_1]}_{=: N(e_1)} - t^2\partial_{tt} (v_1(r\lambda,t)), 
\end{align*}
where $v_1(R,t) = \lambda^{\frac{3}{2}}(t\lambda)^{-2} V_1(R)$, $V_1(R)$ is as in (\ref{v_1 formula}), $R = rt^{-1-\nu}$, and
\begin{align*}
    t^2\partial_{tt} (v_1(r\lambda,t)) &= (t^2\partial_{tt}v_1)(R,t) + \left(-\frac{3}{2}-\frac{3}{2}\nu \right)\left(-\frac{3}{2}-\frac{3}{2}\nu -1 \right) (R\partial_Rv_1)(R,t) \\
    &\phantom{=}+ 2\left(-\frac{3}{2}-\frac{3}{2}\nu \right) (R\partial_R t\partial_t v_1)(R,t) + \left(-\frac{3}{2}-\frac{3}{2}\nu \right)^2 (R^2\partial_{RR} v_1)(R,t) \\
    &\in \frac{ \lambda^{\frac{3}{2}} }{ (t \lambda)^2 } S^2(R^0,\log(R)) \\
    t^2F(u_0) &\in \frac{ \lambda^{\frac{3}{2}} }{ (t \lambda)^{-2} } S^0(R^{-7}), \quad t^2F'(u_0)v_1 \in \lambda^{\frac{3}{2}}  S^2(R^{-4},\log(R)).
\end{align*}
Moreover, we note that $t^2\partial_{tt} (v_1(r\lambda,t))$ has a constant dominant term 
$$
\frac{1}{4}(\nu-3)(\nu-5) C_1(\nu)R^0 =: C_2(\nu)R^0
$$
at $R \to +\infty$ and that $t^2\partial_{tt} (v_1(r\lambda,t))$ is the dominant component of the error near the tip of the cone. Using Corollary \ref{nonlinear rules for V,E}, we conclude that
$$
N(e_1) \cdot (1-\chi_{[1/m,+\infty)}) \left(\frac{(t\lambda)^{\frac{2}{3}}}{R} \right) \in E_{\mathrm{ori},2} + E_{\mathrm{tip},2},
$$
and one checks separately that
\begin{align*}
 t^2\partial_{tt} (v_1(r\lambda,t)) \cdot (1-\chi_{[1/m,+\infty)}) &\in E_{\mathrm{ori},1} + E_{\mathrm{tip},1} \\
 t^2F'(u_0)v_1 \cdot (1-\chi_{[1/m,+\infty)}) &\in E_{\mathrm{ori},2} + E_{\mathrm{tip},2} \\
  t^2F(u_0) \cdot (1-\chi_{[1/m,+\infty)}) &\in E_{\mathrm{ori},2} + E_{\mathrm{tip},2}
\end{align*}
by applying Proposition \ref{inclusion S^2n into V_2k-1} and the inclusion $V_{2k-1} \subset E_{\mathrm{ori},k} + E_{\mathrm{tip},k}$. It remains to analyze the term:
$$
\left( N(e_1) - t^2 \partial_{tt}v_1 \right) \cdot \chi_{[1/m,+\infty)} \left(\frac{(t\lambda)^{\frac{2}{3}}}{R} \right),
$$
which contributes to the origin part of an $E_{\mathrm{ori},1}$ element. To see this, we perform a binomial expansion around $u_0$:
\begin{align*}
   N(e_1) \cdot \chi_{[1/m,+\infty)} \simeq  \chi_{[1/m,+\infty)}  \cdot \sum_{n=2}^{\infty} \binom{p}{n} t^2 u_0^p T \left[ \left( \frac{v_1}{u_0} \right)^{n} \right],
\end{align*}
where we recall that $T$ is the “truncation” operator defined in equation (\ref{truncation}), and
\begin{align*}
     \sum_{n=2}^{\infty} \binom{p}{n} t^2 u_0^p T \left[ \left( \frac{v_1}{u_0} \right)^{n} \right] - t^2 \partial_{tt}v_1  &\in  \sum_{n=2}^{\infty}  \dfrac{\lambda^{\frac{3}{2}}}{(t\lambda)^{2(n-1)}} S^{2n}(R^{-7+3n},\log(R)^{3N_0}) \\
    &\phantom{\in}-  \dfrac{ \lambda^{\frac{3}{2}} }{ (t \lambda)^2 } S^2(R^0,\log(R)).
\end{align*}
Hence, the sum has an appropriate form for a $\tilde{e}_{\mathrm{ori}}$ component.

\subsubsection{Construction of $v_{2k}$ from $t^2e_{2k-1}$}

For each term coming from the finite sum of $t^2e_{2k-1}^0$ on $R \geq 2(t\lambda)^{\frac{2}{3}+\varepsilon}$, i.e. each term of the form:
\begin{equation} \label{t^2 tilde(e)_(2k-1)^0}
    t^2\tilde{e}_{2k-1}^0(R,a,t):=  \frac{\lambda^{\frac{3}{2}}}{(t\lambda)^{\alpha}}q_{\alpha,i,j_1,j_2}(a)R^i \log(R)^{j_1}, \quad R > 0, a \in (0,1), 0 < t \leq t_0,
\end{equation}
where $(\alpha,i)$ is $k$-admissible, $j_1, j_2 \geq 0$ and $q_{\alpha,i,j_1,j_2}(a) \in \mathcal{Q}_{\frac{1}{2}+\frac{1}{2}\nu}$, we solve
$$t^2(-\partial_t^2 + \partial_r^2 + \frac{4}{r}\partial_r) \tilde{v}_{2k} = -t^2\tilde{e}_{2k-1}^0$$
and then apply back the omitted cutoff 
$$\chi_{[1,+\infty)}\left( \frac{R}{(t\lambda)^{\frac{2}{3}+\varepsilon}}\right)^{1+j_2}$$
to the solution $\tilde{v}_{2k}$. Summing all these solutions, we obtain the correction $v_{2k}$. As shown below in Theorem \ref{hypergeometric ode with error forcing term}, the solution takes the form:
$$
v_{2k} = \fsum_{\substack{(\alpha,i) \kadm \\ j_1,j_2 \geq 0}} \frac{\lambda^{\frac{3}{2}}}{(t\lambda)^{\alpha}} R^i \left( \sum_{0 \leq l \leq j_1} Q_{\alpha,i,j_1,j_2,l}(a) \log(R)^{l} \right) \chi_{[1,+\infty)}\left( \frac{R}{(t\lambda)^{\frac{2}{3}+\varepsilon}}\right)^{j_2},$$
where $Q_{\alpha,i,j_1,j_2,l}(a) \in a^2\mathcal{Q}_{\frac{1}{2}+\frac{1}{2}\nu}$. The correction term $v_{2k}$ has comparable size to $v_{2k-1}$ near the tip of the cone $a \sim 1$. However, we will obtain a smaller error near the tip of the cone. Using an appropriate change of variables, we reduce the problem to solving an hypergeometric equation, which we first study in the following lemma.



\begin{lemma}[Hypergeometric ODE with $\tilde{\mathcal{Q}}$ forcing term] \label{Q inhomogeneous hypergeom}
   Write $z = a \in \mathbb C$. Let $\alpha, \beta, \gamma \in \mathbb R$ such that $\gamma \notin \mathbb Z_{\leq 0}$ and $\gamma - \alpha - \beta > 0$. Let $b,r \in \mathbb R$ and $q(z) \in z^r\tilde{\mathcal{Q}}_{b}$. Then the inhomogeneous hypergeometric equation
    \begin{equation}
    z(1-z)w''(z) + (\gamma - (\alpha + \beta + 1)z)w'(z) - \alpha \beta w(z) = q(z), \quad 0 < z < 1 \label{eq:hypergoemetric equation with Q tilde forcing term}
\end{equation} 
has a particular solution $w(z) \in z^{r+1}\tilde{\mathcal{Q}}_{\min\{b+1, \gamma-\alpha-\beta\}}$. 

Moreover, if $-r-1 \notin \mathbb N_{\geq 0}$ and $-\gamma - r \notin \mathbb N_{\geq 0}$ (e.g. if $r = 0$) and the worst logarithmic singularity of $q(z)$ near $z = 0$ is bounded by $\log(z)^J$, $J \in \mathbb N_{\geq 0}$, then so is the worst logarithmic singularity of the solution $w(z)$.
\end{lemma}

\begin{proof}
Around $|z| \leq a_0$, we expand:
$$q(z) = z^r \left( \sum_{j=0}^{L}q_{j}(z) \log(z)^j \right)$$
Equation (\ref{eq:hypergoemetric equation with Q tilde forcing term}) near zero becomes
$$
w''(z) + \frac{(\gamma - (\alpha + \beta + 1)z)}{z(1-z)}w'(z) - \frac{\alpha \beta}{z(1-z)} w(z) = z^{-1}(1-z)^{-1}q(z).
$$
Hence, we must solve a finite number of hypergeometric equations of the form:
$$
w''(z) + \frac{(\gamma - (\alpha + \beta + 1)z)}{z(1-z)}w'(z) - \frac{\alpha \beta}{z(1-z)} w(z) = z^{(r+1)-2} \tilde{q}_j(z) \log(z)^j,
$$
where $\tilde{q}_j(z) = (1-z)^{-1}q_j(z)$ is holomorphic around $|z| \leq a_0$. The indicial roots at zero are $\{0,1-\gamma\}$. Using Theorem \ref{thm:inhomogeneous fuchs ode}, there exists a particular solution of the form 
$$w(z) = z^{r+1}\left( \sum_{j=0}^{L} \sum_{k=0}^{j+2} Q_{j,k}(z) \log(z)^k \right)$$
If $0 - (r+1) = -r-1 \notin \mathbb N_{\geq 0}$ and $(1-\gamma) - (r+1) = -\gamma-r \notin \mathbb N_{\geq 0}$, Theorem \ref{thm:inhomogeneous fuchs ode} also ensures that the sum over the indices $k$ only goes from $k = 0$ to $j$, i.e., $Q_{j,k} = 0$ for $k \in \{j+1,j+2\}$ and the logarithmic singularity cannot increase.

Since $q(z)$ is holomorphic on $(0,1) \subset U \subset B(0,1)$, we can use regular ODE theory to get a holomorphic extension of $w(z)$ on $U$ solving the equation. On $|z-1| < a_0$, $z^r$ is analytic, so we can write
$$q(z) = q_0(1-z) + \sum_{i=1}^{+\infty}(1-z)^{\beta(i)} \sum_{j=0}^{iL}q_{i,j}(1-z)\log(1-z)^j,$$
where either the sum is finite or $\beta(i) \geq c(i-1) + b$ for $c > 0$ small enough and the growth condition
$$||q_{i,j}||_{L^{\infty}(|z-1| < z_0)} \leq C^i$$
holds for some $C > 0$. Using Theorem \ref{thm:inhomogeneous fuchs ode} once again (the indicial roots at $z = 1$ are $\{0,\gamma-\alpha-\beta\}$), a particular solution is given by
\begin{align*}
    Q_{0}(1-z) &+ \sum_{i=1}^{+\infty}(1-z)^{\beta(i)+1} \sum_{j=0}^{L_i} \sum_{k=0}^{j+1}Q_{i,j,k}(1-z)\log(1-z)^k \\
    = Q_{0}(1-z) &+ \sum_{i=1}^{+\infty}(1-z)^{\beta(i)+1} \sum_{j=0}^{L_i+1}\tilde{Q}_{i,j}(1-z) \log(1-z)^j
\end{align*}
meaning that $w(z)$ must match this particular solution modulo some linear combination of the fundamental system (see (\ref{fuchs fundamental system}) from Section \ref{section:appendix, ode}) of the ODE. The fundamental system introduces a $(1-z)^{\gamma-\alpha-\beta}\log(1-z)$ singularity in the solution. We note that the hypergeometric equation near $z = 1$ with analytic forcing term $q_0(1-z)$ can always yield a logarithm-free analytic solution $Q_0(1-z)$ thanks to Remark \ref{remark:analytic sol inhomogeneous fuchs equation}.

If the expansion for $q(z)$ is finite, then so is the resulting sum for the solution. If it is infinite, then one must verify that the boundedness condition is still verified by the $\tilde{Q}_{i,j}$ when $i$ is sufficiently large. This holds due to the estimates from Theorem \ref{thm:inhomogeneous fuchs ode}, as because the growth of the $\beta(i)$ exponents is at least linear in $i$, while the logarithmic exponents growth is at most linear. In other words, when $i \in \mathbb N$ is large enough, one has 
\begin{enumerate}
    \item $\beta(i) - \max\{|r_1|,|r_2|\} > (c/2)i > 1$, where $\{r_1, r_2\}$ are the indicial roots of the equation at $z = 1$.
    \item and an exponential upper bound $$\left( \frac{j}{\beta(i)-r_k} \right)^j \lesssim C(c,L)^i, \quad k \in \{1,2\},$$
because $0 \leq j \leq iL$, $\beta(i) - \max\{|r_1|,|r_2|\} \geq (c/2)i$.
\end{enumerate}
\end{proof}

\begin{theorem}[Particular solution to (\ref{even v})]\label{hypergeometric ode with error forcing term}
Let $d \geq 1$. Let $e(R,a,t) = t^{s}q(a)R^i\log(R)^k$ where $s,i \in \mathbb R$, $s - \nu i > -(d-1)/2$, $k \in \mathbb N$, $q \in a^{\delta}\mathcal{Q}_{\beta}$, $\beta, \delta \in \mathbb R$. Then one can find a solution to
$$t^2\left(-\partial_t^2 + \partial_r^2 + \frac{d-1}{r}\partial_r\right) v = e(r\lambda, r/t, t)$$
of the form $$v(R,a,t) = t^sR^i\sum_{0 \leq l \leq k} Q_l(a)\log(R)^l$$ where $Q_l \in a^{\delta+2}\mathcal{Q}_{\min\{\beta+1,s-\nu i + \frac{d-1}{2}\}}$, $R = r \lambda$, $a = r/t$. Moreover, if $\delta - i = 0$ and $q(a)$ has no logarithmic singularity at  $a= 0$, then so does $a^{-2}Q_l(a)$ for any $l$.
\end{theorem}

\begin{proof}
Writing $R^i = a^i(t\lambda)^{i}$, we can assume without loss of generality that $i = 0$ and $s > -(d-1)/2$. Plugging $v$ in the equation and matching powers of $\log(R)$, we find the following system of recursive equations for $Q_k$, $Q_{k-1}$, ..., $Q_0$:
\begin{align}
    t^2(-\partial_t^2 + \partial_r^2 + \frac{d-1}{r}\partial_r) \left[ t^{s} Q_k \right]  &= t^{s} q(a),  \notag \\
    t^2(-\partial_t^2 + \partial_r^2 + \frac{d-1}{r}\partial_r)\left[ t^{s} Q_{k-1} \right]  &= (d-2)kt^{s}Q_k(a)a^{-2} +(2s-1)(\nu+1)kt^{s}Q_k(a), \notag  \\
    t^2(-\partial_t^2 + \partial_r^2 + \frac{d-1}{r}\partial_r)\left[ t^{s} Q_{l-2} \right] &= (d-2)(l-1)t^{s}Q_{l-1}(a)a^{-2}  \notag \\
    &\phantom{=}+(2s-1)(\nu+1)(l-1)t^{s}Q_{l-1}(a)  \notag \\
    &\phantom{=}- l(l-1)t^{s}Q_{l}(a)\left[(\nu + 1)^2 - a^{-2}\right] \notag \\
    &\phantom{=}-2(l-1) t^{s} Q_{l-1}'(a)\left[ (\nu + 1) a - a^{-1} \right], \label{eq: system of equation near tip of the cone}
\end{align}
where $a = r/t$, $0 < r < t < t_0$. Hence, we must solve equations of the form:
$$t^2\left(-\partial_t^2 + \partial_r^2 + \frac{d-1}{r}\partial_r\right) \left[ t^{s} w \right]  = t^{s} f \in t^{s} a^{\delta} \mathcal{Q}_{\beta},$$
which is equivalent to
$$t^2 \left( - \left( \partial_t + \frac{s}{t} \right)^2 + \partial_r^2 + \frac{d-1}{r}\partial_r \right)w(a)  = f(a) \in a^{\delta} \mathcal{Q}_{\beta},$$
or $L_{s} w(a) = f(a)$, $0 < a < 1$, where
$$L_{s} = (1-a^2)\partial_{aa} + ((d-1)a^{-1} + 2as - 2a)\partial_a + (s - s^2).$$
Finally, writing $f(a) = a^{\delta}F(a^2)$ and looking for a solution of the form $w(a) = W(a^2)$, we reduce to an hypergeometric equation for $W(z)$:
$$z(1-z)W''(z) + \left( \frac{d}{2} + z \left(s - \frac{3}{2} \right) \right) W'(z) + \frac{s - s^2}{4} = z^{\frac{\delta}{2}}F(z), \quad 0 < z < 1,$$
\index{z@$z$, a complex variable which replaces $R$ or $a$ depending on the context}
whose parameters are
$$
\tilde{\alpha} = -\frac{s}{2}, \tilde{\beta} = -\frac{s}{2} + \frac{1}{2}, \tilde{\gamma} = \frac{d}{2}.
$$
Since $d \geq 1$ and $s > -(d-1)/2$, we can use Lemma \ref{Q inhomogeneous hypergeom} to obtain a solution $$
W(z) \in a^{\frac{\delta}{2}+1}\tilde{\mathcal{Q}}_{\min\{\beta+1,s+\frac{d-1}{2}\}}.
$$
This yields a solution $w(a) \in a^{\delta+2}\mathcal{Q}_{\min\{\beta+1,s+\frac{d-1}{2}\}}$. If $\delta = 0$ and $f(z)$ has no logarithmic singularity at $z = 0$, then the Lemma also implies that the $\mathcal{Q}$ part of $w(a)$ has no logarithmic singularity at $z = 0$ either.

 Since $q(a) \in a^{\delta}\mathcal{Q}_{\beta}$, we find $Q_k \in a^{\delta+2}\mathcal{Q}_{\min\{\beta+1,s + \frac{d-1}{2}\}}$. Solving for $Q_{k-1}$ using $Q_k$ leads to $$Q_{k-1} \in a^{\delta+2}\mathcal{Q}_{\min\{ \min\{\beta+1,s + \frac{d-1}{2}\} + 1, s+ \frac{d-1}{2}\}} = a^{\delta+2}\mathcal{Q}_{\min\{\beta+2,s + \frac{d-1}{2}\}} $$
Furthermore, observe that
$$
Q_{l-1}(a), a^{-2}Q_{l-1}(a), aQ'_{l-1}(a), a^{-1}Q'_{l-1}(a) = a^{-2} a \partial_a Q_{l-1}(a) \in a^{\delta} \mathcal{Q}_{\min\{\beta+2,s + \frac{d-1}{2}\}-1},
$$
and solving for the other $Q_{l-2}$'s leads to 
$$
Q_{l-2} \in a^{\delta+2}\mathcal{Q}_{\min\{ \min\{\beta+2,s+ \frac{d-1}{2}\}, s + \frac{d-1}{2}\}} = a^{\delta+2}\mathcal{Q}_{\min\{\beta+2,s+ \frac{d-1}{2}\}}.
$$
Finally, if $\delta = 0$ and $q(a)$ has no logarithmic singularity at $a = 0$, then so does $a^{-2}Q_k(a)$ and this property propagates to $Q_{k-1}$ and all $Q_{l-2}$ by induction.
\end{proof}

\begin{corollary}  \label{hypergeometric ode with error forcing term, dim 5}
In dimension $d = 5$, for a forcing term of the form: $$\frac{\lambda^{\frac{3}{2}}}{(t\lambda)^{\alpha}}q(a)R^i \log(R)^j,$$ 
with $\alpha \geq i + 2$ and $q(a) \in \mathcal{Q}_{\frac{1}{2}+\frac{1}{2}\nu}$, one can apply the theorem and obtain coefficients $Q_l \in a^2\mathcal{Q}_{\frac{1}{2} + \frac{1}{2}\nu}$ in the system (\ref{eq: system of equation near tip of the cone}).
\end{corollary}

\begin{proof}
    This is a direct application of Theorem \ref{eq:hypergoemetric equation with Q tilde forcing term} with 
\begin{align*}
        \beta &= \frac{1}{2}\nu + \frac{1}{2} \\
        s - \nu i &= \frac{3}{2}(-1-\nu) + \alpha \nu - \nu i = \nu \left( \alpha - i - \frac{3}{2} \right) - \frac{3}{2} \geq \frac{1}{2}\nu - \frac{3}{2} \\
        s-\nu i + \frac{(d-1)}{2} &\geq \frac{1}{2}\nu + \frac{1}{2}.
\end{align*}
\end{proof}

\begin{remark}[Loss of regularity in higher dimensions]\label{rmk:loss of regularity high dimension}
One of the main difficulties in generalizing this method of constructing blow-ups in higher dimensions is that solving the ODE with a forcing term $$\frac{\lambda^{\frac{d-2}{2}}}{(t\lambda)^2}$$
introduces a $(1-a)^{\frac{1}{2}+ \nu \left( \frac{6-d}{2} \right)}$ singularity, meaning that for $d > 6$, the obtained correction term is not even continuous at $a = 1$ unless $\nu > 0$ is small. 
\end{remark}

When $k = 1$, the dominant component of $-t^2\tilde{e}_1^0$ on the interval $2(t\lambda)^{\frac{2}{3}+\varepsilon} \leq R \leq (t\lambda)$ is of the form:
$$
\frac{\lambda^{\frac{3}{2}}}{(t\lambda)^{2}}C_2(\nu),
$$ 
and arises from $-t^2\partial_{tt}(v_1(r\lambda,t)) \in  \frac{\lambda^{\frac{3}{2}}}{(t\lambda)^{2}}S^{2}(R^0,\log(R))$. The remainder of the error is negligible compared to this term.

One gets a correction term $v_2$ whose dominant component is of the form:
\begin{align}
    \frac{\lambda^{\frac{3}{2}}}{(t\lambda)^{2}}H(a^2) \cdot \chi_{[1,+\infty)}\left( \frac{R}{(t\lambda)^{\frac{2}{3}+\varepsilon}}\right), \label{v_2 dominant component}
\end{align}
 where $H(z)$, $H(0) = 0$, solves
$$
z(1-z)H''(z) + \left( \frac{5}{2} + z \left(s - \frac{3}{2} \right) \right) H'(z) + \frac{s - s^2}{4} = C_2(\nu), \quad 0 < z < 1,
$$
with $s = \frac{3}{2}(-1-\nu) + 2\nu$ and $\nu > 3$. Explicitly, $H(z)$ is given by
\begin{equation*} 
    H(z) = \frac{C_2(\nu)}{\tilde{\alpha}\tilde{\beta}}\left( F(\tilde{\alpha},\tilde{\beta};\tilde{\gamma},z)-1 \right) = 4C_1(\nu)\left( F(\tilde{\alpha},\tilde{\beta};\tilde{\gamma},z)-1 \right), \quad 0 \leq z < 1,
\end{equation*}
where
$$\tilde{\alpha} = -\frac{s}{2}, \tilde{\beta} = -\frac{s}{2} + \frac{1}{2}, \tilde{\gamma} = \frac{5}{2}.$$
This is exactly the function from $\tilde{\mathcal{Q}}_{\frac{1}{2} + \frac{1}{2}\nu}$ defined earlier in (\ref{H(z) definition}) and for which $(H(z^2) + C_1(\nu))^e \in \tilde{\mathcal{Q}}_{\frac{1}{2} + \frac{1}{2}\nu}$ for any exponent $e \in \mathbb R$. 

Finally, observe that $u_0 + v_1 + v_2$, which is equal to $u_0+v_1$ on $0 \leq R \leq (t\lambda)^{\frac{2}{3}+\varepsilon}$ because of the cutoff, is also positive on $(t\lambda)^{\frac{2}{3}+\varepsilon} \leq R \leq (t\lambda)$, $0 < t \leq t_0 \ll 1$, by positivity of $H(z)$.

\subsubsection{Computation of $t^2e_{2k}$ from $v_{2k}$}

The error $t^2e_{2k}$ is given by
$$t^2e_{2k} \simeq E^t(v_{2k}) + t^2e_{2k-1}^1 + t^2[F(v_{2k} + u_{2k-1}) - F(u_{2k-1})],$$
where $E^t(v_{2k})$ denotes the components of $t^2 \square v_{2k}$ where at least one derivative falls in a cutoff $\chi_{[1,+\infty)}$. By construction, $t^2e_{2k-1}^1 \in E_{\mathrm{ori},k}$. Moreover, we prove that $E^t(v_{2k}) \in E_{\mathrm{ori},k}$ as well and 
$$
N(R,t) := t^2[F(v_{2k} + u_{2k-1}) - F(u_{2k-1})] \in E_{\mathrm{ori},k+1} + E_{\mathrm{tip},k+1}.
$$
Assume for simplicity and by linearity that we only have one term:
$$v_{2k} = \frac{\lambda^{\frac{3}{2}}}{(t\lambda)^{\alpha}} R^i \left( \sum_{0 \leq l \leq j} Q_{\alpha,i,j,l}(a) \log(R)^{l} \right) \chi_{[1,+\infty)}\left( \frac{R}{(t\lambda)^{\frac{2}{3}+\varepsilon}}\right) = \tilde{v}_{2k} \cdot \chi_{[1,+\infty)} \in V_{2k} \subset E_{\mathrm{tip},k}.$$
To prove that $E^t(v_{2k}) \in E_{\mathrm{ori},k}$, observe that
\begin{align} \label{Et(v_2k)}
    E^t(v_{2k}) &= \nu\left(\frac{2}{3}+\varepsilon\right) \chi_{[1,+\infty)}'\left( \frac{R}{(t\lambda)^{\frac{2}{3}+\varepsilon}}\right) \frac{R}{(t\lambda)^{\frac{2}{3}+\varepsilon}} t\partial_t(\tilde{v}_{2k}(r\lambda,r/t,t)) \\
    &\phantom{=}+ \nu\left(\frac{2}{3}+\varepsilon\right)\left[\nu \left(\frac{2}{3}+\varepsilon\right)-1\right] \chi_{[1,+\infty)}''\left( \frac{R}{(t\lambda)^{\frac{2}{3}+\varepsilon}}\right) \left(\frac{R}{(t\lambda)^{\frac{2}{3}+\varepsilon}} \right)^2 \tilde{v}_{2k}(r\lambda,r/t,t) \notag \\
    &\phantom{=}- a^{-2} \frac{R} {(t\lambda)^{\frac{2}{3}+\varepsilon}} \chi_{[1,+\infty)}'\left( \frac{R}{(t\lambda)^{\frac{2}{3}+\varepsilon}}\right) \tilde{v}_{2k}(r\lambda,r/t,t) \notag \\
    &\phantom{=}- a^{-2} \frac{R} {(t\lambda)^{\frac{2}{3}+\varepsilon}} \chi_{[1,+\infty)}'\left( \frac{R}{(t\lambda)^{\frac{2}{3}+\varepsilon}}\right) r\partial_r (\tilde{v}_{2k}(r\lambda,r/t,t)) \notag \\
    &\phantom{=}- a^{-2} \left( \frac{R} {(t\lambda)^{\frac{2}{3}+\varepsilon}} \right)^2 \chi_{[1,+\infty)}''\left( \frac{R}{(t\lambda)^{\frac{2}{3}+\varepsilon}}\right) \tilde{v}_{2k}(r\lambda,r/t,t). \notag
\end{align}
is an element of $E_{\mathrm{ori},k}$ supported on the region $(t\lambda)^{\frac{2}{3}+\varepsilon} \leq R \leq 2(t\lambda)^{\frac{2}{3}+\varepsilon}$. Approximating  each $Q_{\alpha,i,j}(a) \in a^2\mathcal{Q}_{\frac{1}{2}+\frac{1}{2}\nu}$ by a finite sum as in the proof of (\ref{prop:product rule 2}) from Proposition \ref{product rules for V, E}, we obtain a finite sum of $\tilde{e}_{\mathrm{mid}}$ components, together with some approximation error belonging to $\mathcal{E}_{N_0,\nu}$. For this term $E^t(v_{2k})$, there is no apparent gain of smallness compared to $t^2e_{2k-1}$ or $v_{2k}$, but the support is now near the origin.

Finally, we deal with the nonlinear part $N(R,t)$ of $t^2e_{2k}$. This part is supported on $(t\lambda)^{\frac{2}{3}+\varepsilon} \leq R \leq 2(t\lambda)^{\frac{2}{3}+\varepsilon}$. Therefore, we can introduce a harmless cutoff $(1-\chi_{[1/m,+\infty)})$, i.e.,
$$
N(R,t) = N(R,t) \cdot (1-\chi_{[1/m,+\infty)}) \left( \frac{(t\lambda)^{\frac{2}{3}}}{R}\right),
$$
and apply Corollary \ref{nonlinear rules for V,E} with $$w_1 = v_1 + \sum_{i=2}^{k}v_{2i-1}, \quad w_2 = v_2 + \sum_{i=2}^{k-1}v_{2i}, \quad w_k = v_{2k}$$ if $k > 1$ to conclude that $N(R,t) \in E_{\mathrm{ori},k+1} + E_{\mathrm{tip},k+1}$. When $k = 1$, we apply the second part of Corollary \ref{nonlinear rules for V,E} with $w_1 = v_1, w_2 = v_2$ to obtain separately $t^2 F(u_0 + v_1 + v_2)$, $t^2F(u_0+v_1) \in E_{\mathrm{ori},2} + E_{\mathrm{tip},2}$ by choosing $w_1 = v_1, w_2 = v_2$.

\subsubsection{Construction of $v_{2k+1}$ from $t^2e_{2k}$}\label{subsubsec: construction v_2k+1}

We solve (\ref{odd v}) again. As in Secion \ref{section:renormalization step, first iterate}, we are led to solve
\begin{equation}
    (t\lambda)^2\mathcal{L}v_{2k+1}(R,t) = t^2e_{2k}^0(R,t), \quad R \geq 0, \quad \mathcal{L} = -\partial_R^2 - \frac{4}{R} \partial_R - pW(R)^{p-1}, \label{eq:v_2k+1 L operator def}
\end{equation}
where $t$ is treated as a parameter and $t^2e_{2k}^0(R,t) \in E_{\mathrm{ori},k}$ is supported on $0 \leq R \leq 2(t\lambda)^{\frac{2}{3}+\varepsilon}$. We solve the equation separately for each part $t^2e_{2k,\mathrm{ori}}^0$, $t^2e_{2k,\mathrm{mid}}^0$ and $t^2e_{2k,\mathrm{tip}}^0$ and, as stated earlier, the negligible part $\eta$ can be ignored and absorbed into $t^2e_{2k}^1$ and subsequent error terms. On the domain $0 \leq R \leq m(t\lambda)^{\frac{2}{3}}$, we can solve the equation for each element 
$$\frac{\lambda^{\frac{3}{2}}}{(t\lambda)^{\alpha+2n}} w_{n}^{\alpha,I,j}(R), n \in \mathbb N,$$
where $w_{n}^{\alpha,I,j} \in S^{2(k-1)}(R^{I+3n}, \log(R)^{J})$ for some common $J \in \mathbb N_{\geq 0}$, as given in the decomposition of $t^2e_{2k}^0(R,t)$. After applying back the cutoffs $\chi_{[1/m,+\infty)}^{1+j}$ to the solution, the resulting solution is
$$v_{2k+1,\mathrm{ori}}(R,t) = \fsum_{ \substack{(\alpha,I) \kadm \\ j \geq 0} }  \sum_{n = 0}^{+\infty} \ \frac{\lambda^{\frac{3}{2}}}{(t\lambda)^{\alpha+2+2n}} W_{n}^{\alpha,I,j}(R)\chi_{[1/m,+\infty)}\left(\frac{(t\lambda)^{\frac{2}{3}}}{R} \right)^{1+j},$$
where $W_n^{\alpha,I,j} \in S^{2k}(R^{I+2+3n}, \log(R)^{J+2})$. In particular, if $(\alpha,I)$ was $k$-admissible on $C_{\mathrm{ori}}$, then the new pair $(\alpha + 2, i + 2)$ is $(k+1)$-admissible on $C_{\mathrm{ori}}$.

\begin{theorem}  \label{usual step 1}
Let $w(R,t) \in \dfrac{\lambda^{\frac{3}{2}}}{(t\lambda)^{\alpha}}S^{2n}(R^{I}, \log(R)^J)$ for some $\alpha \in \mathbb R$, $I,J,n \in \mathbb N_{\geq 0}$. Let $\mathcal{L}$ be as in (\ref{eq:v_2k+1 L operator def}). Then the correction term $v(R,t)$ defined by the equation
$$(t\lambda)^2 \mathcal{L}v(R,t) = w(R,t), \quad v(0,t) = v'(0,t) = 0,$$
 belongs to $\dfrac{\lambda^{\frac{3}{2}}}{(t\lambda)^{\alpha+2}}S^{2n+2}(R^{I+2} \log(R)^{J+2})$, and satisfies the estimates
\begin{align*}
    ||v||_{S, \mathrm{ori}} &\leq C(\mathcal{L},R_0)||w||_{S, \mathrm{ori}}, \\
    ||v||_{S,\infty} &\leq C(\mathcal{L},R_0,J,m)||w||_{S,\infty},
\end{align*}
for some constant $C$ which is independent of $n$, $I$ and $\alpha$.
\end{theorem}

\begin{proof}
Let $w(z,t) = \lambda^{\frac{3}{2}} (t\lambda)^{-\alpha} w(z)$ and $v(z,t) =  \lambda^{\frac{3}{2}} (t\lambda)^{-\alpha-2} v(z)$. On a neighbourhood of $|z| \leq \sqrt{15}/2$, $w(z)$ is holomorphic with a zero of order $2n$. Using a logarithmic-free fundamental system $\{u_1, u_2\}$ found as in (\ref{fuchs fundamental system}) ($r_1 = 0, r_2 = -3$, $u_1(z) = o(1)$, $u_2(z) = o(z^{-3})$, $W(u_1,u_2)(z)^{-1} = o(z^4)$) on this neighbourhood, extended smoothly to $\mathbb R$ using regular ODE theory, one obtains the solution
\begin{align*}
    v(z) &= \int_{[0,z]} [u_2(z)u_1(y)-u_1(z)u_2(y)]W(u_1,u_2)(y)^{-1} w(y)dy, \\
     v'(z) &= \int_{[0,z]} [u_2'(z)u_1(y)-u_1'(z)u_2(y)]W(u_1,u_2)(y)^{-1} w(y)dy,
\end{align*}
which has the desired regularity and a zero of order $2n+2$ at $z = 0$. Moreover, 
\begin{align*}
    ||v(z)||_{L^{\infty}(F)}  + ||v'(z)||_{L^{\infty}(F)} \leq C(\mathcal{L},R_0) \cdot ||w(z)||_{L^{\infty}(F)} 
\end{align*}
where $F = \{z \in \mathbb C: |z| \leq \sqrt{15}/2 \} \cup \{z \in \mathbb R:  \sqrt{15}/2 \leq z \leq R_0 \}$, because the singularity of the integrand at the origin (which comes from $u_1,u_2$) is removable.

On a neighbourhood of $\Re(z) > 0, |z| \geq R_0$, there is an expansion
$$w(z) = z^I \sum_{j=0}^J w_{j} \log(z)^j$$ 
and we can use Theorem \ref{thm:inhomogeneous fuchs ode} ($r_1 = 3, r_2 = 0, \beta = -I-2$, see also (\ref{step 1 at infinity}) for the ODE at infinity) to find a particular solution of the form 
$$\tilde{v}(z) = z^{I+2}\sum_{j=0}^{J}\sum_{k=0}^{j+2} v_{j,k}(z^{-1})  \log(z)^k = z^{I+2} \sum_{k=0}^{J+2} \underbrace{\left( \sum_{j = \max\{0,k-2\}}^{J} v_{j,k}(z^{-1}) \right)}_{=: v_k} \log(z)^k$$
with the coefficient estimate
$$||v_{j,k}||_{A(|y| \leq R_0^{-1})} \leq C(\mathcal{L},R_0,J) ||w_{j}||_{L^{\infty}(|y| \leq R_0^{-1})}.$$
This implies that
\begin{align*}
    ||\tilde{v}||_{S,\infty} &= m^{I+2} \max_{0 \leq k \leq J+2} ||v_{k}||_{L^{\infty}(|y| \leq R_0^{-1})} \\ 
    &\leq m^{-2} J C(\mathcal{L},R_0,J) m^I \max_{0 \leq j \leq J} ||w_{j}||_{L^{\infty}(|y| \leq R_0^{-1})} \\
&\leq m^{-2} J C(\mathcal{L},J) ||w||_{S,\infty}. 
\end{align*}
Hence, $v(z)$ must match this particular solution modulo some linear combination 
$$c_1 U_1(z^{-1}) + c_2U_2(z^{-1}) = z^{I} \left( c_1 z^{-I}U_1(z^{-1}) + c_2z^{-I}U_2(z^{-1}) \right) $$
of the fundamental system $\{U_1(y), U_2(y) = \tilde{U}_2(y) + c \cdot U_1(y) \log(y) \}$ at infinity. It remains to check that 
$$||c_1 y^{I} U_1(y)||_{A^{\infty}(|y| \leq R_0^{-1})} +||c_2y^{I} \tilde{U}_2(y)||_{A^{\infty}(|y| \leq R_0^{-1})}  \leq C(\mathcal{L},R_0,J,m)||w||_{S,\infty}$$
to conclude the proof. It is sufficient to prove that
$$
|c_1| + |c_2| \leq R_0^I \cdot C(\mathcal{L},R_0,J,m) ||w||_{S,\infty}.
$$
To this end, we first observe that
$$|v(R_0)| + |v'(R_0)| \leq C(\mathcal{L},R_0) ||w||_{S, \mathrm{ori}}$$
and
$$|\tilde{v}(R_0)| + |\tilde{v}'(R_0)| \leq R_0^{I} \cdot C(\mathcal{L},R_0,J) ||w||_{S,\infty}$$
using our estimates on $v$, $v'$,$\tilde{v}$, $\tilde{v}'$. Since $c_1$, $c_2$ solves the linear system:
\begin{align*}
    c_1U_1(R_0^{-1}) + c_2U_2(R_0^{-1}) + \tilde{v}(R_0) &= v(R_0), \\
    c_1U_1'(R_0^{-1}) + c_2U_2'(R_0^{-1}) + \tilde{v}'(R_0) &= v'(R_0), 
\end{align*}
one finds that
$$
|c_1| + |c_2| \leq R_0^{I} \cdot C(\mathcal{L},R_0,J,m).
$$
This finishes the proof.
\end{proof}

 Applying back the cutoffs $$\chi_{[1/m,+\infty)}\left( \frac{(t\lambda)^{\frac{2}{3}}}{R} \right)^{1+j}$$ creates an additional error $E^t(v_{2k+1,\mathrm{ori}}) \in E_{\mathrm{ori},k+1}$, supported on $\frac{m}{2} (t\lambda)^{\frac{2}{3}} \leq R \leq m(t\lambda)^{\frac{2}{3}}$, made of those terms in $\mathcal{L}v_{2k+1,\mathrm{ori}}$ where at least one derivative hits the cutoff. Indeed, assuming for simplicity and by linearity that we have only one sum
$$v_{2k+1,\mathrm{ori}}(R,t) = \chi_{[1/m,+\infty)} \left(\frac{(t\lambda)^{\frac{2}{3}}}{R} \right) \underbrace{\sum_{n = 0}^{+\infty} \ \frac{\lambda^{\frac{3}{2}}}{(t\lambda)^{\alpha+2+2n}} W_{n}^{\alpha,I}(R)}_{=: \tilde{v}_{2k+1,\mathrm{ori}}},$$
this leads to an error term:
\begin{align*}
    E^t(v_{2k+1,\mathrm{ori}}) &=  +\frac{1}{R^2} \left( \frac{(t\lambda)^{\frac{2}{3}}}{R} \right)^2 \chi''_{[1/m,2/m]}\left( \frac{(t\lambda)^{\frac{2}{3}}}{R} \right) \tilde{v}_{2k+1,\mathrm{ori}} \\
    &-  \frac{1}{R^2} \left( \frac{(t\lambda)^{\frac{2}{3}}}{R} \right) \chi'_{[1/m,2/m]}\left( \frac{(t\lambda)^{\frac{2}{3}}}{R} \right) R\partial_R \tilde{v}_{2k+1,\mathrm{ori}} \\
    &- \frac{2}{R^2}  \left( \frac{(t\lambda)^{\frac{2}{3}}}{R} \right) \chi'_{[1/m,2/m]}\left( \frac{(t\lambda)^{\frac{2}{3}}}{R} \right) \tilde{v}_{2k+1,\mathrm{ori}}.
\end{align*}
It follows from Proposition \ref{prop:stability under differentiation} and point (\ref{prop:product rule 2}) from Proposition \ref{product rules for V, E}, as well as the gain of smallness of order $R^{-2} \sim (t\lambda)^{-\frac{4}{3}}$, that $E^t(v_{2k+1,\mathrm{ori}}) \in E_{\mathrm{ori},k+1}$. Hence, this error is absorbed into the next error $t^2e_{2k+1}$.

Next, we solve (\ref{odd v}) with each term 
 $$w(R,t) = \frac{\lambda^{\frac{3}{2}}}{(t\lambda)^{\alpha}}R^i \log(R)^{j_1} g_{\alpha,i,j_1,j_2} \left( \frac{(t\lambda)^\frac{2}{3}}{R} \right) (1-\chi_{[1,+\infty)}) \left( \frac{R}{(t\lambda)^{\frac{2}{3}+\varepsilon}} \right)^{1+j_2}$$
 coming from $t^2e_{2k,\mathrm{mid}}^0$. As before, we first solve the ODE without the cutoff, then apply the cutoff afterward. The error $E^t(v_{2k+1,\mathrm{mid}})$ caused by this simplification, which is supported on $(t\lambda)^{\frac{2}{3}+\varepsilon}\leq R \leq 2(t\lambda)^{\frac{2}{3}+\varepsilon}$, can be included in $t^2e_{2k+1}$ similarly to $E^t(v_{2k+1,\mathrm{ori}})$. In order to show this, we start with a lemma concerning primitives for our $g_{\alpha,i,j_1,j_2}(y)$ functions.

\begin{lemma} \label{lemma:primitive of z^ig(z)log(z)^j}
    Assume that $g(z)$ is smooth on $(0,+\infty)$, zero on $(2/m,+\infty)$ and expands as a finite sum of holomorphic functions and logarithms on $|z| < z_0$,
    $$g(z) = \sum_{j=0}^J g_j(z) \log(z)^j, \quad |z| < z_0, z \notin \mathbb R_{\leq 0}.$$
    
    If $I \in \mathbb Z$, there exists a primitive of $z^I g(z)$ of the form $z^{\min\{I+1,0\}}G(z)$, where $G(z)$ is smooth on $(0,+\infty)$, zero on $(2/m,+\infty)$ and expands as a finite sum of holomorphic functions and logarithms on $|z| < z_0$. Explicitly, 
    $$z^{\min\{I+1,0\}}G(z) = \int_{z}^{2/m} y^I g(y) dy$$
    when $z \in (0,+\infty)$.
\end{lemma}

\begin{proof}
Without loss of generality, we assume that $g(z) = g_{J}(z)\log(z)^J$ near $z = 0$,
$$g_J(z) = \sum_{n=0}^{\infty}g_nz^n, \quad |z| < z_0$$
instead of having a sum with different logarithmic exponents. The primitive
    $$\tilde{G}(z) = \int_{z}^{2/m} y^I g(y) dy, \quad z \in (0,+\infty)$$
     of $z^I g(z)$ is smooth on $(0,+\infty)$, zero on $(2/m,+\infty)$ and extends holomorphically on $|z| < z_0$, $z \notin \mathbb R_{\leq 0}$, via:
     $$G(z) = \int_{[z_0/2,z]} y^I g_J(y) \log(y)^J dy + \underbrace{\int_{z_0/2}^{2/m} y^I g(y) dy}_{\text{cst}}, \quad |z| < z_0, z \notin \mathbb R_{\leq 0}.$$
     We compute
     \begin{align*}
\int_{[z_0/2,z]} y^I g_J(y)\log(y)^J dy &= \int_{[z_0/2,z]} \sum_{n=0}^{\infty} g_n y^{n+I} \log(y)^J dy \\
    &= \sum_{n=0}^{\infty} \int_{[z_0/2,z]} g_n y^{n+I} \log(y)^J dy, \quad |z| < z_0, z \notin \mathbb R_{\leq 0}.
\end{align*}
If $\delta \in \mathbb R \setminus \{-1\}$, $J \in \mathbb N_{\geq 0}$, then a primitive of $y^{\delta} \log(y)^J$ is given by 
$$z^{\delta+1} \sum_{k=0}^J \frac{(-1)^k J!}{(J-k)!} \cdot \frac{\log(z)^{J-k}}{(\delta+1)^{k+1}},$$
and if $\delta = -1$, a primitive is given by $(J+1)^{-1}\log(y)^{J+1}$. For $k \in \{0,...,J\}$, let
\begin{align*}
    A_{J,k}(z) =  \frac{(-1)^k J!}{(J-k)!} \sum_{\substack{n=0 \\ n + I \neq -1}}^{+\infty} \frac{g_{n}}{(n+I+1)^{k+1} }z^n,
\end{align*}
which is holomorphic on $|z| < |z_0|$. Then, for $|z| < z_0, z \notin \mathbb R_{\leq 0}$, we obtain
$$\tilde{G}(z) = z^{I+1} \left( \sum_{k=0}^J A_{J,k}(z) \log(z)^{J-k} + z^{-I-1}g_{-1-I}(J+1)^{-1}\log(y)^{J+1}\right) + C,$$
     with the convention that $g_{-1-I} = 0$ if $I \notin \mathbb Z_{\leq -1}$. If $I + 1 \in \mathbb N_{\geq 0}$, we set $G(z) = \tilde{G}(z)$. Otherwise, we set $G(z) = z^{-I-1}\tilde{G}(z)$.
\end{proof}

\begin{theorem}  \label{second step 1}
Let $w(R,t) = \dfrac{\lambda^{\frac{3}{2}}}{(t\lambda)^{\alpha}}R^I \log(R)^J g\left(\dfrac{(t\lambda)^{\frac{2}{3}}}{R}\right)$ for some $k$-admissible pair $(\alpha, I)$ on $C_{\mathrm{mid}}$, $J \in \mathbb N_{\geq 0}$ and a smooth $g(y)$ on $(0,+\infty)$ which is zero on $(2/m,+\infty)$ and extends as a finite sum of holomorphic functions and logarithms near $0$. Let $\mathcal{L}$ be as in (\ref{eq:v_2k+1 L operator def}). Then the correction term $v(R,t)$ obtained by solving
$$(t\lambda)^2 \mathcal{L}v(R,t) = w(R,t), \quad v(0,t) = v'(0,t) = 0$$
is of the form:
$$v(R,t) = \fsum_{\substack{(\tilde{\alpha},i) \ \mathrm{(k+1)}\text{-adm} \\ 0 \leq j \leq J+2}} \frac{\lambda^{\frac{3}{2}}}{(t\lambda)^{\tilde{\alpha}}}R^{i} \log(R)^j G_{\tilde{\alpha},i,j} \left( \frac{(t\lambda)^{\frac{2}{3}} }{R}\right) + \eta,$$
where $(\tilde{\alpha},i)$ is $(k+1)$-admissible on $C_{\mathrm{mid}}$, $\eta \in \mathcal{E}_{N_0,\nu}$, $G_{\tilde{\alpha},i,j}(y)$ is smooth on $(0,+\infty)$, zero on $(2/m,+\infty)$ and extends as a finite sum of holomorphic functions and logarithms near $y = 0$.
\end{theorem}

\begin{remark}
    The specific initial condition for the ODE is not critical, as we will multiply the solution by the cutoffs we ignored, which will make it zero at $R = 0$.
\end{remark}

\begin{proof} 
 The solution is explicitly given by
\begin{align*}
    v(R,t) = \frac{1}{(t\lambda)^2} R^{-2}\int_{\frac{m}{2}(t\lambda)^{\frac{2}{3}}}^{R}\left[ \theta(R)\phi(s) - \theta(s)\phi(R) \right]s^2w(s,t)ds,
\end{align*}
where $\{\phi, \theta\}$ is the fundamental system from Section 4 (see (\ref{phi}) and (\ref{theta})). When $\Re(z) > 0$, $|z| \geq R_0$, we have:
\begin{align*}
    \phi(z) &= \underbrace{z^{-1}\sum_{n=0}^{3N_0+1}a_iz^{-i}}_{\tilde{\phi}(z)} + z^{-1}\sum_{n=3N_0+2}^{+\infty}a_iz^{-i} \\
        \theta(z) &= \underbrace{c \log(z)\tilde{\phi}(z) + z^{2}\sum_{n=0}^{3N_0+3}b_iz^{-i}}_{\tilde{\theta}(z)} + z^{-1}\sum_{n=3N_0+4}^{+\infty}b_iz^{-i} + c \log(z)(\phi(z)-\tilde{\phi}(z)).
\end{align*}
Plugging into the formula for $v$, we obtain
\begin{align*}
    v(R,t) &= \frac{1}{(t\lambda)^2} R^{-2}\int_{\frac{m}{2}(t\lambda)^{\frac{2}{3}}}^{R}\left[ \tilde{\theta}(R)\tilde{\phi}(s) - \tilde{\theta}(s)\tilde{\phi}(R) \right]s^2w(s,t)ds + \eta,
\end{align*}
where $\eta \in \mathcal{E}_{N_0,\nu}$. Expanding the integrand, we are led to analyze a finite number of terms of the form:
$$\frac{1}{(t\lambda)^{2+\alpha}} \log(R)^{\delta} R^{l} \int_{\frac{m}{2}(t\lambda)^{\frac{2}{3}}}^{R}s^{i} \log(s)^{j} g \left( \frac{(t\lambda)^{\frac{2}{3}}}{s} \right) ds,$$
where $\delta \in \{0,1\}$, $(\alpha,I)$ is $k$-admissible on $C_{\mathrm{mid}}$ by hypothesis, $0 \leq j \leq J + 1$ and either ($i \leq  I + 4$, $l \leq -3$) or ($i \leq  I + 1$, $l \leq 0$) depending on whether we deal with $\tilde{\theta}(R)\tilde{\phi}(s)$ or $\tilde{\theta}(s)\tilde{\phi}(R)$. We do a change of variables $s = (t\lambda)^{\frac{2}{3}}y^{-1}$ in the integral and obtain
$$
(t\lambda)^{-2-\alpha+\frac{2}{3}(i+1)} \log(R)^{\delta} R^{l} \int_{\frac{(t\lambda)^{\frac{2}{3}}}{R}}^{m/2}y^{-i-2} \left( - \log(y) + \log((t\lambda)^{\frac{2}{3}}) \right)^{j} g\left( y \right) dy.
$$
We are thus left to analyze a finite number of integrals of the form:
$$
(I) = (t\lambda)^{-2-\alpha+\frac{2}{3}(i+1)} \log((t\lambda)^{\frac{2}{3}})^{j_1}  \log(R)^{\delta} R^{l} \int_{\frac{(t\lambda)^{\frac{2}{3}}}{R}}^{m/2} y^{-i-2}  \log(y)^{j_2} g \left( y \right) dy
$$
with $j = j_1 + j_2$.
Let
$$y^{\min\{-i-1,0\}}G_{i,j_2}(y) = \int_{y}^{m/2} s^{-i-2}  \log(s)^{j_2} g \left( s \right) ds$$
be the primitive from Lemma \ref{lemma:primitive of z^ig(z)log(z)^j}, where $G_{i,j_2}(y)$ a smooth function on $(0,+\infty)$ being zero when $y > 2/m$ and which expands as a finite sum of holomorphic functions and logarithms near $y = 0$. Plugging the primitive into (I), we finally get
$$
(I) = (t\lambda)^{-2-\alpha+\frac{2}{3}(i+1)} \left[ 
 \log \left( \frac{(t\lambda)^{\frac{2}{3}}}{R} \right) + \log(R) \right]^{j_1} \log(R)^{\delta} R^{l} \left( \frac{(t\lambda)^{\frac{2}{3}}}{R} \right)^{\min\{-i-1,0\}} G_{i,j_2} \left( \frac{(t\lambda)^{\frac{2}{3}}}{R} \right).
 $$
Finally, we check the smallness. If $-i-1 \geq 0$, then the smallness is determined by
$$(t\lambda)^{-2-\alpha+\frac{2}{3}(i+1)} R^l$$
 where $l \leq 0$ and $i+l+1 \leq I+2$. In the middle region $m(t\lambda)^{\frac{2}{3}} \leq R \leq (t\lambda)^{\frac{2}{3}+\varepsilon}$, this is bounded by
\begin{align*}
    (t\lambda)^{-2-\alpha+\frac{2}{3}(i+1)+\frac{2}{3}l} &= \frac{(t\lambda)^{\frac{2}{3}(i+1+l)}}{(t\lambda)^{\alpha+2}} =  \frac{(t\lambda)^{\frac{2}{3}(I+2)}}{(t\lambda)^{\alpha+2}} \cdot (t\lambda)^{\frac{2}{3}(i+l+1-I-2)} \\
    &\lesssim \frac{(t\lambda)^{\frac{2}{3}(I+2)}}{(t\lambda)^{\alpha+2}} \lesssim \frac{|R|^I}{(t\lambda)^{\alpha}} \cdot \frac{1}{(t\lambda)^{\frac{2}{3}}}
\end{align*}
 meaning that $(\alpha + 2 - \frac{2}{3}(i+1), l)$ is $(k+1)$-admissible on $C_{\mathrm{mid}}$ if $(\alpha,I)$ is $k$-admissible. Otherwise, $-i-1<0$ and the smallness is determined by
$$(t\lambda)^{-2-\alpha} R^{l+i+1} = \frac{R^{I+2}}{(t\lambda)^{\alpha+2}} \cdot R^{l+ i-I-1}$$
 where $l+i-I-1 \leq 0$. In that case, there is also a gain of smoothness in the middle region and  $(\alpha + 2, l+i+1)$ is $(k+1)$-admissible.
\end{proof}

Finally, we handle the contribution from  $t^2e_{2k,\mathrm{tip}}^0$, with a result structurally similar to the one for $t^2e_{2k,\mathrm{mid}}^0$. However, there is a subtle difference regarding the admissibility of the resulting exponent pairs, which are admissible on $C_{\mathrm{tip}}$ while the initial pairs were only admissible on $C_{\frac{2}{3}+\varepsilon}$. Specifically, the new $R^i$-factors for the solution all have negative exponents $i \leq 0$. For such a term, admissibility on $C_{\frac{2}{3}+\varepsilon}$ and $C_{\mathrm{tip}}$ are equivalent.

 \begin{theorem} \label{third step 1}
Let $w(R,t) = \dfrac{\lambda^{\frac{3}{2}}}{(t\lambda)^{\alpha}}R^I \log(R)^J h\left(\dfrac{R}{(t\lambda)^{\frac{2}{3}+\varepsilon}}\right)$ for some $k$-admissible pair $(\alpha,I)$ on $C_{\frac{2}{3}+\varepsilon}$ and a smooth $h(y)$ which has compact support in $(1,2)$. Let $\mathcal{L}$ be as in (\ref{eq:v_2k+1 L operator def}). Then the correction term $v(R,t)$ obtained by solving
$$(t\lambda)^2 \mathcal{L}v(R,t) = w(R,t), \quad v(0,t) = v'(0,t) = 0$$
is of the form:
$$v(R,t) = \fsum_{\substack{(\tilde{\alpha},i) \ \mathrm{(k+1)}\text{-adm} \\ 0 \leq j_1 \leq 1 \\
0 \leq j_2 \leq J}} \frac{\lambda^{\frac{3}{2}}}{(t\lambda)^{\tilde{\alpha}}}R^{i} \log(R)^{j_1} \log(t\lambda)^{j_2} H_{\tilde{\alpha},i,j} \left( \frac{R}{(t\lambda)^{\frac{2}{3}+\varepsilon}} \right) + \eta,$$
where $i \leq 0$, $(\tilde{\alpha},i)$ is $(k+1)$-admissible on $C_{\mathrm{tip}}$, $\eta \in \mathcal{E}_{N_0,\nu}$ has no singularity at the boundary $R = (t\lambda)$ and $H_{\tilde{\alpha},i,j}(y)$ is smooth, constant outside $[1,2]$ and vanishes in a neighbourhood of $y = 1$.
 \end{theorem}

\begin{proof}
    As in the proof Theorem \ref{second step 1}, the solution $v(R,t)$ is composed of a negligible term and a finite sum of elements of the form:
    $$\frac{1}{(t\lambda)^{2+\alpha}} \log(R)^{\delta} R^{l} \int_{(t\lambda)^{\frac{2}{3}+\varepsilon}}^{R}s^{i} \log(s)^{j} h \left( \frac{s}{(t\lambda)^{\frac{2}{3}+\varepsilon}} \right) ds,$$
where $\delta \in \{0,1\}$, $(\alpha,I)$ is $k$-admissible on $C_{\mathrm{mid}}$ by hypothesis, $0 \leq j \leq J + 1$ and either ($i \leq  I + 4$, $l \leq -3$) or ($i \leq  I + 1$, $l \leq 0$) depending on whether we deal with $\tilde{\theta}(R)\tilde{\phi}(s)$ or $\tilde{\theta}(s)\tilde{\phi}(R)$.  We do a change of variables $s = (t\lambda)^{\frac{2}{3}+\varepsilon}y$ in the integral, which leads to a finite number of integrals of the form:
$$
(I) = (t\lambda)^{-2-\alpha+\left( \frac{2}{3} + \varepsilon \right) (i+1)} \log((t\lambda)^{\frac{2}{3}+\varepsilon})^{j_1}  \log(R)^{\delta} R^{l} \int_{1}^{\frac{R}{(t\lambda)^{\frac{2}{3}+\varepsilon}}} y^{i}  \log(y)^{j_2} h \left( y \right) dy
$$
with $j = j_1 + j_2$.
Let
$$H_{i,j_2}(y) = \int_{1}^{y} s^{i}  \log(s)^{j_2} h \left( s \right) ds.$$
Since $y^j \log(y)^{j_2} h(y)$ is smooth with compact support in $(1,2)$, the primitive $H_{i,j_2}(y)$ is a smooth function on $(0,+\infty)$, vanishing when $y < 1$ and constant when $y > 2$. We finally get the desired form:
$$
(I) =(t\lambda)^{-2-\alpha+\left( \frac{2}{3} + \varepsilon \right) (i+1)}  \log((t\lambda)^{\frac{2}{3}+\varepsilon})^{j_1}   \log(R)^{\delta} R^{l} H_{i,j_2} \left( \frac{R}{(t\lambda)^{\frac{2}{3}+\varepsilon}} \right).
 $$
Finally, we check the smallness, which is determined by
$$(t\lambda)^{-2-\alpha+\left( \frac{2}{3} + \varepsilon \right) (i+1)}R^l,$$
 where $l \leq 0$ and $i+l+1 \leq I+2$. On $R \sim (t\lambda)^{\frac{2}{3}+\varepsilon}$ or on the tip of the cone $(t\lambda)^{\frac{2}{3}+\varepsilon} \leq R \leq (t\lambda)$, this is bounded by
\begin{align*}
    (t\lambda)^{-2-\alpha+\left( \frac{2}{3} + \varepsilon \right) (i+1+l)} &= \frac{(t\lambda)^{\left( \frac{2}{3} + \varepsilon \right) (i+1+l)}}{(t\lambda)^{\alpha+2}} =  \frac{(t\lambda)^{\left( \frac{2}{3} + \varepsilon \right) (I+2)}}{(t\lambda)^{\alpha+2}} \cdot (t\lambda)^{\left( \frac{2}{3} + \varepsilon \right) (i+l+1-I-2)} \\
    &\lesssim \frac{(t\lambda)^{\left( \frac{2}{3} + \varepsilon \right) (I+2)}}{(t\lambda)^{\alpha+2}} \lesssim \frac{|R|^I}{(t\lambda)^{\alpha}} \cdot \frac{1}{(t\lambda)^{\frac{2}{3}-2\varepsilon}},
\end{align*}
 meaning that $\left[ \alpha + 2 - \left(\frac{2}{3}+\varepsilon \right) (i+1), l\right]$ is $(k+1)$-admissible on $C_{\frac{2}{3}+\varepsilon}$ if $(\alpha,I)$ is $k$-admissible. Since $l \leq 0$, admissibility on $C_{\frac{2}{3}+\varepsilon}$ and on $C_{\mathrm{tip}}$ are equivalent.
\end{proof}

All of this leads to the desired correction term $v_{2k+1} \in V_{2k+1}$.
\begingroup
\allowdisplaybreaks
\subsubsection{Computation of $t^2e_{2k+1}$ from $v_{2k+1}$}

The error term $t^2e_{2k+1} = t^2e_{2k+1}^0 + t^2e_{2k+1}^1\in E_{\mathrm{tip},k+1} + E_{\mathrm{ori},k+1}$ is given by
\begin{align*}
     t^2e_{2k+1}^0 &\simeq
     t^2[F(v_{2k+1} + u_{2k}) - F(u_{2k}) - F'(u_0)v_{2k+1}] \\
     &+ t^2\partial_{tt} v_{2k+1} + t^2e_{2k}^1 + E^t(v_{2k+1}),
\end{align*}
    where $E^t(v_{2k+1})$ consists of those components in $\mathcal{L} v_{2k+1,\mathrm{ori}}$ and $\mathcal{L} v_{2k+1,\mathrm{mid}}$ where at least one derivative hits the cutoff. It has already been established in Section \ref{subsubsec: construction v_2k+1} that $E^t(v_{2k+1}) \in E_{\mathrm{ori},k+1}$. A straightforward computation also shows that $t^2\partial_{tt} v_{2k+1} \in V_{2k+1}$ and we know, by Proposition \ref{product rules for V, E}, that
    $$V_{2k+1} \subset E_{\mathrm{ori},k+1} + E_{\mathrm{tip},k+1}.$$ 
    Hence, it remains to verify that
\begin{align*}
    N(e_{2k+1}) &= t^2[F(v_{2k+1} + u_{2k}) - F(u_{2k}) -F'(u_0)v_{2k+1}  ] \\
&\in E_{\mathrm{ori},k+1} + E_{\mathrm{tip},k+1}.
\end{align*}
Using Corollary \ref{nonlinear rules for V,E} with 
 $$w_1 = v_1 + \sum_{i=2}^{k}v_{2i-1}, \quad w_2 = v_2 + \sum_{i=2}^{k}v_{2i}, \quad w_k = v_{2k+1},$$ 
 we only need to handle
 $$
 N(e_{2k+1})_{\mathrm{ori}}:=  \chi_{[1/m,+\infty)} \left(\frac{(t\lambda)^{\frac{2}{3}}}{R} \right)\cdot N(e_{2k+1}) \in E_{\mathrm{ori},k+1}.
 $$
On the region $0 \leq R \leq m(t\lambda)^{\frac{2}{3}}$, write
$$N(e_{2k+1}) = t^2[F(v_{2k+1} + u_{2k}) - F(u_{2k}) - F'(u_{2k})v_{2k+1} + (F'(u_{2k})-F'(u_0))v_{2k+1}  ], $$
so that
\begin{align}
N(e_{2k+1})_{\mathrm{ori}} &\simeq  \chi_{[1/m,+\infty)} \cdot \sum_{ \substack{ i \geq 0 \\ 3N_0+1 \geq j \geq 0 \\
3N_0+1 \geq l \geq 2}} \binom{p}{i,j,l} t^2 u_0^p T \left[ \left( \frac{v_1}{u_0} \right)^{i} \right] \left( \frac{u_{2k}-u_0-v_1}{u_0} \right)^{j} \left( \frac{v_{2k+1}}{u_0} \right)^{l} \notag \\
&+ \chi_{[1/m,+\infty)}  \cdot \frac{v_{2k+1}}{u_0} \sum_{ \substack{i,j \geq 0 \\ i+j \geq 1 \\ 3N_0+1 \geq j \geq 0 }} \binom{p-1}{i,j} t^2 u_0^{p}  T \left[ \left( \frac{v_1}{u_0} \right)^{i} \right]  \left( \frac{u_{2k}-u_0-v_1}{u_0} \right)^{j}, \label{N(e_2k+1)^1_ori}
\end{align}
where $T$ is the truncation operator from (\ref{truncation}). Observe the following additional facts:
\begin{enumerate}
\item For any $k \geq 1$, 
\begin{align*}
    t^2 u_0^p V_{2k-1} \subset \lambda^{\frac{3}{2}} (t\lambda)^2 V_{2k-1}, \quad
    t^2 u_0^p \left( \frac{v_1}{u_0} \right) V_{2k-1} \subset \lambda^{\frac{3}{2}} V_{2k-1}.
\end{align*}
    \item \begin{align*}
    \chi_{[1/m,+\infty)} \cdot \sum_{i \geq 2} c_i  t^2 u_0^p T \left[ \left( \frac{v_1}{u_0} \right)^{i} \right] &\in V_1
\end{align*}
for any $(c_i)_{i\geq 0} \in \ell^{\infty}$.
\item For any $k \geq 1$, 
$$\chi_{[1/m,+\infty)}  u_0^{-1} V_{2k-1}  \cdot \subset \chi_{[1/m,+\infty)} \cdot \left[ \frac{ \lambda^{\frac{3}{2}} } {(t\lambda)^2} \right]^{-1} V_{2k-1} \subset V_{2k-1}.$$
\end{enumerate}

Examining first sum in (\ref{N(e_2k+1)^1_ori}), we must distinguish three cases : $i = 0, i = 1$ and $i \geq 2$ cases. Using the product rules from Proposition \ref{product rules for V, E}, we find that we have a finite sum of the form:
\begin{align*}
     \sum_{ \substack{ 3N_0+1 \geq j \geq 0 \\
3N_0+1 \geq l \geq 2}}  \left[ t^2u_0^p + t^2 u_0^p \left( \frac{v_1}{u_0} \right) + V_{1} \right] \cdot V_{2 \cdot 2 - 1}^j \cdot V_{2(k+1)-1}^l \cdot  \left[ \frac{ \lambda^{\frac{3}{2}} } {(t\lambda)^2} \right]^{-j-l}\\
\subset \left( 1 + (t\lambda)^2 + (t\lambda)^4 \right) E_{\mathrm{ori},k+3} \subset  E_{\mathrm{ori},k+1}.
\end{align*}
The other sum is treated analogously.
\endgroup

\section{The spectral theory of the linearized operator} \label{section:spectral theory of the linearized operator}
At this point, we pivot from constructing the approximate solution to developing the analytical tools required to find the final correction term needed to obtain an exact solution of (NLW) on the light cone. When solving equation (\ref{space-variable equation epsilon}) for $\varepsilon$ in Section \ref{section:exact solution fourier method}, the following Sturm-Liouville operator arises
\begin{align*}
    \mathcal{L} &= -\partial_{RR} - pW(R)^{p-1} + \frac{1}{R^2} \cdot \left( \frac{(d-3)(d-1)}{4} \right) \\
    &= -\partial_{RR} - \frac{d^2(d+2)(d-2)}{(R^2 + d(d-2))^2} + \frac{1}{R^2} \cdot \left( \frac{(d-3)(d-1)}{4} \right) \\
    &= -\partial_{RR} + V(R)
\end{align*}
\index{L@$\mathcal{L}, \mathcal{L}_R$, the perturbed Schrödinger operator}
on $(0,+\infty)$, which is self-adjoint on 
$$\dom(\mathcal{L}) = \{ f \in L^2((0,+\infty)): f,f' \in AC_{loc}((0,+\infty)), \ \mathcal{L}f \in L^2((0,+\infty))\}.$$
This section is dedicated to studying this perturbed Schrödinger operator. We will characterize its spectrum, construct its generalized eigenfunctions $\phi(R,\xi)$ as well as a spectral measure $d\rho(\xi)$ which will be key tools to construct a generalized Fourier transform. When expressing the equation for $\varepsilon$ in the generalized Fourier space, $\mathcal{L}$ is transformed into a multiplication by $\xi$, which will make the equation easier to solve using a fixed point argument.

We will study this operator for general $d \in \mathbb N_{\geq 4}$ and will later restrict to $d \in \{4,5\}$. We aim to obtain precise asymptotic estimates on the spectral measure, the eigenfunctions and the Jost solution. A fundamental system $\{\phi(R), \theta(R)\}$ of $\mathcal{L}f = 0$ with $W(\theta,\phi) = 1$ is given by
\begin{equation}
    \phi(R) = R^{\frac{d-1}{2}} \frac{d}{d \lambda}\Bigr|_{\lambda = 1} \left( \lambda^{\frac{d-2}{2}}W(\lambda R) \right) =  \frac{R^{\frac{d-1}{2}}(R^2-d(d-2))}{(R^2+d(d-2))^{\frac{d}{2}}}
    \label{phi}
\end{equation}
\begin{equation}
    \theta(R) = \frac{R^{\frac{d-1}{2}}(R^2-d(d-2))}{(R^2+d(d-2))^{\frac{d}{2}}} \cdot \int^R \frac{1}{\phi(s)^2}ds, 
    \label{theta}
\end{equation}
\index{phi@$\phi(R)$, an explicit element of the fundamental system for $\mathcal{L}u = 0$}
\index{theta@$\theta(R)$, an explicit element of the fundamental system for $\mathcal{L}u = 0$}
where 
$$
\int^R \frac{1}{\phi(s)^2}ds = \begin{dcases}
a \log(R) + \frac{b}{R^2-d(d-2)} + \sum_{\subalign{i&=-(d-2) \\ i &= 0 \mod 2}}^{d-2} c_i R^i &\text{ if } d = 0 \mod 2 \\
\phantom{a \log(R) +} \frac{bR}{R^2-d(d-2)} + \sum_{\subalign{i&=-(d-2) \\ i &= 1 \mod 2}}^{d-2} c_i R^i &\text{ if } d = 1 \mod 2,
\end{dcases}
$$
and the following asymptotics hold for $\phi$ and $\theta$:
\begin{equation}
    \phi(R) \asymp \begin{cases}
    R^{\frac{d-1}{2}}, \quad R \to 0 \\
    R^{\frac{3-d}{2}}, \quad R \to +\infty \\
    \end{cases} \quad 
    \theta(R) \asymp \begin{cases}
    R^{\frac{3-d}{2}}, \quad R \to 0 \\
    R^{\frac{d-1}{2}}, \quad R \to +\infty \\
    \end{cases} 
    \label{phi theta asympotics}
\end{equation}
with symbol-type behaviour of the derivatives.

From the behaviour of $\mathcal{L}$ at the endpoints and the number of zeros of $\phi(R)$ on $(0,+\infty)$, one deduces the following spectral properties for $\mathcal{L}$:

\begin{proposition}[Properties of the Sturm-Liouville operator]
For $d \in \mathbb N_{\geq 4}$, $\mathcal{L}$ is limit point at zero and limit point at infinity. Moreover, the spectrum of $\mathcal{L}$ decomposes as follows:
\begin{enumerate}
    \item Essential spectrum: $\spec_{ess}(\mathcal L)= [0,+\infty)$.
    \item Absolutely continuous spectrum: $\spec_{ac}(\mathcal L) = [0,+\infty)$.
    \item Singularly continuous spectrum: $\spec_{sc}(\mathcal L) = \emptyset$.
    \item Pure point spectrum: $\spec_{pp}(\mathcal L) = \{\xi_4\}$ if $d = 4$ and $\spec_{pp}(\mathcal L) = \{\xi_d,0\}$ otherwise for some $\xi_d < 0$.
\end{enumerate}
\end{proposition}

\begin{proof}
All of these properties follow only from the behaviour of $\phi(R)$ and the potential $V(R)$. One can look at \cite{Fulton1} for the limit-point, limit-circle dichotomy and at \cite[Lemma 9.35]{Teschl}, \cite[Theorem 15.3]{Weidmann}, \cite[Chapter XIII.7, Theorem 40, Theorem 55 and Corollary 56]
{Dunford} for the remainder of the statement.
\end{proof}

More generally, one can find a fundamental system $\{\phi(R,z), \theta(R,z)\}$ for $\mathcal{L}u = zu$ of the following form:
\index{phi z@$\phi(R,z)$, part the fundamental system for $\mathcal{L}u = zu$ obtained in Proposition \ref{fundamental system Lu = zu}}
\index{theta z@$\theta(R,z)$, part the fundamental system for $\mathcal{L}u = zu$ obtained in Proposition \ref{fundamental system Lu = zu}}

\begin{proposition}[Expansion for $\phi(R,z)$]\label{fundamental system Lu = zu}
For any $z \in \mathbb C$, there exists a fundamental system $\{\phi(R,z), \theta(R,z)\}$, $W(\theta,\phi) = 1$, for $\mathcal{L}u = zu$, real-valued whenever $z \in \mathbb R$ and such that $\phi(R,z)$ is given by the following absolutely convergent series:
\begin{equation}
    \phi(R,z) = \phi(R) + R^{\frac{d-1}{2}} \sum_{j=1}^{\infty}(R^2z)^j \phi_j(R^2),
    \label{phi(R,z)}
\end{equation}
where $\phi_j$ is holomorphic on $U = \{z \in \mathbb C: \Re(z) > -\frac{1}{2}d(d-2)\}$,
\begin{align*}
|\phi_j(u)| &\leq \frac{C^j}{(j-1)!} (1 + |u|)^{-1}  \cdot \log(1 + |u|)^{\delta_{d = 4}}, \quad \forall u \in U, \\
|\phi_1(u)| &\gtrsim (1 + |u|)^{-1}  \cdot \log(1 + |u|)^{\delta_{d = 4}}, \quad \forall u \in U, |u| \gtrsim 1,
\end{align*}
for some constant $C > 0$ independent of $j$ and $u$ (the logarithm appears only in dimension $d = 4$). Moreover, $\theta(R,z)$ is entire with respect to $z$, $R^{\frac{d-3}{2}}\theta(R,z) \in C^0([0,+\infty) \times \mathbb C)$ and it is a Frobenius type solution in the following sense:
$$\lim \limits_{R \to 0} R^{-(l+1)} \frac{d^{n_l + 1}}{dz^{n_l + 1}}\theta(R,z) = 0, \quad l = \frac{d-3}{2}, n_l = \floor{l + 1/2}.$$
\end{proposition}

\begin{proof}
The exact form of $\theta(R,z)$ does not matter. We only need its existent, which is proved in \cite{Kostenko2}. As for $\phi(R,z)$, we try to look for a solution of the form: $$\phi(R,z) = R^{-\frac{d-1}{2}}\sum_{j=0}^{\infty}z^jf_j(R), \quad \phi(0,z) = \phi'(0,z) = 0,$$
where $f_0(R) = R^{\frac{d-1}{2}}\phi(R)$. Then, $f_j$ must solve
$$\mathcal{L}(R^{-\frac{d-1}{2}}f_j) = R^{-\frac{d-1}{2}}f_{j-1}, \quad f_j(0) = f_j'(0) = 0.$$
By induction, we find such solutions $f_j$ and prove that $f_j$ has a zero of order $R^{(d-1)+2j}$ at $R = 0$, $R^{-(d-1)}f_j(R) = g_j(R^2)$ for some $g_j$ analytic on $U$ and for $z \in U$, $j \geq 1$:
$$|f_j(z)| \leq \frac{C^j}{(j-1)!} (1+|z|)^{(d-1)+2(j-1)}(\cdot \log(1 + |z|) \text{ if } d = 4).$$
Using the variation of parameters, one gets
\begin{align*}
    f_j(R) &= -\int_0^R \frac{ R^{\frac{d-1}{2}} }{ s^{\frac{d-1}{2}} } \left[\phi(R)\theta(s) - \phi(s)\theta(R) \right] f_{j-1}(s)ds.
\end{align*}
Let $c(R) = \int^s \frac{1}{\phi(s)^2} ds$ as in (\ref{theta}) and
$$F_{j-1}(R) = \int_{0}^R s^{-\frac{d-1}{2}}\phi(s)f_{j-1}(s)ds = \int_{0}^R \frac{s^2-d(d-2)}{(s^2+d(d-2))^{\frac{d}{2}}}f_{j-1}(s)ds.$$
Then $F_{j-1}(R)$ has a zero of order $R^{d + 2(j-1)}$ at $R = 0$ and $R^{-d}F_{j-1}(R) = G_{j-1}(R^2)$ where $G_{j-1}$ is analytic on $U$. Hence, for $R < d(d-2)$, we can use integration by parts and write 
\begin{align}
    f_j(R) &= -R^{\frac{d-1}{2}} \phi(R) \int_0^R s^{-\frac{d-1}{2}}\phi(s) \left[c(s)-c(R) \right] f_{j-1}(s)ds \notag \\
    &= R^{\frac{d-1}{2}} \phi(R) \int_0^R F_{j-1}(s)c'(s)ds \notag \\
    &= \frac{R^{d-1}(R^2-d(d-2))}{(R^2+d(d-2))^{\frac{d}{2}}} \int_0^R G_{j-1}(s^2)\frac{s(s^2 + d(d-2))^d}{(s^2-d(d-2))^2}ds.
    \label{f_j at zero}
\end{align}
From this last formula (\ref{f_j at zero}), we deduce that $f_j$ has a zero of order $R^{(d-1)+2j}$ at $R = 0$, $f_j$ extends as an holomorphic function on $U \cap B_{d(d-2)}(0)$ and $R^{-(d-1)}f_j$ has an even expansion around $R = 0$. In fact, $f_j$ extends holomorphically on $U$. If we let $c_1(R) = c(R)$ modulo the part which is singular at $R = d(d-2)$ and $c_2(R) = c(R) - c_1(R)$, then we can write
\begin{align}
    f_j(R) &= -R^{\frac{d-1}{2}} \phi(R) \int_0^R s^{-\frac{d-1}{2}}\phi(s) \left[c_2(s)-c_2(R) \right] f_{j-1}(s)ds \notag \\
    &\phantom{=}+ R^{\frac{d-1}{2}} \phi(R) \int_0^R F_{j-1}(s)c_1'(s)ds, \label{f_j}
\end{align}
and this extends holomorphically on $U$ since multiplication by $\phi(R)$ and $\phi(s)$ removes the singularity. Then one can bound $f_j(z)$ as follows:
\begin{align*}
    |f_j(z)| \leq C \left( (1+|z|) \int_0^{|z|} (1+|s|)^{1-d} |f_{j-1}(s)|ds +  \int_0^{|z|} (1+|s|)^{2-d} |f_{j-1}(s)|ds \right. \\
    \phantom{=} \left. + (1+|z|) \int_0^{|z|} |s|^{2d-3} |s^{-d}F_{j-1}(s)| ds\right) \quad \forall z \in U
\end{align*}
using (\ref{f_j}). Moreover, one has
\begin{align*}
    |f_0(z)| \leq C(1+|z|), \quad |z^{-d} F_0(z)| \leq C (1+|z|)^{-d}\cdot \log(1 + |z|)^{\delta_{d=4}},
\end{align*}
so we deduce
\begin{align*}
    |f_1(z)| &\leq 3C^2 (1+|z|)^{d-1}\cdot \log(1 + |z|)^{\delta_{d=4}}, \\
    |z^{-d}F_1(z)| &\leq 3C^2 (1+|z|)^{2-d}\cdot \log(1 + |z|)^{\delta_{d=4}}.
\end{align*}
It follows by induction that
\begin{align*}
    |f_j(z)| &\leq \frac{3^jC^{j+1}}{(j-1)!} (1+|z|)^{(d-1)+2(j-1)}\cdot \log(1 + |z|)^{\delta_{d=4}}, \\
    |z^{-d}F_j(z)| &\leq \frac{3^jC^{j+1}}{j!} (1+|z|)^{2j-d}\cdot \log(1 + |z|)^{\delta_{d=4}}.
\end{align*}
Finally, set $\phi_j(R^2) = R^{-(d-1)}R^{-2j}f_j(R)$. Note that in dimension $d = 4$, one explicitly has
$$F_0(z) = -\frac{\splitfrac{ z^6(7+\log(512))+24z^4(2+\log(512))+192z^2(1+\log(512))}{-3\left(z^2+8\right)^3\log\left(z^2+8\right)+1536\log(8)}}{6\left(z^2+8\right)^3},$$
so that we get the lower bound
\begin{align*}
|f_1(z)| &\gtrsim \left| z^{\frac{d-1}{2}} \phi(z) \int_1^z F_{j-1}(y)c_1'(y)dy \right| \gtrsim |z| \int_1^{|z|} |F_j(y)| \cdot y^{d-3} dy \\
&\gtrsim |z| \int_1^z \log(y)y^{d-3} \gtrsim |z|^{d-1} \log(|z|)  \quad \forall z \in U, |z| \gtrsim 1.
\end{align*}
In dimension $d > 4$,
$$F_0(z) = \int_0^z \phi(s)^2ds \geq \min_{y \in B_1(0)} \int_0^y \phi(s)^2ds > 0 \quad \forall z \in U, |z| \geq 1$$
and we deduce the lower bound
\begin{align*}
|f_1(z)| &\gtrsim \left| z^{\frac{d-1}{2}} \phi(z) \int_1^z F_{j-1}(y)c_1'(y)dy \right| \gtrsim |z| \int_1^{|z|} |F_j(y)| \cdot y^{d-3} dy \\
&\gtrsim |z|^{d-1}  \quad \forall z \in U, |z| \gtrsim 1.
\end{align*}
\end{proof}

\begin{remark}
    In particular, $\phi(R,z) \in C^{\infty}((0,+\infty) \times U)$, $U =  \{z \in \mathbb C:  \Re(z) > -\frac{1}{2}d(d-2)\}$. Indeed, fix $R_0 > 0, z_0 \in U$. Around $z = R_0^2$,
    $$
    |\phi_j(z)| \leq \frac{C^j}{(j-1)!}, \quad
 |\phi_j^{(n)}(z)| \lesssim_{n,R_0} \frac{C^j}{(j-1)!}
    $$
    using Cauchy's Integral Formula. Next, observe that
    \begin{align*}
        \partial_{R}^{k_1}\partial_z^{k_2} \left[(R^2z)^j \phi_j(R^2) \right] &= \frac{j!}{(j-k_2)!}z^{j-k_2} \sum_{l_1+l_2=k_1} \binom{k_1}{l_1} \frac{(2j)!}{(2j-l_1)!}R^{2j-l_1} \partial_R^{l_2} \left[ \phi_j(R^2)\right] \\
        &= \fsum_{\substack{l_1 \in \mathbb Z \\ 0 \leq l_2 \leq k_1}} c_{l_1,l_2,k_1,k_2} z^{-k_2} R^{-l_1} (R^2z)^{j-k_2} \phi_j^{(l_2)}(R^2),
    \end{align*}
    where the sum is zero if $k_2 > j$. Hence, the sum
    $$
\fsum_{\substack{l_1 \in \mathbb Z \\ 0 \leq l_2 \leq k_1}} \sum_{j=k_2+1}^{\infty} c_{l_1,l_2,k_1,k_2} z^{-k_2} R^{-l_1} (R^2z)^{j-k_2} \phi_j^{(l_2)}(R^2)
    $$
    converges uniformly around $(R,z) = (R_0,z_0)$ since
$$
|c_{l_1,l_2,k_1,k_2} z^{-k_2} R^{-l_1} (R^2z)^{j-k_2} \phi_j^{(l_2)}(R^2)| \lesssim_{l_1,l_2,k_1,k_2,R_0,z_0} \frac{C(R_0,z_0)^j}{(j-1)!}.
$$
Therefore, the sum of derivatives
    $$
    \sum_{j=k_2+1}^{\infty}  \partial_{R}^{k_1}\partial_z^{k_2} \left[(R^2z)^j \phi_j(R^2) \right]
    $$
    converges uniformly for any $k_1, k_2 \in \mathbb N$, which is sufficient to prove the smoothness of $\phi(R,z)$.
\end{remark}
\begin{remark}
There exists $0 < \xi_0 < \delta_0 < \delta_1 \ll 1$ (depending on the absolute constants from the Proposition \ref{fundamental system Lu = zu}) such that for all $0 < \xi < \xi_0$, for all $\delta \in [\delta_0, \delta_1]$, one has
\begin{equation}
    \phi(R,\xi) \gtrsim R^{\frac{d-5}{2}} \cdot \log(1+R^2)^{\delta_{d=4}}  \label{phi lower bound}
\end{equation}
when $R = \delta \xi^{-\frac{1}{2}}$.
\end{remark}

\begin{corollary} \label{phi upper bound R^2z < 1}
When $|R^2z| \lesssim 1$, we have the following pointwise estimate on $\phi(R,z)$ for any $k \geq 0$ and $l \geq 1$:
\begin{align*}
|(R^k\partial_R^k)\phi(R,z)| &\lesssim R^{\frac{d-1}{2}}\langle R \rangle^{-(d-2)} + (R^2z) R^{\frac{d-1}{2}} \langle R^2 \rangle^{-1}  \cdot \log(1+R^2)^{\delta_{d=4}}   \\
|(R^k\partial_R^k)(z^l\partial_z^l)\phi(R,z)| &\lesssim (R^2z)^l R^{\frac{d-1}{2}} \langle R^2 \rangle^{-1}\log(1+R^2)^{\delta_{d=4}}.
\end{align*}
In particular, 
$$|(R^k\partial_R^k)\phi(R,z)| + |(R^k\partial_R^k)(z^l\partial_z^l)\phi(R,z)| \lesssim R^{\frac{d-1}{2}} \langle R \rangle^{-2} \cdot \log(1+R^2)^{\delta_{d=4}}  \lesssim 1$$
is uniformly bounded on $|R^2z| \lesssim 1$ if $d \in \{4,5\}$.
\end{corollary}

\begin{proof}
It suffices to differentiate the series (\ref{phi(R,z)}) and then distinguish the cases $R \leq 1$, $R \geq 1$.
\end{proof}

    


Now, we define the singular $m$-function and state the spectral theorem for the self-adjoint operator $\mathcal{L}$:

\begin{definition}[Singular $m$-function]
Let $m(z): \mathbb C \setminus \mathbb R \rightarrow \mathbb C$, $m(\overline{z}) = \overline{m(z)}$, be the singular Weyl-Titchmarsh $m$-function. It is the unique function for which $\theta(R,z) + m(z)\phi(R,z)$ belongs to $L^2([1,+\infty))$ and solves $\mathcal{L}u = zu$ on $(0,+\infty)$. 
\end{definition}

\begin{theorem}[Spectral theorem] \label{spectral theorem}
The singular $m$-function is a generalized Nevanlinna function which defines a non-negative spectral density $d\rho$ on $\mathbb R$ via
$$\frac{1}{2}(d \rho((a,b)) + d \rho([a,b]))= \lim \limits_{\varepsilon \to 0} \frac{1}{\pi} \int_a^b \Im m(t + i \varepsilon) dt$$
such that the generalized Fourier transform
\begin{align*}
\mathcal{F}: L^2((0,+\infty)) &\rightarrow L^2(\mathbb R, d \rho) \\
f &\mapsto \hat{f}(\xi):= \lim \limits_{r \to +\infty} \int_0^r \phi(s,\xi)f(s)ds 
\end{align*}
\index{F@$\mathcal{F}$, the generalized Fourier transform}
is a unitary operator with inverse
\begin{align*}
\mathcal{F}^{-1}: L^2(\mathbb R, d \rho) &\rightarrow  L^2((0,+\infty)) \\
F &\mapsto \check{F}(R):= \lim \limits_{r \to +\infty} \int_{-r}^r \phi(R,\xi)F(\xi)d\rho(\xi) .
\end{align*}
Here, the limits must be understood as limits of functions in their respective $L^2$-space.

Moreover, if $E$ is the unique spectral family associated to the self-adjoint operator $\mathcal{L}$ on $\dom(\mathcal{L})$ and, for $f \in L^2((0,+\infty))$, $d \mu_{f}$ is its spectral measure, then
$$d \mu_f = |\hat{f}|^2 d \rho$$
\end{theorem}

\begin{proof}
See \cite[Lemma 3.4]{Gesztesy}, \cite[Theorem 3.4, Corollary 3.5]{Kostenko1} and \cite[Theorem 4.5]{Kostenko2}. 
\end{proof}

\begin{remark}
     If $\xi^*$ is an eigenvalue, the inverse Fourier transform of $\delta_{\xi^*}$ is a multiple of the eigenfunction $\phi(R,\xi^*)$. In other words, the $L^2(\mathbb R, d\rho)$-limit
    $$\lim \limits_{r \to +\infty} \int_{0}^r \phi(R,\xi)\phi(R,\xi^*)dR$$
    is $0$ d$\rho$-almost everywhere on $\{\xi \in \mathbb R: \xi \neq \xi^*\}$.
\end{remark}

\begin{notation}\label{not: notation of inner product for L^2 limit}
    We will write
    $$\innerp[\Big]{ f(s) }{ \phi(s,\xi)}_{L^2((0,+\infty))}:= \lim\limits_{r \to +\infty} \int_0^r \phi(s,\xi)f(s)ds$$
   as a limit of functions in $L^2(\mathbb R, d\rho)$,  even though
    $$\int_0^{+\infty} |\phi(s,\xi)f(s)|ds$$
    need not to be finite.
\end{notation}

\begin{remark}[On the decomposition of $d\rho$]
Since $\spec(\mathcal L) = \{\xi_d, 0\} \cup [0,+\infty) = \spec_{pp}(\mathcal{L}) \cup \spec_{ac}(\mathcal{L})$, we can write
$$
d\rho(\xi) = \frac{1}{||\phi(R,\xi_d)||_{L^2}} \delta_{\xi_d}(\xi) + \frac{1}{||\phi(R)||_{L^2}} \delta_{0}(\xi) + \rho(\xi) \chi_{(0,+\infty)}(\xi) d\xi
$$ for some $\rho(\xi) \in L^1_{loc}([0,+\infty))$.
\index{rho@$\rho(\xi)d\xi$, the absolutely continuous part of the spectral measure $d\rho(\xi)$}
\end{remark}

Next, we introduce the Jost solution which will be useful in the computation of asymptotics for $\rho(\xi)$. Note that in the following definition, we use the principal branch of the complex logarithm in order to define roots.

\begin{definition}[Jost solution]
For $z \in \mathbb C \setminus \mathbb R_{\leq 0}$, $\Im z \geq 0$, let $\psi^{+}(R,z)$ denote the Jost solution to $\mathcal{L}u = zu$ at $R = +\infty$ normalized so that 
$$\psi^+(R,z) \sim z^{-\frac{1}{4}}e^{iR\sqrt{z}}, \quad R|\sqrt{z}| \to +\infty.$$
It is given by $z^{-\frac{1}{4}}f_+(R,z)$, where $f_+(R,z)$ is the unique fixed point of 
\begin{equation}
    f_+(R,z) = e^{iR\sqrt{z}} - \int_R^{+\infty}\frac{\sin(\sqrt{z}(R-R'))}{\sqrt{z}}V(R')f_+(R',z)dR', \quad R > 0, \Im z \geq 0, z \neq 0, \label{f_+ definition}
\end{equation}
where $\mathcal{L} = -\partial_{RR} + V(R)$. 
For its construction, see \cite[Section 12.1.1]{Newton}.
\index{psi plus@$\psi^+(R,\xi)$, Jost solution for $\mathcal{L}u=zu$}
\end{definition}


Next, we give an approximation for $\psi^+(R,\xi)$ which is useful when $R^2 \xi \gtrsim 1$:

\begin{proposition}\label{prop:psi^+ in terms of sigma}
For $\xi > 0$, $R^2 \xi \gtrsim 1$, $\psi^+(R,\xi)$ is of the form
$$\psi^+(R, \xi) = \xi^{-\frac{1}{4}}e^{iR\xi^{\frac{1}{2}}} \sigma(R\xi^{\frac{1}{2}},R)$$
where $\sigma(q,r)$ is well-approximated by the series
$$\sigma(q,r) \approx \sum_{j=0}^{+\infty}q^{-j}\psi_j^+(r)$$
for some zeroth-order symbol $\psi_j^+$ being analytic on $(0,+\infty]$, i.e.,
$$\sup_{r > 0}|(r\partial_r)^k \psi_j^+(r)| < +\infty \ \forall k \in \mathbb N_{\geq 0},$$
in the following sense:
$$\sup_{r > 0} \left| (r \partial_r)^{\alpha} (q \partial_q)^{\beta} \left[ \sigma(q,r) - \sum_{j=0}^{j_0}q^{-j}\psi_j^+(r) \right] \right| \leq c_{\alpha, \beta, j_0} q^{-j_0-1}, \quad \forall q \geq 1,$$
for any $\alpha, \beta \in \mathbb N_{\geq 0}$, for any $j_0$ large enough.
\end{proposition}

\begin{proof}
See \cite[Proposition 4.6]{Krieger_2007}.
\end{proof}

\begin{remark}\label{rmk: boundedness of sigma}
    In particular,
    $$\sup_{r > 0, q > 1}| (r \partial_r)^{\alpha} (q \partial_q)^{\beta}\sigma(q,r)| < +\infty \ \forall \alpha,\beta \in \mathbb N_{\geq 0}.$$
\end{remark}

\begin{corollary}\label{cor: psi^+ estimate}
For all $R, \xi > 0$, $R^2 \xi \gtrsim 1$, the following pointwise estimates hold for $\psi^+(R,\xi)$:
\begin{align}
    |(\xi^l\partial_{\xi}^l) (R^k \partial_R^k) \psi^+(R,\xi)| &\lesssim \xi^{-\frac{1}{4}} (R\xi^{\frac{1}{2}})^{l+k} \quad \forall k, l \geq 0. \label{psi_+ upper bound}
\end{align}
\end{corollary}

\begin{proof}
Observe that
\begin{align*}
     (\xi\partial_{\xi})F(R\xi^{\frac{1}{2}},R) &= \frac{1}{2} R\xi^{\frac{1}{2}}\partial_q F(R\xi^{\frac{1}{2}},R) = G(R\xi^{\frac{1}{2}},R), \quad G(q,r) = \frac{1}{2} q \partial_q F, \\
  (R\partial_{R})F(R\xi^{\frac{1}{2}},R) &= R\xi^{\frac{1}{2}}\partial_q F(R\xi^{\frac{1}{2}},R) + R \partial_R F(R\xi^{\frac{1}{2}},R) = H(R\xi^{\frac{1}{2}},R), \quad H = (q \partial_q  + r \partial_r) F,
\end{align*}
so that, by induction,
$$
 |(\xi\partial_{\xi})^l (R \partial_R)^k \sigma(R\xi^{\frac{1}{2}},R)| \lesssim  \sup_{\alpha \leq l, \beta \leq k} \sup_{r > 0, q > 1}| (r \partial_r)^{\alpha} (q \partial_q)^{\beta}\sigma(q,r)| < +\infty
$$
The same inequality holds with $(\xi^l \partial_{\xi}^l) (R^k \partial_R^k)$, as it is a linear combination of the differential operators $(\xi\partial_{\xi})^i (R \partial_R)^j$ for $i \leq l, j \leq k$. One also checks that
\begin{align*}
    (R^k\partial_R^k)(\xi^{-\frac{1}{4}}e^{iR\xi^{\frac{1}{2}}}) &= i^ke^{iR\xi^{\frac{1}{2}}} \xi^{-\frac{1}{4}}(R\xi^{\frac{1}{2}})^k, \\
    \left| (\xi^l \partial_{\xi}^l) (\xi^{\alpha} e^{iR\xi^{\frac{1}{2}}}) \right| &\lesssim \xi^{\alpha}(R\xi^{\frac{1}{2}})^l, \\
    \left|(\xi^l \partial_{\xi}^l)  [(R^k\partial_R^k)(\xi^{-\frac{1}{4}}e^{iR\xi^{\frac{1}{2}}})] \right| &\lesssim \xi^{-\frac{1}{4}}(R\xi^{\frac{1}{2}})^{k+l},
\end{align*}
which yields the result $\psi^+(R,\xi) = \xi^{-\frac{1}{4}}e^{iR\xi^{\frac{1}{2}}} \sigma(R\xi^{\frac{1}{2}},R)$ using the product rule.
\end{proof}

Now, we are ready to give growth estimates on the spectral density $\rho(\xi)$.

\begin{proposition}\label{prop:spectral density properties}
For $R > 0$ and $\xi > 0$, we have
\begin{equation*}
    \phi(R,\xi) = a(\xi)\psi^{+}(R,\xi) + \overline{a(\xi)\psi^{+}(R,\xi)},
\end{equation*}
where $a(\xi)$ is smooth, non-zero, has asymptotics
$$|a(\xi)| \asymp \begin{cases} \xi^{\frac{6-d}{4}} \cdot |\log(\xi)|^{\delta_{d=4}}, \ &\xi \ll 1 \\
\xi^{\frac{2-d}{4}}, \ &\xi \gg 1
\end{cases}$$
and symbol-type upper bounds
$$|(\xi \partial_{\xi})^k a(\xi)| \leq c_k|a(\xi)| \ \forall \xi > 0.$$
Moreover,
$$\rho(\xi) = \frac{1}{\pi |a(\xi)|^2}$$
and the corresponding asymptotics are
\begin{equation}
    |\rho(\xi)| \asymp \begin{cases} \xi^{\frac{d}{2} - 3} \cdot |\log(\xi)|^{-2 \delta_{d = 4}}, \ &\xi \ll 1 \\
 \xi^{\frac{d}{2} - 1}, \ &\xi \gg 1
\end{cases}
\end{equation}
with symbol-type upper bounds 
$$
|(\xi \partial_{\xi})^k \rho(\xi)| \leq \tilde{c}_k|\rho(\xi)| \ \forall \xi > 0.
$$
\index{rho@$\rho(\xi)d\xi$, the absolutely continuous part of the spectral measure $d\rho(\xi)$}
\end{proposition}

\begin{proof}
Following \cite[Proposition 4.7]{Krieger_2007}, we find that
\begin{align}
    a(\xi) &= -\frac{i}{2}W(\phi(\cdot, \xi), \psi^{+}(\cdot, \xi)), \label{a(xi) upper bound} \\
    |a(\xi)| &\geq \frac{|\partial_R \phi(R,\xi)|}{2|\partial_R \psi^{+}(R,\xi)|}, \label{a(xi) lower bound} \\
    \rho(\xi) &= \frac{1}{\pi |a(\xi)|^2}. \notag
\end{align}
The behaviour of $\rho(\xi)$ for large $\xi$ is well-known: see \cite[Theorem 2.1]{Kostenko_2013}.  

For small $\xi$, we proceed as in \cite[Proposition 4.7]{Krieger_2007}. We take $R = \delta \xi^{-\frac{1}{2}}$ as in (\ref{phi lower bound}) ($\delta$ is fixed and $\xi \to 0^+$) so that one gets
\begin{align*}
    |(R \partial_R)^i \phi(R, \xi)| &\sim \xi^{\frac{5-d}{4}} \cdot \log(1 + \xi^{-1})^{\delta_{d = 4}}, \\
    |(R \partial_R)^i \psi^+(R,\xi)| &\lesssim \xi^{-\frac{1}{4}},
\end{align*}
for $i = 0, 1$ using (\ref{phi(R,z)}), (\ref{phi lower bound}) and (\ref{psi_+ upper bound}). We conclude by applying these estimates with (\ref{a(xi) upper bound}) and (\ref{a(xi) lower bound}).
\end{proof}

\begin{corollary} \label{phi upper bound R^2z > 1}
When $R^2 \xi \gtrsim 1$, the following pointwise estimates hold for $\phi(R,\xi)$ for any $k, l \geq 0$:
\begin{align*}
    |\xi^l\partial_{\xi}^l R^k \partial_R^k \phi(R,\xi)| &\lesssim \xi^{\frac{5-d}{4}} (R\xi^{\frac{1}{2}})^{l+k} \langle \xi \rangle^{-1} \cdot \left( 1_{0 < \xi < 1/2}(\xi) |\log(\xi)| \right)^{\delta_{d = 4}}. 
\end{align*}
The logarithm appears only in dimension $d = 4$ for small $\xi > 0$. In particular, in dimension $d \in \{4,5\}$, one has
$$|\phi(R,\xi)| \lesssim \langle \xi \rangle^{\frac{1-d}{4}}, \quad R^2 \xi \gtrsim 1.$$
\end{corollary}

\begin{proof}
    Write $\phi = 2 \Re( a(\xi) \psi^+(R,\xi))$ and use the estimates from Corollary \ref{cor: psi^+ estimate}, as well as 
    $$
    |(\xi \partial_{\xi})^k a(\xi)| \leq c_k|a(\xi)| \lesssim \xi^{\frac{6-d}{2}} \cdot \left( 1_{0 < \xi < 1/2}(\xi) |\log(\xi)| \right)^{\delta_{d = 4}}, \quad \xi > 0
    $$
    from Proposition \ref{prop:spectral density properties}.
\end{proof}



\begin{corollary} \label{phi upper bound}
Assume $d \in \{4,5\}$. Fix $0 < \xi_0 < 1$. For all $R, \xi > 0$, the following pointwise estimates hold for $\phi(R,\xi)$ when $l \geq 1$:
\begin{align*}
|\phi(R,\xi)| &\lesssim \langle \xi \rangle^{\frac{1-d}{4}}, \\
|R\partial_R \phi(R,\xi)| &\lesssim \min\{R\xi^{\frac{3-d}{2}},R^{\frac{d-1}{2}}\} \text{ if } \xi > 1,\\
|\partial_{\xi}^l \phi(R,\xi)| &\lesssim \min\{R^{\frac{d-1}{2}+2l},\xi^{\frac{1-d}{4}}(R\xi^{-\frac{1}{2}})^l\} \text{ if } \xi > \xi_0,  \\
    |\partial_{\xi}^l \phi(R,\xi)| &\lesssim \min\{ R^{\frac{d-1}{2}+2(l-1) }\log(1+R^2)^{\delta_{d = 4}},\xi^{\frac{5-d}{4}}(R\xi^{-\frac{1}{2}})^l |\log(\xi)|^{\delta_{d = 4}} \} \text{ if } \xi < \xi_0.
\end{align*}
\end{corollary}

\begin{proof}
    This is a combination of Corollary \ref{phi upper bound R^2z < 1} and Corollary \ref{phi upper bound R^2z > 1}.
\end{proof}

\begin{corollary} \label{F(xi,eta) estimates}
Assume $d \in \{4,5\}$. Let 
\begin{align*}
    W_0(R) &= [\mathcal{L},R\partial_R]-2\mathcal{L} = [V,R\partial_R]-2V(R) = -2V(R) - R\partial_RV(R) \\
    &= \frac{2d^2\left(d^2-4\right)\left((d-2)d-R^2\right)}{\left((d-2)d+R^2\right)^3},
\end{align*}
where $\mathcal{L} = -\partial_{RR} + V(R)$ (not to be confused with the ground state $W(x)$). The symmetric function
    $$|F(\xi,\eta)| = \innerp[\Big]{ W_0(R) \phi(R,\xi) }{ \phi(R,\eta) }_{L^2((0,+\infty))}$$
    is of class $C^1([0,+\infty)\times[0,+\infty)) \cap C^2((0,+\infty)\times(0,+\infty))$ and satisfies:
    \begin{align*}
        |F(\xi,\eta)| &\lesssim \begin{cases}
            \xi + \eta \quad &\text{ if } \xi + \eta \leq 1 \\
             (\xi+\eta)^{\frac{1-d}{2}}(1+|\xi^{\frac{1}{2}}-\eta^{\frac{1}{2}}|)^{-N} \quad &\text{ if } \xi + \eta \geq 1 
        \end{cases} \\
        |\partial_{\xi}F(\xi,\eta)| + |\partial_{\eta}F(\xi,\eta)|&\lesssim \begin{cases}
            1 \quad &\text{ if } \xi + \eta \leq 1 \\
             (\xi+\eta)^{\frac{-d}{2}}(1+|\xi^{\frac{1}{2}}-\eta^{\frac{1}{2}}|)^{-N} \quad &\text{ if } \xi + \eta \geq 1 
        \end{cases} \\
        |\partial_{\xi}\partial_{\eta}F(\xi,\eta)| &\lesssim \begin{cases}
            |\log(\xi + \eta)|^3 \quad &\text{ if } \xi + \eta \leq 1, \ d = 4 \\
           |\xi + \eta|^{-\frac{1}{2}}(1 + |\log(\xi/\eta)|) \quad &\text{ if } \xi + \eta \leq 1, \ d = 5 \\
             (\xi+\eta)^{\frac{-1-d}{2}}(1+|\xi^{\frac{1}{2}}-\eta^{\frac{1}{2}}|)^{-N} \quad &\text{ if } \xi + \eta \geq 1 
        \end{cases} 
    \end{align*} 
    for any fixed $N \in \mathbb N$.
    \index{F(xi,eta)@$F(\xi,\eta)$, a symmetric function which contributes to $\mathcal{K}$}
\end{corollary}

\begin{proof}
    The bounds from Corollary \ref{phi upper bound} imply
    \begin{align*}
        |F(\xi,\eta)| &\lesssim \langle \xi \rangle^{\frac{1-d}{4}} \langle \eta \rangle^{\frac{1-d}{4}}, \\
        |\partial_{\xi}F(\xi,\eta)| &\lesssim \langle \xi \rangle^{\frac{-1-d}{4}} \langle \eta \rangle^{\frac{1-d}{4}}, \\
        |\partial^2_{\xi \eta}F(\xi,\eta)| &\lesssim \xi^{\frac{-1-d}{4}}\eta^{\frac{-1-d}{4}} \text{ if } \xi > \frac{1}{2}, \eta > \frac{1}{2},
    \end{align*}
    which yield the desired bounds when $\xi \sim \eta$, $\xi + \eta \geq 1$. The same bounds from Corollary \ref{phi upper bound} combined with Dominated Convergence proves the $C^2$-regularity of $F(\xi,\eta)$, as well as the right-continuity at zero of $F(\xi,\eta)$, $\partial_{\xi}F$ and $\partial_{\eta}F$.

If $\xi + \eta \geq 1$ and $\xi, \eta$ are separated, then we proceed as in \cite[Theorem 5.1]{Krieger_2007}: Writing the integrals as limits, an integration by parts argument shows that 
\begin{equation*}
   \eta F(\xi,\eta)  = \innerp[\Big]{ W_0(R)\phi(R,\xi)}{\mathcal{L}\phi(R,\eta) }_{L^2_R} = \innerp[\Big]{ [\mathcal{L},W_0(R)]\phi(R,\xi)}{\phi(R,\eta) }_{L^2_R} + \xi F(\xi,\eta),
\end{equation*}
i.e.,
\begin{equation*}\label{integration by parts F(xi,eta)}
   (\eta-\xi) F(\xi,\eta)  = - \innerp[\Big]{ (2W_{0,R}(R)\partial_R - W_{0,RR}(R))\phi(R,\xi)}{\phi(R,\eta) }_{L^2_R}.
\end{equation*}
For fixed $\xi, \eta \geq 0$, the integration by parts is justified because $\phi(R,\eta)$, $\partial_R\phi(R,\eta)$ are bounded and $W_0(R)\phi(R,\xi)$, $\partial_{R}(W_0(R)\phi(R,\xi))$ vanish at zero and infinity.

By iteration, for arbitrary $k \in \mathbb N$, there exists rational functions $W_j^{\mathrm{odd}}(R)$, $W_j^{\mathrm{even}}(R) \in C^{\omega}([0,+\infty))$, $0 \leq j \leq k$, decaying as
$$\langle R \rangle |W_j^{\mathrm{odd}}(R)| + |W_j^{\mathrm{even}}(R)| \lesssim \langle R \rangle^{4-2k},$$
and respectively having odd/even expansions at $R = 0$, for which
\begin{equation}\label{iteration of integration by parts F(xi,eta)}
    (\eta-\xi)^{2k}F(\xi,\eta) = \innerp[\Big]{ \Bigl(\sum_{j=0}^{k-1}\xi^j W_j^{\mathrm{odd}}(R)\partial_R + \sum_{j=0}^{k}\xi^j W_j^{\mathrm{even}}(R) \Bigr) \phi(R,\xi) }{ \phi(R,\eta) }_{L^2_R}.
\end{equation}
This implies the desired bound for $F(\xi,\eta)$ and then differentiating with respect to $\xi$ and/or $\eta$ implies the other ones. If $\xi + \eta \lesssim 1$, one observes that $F(0,0) = 0$ because $$W_0(R) \phi(R,0) = \left([\mathcal{L},R\partial_R]-2\mathcal{L}\right)\phi(R,0) = -\mathcal{L}R\partial_R \phi(R,0),$$
    thus integration by parts yields
    $$F(0,0) = \innerp[\Big]{ -R\partial_R\phi(R,0) }{ \mathcal{L}\phi(R,0) }_{L^2_R} = 0.$$
    Combining this with the differentiability at zero, we obtain the bounds for $F(\xi,\eta)$, $\partial_{\xi}F$ and $\partial_{\eta}F$. It remains to prove the estimates for the second derivative. Since the case $d = 4$ has been treated in \cite[Theorem 5.1]{Krieger_2007}, we assume $d = 5$ but the strategy is the same. Using Corollary \ref{phi upper bound}, for $0 < \eta \leq \xi < 1/2$, one gets 
    \begin{align*}
        |\partial_{\xi}\partial_{\eta}F(\xi,\eta)| &\lesssim \int_{0}^{+\infty}\langle R \rangle^{-4} |\partial_{\xi} \phi(R,\xi)\partial_{\eta} \phi(R,\eta)| dR \\
        &\lesssim \int_{0}^{\xi^{-\frac{1}{2}}}\langle R \rangle^{-4} R^{4}dR +  \int_{\xi^{-\frac{1}{2}}}^{\eta^{-\frac{1}{2}}}\langle R \rangle^{-4} R^2 (R\xi^{-\frac{1}{2}}) dR  \\
        &\phantom{\lesssim}+ \int_{\eta^{-\frac{1}{2}}}^{+\infty} \langle R \rangle^{-4}(R\xi^{-\frac{1}{2}})  (R\eta^{-\frac{1}{2}})dR \\
        &\lesssim \xi^{-\frac{1}{2}} \cdot (1 +|\log(\xi/\eta)|) \lesssim (\xi+\eta)^{-\frac{1}{2}} \cdot (1 + |\log(\xi/\eta)|) 
    \end{align*}
    and we conclude by symmetry. Similarly, we can show that the derivatives $\partial_{\xi}^2F$ and $\partial_{\eta}^2F$ exist if $0 < \eta, \xi \leq 1/2$, but we do not need to estimate them.
\end{proof}

\section{The Transference Identity}\label{section:transferrance identity}

Assume $d \in \{4,5\}$ for the remainder of this paper. We are interested in studying the error one makes when passing from $\mathcal{F}(R \partial_R u)$ to $-2\xi\partial_{\xi}\mathcal{F}(u)$, as the operator $R \partial_R$ appears naturally when solving for $\varepsilon$ in Section \ref{section:exact solution fourier method}. This will allow us to translate the wave equation from physical to (generalized) Fourier space where the operator $\mathcal{L}$ is replaced by a multiplication by $\xi$ and one can ``interchange'' the operators $R \partial_R$ and $\mathcal{F}$ up to a controllable error. 
 
Let $L^2(\mathbb R, d\rho)$ denote the set of $d\rho$-measurable functions $f(\xi)$ that are square-integrable. These functions admit the following representation $d\rho$ almost-everywhere:
$$
f(\xi) = f(\xi_d) \delta_{\xi_d} + f(0) \delta_{0, d = 5} + f(\xi)1_{\xi > 0}(\xi),
$$
where $f_c(\xi):= f(\xi)1_{\xi > 0}(\xi) \in L^2((0,+\infty),\rho(\xi)d\xi)$. Let $C_c^{\infty}(\spec_{pp}(\mathcal{L}) \cup (0,+\infty))$ denote the subset of $L^2(\mathbb R, d\rho)$ for which $f_c(\xi) = f(\xi)1_{\xi > 0}(\xi) \in C_c^{\infty}((0,+\infty))$.  Our goal is to study the difference operator $\mathcal{K}$:
\begin{align*}
    \mathcal{K}:   C_c^{\infty}(\spec_{pp}(\mathcal{L}) \cup (0,+\infty)) &\mapsto L^2(\mathbb R,d\rho) \\
    f(\xi)  &\mapsto \mathcal{K}(f):=  \mathcal{F}(R \partial_R \mathcal{F}^{-1}f) + \mathcal{F}(\mathcal{F}^{-1} 2\xi \partial_{\xi} f),
\end{align*}
\index{K@$\mathcal{K}$, the transference operator}
where $\xi\partial_{\xi}$ acts as zero on the discrete component, and show that this is a well-defined bounded operator inbetween some weighted $L^2$-spaces.  As a first step, we are going to show that $\mathcal{K}$ is well-defined. 
\begingroup
\allowdisplaybreaks
\begin{proposition}
    The operator $\mathcal{K}$ is well-defined.
\end{proposition}

\begin{proof}
Denoting by $\phi_d(R)$ and $\phi_0(R)$ the normalized eigenfunctions, observe that
\begin{align*}
    \mathcal{F}^{-1} f &= f(\xi_d)  \mathcal{F}^{-1} \delta_{\xi_d} + f(0) \mathcal{F}^{-1} \delta_{0, d = 5} + \mathcal{F}^{-1} f_c  \\
    &= f(\xi_d) \phi_d(R) + f(0) \phi_0(R) \delta_{d = 5} + \mathcal{F}^{-1} f_c, \\
     \mathcal{F} (R \partial_R \mathcal{F}^{-1} f) &= f(\xi_d) \mathcal{F} (R \partial_R \phi_d) + f(0) \mathcal{F} (R \partial_R \phi_0)\delta_{d = 5} + \mathcal{F} (R \partial_R \mathcal{F}^{-1} f_c), \\
     \mathcal{F}(\mathcal{F}^{-1} 2\xi \partial_{\xi} f) &=  2\xi \partial_{\xi} f_c,
\end{align*}
where we recall that we write
    $$
    (\mathcal{F}f)(\xi) = \innerp[\Big]{ f(s) }{ \phi(s,\xi)}_{L^2((0,+\infty))}:= \lim\limits_{r \to +\infty} \int_0^r \phi(s,\xi)f(s)ds
    $$
   as a limit of functions in $L^2(\mathbb R, d\rho)$,  even though
    $$\int_0^{+\infty} |\phi(s,\xi)f(s)|ds$$
    need not to be finite.  In order to show that $\mathcal{K}$ is well-defined, it is necessary to show that $\mathcal{F}(R\partial_R \phi_d(R))$, $\mathcal{F}(R\partial_R \phi_0(R))$ and $\mathcal{F}(R\partial_R \mathcal{F}^{-1})$ are well-defined. To this end, it suffices to show that $R\partial_R \phi_d(R)$, $R\partial_R \phi_0(R)$ and $R\partial_R \mathcal{F}^{-1} f$ are $L^2((0,+\infty))$-functions in Lebesgue sense (Theorem \ref{spectral theorem}). Jost solution theory shows that $\phi_d(R)$ decays exponentially as $R \to +\infty$ and, by (\ref{phi theta asympotics}), $\phi_0(R)$ decays as $R^{-1}$ in dimension $d = 5$ where it appears with symbol-type behaviour of the derivatives. Hence, $R\partial_R \phi_d(R)$ and  $R\partial_R \phi_0(R)$  are $L^2$-functions. If $f \in C_c^{\infty}((0,+\infty))$, then our estimates on $R\partial_R\phi(R,\xi)$ from Corollary \ref{phi upper bound} shows that $R\partial_R\mathcal{F}^{-1}f$ is, a priori, only bounded. Yet, $\mathcal{F}^{-1}f$ decays like a Schwartz function: one can get rid of powers of $R$ using successive integration by parts as we show in Lemma \ref{lemma:F^-1 f is schwartz if f is smooth with compact support}.
\end{proof}
\endgroup
\begin{lemma}\label{lemma:F^-1 f is schwartz if f is smooth with compact support}
    Let $f \in C_c^{\infty}((0,+\infty))$. Then $\mathcal{F}^{-1}f(R)$ has an arbitrary fast polynomial decay at infinity. 
\end{lemma}

\begin{proof}
    Let $K \subset (0,+\infty)$ be the compact support of $f$ and
    $$
\mathcal{F}^{-1} f(R) = \int_{0}^{+\infty} \phi(R,\xi)f(\xi)d\rho(\xi) = \int_{K} \phi(R,\xi)f(\xi)\rho(\xi)d\xi,
    $$
    where $\rho(\xi)$ is smooth on $(0,+\infty)$ (Corollary \ref{prop:spectral density properties}). As $\phi(R,\xi)$ is smooth on $(0,+\infty) \times (0,+\infty)$, $K \subset (0,+\infty)$ is compact, one can interchange derivative and integral thanks to Dominated Convergence, i.e.,
    \begin{align*}
            (R\partial_R)^n\mathcal{F}^{-1} f(R) &=  \int_{K} (R\partial_R)^n\phi(R,\xi)f(\xi)\rho(\xi)d\xi \\
            &= 2 \Re \left(  \int_{K} a(\xi) (R\partial_R)^n\left[ \psi^+(R,\xi) \right] f(\xi)\rho(\xi)d\xi \right). 
    \end{align*}
    Observe that
 \begin{align*}
           (R\partial_{R})e^{iR\xi^{\frac{1}{2}}} F_0(R\xi^{\frac{1}{2}},R) &= e^{iR\xi^{\frac{1}{2}}} \left[ iR \xi^{\frac{1}{2}}F(R\xi^{\frac{1}{2}},R) + R\xi^{\frac{1}{2}}\partial_q F(R\xi^{\frac{1}{2}},R) + R \partial_R F(R\xi^{\frac{1}{2}},R) \right] \\
      &= e^{iR\xi^{\frac{1}{2}}} H(R\xi^{\frac{1}{2}},R), \quad H = (iq + q \partial_q  + r \partial_r) F_0.
 \end{align*}
If $F_0(q,r) = \sigma(q,r)$ is the function coming from $\psi^+$, then $(q \partial_q)^i(r \partial_r)^{j} \sigma$ is bounded on $(q,r) \in [1,+\infty) \times (0,+\infty)$ for all fixed $i,j \geq 0$ (Remark \ref{rmk: boundedness of sigma}). Hence,
 \begin{align*}
           (R\partial_{R})e^{iR\xi^{\frac{1}{2}}} \sigma(R\xi^{\frac{1}{2}},R) &= e^{iR\xi^{\frac{1}{2}}} H(R\xi^{\frac{1}{2}},R), \quad H = iqF_0 + F_1, \quad F_1 = (q \partial_q  + r \partial_r) F_0,
 \end{align*}
where $(q \partial_q)^i(r \partial_r)^{j} F_1$ is bounded on $(q,r) \in [1,+\infty) \times (0,+\infty)$. By induction, it holds that 
\begin{align*}
           (R\partial_{R})^ne^{iR\xi^{\frac{1}{2}}} \sigma(R\xi^{\frac{1}{2}},R) &= e^{iR\xi^{\frac{1}{2}}} H(R\xi^{\frac{1}{2}},R), \\
           H &= (iq)^nF_0 + (iq)^{n-1}F_1 + ... + (iq)^{n-1}F_{n-1} + F_n,
 \end{align*}
 where $(q \partial_q)^i(r \partial_r)^{j} F_k$ is bounded on $(q,r) \in [1,+\infty) \times (0,+\infty)$ for all $i,j \geq 0$ and all $0 \leq k \leq n$.  It is now sufficient to prove that for any smooth $F$ as above, there exists $C = C(F,f,a,\rho,K,n) > 0$ for which
 $$
\sup_{R > 0} \left| \int_{K} e^{iR\xi^{\frac{1}{2}}} (iR\xi^{\frac{1}{2}})^n F(R\xi^{\frac{1}{2}},R) a(\xi)f(\xi)\rho(\xi)d\xi \right| \leq C.
 $$
 Write $a(\xi) f(\xi) \rho(\xi) = \tilde{f}(\xi)$. First, observe that
 $$
(\xi \partial_{\xi})^n \left[ G(R\xi^{\frac{1}{2}}) \right] = \frac{1}{2^n} [(q\partial_q)^n G](R\xi^{\frac{1}{2}}).
 $$
 Applying this identity with $G(q) = e^{iq}$, one obtains
 $$
(\xi \partial_{\xi})^n \left[  e^{iR\xi^{\frac{1}{2}}} \right] = \frac{i^n}{2^n} (R\xi^{\frac{1}{2}})^n e^{iR\xi^{\frac{1}{2}}}. 
 $$
 Hence,
\begin{align*}
  \int_{K} e^{iR\xi^{\frac{1}{2}}} (iR\xi^{\frac{1}{2}})^{n+1} F(R\xi^{\frac{1}{2}},R) \tilde{f}(\xi) d\xi \\
  = c_n \xi \int_{0}^{+\infty} (\xi \partial_{\xi})^{n+1} \left[  e^{iR\xi^{\frac{1}{2}}} \right] F(R\xi^{\frac{1}{2}},R)  \tilde{f}(\xi) d\xi \\
  =  c_n \xi \int_{0}^{+\infty} e^{iR\xi^{\frac{1}{2}}} (-1-\xi\partial_{\xi})^{n+1} \left[ F(R\xi^{\frac{1}{2}},R)  \tilde{f}(\xi) \right] d\xi.
\end{align*}
As 
$$
(\xi\partial_{\xi})^i \left[ F(R\xi^{\frac{1}{2}},R)  \tilde{f}(\xi) \right] = \sum_{i_1 + i_2 = i} \binom{i}{i_1,i_2} (\xi\partial_{\xi})^{i_1} \tilde{f}(\xi) \frac{1}{2^{i_2}} [(q\partial_q)^{i_2} F](R\xi^{\frac{1}{2}},R)
$$
is bounded for $\xi \in K, R  >0$, this finishes the proof.
\end{proof}

  Our next goal is to prove boundedness of $\mathcal{K}$ on some appropriate weighted $L^2$-spaces by finding and analysing its kernel. Representing $L^2(\mathbb R, d\rho)$ as $\mathbb R^{d-3} \times L^2(\mathbb R, \rho(\xi)d\xi)$ using the natural map
$$
f(\xi) = f(\xi_d) \delta_{\xi_d} + f(0) \delta_{0, d = 5} + f_c(\xi) \mapsto (f(\xi_d),f(\xi_0),f_c(\xi)),
$$
observe that
\begin{align*}
     \mathcal{K}f &= f(\xi_d) \mathcal{F} (R \partial_R \phi_d) + f(0) \mathcal{F} (R \partial_R \phi_0)\delta_{d = 5} + \mathcal{F} (R \partial_R \mathcal{F}^{-1} f_c)  + 2\xi\partial_{\xi}f_c \\
     &= f(\xi_d) \mathcal{K}_d + f(0) \mathcal{K}_0 \delta_{d = 5} + \mathcal{K}_c(f_c) \\
     &= (\mathcal{K}_d, \mathcal{K}_0, \mathcal{K}_c) \cdot (f(\xi_d), f(0), f_c),
\end{align*}  
where ``d'' represents the negative discrete eigenvalue $\xi_d$, ``0'' stands for the 0-eigenvalue when $d = 5$ and ``c'' represents the continuous part of the spectrum. Extracting the discrete and continuous components from each Fourier transform (the discrete component is obtained by evaluating the Fourier transform at $\xi = \xi_d$ and $\xi = 0$), the operator $\mathcal{K}$ admits, respectively in $d = 4$ and $d = 5$, the following matrix representation
$$\begin{pmatrix} \mathcal{K}_{dd} & \mathcal{K}_{cd} \\
\mathcal{K}_{dc} & \mathcal{K}_{cc} \\
\end{pmatrix}, \quad \begin{pmatrix} \mathcal{K}_{dd} & \mathcal{K}_{0d} & \mathcal{K}_{cd} \\
\mathcal{K}_{d0} & \mathcal{K}_{00} & \mathcal{K}_{c0} \\
\mathcal{K}_{dc} & \mathcal{K}_{0c} & \mathcal{K}_{cc}
\end{pmatrix},$$
where for $x \in \{ ``d``,  ``c`` ,  ``0``\}$, $K_{xd}(\cdot)$ represents the evaluation of $\mathcal{K}_x(\cdot)(\xi)$ at $\xi = \xi_d$, $K_{x0}$ represents the value at $\xi = 0$ and $K_{xc}$ represents the continuous component of the transform. More precisely, 
\begin{alignat*}{2}
    \mathcal{K}_{dd} &= \Bigl< R\partial_R \phi_d(R),\phi_d(R) \Bigr>_{L^2_R}, \quad \mathcal{K}_{0d} = \Bigl< R\partial_R \phi_0(R),\phi_d(R) \Bigr>_{L^2_R}, \\
    \mathcal{K}_{cd} &= \Bigl< \int_0^{+\infty}f(\xi)R\partial_R\phi(R,\xi)\rho(\xi)d\xi ,\phi_d(R) \Bigr>_{L^2_R}, \\
     \mathcal{K}_{dc} &= \Bigl< R\partial_R \phi_d(R),\phi(R,\eta) \Bigr>_{L^2_R}, \quad \mathcal{K}_{0c} = \Bigl< R\partial_R \phi_0(R),\phi(R,\eta) \Bigr>_{L^2_R}, \\
    \mathcal{K}_{cc} &= \Bigl< \int_0^{+\infty}\left[ f(\xi)R\partial_R\phi(R,\xi) + 2 \xi\partial_\xi f(\xi) \phi(R,\xi)\right] \rho(\xi)d\xi ,\phi(R,\eta) \Bigr>_{L^2_R},
\end{alignat*}
and similarly for $\mathcal{K}_{d0}, \mathcal{K}_{00}, \mathcal{K}_{c0}$. We remark once again that some of these inner products only make sense as a limit of $L^2(\mathbb R, d\rho)$-functions (see Theorem \ref{spectral theorem} and Notation \ref{not: notation of inner product for L^2 limit}). We start by making $\mathcal{K}_{d}$ and $\mathcal{K}_0$ more explicit.

\begin{proposition}
    One has
    $$\mathcal{K}_{dd} = \mathcal{K}_{00} = -\frac{1}{2}, \quad \mathcal{K}_{0d} = \mathcal{K}_{d0} \in \mathbb R,$$
    as well as,
    $$
\mathcal{K}_{0c}(\eta) = \frac{F(0,\eta)}{\eta}, \quad \mathcal{K}_{dc}(\eta) = ||\phi(R,\xi_d)||_{L^2_R} \cdot \frac{F(\xi_d,\eta)}{\eta-\xi_d}, \quad \eta \geq 0.
    $$
    where $F$ is as in Corollary \ref{F(xi,eta) estimates}. Moreover, $\mathcal{K}_{0c}$, $\eta\partial_{\eta}\mathcal{K}_{0c}$, $\mathcal{K}_{dc}$ and $\eta\partial_{\eta}\mathcal{K}_{dc}$ are continuous on $[0,+\infty)$ and have an arbitrary fast polynomial decay at infinity.
\end{proposition}

    \begin{proof}

Jost solution theory shows that $\phi_d(R)$ decays exponentially as $R \to +\infty$ and $\phi_0(R)$ decays as $R^{-1}$ in dimension $d = 5$ by (\ref{phi theta asympotics}) . Hence,
$$\mathcal{K}_{dd} = \mathcal{K}_{00} = -\frac{1}{2}, \quad \mathcal{K}_{0d} = \mathcal{K}_{d0} \in \mathbb R.$$
 An integration by parts shows that
    $$\mathcal{K}_{0c}(\eta) = \frac{F(0,\eta)}{\eta}, \quad \eta \geq 0,$$
    since
     $$
     W(R) \phi_0(R) = \left([\mathcal{L},R\partial_R]-2\mathcal{L} \right)\phi_0(R) = \mathcal{L}R\partial_R \phi_0(R).
     $$
For $\eta \geq 0$ fixed, the integration by parts is justified because $\phi(R,\eta)$, $\partial_R\phi(R,\eta)$ are bounded and $R\partial_R \phi_0(R)$, $\partial_{R}(R\partial_R \phi_0(R))$ vanish at zero and infinity.  Similarly, $\mathcal{K}_{dc}(\eta)$ is given by
     $$\mathcal{K}_{dc}(\eta) = \frac{\innerp[\big]{ W(R) \phi_d(R) }{ \phi(R,\eta) }_{L^2_R}}{\eta-\xi_d} = ||\phi(R,\xi_d)||_{L^2_R} \cdot \frac{F(\xi_d,\eta)}{\eta-\xi_d}, \quad \eta \geq 0,$$
    because for $\eta \geq 0$,
    $$(\eta-\xi_d)\mathcal{K}_{dc}(\eta) = \Bigl< [\mathcal{L},R\partial_R] \phi_d(R), \phi(R,\eta) \Bigr>_{L^2_R}, \quad \Bigl< \phi_d(R), \phi(R,\eta) \Bigr>_{L^2_R} \propto \delta_{\xi_d}(\eta) = 0.$$
    Moreover, since formula (\ref{iteration of integration by parts F(xi,eta)}) also holds for $\xi = 0$ and $\xi = \xi_d$,  it follows that $\mathcal{K}_{0c}$, $\eta\partial_{\eta}\mathcal{K}_{0c}$, $\mathcal{K}_{dc}$ and $\eta\partial_{\eta}\mathcal{K}_{dc}$ are continuous on $[0,+\infty)$ and have an arbitrary fast polynomial decay at infinity.
\end{proof}
    
    It remains to study the components of $\mathcal{K}_c$, i.e.  $\mathcal{K}_{c0}$, $\mathcal{K}_{cd}$ and $\mathcal{K}_{cc}$. 
    \begingroup
\allowdisplaybreaks
    \begin{proposition}
        If $f \in C_c^{\infty}((0,+\infty))$, then 
        \begin{align*}
            \mathcal{K}_{c0}f &= - \int_0^{+\infty}f(\xi) \rho(\xi) K_{0c}(\xi)d\xi, \\
            \mathcal{K}_{cd}f &= - \int_0^{+\infty}f(\xi) \rho(\xi) K_{dc}(\xi)d\xi.
        \end{align*}
    \end{proposition}
    \begin{proof}
    One has
    $$\mathcal{K}_{c0} = \lim \limits_{r \to +\infty}\int_0^r \int_0^{+\infty} f(\xi)R\partial_R\phi(R,\xi)\rho(\xi) \phi_0(R)dR d\xi.$$
    Since $f(\xi)$ has a compact support, we can interchange the order of integration and obtain
    $$\mathcal{K}_{c0} = \lim \limits_{r \to +\infty}  \int_0^{+\infty}f(\xi) \rho(\xi) \left( \int_0^r R\partial_R\phi(R,\xi) \phi_0(R) dR \right) d\xi. $$
    Integrating by parts in the $R$-integral, we get
    \begin{align*}
        \int_0^r R\partial_R\phi(R,\xi) \phi_0(R) dR &= r\phi(r,\xi)\phi_0(r) - \int_0^r \phi(R,\xi) R\partial_R\phi_0(R) dR \\
        &\phantom{=}- \int_0^r \phi(R,\xi)\phi_0(R) dR.
    \end{align*}
    The second and third term converges in $L^2(\mathbb R, d\rho)$. Hence,
    \begin{align*}
        \mathcal{K}_{c0}f &= \lim \limits_{r \to +\infty} r\phi_0(r) \int_0^{+\infty}f(\xi) \phi(r,\xi)\rho(\xi)d\xi - \int_0^{+\infty}f(\xi) \rho(\xi) K_{0c}(\xi)d\xi \\
        &\phantom{=}- \int_0^{+\infty}f(\xi) \rho(\xi) \delta_0(\xi) d\xi \\
        &= \left( \lim \limits_{r \to +\infty} r\phi_0(r) \right)  \left(\lim \limits_{r \to +\infty}  \mathcal{F}^{-1}(f)(r) \right) - \int_0^{+\infty}f(\xi) \rho(\xi) K_{0c}(\xi)d\xi \\
        &= - \int_0^{+\infty}f(\xi) \rho(\xi) K_{0c}(\xi)d\xi,
    \end{align*}
    because $\mathcal{F}^{-1}f(R)$ decays like a Schwartz function. Similarly, one computes
    $$\mathcal{K}_{cd}f = - \int_0^{+\infty}f(\xi) \rho(\xi) K_{dc}(\xi)d\xi.$$
    \end{proof}
        
    Finally, for $\mathcal{K}_{cc}$, we integrate by parts with respect to $\xi$ in the component
\begin{align*}
    \int_0^{+\infty} 2 \xi\partial_\xi f(\xi) \phi(R,\xi) \rho(\xi)d\xi &= - \int_0^{+\infty} 2 f(\xi) \phi(R,\xi) \xi\partial_\xi \rho(\xi)d\xi \\
    &\phantom{=}- \int_0^{+\infty} 2 f(\xi) \phi(R,\xi) \rho(\xi)d\xi \\
    &\phantom{=} - \int_0^{+\infty} 2 f(\xi)  \xi\partial_\xi\phi(R,\xi) \rho(\xi)d\xi \\
    &= -\int_0^{+\infty} 2 f(\xi) \left(1 + \frac{\xi\rho'(\xi)}{\rho(\xi)} \right) \phi(R,\xi)  \rho(\xi) d\xi \\
    &\phantom{=} - \int_0^{+\infty} 2 f(\xi)  \xi\partial_\xi\phi(R,\xi) \rho(\xi)d\xi,
\end{align*}
    yielding
    $$\mathcal{K}_{cc}f(\eta) = \Bigl< \int_0^{+\infty}f(\xi) \left[ R\partial_R - 2 \xi\partial_\xi \right]\phi(R,\xi) \rho(\xi)d\xi ,\phi(R,\eta) \Bigr>_{L^2_R} - 2 \left(1 + \frac{\eta\rho'(\eta)}{\rho(\eta)} \right)f(\eta)$$
    when $\eta \geq 0$. Again, the function
    $$u(R) = \int_0^{+\infty}f(\xi) \left[ R\partial_R - 2 \xi\partial_\xi \right]\phi(R,\xi) \rho(\xi)d\xi $$
    decays like a Schwartz function when $f \in C^{\infty}_c((0,+\infty))$ and the bounds one can get on $\sup_{R\geq 1}|R^nu^{(m)}(R)|$ will depend only on $n$, $m$, $||f||_{\infty}$ and $\text{supp}(f)$. Since $\phi(R,\eta)$ is uniformly bounded, it follows that $\mathcal{K}_{cc}$ is continuous if $C^{\infty}_c$ and $C^{\infty}$ are respectively endowed with the test function topology and the $L^{\infty}$ topology. Hence, the Schwartz kernel theorem (\cite[Chapter V]{Hormander}) shows that one can write
    $$\mathcal{K}_{cc}f(\eta) = \int_0^{+\infty}K(\eta,\xi)f(\xi)d\xi$$
    for some distribution-valued kernel $\eta \mapsto K(\eta,\xi)$ which is made more explicit in the following theorem.
 \endgroup
    \begin{theorem}
        The operator $\mathcal{K}_{cc}$ admits the following representation
        $$\mathcal{K}_{cc} = -\left( \frac{3}{2} + \frac{\eta \rho'(\eta)}{\rho(\eta)} \right) \delta(\xi-\eta) + \mathcal{K}_{cc}^0$$
        where $\mathcal{K}_{cc}^0$ has kernel 
        $$K_{cc}^0(\eta,\xi) = \frac{\rho(\xi)}{\eta-\xi}F(\xi,\eta)$$
        and $F(\xi,\eta)$ is as in Corollary \ref{F(xi,eta) estimates}.
        \index{F(xi,eta)@$F(\xi,\eta)$, a symmetric function which contributes to $\mathcal{K}$}
    \end{theorem}

    \begin{proof}
        This follows from another integration by parts and a change of integration order. This is proved exactly as in \cite[Theorem 5.1]{Krieger_2007}.
    \end{proof}

    \begin{definition} \label{def:l^2,alpha_rho}
        Let $L^{2,\alpha}_{\rho}$, $\alpha \in \mathbb R$, be the set of d$\rho$-measurable functions $f(\xi)$ for which the following norm is finite:
        $$||f||_{L^{2,\alpha}_{\rho}}^2 = |f(\xi_d)|^2 + |f(0)|^2 \cdot \delta_{d = 5} + \int_0^{+\infty}|f(\xi)|^2 \langle \xi \rangle^{2\alpha}\rho(\xi)d\xi.$$
        This space can be represented as $\mathbb R^{d-3} \times L^2((0,+\infty), \langle \xi \rangle^{2\alpha}\rho(\xi)d\xi)$.
        \index{L 2 alpha rho@$L^{2,\alpha}_{\rho}$, a normed space on Fourier side}
    \end{definition}

    \begin{theorem} \label{thm:boundedness of K operator}
        For any $\alpha \in \mathbb R$, the operators $\mathcal{K}_{cc}^0$, $\mathcal{K}$, $[\mathcal{K},\xi\partial_{\xi}]$ maps
        $$\mathcal{K}_{cc}^0: L^{2,\alpha}_{\rho} \rightarrow L^{2,\alpha + 1/2}_{\rho}, \quad \mathcal{K}: L^{2,\alpha}_{\rho} \rightarrow L^{2,\alpha}_{\rho}, \quad [\mathcal{K},\xi\partial_{\xi}]: L^{2,\alpha}_{\rho} \rightarrow L^{2,\alpha}_{\rho}$$
        continuously, where $\xi \partial_{\xi}$ acts as zero on the discrete component.
    \end{theorem}

    \begin{proof}
    Recall that if $f \in \mathbb R^{2} \times C^{\infty}_c((0,+\infty))$ (in dimension 5, but the same can be said in dimension 4), then
    $$\mathcal{K}f(\eta) = \begin{pmatrix} k_{1,1}f(\xi_d) + k_{1,2}f(0) -\langle f \cdot \rho, K_{dc} \rangle_{L^2((0,+\infty))} \\
k_{2,1}f(\xi_d) + k_{2,2}f(0) -\langle f \cdot \rho, K_{0c} \rangle_{L^2((0,+\infty))}  \\
\mathcal{K}_{dc}(\eta)f(\xi_d) + \mathcal{K}_{0c}(\eta)f(0) + \mathcal{K}_{cc}f (\eta)
\end{pmatrix}.
$$
Similarly, we have that
\begin{align*}
    [\mathcal{K},\xi\partial_{\xi}]f(\eta) &= \begin{pmatrix} -\langle \xi\partial_{\xi} f, \rho \cdot K_{dc} \rangle_{L^2((0,+\infty))} \\
-\langle \xi\partial_{\xi}  f, \rho \cdot K_{0c} \rangle_{L^2((0,+\infty))}  \\
 \mathcal{K}_{cc}(\xi\partial_{\xi}f) (\eta) - \eta \partial_{\eta} \mathcal{K}_{cc}f (\eta) - \eta \partial_{\eta} \mathcal{K}_{dc}(\eta)f(\xi_d) - \eta \partial_{\eta} \mathcal{K}_{0c}(\eta)f(0)
\end{pmatrix} \\
&= \begin{pmatrix} \langle f , (\rho + \xi\partial_{\xi}\rho) \cdot K_{dc} + \rho \cdot \xi\partial_{\xi} K_{dc} \rangle_{L^2((0,+\infty))} \\
\langle f, (\rho + \xi\partial_{\xi} \rho) \cdot K_{0c} + \rho \cdot \xi\partial_{\xi} K_{0c} \rangle_{L^2((0,+\infty))}  \\
 \mathcal{K}_{cc}(\xi\partial_{\xi}f) (\eta) - \eta \partial_{\eta} \mathcal{K}_{cc}f (\eta) - \eta \partial_{\eta} \mathcal{K}_{dc}(\eta)f(\xi_d) - \eta \partial_{\eta} \mathcal{K}_{0c}(\eta)f(0)
\end{pmatrix},
\end{align*}

where $k_{i,j} \in \mathbb R$ and $\mathcal{K}_{dc}$, $\xi\partial_{\xi} \mathcal{K}_{dc}$, $\mathcal{K}_{0c}$, $\xi\partial_{\xi} \mathcal{K}_{0c}$ are continuous on $[0,+\infty)$ and fast-decaying functions. Hence, it suffices to study the mapping properties of $$\mathcal{K}_{cc} = -\left( \frac{3}{2} + \frac{\eta \rho'(\eta)}{\rho(\eta)} \right) \delta(\xi-\eta) + \mathcal{K}_{cc}^0, \quad [\mathcal{K}_{cc}, \xi\partial_{\xi}] = \eta\partial_{\eta}\left( \frac{\eta \rho'(\eta)}{\rho(\eta)} \right) \delta(\xi-\eta) +  [\mathcal{K}_{cc}^0, \xi\partial_{\xi}]$$
for which the dirac-delta contribution causes no issue.

\item[\textbf{Boundedness of $\mathcal{K}_{cc}^0$}:]

The boundedness of $\mathcal{K}_{cc}^0$ is equivalent to proving that the kernel 
$$\tilde{K}^0_{cc}(\eta,\xi) = \rho(\eta)^{\frac{1}{2}} \langle \rho \rangle^{\alpha + 1/2}K_{cc}^0(\eta,\xi)  \langle \eta \rangle^{-\alpha} \rho(\xi)^{-\frac{1}{2}}: L^2((0,+\infty)) \rightarrow L^2((0,+\infty))$$
acts, as a principal value integral, continuously. Write
\begin{align*}
\tilde{K}^0_{cc}(\eta,\xi) &= \frac{\sqrt{\rho(\eta)\rho(\xi)} \langle \eta \rangle^{\alpha + 1/2} \langle \xi \rangle^{-\alpha} F(\xi,\eta)}{\eta-\xi} = \frac{\tilde{F}(\xi,\eta)}{\eta-\xi}.
\end{align*}
First, we split the kernel into two regions: the diagonal
$$D =\{(\eta,\xi) \in (0,+\infty)^2: \frac{1}{4}\xi \leq \eta \leq 4\xi\},$$
where $\xi \sim \eta$, and its complementary $(0,+\infty)^2 \setminus D$, where one always has $|\xi-\eta| \geq \frac{5}{3}(\xi + \eta)$. Our estimates on $F(\xi,\eta)$ from Corollary \ref{F(xi,eta) estimates} show that
$$\left| \frac{\tilde{F}(\xi,\eta)}{\xi-\eta}  \right| \lesssim \begin{cases}
    \eta^{\frac{d-6}{4}} \xi^{\frac{d-6}{4}} (\cdot |\log(\xi)\log(\eta) |^{-1} \text{ if } d = 4) &\text{ if } (\xi,\eta) \notin D, \xi + \eta \lesssim 1 \\
    (1+\xi)^{-N}(1+\eta)^{-N} &\text{ if } (\xi,\eta) \notin D, \xi + \eta \gtrsim 1,
\end{cases}$$
meaning that $1_{(0,+\infty)^2 \setminus D} \cdot \tilde{K}^{0}_{cc}$ is a Hilbert-Schmidt kernel on $(0,+\infty)^2$. It remains to treat the diagonal part. In that region, one has the following estimates
    \begin{align*}
        |\tilde{F}(\xi,\eta)| &\lesssim \begin{cases}
           \eta^{\frac{d-4}{2}} (\cdot |\log(\eta)|^{-2} \text{ if } d = 4) \quad &\text{ if } \xi + \eta \leq 1 \\
            (1+|\xi^{\frac{1}{2}}-\eta^{\frac{1}{2}}|)^{-N} \quad &\text{ if } \xi + \eta \geq 1 
        \end{cases} \\
        |\partial_{\xi}\tilde{F}(\xi,\eta)| + |\partial_{\eta}\tilde{F}(\xi,\eta)|&\lesssim \begin{cases}
            \eta^{\frac{d-6}{2}} (\cdot |\log(\eta)|^{-2} \text{ if } d = 4) \quad &\text{ if } \xi + \eta \leq 1 \\
             \eta^{-\frac{1}{2}}(1+|\xi^{\frac{1}{2}}-\eta^{\frac{1}{2}}|)^{-N} \quad &\text{ if } \xi + \eta \geq 1 .
        \end{cases} 
    \end{align*} 
We write
\begin{align*}
1_{D} \cdot \tilde{K}^0_{cc}(\eta,\xi) &= 1_{D} \cdot \left( \frac{\tilde{F}(\xi,\xi)}{\eta-\xi} + \frac{\tilde{F}(\xi,\eta)-\tilde{F}(\xi,\xi)}{\eta-\xi} \right) = (A) + (B)
\end{align*}
and the $L^2$-boundedness of $(A)$ follows from the boundedness of the Hilbert transform. 

If we further split the diagonal $D$ into $D \cap [0,1]^2$ and $D \setminus [0,1]^2$, then $(B) \cdot 1_{D \cap [0,1]^2}$ is Hilbert-Schmidt on $(0,+\infty)^2$. As for $(B) \cdot 1_{D \setminus [0,1]^2}$, we use Schur test: it suffices to prove that
$$\sup_{\xi \geq 0} \int_{\eta \geq 0} 1_{D \setminus [0,1]^2} \cdot \left| \frac{\tilde{F}(\xi,\eta)-\tilde{F}(\xi,\xi)}{\eta-\xi} \right| d\eta + \sup_{\eta \geq 0} \int_{\xi \geq 0} 1_{D \setminus [0,1]^2} \cdot \left| \frac{\tilde{F}(\xi,\eta)-\tilde{F}(\xi,\xi)}{\eta-\xi} \right| d\xi < +\infty$$
to get the $L^2((0,+\infty))$-boundedness. Given $\eta \sim \xi$ on $D$ and our bounds on $\nabla F$, it is enough to prove
$$\sup_{\xi \geq 1} \int_{\eta \geq 1} \eta^{-\frac{1}{2}}(1 + |\xi^{\frac{1}{2}}-\eta^{\frac{1}{2}}|)^{-N} d\eta < +\infty,$$
which is the case because
\begin{align*}
    \int_{\eta \geq 1} \eta^{-\frac{1}{2}}(1 + |\xi^{\frac{1}{2}}-\eta^{\frac{1}{2}}|)^{-N} d\eta &= \int_{1}^{\xi} \eta^{-\frac{1}{2}}(1 + \xi^{\frac{1}{2}}-\eta^{\frac{1}{2}})^{-N} d\eta + \int_{\xi}^{+\infty} \eta^{-\frac{1}{2}}(1 - \xi^{\frac{1}{2}}+\eta^{\frac{1}{2}})^{-N} d\eta \\
    &\lesssim \left| \left[ (1 + \xi^{\frac{1}{2}}-\eta^{\frac{1}{2}})^{-N+1} \right]^{\eta = \xi}_{1} \right| + \left| \left[ (1 - \xi^{\frac{1}{2}}+\eta^{\frac{1}{2}})^{-N+1} \right]^{\eta = +\infty}_{\xi} \right| \\
    &\lesssim  \left( \xi^{\frac{1}{2}} \right)^{-N+1} + 1 \\
    &\lesssim 1 \quad \forall \xi \geq 1.
\end{align*}
\begin{remark}
    In the above proof, we remark that one can also choose to write 
\begin{align*}
1_{D} \cdot \tilde{K}^0_{cc}(\eta,\xi) &= 1_{D} \cdot \left( \frac{\tilde{F}(\eta,\eta)}{\eta-\xi} + \frac{\tilde{F}(\xi,\eta)-\tilde{F}(\eta,\eta)}{\eta-\xi} \right)
\end{align*}
and use the same analysis. This is useful for the commutator estimate, because we can choose if we wish to differentiate $\tilde{F}$ with respect to $\xi$ or $\eta$. We do not need to have an estimate on the full gradient of $\tilde{F}$.
\end{remark}


\item[\textbf{Boundedness of $[\mathcal{K}_{cc}^0,\xi\partial_{\xi}]$}:]

An integration by parts shows that the commutator $[\xi\partial_{\xi},\mathcal{K}_{cc}^0]$ has kernel
\begin{align*}
    K^{0,com}_{cc} &= (\eta\partial_{\eta} + \xi\partial_{\xi})K^0_{cc}(\eta,\xi) + K_{cc}^0(\eta,\xi) \\
    &= \frac{\rho(\xi)}{\eta-\xi}\left( \frac{\xi\rho'(\xi)}{\rho(\xi)}F(\xi,\eta) + (\eta\partial_{\eta} + \xi\partial_{\xi})F(\xi,\eta) \right).
\end{align*}
Then $\eta\partial_{\eta} F$, $\xi \partial_{\xi}F$ and  $\xi\rho'(\xi)\rho(\xi)^{-1}F(\xi,\eta)$ all satisfy the same estimates as $F(\xi,\eta)$, the only difference being on the diagonal, away from zero, where we lose a factor $\eta^{-\frac{1}{2}}$, but we gain it back by considering the $L^{2,\alpha}_{\rho} - L^{2,\alpha}_{\rho}$ boundedness instead of $L^{2,\alpha}_{\rho} - L^{2,\alpha+1/2}_{\rho}$. 

We remark that for $\eta \partial_{\eta} F$ (resp. $\xi \partial_{\xi} F$), we only need the estimate for the derivative with respect to $\xi$ (resp. $\eta$) to be the same as the one for $\nabla F$.
\end{proof}

\section{Exact solutions by means of Fourier method} \label{section:exact solution fourier method}
In this section, we construct the final piece of the solution $\varepsilon$, which corrects the approximate solution $u_k$ from Theorem \ref{Thm:Approximate Solution} and Theorem \ref{Thm:Approximate Solution d = 4} to an exact solution within the cone. We first transform the evolution equation for $\varepsilon$ into the generalized Fourier space and we formulate the equation as a fixed-point problem given by (\ref{contraction operator}). The core of the section is to prove, by using the properties of the transference operator, that this map is a contraction on a carefully chosen Banach space. 

Let $\varepsilon(r,t)$ be a solution of (\ref{space-variable equation epsilon}). Substituting
$$
\tilde{\varepsilon}(\tau,R) = R^{\frac{d-1}{2}}\varepsilon(t(\tau),r(\tau,R)), \quad \tau = \nu^{-1}t^{-\nu}, \ R = \lambda(t)r, \ \lambda(t) = t^{-1-\nu},
$$
\index{tau@$\tau = \nu^{-1}t^{-\nu}$}
we get
\begin{equation} \label{new space-variable equation epsilon}
    \mathcal{D}^2 \tilde{\varepsilon} - (d-2) \beta(\tau) \mathcal{D} \tilde{\varepsilon} - \frac{(d-1)(d-3+d\nu - \nu)}{4\nu} \dot{\beta}(\tau) \tilde{\varepsilon} + \mathcal{L} \tilde{\varepsilon} = \mathcal{N}\tilde{\varepsilon},
\end{equation}
where the following notations were used:
\begin{align*}
\lambda(\tau) &= \lambda(t(\tau)) = (\nu \tau)^{\frac{1+\nu}{\nu}}, \ \dot{\lambda}(\tau) = \partial_{\tau} \lambda(t(\tau)), \\
    \beta(\tau) &= \dot{\lambda}(\tau)\lambda(\tau)^{-1} = \frac{1+\nu}{\tau \nu}, \\
    \dot{\beta}(\tau) &= \partial_{\tau}\beta = -\frac{1+\nu}{\tau^2 \nu}, \\
    \mathcal{D} &= \partial_{\tau} + \beta(\tau)R \partial_R,  \\
    \mathcal{L} &= -\partial_{RR} - pW(R)^{p-1} + \frac{1}{R^2} \cdot \left( \frac{(d-3)(d-1)}{4} \right),  \\
    \mathcal{N} &= \lambda(\tau)^{-2}R^{\frac{d-1}{2}}\left[e_{k-1} + F(u_{k-1}+\chi\left(R \tau^{-1}\right)R^{-\frac{d-1}{2}}\tilde{\varepsilon}) \right.  \\
     &\phantom{=} \left.-F(u_{k-1}) - F'(u_0)\chi\left(R \tau^{-1}\right) R^{-\frac{d-1}{2}}\tilde{\varepsilon}\right].
\end{align*}
We assume until the end of the paper that $u_0, u_{k-1}$, $e_{k-1}$ are always extended on $0 \leq R < +\infty$, $\tau \geq \tau_0$, with the same size and regularity as well as being supported in $0 \leq R < 2\tau$. We remark that we have added a term $\chi(R\tau^{-1})$ in front of $R^{-\frac{d-1}{2}}\tilde{\varepsilon}$, where $0 \leq \chi = 1-\chi_{[1,+\infty)} \leq 1$ is a smooth transition function which is $1$ on $|x| \leq 1$ and $0$ on $|x| \geq 2$, in the definition of $\mathcal{N}$. All of this does not change the equation on the cone $0 \leq R < \tau$, $\tau \geq \tau_0$, of interest.
\index{chi2@$\chi(R\tau^{-1})$, another transition function}

\begin{remark}[On extending $u_k, e_k$ outside the cone] \label{extension of v_k, e_k}
In the following, we briefly describe how one can, for fixed $t$, extend $u_k, e_k$ on the whole $\mathbb R^d$. Multiplying this extension by $\chi(R\tau^{-1})$ restricts its support to the desired region $0 \leq R < 2\tau$.

Note that $u_0, v_1, e_1$ are already defined on the whole $0 \leq R < +\infty$, $0 < t \leq t_0$. For these terms, we only need to apply the cutoff. Otherwise, write $v(R,t)$ or $e(R,t)$ as a function $f(a,t)$, $a = R/(t\lambda) \in [0,1]$, and for fixed $t$, extend $f(a,t)$ on $[0,+\infty)$ while keeping the same Hölder regularity and a comparable Hölder constant, which gives the $(t\lambda)$ smallness. This can be done via Whitney's Extension Theorem (\cite[Chapter VI, Theorem 4]{Stein_Whitney_Extension}). However, since we are working on an interval, we can use the following simpler construction:
$$
\tilde{f}(a,t) =  \int_1^{+\infty} \phi(y)\psi(1-a(1-y))f(1-a(1-y),t)dy, \quad a > 1,
$$
where $\psi \in C^{\infty}(\mathbb R)$ is a smooth transition function which is $1$ on $(-\infty,1/2]$ and $0$ on $[3/4,+\infty)$ and $\phi \in C^{0}([1,+\infty))$ is a continuous function satisfying
   $$\int_{1}^{+\infty}\phi(s)ds = 1, \quad \int_{1}^{+\infty}s^n\phi(s)ds = 0 \ \forall n \geq 1, \quad \lim \limits_{s \to +\infty}s^n\phi(s) = 0, \quad \forall n \geq 0.$$
\end{remark}

Now, we translate equation (\ref{new space-variable equation epsilon}) to the Fourier side. Observe that
$$\mathcal{F} \left( \partial_{\tau} + \beta(\tau)R \partial_R \right) = \left(\partial_{\tau}-2\beta(\tau)\xi\partial_{\xi} \right)\mathcal{F} + \beta(\tau)\mathcal{K}\mathcal{F}$$
and
\begin{align*}
  \mathcal{F} \left( \partial_{\tau} + \beta(\tau)R \partial_R \right)^2 &= \left[\left(\partial_{\tau}-2\beta(\tau)\xi\partial_{\xi}\right)\mathcal{F} + \beta(\tau)\mathcal{K}\mathcal{F}\right] \left( \partial_{\tau} + \beta(\tau)R \partial_R \right) \\
  &= \left(\partial_{\tau}-2\beta(\tau)\xi\partial_{\xi} \right)^2\mathcal{F} + \left(\partial_{\tau}-2\beta(\tau)\xi\partial_{\xi} \right) \beta(\tau)\mathcal{K}\mathcal{F} \\
  &\phantom{=}+ \beta(\tau)\mathcal{K}\left(\partial_{\tau}-2\beta(\tau)\xi\partial_{\xi} \right)\mathcal{F} + \beta(\tau)^2\mathcal{K}^2\mathcal{F},
\end{align*}
i.e.,
\begin{align*}
\mathcal{D}_{\tau} &= \partial_{\tau}-2\beta(\tau)\xi\partial_{\xi} \\
\mathcal{F} \mathcal{D} &= \mathcal{D}_{\tau}\mathcal{F} + \beta(\tau)\mathcal{K}\mathcal{F} \\
\mathcal{F} \mathcal{D}^2 &= \mathcal{D}_{\tau}^2\mathcal{F} + \beta(\tau) \mathcal{K} D_{\tau}\mathcal{F} + D_{\tau} \mathcal{K}\beta(\tau)\mathcal{F} + \beta(\tau)^2\mathcal{K}^2\mathcal{F} \\
&= \mathcal{D}_{\tau}^2\mathcal{F} + \beta(\tau) \mathcal{K} D_{\tau}\mathcal{F} + \beta(\tau) D_{\tau} \mathcal{K} \mathcal{F} + \dot{\beta}(\tau)\mathcal{K} \mathcal{F} + \beta(\tau)^2\mathcal{K}^2\mathcal{F} \\
&= \mathcal{D}_{\tau}^2\mathcal{F} + 2\beta(\tau) \mathcal{K} D_{\tau}\mathcal{F} + \beta(\tau) [D_{\tau},\mathcal{K}]  \mathcal{F} + \dot{\beta}(\tau)\mathcal{K} \mathcal{F}+ \beta(\tau)^2\mathcal{K}^2\mathcal{F} \\
&= \mathcal{D}_{\tau}^2\mathcal{F} + 2\beta(\tau) \mathcal{K} D_{\tau}\mathcal{F} - 2 \beta(\tau)^2 [\xi\partial_{\xi},\mathcal{K}]  \mathcal{F} + \dot{\beta}(\tau)\mathcal{K} \mathcal{F}+ \beta(\tau)^2\mathcal{K}^2\mathcal{F}.
\end{align*}
Therefore, (\ref{new space-variable equation epsilon}) rewrites as
\begin{equation}\label{fourier-variable equation epsilon}
\begin{aligned}
    (\mathcal{D}^2_{\tau} + \beta(\tau) \mathcal{D}_{\tau} + \xi ) \mathcal{F}\tilde{\varepsilon} &= (d-1) \beta(\tau) \mathcal{D}_{\tau} \mathcal{F} \tilde{\varepsilon} + (d-2) \beta(\tau)^2\mathcal{K}\mathcal{F} \tilde{\varepsilon}\\
    &\phantom{=}-2\beta(\tau) \mathcal{K} D_{\tau}\mathcal{F} \tilde{\varepsilon} + 2\beta(\tau)^2[\xi\partial_{\xi},\mathcal{K}] \mathcal{F} \tilde{\varepsilon} - \dot{\beta}(\tau)\mathcal{K}\mathcal{F}\tilde{\varepsilon} \\
    &\phantom{=}- \beta(\tau)^2\mathcal{K}\mathcal{F}\tilde{\varepsilon} + \frac{(d-1)(d-3+d\nu - \nu)}{4\nu} \dot{\beta}(\tau) \mathcal{F}\tilde{\varepsilon}  +  \mathcal{F}\mathcal{N}\tilde{\varepsilon}.
\end{aligned}
\end{equation}
Since
$$S^{-1}\partial_{\tau}S = \mathcal{D}_{\tau}, \quad S^{-1}\lambda(\tau)^{-2}\xi S = \xi, \quad (Sg)(\tau,\xi):= g(\tau, \lambda(\tau)^{-2}\xi),$$
the operator on the left-hand side of (\ref{fourier-variable equation epsilon}) can be inverted as was shown in \cite[Section 3]{Krieger_2012}. If
$$\underline{\mathbf{x}}(\tau,\xi) = \begin{pmatrix}
    x_d(\tau) \\
    x_0(\tau) \\
    x(\tau,\xi)
\end{pmatrix} \quad \underline{\mathbf{f}}(\tau,\xi) = \begin{pmatrix}
    f_d(\tau) \\
    f_0(\tau) \\
    f(\tau,\xi)
\end{pmatrix} \quad \underline{\mathbf{\xi}} = \xi \cdot \begin{pmatrix}
    1_{\xi = \xi_d} & 0 & 0 \\
    0 & 1_{\xi = 0} & 0 \\
    0 & 0 & 1_{\xi > 0}
\end{pmatrix}$$
in dimension 5 (and similarly in dimension 4), then the inhomogeneous problem
$$(\mathcal{D}^2_{\tau} + \beta(\tau) \mathcal{D}_{\tau} + \underline{\mathbf{\xi}} ) \underline{\mathbf{x}}(\tau,\xi) = \underline{\mathbf{f}}(\tau,\xi), \quad  \tau > 0, \xi \in \{\xi_d\} \cup [0,+\infty) \label{eq:inhomogeneous problem with D_tau}$$ 
is solved as
\begin{equation} \label{fourier-variable equation inversion}
\begin{aligned}
    x(\tau, \xi) &= \int_{\tau}^{+\infty} H(\sigma,\tau, \lambda(\tau)^2 \xi), f\left(\sigma,\frac{\lambda(\tau)^2}{\lambda(\sigma)^2} \xi\right)d\sigma, \\
    H(\sigma,\tau,\xi) &= \xi^{-\frac{1}{2}}  \sin \left[ \xi^{\frac{1}{2}} \int_{\tau}^{\sigma} \lambda(u)^{-1}du \right], \\
    \mathcal{D}_{\tau} x(\tau, \xi) &= \int_{\tau}^{+\infty} (\partial_{\tau}H)(\sigma,\tau, \lambda(\tau)^2 \xi) f\left(\sigma,\frac{\lambda(\tau)^2}{\lambda(\sigma)^2} \xi\right)d\sigma, \\
    x_0(\tau) &= \int_{\tau}^{+\infty}H_0(\tau,\sigma) f_0(\sigma)d\sigma, \quad H_0(\tau,\sigma) = \nu \sigma^{\frac{1+\nu}{\nu}} \left(\tau^{-\frac{1}{\nu}} - \sigma^{-\frac{1}{\nu}} \right), \\
    x_d(\tau) &= \int_{\tau}^{+\infty}H_d(\tau,\sigma) \tilde{f}_d(\sigma)d\sigma, \quad H_{d}(\tau,\sigma) = -\frac{1}{2}|\xi_d|^{-\frac{1}{2}} \exp \left( -\frac{1}{2}|\xi_d|^{\frac{1}{2}} |\tau-\sigma| \right), \\
    \tilde{f}_d(\tau) &= f_d(\tau)-\beta(\tau)x_d(\tau).
\end{aligned}
\end{equation}

\begin{definition}\label{definition banach space}
        For $\alpha \in \mathbb R$, $N \in \mathbb N$, $\tau_0 \geq 1$, let $L^{\infty,N} L^{2,\alpha}_{\rho}$ be the set of measurable functions $f(\tau, \xi)$ for which the following norm is finite
        $$||f||_{L^{\infty,N} L^{2,\alpha}_{\rho}} = \sup_{\tau \geq \tau_0} \tau^{N} ||f(\tau, \cdot)||_{L^{2,\alpha}_{\rho}} < +\infty$$
        and $L^{2,\alpha}_{\rho} = \mathbb R^{d-3} \times L^2((0,+\infty), \langle \xi \rangle^{2\alpha}\rho(\xi)d\xi)$ as in Definition \ref{def:l^2,alpha_rho}.
        \index{L infty N, 2 alpha rho@$L^{\infty,N} L^{2,\alpha}_{\rho}$, a mixed-normed space on Fourier side}
\end{definition}

Our last goal in this section is to prove that the fixed-point iteration $(\underline{\mathbf{x}}_{n}, \mathcal{D}_{\tau}\underline{\mathbf{x}}_{n})$ defined via 
\begin{align} \label{contraction operator}
    \underline{\mathbf{x}}_{n+1} &= (\mathcal{D}^2_{\tau} + \beta(\tau) \mathcal{D}_{\tau} + \xi )^{-1} \left[ (d-1) \beta(\tau) \mathcal{D}_{\tau} \underline{\mathbf{x}}_{n} + (d-2) \beta(\tau)^2\mathcal{K}\underline{\mathbf{x}}_{n} -2\beta(\tau) \mathcal{K} D_{\tau}\underline{\mathbf{x}}_{n} \right. \notag \\
    &\phantom{=}\left. + 2\beta(\tau)^2[\xi\partial_{\xi},\mathcal{K}]\underline{\mathbf{x}}_{n} - \dot{\beta}(\tau)\mathcal{K}\underline{\mathbf{x}}_{n} - \beta(\tau)^2\mathcal{K}\underline{\mathbf{x}}_{n} \right. \notag \\
    &\phantom{=}\left.+ \frac{(d-1)(d-3+d\nu - \nu)}{4\nu} \dot{\beta}(\tau) \mathcal{F}\tilde{\varepsilon}  +  \mathcal{F}\mathcal{N} (\mathcal{F}^{-1}\underline{\mathbf{x}}_{n}) \right] \\
    &= (\mathcal{D}^2_{\tau} + \beta(\tau) \mathcal{D}_{\tau} + \xi )^{-1} \left[ T (\underline{\mathbf{x}}_{n}, \mathcal{D}_{\tau}\underline{\mathbf{x}}_{n}) + \tilde{\mathcal{N}} \underline{\mathbf{x}}_{n} + \lambda(\tau)^{-2}\mathcal{F}R^{\frac{d-1}{2}}e_{k-1} \right], \notag \\
    \underline{\mathbf{x}}_0 &= (\mathcal{D}^2_{\tau} + \beta(\tau) \mathcal{D}_{\tau} + \xi )^{-1} \left[ \lambda(\tau)^{-2}\mathcal{F}R^{\frac{d-1}{2}}e_{k-1}\right], \notag
\end{align}
with $T$ linear and $\tilde{\mathcal{N}}$ nonlinear, is a Cauchy sequence in the Banach space
$$
\underline{\mathbf{x}} = (x_d, x_0, x) \in  L^{\infty,N-2}L^{2,\alpha+\frac{1}{2}}_{\rho}
$$
for an appropriate choice of $\alpha, N, \tau_0$. To this end, we proceed as follows. First, we have a look at the boundedness of the inverse operator
$$
\begin{pmatrix}
    (\mathcal{D}^2_{\tau} + \beta(\tau) \mathcal{D}_{\tau} + \xi )^{-1} \\
    D_{\tau} (\mathcal{D}^2_{\tau} + \beta(\tau) \mathcal{D}_{\tau} + \xi )^{-1}
\end{pmatrix} :  L^{\infty,N}L^{2,\alpha}_{\rho} \rightarrow  L^{\infty,N-2}L^{2,\alpha+\frac{1}{2}}_{\rho} \times  L^{\infty,N-1}L^{2,\alpha}_{\rho}.
$$
We prove that for arbitrarily small $\kappa > 0$, this map is bounded with norm $\leq \kappa$ if $N$ is large enough (depending on $\kappa$, $\nu$, $\alpha$) and $\tau_0$ is small enough (depending on $\kappa, \nu, \alpha, N$) in Definition \ref{definition banach space}. This is the content of Theorem \ref{thm:D_tau boundedness}. Then, we observe that the linear part $T$ of (\ref{contraction operator}) is a bounded operator inbetween
$$
T : L^{\infty,N-2}L^{2,\alpha+\frac{1}{2}}_{\rho} \rightarrow L^{\infty,N}L^{2,\alpha}_{\rho},
$$
thanks to Theorem \ref{thm:boundedness of K operator} and the gain of smallness coming from $\beta(\tau)$. Moreover, the operator norm of $T$ depends only on $\nu$, $\alpha$ and is independent of $N$ or how $\tau_0$ is chosen in Definition \ref{def:l^2,alpha_rho}. Similarly, we prove in Proposition \ref{prop:local lipschitz} that the nonlinear part $\mathcal{N}$ is locally Lipschitz (with local constant independent of $N$ and $\tau_0$) from $L^{\infty,N-2}L^{2,\alpha+\frac{1}{2}}_{\rho}$ to $L^{\infty,N}L^{2,\alpha}_{\rho}$, and that the (non-inverted) forcing term $\lambda(\tau)^{-2}\mathcal{F}R^{\frac{d-1}{2}}e_{k-1}$ belongs  $L^{\infty,N}L^{2,\alpha}_{\rho}$. This is enough to prove that the fixed-point iteration (\ref{contraction operator}) converges for an appropriate choice of $N$ and $\tau_0$.

\begin{theorem}\label{thm:D_tau boundedness}
    Let $\underline{\mathbf{f}} = (f,f_0,f_d)  \in L^{\infty,N} L^{2,\alpha}_{\rho}$. If
    $$
    (x,x_0,x_d) =  (\mathcal{D}^2_{\tau} + \beta(\tau) \mathcal{D}_{\tau} + \xi )^{-1}(f,f_0,f_d)
    $$
    is inverted as in (\ref{fourier-variable equation inversion}), then
    \begin{align*}
        ||x||_{L^{\infty,N-2} L^{2,\alpha+\frac{1}{2}}_{\rho}} + ||\mathcal{D}_{\tau} x||_{L^{\infty,N-1} L^{2,\alpha}_{\rho}} &\leq C(\nu,\alpha) \frac{1}{N} ||f||_{L^{\infty,N} L^{2,\alpha}_{\rho}}, \\
        ||x_0||_{L^{\infty,N-2}} + ||\partial_{\tau} x_0||_{L^{\infty,N-1}} &\leq C(\nu) \frac{1}{N} ||f_0||_{L^{\infty,N}}, \\ 
        ||x_d||_{L^{\infty,N}} + ||\partial_{\tau}x_d||_{L^{\infty,N}} &\leq C(\nu,N)||f_d||_{L^{\infty,N}}
    \end{align*}
    for all $N \geq N_0(\nu,\alpha)$, independently of how $\tau_0 \geq 1$ is chosen in Definition \ref{def:l^2,alpha_rho}.
\end{theorem}
\begingroup
\allowdisplaybreaks
\begin{proof}
    This is a consequence of the following bounds
    \begin{align*}
        |H(\sigma,\tau,\xi)| &\leq C(\nu) \cdot \tau \langle \xi \rangle^{-\frac{1}{2}}, \\
        |\partial_{\tau}H(\sigma,\tau,\xi)| &\leq C(\nu) \cdot 1, \\ 
        \Bigl| \Bigl| f \left( \sigma, \frac{\lambda(\tau)^2}{\lambda(\sigma)^2} \cdot \right) \langle \cdot \rangle^{\beta}  \Bigr| \Bigr|_{L^{2,\alpha}_{\rho}} &\leq \left( \frac{\sigma}{\tau} \right)^{C(\nu,\alpha,\beta)} ||f(\sigma, \cdot)||_{L^{2,\alpha+\beta}_{\rho}}, \\
        |H_0(\tau,\sigma)| &\leq C(\nu) \cdot \tau  \left( \frac{\sigma}{\tau} \right)^{2}, \\
        |\partial_{\tau}H_0(\tau,\sigma)|  &\leq C(\nu) \cdot \left( \frac{\sigma}{\tau} \right)^{2},
    \end{align*}
    for some large constants $C > 0$ and the exponential decay of $H_d$, $\partial_{\tau}H_d$. For example, one deduces
    \begin{align*}
        ||x(\tau,\cdot)||_{L^{2,\alpha+1/2}_{\rho}} &\leq \int_{\tau}^{+\infty} \tau \Bigl| \Bigl| f \left( \sigma, \frac{\lambda(\tau)^2}{\lambda(\sigma)^2} \cdot \right) \langle \cdot \rangle^{-\frac{1}{2}}  \Bigr| \Bigr|_{L^{2,\alpha+1/2}_{\rho}}d\sigma \\
        &\leq \int_{\tau}^{+\infty} \tau \left( \frac{\sigma}{\tau} \right)^{C(\nu,\alpha)} ||f(\sigma, \cdot)||_{L^{2,\alpha+1/2}_{\rho}} d\sigma \\
        &\leq \int_{\tau}^{+\infty} \tau^{1-C(\nu,\alpha)} \left( \frac{\sigma^{C(\nu,\alpha)}}{\sigma^{N}} \right) \sigma^N ||f(\sigma, \cdot)||_{L^{2,\alpha+1/2}_{\rho}} d\sigma \\
        &\leq \tau^{1-C(\nu,\alpha)} \left( \int_{\tau}^{+\infty} \sigma^{C(\nu,\alpha)-N} d\sigma \right) ||f||_{L^{\infty,N} L^{2,\alpha}_{\rho}} \\
        &\leq \frac{1}{N-C(\nu,\alpha)-1}\tau^{2-N} ||f||_{L^{\infty,N} L^{2,\alpha}_{\rho}} 
    \end{align*}
    if $N > C(\nu,\alpha) + 1$.
\end{proof}
\endgroup
As a corollary, when $\nu$ and $\alpha$ are fixed, for any arbitrarily small constant $\kappa > 0$, one can fix $N$ large enough (depending on $\kappa,\nu,\alpha$) and then $\tau_0 \geq 1$ large enough (depending on $\kappa,\nu, N$) in the definition of $L^{\infty,N} L^{2,\alpha}_{\rho}$ so that 
    \begin{align}\label{contraction constant kappa}
        ||x||_{L^{\infty,N-2} L^{2,\alpha+\frac{1}{2}}_{\rho}} + ||\mathcal{D}_{\tau} x||_{L^{\infty,N-1} L^{2,\alpha}_{\rho}} &\leq \kappa ||f||_{L^{\infty,N} L^{2,\alpha}_{\rho}}, \\
        ||x_0||_{L^{\infty,N-2}} + ||\partial_{\tau} x_0||_{L^{\infty,N-1}} &\leq \kappa ||f_0||_{L^{\infty,N}}, \notag \\ 
        ||x_d||_{L^{\infty,N-2}} + ||\partial_{\tau}x_d||_{L^{\infty,N-1}} &\leq  \kappa ||f_d||_{L^{\infty,N}}. \notag
    \end{align}

\begin{lemma} \label{lemma:isometry L, H}
    For $\alpha \geq 0$ fixed, we have the following equivalence of norms 
    $$
    ||\underline{\mathbf{x}}||_{L^{2,\alpha}_{\rho}} \asymp ||R^{-\frac{d-1}{2}}\mathcal{F}^{-1}\underline{\mathbf{x}}||_{H^{2\alpha}_{\mathrm{rad}}(\mathbb R^d)}.
    $$
\end{lemma}

\begin{proof}
    See \cite[Lemma 6.6]{Krieger_2007_waveEq}.
\end{proof}

    

\begin{theorem}[Strauss estimates] \label{Strauss Estimates Thm}
    Let $u(x) \in H_{\mathrm{rad}}^s(\mathbb R^d)$, $d \geq 2$, $1 < 2s < d$. Then $u(x)$ is continuous a.e. on $x \neq 0$. Moreover, there exists some universal constant $C(d,s) > 0$ for which
    \begin{align*}
        |u(x)| &\leq C |x|^{\frac{1}{2}-\frac{d}{2}}||u||_{H^s(\mathbb R^d)} \quad |x| \geq 1, \\
        |u(x)| &\leq C |x|^{s-\frac{d}{2}}||u||_{H^s(\mathbb R^d)} \quad |x| \leq 1.
    \end{align*}
\end{theorem}

\begin{proof}
    See \cite[Theorem 10 and Theorem 13]{strauss_estimates}.
\end{proof}

\begin{corollary} \label{Strauss Estimates Corollary}
     Let $u(x) \in H_{\mathrm{rad}}^s(\mathbb R^d)$, $d \geq 2$, $s > d/2$. Consider
     $$v = \prod_{l=1}^n  \partial_r^{i_l}u$$
     and assume that $i_l \in \mathbb N_{\geq 0}$, $s-i_l > 1/2$ and $s-i_1 - ... - i_n \geq 0$. Then $v(x) \in L^2(\mathbb R^d)$ with
     $$||v||_{L^2(\mathbb R^d)} \lesssim ||u||_{H^s(\mathbb R^d)}^n.$$
     Similarly, it holds that
     $$w = \partial_r^i u \cdot u^{n-1} \in L^2(\mathbb R^d), \quad ||w||_{L^2(\mathbb R^d)} \lesssim ||u||_{H^s(\mathbb R^d)}^n,$$
     for any $n \geq 1$ and $s \geq i \geq 0$.
\end{corollary}

\begin{proof}
We start with $w$ which holds true because $\partial_r^i u \in L^2$ and $u^{n-1} \in L^{\infty}$. As for $v$, the case $n = 1$ is trivial. Assume $n \geq 2$. Using Strauss Estimates, each term in the product is $L^{\infty} \cap L^2$ away from the origin and at the origin, the worse singularity that can happen for $\partial_r^{i_l}u$ is $|r|^{\min\{s-i_l-d/2-\delta,0\}}$ where we fix $$0 < \delta < \min\{n^{-1}(n-1)(s-d/2), s-i_l-d/2: l \in \{1,...,n\}, s-i_l-d/2 > 0\}.$$ 
Hence, the product (\ref{nonlinear sobolev estimate, square term}) is in $L^2(\mathbb R^d)$ away from the origin and has a singularity at worse 
$$|r|^{\min\{s-\sum i_l - \frac{d}{2} - n\delta + (n-1) (s-d/2),0\}} \leq |r|^{-\frac{d}{2}+}$$
at the origin which is also square-integrable.
\end{proof}

\begin{proposition}\label{prop:local lipschitz}
    Let $d \in \{4,5\}$ and $d/2 < 1 + 2\alpha < 1 + (6-d)\nu/2$. In particular, $\nu > 3$ if $d = 5$ and $\nu > 1$ if $d = 4$. Let also $N, \tau_0 \gg 1 + \nu$ and consider a pair $(u_{k-1},e_{k-1})$ with $e_{k-1}$ having smallness of order $\tau^{-2N}$, obtained from Theorem \ref{Thm:Approximate Solution} or Theorem \ref{Thm:Approximate Solution d = 4}, extended outside the cone as functions having support in $0 \leq R < 2\tau$, as well as the same regularity and smallness (see Remark \ref{extension of v_k, e_k}).

    Then the forcing term in (\ref{contraction operator}) satisfies $\lambda(\tau)^{-2}\mathcal{F}R^{\frac{d-1}{2}}e_{k-1} \in L^{\infty,N} L^{2,\alpha}_{\rho}$ and the non-linear map $\tilde{\mathcal{N}}$ given by
    \begin{align*}
        \underline{\mathbf{x}} \mapsto &\lambda(\tau)^{-2}\mathcal{F} R^{\frac{d-1}{2}}\left[F\left( u_{k-1}+\chi\left(R\tau^{-1}\right) R^{-\frac{d-1}{2}}\mathcal{F}^{-1}\underline{\mathbf{x}} \right) \right.
    \\
    &\phantom{=}\left. -F(u_{k-1})-F'(u_0) \chi 
  \left(R\tau^{-1}\right)R^{-\frac{d-1}{2}}\mathcal{F}^{-1}\underline{\mathbf{x}}\right]
    \end{align*}
    is locally Lipschitz from $L^{\infty,N-2} L^{2,\alpha+1/2}_{\rho}$ to $L^{\infty,N} L^{2,\alpha}_{\rho}$, with Lipschitz constants independent of $N$ or $\tau_0$.
\end{proposition}

\begin{proof}
The forcing term $e_{k-1}$ has regularity $e_{k-1} \in L^{\infty,2N} H^{2\alpha}_{\mathrm{rad}}$ by construction, hence \sloppy $\lambda(\tau)^{-2}\mathcal{F}R^{\frac{d-1}{2}}e_{k-1} \in L^{\infty,N} L^{2,\alpha}_{\rho}$ by Proposition \ref{lemma:isometry L, H}. Using the same Proposition \ref{lemma:isometry L, H}, it suffices to prove that
    $$y \mapsto \lambda(\tau)^{-2}\left[F(u_{k-1}+\chi \left(R\tau^{-1}\right) y)-F(u_{k-1})-F'(u_0) \chi \left(R\tau^{-1}\right) y\right]$$
    is locally Lipschitz from $L^{\infty,N-2} H^{2\alpha+1}_{\mathrm{rad}}$ to $L^{\infty,N} H^{2\alpha}_{\mathrm{rad}}$. We note that the Lipschitz constants do not depend on $N$ or $\tau_0$ because if the map is locally Lipschitz from $L^{\infty,N^*-2} L^{2,\alpha+1/2}_{\rho}$ to $L^{\infty,N^*} L^{2,\alpha}_{\rho}$ for some specific $N^* = N^*(\nu)$, then it is also locally Lipschitz with the same constant for any $N \geq N^*$. The same holds if we take a bigger $\tau_0$ in the definition of these spaces. 

    If $1+2\alpha > d/2$, then $H^{2\alpha+1}(\mathbb R^d)$ is an algebra and $H^{2\alpha+1} \hookrightarrow L^{\infty}(\mathbb R^d) \cap C^0(\mathbb R^d)$. Hence, $y$ decays as $\tau^{-N} \lesssim \tau^{-N/2}R^{-N/2}$ and is negligible compared to $u_{k-1}$, whose dominant component is $u_{0} \cdot \chi(R\tau^{-1})$ (see the extension from Remark \ref{extension of v_k, e_k}). In particular, $u_{k-1} + sy$, $s \in [0,1]$, stays non-negative on $0 \leq R < +\infty, \tau \geq \tau_0$. From now on, we ignore the cutoff and simply assume that $y$ is supported on $0 \leq R < 2\tau$ and negligible compared to $u_{k-1}$. 
    
    First, we prove that the mapping has the correct range. Write
    \begin{equation}\label{nonlinear sobolev estimate, integral form}
    \begin{aligned}
        \lambda^{-2}\left[ F(u_{k-1}+ y)-F(u_{k-1})-F'(u_0)y \right] &= \lambda^{-2}\left[ F(u_{k-1}+y)-F(u_{k-1})-F'(u_{k-1})y \right. \\
        &+ \left.F'(u_{k-1})y - F'(u_0)y \right]\\
        &= \lambda^{-2} y^2 \int_0^1\int_0^1 s_1F''(u_{k-1}+s_1s_2y)ds_2ds_1 \\
        &+ \lambda^{-2} y(u_{k-1}-u_0)\int_0^1 F''(u_0 + s(u_{k-1}-u_0)) ds.
    \end{aligned}
    \end{equation}
    
    Then we need to estimate the $L^{\infty,N}L^2$ and $L^{\infty,N}\dot{H}^{2\alpha}$ norms of these products.    Combining $H^{2\alpha+1} \hookrightarrow L^2(\mathbb R^d) \cap L^{\infty}(\mathbb R^d)$ with
    $$||u_{k-1}||_{L^{\infty}(\mathbb R^d)} \asymp \lambda^{\frac{d-2}{2}}, \quad ||u_{k-1}-u_0||_{L^{\infty}(\mathbb R^d)} \asymp \lambda^{\frac{d-2}{2}}\tau^{-2},$$
    the $L^{\infty,N}L^2$ bound will follow from the inequality $|a+b|^{p-2} \leq |a|^{p-2} + |b|^{p-2}$ and a simple application of Hölder's inequality.

    As for the $L^{\infty,N}\dot{H}^{2\alpha}$ bound, for $\tau \geq \tau_0$ fixed, we can (weakly) differentiate $\ceil{2\alpha} \in [2,1+2\alpha)$ times both products in (\ref{nonlinear sobolev estimate, integral form}) and estimate them. We need to be careful with the $\ceil{2\alpha}$-th derivative which can introduce a singularity for $\partial_R^{\ceil{2\alpha}} u_{k-1}$ at $R = \tau$ or make $\partial_R^{\ceil{2\alpha}}y$ of low regularity $L^2 \cap H^{\frac{1}{2}}$ where Strauss estimates are unavailable.
    
    First, we treat the product which is nonlinear with respect to $y$. Using Faà di Bruno formula, one has
    $$\left| \partial_R^n   F''\left( u_{k-1}+s_1s_2y \right) \right| \lesssim  \sum_{1 \cdot m_1 + ... + n \cdot m_n = n} \left| F^{(2 + m_1 + ... + m_n)} \left( u_{k-1}+s_1s_2y \right)  \right| \prod_{l=1}^n \left| \partial_R^l (u_{k-1} + s_1s_2y) \right|^{m_l}.$$
    Since
    \begin{align}
        |u_{k-1}+s_1s_2y| &\asymp |u_{k-1}|,  \\
        | \partial_R^{l} (u_{k-1}+s_1s_2y)| &\lesssim \frac{|u_{k-1}|}{1+R^l} + s_1s_2|\partial_R^l y|, \quad \text{ if } l < \ceil{2\alpha}, \\
         | \partial_R^{\ceil{2\alpha}} (u_{k-1}+s_1s_2y)| &\lesssim \frac{|u_{k-1}|}{1+R^{\ceil{2\alpha}}} \cdot \left(1 + \chi(R\tau^{-1}) \cdot \left|1-\frac{R}{\tau} \right|^{-\frac{1}{2}+}\right) +  s_1s_2|\partial_R^{\ceil{2\alpha}} y|,
    \end{align}
    (see Theorem \ref{Thm:Approximate Solution} and Theorem \ref{Thm:Approximate Solution d = 4}) in order to estimate the $L^2(\mathbb R^d)$ norm of
    $$\lambda^{-2} \partial_{R}^{\ceil{2\alpha}} \left( y^2 \int_0^1\int_0^1 s_1F''(u_{k-1}+s_1s_2y)ds_2ds_1 \right),$$
    it is enough to estimate the norm of
    \begin{equation} \label{nonlinear sobolev estimate, square term}
         \lambda^{-2} \lvert \partial_R^i y \rvert \cdot \lvert \partial_R^j y \rvert \cdot |u_{k-1}|^{p-2} \left(1 + \prod_{l=1}^n \left| \frac{\partial_R^l y}{u_{k-1}} \right|^{m_l} \right) 
    \end{equation}
    when $i + j + n = \ceil{2\alpha}$, $1 \cdot m_1 + ... + n \cdot m_n = n$, and to treat separately the term
    \begin{equation} \label{nonlinear sobolev estimate, square term, special case}
        \lambda^{-2} y^2 \cdot |u_{k-1}|^{p-2} \left(1 + \chi(R\tau^{-1}) \cdot \left|1-\frac{R}{\tau} \right|^{-\frac{1}{2}+}  \right).
    \end{equation}
    In both cases, we can ignore the $\lambda^{-2}$ and $u_{k-1}$ because
     $$||u_{k-1}||_{L^{\infty}(\mathbb R^d)} \asymp  \lambda^{\frac{d-2}{2}}, \quad ||u_{k-1}^{-1}||_{L^{\infty}(R \sim \tau)} \lesssim \frac{\tau^2}{\lambda^{\frac{d-2}{2}}},$$
    which cause no issue since we will get a $\tau^{-2N}$, $N \gg \ceil{2\alpha} \geq 2$, factor from the $y$ products. Moreover, 
    \begin{align*}
        ||y^2||_{L^2(\mathbb R^d)} +  \left| \left| y^2 \cdot \partial_R^{\ceil{2\alpha}} y \right| \right|_{L^2(\mathbb R^d)} + \left| \left|  \partial_R^i y \cdot \partial_R^j y \cdot \left( \prod_{l=1}^n \partial_R^l y \right)^{m_l} \right| \right|_{L^2(\mathbb R^d)} \\
        \lesssim ||y||_{H^{1+2\alpha}}^2 + ||y||_{H^{1+2\alpha}}^3 + ||y||_{H^{1+2\alpha}}^{\ceil{2\alpha}}
    \end{align*}
    using Corollary \ref{Strauss Estimates Corollary}, which allows treating (\ref{nonlinear sobolev estimate, square term}). It remains to estimate (\ref{nonlinear sobolev estimate, square term, special case}). In that case,
    \begin{align*}
    \left| \left|  y^2 \cdot \chi(R\tau^{-1}) \cdot \left(1-\frac{R}{\tau} \right)^{-\frac{1}{2}+}  \right| \right|_{L^2(R \nsim \tau)} &\lesssim ||y^2||_{L^2} \lesssim \tau^{-2(N-2)
} ||y||_{L^{\infty,N-2}H^{1+2\alpha}}^2,  \\
\left| \left|  y^2 \cdot \chi(R\tau^{-1}) \cdot \left(1-\frac{R}{\tau} \right)^{-\frac{1}{2}+}  \right| \right|_{L^2(R \sim \tau)} 
     &\lesssim ||y||_{L^{\infty}}^2  \left| \left| \left( 1-\frac{R}{\tau} \right)^{-\frac{1}{2}+} \right| \right|_{L^2(R \sim \tau)}  \\
        &\lesssim \tau^{-2(N-2) + d/2
} ||y||_{L^{\infty,N-2}H^{1+2\alpha}}^2.
    \end{align*}

    Now, we deal with the term in (\ref{nonlinear sobolev estimate, integral form}) which is linear in $y$. Similarly to the nonlinear case, it is enough to estimate
    \begin{equation} \label{linear sobolev estimate, square term}
         \lambda^{-2} \lvert \partial_R^i y \rvert \cdot \lvert \partial_R^j (u_{k-1}-u_0) \rvert \cdot |u_{k-1}|^{p-2}
    \end{equation}
    when $i + j \leq \ceil{2\alpha}$ and to treat separately the case
    \begin{equation} \label{linear sobolev estimate, square term, special case}
                 \lambda^{-2} \lvert y \rvert \cdot \lvert u_{k-1}-u_0 \rvert \cdot |u_{k-1}|^{p-2} \cdot \left( 1 + \chi(R\tau^{-1}) \cdot R^{-\ceil{2\alpha}} \cdot \left|1-\frac{R}{\tau} \right|^{-\frac{1}{2}+} \right).
    \end{equation}
    Using 
$$||u_{k-1}||_{L^{\infty}(\mathbb R^d)}^{p-2} \asymp  \lambda^{\frac{6-d}{2}}, \quad ||u_{k-1}-u_0||_{L^{\infty}(\mathbb R^d)} \lesssim \frac{\lambda^{\frac{d-2}{2}}}{\tau^2},$$
    and noticing that the $L^2$-contribution of the singularity on $R \sim \tau$ is cancelled by the $R^{-\frac{d+3}{2}}$ pointwise decay, which comes from the $\ceil{2\alpha} \geq 2$ derivative of $u_{k-1}-u_0, u_{k-1}$ and Strauss estimates applied to $y$, we observe a smallness gain of $\tau^{-2}$ and conclude as before.

    As for the local Lipschitz bound, one needs to estimate
    \begin{align*}
       \lambda^{-2}&\left[F(u_{k-1}+ y_1)-F(u_{k-1} + y_2)-F'(u_{k-1})(y_1-y_2) \right] 
       \\ +  \lambda^{-2}&\left[F'(u_{k-1})(y_1-y_2)-F'(u_0)(y_1-y_2) \right] \\
        = \lambda^{-2} &(y_1-y_2)^2 \int_0^1\int_0^1 s_1F''(u_{k-1} + s_1s_2(y_1-y_2)]) ds_2ds_1 \\
        +  \lambda^{-2} &(y_1-y_2)(u_{k-1}-u_0)\int_0^1F''(u_0 + s(u_{k-1}-u_0))ds
    \end{align*}
    for $y = y_1-y_2 \in L^{\infty,N-2}H^{2\alpha+1}_{\mathrm{rad}}$ supported on $0 \leq R < 2\tau$. The proof is exactly the same as before.
\end{proof}

For $\nu, \alpha$ as in Proposition \ref{prop:local lipschitz}, the strategy is to choose a threshold $N^*(\nu), \tau^*(\nu)$ for which the proposition applies when $N \geq N^*(\nu), \tau_0 \geq \tau^*(\nu)$ with some Lipschitz constant $C^*$ near $0$. Now, we can choose $\kappa$ small enough in (\ref{contraction constant kappa}) (depending on $C^*$, $\mathcal{K}$, $\nu$, $\alpha$) and then $N$, $\tau_0$ large enough so that the right-hand side operator of (\ref{contraction operator}) becomes a contraction on a small closed ball centered at 
\begin{align*}
\underline{\mathbf{x}}_0 &= (\mathcal{D}^2_{\tau} + \beta(\tau) \mathcal{D}_{\tau} + \xi )^{-1} \left[ \lambda(\tau)^{-2}\mathcal{F}R^{\frac{d-1}{2}}e_{k-1}\right] \in L^{\infty,N-2}L^{2,\alpha+\frac{1}{2}}_{\rho} ,
\end{align*}
which proves that the fixed-point iteration (\ref{contraction operator}) converges.

 \section{End of the proof} \label{section:end of proof}
 This final section concludes the proof of the main theorem by showing that the exact solution $u$, which has been rigorously constructed inside the light cone, extends as an exact solution on $\mathbb R^d$. The argument proceeds in three steps. First, we use the constructed solution $u$ to define initial data at a small time $t_0$ and invoke local well-posedness theory to guarantee the existence of an exact solution $v$ evolving backward from this data. Second, we apply the principle of finite speed of propagation to the difference $w = u - v$ to prove that $w$ must be zero inside the cone. Finally, we rely on the small-data global well-posedness theory to show that $v$ does not blow up before time zero.
 
 Let $d \in \{4,5\}$ and $\nu > (d-2)/(6-d)$. Theorems \ref{Thm:Approximate Solution} and \ref{Thm:Approximate Solution d = 4} alongside the fixed point argument from Section \ref{section:exact solution fourier method} show that there exists a radial function $u(x,t)$ on $\mathbb R^d \times [0,t_0]$, $t_0 \ll 1$, which solves a nonlinear equation
\begin{align*}
    \square u &= \square u_{k-1} + e_{k-1} + F\left[ \left( 1-\chi\left(R \tau^{-1}\right) \right) u_{k-1} + \chi \left(R \tau^{-1}\right) u \right] \\
    &\phantom{=}-F(u_{k-1}) + F'(u_0)\left( 1-\chi\left(R \tau^{-1}\right) \right) (u-u_{k-1}),
\end{align*}
 where 
 \begin{enumerate}
 \item $F(x) = |x|^{p-1}x$.
     \item $0 \leq \chi \leq 1$ is a smooth cutoff which is $1$ on $|x| \leq 1$ and $0$ on $|x| \geq 2$.
     \item $u_{k-1}, e_{k-1}$ were extended outside the cone (Remark \ref{extension of v_k, e_k}).
     \item Inside the cone $0 < |x| < t, 0 < t \leq t_0$, the relation $e_{k-1} = F(u_{k-1})-\square u_{k-1}$ holds and $\chi(R\tau^{-1}) = 1$, so that $\square u = F(u)$.
 \end{enumerate}
 This solution is of the form
$$
u(x,t) = \lambda(t)^{\frac{d-2}{2}}W(\lambda(t)x) \chi(|x|/t) + \eta(x,t), \quad \eta(x,t) = u^e(x,t) + \varepsilon(x,t), \quad \lambda(t) = t^{-1-\nu},
$$
where $u^e \in C^{\frac{1}{2}+\frac{6-d}{2}\nu-}(\mathbb R^d)$ has support in $0 \leq R < 2\tau$ and
\begin{align*}
\sup_{0 < t < t_0} t^{-\frac{6-d}{2}\nu-1}||u^e||_{H^{1+\frac{6-d}{2}\nu-}(\mathbb R^d)} +  t^{-\frac{6-d}{2}\nu} ||\partial_t u^e||_{H^{\frac{6-d}{2}\nu-}(\mathbb R^d)} &< +\infty, \\
\sup_{0 < t < t_0} t^{-N_0}||\varepsilon||_{H^{1+\frac{6-d}{2}\nu-}(\mathbb R^d)}+ t^{-N_0+1}||\partial_t \varepsilon||_{H^{\frac{6-d}{2}\nu-}(\mathbb R^d)} &< +\infty,
\end{align*}
for an arbitrarily large $N_0 \gg 1+\nu$, as well as
\begin{align*}
 \lim \limits_{t \to 0} \int_{|x| < c t} \left| \partial_{x_i} \left[ \lambda(t)^{\frac{d-2}{2}}W(\lambda(t)x) \right] \right|^2 dx &= \lim \limits_{t \to 0} \int_{|x| < c \lambda t} \left| \partial_{x_i} W(x) \right|^2 dx = ||\partial_{x_i} W||_{L^2}^2, \\
    \lim \limits_{t \to 0} \int_{|x| > c t} \left| \partial_{x_i} \left[ \lambda(t)^{\frac{d-2}{2}}W(\lambda(t)x) \right] \right|^2 dx &= \lim \limits_{t \to 0} \int_{|x| > c \lambda t} \left| \partial_{x_i} W(x) \right|^2 dx = 0, \\
    \int_{|x| < c t} \left| \partial_{t} \left[ \lambda(t)^{\frac{d-2}{2}}W(\lambda(t)x) \right] \right|^2 dx &= \mathcal{O}\left((t\lambda)^{-\frac{1}{2}} \right), \\
    \int_{|x| < c t} \left| \lambda(t)^{\frac{d-2}{2}}W(\lambda(t)x) \right|^{p+1} dx &= \mathcal{O}\left((t\lambda)^{-\frac{d-2}{2}}\log(t\lambda) \right),
\end{align*}
for any constant $c > 0$. These estimates hold true with $u_0 = \lambda(t)^{\frac{d-2}{2}}W(\lambda(t)x) \chi(|x|/t)$ instead of $\lambda(t)^{\frac{d-2}{2}}W(\lambda(t)x)$. It suffices to note that the cutoff $\chi(|x|/t)$ creates harmless additional terms
$$
\lambda(t)^{\frac{d-2}{2}}W(\lambda(t)x) \chi'(|x|/t) \frac{x_i}{t^2}, \quad \lambda(t)^{\frac{d-2}{2}}W(\lambda(t)x) \chi'(|x|/t) \frac{1}{t},
$$
when taking derivatives. On the region of support $|x| \sim t$, one has 
$$
||\lambda(t)^{\frac{d-2}{2}}W(\lambda(t)x) \chi'(|x|/t) t^{-1} ||_{L^2} \lesssim t^{\frac{1}{2}} \left| \left|\nabla_x \left[ \lambda(t)^{\frac{d-2}{2}}W(\lambda(t)x) \right] \right| \right|_{L^2} \lesssim t^{\frac{1}{2}} ||W||_{\dot{H}^1}
$$
thanks to Hölder's inequality and Hardy's inequality. 

Finally, we remark that the regularity of $(\eta, \partial_t \eta)$ is at least $H^2 \times H^1$. Our final goal is to construct a solution to (NLW) on $\mathbb R^d \times (0,t_0]$ that coincides with $u$ inside the cone $0 < |x| \leq t, 0 < t \leq t_0$. Recall now the following local well-posedness theorem for (NLW).
\begin{theorem}[Local well-posedness for (NLW)]\label{thm:local well-posedness for nlw}
    Let $S(t)((u_0, u_1))$ denote the solution operator for the linear wave equation with initial conditions $(u_0, u_1) \in \dot{H}^1 \times L^2$ at $t = 0$. Let $0 \in I$ be any interval. If
    $$
    ||(u_0, u_1)||_{\dot{H}^1 \times L^2} \leq A,
    $$
    there exists $\delta(A) > 0$ for which 
    $$
    ||S(t)(u_0, u_1)||_{L_{t,x}^{\frac{2(d+1)}{(d-2)}}(\mathbb R^d \times I)} \leq \delta
    $$
    implies the existence of a unique solution $(u, \partial_t u) \in C^0(I, \dot{H}^1 \times L^2)$ of (NLW) with initial data $(u_0, u_1)$ at $t = 0$.
\end{theorem}

\begin{proof}
   See \cite[Theorem 2.7]{Kenig_Merle_local_wellposed_low_dim}.
\end{proof}

\begin{remark}[Additional properties]\label{rmk:strichartz,global well-posedness,regularity for nlw}
    Strichartz estimates (\cite[Lemma 2.1]{Kenig_Merle_local_wellposed_low_dim}) show that
$$
||S(t)(u_0, u_1)||_{L_{t,x}^{\frac{2(d+1)}{(d-2)}}(\mathbb R^d \times I)} \leq C ||(u_0, u_1)||_{\dot{H}^1 \times L^2}
$$
for some constant independent of $I$, meaning that Theorem \ref{thm:local well-posedness for nlw} applies if we choose $I$ sufficiently small.

 Moreover, if $A$ is small enough, the conclusion of the theorem always holds and we obtain existence of a global solution (see the proof of Theorem 2.7 in \cite{Kenig_Merle_local_wellposed_low_dim}, as well as Remark 2.10 in the same paper).

 One also has persistence of regularity: if $(u_0,u_1) \in \dot{H}^1 \cap \dot{H}^{1+\mu} \times H^{\mu} = \mathcal{H}$ for some $0 \leq \mu \leq 1$, then $(u,\partial_t u) \in C^0(I,\mathcal{H})$ (\cite[Remark 2.9]{Kenig_Merle_local_wellposed_low_dim}). In the $(u_0,u_1) \in H^2 \times H^1$ case, it follows from Duhamel formula that $\partial_{tt} u \in C^0(I,L^2)$. 
\end{remark}


\begin{theorem}[Finite Speed of Propagation] \label{thm:finite propagation speed}
    Let $x_0 \in \mathbb R^d$, $t_0 \geq 0$ and
    $$
    K = \{(x,t): 0 \leq |x-x_0| \leq t_0-t, 0 \leq t \leq t_0\}.
    $$ 
    Let $u(x,t)$ be a ``strong'' solution of a nonlinear wave equation
    $$
    \square u = f(u), \quad t \in I = [0,T]
    $$
    for which 
    $$
    u(x,0) = \partial_t u(x,0) = 0, \quad x \in B(x_0,t_0), \quad x_0 \in \mathbb R^d, t_0 > 0.
    $$
    By “strong” solution, we mean that $u(x,t)$ has the following smoothness
    $$
    u \in C^0([0,T], H^2), \quad \partial_t u \in C^0( [0,T], H^1), \quad \partial_{tt} u \in C^0( [0,T], L^2),
    $$
    and $u(x,t) \in L^{\infty}(K \cap \mathbb (\mathbb R^d \times [0,t_1]))$ for any $0 \leq t_1 < \min\{t_0,T\}$. Moreover, assume that, given this regularity, $f(u)$ is a measurable function of $(x,t)$ for which
    $$
|f(u)| \leq C(||u||_{L^{\infty}(K \cap (\mathbb R^d \times [0,t_1]))})|u|
    $$
    almost everywhere on $K \cap (\mathbb R^d \times [0,T])$. Then $u = 0$ on $K \cap (\mathbb R^d \times [0,T])$.
\end{theorem}


\begin{proof}
    We follow the usual energy argument (\cite[Section 12.1.2]{evans10}). For any $0 \leq t_1 < \min\{t_0,T\}$, we prove that $u = 0$ on $\{(x,t): 0 \leq |x-x_0| \leq t_0-t, 0 \leq t \leq t_1 \}$. Let
    $$
    E(t) = \frac{1}{2} \int_{B(x_0,t_0-t)} u_t^2 + |\nabla_x u|^2 + u^2 dx, \quad 0 \leq t_1.
    $$
    Differentiation yields
    $$
    E'(t) =  \int_{B(x_0,t_0-t)} u_t u_{tt} + \nabla_x u \cdot \nabla_x u_t  + u u_t dx -  \frac{1}{2} \int_{|x-x_0| = t_0-t} u_t^2 + |\nabla_x u|^2 + u^2 dS.
    $$
    The differentiation formula 
    $$
    \frac{d}{dt} \int_{B(x_0,t_0-t)} g(t,x) dx = \int_{B(x_0,t_0-t)} g_t(t,x) dx - \int_{|x-x_0| = t_0-t} g(t,x) dS
    $$
    is justified for classical $g(t,x) \in C^1(I \times \mathbb R^d)$ functions. It also holds in the vector-valued setting $g \in C^1(I, L^2(\mathbb R^d))$ by approximating 
    $$
    g(t,x) = \sum_{n=1}^{\infty}\langle g(t,\cdot), v_k(\cdot) \rangle_{L^2} \cdot v_k(x)
    $$
    using a complete orthonormal basis of smooth functions $\{v_1,v_2,...\}$ of $L^2(\mathbb R^d)$. Uniform convergence in the $C^1([0,T],L^2(\mathbb R^d))$-norm follows from Dini's theorem: the decreasing sequence of continuous functions
    $$
    t \mapsto \left| \left| \sum_{k=n}^{\infty}\langle g(t,\cdot), v_k(\cdot) \rangle_{L^2} \cdot v_k(x) \right| \right|_{L^2(\mathbb R^d)}^2 = \sum_{k=n}^{\infty}\langle g(t,\cdot), v_k(\cdot) \rangle_{L^2}^2, \quad t \in [0,T],
    $$
    converges pointwisely to $0$ as $n \to +\infty$, hence uniformly by Dini's theorem. Uniform convergence of the sequence 
    $$
    g_t(t,x) = \sum_{n=1}^{\infty}\langle g_t(t,\cdot), v_k(\cdot) \rangle_{L^2} \cdot v_k(x)
    $$
    holds in the same fashion. An integration by parts, which is valid for Sobolev functions (\cite{evans_measure_theory}, Section 4.6), yields 
    \begin{align*}
    E'(t) =  \int_{B(x_0,t_0-t)} u_t (u_{tt} + \Delta_x u + u ) dx &-  \frac{1}{2}\int_{|x-x_0| = t_0-t} u_t^2 + |\nabla_x u|^2 + u^2 dS \\
    &+ \int_{|x-x_0| = t_0-t} (\nabla_x u \cdot \nu) u_t dS, 
    \end{align*}
    where $\nu$ is the normal outward-pointing vector of the surface $|x-x_0| = t_0-t$. Hence,
    $$
    E'(t) \leq \int_{B(x_0,t_0-t)} u_t (-F(u) + u ) dx \leq C   \int_{B(x_0,t_0-t)} |u_tu| dx \leq C E(t)
    $$ 
    since $ab \leq a^2/2 + b^2/2$ for $a,b \geq 0$. Since $E(0) = 0$, Grönwall's Lemma implies that $E(t) = 0$.
\end{proof}

Let $0 < A \ll 1$ be small enough so that the global well-posedness result from Remark \ref{rmk:strichartz,global well-posedness,regularity for nlw} holds for initial data $||(u_0,u_1)||_{\dot{H}^1 \times L^2} \leq 3A$. Let $0 < t_0 \ll 1$ be small enough (depending on $A$ and $||W||_{\dot{H^1}}$) so that for all $0 < t \leq t_0$,
\begin{align} 
\left| \int_{|x| < \frac{1}{2} t} |\nabla_{x} u(t)|^2 dx  - ||W||_{\dot{H^1}}^2 \right| &\ll A, \label{estimate on nabla_x u inside cone}  \\
\int_{|x| > t} |\nabla_{x} u(t)|^2dx + \int_{x \in \mathbb R^d} |\partial_t u(t)|^2  + |u(t)|^{p+1} dx &\ll A, \label{negligble estimate on u}  
\end{align}
where $u_0 = \lambda(t)^{\frac{d-2}{2}}W(\lambda(t)x) \chi(|x|/t)$.

Let $v$ be the local solution of (NLW) constructed at $t_0$ with initial data $(u(t_0), \partial_tu(t_0))$ (we solve (NLW) backwards in time). Assume that $v$ exists on $I = [T, t_0]$ with $0 < T < t_0$. Then the difference $w = u-v$ solves a nonlinear equation
$$
\square w = f(w)
$$
with initial conditions
    $$
    u(x,0) = \partial_t u(x,0) = 0, \quad x \in B(0,t_0),
    $$
and where $f(w) = F(u)-F(u-w) = w \int_{-1}^0 F'(u + s(u-v))ds$ on the cone $0 \leq |x| \leq t, 0 < t \leq t_0$. By Theorem \ref{thm:finite propagation speed} ($w$ has the desired regularity and the local boundedness properties hold thanks to radiality of $u$, $v$), $w = 0$ and $v = u$ on the section of the cone $0 \leq |x| \leq t, 0 < T \leq t \leq t_0$.

Next, we prove that the $\dot{H}^1 \times L^2$-norm of $v$ stays small outside the cone via conservation of energy. Let 
\begin{align*}
    E(v(t),\partial_t v(t)) = E(v(t_0), \partial_t v(t_0)) &= 
    E(u(t_0),\partial_t u(t_0)) \\
    &= \int_{\mathbb R^d} \frac{1}{2} |\nabla_{t,x}u(t_0)|^2 - \frac{1}{p+1}|u(t_0)|^{p+1}dx \\
    &= \frac{1}{2}||W||_{\dot{H}^1}^2  \pm A/8, \quad t \in I
\end{align*}
using (\ref{estimate on nabla_x u inside cone}) and (\ref{negligble estimate on u}. Inside the cone, $v = u$ so we also have 
$$
\int_{|x| < \frac{1}{2} t} \frac{1}{2} |\nabla_{t,x}v(t)|^2 + \frac{1}{p+1}|v(t)|^{p+1}dx = \frac{1}{2}||W||_{\dot{H}^1}^2 \pm A/8, \quad t \in I.
$$
Hence,
\begin{align*}
        \int_{|x| \geq \frac{1}{2} t} \frac{1}{2} |\nabla_{t,x}v(t)|^2 + \frac{1}{p+1}|v(t)|^{p+1}dx &\leq A/4 + 2 \int_{|x| \geq \frac{1}{2} t} \frac{1}{p+1}|v(t)|^{p+1}dx \\
    &\leq A/4 + C \left( \int_{|x| \geq \frac{1}{2} t} |\nabla_x v(t)|^{2}dx  \right)^{\frac{p+1}{2}}, \quad t \in I,
\end{align*}
where $C$ is the norm of the Sobolev embedding $\dot{H}^1(\mathbb R^d) \hookrightarrow L^{p+1}(\mathbb R^d)$. In other words, the continuous function
$$
t \mapsto h(t) =  \int_{|x| \geq \frac{1}{2} t} |\nabla_{t,x}v(t)|^2, \quad t \in I,
$$
satisfies
\begin{align*}
    h(t) &\leq A/2 + 2Ch(t)^{\frac{p+1}{2}}, \quad t \in I, \\
    h(t_0) &\ll A.
\end{align*}
Assume for a contradiction that there is some $t \in I$ for which $h(t) > A$. By continuity, there is a $t^* \in I$  for which $h(t^*) = A$, meaning that
$$
A \leq A/2 + 2C A^{\frac{p+1}{2}} \iff A \geq \left( \frac{1}{4C} \right)^{\frac{2}{p-1}}.
$$
The quantity on the right-hand side depends only on $d$ and $p$, therefore, we can choose $A$ to be sufficiently small at the outset to prevent this from occurring. Thus, $h(t) < A$ for all $t \in I$.

Finally, we prove that $v$ exists up to time $t = 0$ (hence $v$ extends $u$ outside the cone). Assume that $v$ exists for time $t \in [T_{-},t_0)$ with $0 < T_{-} \leq t_0$. Let $0 \leq \psi \leq 1$ be a smooth cutoff which is $1$ on $|x| \leq 1/2$ and $0$ on $|x| \geq 3/4$. Consider the solution $w$ obtained by solving (NLW) with initial data
$$
(w_0, w_1) = \left(1-\psi\left( \frac{|x|}{T_{-}} \right) \right) v(T), \partial_t v(T) \in H^2 \times H^1.
$$
For $|x| \geq 3T_{-}/4$, this coincides with $v(T_{-}), \partial_t v(T_{-})$ and we have small energy
$$
\int_{|x| \geq 3T_{-}/4} |\nabla_x w_0|^2 + |w_1|^2 dx \leq A.
$$
If $t_0$ was chosen sufficiently small in the first place (depending on $\psi$, $d$, $p$, which are independent of $v$ and its interval of existence), then
$$
\int_{\frac{1}{2}T_{-} \leq |x| \leq 3T_{-}/4} |\nabla_x w_0|^2 + |w_1|^2 dx \leq 2A,
$$
as well using Hardy's inequality to estimate the components where the derivative falls in the cutoff. The small-energy global well-posedness theory implies that $w$ is a global solution and finite speed of propagation implies that $v = w$ on some neighbourhood of $\{x \in \mathbb R^d: |x| \geq 3T_{-}/4\} \times \{T_{-}\} \subset \mathbb R^d \times [T_{-},t_0]$. Hence, $v$ can be extended as
$$
v(x,t) = \begin{cases}
    u(x,t), \quad (x,t) \in \{(x,t): 0 < |x| < t, 0 < t \leq T_{-} ]\} \\
    w(x,t), \quad (x,t) \in \{(x,t): |x| \geq t, 0 < t \leq T_{-} ]
\end{cases}
$$
which concludes the proof.

\section{Renormalization Step in dimension 4} \label{section:appendix, dimension 4}
We perform the main inductive argument of the renormalization procedure in dimension $d = 4$, explaining how to construct the even correction terms $v_{2k}$ from the error $e_{2k-1}$ by solving a wave-like equation in self-similar coordinates, and the odd correction terms $v_{2k+1}$ from the error $e_{2k}$ using an elliptic-like equation. We prove that at each step, there is a systematic decrease in the error, finishing the proof of Theorem \ref{Thm:Approximate Solution}. The formalism in dimension $d = 4$ is slightly different. In this situation, we can use the simpler framework from Krieger-Schlag-Tataru (\cite{Krieger_2007_waveEq}, the 3D case) to construct our approximate solution as the exponent $p = 3$ in the nonlinearity is an integer. Let $\nu > 0$ in this section. The restriction $\nu > 1$ in Theorem \ref{thm:blow-up, main thm} arises only in Proposition \ref{prop:local lipschitz}.

\begin{definition}
Let $\tilde{\mathcal{Q}}_{\beta}, \mathcal{Q}_{\beta}$ be defined as in Definition \ref{definition of tilde Q} but with no logarithmic singularity $\log(a)^j$ in the expansion at $a = 0$. In other words, we restrict to functions which are holomorphic at $a = 0$. We define $\mathcal{Q} = \mathcal{Q}_{\nu + 1/2}$ and $\mathcal{Q}' = \mathcal{Q}_{\nu - 1/2}$. This new family $\mathcal{Q}'$ is obtained from $\mathcal{Q}$ by applying $a \partial_a$, $a^{-1} \partial_a$ or $(1-a^2) \partial_{aa}$ and $\mathcal{Q}$ is obtained from $\mathcal{Q'}$ by applying $(1-a^2)$. Moreover, $\mathcal{Q} \subset \mathcal{Q}'$.
\index{Q@$\tilde{\mathcal{Q}}$, $\mathcal{Q}$, families of holomorphic functions with respect to the self-similar variable $a$}
\end{definition}

\begin{definition}[Space $S^m(R^k \log (R)^l, \mathcal{Q})$]
$S^m(R^k \log (R)^l)$ is the class of real-analytic functions $w(R): [0,\infty)  \rightarrow \mathbb R$ for which
\index{S-space@$S^{2n}(R^I,\log(R)^J)$, a vector space of smooth functions with respect to the radial variable $R$}
\begin{enumerate}
    \item $w$ has a zero of order $m$ at $R = 0$ and $R^{-m}w(R)$ has an even Taylor expansion at $R = 0$.
    \item $w(R)$ has the following expansion at $R = +\infty$
    $$w(R,a,b) = R^k \sum_{j=0}^l w_j(R^{-1}) \log(R)^{j},$$
    where $w_j$ has an even Taylor expansion at $R = 0$.
\end{enumerate}

$IS^m(R^k \log (R)^l, \mathcal{Q})$ will denote the space of analytic functions $u(r,t)$ on the cone $C_0 = \{(r,t): 0 \leq r < t, 0 < t < t_0\}$ given by a finite sum
$$u(r,t) = \fsum_{i} P_i((t\lambda)^{-2}) Q_i(r/t) w(r \lambda) = \fsum_{i} P_i(b) Q_i(a) w_i(r \lambda)$$
on the cone, for some polynomials $P_i(b)$, $Q_i \in \mathcal{Q}$ and $w_i \in S^m(R^k \log (R)^l)$. We have a similar definition with $\mathcal{Q}'$ instead of $\mathcal{Q}$.
\index{b@$b = (t\lambda)^{-2}$, a useful variable in dimension $d = 4$}
\index{IS-space@$IS^m(R^k \log (R)^l, \mathcal{Q})$, a vector space of smooth functions with respect to the variables $a$,$b$,$R$}
\end{definition}

\begin{proposition}
The following simple rules of calculations will be used throughout the proof:

\begin{enumerate}
    \item $(t\lambda)^{-2} = b = a^2R^{-2}$
    \item $IS^{m_1}(R^{k_1} \log (R)^{l_1}, \mathcal{Q}) IS^{m_2}(R^{k_2} \log (R)^{l_2}, \mathcal{Q}) \subset IS^{m_1+m_2}(R^{k_1+k_2} \log (R)^{l_1+l_2}, \mathcal{Q})$
    \item $P(b,a^2) IS^m(R^k \log (R)^l, \mathcal{Q}) \subset IS^m(R^k \log (R)^l, \mathcal{Q})$ for any bivariate polynomial $P(x,y)$
    \item $IS^m(R^k \log (R)^l, \mathcal{Q}) = R^i IS^{m-i}(R^{k-i} \log(R)^l, \mathcal{Q})$ for any $i \in \mathbb Z_{\leq m}$
    \item $IS^m(R^k \log (R)^l, \mathcal{Q}) = (1+R^2)^{i/2} IS^{m}(R^{k-i} \log(R)^l, \mathcal{Q})$ for any $i \in \mathbb Z$
    \item $b^i(1+R^2)^i IS^m(R^k \log (R)^l, \mathcal{Q}) = (b+a^2)^i  IS^m(R^k \log (R)^l, \mathcal{Q}) \subset IS^m(R^k \log (R)^l, \mathcal{Q})$ for any $i \in \mathbb N$
\end{enumerate}
The same rules hold with $\mathcal{Q}'$ instead of $\mathcal{Q}$. Moreover, any differential operator mapping $\mathcal{Q}$ to $\mathcal{Q}'$ (such as $a \partial_a$) maps $IS^m(R^k \log (R)^l, \mathcal{Q})$ to  $IS^m(R^k \log (R)^l, \mathcal{Q}')$. The same statement holds when exchanging the roles of $\mathcal{Q}$ and $\mathcal{Q}'$.
\end{proposition}

\begin{proposition} \label{b-decomposition}
Let $w(R,a,b) \in S^m(R^k \log (R)^l, \mathcal{Q})$. Then 
$$
w(R,a,b) - w(R,a,0) \in bS^m(R^k \log (R)^l, \mathcal{Q}).
$$
\end{proposition}

\begin{proof}
Write
$$w(R,a,b) - w(R,a,0) = b \int_0^1 \partial_b w(R,a,tb) dt.$$
\end{proof}

\begin{theorem}\label{Thm:Approximate Solution d = 4}
In dimension $d = 4$, we prove that
\begin{align}
    v_{2k-1} &\in \frac{\lambda}{(t\lambda)^{2k}} IS^2(R^0 \log(R)^{m_k}, Q) ,\label{v_2k-1}\\
    t^2e_{2k-1} &\in \frac{\lambda}{(t\lambda)^{2k}} IS^{2}(R^0 \log(R)^{p_k},Q'), \label{e_2k-1}\\
    v_{2k} &\in \frac{\lambda}{(t\lambda)^{2k}} a^2IS^0(R^0 \log(R)^{p_k}, Q) \subset \frac{\lambda}{(t\lambda)^{2k+2}} IS^2(R^2 \log(R)^{p_k}, Q), \label{v_2k}\\
    t^2e_{2k} &\in \frac{\lambda}{(t\lambda)^{2k}} ( IS^0(R^{-2} \log(R)^{q_k},Q) + b IS^{2}(R^0 \log(R)^{q_k},Q') ). \label{e_2k}
\end{align}

for some increasing sequences of non-negative integers $m_k$, $p_k$, $q_k$, where $p_0 = q_0 = 0$, $m_1 = 1$. Moreover, for the $IS(\cdot, \cdot)$ part of $v_{2k-1}$ and $v_{2k}$, one can find representatives which do not depend on $b$. Additionally, $v_1$ and $t^2e_1$ have no $\mathcal{Q}$ element in their definition and the dominant components of $v_1, t^2e_1, v_2$ have no logarithm.
\end{theorem}

\subsubsection{Initialization}

One checks that $u_0, t^2 e_0 \in \lambda IS^0(R^{-2})$. We also define
\begin{align*}
M_k(v) &= v (3u_k^2 + 3 u_k v + v^2), \\
N_{2k-1}(v) &= M_{2k-2}(v) - pu_0^{p-1}v, \quad N_{2k}(v) = M_{2k-1}(v).
\end{align*}

\subsubsection{Construction of $v_{2k-1}$ from $e_{2k-2}$}

We write $$t^2 e_{2k-2} = t^2 \hat{e}_{2k-2}^0 + t^2 \hat{e}_{2k-2}^1 \in \frac{\lambda}{(t\lambda)^{2k-2}}( IS^0(R^{-2} \log(R)^{q_{k-1}},Q) + b IS^{2}(R^0 \log(R)^{q_{k-1}},Q') )$$
and further split $t^2 \hat{e}_{2k-2}^0$ into $\lambda (t\lambda)^{-(2k-2)}(w^0 +bw^1)$, where $w^0$ does not depend on $b$, as in Proposition \ref{b-decomposition}. We then set $$t^2 e^0_{2k-1} = \frac{\lambda}{(t\lambda)^{2k-2}}w^0(R,a).$$
In radial coordinates, (\ref{odd v}) reads as 
$$t^2 \mathcal{L}_R v_{2k-1}(r,t) = t^2 e^0_{2k-2}(r,t),$$
where 
$\mathcal{L}_r = -\partial_r^2 - \frac{3}{r} \partial_r - 3u_0^2$ and $t$ is a parameter. We do the change of variables $R = r \lambda(t)$ and get
$$(t \lambda)^2 \mathcal{L} v_{2k-1}(R,t) = t^2 e^0_{2k-2}(r,t),$$
where 
$\mathcal{L} = -\partial_R^2 - \frac{3}{R} \partial_R - 3W(R)^2$. We write $t^2 e^0_{2k-2}(r,t) = \lambda (t\lambda)^{-(2k-2)}w^0(R,a)$ and look for a solution $\lambda (t\lambda)^{-2k}v(R,a)$ by treating $a, t$ as parameters. This is the same as solving
$$\mathcal{L} v(R,a) = w^0(R,a), \ a = r/t,$$
where we ignore terms of $\mathcal{L} v$ which involve $\partial_a$ or $\partial_{aa}$. The initial conditions required are $v(0,a) = v'(0,a) = 0$. Then we prove that 
$$v_{2k-1} \in \frac{\lambda}{(t\lambda)^{2k}} IS^2(R^0 \log(R)^{q_{k-1}+2}, \mathcal{Q})$$
as in dimension $5$. We note that for $v_1$, there is no logarithm in the dominant component: the equation at infinity is a regular singular ODE given by (\ref{step 1 at infinity}) with $-2z^{-1}V_1$ replaced by $-z^{-1}V_1$. Hence, applying Theorem \ref{thm:inhomogeneous fuchs ode} ($r_1 = 2, r_2 = 0, \beta = 1$) yields an analytic solution and a logarithmic component $u_1(z)\log(z)$ where $u_1(z) = o(z^2)$, meaning that $v_{1} \approx R^0 + R^{-1} + R^{-2}\log(R)$ as $R \to +\infty$.

\subsubsection{Construction of $e_{2k-1}$ from $v_{2k-1}$}

We have
$$t^2e_{2k-1} = t^2N_{2k-1}(v_{2k-1}) + t^2e^1_{2k-2} - t^2E^tv_{2k-1} - t^2E^av_{2k-1},$$
where $E^tv_{2k-1} = \partial_{tt} [\lambda (t\lambda)^{-2k} v(r \lambda ,r/t)]$ but we ignore the $a$-derivatives and 
$$
E^av_{2k-1} = \square [\lambda (t\lambda)^{-2k} v(r \lambda,r/t)]
$$
but we keep only the terms where at least one $a$-derivative appears. The proof that $t^2 N_{2k-1}(v_{2k-1})$ belongs to the right space is an algebraic computation using the fact that 
\begin{align}
\begin{split}
    u_{k}-u_0 = v_1 + \sum_{i=2}^{2k-2}v_j &\in \frac{\lambda}{(t \lambda)^2}IS^2(R^0 \log(R)^{n_k},\mathcal{Q}) + \frac{\lambda}{(t \lambda)^4}IS^2(R^2 \log(R)^{n_k},\mathcal{Q}) \\
    &\subset \frac{\lambda}{(t \lambda)^2}(IS^2(R^0 \log(R)^{n_k},\mathcal{Q}) + a^2R^{-2} IS^2(R^2 \log(R)^{n_k},\mathcal{Q})) \\
    &\subset \frac{\lambda}{(t \lambda)^2} IS^0(R^0 \log(R)^{n_k},\mathcal{Q}) \\
    &\subset \lambda b (1+R^2) IS^0(R^{-2} \log(R)^{n_k},\mathcal{Q}) \\
    &\subset \lambda IS^0(R^{-2} \log(R)^{n_k},\mathcal{Q})
\end{split}
\end{align}
and $u_0 \in \lambda IS^0(R^{-2})$, so that $u_{k} \in \lambda IS^0(R^{-2} \log(R)^{n_k},\mathcal{Q})$ as well.

For $t^2 E^tv_{2k-1}$, we observe that
$$t^2 \partial_{tt} \left( \frac{\lambda}{(t\lambda)^{2k}} S^2(R^0 \log(R)^{m_k}) \right) \subset \frac{\lambda}{(t\lambda)^{2k}} S^2(R^0 \log(R)^{m_k}).$$
Finally, for $t^2 E^av_{2k-1}$, we write $v_{2k-1}(r,t) = \lambda (t\lambda)^{-2k}v(R,a)$ and observe that
\begin{align*}
    t^2 E^av_{2k-1} &= 2 t^2 \partial_{t} \left( \frac{\lambda}{(t \lambda)^{2k}} \right) v_{a}(R,a) \frac{-r}{t^2} + t^2 \frac{\lambda}{(t \lambda)^{2k}} \left( 2v_{aR}(R,a) \partial_{t}(\lambda)  \frac{-r^2}{t^2}  \right. \\
    &+ \left. v_{a}(R,a) \frac{2r}{t^3} +  v_{aa}(R,a) \frac{r^2}{t^4} - v_{a}(R,a) \frac{3}{rt} - v_{aa}(R,a) \frac{1}{t^2} - v_{aR}(R,a) \frac{2 \lambda }{t} \right) \\
    &= -2 t \partial_{t} \left( \frac{\lambda}{(t \lambda)^{2k}} \right) a v_{a}(R,a) + \frac{\lambda}{(t \lambda)^{2k}} \left( 2(\nu + 1)aR v_{aR}(R,a) \right. \\
    &+ \left. 2av_{a}(R,a) - (1 - a^2)v_{aa}(R,a) - 3a^{-1}v_{a}(R,a) - 2a^{-1}Rv_{aR}(R,a)\right) \\
    &\in \frac{\lambda}{(t\lambda)^{2k}} IS^2(R^0 \log(R)^{q_{k-1}}, \mathcal{Q}')
\end{align*}
and we note that, since the dominant component of $v_1$ contains no logarithm, the same is true for the dominant component of $t^2e_1$.

\subsubsection{Construction of $v_{2k}$ from $e_{2k-1}$}

As in Step 1, we keep from $t^2 e_{2k-1}$ the part $t^2 e_{2k-1}^0$ whose $IS(\cdot, \cdot)$ component is independent of $b$. We consider the main asymptotic component of $t^2e_{2k-1}^0$, i.e.,
$$t^2 \tilde{e}^0_{2k-1} = \frac{\lambda}{(t\lambda)^{2k}} \sum_{j=0}^{p_k} q_j(a) \log(R)^j, \ q_j \in \mathcal{Q}',$$
and solve the equation
$$t^2(-\partial_t^2 + \partial_r^2 + \frac{3}{r}\partial_r) \tilde{v}_{2k} = -t^2 \tilde{e}^0_{2k-1}.$$
Using Theorem \ref{hypergeometric ode with error forcing term}, we find a solution $\tilde{v}_{2k}$ of the following form
$$\tilde{v}_{2k} = -\frac{\lambda}{(t\lambda)^{2k}} \sum_{j=0}^{p_k} W_{j}(a) \log(R)^j, \quad W_{j}(a) \in a^2\mathcal{Q},$$
where we note that no logarithmic singularity has been created at $a = 0$. Then, we define
$$v_{2k} =  -\frac{\lambda}{(t\lambda)^{2k}} \sum_{j=0}^{p_k} W_{2k,j}(a) \frac{1}{2^j} \log(1+R^2)^j \in \frac{\lambda}{(t\lambda)^{2k}} a^2IS^2(R^0 \log(R)^{p_k}, Q),$$
where we use that
\begin{align}
\begin{split}
    \log(1+R^2)^j &= \left(\log\left(1 + \frac{1}{R^2}\right) + 2 \log(R)\right)^j \\
    &= \left( \sum_{n=1}^{\infty} \frac{(-1)^{n+1}}{n} R^{-2n} + 2\log(R) \right)^j, \quad R > 1 \label{log decomposition}
\end{split}
\end{align}
to get the expansion at infinity. As for $v_1$ and $t^2e_1$, the dominant component of $v_2$ contains no logarithm.

\subsubsection{Construction of $e_{2k}$ from $v_{2k}$}

We define 
$$t^2 e^0_{2k-1} = \frac{\lambda}{(t\lambda)^{2k}} \sum_{j=0}^{p_k} q_j(a) \left(\frac{1}{2} \log(1+R^2)\right)^j$$
and write
\begin{align*}
    t^2e_{2k} &= t^2(e_{2k-1} - \square v_{2k} + N_{2k}(v_{2k})) \\
    &= t^2(e_{2k-1} - e^0_{2k-1}) + t^2(e^0_{2k-1} - \square v_{2k}) + t^2 N_{2k}(v_{2k}) 
\end{align*}
and we prove that each part belongs to the right space. The proof of $$t^2 N_{2k}(v_{2k}) \in \frac{\lambda}{(t\lambda)^{2k}}  IS^2(R^{-2} \log(R)^{p_k}, \mathcal{Q})$$ is just simple algebra. By construction, one has 
$$t^2(e_{2k-1} - e^0_{2k-1}) \in \frac{\lambda}{(t\lambda)^{2k}} IS^0(R^{-2} \log(R)^{p_k}, \mathcal{Q}'),$$
where everything is straightforward except the expansion at infinity which follows from formula (\ref{log decomposition}) and the fact that the main asymptotic component of $t^2(e_{2k-1} - e^0_{2k-1})$ is of the form
$$\frac{\lambda}{(t\lambda)^{2k}} \sum_{j=0}^{p_k} q_j(a) \left( \log(R)^j - \frac{1}{2^j} \log(1+R^2)^j \right).$$
This finishes the study of the term $t^2(e_{2k-1} - e^0_{2k-1})$ since one has the inclusion
$$IS^0(R^{-2} \log(R)^{p_k}, \mathcal{Q}') \subset IS^0(R^{-2} \log(R)^{p_k}, \mathcal{Q}) + bIS^2(R^0 \log(R)^{p_k}, \mathcal{Q}')$$
by writing $w(R,a,b) = (1-a^2)w(R,a,b) + bR^2w(R,a,b)$.

It remains to prove that 
$$f = t^2(e^0_{2k-1} - \square v_{2k}) \in \frac{\lambda}{(t\lambda)^{2k}} IS^0(R^{-2} \log(R)^{p_k}, \mathcal{Q}').$$
Again, the only part which is not straightforward is the behaviour at infinity. For this, we write
\begin{align*}
    f &= t^2(\tilde{e}^0_{2k-1} - \square \tilde{v}_{2k}) + t^2(e^0_{2k-1} - \tilde{e}^0_{2k-1}) - t^2 \square(v_{2k} - \tilde{v}_{2k}) \\
    &= 0 + t^2(e^0_{2k-1} - \tilde{e}^0_{2k-1}) + t^2\square(v_{2k} - \tilde{v}_{2k}).
\end{align*}
Note that $t^2(e^0_{2k-1} - \tilde{e}^0_{2k-1})$ is the main asymptotic component of $-t^2(e_{2k-1} - e^0_{2k-1})$ which was already studied earlier. So it remains to deal with $t^2\square(v_{2k} - \tilde{v}_{2k})$, where $$(v_{2k} - \tilde{v}_{2k}) = \frac{\lambda}{(t\lambda)^{2k}}\sum_{j=1}^{p_k}a^2W_j(a) [\log(R)^j - \log(1+R^2)^j], \quad W_j \in \mathcal{Q}.$$
One computes explicitly
\begin{align*}
t^2 \square \left( \frac{\lambda}{(t\lambda)^{2k}}w(R,a) \right) = \frac{\lambda}{(t\lambda)^{2k}} &\left[ 1 + a^{-1}\partial_a + a \partial_a + (a^2-1)\partial_{aa} + a^{-2}R\partial_R + R\partial_R  \right. \\
&\left. \phantom{[1} + a^{-1} \partial_a R \partial_R + a \partial_a R \partial_R + R^2 \partial_{RR} + a^{-2} R^2 \partial_{RR} \right] w(R,a)
\end{align*}
up to some omitted multiplicative constants depending on $k,\nu$. Hence, for terms of the form
$$w(R,a) = \frac{\lambda}{(t\lambda)^{2k}}a^2W_j(a) [\log(R)^j - \log(1+R^2)^j], \quad W_j \in \mathcal{Q},$$
we get the desired results. 

\section{Some results about regular singular ODEs} \label{section:appendix, ode}

In this last section, we consider a linear ordinary differential equation
\begin{equation} \label{regular singular fuchs ode}
    u''(z) + p(z)u'(z) + q(z)u(z) = 0, \quad z \in \mathbb C,
\end{equation}
around a regular singular point $0 \in \mathbb C$, meaning that
$$p(z) = \frac{1}{z}\sum_{n=0}^{+\infty}p_n z^n, \quad q(z) = \frac{1}{z^2}\sum_{n=0}^{+\infty}q_n z^n, \quad |z| < R.$$
The standard method for finding solutions for this equation (with zero or analytical forcing term) is to make a power series Ansatz. This is called the Frobenius method. The goal of this section is to generalize the Frobenius method to solve the equation with power and logarithmic forcing terms. These results are used extensively in the previous sections.

We recall (see \cite[Chapter 4]{Teschl_ode}) that if $\{r_1, r_2\}$, $\Re(r_1) \geq \Re(r_2)$, are the roots of the indicial equation $\alpha^2 + (p_0-1)\alpha+q_0 = 0$, then one can find a fundamental system of the form
\begin{equation} \label{fuchs fundamental system}
    u_1(z) = z^{r_1} h_1(z), \quad u_2(z) = \underbrace{z^{r_2}h_2(z)}_{=: \tilde{u}_2(z)} + c \cdot \log(z) u_1(z),
\end{equation}
where $h_i(z)$ is analytic at $0$ with $h_i(0) = 1$ and radius of convergence at least equal to the distance between $0$ and the next singularity of $p(z)$ and $q(z)$. Moreover, if $r_1 - r_2 \notin \mathbb N_{\geq 0}$, then the constant $c \in \mathbb C$ is necessarily $0$ and if $r_1 = r_2$, then $c \in \mathbb C$ is necessarily non-zero.

Finally, one observes that $r_1 + r_2 = 1 - p_0$, $r_1r_2 = q_0$ and the Wronskian is of the form
\begin{equation} \label{fuchs wronskian}
    W(z) = Cz^{-p_0}\exp\left(-\sum_{n=1}^{+\infty}\frac{p_n}{n} z^{n} \right) = z^{-p_0}h_3(z),
\end{equation}
where $h_3(z)$ is analytic and non-zero on $|z| < R$.

In the following, we are interested in solving inhomogeneous regular singular problems where the forcing term can be a combination of powers and logarithms.



\begin{proposition}[Parseval Identity]
    Let $f(z)$ be holomorphic on $B(0,R)$ with
    $$f(z) = \sum_{n=0}^{+\infty}f_nz^n, \quad |z| < R.$$
    Then for any $0 < r < R$,
    $$\sum_{n=0}^{\infty}|f_n|^2r^{2n} = \frac{1}{2\pi}\int_0^{2\pi} |f(re^{i\theta})|^2 d\theta.$$
\end{proposition}

\begin{proof}
    Write 
    $$ \frac{1}{2\pi}\int_0^{2\pi} |f(re^{i\theta})|^2 d\theta =  \frac{1}{2\pi}\int_0^{2\pi} f(re^{i\theta}) \overline{f(re^{i\theta})} d\theta$$
    and expand both $f(re^{i\theta})$, $\overline{f(re^{i\theta})}$ around zero.
\end{proof}

\begin{definition}[Wiener Space]\label{wiener space}
    The Wiener Space $A(|z| < R)$ is the normed vector space of holomorphic functions 
    $$f(z) = \sum_{n=0}^{\infty} f_nz^n, \quad  |z| < R,$$
    with coefficients $(f_n R^n)_{n = 0}^{\infty} \in \ell^1(\mathbb N)$. The norm on this space is defined as
    $$||f||_{A(|z| < R)}:= ||(f_n R^n)||_{\ell^1(\mathbb N)}$$
    and we observe that
    $$||\sum_{n=m_1}^{m_2}f_nz^n ||_{L^{\infty}(|z| < R)} \leq ||f||_{A(|z| < R)} \quad \forall 0 \leq m_1 \leq m_2 \leq +\infty,$$
    as well as
    $$
||f^{(k)}||_{A(|z| < \delta R)} = \left( \delta R \right)^{-k}  \sum_{n=0}^{+\infty} \frac{n!}{(n-k)!} \delta^{n} |f_n| R^n \lesssim_{\delta,R} ||f||_{A(|z| < R)} \quad \forall k \in \mathbb N_{\geq 0}
    $$
    whenever $0 < \delta < 1$.
    
    This space is an algebra and for $g(z) \in A(|z| < R)$ fixed, the multiplication operator $T_g: f \mapsto f \cdot g$ is a bounded operator with norm $||T_f|| = ||g||_{A(|z| < R)}$.
    \index{A(z < R)@$A(\lvert z \rvert < R)$, the Wiener space}
\end{definition}




\begingroup
\allowdisplaybreaks
\begin{theorem}[Nonhomogeneous ODE with singular forcing term] \label{thm:inhomogeneous fuchs ode}
Consider a regular singular ODE (\ref{regular singular fuchs ode}) around $0 \in \mathbb C$ with indicial roots $\{r_1,r_2\}$, $\Re(r_1) \geq \Re(r_2)$, and assume that $p(z), q(z)$ have radius convergence $R + \varepsilon$. Let $u_1(z) = z^{r_1}h_1(z), u_2(z) = z^{r_2}h_2(z) + c \cdot u_1(z)\log(z)$ be the fundamental system from (\ref{fuchs fundamental system}). Let $\beta \in \mathbb C$, $j \in \mathbb N_{\geq 0}$ and 
$$g(z) = \sum_{n=0}^{\infty}g_nz^n, \quad |z| < R+\varepsilon.$$ 
In the following, we define $g_{r_i-\beta}$ to be the $(r_i-\beta)$-th order term of $g(z)$ series expansion at $z = 0$ if $r_i - \beta \in \mathbb N_{\geq 0}$ and $g_{r_i-\beta} = 0$ otherwise.

The inhomogeneous equation
\begin{equation} \label{eq:inhomogeneous fuchs ode}
    w''(z) + p(z)w'(z) + q(z)w(z) = z^{\beta-2} g(z) \log(z)^j, \quad |z| < R+\varepsilon, z \notin \mathbb R_{\leq 0},
\end{equation} 
has a particular solution $w(z)$ given by 
\begin{equation} \label{eq:inhomogeneous fuchs power series solution}
    w(z) = z^{\beta} \sum_{k=0}^{j+2} w_k(z)\log(z)^k,
\end{equation} 
where 
\begin{enumerate}
    \item Each $w_k(z)$ is holomorphic on $|z| < R+\varepsilon$.
    \item $w_{j+2}(z) = 0$ if $r_1-r_2 \notin \mathbb N_{\geq 0}$ or $r_2-\beta \notin \mathbb N_{\geq 0}$. In other words, a non-trivial $\log(z)^{j+2}$ factor can only occur when both $r_1-\beta, r_2-\beta \in \mathbb N_{\geq 0}$. 
    \item $w_{j+1}(z) = 0$ if both $r_1-\beta, r_2-\beta \notin \mathbb N_{\geq 0}$. In other words, a non-trivial $\log(z)^{j+1}$ factor can only occur when $r_1-\beta \in \mathbb N_{\geq 0}$ or $r_2-\beta \in \mathbb N_{\geq 0}$.
    \item There exists $C(u_1(z),u_2(z),r_1,r_2,R)$ such that
    \begin{align*}
    ||w_{j+2}(z)||_{A(|z| < R)} &\leq C \cdot (j+1)^{-2} \cdot ||g(z)||_{L^{\infty}(|z| < R)},  \\   
    ||w_{j+1}(z)||_{A(|z| < R)}  &\leq C \cdot (j+1)^{-1} \cdot ||g(z)||_{L^{\infty}(|z| < R)}.
    \end{align*}    
    If both $r_1-\beta < 0, r_2-\beta < 0$, then one also has 
   \begin{equation}
          ||w_{k}(z)||_{A(|z| < R)}  \leq C \cdot \left( \frac{j}{\min_{i\in\{1,2\}} \{|\beta-r_i| \} } \right)^j \cdot ||g(z)||_{L^{\infty}(|z| < R)} \label{exponential upper bound on w_k}
   \end{equation}
    for all $k \in \{0,1,...,j\}$. Otherwise, 
    $$||w_{k}(z)||_{A(|z| < R)} \leq C \cdot j^j \cdot \min_{i \in \{1,2\}} \text{dist}(r_i-\beta,\mathbb Z \setminus \{r_i-\beta\} )^{-(j+1)} \cdot ||g(z)||_{L^{\infty}(|z| < R)}.$$
    \item If we write
   \begin{align*}
   w'(z) = z^{\beta-1} \sum_{k=0}^{j+2} w_{1,k}(z)\log(z)^k, \\
   w''(z) = z^{\beta-2} \sum_{k=0}^{j+2} w_{2,k}(z)\log(z)^k,
   \end{align*}
    then the estimates in (4) also hold with $w_{i,k}(z)$ instead of $w_k(z)$. For the second derivative, the constant $C$ will also depend on $p(z), q(z)$. 
    \item For $n > 2$, if we write
   $$w^{(n)}(z) = z^{\beta-n} \sum_{k=0}^{j+2} w_{n,k}(z)\log(z)^k,$$
   then the estimates in (4) also hold with $w_{n,k}(z)$ instead of $w_k(z)$ if we replace 
   \begin{align*}
       C(u_1,u_2,r_1,r_2,R) &\mapsto C(n,p,q,u_1,u_2,r_1,r_2,R) \cdot j^{n-2} \cdot \prod_{i=1}^{n-2}|\beta-i-1| \\
       ||g(z)||_{L^{\infty}(|z| < R)} &\mapsto \max_{0 \leq i \leq n-2} ||g^{(i)}(z)||_{L^{\infty}(|z| < R)}
   \end{align*}
\end{enumerate}
\end{theorem}

\begin{proof}
We use the fundamental system $\{u_1,u_2\}$ given by (\ref{fuchs fundamental system}). By variation of parameters, a particular solution is given by any 
\begin{align*}
    \tilde{w}(z) = \int_{[R/2,z]} [u_2(z)u_1(y)-u_1(z)u_2(y)]W(u_1,u_2)(y)^{-1}y^{\beta-2} g(y) \log(y)^j dy, \\
    \tilde{w}'(z) = \int_{[R/2,z]} [u_2'(z)u_1(y)-u_1'(z)u_2(y)]W(u_1,u_2)(y)^{-1}y^{\beta-2} g(y) \log(y)^j dy,
\end{align*}
when $|z| < R+\varepsilon, z \notin \mathbb R_{\leq 0}$, modulo some linear combination of $\{u_1,u_2\}$. Write
\begin{align*}
    u_1(y) \cdot W(u_1,u_2)(y)^{-1} \cdot g(y) &= y^{r_1 + p_0 }\underbrace{\sum_{n=0}^{\infty}a_n y^n}_{=: a(y)}, \quad |y| < R+\varepsilon, \\
    u_2(y) \cdot W(u_1,u_2)(y)^{-1} \cdot g(y) &= y^{r_2 + p_0} \underbrace{\sum_{n=0}^{\infty}b_n y^n}_{=: b(y)} + c \cdot y^{r_1 + p_0} \log(y) a(y) , \quad |y| < R+\varepsilon.
\end{align*}
Expanding everything in the variation of parameters, we find
\begin{align*}
    \tilde{w}(z) &= z^{r_2} h_2(z) \sum_{n=0}^{+\infty} a_{n} \int_{\frac{R}{2}}^z y^{r_1 + p_0 + \beta - 2+ n} \log(y)^j dy \\
    &+ c \cdot z^{r_1}h_1(z) \log(z) \sum_{n=0}^{+\infty} a_{n} \int_{\frac{R}{2}}^z y^{r_1 + p_0 + \beta - 2 + n} \log(y)^j dy \\
    &- z^{r_1}h_1(z) \sum_{n=0}^{+\infty} b_{n} \int_{\frac{R}{2}}^z y^{r_2 + p_0 + \beta - 2 + n} \log(y)^jdy \\
    &-c \cdot z^{r_1}h_1(z) \sum_{n=0}^{+\infty} a_{n} \int_{\frac{R}{2}}^z y^{r_1 + p_0 + \beta -2 + n} \log(y)^{j+1} dy.
\end{align*}
Using  $r_1 + r_2 = 1-p_0$, this rewrites as
\begin{align*}
\tilde{w}(z) &=  z^{r_2} h_2(z) \sum_{n=0}^{+\infty} a_{n} \int_{\frac{R}{2}}^z y^{\beta-r_2 + n - 1} \log(y)^j dy \\
&+ c \cdot z^{r_1}h_1(z) \log(z) \sum_{n=0}^{+\infty} a_{n} \int_{\frac{R}{2}}^z y^{\beta-r_2 + n - 1} \log(y)^j dy \\
    &- z^{r_1}h_1(z) \sum_{n=0}^{+\infty} b_{n} \int_{\frac{R}{2}}^z y^{\beta-r_1 + n - 1} \log(y)^j dy  \\
    &- c \cdot z^{r_1}h_1(z) \sum_{n=0}^{+\infty} a_{n} \int_{\frac{R}{2}}^z y^{\beta-r_2 + n - 1} \log(y)^{j+1}dy.
\end{align*}
Replacing $z^{r_i}h_i(z)$ above by 
$$\frac{d}{dz}(z^{r_i}h_i(z)) = z^{r_i-1}(r_ih_i(z) + zh_i'(z)) = z^{r_i-1}\tilde{h}_i(z),$$ 
we get a similar formula for $\tilde{w}'(z)$. If $\delta \in \mathbb R \setminus \{-1\}$, $j \in \mathbb N_{\geq 0}$, then a primitive of $y^{\delta} \log(y)^j$ is given by 
$$z^{\delta+1} \sum_{k=0}^j \frac{(-1)^k j!}{(j-k)!} \cdot \frac{\log(z)^{j-k}}{(\delta+1)^{k+1}}$$
and if $\delta = -1$, then a primitive is given by $(j+1)^{-1}\log(y)^{j+1}$. For $k \in \{0,...,j\}$, let
\begin{align*}
    A_{j,k}(z) =  \frac{(-1)^k j!}{(j-k)!} \sum_{\substack{n=0 \\ n \neq r_2-\beta}}^{+\infty} \frac{a_{n}}{(\beta-r_2+n)^{k+1} }z^n, \quad B_{j,k}(z) =  \frac{(-1)^k j!}{(j-k)!} \sum_{\substack{n=0 \\ n \neq r_1-\beta}}^{+\infty} \frac{b_{n}}{(\beta-r_1+n)^{k+1} }z^n.
\end{align*}
Ignoring the integration constants (which are zero modulo the fundamental system), we can find a solution of the form
\begin{align*}
    w(z) &= z^{\beta} \left[ h_2(z) \sum_{k=0}^j A_{j,k}(z)\log(z)^{j-k} - h_1(z)\sum_{k=0}^j B_{j,k}(z)\log(z)^{j-k} \right. \\
    &\phantom{=}+   c \cdot z^{r_1-r_2}h_1(z) \left( \sum_{k=1}^{j} A_{j,k}(z)\log(z)^{j+1-k}- \sum_{k=1}^{j+1} A_{j+1,k}(z)\log(z)^{j+1-k} \right)  \\
    &\phantom{=} + \left( h_2(z) a_{r_2-\beta} z^{r_2-\beta} - h_1(z)b_{r_1-\beta}z^{r_1-\beta}\right) \frac{\log(z)^{j+1}}{(j+1)} \\
     &\phantom{=} +  \left. c \cdot z^{r_1-r_2}h_1(z)a_{r_2-\beta}z^{r_2-\beta} \frac{\log(z)^{j+2}}{(j+1)(j+2)}  \right],
\end{align*}
where we recall that $c = 0$ if $r_1-r_2 \notin \mathbb N_{\geq 0}$ and we used the convention that $a_{r_i-\beta} = b_{r_i-\beta} = 0$ if $r_i - \beta \notin \mathbb N_{\geq 0}$. A similar expression for $w'(z)$ can be obtained by replacing $\beta$ by $\beta-1$ and $h_i(z)$ by $r_ih_i(z) + zh_i'(z)$.

\item[\textbf{Estimates on the solution: }] It is clear from the definition that
$$||a(z)||_{L^{\infty}(|z| < R)} + ||b(z)||_{L^{\infty}(|z| < R)} \lesssim_{u_1(z),u_2(z),R} ||g||_{L^{\infty}(|z| < R)}.$$
Moreover, it follows from Cauchy-Schwarz and Parseval identity that
$$||A_{j,k}(z)||_{A(|z| < R)} + ||B_{j,k}(z)||_{A(|z| < R)} \lesssim_{u_1(z),u_2(z),R} S_{r_1,r_2,j,k,\beta} \cdot ||g||_{L^{\infty}(|z| < R)},$$
where
$$S_{r_1,r_2,j,k,\beta} = \frac{j!}{(j-k)!} \sum_{i=1}^2 \left( \sum_{\substack{n=0 \\ n \neq r_i-\beta}}^{+\infty} \frac{1}{|\beta-r_i+n|^{2(k+1)} }  \right)^{\frac{1}{2}}.$$
If both $r_1-\beta, r_2-\beta \notin \mathbb R_{\geq 0}$, one has 
$$S_{r_1,r_2,j,k,\beta} \lesssim \left( \frac{j}{\min_{i\in\{1,2\}} \{|\beta-r_i| \} } \right)^j.$$
Otherwise, 
$$S_{r_1,r_2,j,k,\beta} \lesssim  j^j \cdot \min\{|\beta-r_i+n|: i \in \{1,2\}, n \in \mathbb Z, n \neq r_i-\beta \}^{-(j+1)}.$$
Finally, one has Cauchy's inequality
$$|a_{r_i-\beta}R^{r_i-\beta}| \leq ||a(z)||_{L^{\infty}(|z| < R)}, \quad |b_{r_i-\beta}R^{r_i-\beta}| \leq ||b(z)||_{L^{\infty}(|z| < R)}.$$
Combining everything, together with the boundedness of the multiplication operator on $A(|z| < R)$, we get the desired estimate for $w(z)$ and $w'(z)$. The estimates for $w^{(2+n)}(z)$, $n \geq 0$, follow by differentiating the ODE $n$ times.
\end{proof}
\endgroup
\begin{remark}[On analytic solutions to the inhomogeneous equation]\label{remark:analytic sol inhomogeneous fuchs equation}
    If $\beta \in \mathbb N_{\geq 2}$ and $j = 0$, i.e., the inhomogeneous equation (\ref{eq:inhomogeneous fuchs ode}) has an analytic forcing term, and both $r_1-\beta, r_2-\beta \notin \mathbb N_{\geq 0}$, we observe that the solution (\ref{eq:inhomogeneous fuchs power series solution}) is analytic. In this case, the solution takes the form
    \begin{align*}
        w(z) &= z^{\beta} \left( h_2(z) A_{0,0}(z) - h_1(z) B_{0,0}(z) \right) \\
        &= z^{\beta-r_2} u_2(z) A_{0,0}(z) - z^{\beta-r_1}u_1(z) B_{0,0}(z).
    \end{align*}
    If $\beta = 2$, $j = 0$, $r_2-\beta \notin \mathbb N_{\geq 0}$ but $r_2 \in \mathbb N_{\geq 0}$ and $r_1-\beta \in \mathbb N_{\geq 0}$, then the solution  (\ref{eq:inhomogeneous fuchs power series solution}) has a logarithmic component proportional to $u_1(z)\log(z)$ that can be removed by adding an appropriate multiple of $u_2(z)$, thus yielding an analytic solution.
\end{remark}

\addcontentsline{toc}{section}{Index of Notations}
\printindex

\chapter{Estimates for solutions of parabolic equations}

The content of this chapter is taken from the appendices of my article \cite{samuelian2025blowuptrees}, to appear soon in \textit{Calculus of Variations and Partial Differential Equations}. For the sake of completeness, proofs which were omitted in the original paper have been added.

\section{Calderon-Zygmund estimates ($L^p$-estimates) for the heat equation}
In this appendix, we present the so-called Calderon-Zygmund estimates for the linear heat equation. The goal is to estimate $||D^2_x u||_{L^p_{x,t}}$ in terms of $u$ and $u_t - \Delta u$. The global estimate can be proved using the generalized theory of Calderon-Zygmund operators in $\mathbb R^{n+1}$, as in \cite[Chapter 7]{LemariRieusset2002RecentDI}, which we then localize (Theorem \ref{thm:Interior L^p Estimates for the Heat Equation}). Before that, we recall the concept of distributional solution to the heat equation, which will allow us to derive the $L^p$ estimates without assuming any a priori differentiability on the solution. We also study what kind of regularity is needed on a global solution to ensure that it is given by convolution with the Gaussian heat kernel
$$
G(x,t) = \frac{1}{(4\pi t)^{n/2}}e^{-\frac{|x|^2}{4t}}, \quad x \in \mathbb R^n, t  > 0.
$$
The usual criterion (given for example in \cite{evans10}) is that the solution should not grow faster than some exponential. We prove that boundedness of a mixed norm $L^{\infty}_{t,x}((0,T);L^p(\mathbb R^n))$ (and other related criterions) also ensures the same result (Lemma \ref{lemma:uniqueness heat equation}).

\begin{lemma}[Gaussian Kernel Pointwise Estimates]
    Let $l \in \mathbb N_{\geq 0}$ and $\alpha \in \mathbb N_{\geq 0}^n$ be a multi-index. There is a constant $C_{n,l,|\alpha|} > 0$ for which the $n$-th dimensional Gaussian heat kernel satisfies 
    \begin{align}
        |\partial_t^l \partial_{x}^{\alpha}  G(x,t)| &\leq C_{n,l,|\alpha|} \left( \frac{1}{4t} \right)^{(|\alpha|+2l+n)/2}  \left( 1 + \frac{|x|^2}{4t} \right)^{(|\alpha|+2l)/2} e^{-\frac{|x|^2}{4t}}. \label{Gaussian heat kernel upper bounds}
    \end{align}
    when $x \in \mathbb R^n$ and $t > 0$. In particular,
    \begin{align}
        |\partial_t^l \partial_{x}^{\alpha}  G(x,t)| &\leq C_{n,l,|\alpha|} \left( \frac{1}{4t} \right)^{(|\alpha|+2l+n)/2} \label{Gaussian heat kernel upper bounds, t only} \\
        |\partial_t^l \partial_{x}^{\alpha}  G(x,t)| &\leq C_{n,l,|\alpha|} \left( \frac{1}{|x|} \right)^{|\alpha|+2l+n} \label{Gaussian heat kernel upper bounds, x only},
    \end{align}
    when $x \in \mathbb R^n$ and $t > 0$.
\end{lemma}

\begin{proof}
    Assume first that $n = 1$. Let $F(x) = e^{-x^2}$. A simple proof by induction shows that $F^{(k)}(x) = p_k(x) F(x)$ for some polynomial of degree at most $k$. Then
\begin{align*}
    G(x,t) &= \left( \frac{1}{4\pi t} \right)^{1/2} F \left( \frac{x}{2 \sqrt{t}}\right) \\
    \partial_x^k G(x,t) &= \pi^{-n/2} \left( \frac{1}{4t} \right)^{(k+1)/2} p_k \left( \frac{x}{2 \sqrt{t}}  \right) F \left(\frac{x}{2 \sqrt{t}} \right) \\
    \partial_t \partial_x^k G(x,t) &\sim \pi^{-n/2} \left( \frac{1}{4t} \right)^{(k+3)/2} p_k \left( \frac{x}{2 \sqrt{t}}  \right) F \left(\frac{x}{2 \sqrt{t}} \right) \\
    &+ \pi^{-n/2} \left( \frac{1}{4t} \right)^{(k+1)/2} p_k' \left( \frac{x}{2 \sqrt{t}}  \right) \left( \frac{x}{t^{1+1/2}}  \right) F \left(\frac{x}{2 \sqrt{t}} \right)   \\ 
    &+ \pi^{-n/2} \left( \frac{1}{4t} \right)^{(k+1)/2} p_k \left( \frac{x}{2 \sqrt{t}}  \right) p_1\left( \frac{x}{2 \sqrt{t}}  \right) F \left(\frac{x}{2 \sqrt{t}} \right) \frac{x}{t^{1+1/2}} \\
    &= \pi^{-n/2} \left( \frac{1}{4t} \right)^{(k+3)/2} \tilde{p}_{k+2} \left( \frac{x}{2 \sqrt{t}}  \right) F \left(\frac{x}{2 \sqrt{t}} \right) \\
    \partial_t^l \partial_x^k G(x,t) &= \pi^{-n/2} \left( \frac{1}{4t} \right)^{(k+2l+1)/2} \tilde{p}_{k+2l} \left( \frac{x}{2 \sqrt{t}}  \right) F \left(\frac{x}{2 \sqrt{t}} \right)
\end{align*}
for some polynomial $\tilde{p}_{k+2l}$ of degree at most $k+2l$. In the multidimensional case, if $G_1(x,t)$ denotes the above one-dimensional kernel, one finds
\begin{align*}
    G(x,t) &= \prod_{i=1}^n G_1(x_i,t) \\
    \partial_{x}^{\alpha}  G(x,t) &= \prod_{i=1}^n \partial_{x_i}^{\alpha_i}G_1(x_i,t), \quad \alpha = (\alpha_1,...,\alpha_n) \\
    \partial_t^l \partial_{x}^{\alpha}  G(x,t) &\sim \sum_{\substack{\beta, |\beta| = l \\ \beta_i \geq 0}} \prod_{i=1}^n \partial_t^{\beta_i} \partial_{x_i}^{\alpha_i}G_1(x_i,t), \quad \beta = (\beta_1,...,\beta_n)
\end{align*}
and one immediately gets the upper bound
$$
|\partial_t^l \partial_{x}^{\alpha}  G(x,t)| \lesssim \left( \frac{1}{4t} \right)^{(|\alpha|+2l+n)/2}  \left( 1 + \frac{|x|^2}{4t} \right)^{(|\alpha|+2l)/2} e^{-\frac{|x|^2}{4t}}.
$$
Using the boundedness of $P(x)e^{-x}, x \in [0,+\infty)$, for any real-valued polynomial $P(x)$, we also find
$$
|\partial_t^l \partial_{x}^{\alpha}  G(x,t)| \lesssim \left( \frac{1}{4t} \right)^{(|\alpha|+2l+n)/2}.
$$
Similarly, 
$$
|\partial_t^l \partial_{x}^{\alpha}  G(x,t)| \lesssim \left( \frac{1}{|x|} \right)^{|\alpha|+2l+n} 
$$
by writing 
$$
\left( \frac{1}{4t} \right)^{(|\alpha|+2l+n)/2} = \left( \frac{1}{|x|^2} \right)^{(|\alpha|+2l+n)/2} \left( \frac{|x|^2}{4t} \right)^{(|\alpha|+2l+n)/2}
$$
in (\ref{Gaussian heat kernel upper bounds}).

\end{proof}

\begin{definition}[Distributional Solution]\label{def:distributional solution}
    Let $\Omega \subset \mathbb R^{n+1}$ be open and $f \in \mathcal{D}'(\Omega)$ be a distribution. We say that $u \in \mathcal{D}'(\Omega)$ is a distributional solution to the linear heat equation
\begin{equation}
u_t - \Delta u = f(x,t), \quad (x,t) \in \Omega, \label{eq:heat equation distribution sense}
\end{equation}
if $\langle u_t - \Delta u - f, \phi \rangle = 0$ for any test function $\phi(x,t) \in C^{\infty}_c(\Omega)$.
\end{definition}

\begin{theorem}[Hypoellipticity of the heat operator]\label{hypoellipticity}
Let $\Omega \subset \mathbb R^{n+1}$ be open. Assume that $u \in \mathcal{D}'(\Omega)$ is a distributional solution to the linear heat equation
$$
u_t - \Delta u = f(x,t), \quad (x,t) \in \Omega,
$$
where $f \in C^{\infty}(\Omega)$, then $u \in C^{\infty}(\Omega)$. 
\end{theorem}

\begin{proof}
    This follows from the fact that the Gaussian heat kernel (extended as zero for negative times) is a $C^{\infty}(\mathbb R^{n+1} \setminus \{(0,0)\})$ fundamental solution for the heat equation (\cite{folland}, Chapter 1, Theorem 1.58).
\end{proof}

\begin{lemma}[Gaussian Convolution Estimates for Inhomogeneous Problem] \label{lemma:gaussian convolution estimate}
    Let $1 \leq p \leq +\infty, 0 < T < +\infty$, $f \in L^p(\mathbb R^n \times (0,T))$. The convolution with the heat kernel
   \begin{equation}
           u(x,t) = \int_{0}^t \int_{\mathbb R^n} G(y,s)f(x-y,t-s)dyds, \quad x \in \mathbb R^n, t \in [0,T) \label{gaussian convolution solution}
   \end{equation}
    is a well-defined $C^0([0,T);L^p_x(\mathbb R^n))$-function which solves the heat equation (\ref{eq:heat equation distribution sense}) in the sense of distributions. Moreover, there exists $C = C(n)$ such that 
\begin{align}
     ||u||_{L^p(\mathbb R^n \times (0,t))}&\leq C t ||f||_{L^p(\mathbb R^n \times (0,t))}, \label{Gaussian heat kernel convolution, young inequality} \\
     ||\nabla_xu||_{L^p(\mathbb R^n \times (0,t))} &\leq C \sqrt{t} ||f||_{L^p(\mathbb R^n \times (0,t))}, \label{Gaussian heat kernel convolution, young inequality gradient} 
 \end{align}
 for all $t \in (0,T]$, where $\nabla_x u$ is a distributional (hence weak) derivative, and
 \begin{equation}
     ||u(\cdot,t)||_{L^p(\mathbb R^n)} \leq Ct^{1-\frac{1}{p}}||f||_{L^p(\mathbb R^n \times (0,t))} \label{eq: weak continuity of u(t,x), bound}
 \end{equation}
 for almost every $t \in (0,T)$. If $p = +\infty$, one actually has $u(x,t) \in C^0([0,T) \times \mathbb R^n)$.
\end{lemma}

\begin{proof}
Observe that
$$
 u =  [G \cdot 1_{(0,T)}] *_{x,t} [f \cdot 1_{(0,T)}]
 $$
 is well-defined and measurable. It follows from Young's Convolution inequality that
 \begin{align*}
     ||u||_{L^p(\mathbb R^n \times (0,T))}&\leq 
||G||_{L^1(\mathbb R^n \times (0,T))} ||f||_{L^p(\mathbb R^n \times (0,T))}  \leq ||f||_{L^p(\mathbb R^n \times (0,T))}.
\ \end{align*}
One checks that $u(x,t)$ solves the heat equation in the sense of distributions using Fubini's theorem.
 
 As for the weak derivative, the convolution makes sense as a convolution of distributions if $f$ has compact support. Hence, we can differentiate and apply Young's inequality once again to obtain 
 $$
 ||\nabla_x u||_{L^p(\mathbb R^n \times (0,T))} \leq ||\nabla_x G||_{L^1(\mathbb R^n \times (0,T))} ||f||_{L^p(\mathbb R^n \times (0,T))}
$$
since
 $$
 ||\nabla_x G(\cdot,t)||_{L^p(\mathbb R^n)} \sim t^{-(n/2)(1-1/p)-1/2}, \quad 1 \leq p \leq +\infty.
 $$
 If $f$ does not compact support, one can approximate $f$ by a sequence of $L^p(\mathbb R^n \times (0,T))$-functions $(f_n)_{n \geq 1}$ with compact support. The sequence of corresponding solutions $(u_n)_{n \geq 1}$ converge to $u$ in $L^p$, as
 $||u_n-u||_{L^p(\mathbb R^n \times (0,T))} \lesssim T||f_n-f||_{L^p(\mathbb R^n \times (0,T))}$
 thanks to (\ref{Gaussian heat kernel convolution, young inequality}) The sequence of weak derivative is Cauchy in $L^p$ thanks to (\ref{Gaussian heat kernel convolution, young inequality gradient}) and must converge to $\nabla u$ in the sense of distributions, hence in $L^p$ as well.

Next, we prove (\ref{eq: weak continuity of u(t,x), bound}). If $f(x,t) \in C^{\infty}_c(\mathbb R^n \times (0,T))$, one can first apply Young Convolution Inequality with respect to $x$ for fixed $0 < s < t < T$,
 $$
\bigg| \bigg| \int_{\mathbb R^n}G(y,s)f(x-y,t-s)dy \bigg| \bigg|_{L^p(\mathbb R^n)} \leq ||G(\cdot,s)||_{L^{1}(\mathbb R^n)} ||f(\cdot,t-s)||_{L^p(\mathbb R^n)},
 $$
 and then for fixed $0 <t < T$, use Minkowski's inequality,
\begin{align}
   ||u(\cdot,t)||_{L^p(\mathbb R^n)} &\leq  \bigg| \bigg|\int_{0}^{T}   \int_{\mathbb R^n}G(y,s)f(x-y,t-s)1_{(0,T)}(t-s)dy  ds \bigg| \bigg|_{L^p(\mathbb R^n)}\notag  \\
   &\leq \int_{0}^{T} \bigg| \bigg| \int_{\mathbb R^n}G(y,s)f(x-y,t-s)1_{(0,T)}(t-s)dy \bigg| \bigg|_{L^p(\mathbb R^n)} ds \notag \\
    &\leq \int_0^{T} ||f(\cdot,t-s)1_{(0,T)}(t-s)||_{L^p(\mathbb R^n)} ds = \int_0^{t} ||f(\cdot,t-s)||_{L^p(\mathbb R^n)} ds \notag \\
    &\leq t^{1-\frac{1}{p}} ||f||_{L^p(\mathbb R^n \times (0,t))}. \label{eq: heat equation, inhomogeneous problem, pointwise bound for fixed time}
\end{align}
One also finds
\begin{align}
     ||u(\cdot,t_1)-u(\cdot,t_2)||_{L^p(\mathbb R^n)} \leq \int_0^{T} ||f(\cdot,t_1-s)1_{(0,T)}(t_1-s) - f(\cdot,t_2-s)1_{(0,T)}(t_2-s)||_{L^p(\mathbb R^n)} ds. \label{eq: heat equation, inhomogeneous problem, time continuity in L^p}
\end{align}

More generally, if $f(x,t) \in L^p(\mathbb R^n \times (0,T))$, one has a solution $u(x,t) \in L^p(\mathbb R^n \times (0,T))$ as before. In particular, $u(\cdot,t) \in L^p(\mathbb R^n)$ is measurable and $L^p$ for almost every $t \in (0,T)$ by Fubini. One can approximate $f$ in $L^p$ by a sequence of $C^{\infty}_c(\mathbb R^n \times (0,T))$-functions $(f_n)_{n \geq 1}$ and the sequence of corresponding solutions $(u_n)_{n \geq 1}$ converge to $u$ in $L^p$. In other words,
$$
||u_n-u||_{L^p(\mathbb R^n \times (0,T))} = \big|\big| ||u_n-u||_{L^p_x(\mathbb R^n)} \big|\big|_{L^p_t(\mathbb R^n)} \to 0,
$$
meaning that, if $p < +\infty$, along a subsequence (still denoted with index $n$), $||u_n(\cdot,t)|| \to ||u(\cdot,t)||_{L^p_x(\mathbb R^n)}$ for almost every $t \in (0,T)$. Hence, one can pass to the limit in (\ref{eq: heat equation, inhomogeneous problem, pointwise bound for fixed time}) for almost every $t$. 

Similarly, one can pass to the limit in (\ref{eq: heat equation, inhomogeneous problem, time continuity in L^p}) for almost every $0 < t_1, t_2 < T$. As the translation operator in $L^p$ is continuous, one deduces that, given $\varepsilon > 0$, there is a $\delta > 0$ such that for almost every $0 < t_1, t_2 < T$, if $|t_1-t_2| < \delta$, then 
$$
||u(\cdot,t_1)-u(\cdot,t_2)||_{L^p(\mathbb R^n)} \leq \varepsilon,
$$
which proves that $u(x,t)$ has a $C^0((0,T);L^p(\mathbb R^n))$ representative.

If $p = +\infty$, the bound (\ref{eq: weak continuity of u(t,x), bound}) follows from (\ref{Gaussian heat kernel convolution, young inequality}) and continuity follows from Dominated Convergence.


\end{proof}

\begin{lemma}[Gaussian Convolution Estimates for Homogeneous Problem] \label{lemma:gaussian convolution estimate, homogeneous problem}
    Let $1 \leq p \leq +\infty$, $g \in L^p(\mathbb R^n)$. The convolution with the heat kernel
   \begin{equation}
           u(x,t) = \int_{\mathbb R^n} G(y,t)g(x-y)dy, \quad x \in \mathbb R^n, t \in [0,+\infty) \label{gaussian convolution solution, homogeneous problem}
   \end{equation}
    is a well-defined $C^0([0,+\infty);L^p_x(\mathbb R^n))$-function which solves the heat equation (\ref{eq:heat equation distribution sense}) with zero forcing term $f = 0$ in the sense of distributions. Moreover, there exists $C = C(n)$ such that 
\begin{align*}
     ||u(\cdot,t)||_{L^p(\mathbb R^n)}&\leq C ||g||_{L^p(\mathbb R^n)},\\
     ||\nabla_x u(\cdot,t)||_{L^p(\mathbb R^n)}&\leq C t^{-\frac{1}{2}} ||g||_{L^p(\mathbb R^n)},
 \end{align*}
 for all $t > 0$, where $\nabla_x u$ is a distributional (hence weak) derivative.
 
Moreover, if $p < +\infty$ or if $g \in C^0_b(\mathbb R^n)$ in case $p = +\infty$, one has
 \begin{equation}
    \lim \limits_{t \to 0} ||u(\cdot,t)-g(\cdot)||_{L^p(\mathbb R^n)} = 0.
\end{equation}
\end{lemma}

\begin{proof}
    The estimates follow from Young Convolution's inequality, as in Lemma \ref{lemma:gaussian convolution estimate}. The continuity and convergence $u(x,t) \to g(x)$ as $t \to 0$ follows from the fact that $G(x,t)$ is an approximate identity.
\end{proof}

\begin{lemma}[Uniqueness for the initial-value problem] \label{lemma:uniqueness heat equation}
    Let $0 < T < +\infty$, $f \in \mathcal{D}'(\mathbb R^n \times (0,T))$, $g \in \mathcal{D}'(\mathbb R^n)$. The initial-value problem
    \begin{align*}
        u_t - \Delta u &= f(x,t), \quad (x,t) \in \mathbb R^n \times (0,T),  \\
        u(x,0) &= g(x), \quad x \in \mathbb R^n, \notag 
    \end{align*}
    has exactly one solution (in the sense of distributions) of class $L^{\infty}_{t,x}((0,T);L^p(\mathbb R^n))$ with $1 < p < +\infty$, of class $L^p(\mathbb R^n \times (0,T))$ with $1 < p < +\infty$ and $u(x,t) = g(x)$ for all $t$ close to $t = 0$, and of class $C^{0}_b(\mathbb R^n \times [0,T))$, where the initial condition is understood as follows
    $$
    \lim \limits_{t \to 0^+} \langle u(\cdot,t), \phi \rangle = \langle g, \phi \rangle \quad \forall \phi \in C^{\infty}_c(\mathbb R^n).
    $$
    If $f \in L^p(\mathbb R^n \times (0,T))$ and $g \in L^p(\mathbb R^n)$, resp. $f \in C^0_b(\mathbb R^n \times (0,T))$ and $g \in C^0_b(\mathbb R^n)$, such solution is given by convolution of $f$ and $g$ with the heat kernel
      \begin{equation}
           u(x,t) = \int_{\mathbb R^n} G(y,t)g(x-y)dy + \int_{0}^t \int_{\mathbb R^n} G(y,s)f(x-y,t-s)dyds, \quad x \in \mathbb R^n, t \in [0,T). \label{gaussian convolution solution, full problem}
   \end{equation}
\end{lemma}

\begin{proof}
We start with existence. One checks that $u(x,t) \in L^{\infty}_{t,x}((0,T);L^p(\mathbb R^n))$ given by (\ref{gaussian convolution solution, full problem}) solves the heat equation in the sense of distributions using Fubini's Theorem. As for the initial condition, write
$$
u(x,t) = \int_{\mathbb R^n} G(y,t)g(x-y)dy + \int_{0}^t \int_{\mathbb R^n} G(y,s)f(x-y,t-s)dyds = u_0(x,t) + u_1(x,t)
$$
and let $\phi \in C^{\infty}_c(\mathbb R^n)$. As in Lemma \ref{lemma:gaussian convolution estimate}, it follows from Fubini's theorem that $u_1(\cdot,t) \in L^p(\mathbb R^n)$ for almost every $t \in (0,T)$ and
$$
| \langle u_1(\cdot,t), \phi \rangle | \leq ||u_1(\cdot,t)||_{L^p(\mathbb R^n)} ||\phi||_{L^q(\mathbb R^n)} \leq C t^{1-\frac{1}{p}} ||\phi||_{L^q(\mathbb R^n)},
$$
leading to $\lim \limits_{t \to 0^+} | \langle u_1(\cdot,t), \phi \rangle | = 0$. Similarly, $\lim \limits_{t \to 0^+} \langle u_0(\cdot,t)-g, \phi \rangle  = 0$ as $u_0(\cdot,t) \to g$ in $L^p$ (resp. $C^0_b$) when $t \to 0$ (see Lemma \ref{lemma:gaussian convolution estimate, homogeneous problem}).

In the continuous case, after a change of variables,
\begin{equation}
         u(x,t) \sim \int_{\mathbb R^n} e^{-\frac{|y|^2}{4}}g(x-t^{\frac{1}{2}}y)dy + \int_{0}^t \int_{\mathbb R^n} e^{-\frac{|y|^2}{4}}f(x-s^{\frac{1}{2}}y,t-s)dyds, \quad x \in \mathbb R^n, t \in [0,T).  \label{eq:heat convolution, change of variable, continuity}
\end{equation}
     Hence, if $f$, $g$ are continuous and bounded, then one can apply Dominated Convergence to show that $u(t,x) \in C^0_b( [0,T) \times \mathbb R^n)$ as well.

     Finally, if there were two distributional solutions $u, v$ of class $L^p(\mathbb R^n \times (0,T))$, $1 < p < +\infty$, (resp. $C^0_b$) to the initial-value problem, then $w = u-v \in L^p(\mathbb R^n \times (0,T))$ (resp. $C^0_b$) would be a solution to the initial-value problem with zero forcing term and zero initial condition. This solution is $C^{\infty}(\mathbb R^n \times (0,T))$ by hypoellipticity. Now, let $\phi_{\varepsilon}(x) \in C^{\infty}_c(\mathbb R^n)$ be a standard mollifier. Then $w_{\varepsilon} = \phi_{\varepsilon} *_x w \in C^{\infty}(\mathbb R^n \times (0,T))$ is a solution to the linear heat equation which is weakly continuous at zero with $w_{\varepsilon}(x,0) = 0$. Moreover, it is continuous at zero in the classical sense with $w_{\varepsilon}(x,0) = 0$. If $u(x,t) = v(x,t) = g(x)$ for all $t$ close to $0$, then this is verified as $w_{\varepsilon}(x,t) = 0$ for all $t$ close to $0$. Else, if $u$, $v \in L^{\infty}_{t,x}((0,T);L^p(\mathbb R^n))$, then $w_{\varepsilon}$ is bounded (by Young's inequality). Moreover, for fixed $x_0 \in \mathbb R^n$ and $\varepsilon > 0$,
     \begin{align*}
         w_{\varepsilon}(x,t) &= \int_{\mathbb R^n} w(y,t) \phi_{\varepsilon}(x-y)dx = \langle w(\cdot,t), \phi_{\varepsilon}(x-\cdot) \rangle \\
         &= \langle w(\cdot,t), \phi_{\varepsilon}(x_0-\cdot) \rangle  + \langle w(\cdot,t), \phi_{\varepsilon}(x-\cdot) - \phi_{\varepsilon}(x_0-\cdot) \rangle  = A + B.
     \end{align*}
     One has $A \to 0$ as $x \to x_0, t \to 0$ since $\phi_{\varepsilon}(x_0-\cdot) \in C^{\infty}_c(\mathbb R^n)$ is a test function. Similarly, we have the inequality
     \begin{align*}
         |B| \leq ||w_{\varepsilon}||_{L^{\infty}_{x,t}(\mathbb R^n \times (0,T))}  \int_{\mathbb R^n} |\phi_{\varepsilon}(x-y) - \phi_{\varepsilon}(x_0-y)| dy,
     \end{align*}
     which also tends to $0$ as $x \to x_0, t \to 0$ by continuity of the translation operator in $L^1$.
     
     The classical theory for the heat equation (see \cite[Chapter 2.3.1, Theorem 7]{evans10}, which we apply on $\mathbb R^n \times [0,T']$ for any $0 < T' < T$) proves that $w_{\varepsilon} = 0$ is the only solution for all $\varepsilon > 0$. 
     
     As $w \in L^p(\mathbb R^n \times (0,T))$, $w(\cdot,t) \in L^p(\mathbb R^n)$ for almost every $t \in (0,T)$ by Fubini's Theorem. In particular, $w(\cdot,t) \in C^0(\mathbb R^n) \cap L^p(\mathbb R^n)$ for almost every $t$. For such $t \in (0,T)$ fixed, $w_{\varepsilon}(\cdot,t) \to w(\cdot,t)$ in $L^p(\mathbb R^n)$ which implies $w(\cdot,t) = 0$ almost everywhere. We conclude that $w = 0$ everywhere by smoothness. One proceeds in a similar fashion for the $C^0_b$ case.
\end{proof}

Next, we introduce the notion of generalized Calderon-Zygmund operators.

\begin{definition}[Generalized Balls] \label{def: generalized open balls}
    Let $\mathbb R^n$ be endowed with the euclidean topology and let $\mu$ be a Borel measure on $\mathbb R^n$ with $0 < \mu(\mathbb R^n) \leq +\infty$. A family of generalized open balls $\{B(x,\delta)\}_{x \in \mathbb R^n, \delta > 0}$ is a collection of sets satisfying
    \begin{enumerate}
    \setcounter{enumi}{-1}
        \item $B(x,\delta)$ is open, bounded, non-empty and $B(x,\delta_1) \subset B(x,\delta_2)$ whenever $0 < \delta_1 < \delta_2$.
        \item There exists a constant $c_1 > 1$ for which $B(x,\delta) \cap B(y,\delta) \neq \emptyset \Longrightarrow B(y,\delta) \subset B(x,c_1\delta)$ for any $x,y \in \mathbb R^n,\delta >0$.
        \item There  exists a constant $c_2 > 1$ for which $\mu(B(x,c_1\delta)) \leq c_2 \mu(B(x,\delta))$ for any $x \in \mathbb R^n, \delta >0$.
        \item $\bigcap_{\delta > 0} \overline{B}(x,\delta) = \{x\}$, $\bigcup_{\delta > 0} B(x,\delta) = \mathbb R^n$ for all $x \in \mathbb R^n$.
        \item For all open set $U \subset \mathbb R^n$ and $\delta > 0$, the function $x \mapsto \mu(\{B(x,\delta) \cap U\})$ is continuous.
    \end{enumerate}
\end{definition}

\begin{definition}[Calderon-Zygmund Operator] \label{def: Calderon-Zygmund operator}
    Fix a family of generalized open balls $\{B(x,\delta)\}$ as in Definition \ref{def: generalized open balls}. A Calderon-Zygmund operator $T$ is a bounded linear operator $T : L^q(\mathbb R^n,d\mu) \rightarrow L^q(\mathbb R^n,d\mu)$, $q \in (1,+\infty]$, for which:
    \begin{enumerate}
        \item There exists a measurable function $K : \mathbb R^n \times \mathbb R^n \rightarrow \mathbb R$, called the kernel of $T$, such that for any $f \in L^q(\mathbb R^n,d\mu)$ of compact support, for almost every $x \in \mathbb R^n \setminus \supp(f)$,
        $$
        Tf(x) = \int_{\mathbb R^n} K(x,y)f(y)d\mu(y), \quad  \int_{\mathbb R^n} |K(x,y)f(y)|d\mu(y) < + \infty.
        $$
        \item There exists $c > 1$ and $A  > 0$ for which
\begin{equation}
            \int_{\mathbb R^n \setminus B(y,c\delta)} |K(x,y) - K(x,z)| d\mu(x) \leq A \quad \forall y \in \mathbb R^n, \forall \delta > 0, \forall z \in B(y,\delta). \label{def: Calderon-Zygmund operator, integral condition}
\end{equation}
    \end{enumerate}
\end{definition}

\begin{theorem}[Boundedness of Calderon-Zygmund Operators]
    Any Calderon-Zygmund operator which is $L^q$-bounded extends uniquely to a bounded operator $T : L^p(\mathbb R^n,d\mu) \rightarrow L^p(\mathbb R^n,d\mu)$ for any $1 < p < q$.
\end{theorem}

\begin{proof}
     See \cite{Stein}, Chapter 1.5, Theorem 3.
\end{proof}
  
\begin{theorem}[Global $L^p$ Estimates for the Heat Equation]\label{thm: global lp estimate}
    Let $1 < p < +\infty$ and $0 < T < +\infty$. Assume that $u \in L^{p}(\mathbb R^n \times (0,T))$ is the solution (\ref{gaussian convolution solution}) given by convolution with the heat kernel to the linear heat equation
    \begin{align*}
    u_t - \Delta u &= f(x,t), \quad (x,t) \in \mathbb R^n \times (0,T),  \\
    u(x,0) &= 0, \quad x \in \mathbb R^n,
    \end{align*}
    where $f \in L^p(\mathbb R^n \times (0,T))$.
    
    There exists $C = C(n,p)$ (independent of $T$) for which
    $$
    ||D^2_xu||_{L^p_{x,t}(\mathbb R^n\times (0,T))} \leq C ||f||_{L^p_{x,t}(\mathbb R^n\times (0,T))},
    $$
    where $D^2_x u$ is a distributional (hence weak) derivative.
\end{theorem}

\begin{proof}
 This uses the generalized theory of Calderon-Zygmund operators, following the proof in \cite[Chapter 7]{LemariRieusset2002RecentDI}. We first prove a global estimate on $||D^2_{x,t} u||_{L^p_{x,t}}$ by showing that the composition of $D^2_x$ and the solution operator is a $L^2$-bounded Calderon-Zygmund operator. We then get $L^p$-boundedness for $p > 2$ using a duality argument.

  We will assume without loss of generality that $f \in C^{\infty}_c(\mathbb R^n \times (0,T))$. If the theorem is proved for such class of forcing terms, one can approximate $f \in L^p(\mathbb R^n \times (0,T))$ by a sequence of $C^{\infty}_c(\mathbb R^n \times (0,T))$-functions $(f_n)_{n \geq 0}$, $||f_n-f||_{L^{p}(\mathbb R^n \times (0,T))} \to 0$. The corresponding sequence of solutions $u_n$ given by convolution with the heat kernel must converge to $u$ in $L^{p}(\mathbb R^n \times (0,T))$ because of (\ref{Gaussian heat kernel convolution, young inequality}). The corresponding sequence of derivatives $D^2_xu_n$ is uniformly bounded in $L^p(\mathbb R^n \times [0,T])$-norm by $2C(n,p)||f||_{L^p(\mathbb R^n \times [0,T])}$ for large $n$ and converge to $D^2_xu$ in the sense of distributions. Using test functions, we deduce easily that $||D^2_xu||_{L^p} \leq 2C(n,p)||f||_{L^p(\mathbb R^n \times [0,T])}$ as well.
 
For $x \in \mathbb R^n$, $t \in [0,T]$, one has 
\begin{align*}
    \partial_{x_i}\partial_{x_j} Sf(x,t) &=  \int_{0}^t \int_{\mathbb R^n}  G(y,s)(\partial_{y_i}\partial_{y_j}  f)(x-y,t-s)dyds = S (\partial_{y_i}\partial_{y_j}  f)(x,t)\\
    |\partial_{x_i}\partial_{x_j}  Sf(x,t)| &\lesssim t ||\Delta_x f ||_{L^{\infty}_{x,t}(\mathbb R^n \times [0,T])}.
\end{align*}
For fixed $T \geq t \geq s > 0$, an integration by parts and a change of variables show that 
\begin{align} \label{calderon zygmund integral of Delta S f}
    \partial_{x_i}\partial_{x_j} Sf(x,t)&=  \int_{0}^t \int_{\mathbb R^n}  ( \partial_{y_i}\partial_{y_j} G)(x-y,t-s)f (y,s)dyds \notag \\
    &=  \int_{\mathbb R} \int_{\mathbb R^n}  ( \partial_{y_i}\partial_{y_j} G)(x-y,t-s)f (y,s)dyds, 
\end{align}
where $(\partial_{x_i}\partial_{x_j} G)(x,t)$ is defined as $0$ for negative times and $f$ is extended as $0$ for $t \in [0,T]^c$ if $T < +\infty$ (it does not matter if $f$ is not smooth on the whole space).

We show that the convolution operator $\partial_{x_i}\partial_{x_j} \circ S$ given by (\ref{calderon zygmund integral of Delta S f}) is a bounded operator from $L^{p}_{x,t}(\mathbb R^{n+1})$ to itself for $1 < p \leq 2$. 

Let us prove first that $\Delta_x \circ S$ is a $L^2$-bounded operator as can be seen when writing (\ref{calderon zygmund integral of Delta S f}) on Fourier side. The Fourier transform with respect to $x$ of $(\Delta_x G)(x,t)$ is
$$
\left( \mathcal{F}_x \Delta_x G\right)(\xi,t) = -|\xi|^2 \left( \mathcal{F}_x G\right)(\xi,t) = -|\xi|^2 \exp \left( - |\xi|^2t \right), \quad \xi \in \mathbb R^n, t > 0
$$
and $0$ for negative times. Taking now the Fourier transform with respect to time, we find
\begin{align*}
    \left( \mathcal{F}_{x,t} \Delta_x G\right)(\xi,\tau) &= \int_0^{+\infty} -|\xi|^2e^{-|\xi|^2t -it\tau}dt \\
    &= \int_0^{+\infty} (-|\xi|^2-i\tau)e^{(-|\xi|^2 -i\tau) t}dt  +  i\tau \int_0^{+\infty} e^{(-|\xi|^2 -i\tau)t}dt \\
    &= -\frac{|\xi|^2}{|\xi|^2 + i \tau}, \quad \xi \in \mathbb R^n, \tau \in \mathbb R,
\end{align*}
meaning that $|\left( \mathcal{F}_{x,t} \Delta_x G\right)| \leq 1$, from which the $L^2$-boundedness follows:
$$
||\Delta_x S g||_{L^2_{x,t}(\mathbb R^n\times [0,+\infty))}  \leq C ||g||_{L^2_{x,t}(\mathbb R^n\times [0,+\infty))}.
$$
One also deduces
$$
|| \partial_{x_i}\partial_{x_j}  S g||_{L^2_{x,t}(\mathbb R^n\times [0,+\infty))} \lesssim  ||\Delta_x S g||_{L^2_{x,t}(\mathbb R^n\times [0,+\infty))}  \lesssim  ||g||_{L^2_{x,t}(\mathbb R^n\times [0,+\infty))} 
$$
by passing to Fourier side again.

It remains to show that $ \partial_{x_i}\partial_{x_j} \circ S$ is a Calderon-Zygmund operator by studying its kernel $1_{(0,+\infty)}(t-s) \partial_{x_i}\partial_{x_j} G(x-y,t-s)$.

Define the family of anisotropic balls 
$$
B(x,t,\delta) =  (x,t) + \{(y,s) : |y| < \delta, |s|^{1/2} < \delta \}.
$$
Because we are dealing with a convolution operator, (\ref{def: Calderon-Zygmund operator, integral condition}) rewrites as
$$
I = \int_{\mathbb R^{n+1} \setminus B(0,0,2\delta)} |  1_{(0,+\infty)}(t)\Delta_x G(x,t) - 1_{(0,+\infty)}(t+s)\Delta_x G(x+y,t+s)| dxdxt \leq A
$$
for $c = 2$ and any fixed $(y,s) \in B(0,0,\delta)$. Let $D_1 = \mathbb R^{n+1} \setminus B(0,0,2\delta) \cap \{t > 0, t+s > 0\}$, $D_2 = \mathbb R^{n+1} \setminus B(0,0,2\delta) \cap \{t > 0, t+s < 0\}$, $D_3 =  \mathbb R^{n+1} \setminus B(0,0,2\delta) \cap \{t < 0, t+s > 0\}$ and $I_1, I_2, I_3$ be the corresponding integrals. 

On $D_1$, we use the mean-value theorem. 
\begin{align*}
    I_1 &\lesssim \int_{\mathbb R^{n+1} \setminus B(0,0,2\delta)} |y \nabla_x \partial_{x_i}\partial_{x_j} G(z,r)| +  |t \partial_t \partial_{x_i}\partial_{x_j} G(z,r)|dxdt \\
& \lesssim   \int_{\mathbb R^{n+1} \setminus B(0,0,2\delta)} \delta |\nabla_x \partial_{x_i}\partial_{x_j} G(z,r)| +  \delta^2 |\partial_t \partial_{x_i}\partial_{x_j} G(z,r)|dxdt
\end{align*}
for some $(z,r) \in B(x,t,\delta)$ depending on $(x,t), (y,s)$. Since $(x,t) \in B(0,0,2\delta)^c$, $|x| \gtrsim |z|$ and $|t| \gtrsim |r|$. On $|x| \geq t^{1/2} \geq 2\delta$, we use (\ref{Gaussian heat kernel upper bounds, x only}). Otherwise, use (\ref{Gaussian heat kernel upper bounds, t only}). We find
\begin{align*}
    I_1 &\lesssim \delta \int_{4\delta^2}^{+\infty} dt \int_{|x| \geq |t|^{1/2}} \frac{1}{|z|^{n+3}} dx + \delta^2 \int_{|x| \geq 2\delta} dx \int_{|x|^2}^{+\infty} \frac{1}{|r|^{n/2+2}}dt  \\
    &\lesssim  \delta \int_{4\delta^2}^{+\infty} dt \int_{|x| \geq |t|^{1/2}} \frac{1}{|x|^{n+3}} dx + \delta^2 \int_{|x| \geq 2\delta} dx \int_{|x|^2}^{+\infty} \frac{1}{|t|^{n/2+2}}dt  \\
    &\lesssim 1.
\end{align*}
On both $D_2$ and $D_3$, $|t| \leq |s| \leq \delta^2$. Since $(x,t) \in B(0,0,2\delta)^c$, one must have $|x| \geq \delta \geq |t|^{1/2}$. We focus on $D_3$, $D_2$ being estimated similarly. We use (\ref{Gaussian heat kernel upper bounds, x only}) again and notice that $|x| \gtrsim |x+y|$. Hence,
\begin{align*}
    I_3 &\lesssim \int_{|t| \leq \delta^2} dt \int_{|x| \geq \delta} |\Delta_x G(x+y,t+s)| dx \lesssim \int_{|t| \leq \delta^2} dt \int_{|x| \geq \delta} \frac{1}{|x+y|^{n+2}} dx \\
&\lesssim \int_{|t| \leq \delta^2} dt \int_{|x| \geq \delta} \frac{1}{|x|^{n+2}} dx  \lesssim 1,
\end{align*}
which concludes our proof for $1 < p \leq 2$.

If $2 < p < +\infty$ is fixed, $q$ is the dual exponent, $g, \phi \in C^{\infty}_c(\mathbb R^{n+1})$, $\phi^*(x) = \phi(-x)$ and $||\phi||_{L^q_{x,t}} \leq 1$, one has
\begin{align*}
    \langle \partial_{x_i}\partial_{x_j} S g, \phi \rangle &= \langle G * g, \partial_{x_i}\partial_{x_j} \phi \rangle = \langle g, G * \partial_{x_i}\partial_{x_j} \phi^* \rangle = \langle g, \partial_{x_i}\partial_{x_j} S \phi^* \rangle  \\
    |\langle \partial_{x_i}\partial_{x_j} S g, \phi \rangle | &= |\langle g, \partial_{x_i}\partial_{x_j} S \phi^* \rangle| \leq ||g||_{L^p_{x,t}} \cdot ||\partial_{x_i}\partial_{x_j} S \phi^*||_{L^q_{x,t}} \lesssim ||g||_{L^p_{x,t}}.
\end{align*}
To summarize, we have proved that for $1 < p < +\infty$, there exists $C = C(n,p)$ for which
$$
||\partial_{x_i}\partial_{x_j} Sg||_{L^p_{x,t}(\mathbb R^n\times [0,+\infty))} \leq C ||f||_{L^p_{x,t}(\mathbb R^n\times [0,+\infty))}.
$$
Applying this estimate with $g(x,t) = f(x,t)$ where, if $T < +\infty$, $f(x,t)$ is extended as $0$ for $t \in [0,T]^c$, we find
\begin{align*}
    ||\partial_{x_i}\partial_{x_j} u||_{L^p_{x,t}(\mathbb R^n\times [0,T])} = ||\partial_{x_i}\partial_{x_j} S f||_{L^p_{x,t}(\mathbb R^n\times [0,T])}  &\leq ||\partial_{x_i}\partial_{x_j} Sf||_{L^p_{x,t}(\mathbb R^n\times [0,+\infty))}  \\
    &\leq C ||f||_{L^p_{x,t}(\mathbb R^n\times [0,+\infty))} \\
    &= C ||f||_{L^p_{x,t}(\mathbb R^n\times [0,T])},
\end{align*}
which concludes our proof.
\end{proof}

\begin{lemma}[Interpolation Inequality] \label{interpolation lemma, sobolev version}
    For $1 \leq p < +\infty$, there exists $K = K(n,p)$ for which
    \begin{align}
        ||\nabla_x u||_{L^{p}_{x}(\mathbb R^n)} &\leq K \left[ \varepsilon  ||D^2_x u||_{L^{p}_{x}(\mathbb R^n)}  + \varepsilon^{-1} ||u||_{L^{p}_{x}(\mathbb R^n)} \right],\label{L^p norm of nabla_x u, interpolation}
    \end{align}
    for any $\varepsilon > 0$ and $u \in W^{2,p}(\mathbb R^n)$. 
    
    Similarly, assume that $u \in L^p(\mathbb R^n \times I)$, $I \subset \mathbb R$ open, has distributional (hence weak) derivatives $\nabla_x u, D^2_x u \in L^p(\mathbb R^n \times I)$. Then,
    \begin{align}
        ||\nabla_x u||_{L^{p}_{x}(\mathbb R^n \times I)} &\leq 2K \left[ \varepsilon  ||D^2_x u||_{L^{p}_{x}(\mathbb R^n \times I)}  + \varepsilon^{-1} ||u||_{L^{p}_{x}(\mathbb R^n \times I)} \right] \label{L^p norm of nabla_x u, interpolation with time}
    \end{align}
    for any $\varepsilon > 0$.

    The Lemma still holds if one replaces $\mathbb R^n$ by the half-space $\mathbb R^n \cap \{x_n > 0\}$.
\end{lemma}

\begin{proof}
    See \cite[Lemma 5.5 and Remark 5.7]{adams2003sobolev} for the first part (\ref{L^p norm of nabla_x u, interpolation}) of the theorem. 

    As for the second part, we do the half-space $\mathbb  R^{n,+} = \mathbb R_n \cap \{x_n > 0\}$ case (and the full space one is similar). If $u \in C^{\infty}(\mathbb R^{n,+} \times I)$, one can apply (\ref{L^p norm of nabla_x u, interpolation}) for fixed $t \in I$, integrate and deduce
\begin{align*}
        \int_{I}||\nabla_x u(\cdot,t)||_{L^{p}_{x}(\mathbb  R^{n,+})}^pdt &\leq K^p \int_{I} \left[ \varepsilon  ||D^2_x u(\cdot,t)||_{L^{p}_{x}(\mathbb  R^{n,+})}  + \varepsilon^{-1} ||u(\cdot,t)||_{L^{p}_{x}(\mathbb  R^{n,+})} \right]^pdt \\
         &\leq 2^p K^p \int_{I} \left[ \varepsilon  ||D^2_x u(\cdot,t)||_{L^{p}_{x}(\mathbb  R^{n,+})}  \right]^p + \left[ \varepsilon^{-1} ||u(\cdot,t)||_{L^{p}_{x}(\mathbb  R^{n,+})} \right]^pdt,
\end{align*}
which leads to (\ref{L^p norm of nabla_x u, interpolation with time}) after taking the $p$-th root. 

Write $I = (a,b)$, $a,b \in \mathbb R \cup \{\pm\infty\}$. In the general case, it suffices to approximate $u$ using a mollification. Let $u_{\varepsilon}(x,t) = \phi_{\varepsilon} *_{x,t} u(x,t)$, $(x,t) \in \mathbb R^n \cap \{x_n > \varepsilon\} \times (a+\varepsilon, b-\varepsilon)$ (note that the value of $u(x,t)$ outside $\mathbb R^{n,+} \times I$ does not matter to compute this integral) and $\phi_{\varepsilon} = \varepsilon^{-n-1} \phi(\varepsilon^{-1}x, \varepsilon^{-1}t)$ is an approximate identity defined by some $\phi \in C^{\infty}_c(B_{1/2}(0)) \subset C^{\infty}_c(\mathbb R^{n+1})$, $\phi \geq 0$, $\int_{\mathbb R^{n+1}} \phi dxdt = 1$. 

One verifies that $u_{\varepsilon} \in C^{\infty}(\mathbb R^n \cap \{x_n \geq \varepsilon\} \times [a+\varepsilon, b-\varepsilon])$ and
$$
\nabla_x u_{\varepsilon} = \phi_{\varepsilon} * \nabla_x u, \quad D^2_x u_{\varepsilon} = \phi_{\varepsilon} * D^2_x u, \quad (x,t) \in \mathbb R^n \cap \{x_n \geq \varepsilon\}  \times [a+\varepsilon, b-\varepsilon],
$$
as well as $u_{\varepsilon} \to u$, $\nabla_x u_{\varepsilon} \to \nabla_x u$, $D^2_xu_{\varepsilon} \to D^2_x u$ in $L^p_{loc}(\mathbb R^{n,+} \times I)$. Moreover, they are uniformly bounded in $L^p(\mathbb R^n \cap \{x_n \geq \varepsilon\}  \times [a+\varepsilon, b-\varepsilon])$. For any fixed compact set $S \subset \mathbb R^{n,+} \times (a,b)$, the smooth case leads to
\begin{align*}
    ||\nabla_x u||_{L^{p}_{x}(S)} &\leq  ||\nabla_x u_{\varepsilon}||_{L^{p}_{x}(S)} \leq \lim_{\varepsilon \to 0} ||\nabla_x u_{\varepsilon}||_{L^{p}_{x}(\mathbb R^{n} \cap \{x_n > \varepsilon\} \times (a+\varepsilon, b- \varepsilon))} \\
&\leq \lim_{\varepsilon \to 0} 2K \delta  ||D^2_x u_{\varepsilon}||_{L^{p}_{x}(\mathbb R^{n} \cap \{x_n > \varepsilon\} \times (a+\varepsilon, b- \varepsilon))} \\
&+ 2K\delta^{-1} ||u_{\varepsilon}||_{L^{p}_{x}(\mathbb R^{n} \cap \{x_n > \varepsilon\} \times (a+\varepsilon, b- \varepsilon))} \\
&\leq \lim_{\varepsilon \to 0} 2K \left[ \delta ||D^2_x u||_{L^{p}_{x}(\mathbb R^{n,+} \times I)}  + \delta^{-1} ||u||_{L^{p}_{x}(\mathbb R^{n,+} \times I)} \right] \\
&\leq 2K \left[ \delta  ||D^2_x u||_{L^{p}_{x}(\mathbb R^{n,+} \times I)}  +\delta^{-1} ||u||_{L^{p}_{x}(\mathbb R^{n,+} \times I)} \right], \quad \delta > 0. 
\end{align*}
It suffices to take increasing compact sets to conclude.
\end{proof}

\begin{theorem}[Interior $L^p$ Estimates for the Heat Equation] \label{thm:Interior L^p Estimates for the Heat Equation}
    For $x,t \in \mathbb R^n \times \mathbb R$ and $r > 0$, let $C(x,t;r)$ denote the open cylinder
    $$
    C(x,t;r) := \{(y,s) \in \mathbb R^n \times \mathbb R: |x-y| < r, t-r^2 < s < t\}.
    $$
    Assume that $u \in L^p(C(x,t;r))$, $1 < p < +\infty$, $\nabla_x u \in L^p(C(x,t;r))$, is a distributional solution to the linear heat equation
    \begin{align}
    u_t - \Delta u &= f(x,t), \quad (x,t) \in C(x,t;r)
    \end{align}
    with $f \in  L^p(C(x,t;r))$.
    
    There exists $C = C(n,p)$ (independent of $x$, $t$, $r$) for which
    $$
    ||D^2_x u||_{L^p_{x,t}(C(x,t;r/2))} \leq C r^{-2} \left( ||f||_{L^p_{x,t}(C(x,t;r))} + ||u ||_{L^p_{x,t}(C(x,t;r))} \right).
    $$ 
\end{theorem}

\begin{proof}
We follow the localization argument from \cite{krylov}, which is done for the Schauder estimates (instead of the $L^p$-estimates) but also works in that setting.

After translation and parabolic scaling, one can assume without loss of generality that $x = t = 0$ and $r = 1$. 
    
    We cover $C_{1} = C(0,0;1)$ with an increasing sequence of cylinders
    $$
    C_{1} = \bigcup_{l=1}^{\infty} P_l, \quad P_l = C(0,0;R_l), \quad R_l = \sum_{i=1}^l 2^{-i},
    $$
    and consider smooth cut-off functions $\chi_l(x,t) \in C^{\infty}(\mathbb R^n \times \mathbb R;[0,1])$ for which $\chi_l = 1$ on $P_l$,
    $$
    \supp(\chi_l) \subset \{(y,s) \in \mathbb R^{n+1}: |y| \leq R_{l+1}, -R_{l+1}^2 \leq s \leq R_{l+1}^2\} = -P_{l+1} \cup P_{l+1}
    $$
    and
    $$
    ||\chi_l||_{C^{2,1}_{x,t}} \lesssim 2^{2l} =: \rho^{-l}.
    $$
    Then $v_l = u\chi_l \in L^p(\mathbb R^n \times [-1,0))$ is a solution to the initial value
\begin{align*}
        \partial_t v_l - \Delta v_l &= f\chi_l +u \partial_t \chi_l + u \Delta \chi_l + 2 \nabla_x u \cdot \nabla_x \chi_l, \quad (x,t) \in \mathbb R^n \times [-1,0) \\
        v_l(x,-1) &= 0, \quad x \in \mathbb R^n
\end{align*}
which is identically zero near $t = -1$ and has a $L^p(\mathbb R^n \times [-1,0))$-forcing term. Such a solution is unique and given by convolution with the (time-translated) heat kernel (Lemma \ref{lemma:uniqueness heat equation}). Then Global $L^p$-estimates (Theorem \ref{thm: global lp estimate}) show that $D^2_x (u \chi_l) \in L^{p}$ and
\begin{align}
   ||D^2_x(u\chi_l)||_{L^p_{x,t}(P_{l+1})} &=  ||D^2_x(u\chi_l)||_{L^p_{x,t}(\mathbb R^n \times [-1,0])} \notag  \\
   &\leq C ||f\chi_l +u \partial_t \chi_l + u \Delta \chi_l + 2 \nabla_x u \cdot \nabla_x \chi_l||_{L^p_{x,t}(\mathbb R^n \times [-1,0])} \notag \\
    &= C||f\chi_l +u \partial_t \chi_l + u \Delta \chi_l + 2 \nabla_x u \cdot \nabla_x \chi_l||_{L^p_{x,t}(P_{l+1})} \notag \\
    &\leq \rho^{-l} C\left(  ||f||_{L^p_{x,t}(P_{l+1})}  + ||u||_{L^p_{x,t}(P_{l+1})}  + ||\nabla_x u||_{L^p_{x,t}(P_{l+1})}  \right) \label{lp estimates proof, upper bound on D^2 u chi}
\end{align}
for all $l \geq 1$. We introduce back $\chi_{l+1}$
    \begin{align*}
    ||D^2_x(u\chi_l)||_{L^p_{x,t}(P_{l+1})} &\leq \rho^{-l} C\left(  ||f||_{L^p_{x,t}(P_{l+1})}  + ||u||_{L^p_{x,t}(P_{l+1})}  + ||\nabla_x u||_{L^p_{x,t}(P_{l+1})}  \right)\\
    &=\rho^{-l} C\left(  ||f||_{L^p_{x,t}(P_{l+1})}  + ||u \chi_{l+1}||_{L^p_{x,t}(P_{l+1})}  + ||\nabla_x (u\chi_{l+1})||_{L^p_{x,t}(P_{l+1})}  \right) \\
    &\leq \rho^{-l} C\left(  ||f||_{L^p_{x,t}(C_1)}  + ||u \chi_{l+1}||_{L^p_{x,t}(\mathbb R^n \times (-\infty,0])}  \right. \\
    &\left. \phantom{\leq}+ ||\nabla_x (u\chi_{l+1})||_{L^p_{x,t}(\mathbb R^n \times (-\infty,0])}  \right) 
    \end{align*}
    and from now on, $C = C(n,p) > 1$ is fixed and $\tilde{C} = \tilde{C}(n,p) > 1$ denotes the constant from the interpolation Lemma \ref{interpolation lemma, sobolev version}. One can apply Lemma \ref{interpolation lemma, sobolev version} with $\varepsilon = (2C\tilde{C})^{-1}\rho^2 \rho^l$ and obtain
   \begin{align*}
||D^2_x(u\chi_l)||_{L^p_{x,t}(P_{l+1})}  &\leq \rho^{-l} 4C\tilde{C}\left( ||f||_{L^p_{x,t}(C_1)}   + \varepsilon ||D^2_x (u\chi_{l+1})||_{L^p_{x,t}(\mathbb R^n \times (-\infty,0])}  \right. \\
& \left.\phantom{\leq}+ \varepsilon^{-1} ||u \chi_{l+1}||_{L^p_{x,t}(\mathbb R^n \times (-\infty,0])}  \right) \\
    &\leq \rho^{-l} 4C \tilde{C} \left( ||f||_{L^p_{x,t}(C_1)}  + \varepsilon ||D^2_x (u\chi_{l+1})||_{L^p_{x,t}(P_{l+2})} + \varepsilon^{-1} ||u||_{L^{p}_{x,t}(C_1)} \right) \\
    &\leq \rho^{-2l} \underbrace{(4C \tilde{C})^2}_{=: \hat{C}} \left( ||f||_{L^p_{x,t}(C_1)}  + ||u ||_{L^p_{x,t}(C_1)} \right) + \frac{\rho^2}{2}|| D^2_x (u\chi_{l+1})||_{L^p_{x,t}(P_{l+2})}.
  \end{align*}
  Starting with $l = 1$ and iterating, we obtain
\begin{align}
||D^2_x(u\chi_1)||_{L^p_{x,t}(P_{2})} &\leq  \rho^{-2}\hat{C} \left( ||f||_{L^p_{x,t}(C_1)}  + ||u ||_{L^p_{x,t}(C_1)} \right)  + \frac{\rho^2}{2}|| D^2_x (u\chi_{2})||_{L^p_{x,t}(P_{3})} \notag \\
&\leq \left( \rho^{-2}\hat{C} + \frac{\rho^2}{2} \cdot \rho^{-4}\hat{C} \right) \left( ||f||_{L^p_{x,t}(C_1)}  + ||u ||_{L^p_{x,t}(C_1)} \right)  \notag \\
&\phantom{\leq}+ \left(\frac{\rho^2}{2}\right)^2 || D^2_x (u\chi_{3})||_{L^p_{x,t}(P_{4})} \notag \\
&\leq \rho^{-2} \hat{C} \left( 1 + \frac{1}{2} + ... + \frac{1}{2^{k-2}} \right) \left( ||f||_{L^p_{x,t}(C_1)}  + ||u ||_{L^p_{x,t}(C_1)} \right) \notag \\
&\phantom{\leq}+ \left(\frac{\rho^2}{2}\right)^{k-1} || D^2_x (u\chi_{k})||_{L^p_{x,t}(P_{k+1})} \label{lp estimate proof: final inequality}
  \end{align}
  for any $k > 1$. Since $\rho < 1$ and
\begin{align*}
      || D^2_x (u\chi_{k})||_{L^p_{x,t}(P_{k+1})}  &\leq \rho^{-k} C \left(  ||f||_{L^p_{x,t}(C_{1})}  + ||u||_{L^p_{x,t}(C_{1})}  + ||\nabla_x u||_{L^p_{x,t}(C_{1})}  \right)  
\end{align*}
by  (\ref{lp estimates proof, upper bound on D^2 u chi}), we can let $k \to +\infty$ in (\ref{lp estimate proof: final inequality}). Observing that $u = u \chi_1$ on $C(0,0;1/2) \subset P_2$, this concludes the proof. 
\begingroup
\allowdisplaybreaks
We note that we only need an a priori estimate on $u, \nabla_x u$ but not on $D^2_x u$ to take $k \to +\infty$. The same will be true for the Schauder estimates (Theorem \ref{interior schauder estimates}).
\end{proof}

We prove a last lemma which will be useful in the context of the Parabolic Sobolev Embedding.

\begin{lemma}[Global $L^p$ Estimates for compactly supported solutions of parabolic equations]\label{lemma:global l^p estimates for parabolic operator, compact support}
    Let $\alpha \in (0,1)$, $0 < T < +\infty$ and $1 < p < +\infty$. Consider a parabolic operator
    $$
    L = c(x,t) + \sum_{i=1}^n b^i(x,t) D_{x_i} + \sum_{i,j=1}^n a^{i,j}(x) D_{x_i}D_{x_j}
    $$
    defined on $\overline{\Omega}$, $\Omega = B_{\kappa}(x_1) \times (0,T)$, satisfying
    \begin{align*}
    ||a,b,c||_{\alpha,\alpha/2;\overline{\Omega}} &\leq \Lambda, \quad \Lambda > 0, \\
        \Lambda^{-1}|\xi|^2 \leq \sum_{i,j=1}^n a^{i,j}(x,t)\xi_i\xi_j &\leq \Lambda |\xi|^2, \quad (x,t) \in \overline{\Omega}, \xi \in \mathbb R^n,
    \end{align*}
    and where $a^{i,j}(x)$ does not depend on time.
    
    There exists $\kappa_0 = \kappa_0(n,\alpha,\Lambda) > 0$ and $C = C(n,p,\alpha,\Lambda) > 0$ such that if $\kappa \leq \kappa_0$, for any $u \in L^{\infty}_{t,x}((0,T);L^p(\mathbb R^n))$ with weak derivatives $\nabla_x u, (\partial_t - L)u, D^2_xu \in L^p(B_{\kappa}(x_1) \times (0,T))$, that is compactly supported in the sense that $\overline{\supp(u(\cdot,t))} \subset B_{\kappa}(x_1)$ for all $t \in (0,T)$ and which is weakly continuous at zero in the sense that
    $$
    \lim \limits_{t \to 0^+} \langle u(\cdot,t), \phi \rangle = 0 \quad \forall \phi \in C^{\infty}_c(\mathbb R^n),
    $$
    then one has
    $$
    ||D^2_x u||_{L^p(\mathbb R^n \times (0,T))} \leq C\left( ||(\partial_t - L)u||_{L^p(\mathbb R^n \times (0,T))} + ||u||_{L^p(\mathbb R^n \times (0,T))} \right).
$$
\end{lemma}

\begin{proof}
    Without loss of generality, we translate the problem and assume that $x_1 = 0$. We also assume that $b^i = c = 0$ and only deal with the second order terms of $L$, as one can use Interpolation (Lemma \ref{interpolation lemma, sobolev version}) with a sufficiently small $\varepsilon$ for the lower order terms. 
    
    Write
    $$
\partial_t u - \sum_{i,j=1}^n a^{i,j}(x) D_{x_i}D_{x_j}u = \partial_t u - \sum_{i,j=1}^n a^{i,j}(x_1) D_{x_i}D_{x_j}u + \sum_{i,j=1}^n \left[ a^{i,j}(x_1) - a^{i,j}(x) \right] D_{x_i}D_{x_j}u.
    $$
    Let $P$ be an orthogonal matrix for which $P^T (a_{i,j}(x_1)) P = I_n$. Let $v(x,t) = u(Px,t)$ (in other words, we are rotating $\mathbb R^n$). Then 
    $$(\partial_t - \Delta)v = \partial_t u - \sum_{i,j=1}^n a^{i,j}(x_1) D_{x_i}D_{x_j}u = (\partial_t - L_0)u.
    $$
    As
\begin{align*}
    ||(\partial_t - L_0)u||_{L^{p}_{x,t}(\mathbb R^n \times (-\infty,0])} &= ||(\partial_t - L_0)u||_{L^{p}_{x,t}(\Omega)} \leq ||(\partial_t - L)u||_{L^{p}_{x,t}(\Omega)}  \\
    &\phantom{=}+ \bigg|\bigg| \sum_{i,j=1}^n \left[ a^{i,j}(x_1) - a^{i,j}(x) \right] D_{x_i}D_{x_j}u \bigg|\bigg|_{L^{p}_{x,t}(\Omega)}  \\
    &\lesssim  ||(\partial_t - L)u||_{L^p(\mathbb R^n \times (-\infty,0])} + \kappa^{\alpha} \Lambda ||D^2_x u||_{L^p(\mathbb R^n \times (-\infty,0])},
 \end{align*}
 the function $v \in L^{\infty}((0,T);L^p_x(\mathbb R^n))$ must be the unique solution to the heat equation
     \begin{align*}
    v_t - \Delta v &= (\partial_t - L_0)u \in L^p(\mathbb R^n \times (0,T)),  \\
    v(x,0) &= 0, \quad x \in \mathbb R^n,
    \end{align*}
    given by convolution with the heat kernel (Lemma \ref{lemma:uniqueness heat equation}). Note that $D^2_xu $, $(\partial_t - L_0) u$ and $(\partial_t - L)u$ were naturally extended as zero for $x \in B_{1}(0)^c$. Hence one can estimate
\begin{align*}
        ||D^2_x u||_{L^p(\mathbb R^n \times (-\infty,0])} 
        &\lesssim ||D^2_x v||_{L^p(\mathbb R^n \times (-\infty,0])} \\
        &\leq C\left( ||(\partial_t - \Delta)v||_{L^{p}_{x,t}(\mathbb R^n \times (-\infty,0])} + ||v||_{L^{p}_{x,t}(\mathbb R^n \times (-\infty,0])} \right) \\
    &\lesssim C\left( ||(\partial_t - L_0)u||_{L^{p}_{x,t}(\mathbb R^n \times (-\infty,0])} + ||u||_{L^{p}_{x,t}(\mathbb R^n \times (-\infty,0])} \right) 
\end{align*}
    using Theorem \ref{thm: global lp estimate}. Choose $\kappa$ small enough and the theorem is proved.
\end{proof}

\endgroup




\section{Parabolic Sobolev Embedding}
In this section, we prove a version of the Sobolev-Morrey embedding into a $C^{1+\alpha, \alpha/2}_{x,t}$ Hölder space from a parabolic Sobolev space where $u$, $\nabla_x u$, $D^2_x u$, $\partial_t u \in L^p(\Omega \times [0,T])$ and $p$ is large enough (Theorem \ref{thm:parabolic sobolev embedding}). The classical reference for this result is \cite[Lemma 3.3, Chapter 2]{ladyzhenskaia1968linear}. However, after verification, the proof is not in this book and could not be found elsewhere. I thank my advisor, Pr. Joachim Krieger, for the proof of Theorem \ref{thm:parabolic sobolev embedding} which uses a simple Littlewood-Paley decomposition. We then generalize this estimate to the half-space setting with Dirichlet boundary condition by reflection, we localize the estimate (Corollary \ref{cor: local parabolic sobolev embedding}) and finally obtain it for general bounded domains with smooth boundary by straightening the boundary and  freezing the coefficients of the newly obtained parabolic operator (Theorem \ref{thm:boundary parabolic sobolev embedding}).

\begin{definition}[Smooth Boundary]\label{def:smooth boundary}
 Let $\Omega \subset \mathbb R^n$ be an open set. We say that $\Omega$ has a $C^{\alpha}$-boundary, $\alpha \in \mathbb R_{> 0} \cup \{\infty\}$, if for any $x_0 \in \partial \Omega$, there is an open neighborhood $U$ of $x_0 \in \mathbb R^n$ and a $C^{\alpha}$-map $\gamma : \mathbb R^{n-1} \rightarrow \mathbb R$ for which, up to an orthogonal change of coordinates,
    $$
    \Omega \cap U = \{x \in U : \gamma(x_1,...,x_{n-1}) < x_n\}.
    $$
    In particular,
    $$
    \partial \Omega \cap U = \{x \in U : \gamma(x_1,...,x_{n-1}) = x_n\}.
    $$

    We say that $\Omega$ satisfies the interior cone condition if there exists a bounded, open, cone $C \subset \mathbb R^n$ of fixed aperture angle $\kappa \in (0,\pi]$ and height $h > 0$ such that for any $x \in \Omega$, there is a cone $C_x \subset \Omega$ congruent to $C$ (i.e. obtained by translation, rotation and symmetry from $C$) with vertex/apex at $x$.
    
    We say that  $\Omega$ satisfies the interior ball condition if there exists $\rho > 0$ such that for any $x \in \Omega$, there exists a ball $x \in B_{\rho}(y) \subset \Omega$.
\end{definition}

\begin{lemma}
    Let $U$ be an open, bounded set. If $U$ has a $C^1$ boundary, then it satisfies the interior cone condition. If $U$ has a $C^2$ boundary, then it satisfies the interior ball condition.
\end{lemma}

\begin{proof}
Consider a boundary point $a = (\tilde{a},a_n) \in \partial \Omega$ and straighten locally the boundary as in the definition \ref{def:smooth boundary}. The mapping $\psi : (x_1,...,x_n) \mapsto (x_1,...,x_{n-1},x_n - \gamma(x_1,...,x_{n-1}))$ is a diffeomorphism which maps a neighborhood $a \in V \subset \mathbb R^n$ 
     to some open set $U = \psi(V)$, $\psi(\Omega \cap V) = U \cap \{y \in \mathbb R^n : y_n > 0\}$ and $\psi(\partial \Omega \cap V) = U \cap \{y \in \mathbb R^n : y_n = 0\}$. 
     
     Take a small ball $B_{\rho}(\psi(a)) \subset U$. For $h = \sqrt{\rho}/4$ and any point $p \in B_{\rho/2}(\psi(a)) \cap \{y_n \geq 0\}$, the open cone 
     $$
     p + C = \{p + (\tilde{y},y_n) \in \mathbb R^n : 0 < |\tilde{y}| < y_n, 0 < y_n < h\} 
     $$
     stays in the upper half-ball $B^+_{\rho}(\psi(a)) = \{y \in B_{\rho}(\psi(a)) : y_n > 0\}$.
     
     Indeed, let $L$ be the Lipschitz constant of $\gamma$ on some compact set $K$ containing the projection of $(\tilde{x},x_n) \in \overline{\Omega} \mapsto \tilde{x} \in K \subset \mathbb R^{n-1}$. Then we show that the cone
     $$
    \psi^{-1}(p) + \tilde{C} = \{\psi^{-1}(p) + (\tilde{x},x_n) \in \mathbb R^n : 0 < |\tilde{x}| < x_n/(2L+1), 0 < x_n < 2h/3\} \subset \Omega
     $$
    is included in $\psi^{-1}(p+C) \subset \Omega$. Indeed, for such a point $\psi^{-1}(p) + (\tilde{x},x_n) \in \psi^{-1}(p) + \tilde{C}$, one has
    \begin{align*}
\psi(\psi^{-1}(p) + (\tilde{x},x_n)) &= \psi(\tilde{p}+\tilde{x},p_n + \gamma(\tilde{p}) + x_n) =
(\tilde{p}+\tilde{x}, p_n + \gamma(\tilde{p})+ x_n - \gamma(\tilde{p}+\tilde{x})) \\
&= p + (\tilde{x},x_n+ \gamma(\tilde{p}) -\gamma(\tilde{p}+\tilde{x})).
    \end{align*}
    Then we find that
    $$
x_n+ \gamma(\tilde{p}) -\gamma(\tilde{p}+\tilde{x}) \in (x_n/2,h) \subset (0,h),
    $$
    because $|\gamma(\tilde{p}) - \gamma(\tilde{p}+\tilde{x})| \leq L|\tilde{x}| \leq Lx_n/(2L+1) < x_n/2$ and $0 < 3x_n/2 < h$, as well as 
    $$
|\tilde{x}| \leq x_n/(2L+1) \leq x_n/2 < x_n+ \gamma(\tilde{p}) -\gamma(\tilde{p}+\tilde{x}),
    $$
    meaning that $\psi(\psi^{-1}(p) + (\tilde{x},x_n)) \in p + C$.

    For each $a \in \partial \Omega$, there is a neighborhood $V_a \subset \mathbb R^n$ for which one has a $C^1$-coordinate chart $\psi : \mathbb R^n \to \mathbb R^n$ as above such that $\psi(V_a)$ contains a ball $B_{\rho(a)}(\psi(a))$, $\psi$ and $\gamma$ are of Lipschitz constant $\leq L(a)$ on $\overline{\Omega}$ and its $(n-1)$-th dimensional projection, and $\psi^{-1}$ of Lipschitz constant $J(a)$ on $\psi(\Omega)$. Let $a_1,...,a_m \in \partial \Omega$ such that $\psi^{-1}_i(B_{\rho_i/(4J_iL_i+1)}(a_i))$, $\rho_i = \rho(a_i)$, $L_i = L(a_i)$, $J_i = J(a_i)$ is a finite cover of $\partial \Omega$, where $\psi_i$ is the coordinate chart associated with $a_i$. 
    
    Then for each $a \in \partial \Omega$, $a \in \psi^{-1}_i(B_{\rho_i/(4J_iL_i+1)}(\psi_i(a_i)))$ for some $i \in \{1,...,m\}$ and one can use the coordinate chart $\psi_i$ to deduce
   \begin{align*}
       |a-a_i| \leq |\psi^{-1}(\psi(a)) - \psi^{-1}(\psi(a_i))| \leq J_i \rho_i/(4J_iL_i+1) &\leq \rho_i/4L_i \\
       |\psi(a)-\psi(a_i)| \leq L_i |a-a_i| \leq  L_i\rho_i/(4L_i) &\leq \rho_i/4,
   \end{align*}
    meaning that $V_{a_i}$ is a neighborhood of $a$ for which $\psi(V_{a_i})$ contains a ball $B_{\rho_i/4}(\psi(a))$ (to be distinguished from the initial $B_{\rho_i}(\psi(a_i))$). Hence, the above proof shows that for any $a \in \psi^{-1}(B_{\rho_i/8}(\psi_i(a_i)) \cap \{y_n \geq 0\})$, up to an orthogonal transformation, $\Omega$ contains an open cone $C$ of vertex $a$, height $\min_i \sqrt{\rho_i}/8$, direction $(0,...0,1)$ and angle depending only on $\{L_1,...,L_m\}$. 

   A simple argument by contradiction finally shows that the set
   $$
\bigcup_{i=1}^m \psi^{-1}(B_{\rho_i/8}(\psi_i(a_i)) \cap \{y_n \geq 0\})
   $$
   covers an $\varepsilon$-neighborhood 
   $$
   \{x \in \overline{\Omega} : \dist(x,\partial \Omega) < \varepsilon\}
   $$
   of $\Omega$. Taking a smaller cone of height $\lesssim \varepsilon$, we can put a congruent cone $C_x$ at all points $x \in \Omega$.

   The argument in case of a $C^2$-boundary is similar.
\end{proof}

\begin{lemma}[\cite{ladyzhenskaia1968linear}, Chapter 2, Lemma 3.1] \label{lemma:ladyzenskaia reduction Holder regularity} 
    Let $\Omega \subset \mathbb R^{n}$ be an open set satisfying the interior cone condition with a cone $C$ of height $h$. Let $I \subset \mathbb R$ be any open interval.

    Assume that $u \in C^0_b(\Omega \times I)$ is $\alpha$-Hölder-continuous with respect to $t$, meaning that 
    $$
    |u(x,t)-u(x,t')| \leq \mu_1 |t-t'|^{\alpha} \quad \forall x \in \Omega, \forall t,t' \in I.
    $$
    Further assume that $\nabla_x u \in C^0_b(\Omega \times I)$ is $\beta$-Hölder continuous with respect to $x$, meaning that 
    $$
    |\nabla_x u(x,t)- \nabla_x u(y,t)| \leq \mu_2 |x-y|^{\beta} \quad \forall x,y \in \Omega, \forall t \in I.
    $$
    
    Then $\nabla_x u$ is $\frac{\alpha \beta}{1+\beta}$-Hölder continuous with respect to $t$ with Hölder constant
    $$
    |\nabla_x u(x,t)- \nabla_x u(x,s)| \leq 2\left( \mu_1+\mu_2+ h^{-\beta}||\nabla_x u||_{L^{\infty}(\Omega \times I)} \right) |t-s|^{\frac{\alpha \beta}{1+\beta}} \quad \forall t,s \in I, \forall x \in \Omega.
    $$
    In particular, if $\Omega = \mathbb R^n$ or $\Omega = \mathbb R^n \cap \{x_n > 0\}$, $h$ can be taken arbitrarily large so that 
    $$
    |\nabla_x u(x,t)- \nabla_x u(x,s)| \leq 2\left( \mu_1+\mu_2 \right) |t-s|^{\frac{\alpha \beta}{1+\beta}} \quad \forall t,s \in I, \forall x \in \Omega.
    $$
\end{lemma}

\begin{proof}
    Let $x_1,x_2 \in \Omega$ be two points such that the segment $[x_1,x_2] \subset \Omega$ stays in $\Omega$. Let $v = x_2-x_1$ and write $\partial_v u(y,s) = (x_2-x_1) \cdot \nabla_x u(y,s)$. For $t_1,t_2 \in I$, the mean-value theorem for integrals shows that
    \begin{align}
        &\left| \int_{0}^1 \partial_v u(x_1 + h(x_2-x_1),t_2) -  \partial_v u(x_1,t_1)dh \right|  = \left| \partial_v u(x^*,t_2) -  \partial_v u(x_1,t_1) \right|, \quad x^* \in [x_1,x_2] \notag  \\
        &= \left| \partial_v u(x^*,t_2) - \partial_v u(x_1,t_2) + \partial_v u(x_1,t_2) -  \partial_v u(x_1,t_1) \right| \notag \\
        &\geq \left| \partial_v u(x_1,t_2) -  \partial_v u(x_1,t_1) \right|-\mu_2|x_2-x_1|^{1+\beta}. \label{parabolic sobolev embedding: intermediate estimate}
    \end{align}
    Write 
    \begin{align*}
        \partial_v u(x_1 + h(x_2-x_1),t_2) -  \partial_v u(x_1,t_1) &= \partial_v u(x_1 + h(x_2-x_1),t_2) - \partial_v u(x_1 + h(x_2-x_1),t_1) \\
        &+ \partial_v u(x_1 + h(x_2-x_1),t_1) -\partial_v u(x_1,t_1), 
    \end{align*}
    and further estimate
    \begin{align*}
        \left| \int_{0}^1 \partial_v u(x_1 + h(x_2-x_1),t_1) -\partial_v u(x_1,t_1)  dh \right| \leq \mu_2|x_2-x_1|^{1+\beta}, 
    \end{align*}
    as well as
    \begin{align*}
        \biggl| \int_0^1 \partial_v u(x_1 + h(x_2-x_1),t_2) &- \partial_v u(x_1 + h(x_2-x_1),t_1)  dh \biggr| \\
        &=  \bigr| \left[u(x_2,t_2)-u(x_1,t_2)\right]- \left[u(x_2,t_1)-u(x_1,t_1)\right] \bigl| \\
        &\leq |u(x_2,t_2)-u(x_2,t_1)| +|u(x_1,t_2)-u(x_1,t_1)| \\
        &\leq 2 \mu_1 |t_1-t_2|^{\alpha}.
    \end{align*}
    Plugging those in (\ref{parabolic sobolev embedding: intermediate estimate}), we deduce that  
    $$
\left| \nabla_x u(x_1,t_2) -  \nabla_x u(x_1,t_1) \right| \leq 2\mu_2|x_2-x_1|^{\beta}  + 2 \mu_1 \frac{|t_1-t_2|^{\alpha}}{|x_2-x_1|}.
    $$
    Now, any $x_1 \in \Omega$ is the vertex of a cone $C_{x_1} \subset \Omega$ of height $h > 0$ (which we may assume to be closed by taking a slightly smaller cone in the definition). If $|t_1-t_2| \leq h^{(1+\beta)\alpha^{-1}}$, we can choose $x_2 \in C_{x_1}$ with $[x_1,x_2] \subset C_{x_1}$ and $|x_1-x_2| = |t_1-t_2|^{(1+\beta)^{-1}\alpha} \leq h$, so that
$$
\left| \nabla_x u(x_1,t_2) -  \nabla_x u(x_1,t_1) \right| \leq 2(\mu_1 + \mu_2) |t_2-t_1|^{\frac{\alpha \beta}{1+\beta}}.
$$
If $|t_1-t_2| > h^{(1+\beta)\alpha^{-1}}$ instead, then
\begin{align*}
    \left| \nabla_x u(x_1,t_2) -  \nabla_x u(x_1,t_1) \right| &\leq 2 ||\nabla_x u||_{L^{\infty}(\Omega \times I)} = 2 ||\nabla_x u||_{L^{\infty}(\Omega \times I)} h^{-\beta} h^{\beta} \\
    &\leq 2 ||\nabla_x u||_{L^{\infty}(\Omega \times I)} h^{-\beta} |t_2-t_1|^{\frac{\alpha \beta}{1+\beta}}.
\end{align*}
\end{proof}

 \begin{proposition}[Extension by reflection]\label{prop:extension by reflection}
     Let 
     $$
     u \in C^{2,1}(\mathbb R^n \cap \{x_n > 0\} \times I) \cap L^p(\mathbb R^n \cap \{x_n > 0\} \times I),$$
     $I \subset \mathbb R$ open, $1 \leq p \leq +\infty$. Assume that $\lim \limits_{\varepsilon \to 0^+}u(x,\varepsilon,t) \to 0$ in $L^p_{loc}(\mathbb R^{n-1} \times I)$ and $\lim \limits_{\varepsilon \to 0^+} \nabla_x u (\tilde{x},\varepsilon,t)$ converges to a limit in $L^p_{loc}(\mathbb R^{n-1} \times I)$. Denote by $\nabla_x u(\tilde{x},0,t)$ its limit. Let 
    $$
\tilde{u}(\tilde{x},x_n,t) = \begin{cases} u(\tilde{x},x_n,t), \quad &x_n > 0, \\
- u(\tilde{x},- x_n,t), \quad &x_n < 0, \\
\end{cases}
$$
be the odd reflection of $u$. 

If the distribution $v \in \{\partial_{x_i} u, \partial_{x_i}\partial_{x_j}u, \partial_t u, (\partial_t - \Delta)u, \partial_{x_n}^2 u: i,j \neq n\}$ belongs to  $L^p(\mathbb R^n \cap \{x_n > 0\} \times I)$, then the corresponding distribution $\tilde{v} \in \{\partial_{x_i} \tilde{u}, \partial_{x_i}\partial_{x_j}\tilde{u}, \partial_t \tilde{u}, (\partial_t - \Delta)\tilde{u}, \partial_{x_n}^2 \tilde{u}: i,j \neq n\}$ is given by the odd reflection of $v$. If $v \in \{\partial_{x_n} u, \partial_{x_i}\partial_{x_n} u: i \neq n\}$ belongs to $ L^p(\mathbb R^n \cap \{x_n > 0\} \times I)$, then the corresponding distribution $\tilde{v}$ is given by the even reflection of $v$ (i.e., $\tilde{v}(\tilde{x},x_n,t) = v(\tilde{x},-x_n,t)$ for $x_n < 0$).

Additionally, if $u \in C^{2+\alpha,1+\alpha/2}(\mathbb R^n \cap \{x_n \geq 0\} \times I)$ for $\alpha \in [0,1)$, then $\tilde{u} \in  C^{2+\alpha,1+\alpha/2}(\mathbb R^n \cap \{x_n \geq 0\} \times I)$ as well.

The same result holds if one replaces $\mathbb R^n \cap \{x_n > 0\}$ by a half-ball $B_r^+(x_0) = B_r(x_0) \cap \{x_n > 0\}$ with $(x_0)_n = 0$ and one assumes $\lim \limits_{\varepsilon \to 0^+} u(\tilde{x},\varepsilon,t) \to 0$ in $L^p_{loc}(B_{r}(\tilde{x}_0) \times I)$, where $x_0 = (\tilde{x}_0, (x_0)_n)$, and $\lim \limits_{\varepsilon \to 0^+}  \nabla_x u (\tilde{x},\varepsilon,t)$ converges in $L^p_{loc}(B_r(\tilde{x}_0) \times I)$.
 \end{proposition}

 \begin{proof}
      First, $\tilde{u} \in L^p(\mathbb R^n \times I)$. Next, let $\phi \in C^{2,1}_c(\mathbb R^{n} \times I)$. Write $(\partial_t - \Delta) u = f \in L^p(\mathbb R^n \cap \{x_n > 0\} \times I)$ and $\tilde{f}$ its odd reflection to the whole plane. For $\varepsilon > 0$, integrating by parts shows
      \begin{align*}
   \int_{\mathbb R^n \cap \{x_n > \varepsilon\} \times I} u \partial_{x_i}^{\alpha_i}\partial_{x_j}^{\alpha_j} \phi dxdt &= (-1)^{\alpha_i + \alpha_j}\int_{\mathbb R^n \cap \{x_n > \varepsilon\} \times I} \partial_{x_i}^{\alpha_i}\partial_{x_j}^{\alpha_j}u \phi dxdt, \quad i,j \neq n, \alpha_k \in \{0,1\},  \\
       \int_{\mathbb R^n \cap \{x_n > \varepsilon\} \times I} u \partial_{t} \phi dxdt &= -\int_{\mathbb R^n \cap \{x_n > \varepsilon\} \times I} \partial_{t}u \phi dxdt, \\
         \int_{\mathbb R^n \cap \{x_n > \varepsilon\} \times I} u \partial_{x_i}^{\alpha_i}\partial_{x_n} \phi dxdt &= (-1)^{\alpha_i}\int_{\mathbb R^{n-1} \times I} \partial_{x_i}^{\alpha_i} u(\tilde{x},\varepsilon,t) \phi(\tilde{x},\varepsilon,t) d\tilde{x}dt \\
         &+ (-1)^{1+\alpha_i} \int_{\mathbb R^n \cap \{x_n > \varepsilon\} \times I} \partial_{x_i}^{\alpha_i} \partial_{x_n}u \phi dxdt, \quad i \neq n, \alpha_i \in \{0,1\},
\end{align*}
as well as
\begin{align*}
     \int_{\mathbb R^n \cap \{x_n > \varepsilon\} \times I} u \partial_{x_n}^2 \phi dxdt &= \int_{\mathbb R^{n-1} \times I} u(\tilde{x},\varepsilon,t) \partial_{x_n}\phi(\tilde{x},\varepsilon,t) d\tilde{x}dt \\
     &- \int_{\mathbb R^{n-1} \times I} \partial_{x_n} u(\tilde{x},\varepsilon,t) \phi(\tilde{x},\varepsilon,t) d\tilde{x}dt +\int_{\mathbb R^n \cap \{x_n > \varepsilon\} \times I} \partial_{x_n}^2u \phi dxdt. \\
\end{align*}
As $u(\tilde{x},\varepsilon,t) \to 0$ in $L^p(\mathbb R^{n-1} \times I)$ and $\phi$ is $C^{2,1}$ with compact support, taking the limit $\varepsilon \to 0$ shows
\begin{align*}
    \int_{\mathbb R^n \cap \{x_n > 0\} \times I} u (-\partial_t - \Delta) \phi dxdt &= \int_{\mathbb R^{n-1} \times I} \partial_{x_n} u(\tilde{x}, 0, t)  \phi(\tilde{x},0,t) d\tilde{x}dt \\
    &+ \int_{\mathbb R^n \cap \{x_n > 0\} \times I } f \phi dx dt
\end{align*}
if $f = (\partial_t - \Delta)u \in L^p(\mathbb R^n \cap \{x_n > 0\} \times I)$.

Denoting $\tilde{\phi}(\tilde{x},x_n,t) = -\phi(\tilde{x},-x_n,t)$, then
\begin{align*}
    \langle (\partial_t - \Delta) \tilde{u}, \phi \rangle &=  \langle \tilde{u}, (-\partial_t - \Delta) \phi \rangle \\
    &= \int_{\mathbb R^n \cap \{x_n > 0\} \times I} u (-\partial_t - \Delta) \phi dxdt + \int_{\mathbb R^n \cap \{x_n > 0\} \times I} u (-\partial_t - \Delta) \tilde{\phi} dxdt \\
    &= \int_{\mathbb R^n \cap \{x_n > 0\} \times I } f \cdot (\phi + \tilde{\phi}) dx dt = \int_{\mathbb R^n \times I } \tilde{f} \cdot \phi dx dt,
 \end{align*}
where $\tilde{f} \in L^p(\mathbb R^n \times I)$ is the odd reflection of $f$.  We proceed similarly for the other derivatives.

If $u, \nabla_x u, \partial_tu, D^2_x u$ are continuous up to the boundary $x_n = 0$, then differentiating the boundary condition yields $\partial_t u = \partial_{x_i} u = \partial_{x_i}\partial_{x_j} u = 0$ on $\{x_n = 0\}$ for $i,j \in \{1,...,n-1\}$. It also follows from $(\partial_t - \Delta) u = 0$ that $\partial_{x_n}^2 u = 0$ on $\{x_n = 0\}$. We also have continuity of $\partial_{x_n} u$, $\partial_{x_i}\partial_{x_n} u$, $i \in \{1,...,n-1\}$, on $\{x_n = 0\}$ and it does not matter whether the value is zero or not as they are extended using an even reflection. Hence, $\partial_t \tilde{u}, \nabla_x \tilde{u}, D^2_x \tilde{u}$, which are of class $C^{0}(\mathbb R^n \setminus \{x_n = 0\} \times (-\infty,0])$, extend continuously on $\{x_n = 0\}$. If $D^2_x u \in C^{\alpha}(\mathbb R^n \cap \{x_n > 0\} \times (-\infty,0])$, then the only issue for the Hölder regularity of $\partial_{x_i} \partial_{x_j} \tilde{u}$ (and similarly for $\partial_t \tilde{u}$) is when one has two points $x, y \in \mathbb R^n$ on opposite half-planes $y_n < 0 < x_n$. 

If $\partial_{x_i} \partial_{x_j} \tilde{u}$ is an even reflection, then
$$
|\partial_{x_i} \partial_{x_j} \tilde{u}(\tilde{x},x_n,t) - \partial_{x_i} \partial_{x_j} \tilde{u}(\tilde{y},y_n,t)| = |\partial_{x_i} \partial_{x_j} \tilde{u}(\tilde{x},x_n,t) - \partial_{x_i} \partial_{x_j} \tilde{u}(\tilde{y},-y_n,t)|
$$
and $(\tilde{x},x_n),(\tilde{y},-y_n)$ are on the same half-plane. If it is an odd reflection, then $\partial_{x_i} \partial_{x_j} \tilde{u} = 0$ on $\{x_n = 0\}$ and
\begin{align*}
|\partial_{x_i} \partial_{x_j} \tilde{u}(\tilde{x},x_n,t) - \partial_{x_i} \partial_{x_j} \tilde{u}(\tilde{y},y_n,t)| &\leq |\partial_{x_i} \partial_{x_j} \tilde{u}(\tilde{x},x_n,t) - \partial_{x_i} \partial_{x_j} \tilde{u}(\tilde{x},0,t)| \\
&+ |\partial_{x_i} \partial_{x_j} \tilde{u}(\tilde{y},y_n,t) - \partial_{x_i} \partial_{x_j} \tilde{u}(\tilde{y},0,t)| \\
&\leq  [D^2_xu]_{C^{\alpha}} \cdot \left( |x_n|^{\alpha} + |y_n|^{\alpha} \right) \\
&\lesssim |x_n - y_n|^{\alpha} \lesssim |x-y|^{\alpha}
\end{align*}
as $y_n < 0 < x_n$.

The same argument works for a ball by symmetry of the domain. The integration by parts have no terms on $\partial B_r(x_0)$ because the test function must be compactly supported inside the ball.
 \end{proof}
 
\begin{theorem}[Global Parabolic Sobolev Embedding]\label{thm:parabolic sobolev embedding}
   Let $n+2 < p < +\infty$ and $I \subset \mathbb R$ be any open interval. Assume that $u(x,t) \in L^{p}(\mathbb R^{n} \times I)$ has weak derivatives $\nabla_x u, (\partial_t - \Delta) u \in L^p(\mathbb R^{n} \times I)$. Then $u \in C^{1+2\gamma,\gamma}_{x,t}(\mathbb R^{n}\times I)$ where
    $$
    0 < \gamma < \frac{(p-n-1)(p-n-2)}{p(2p-n-2)}.
    $$
    Moreover, there exists $C = C(n,p,\gamma)$ for which
\begin{align}
    ||u||_{C^{1+2\gamma,\gamma}_{x,t}(\mathbb R^{n}\times I)}
&\leq C \left( 1 + |I|^{-1-2\gamma-\frac{n+2}{p}}\right) \left( ||u||_{L^p(\mathbb R^n \times I)} + ||(\partial_t - \Delta)u||_{L^p(\mathbb R^n \times I)} \right).  \label{ineq:global parabolic sobolev}
\end{align}
The theorem stills holds if one replaces $\mathbb R^n$ by the half-space $\mathbb R^n \cap \{x_n > 0\}$ and one a priori assumes that $u \in C^{2,1}(\mathbb R^n \cap \{x_n > 0\} \times I)$, $\lim \limits_{\varepsilon \to 0^+} u(x,\varepsilon,t) \to 0$ in $L^p_{loc}(\mathbb R^{n-1} \times I)$ and $\lim \limits_{\varepsilon \to 0^+} \nabla_x u (\tilde{x},\varepsilon,t)$ converges in $L^p_{loc}(\mathbb R^{n-1} \times I)$.
\end{theorem}

\begin{notation}
    In the following, $u \in W^{1,p}(\mathbb R^n \times I)$ means that $u, \nabla_x u$ and $\partial_t u$ are in $L^p(\mathbb R^n \times I)$.
\end{notation}

\begin{remark}
    The convergence of $\lim \limits_{\varepsilon \to 0^+}u(x,\varepsilon,t) \to 0$ and $\lim \limits_{\varepsilon \to 0^+} \nabla_x u (\tilde{x},\varepsilon,t)$ in $L^p_{loc}(\mathbb R^{n-1} \times I)$ is a technicality needed for the boundary estimate. As in Theorem \ref{thm:boundary parabolic sobolev embedding}, one can assume that $u \in C^0([0,T),C^1(\mathbb R^n \cap \{x_n \geq 0\}))$ and $u(\tilde{x},0,t) = 0$ instead, which is the regularity obtained by doing a fixed point argument in Theorem \ref{thm:taylor local existence} if one replaces the half-space by a bounded domain with smooth boundary.
\end{remark}

\begin{proof}
First, observe that $D^2_x u, \partial_t u \in L^p(\mathbb R^n \times I)$ thanks to Global $L^p$ estimates (Theorem \ref{thm: global lp estimate}). Then $u(x,t) \in W^{1,p}(\mathbb R^n \times I)$ extends as a function $u(x,t) \in W^{1,p}(\mathbb R^{n+1}) \hookrightarrow C^{0,1-\frac{n+1}{p}}(\mathbb R^{n+1})$ (see \cite[Theorem 5.28]{adams2003sobolev}). In particular, $u$ is bounded on $\mathbb R^{n} \times I$ and Hölder-continuous with respect to $t$ on $I$, uniformly with respect to $x \in \mathbb R^n$. The Hölder constant $\mu_1$ satisfies
$$
\mu_1 \leq C(n,p,I) ||u||_{W^{1,p}(\mathbb R^n \times I)}.
$$

Similarly, for fixed $t \in I$ and $0 < \varepsilon < 1 - \frac{n}{p}$, $u(x,t) \in W^{2,p}(\mathbb R^n) \hookrightarrow W^{2-\varepsilon,p}(\mathbb R^n) \hookrightarrow C^{1,1-\frac{n}{p}-\varepsilon}(\mathbb R^n)$. The $C^{0,1-\frac{n}{p}-\varepsilon}(\mathbb R^n)$-Hölder-norm of $\nabla_x u(x,t)$ is uniform with respect to $t \in \mathbb R$ provided that $u \in L^{\infty}_{t}(W^{2-\varepsilon,p}_x)$, which we will be proving now for an appropriate $\varepsilon$ to be chosen later. The $L^{\infty}$-norm and the Hölder constant $\mu_2$ will satisfy
$$
||\nabla_x u||_{L^{\infty}(\mathbb R^n \times I)} + \mu_2 \leq C ||u||_{L^{\infty}_{t}(W^{2-\varepsilon,p}_x)} \leq C(n,p,\varepsilon) \left( ||u||_{W^{1,p}(\mathbb R^n \times I)} + ||\Delta u||_{L^{p}(\mathbb R^n \times I)} \right).
$$

First, we recall some facts about Sobolev spaces. For $s \in \mathbb R_{\geq 0}$ and $1  < p <\infty$, recall that
$$
W^{s,p}(\mathbb R^n) = H^{s,p}(\mathbb R^n) = \{u \in S'(\mathbb R^n): (1-\Delta)^{s/2} u \in L^p(\mathbb R^n) \} \ (= F^{s}_{p,2}(\mathbb R^n)),
$$
where $S'(\mathbb R^n)$ denotes the set of Schwartz distributions, $F^{s}_{p,2}$ is a Triebel-Lizorkin space and $(1-\Delta)^{s/2}$ is the (spatial) Fourier multiplier $(1-\Delta)^{s/2} = (1+|D|^2)^{s/2}$. For a reference, see \cite[Corollary 3.1, Chapter 3.2.1.2]{sawano}, for integer $s$ and then use complex interpolation to get the spaces inbetween (Theorem 4.29, Chapter 4.2.3.3 from \cite{sawano}). 

We apply a nonhomogeneous spatial Littlewood-Paley decomposition and write
$$
u = \chi(D)u + \sum_{j \geq 0} \phi(2^{-j}D)u,
$$
where $\chi(x)$, resp. $\phi(x)$,  are smooth and supported on a ball, resp. an annulus, centered at the origin (for fixed $t \in \mathbb R$, one has convergence in the sense of Schwartz distribution). 

For $j \geq 0$, one can write 
$$
(1 - \Delta)^{s/2}  \phi(2^{-j}D) u = (1 - \Delta)^{s/2} \theta(D)  \phi(2^{-j}D) u, \quad s \in \mathbb R,
$$
where $\theta \in C^{\infty}$ is a smooth function which is zero in a neighborhood of the origin. Then
\begin{align}
    \big\| (1 - \Delta)^{(2-\varepsilon)/2} \phi(2^{-j}D) u\big\|_{L^p(\mathbb R^n\times I)} &= \big\| \big\| (1 - \Delta)^{(2-\varepsilon)/2} \phi(2^{-j}D) u \big\|_{L^p_x(\mathbb R^n)}^p \big\|_{L^1_t(I)}^{1/p} \notag \\
    &= \big\| \big\| (1 - \Delta)^{-\varepsilon/2} \phi(2^{-j}D)  (1 - \Delta)  u \big\|_{L^p_x(\mathbb R^n)}^p \big\|_{L^1_t(I)}^{1/p} \notag  \\
    &\leq C(\varepsilon)  \left( 2^j \right)^{-\varepsilon} \big\| \big\| \phi(2^{-j}D)  (1 - \Delta)  u\big\|_{L^p_x(\mathbb R^n)}^p \big\|_{L^1_t(I)}^{1/p} \notag \\
    &\leq  C(\varepsilon) \left( 2^j \right)^{-\varepsilon} \big\| \big\|(1 - \Delta) u\big\|_{L^p_x(\mathbb R^n)}^p \big\|_{L^1_t(I)}^{1/p} \notag  \\
    &=  C(\varepsilon) \left( 2^j \right)^{-\varepsilon}\big\|(1 - \Delta) u\big\|_{L^p(\mathbb R^n\times I)}  \label{application of bernstein lemma}
\end{align}
using Bernstein lemma (\cite[Chapter 2.1, Lemma 2.2]{bahouri}) applied with the multiplier $(1-\Delta)^{-\varepsilon/2} \theta(D)$ for the first inequality and Young Convolution Inequality in the second ($\phi(2^{-j}D)u = \mathcal{F}^{-1}(\phi(2^{-j}\cdot ))(x) * u(x,t)$ and $\mathcal{F}^{-1}(\phi(2^{-j}\cdot))(x)$ has $L^1$-norm independent of $j$). Similarly, it holds that
\begin{align*}
     \big\| (1 - \Delta)^{(2-\varepsilon)/2} \chi(D) u\big\|_{L^p(\mathbb R^n\times I)} &= \big\| \big\| (1 - \Delta)^{(2-\varepsilon)/2} \chi(D) u \big\|_{L^p_x(\mathbb R^n)}^p \big\|_{L^1_t(I)}^{1/p} \\
     &\leq  C(\varepsilon)  \big\| \big\| u \big\|_{L^p_x(\mathbb R^n)}^p \big\|_{L^1_t(I)}^{1/p} =  C(\varepsilon)  ||u||_{L^p(\mathbb R^n \times I)} 
\end{align*}
using Young Convolution Inequality with the function $\mathcal{F}^{-1} \left[ (1 + |\cdot |^2)^{(2-\varepsilon)/2} \chi(\cdot) \right](x)$ which is Schwartz (in particular $L^1(\mathbb R^n)$) since it is the inverse Fourier transform of a smooth, compactly supported function.

At the same time, if $J \subset I$ is bounded, then
\begin{align*}
    \left(2^j \right)^{-p \varepsilon} \big\|(1 - \Delta) u\big\|_{L^p(\mathbb R^n\times I)}  &\gtrsim_{\varepsilon} \big\|(1 - \Delta)^{(2-\varepsilon)/2} \phi(2^{-j}D)u \big\|_{L^p( \mathbb R^n\times J)}^p \\
    &\geq |J|\cdot \inf_{t\in J}\int_{\mathbb R^n}|(1 - \Delta)^{(2-\varepsilon)/2} \phi(2^{-j}D)u(\cdot, t)|^p\,dx.
\end{align*}
Hence, if $J \subset I$ is an interval of length $\leq \left( 2^j \right)^{-p \tilde{\varepsilon}}$ with $0 < \tilde{\varepsilon} < \varepsilon$, there exists $t \in J$ with
$$
 \int_{\mathbb R^n}|(1 - \Delta)^{(2-\varepsilon)/2}  \phi(2^{-j}D)u(\cdot, t)|^p\,dx \leq 1 + \inf_{t\in J}\int_{\mathbb R^n}|(1 - \Delta)^{(2-\varepsilon)/2} \phi(2^{-j}D)u(\cdot, t)|^p\,dx,
$$
so that
$$
\big\|(1 - \Delta)^{(2-\varepsilon)/2} \phi(2^{-j}D) u(\cdot, t)\big\|_{L^p(\mathbb R^n)}\lesssim_{\varepsilon} \left( 2^j \right)^{-p(\varepsilon - \tilde{\varepsilon})}\big\|(1 - \Delta) u\big\|_{L^p(\mathbb R^n\times I)}.  
$$
Similarly, one gets
$$
\big\|(1 - \Delta)^{(2-\varepsilon)/2} \chi(D) u(\cdot, t)\big\|_{L^p(\mathbb R^n)}\lesssim_{\varepsilon} ||u||_{L^p(\mathbb R^n \times I)} 
$$
for an interval of length $|J| \leq 1$ and some $t \in J$.

If $t_1 \in J$ is another time in this interval $J$, we show that the control of 
$$
\big\|(1 - \Delta)^{(2-\varepsilon)/2} \phi(2^{-j}D) u(\cdot, t)\big\|_{L^p(\mathbb R^n)}
$$
can be extended to $t_1$. For this, we write
$$
\phi(2^{-j}D)u(x,t_1) = \mathcal{F}^{-1}(\phi(2^{-j}\cdot ))(x) * \left(u(x,t) + \int_{t}^{t_1} \partial_t u(x,s)ds \right)
$$
and use Minkowski's inequality,
\begin{align*}
    \big\|(1 - \Delta)^{(2-\varepsilon)/2} \phi(2^{-j}D) u(\cdot, t_1) \big\|_{L^p(\mathbb R^n)}
&\leq \big\|(1 - \Delta)^{(2-\varepsilon)/2} \phi(2^{-j}D) u(\cdot, t) 
 \big\|_{L^p(\mathbb R^n)} \\
 &+ \int_t^{t_1} \big\|(1 - \Delta)^{(2-\varepsilon)/2} \phi(2^{-j}D) \partial_t u(\cdot, s) \big\|_{L^p(\mathbb R^n)}\,ds.
\end{align*}

The term at time $t$ in the right-hand side is $\lesssim_{\varepsilon} \left( 2^j \right)^{-p(\varepsilon-\tilde{\varepsilon})}||(1 - \Delta) u||_{L^p(\mathbb R^n\times I)}$ and we only need to estimate the double integral, which is
\begin{align*}
    \int_t^{t_1} \big\|(1 - \Delta)^{(2-\varepsilon)/2} \phi(2^{-j}D) \partial_t &u(\cdot, s) \big\|_{L^p(\mathbb R^n)}\,dt \\
    &\leq |t-t_1|^{1-p^{-1}} \cdot \|(1 - \Delta)^{(2-\varepsilon)/2} \phi(2^{-j}D) \partial_t u(\cdot, s) \|_{L^p(\mathbb R^n\times J)}
\end{align*}
by Hölder's inequality.

Since $t_1\in J$, we have $|t-t_1|^{1-p^{-1}}\leq |J|^{1-p^{-1}} \leq \left( 2^j \right)^{-p\tilde{\varepsilon}(1-p^{-1})} = \left( 2^j \right)^{-\tilde{\varepsilon}(p-1)}$. Similarly to (\ref{application of bernstein lemma}), we have the bound 
$$
\|(1 - \Delta)^{(2-\varepsilon)/2}   \phi(2^{-j}D)  \partial_t u\|_{L^p(\mathbb R^n\times J)}\lesssim_{\varepsilon}  \left( 2^j \right)^{2-\varepsilon} \cdot \|\partial_t u\|_{L^p(\mathbb R^n\times J)}
$$
using Bernstein lemma with the multiplier $(1 - \Delta)^{(2-\varepsilon)/2} \theta(D)$ and Young convolution's inequality. It follows that
\begin{align*}
    |t-t_1|^{1-p^{-1}}\cdot \|(1 - \Delta)^{(2-\varepsilon)/2} \phi(2^{-j}D) &\partial_t u(\cdot, s) \|_{L^p(\mathbb R^n\times J)} \\
    &\lesssim_{\varepsilon}  \left( 2^j \right)^{2- \tilde{\varepsilon} p - (\varepsilon - \tilde{\varepsilon}) }\|\partial_t u\|_{L^p(\mathbb R^n\times J)} \\
    &= C(\varepsilon) \left( 2^j \right)^{2- \varepsilon p + (\varepsilon - \tilde{\varepsilon})(p-1) } \|\partial_t u\|_{L^p(\mathbb R^n\times J)}.
\end{align*}
Recall that $0 < \varepsilon < 1 - \frac{n}{p}$. In addition, if $\varepsilon > \frac{2}{p}$ and $0 < \tilde{\varepsilon} < \varepsilon$ is close enough to $\varepsilon$, then
$$
 |t-t_1|^{1-p^{-1}}\cdot \|(1 - \Delta)^{(2-\varepsilon)/2} \phi(2^{-j}D) \partial_t u(\cdot, s) \|_{L^p(\mathbb R^n\times J)} \lesssim_{\varepsilon}  \left( 2^j \right)^{-\delta} \|\partial_t u\|_{L^p(\mathbb R^n\times J)} 
$$
for some small $\delta > 0$.

Such an $\varepsilon > 0$ is guaranteed to exist if
$$
\frac{2}{p} < 1 - \frac{n}{p} \iff p > n+2
$$
As we can split $I$ into intervals of length $\leq \left( 2^j \right)^{-p \tilde{\varepsilon}}$ and $t_1$ was arbitrary, we have proved that
\begin{align*}
    \sup_{t\in I}\|\phi(2^{-j}D) (1 - \Delta)^{(2-\varepsilon)/2}  u(t)  \|_{L^p(\mathbb R^n)} &\lesssim_{\varepsilon}  \left[  \left( 2^j \right)^{-p(\varepsilon - \tilde{\varepsilon})} + \left( 2^j \right)^{-\delta} \right] \\
    &\phantom{\leq} \cdot \left[ \|(1-\Delta) u\|_{L^p(\mathbb R^n\times J)} 
 +\|\partial_t u\|_{L^p(\mathbb R^n\times J)} 
 \right] \\
 &\lesssim_{\varepsilon}  \left( 2^j \right)^{-\tilde{\delta}} \left[ \|(1-\Delta) u\|_{L^p(\mathbb R^n\times J)} 
 + \|\partial_t u\|_{L^p(\mathbb R^n\times J)} \right].
\end{align*}
Similarly, one can prove
$$
\sup_{t\in  I}\|\chi(D) (1 - \Delta)^{(2-\varepsilon)/2}  u(\cdot, t)  \|_{L^p(\mathbb R^n)} \lesssim_{\varepsilon}  \big\|u\big\|_{L^p(\mathbb R^n\times J)} 
 + \big\|\partial_t u\big\|_{L^p(\mathbb R^n\times J)}.
$$
Finally, going back to the Littlewood-Paley decomposition, we find that
\begin{align*}
    \sup_{t\in  I}\| u(\cdot, t)  \|_{W^{2-\varepsilon,p}(\mathbb R^n)} &= \sup_{t\in I}\|(1 - \Delta)^{(2-\varepsilon)/2}  u(\cdot, t)  \|_{L^p(\mathbb R^n)} \\
    &\lesssim_{\varepsilon}  \left[ 1 + \sum_{j \geq 0} \left( 2^j \right)^{-\tilde{\delta}} \right] \cdot \left[ ||u||_{L^p(\mathbb R^n \times I)} + ||\partial_t u||_{L^p(\mathbb R^n \times I)} + ||\Delta u||_{L^p} \right] \\
    &\leq C(\varepsilon) \left[ ||u||_{L^p(\mathbb R^n \times I)} + ||\partial_t u||_{L^p(\mathbb R^n \times I)} + ||\Delta u||_{L^p(\mathbb R^n \times I)} \right],
\end{align*}
which concludes the proof that $\nabla_x u \in L^{\infty}_t(W^{2-\varepsilon,p}_x)$ for $\frac{2}{p} < \varepsilon < 1 - \frac{n}{p}$.

Finally, we choose $\varepsilon = \frac{2}{p} + \delta$, $0 < \delta \ll 1$ arbitrarily small.  By applying Lemma \ref{lemma:ladyzenskaia reduction Holder regularity}, we find that $\nabla_x u(x,t)$ is $\gamma$-Hölder continuous with respect to $t$, uniformly with respect to $x \in \mathbb R^n$, where 
$$
    \gamma = \frac{(p-n-1)(p-n-2)}{p(2p-n-2)} - \hat{\delta} = \frac{\alpha \beta}{1+\beta}, \quad \alpha = 1 - \frac{n+1}{p}, \ \beta = 1 - \frac{n+2}{p} - \delta.
$$
As $2\gamma \leq \beta$, it follows from (\ref{equivalence Holder norm}) that $\nabla_x u(x,t) \in C^{2\gamma, \gamma}_{x,t}(\mathbb R^{n} \times I )$.

To obtain a final inequality where $||\nabla_x u||_{L^p}$ is absent from the right-hand side, one uses the Interpolation Lemma \ref{interpolation lemma, sobolev version}. Similarly, to replace $||D^2_xu||_{L^p}$ and $||\partial_t u||_{L^p}$ by $||(\partial_t - \Delta)u||_{L^p}$, one uses the Global $L^p$ Estimates (Theorem \ref{thm: global lp estimate}).

We finally note that using translations and parabolic rescaling, it follows that the constant $C =
C(n, p, \gamma, I)$ from the theorem depends only on $|I|$ with $$
C(n,p,\gamma,I) \leq
C(n, p, \gamma,[0,1]) \left( 1 + |I|^{-1-2\gamma-\frac{n+2}{p}} \right).
$$

In the half-space case, if $u \in C^{2,1}(\mathbb R^n \cap \{x_n > 0\} \times I\}$, it suffices to extend $u$ by odd reflection (Proposition \ref{prop:extension by reflection}).
\end{proof}

\begin{corollary}[Local Parabolic Sobolev Embedding]\label{cor: local parabolic sobolev embedding}
   Let $n+2 < p < +\infty$, $\Omega \subset \mathbb R^n$ be open and $I \subset \mathbb R$ be open as well. Assume that $u(x,t) \in L^{p}( \Omega \times I)$ has weak derivatives $\nabla_x u, (\partial_t-\Delta) u \in L^p(\Omega \times I)$. For any ball $\overline{B_r(x_0)} \subset \Omega$, one has $u \in C^{1+2\gamma,\gamma}_{x,t}(B_{r/2}(x_0) \times I)$ where
    $$
    0 < \gamma < \frac{(p-n-1)(p-n-2)}{p(2p-n-2)}.
    $$
    Moreover, there exists $C = C(n,p,\gamma)$ for which
    \begin{align*}
            ||u||_{C^{1+2\gamma,\gamma}_{x,t}(B_{r/2}(x_0) \times I)} &\leq C \left( 1 + (|I|/r^2)^{-1-2\gamma-\frac{n+2}{p}} + r^{-1-2\gamma-\frac{n+2}{p}} \right) \\
            &\phantom{\leq}\cdot\left( ||u||_{L^p(\Omega \times I)} + ||(\partial_t - \Delta) u||_{L^p(\Omega \times I)} \right).
        \end{align*}
        The same result holds if one replaces $B_{r}(x_0)$, resp. $B_{r/2}(x_0)$, by some half-ball $B_r^+(x_0) = B_r(x_0) \cap \{x_n > (x_0)_n\}$ (up to a rotation), where the flat part $T_r^+(x_0) = B_r(x_0) \cap \{x_n = (x_0)_n\} \subset \overline{\Omega}$ can touch the boundary while $\overline{B_r^+(x_0)} \setminus T_r^+(x_0) \subset \Omega$. Moreover, one a priori assumes that $u \in C^{2,1}(B_r^+(x_0) \times I)$, $\lim \limits_{\varepsilon \to 0^+} u(\tilde{x},(x_0)_n+\varepsilon,t) \to 0$ in $L^p_{loc}(B_{r}(\tilde{x}_0) \times I)$, where $x_0 = (\tilde{x}_0, (x_0)_n)$, $\lim \limits_{\varepsilon \to 0^+} \nabla_x u(\tilde{x},(x_0)_n+\varepsilon,t)$ converges in $L^p_{loc}(B_{r}(\tilde{x}_0) \times I)$.
        
\end{corollary}

\begin{proof}
  We use the same localization argument from Theorem \ref{interior schauder estimates}. After translation and parabolic rescaling, assume that $x_0 = 0$ and $r = 1$. 

     We cover $B_{1} = B_1(0)$ with an increasing sequence of balls
    $$
    B_{1} = \bigcup_{l=1}^{\infty} B_{R_l}(0), \quad R_l = \sum_{i=1}^l 2^{-i}
    $$
    and consider smooth cut-off functions $\chi_l(x) \in C^{\infty}(\mathbb R^n \times \mathbb R;[0,1])$ for which $\chi_l = 1$ on $B_{R_l}$,
    $$
    \supp(\chi_l) \subset B_{R_{l+1}}
    $$
    and
    $$
    ||\chi_l||_{C^{2}_{x}} \lesssim 2^{2l} =: \rho^{-l}
    $$
    Let $f = (\partial_t - \Delta) u$. Then $v_l = u\chi_l \in L^p(\mathbb R^n \times I)$ (resp. $\in L^p(\mathbb R^{n} \cap \{x_n > 0\} \times I)$ in the half-ball case) satisfies
\begin{align*}
        \partial_t v_l - \Delta v_l &= f\chi_l + u \Delta \chi_l + 2 \nabla_x u \cdot \nabla_x \chi_l, \quad (x,t) \in \mathbb R^n \times I,
\end{align*}
where the right-hand side is a $L^p(\mathbb R^n \times I)$-function (resp. $\in L^p(\mathbb R^{n} \cap \{x_n > 0\} \times I)$). The following inequalities will also hold if one replaces $\mathbb R^n$ by $\mathbb R^n \cap \{x_n > 0\}$ and $B_{R_l}$ by $B_{R_l}^+$. The Global Parabolic Sobolev Embedding (Theorem \ref{thm:parabolic sobolev embedding}) shows that
\begin{align}
   ||u\chi_l||_{C^{1+2\gamma,\gamma}_{x,t}(B_{R_{l+1}} \times I)} &=  ||u\chi_l||_{C^{1+2\gamma,\gamma}_{x,t}(\mathbb R^n \times I)} \notag  \\
   &\leq C ||f\chi_l  + u \Delta \chi_l + 2 \nabla_x u \cdot \nabla_x \chi_l||_{L^p_{x,t}(\mathbb R^n \times I)} + C ||u \chi||_{L^p_{x,t}(\mathbb R^n \times I)}  \notag \\
    &= C||f\chi_l+ u \Delta \chi_l + 2 \nabla_x u \cdot \nabla_x \chi_l||_{L^p_{x,t}(B_{R_{l+1}} \times I)} + C||u \chi_l||_{L^p_{x,t}(B_{R_{l+1}} \times I)} \notag \\
    &\leq \rho^{-l} C\left(  ||f||_{L^p_{x,t}(B_{R_{l+1}} \times I)}  + ||u||_{L^p_{x,t}(B_{R_{l+1}} \times I)}  + ||\nabla_x u||_{L^p_{x,t}(B_{R_{l+1}} \times I)}  \right) \label{parabolic estimates proof, upper bound on u chi}
\end{align}
    We introduce back $\chi_{l+1}$
    \begin{align*}
    ||u\chi_l||_{C^{1+2\gamma,\gamma}_{x,t}(B_{R_{l+1}} \times I)} &\leq \rho^{-l} C\left(  ||f||_{L^p_{x,t}(B_{R_{l+1}} \times I)}  + ||u||_{L^p_{x,t}(B_{R_{l+1}} \times I)}  + ||\nabla_x u||_{L^p_{x,t}(B_{R_{l+1}} \times I)}  \right)  \\
   &= \rho^{-l} C\left(  ||f||_{L^p_{x,t}(B_{R_{l+1}} \times I)}  + ||u \chi_{l+1}||_{L^p_{x,t}(B_{R_{l+1}} \times I)}  + ||\nabla_x (u \chi_{l+1})||_{L^p_{x,t}(B_{R_{l+1}} \times I)}  \right) \\
    &\leq \rho^{-l} C\left(  ||f||_{L^p_{x,t}(\Omega \times I)}  + ||u \chi_{l+1}||_{L^p_{x,t}(\mathbb R^n \times I)}  + ||\nabla_x (u\chi_{l+1})||_{L^p_{x,t}(\mathbb R^n \times I)}  \right) 
    \end{align*}
    and from now on, $C = C(n,p,\gamma,|I|) > 1$ is fixed and $\tilde{C} = \tilde{C}(n,p) > 1$ denotes the constant from the interpolation Lemma \ref{interpolation lemma, sobolev version}. One can apply Lemma \ref{interpolation lemma, sobolev version} with $\varepsilon = (2C\tilde{C})^{-1}\rho^2 \rho^l$ and obtain
   \begin{align*}
||u\chi_l||_{C^{1+2\gamma,\gamma}_{x,t}(B_{R_{l+1}} \times I)} &\leq \rho^{-l} 4C\tilde{C}\left( ||f||_{L^p_{x,t}(\Omega \times I)}   + \varepsilon ||D^2_x (u\chi_{l+1})||_{L^p_{x,t}(\mathbb R^n \times I)}  + \varepsilon^{-1} ||u \chi_{l+1}||_{L^p_{x,t}(\mathbb R^n \times I)}  \right) \\
    &\leq \rho^{-l} 4C \tilde{C} \left( ||f||_{L^p_{x,t}(\Omega \times I)}  + \varepsilon ||D^2_x (u\chi_{l+1})||_{L^p_{x,t}(B_{R_{l+2}} \times I)} + \varepsilon^{-1} ||u||_{L^{p}_{x,t}(\Omega \times I)} \right) \\
    &\leq \rho^{-2l} \underbrace{(4C \tilde{C})^2}_{=: \hat{C}} \left( ||f||_{L^p_{x,t}(\Omega \times I)}  + ||u ||_{L^p_{x,t}(\Omega \times I)} \right) + \frac{\rho^2}{2}|| D^2_x (u\chi_{l+1})||_{L^p_{x,t}(B_{R_{l+2}} \times I)} 
  \end{align*}
  Starting with $l = 1$ and iterating, we obtain
\begin{align}
||D^2_x(u\chi_1)||_{C^{1+2\gamma,\gamma}_{x,t}(B_{R_{2}} \times I)} &\leq  \rho^{-2}\hat{C} \left( ||f||_{L^p_{x,t}(\Omega \times I)}  + ||u ||_{L^p_{x,t}(\Omega \times I)} \right)  + \frac{\rho^2}{2}|| D^2_x (u\chi_{2})||_{L^p_{x,t}(B_{R_{3}} \times I)} \notag \\
&\leq \left( \rho^{-2}\hat{C} + \frac{\rho^2}{2} \cdot \rho^{-4}\hat{C} \right) \left( ||f||_{L^p_{x,t}(\Omega \times I)}  + ||u||_{L^p_{x,t}(\Omega \times I)} \right) \notag \\
&\phantom{\leq}+ \left(\frac{\rho^2}{2}\right)^2 || D^2_x (u\chi_{3})||_{L^p_{x,t}(B_{R_{4}} \times I)} \notag \\
&\leq \rho^{-2} \hat{C} \left( 1 + \frac{1}{2} + ... + \frac{1}{2^{k-2}} \right) \left( ||f||_{L^p_{x,t}(\Omega \times I)}  + ||u ||_{L^p_{x,t}(\Omega \times I)} \right) \notag \\
&\phantom{\leq}+ \left(\frac{\rho^2}{2}\right)^{k-1} || D^2_x (u\chi_{k})||_{L^p_{x,t}(B_{R_{k+1}} \times I)} \label{parabolic estimate proof : final inequality}
  \end{align}
  for any $k > 1$. Since $\rho < 1$ and
\begin{align*}
      || D^2_x (u\chi_{k})||_{L^p_{x,t}(B_{R_{k+1}} \times I)}  &\leq \rho^{-k} C  ||(|u| + |\nabla_xu| + |D^2_xu|) ||_{L^p_{x,t}(\Omega \times I)} 
\end{align*}
we can let $k \to +\infty$ in (\ref{parabolic estimate proof : final inequality}). Observing that $u = u \chi_1$ on $B_{1/2}(0) \times I \subset B_{R_2} \times I$, this concludes the proof when $D^2_x u \in L^p(\Omega \times I)$.
If one omits the hypothesis $D^2_x u \in L^p(\Omega \times I)$, it suffices to use a mollification. Write $I = (a,b)$ with $-\infty \leq a < b \leq +\infty$, $\Omega_{\varepsilon} = \{x \in \Omega: \dist(x,\partial \Omega) > \varepsilon\}$ and let $u_{\varepsilon}(x,t) = \phi_{\varepsilon} *_{x,t} u(x,t)$, $(x,t) \in \overline{\Omega_{\varepsilon}} \times [a+\varepsilon, b-\varepsilon]$ and $\phi_{\varepsilon} = \varepsilon^{-n-1} \phi(\varepsilon^{-1}x, \varepsilon^{-1}t)$ is an approximate identity defined by some $\phi \in C^{\infty}_c(B_{1/2}(0)) \subset C^{\infty}_c(\mathbb R^{n+1})$, $\phi \geq 0$, $\int_{\mathbb R^{n+1}} \phi dxdt = 1$. 

One verifies that $u_{\varepsilon} \in C^{\infty}(\overline{\Omega_{\varepsilon}} \times [a+\varepsilon, b-\varepsilon])$ and
$$
\nabla_x u_{\varepsilon} = \phi_{\varepsilon} * \nabla_x u, \quad (\partial_t - \Delta)  u_{\varepsilon} = \phi_{\varepsilon} * (\partial_t - \Delta)  u, \quad (x,t) \in \Omega_{\varepsilon} \times (a+\varepsilon, b-\varepsilon),
$$
as well as $u_{\varepsilon} \to u$, $\nabla_x u_{\varepsilon} \to \nabla_x u$, $ (\partial_t - \Delta) u_{\varepsilon} \to  (\partial_t - \Delta) u$ in $L^p_{loc}(\Omega \times I)$. Moreover, they are uniformly bounded in $L^p(\overline{\Omega_{\varepsilon}} \times [a+\varepsilon, b-\varepsilon])$. 

For each $\varepsilon$, $D^2_x u_{\varepsilon} = D^2_x \phi_{\varepsilon} * u \in L^p(\Omega_{\varepsilon} \times (a+\varepsilon,b-\varepsilon))$ as a convolution of a $L^{\infty}$ and $L^p$ function. Let $B_{\delta}$ be a smaller ball $B_{\delta}(0) \subset B_1(0)$. Then one can apply the theorem to $u_{\varepsilon}$ on $B_{\delta}$, leading to
\begin{align*}
        ||u_{\varepsilon}||_{C^{1+2\gamma,\gamma}_{x,t}(B_{\delta/2} \times  (a+\varepsilon,b-\varepsilon))}
&\leq C \left( 1 + |I|^{-1-2\gamma-\frac{n+2}{p}}\right) \\
&\phantom{\leq}\cdot\left( ||u_{\varepsilon} ||_{L^p(\Omega_{\varepsilon} \times  (a+\varepsilon,b-\varepsilon))} + ||(\partial_t - \Delta)u_{\varepsilon}||_{L^p(\Omega_{\varepsilon} \times  (a+\varepsilon,b-\varepsilon)} \right) \\
&\leq C \left( 1 + |I|^{-1-2\gamma-\frac{n+2}{p}}\right) \left( ||u ||_{L^p(\Omega \times I)} + ||(\partial_t - \Delta)u||_{L^p(\Omega \times I)} \right),
\end{align*}
where we used Young Convolution's inequality once again. As $u_{\varepsilon} \to u, \nabla_x u_{\varepsilon} \to \nabla_x u$ almost everywhere on $\Omega \times I$ along a sequence $\varepsilon \to 0^+$, one gets the same estimate for $||u_{\varepsilon}||_{C^{1+2\gamma,\gamma}_{x,t}(K_{\delta/2} \times  (a,b))}$. It then suffices to take $\delta \to 1$.

In the half-space case, it suffices to replace $\Omega$ by $B_r(x_0)$ and consider the odd reflection $\tilde{u}$ of $u$ on this ball instead of $u$ defined only on the half-ball (see Proposition \ref{prop:extension by reflection}). We also note that the Laplacian is rotation invariant.
\end{proof}

The method of freezing coefficients allows to deduce $L^p$ estimates for more general parabolic operators.

\begin{lemma}[Global Parabolic Sobolev Embedding for compactly supported solutions of parabolic equations]\label{lemma:global parabolic estimates for parabolic operator, compact support}
    Let $n+2 < p < +\infty$ and
    $$
    0 < \gamma < \frac{(p-n-1)(p-n-2)}{p(2p-n-2)}.
    $$
    Consider a parabolic operator
    $$
    L = c(x,t) + \sum_{i=1}^n b^i(x,t) D_{x_i} + \sum_{i,j=1}^n a^{i,j}(x) D_{x_i}D_{x_j}
    $$
    defined on $\overline{\Omega}$, $\Omega = B_{\kappa}(x_1) \times (0,T)$, satisfying
    \begin{align*}
    ||a,b,c||_{2\gamma,\gamma;\overline{\Omega}} &\leq \Lambda, \quad \Lambda > 0, \\
        \Lambda^{-1}|\xi|^2 \leq \sum_{i,j=1}^n a^{i,j}(x,t)\xi_i\xi_j &\leq \Lambda |\xi|^2, \quad (x,t) \in \overline{\Omega}, \xi \in \mathbb R^n,
    \end{align*}
    and where $a^{i,j}(x)$ does not depend on time.
    
    There exists $\kappa_0 = \kappa_0(n,\alpha,\Lambda) > 0$ and $C = C(n,p,\alpha,\Lambda) > 0$ such that if $\kappa \leq \kappa_0$, for any $u \in L^{\infty}((0,T);L^p_x(\mathbb R^n))$ with weak derivatives $\nabla_x u, \partial_t u, D^2_xu \in L^p(B_{\kappa}(x_1) \times (0,T))$, that is compactly supported in the sense that $\overline{\supp(u(\cdot,t))} \subset B_{\kappa}(x_1)$ for all $t \in (0,T)$ and which is weakly continuous at zero in the sense that
    $$
    \lim \limits_{t \to 0^+} \langle u(\cdot,t), \phi \rangle = 0 \quad \forall \phi \in C^{\infty}_c(\mathbb R^n),
    $$
    then one has
    $$
    ||u||_{C^{1+2\gamma,\gamma}_{x,t}(\mathbb R^{n}\times  (-\infty,0] )} \leq C\left( ||(\partial_t - L)u||_{L^p(\mathbb R^n \times (-\infty,0])} + ||u||_{L^p(\mathbb R^n \times (-\infty,0])} \right).
$$
\end{lemma}

\begin{proof}
    It follows from Lemma \ref{lemma:global l^p estimates for parabolic operator, compact support} that $D^2_xu \in L^p(\mathbb R^n \times (0,T))$, hence $(\partial_t - \Delta) u = (\partial_t - L)u + (L-\Delta)u \in L^p(\mathbb R^n \times (0,T))$ as well. Then
\begin{align*}
        ||u||_{C^{1+2\gamma,\gamma}_{x,t}(\mathbb R^{n}\times  (-\infty,0] )} &\leq C\left( ||(\partial_t - \Delta)u||_{L^p(\mathbb R^n \times (-\infty,0])} + ||u||_{L^p(\mathbb R^n \times (-\infty,0])} \right) \\
        &\leq  C\left( ||(\partial_t - L)u||_{L^p(\mathbb R^n \times (-\infty,0])} + ||D^2_x u||_{L^p(\mathbb R^n \times (-\infty,0])} + ||u||_{L^p} \right) \\
        &\leq C\left( ||(\partial_t - L)u||_{L^p(\mathbb R^n \times (-\infty,0])} + ||u||_{L^p(\mathbb R^n \times (-\infty,0])} \right) 
\end{align*}
    using Lemma \ref{lemma:global l^p estimates for parabolic operator, compact support} and Theorem \ref{thm:parabolic sobolev embedding}.
\end{proof}

\begin{corollary}[Parabolic Sobolev Embedding for parabolic equations with a priori estimates] \label{cor:interior parabolic estimates for parabolic eq with a priori estimates}
   Let $0 < T < +\infty$, $\Omega \subset \mathbb R^n$ be open, $n+2 < p < +\infty$ and
    $$
    0 < \gamma < \frac{(p-n-1)(p-n-2)}{p(2p-n-2)}.
    $$
    Let $L$ be a parabolic operator on $(x,t) \in \overline{\Omega} \times [0,T]$ as in Lemma \ref{lemma:global parabolic estimates for parabolic operator, compact support} with coefficients $a^{i,j}(x)$ independent of $t$. Assume that $u \in L^{\infty}((0,T);L^p_x(\mathbb R^n))$ has weak derivatives $\nabla_x u, \partial_t u, D^2_xu \in L^p(B_{\kappa}(x_0) \times (0,T))$ and is weakly continuous at zero in the sense that
    $$
    \lim \limits_{t \to 0^+} \langle u(\cdot,t), \phi \rangle = 0 \quad \forall \phi \in C^{\infty}_c(\mathbb R^n).
    $$
    For any ball $\overline{B_r(x_0)} \subset \Omega$, one has $u \in C^{1+2\gamma,\gamma}_{x,t}(B_{r/2}(x_0) \times I)$. Moreover, there exists $C = C(n,p,\gamma,\Lambda)$ for which
    \begin{align*}
            ||u||_{C^{1+2\gamma,\gamma}_{x,t}(B_{r/2}(x_0) \times I)} &\leq C \left( 1 + (|I|/r^2)^{-1-2\gamma-\frac{n+2}{p}} + r^{-1-2\gamma-\frac{n+2}{p}} \right) \\
            &\phantom{\leq}\cdot \left( ||u||_{L^p(\Omega \times I)} + ||(\partial_t - L) u||_{L^p(\Omega \times I)} \right).
        \end{align*}
        The same result holds if one replaces $B_{r}(x_0)$, resp. $B_{r/2}(x_0)$, by some half-ball $B_r^+(x_0) = B_r(x_0) \cap \{x_n > (x_0)_n\}$, where the flat part $T_r^+(x_0) = B_r(x_0) \cap \{x_n = (x_0)_n\} \subset \overline{\Omega}$ can touch the boundary while $\overline{B_r^+(x_0)} \setminus T_r^+(x_0) \subset \Omega$. Moreover, one assumes $u \in C^{2,1}(B_r^+(x_0) \times I)$, $\lim \limits_{\varepsilon \to 0^+} u(\tilde{x},(x_0)_n+\varepsilon,t) \to 0$ in $L^p_{loc}(B_{r}(\tilde{x}_0) \times I)$, where $x_0 = (\tilde{x}_0, (x_0)_n)$, $\lim \limits_{\varepsilon \to 0^+}  \nabla_x u(\tilde{x},(x_0)_n+\varepsilon,t)$ converges in $L^p_{loc}(B_{r}(\tilde{x}_0) \times I)$.
\end{corollary}

\begin{proof}
    After translation and parabolic rescaling, one can assume that $x_0 = 0$ and $r = 1$. Then, for any small ball $B_{\kappa}(x_1) \subset B_1(0)$, one can proceed as in the proof of Corollary \ref{cor: local parabolic sobolev embedding} using cut-off functions to get an estimate on $B_{\kappa/2}(x_1)$, but one needs to use Lemma \ref{lemma:global parabolic estimates for parabolic operator, compact support} instead of Theorem \ref{thm:parabolic sobolev embedding}. Then it suffices to cover $B_{1/2}(0)$ with finitely many such estimates.

    In the half-ball case, restrict $\Omega$ to $B_r^+(x_0)$, extend $u$ by odd reflection (See Proposition \ref{prop:extension by reflection}), extend the coefficients $c$, $b^i$, $a^{i,j}$, $a^{n,n}$, $i,j \neq n$, by odd reflection and the coefficients $b^n, a^{i,n}, a^{n,i}$, $i \neq n$, by even reflection and apply the theorem with the ball $B_{r'}(x_0)$ for $r' \to r$.
\end{proof}

   Finally, one obtains a boundary parabolic estimate for bounded domains with smooth boundary by flattening and covering the boundary.
   
\begin{theorem}[Boundary Parabolic Sobolev Embedding on smooth domains]\label{thm:boundary parabolic sobolev embedding}
 Let $0 < T < +\infty$, $\gamma \in (0,1/2)$ and $\Omega \subset \mathbb R^n$ be open, bounded, with $C^{\infty}$-boundary.

    Assume that $u \in C^0([0,T);C^1_0(\overline{\Omega}))$ is a distribution solution to the Dirichlet problem
\begin{align*}
    u_t - \Delta u &= f(x,t), \quad (x,t) \in \Omega \times (0,T), \\
    u &= 0, \quad (x,t) \in \partial \Omega \times (0,T) \cap \Omega \times \{0\},
\end{align*}
with $f \in L^{\infty}(\Omega \times (0,T))$.

There exists $C =  C(n,\gamma, \Omega)$ and $\varepsilon = \varepsilon(n,\gamma,\Omega)$ such that the estimate
    \begin{align*}
         ||u||_{1+2\gamma,\gamma;\Omega_{\varepsilon}  \times (0,T)} 
         &\leq C \bigg( ||f||_{L^{\infty}_{x,t}(\Omega \times (0,T))} + ||u||_{L^{\infty}_{x,t}(\Omega \times (0,T)) } \bigg)
     \end{align*}
     holds on the $\varepsilon$-neighborhood
     $$
    \Omega_{\varepsilon} = \{x \in \Omega: \dist(x,\partial \Omega) < \varepsilon\}
    $$
    of the boundary. 
\end{theorem}

\begin{proof}
One first notes that a $C^0([0,T) \times \overline{\Omega})$ solution is unique thanks to hypoellipticity and the Maximum Principle (Theorem \ref{maximum principle}). The solution $u$ is given by integration with a kernel $k_t(x,y)$ given by
 $$
k_t(x,y) = \sum_{k=1}^{\infty} e^{-t\lambda_k}w_k(x)\overline{w_k}(y),
 $$
  where, as shown in \cite[Chapter 6.5.1]{evans10}, the $(w_k)_{k \geq 1}$ form a countable orthonormal basis of $L^2(\Omega)$ made of eigenfunctions $w_k(x) \in C^{\infty}(\overline{\Omega}), w_k(\partial \Omega) = 0$, $-\Delta w_k = \lambda_k w_k$ on $\Omega$, $\lambda_k > 0$. The sum converges in $L^2(\Omega^2)$ for all $t > 0$ and the kernel is $C^{\infty}((0,+\infty) \times \Omega^2)$ by hypoellipticity because $k_t(x/2,y/2)$ is a distribution solution to the heat equation
 $$
    (\partial_t - \Delta_{x,y}) u(x,y,t) = 0, \quad (x,y) \in 2\Omega \times 2\Omega, t > 0.
 $$ 
 Actually, the sum converges uniformly as 
 $$
||w_k||_{C^0(\Omega)} \lesssim ||w_k||_{H^{m+2}(\Omega)} \lesssim C \left( \lambda_k ||w_k||_{H^m(\Omega)} + ||w_k||_{L^2(\Omega)}  \right) \lesssim C(1 + \lambda_k)^N
 $$
 by elliptic regularity (\cite[Chapter 6.3.2]{evans10}), for $m, N$ large enough which are fixed and independent of $k$. One deduces that $k_t$ extends continuously as a function $k_t(x,y) \in C^0((0,+\infty) \times \overline{\Omega}^2)$ with $k_t(x,y) = 0$ when $x \in \partial \Omega$ or $y \in \partial \Omega$. For fixed $y_0 \in \Omega$, $u(x,t) = G(x-y_0,t) - k_t(x,y_0)$ solves the heat equation on $(0,+\infty)$ with $u(\partial \Omega,t) > 0$. The Maximum Principle shows that $k_t(x,y_0) < G(x-y_0,t)$ on $(0,+\infty) \times \Omega$. 
 
 Hence, the corresponding heat semi-group maps $L^{\infty}(\Omega)$ to $C^0(\overline{\Omega})$ for $t > 0$ with 
     \begin{align*}
||e^{t \Delta} g||_{L^{\infty}(\Omega)} &\leq C ||g||_{L^{\infty}(\Omega)} \quad \forall t > 0.
\end{align*}
Thus the solution operator maps $f \in L^{\infty}(\Omega \times (0,T))$ to $u \in C^0(\overline{\Omega} \times [0,T])$ by Dominated Convergence, and more generally $L^p(\Omega \times (0,T))$ to $L^p(\Omega \times (0,T))$ by comparison with the Gaussian Heat Kernel, with
\begin{equation}
    ||u||_{L^{p}(\Omega \times (0,T))} \leq C T ||f||_{L^{p}(\Omega \times (0,T))}, \quad p \in [1,+\infty]. \label{eq: L infty global estimate}
\end{equation}

Next, observe that the theorem is a direct consequence of Corollary \ref{cor:interior parabolic estimates for parabolic eq with a priori estimates} when $u \in C^{2,1}_{x,t}(\Omega \times (0,T))$, $\partial_t u, D^2_x u \in L^{p}(\Omega \times (0,T))$, for $p = p(\gamma)$ large enough by locally straightening the boundary. Otherwise, let $f_n \to f$ in $L^{p}(\mathbb R^n \times (0,T))$, $(f_n) \in C_c^{\infty}(\Omega \times (0,T))$. The corresponding sequence of solutions $u_n$ to the Dirichlet problem with zero initial data has regularity $u_n \in C^{\infty}(\overline{\Omega} \times [0,T])$ (see \cite[Chapter 7.1.3, Theorem 7]{evans10}). Moreover, the sequence $u_n$ converges to $u$ in $L^{p}(\Omega \times (0,T))$ thanks to (\ref{eq: L infty global estimate}).

One obtains the estimate for each $u_n$ and $u_n - u_m$, which concludes the proof upon passing to the limit.
\end{proof}

\section{Schauder Estimates for parabolic equations}
In this appendix, we present the so-called Schauder estimates (both interior and boundary) for linear parabolic equations. Those estimates provide an upper bound on the Hölder semi-norm $[D^2_x u]_{\alpha,\alpha/2}$ in terms of the semi-norm $[f]_{\alpha,\alpha/2}$ of the forcing term and $||u||_{L^{\infty}_{x,t}}$. If the parabolic operator has constant coefficient, the proof goes as follows: use an argument by contradiction and the derivative estimates to prove a global estimate, which we then localize (Corollary \ref{interior schauder estimates, ver.1}). Then, freeze the coefficients to deduce the estimate for parabolic operators with varying coefficients, so that one can also obtain a boundary Schauder estimates (Theorem \ref{boundary schauder smooth boundary}). We note that the estimates from \cite{Simon1997} or \cite{krylov} require some a priori smoothness on the solution $u$. For the heat operator and thanks to the Parabolic Sobolev Embedding proved in the previous appendix, we show that this smoothness assumption can be weakened.

    In the following, $L_0$ denotes a parabolic operator
    \begin{equation}
            L_0 =\sum_{i,j=1}^n a^{i,j}D_{x_i}D_{x_j} \label{eq: L_0 constant coefficient parabolic}
    \end{equation}
     with constant coefficients satisfying $a^{i,j} = a^{j,i}$ and 
    \begin{align*}
        \Lambda^{-1}|\xi|^2 \leq \sum_{i,j=1}^n a^{i,j}\xi_i\xi_j &\leq \Lambda |\xi|^2, \quad \xi \in \mathbb R^n,
    \end{align*}
for some $\Lambda > 0$ (i.e., $A = (a^{i,j})$ is a symmetric, positive definite matrix).

One should note that if $u \in C^{2,1}_{x,t}(\mathbb R^n \times (-\infty,0])$ is a solution to the free parabolic equation
$$
    u_t - L_0 u = 0, \quad (x,t) \in \mathbb R^n \times (-\infty,0],
$$
and if $P$ is a matrix for which $P^T (a^{i,j}) P = I_n$ (the matrix $P$ is obtained through a orthogonal diagonalization $Q^T A Q = Diag(\lambda_1, ..., \lambda_n)$ and then $P = QDiag(\lambda_1^{-1/2}, ...,\lambda_n^{-1/2})$), then $v(x,t) =  u(Px,t)$ solves
$$
    v_t - \Delta v = 0, \quad (x,t) \in \mathbb R^n \times (-\infty,0].
$$
Similarly, if $u \in C^{2,1}_{x,t}(\mathbb R^n \cap \{\langle d, x \rangle \geq 0\} \times (-\infty,0])$ is a solution to the free parabolic equation
\begin{align*}
        u_t - L_0 u &= 0, \quad (x,t) \in \mathbb R^n \cap \{\langle d,x \rangle > 0\} \times (-\infty,0], \\
        u(x,t) &= 0, \quad (x,t) \in \mathbb R^{n} \cap \{ \langle d, x \rangle = 0\} \times (-\infty,0],
\end{align*}
 and if $P$ satisfies $P^T (a^{i,j}) P = I_n$, $R = (b_1,...,b_{n-1},P^T d/||P^Td||) \in \mathbb R^{n \times n}$ is an orthogonal matrix, then $v(x,t) =  u(PRx,t)$ solves
\begin{align*}
        v_t - \Delta v &= 0, \quad (x,t) \in \mathbb R^n \cap \{x_n > 0\} \times (-\infty,0], \\
        v(x,t) &= 0, \quad (x,t) \in \mathbb R^{n} \cap \{ x_n  = 0\} \times (-\infty,0].
\end{align*}

\begin{proposition}[Derivative estimates for the free Heat Equation] \label{derivative estimates free heat eq}
    Let $U \subset \mathbb R^n$ be open, $t_0 \in \mathbb R$ and $T > 0$. Let $u \in C^{2,1}_{x,t}(U \times (t_0,t_0+T])$ solve the free heat equation
    $$
    u_t - \Delta u = 0, \quad (x,t) \in U \times (t_0,t_0+T] =: U_T.
    $$
    Then $u \in C^{\infty}(U_T)$. Moreover, for any $k,l \in \mathbb N_{\geq 0}$, there exists a constant $C_{k,l,n} > 0$ (independent of $U$, $t_0$, $T$) such that for any closed cylinder $C(x,t;r) \subset U_T$, one has
    $$
    \sup_{C(x,t;r/2)}|D_x^k D_t^l u| \leq \frac{C_{k,l,n}}{r^{k+2l+n+2}}||u||_{L^1(C(x,t;r))}. 
    $$
\end{proposition}

\begin{proof}
    See \cite{evans10}, Chapter 2.3, Theorems 8 and 9.
\end{proof}

\begin{proposition}[Liouville's Theorem for backwards in time Heat Equation] \label{Liouville theorem}
    Assume that $u \in C^{2,1}_{x,t}(\mathbb R^n \times (-\infty,0])$ is a solution to the free parabolic equation
$$
    u_t - L_0 u = 0, \quad (x,t) \in \mathbb R^n \times (-\infty,0], 
$$
where $L_0$ is as in (\ref{eq: L_0 constant coefficient parabolic}).

Assume that $||D^2_xu||_{L^{\infty}_{x,t}(\mathbb R^n \times (-\infty,0])} < +\infty$ or $[D^2_xu]_{\alpha,\alpha/2; \mathbb R^n \times (-\infty,0]} < +\infty$. Then 
$$
u(x,t) = u(0,0) + x^T v + \frac{1}{2} x^T H x + \Tr(H) t
$$
is a constant-coefficient polynomial. 

The proposition still holds if $u \in C^{2,1}(\mathbb R^n \cap \{\langle d,x \rangle \geq 0\} \times (-\infty,0])$, $|d| = 1$, solves the half-space Dirichlet problem
\begin{align*}
        u_t - L_0 u &= 0, \quad (x,t) \in \mathbb R^n \cap \{\langle d,x \rangle > 0\} \times (-\infty,0] \\
        u(x,t) &= 0, \quad (x,t) \in \mathbb R^{n} \cap \{ \langle d, x \rangle = 0\} \times (-\infty,0]
\end{align*}
with $||D^2_xu||_{L^{\infty}_{x,t}(\mathbb R^n \cap \{\langle d,x \rangle > 0\} \times (-\infty,0])} < +\infty$ or $[D^2_xu]_{\alpha,\alpha/2; \mathbb R^n \cap \{\langle d,x \rangle > 0\} \times (-\infty,0]} < +\infty$. 
\end{proposition}

\begin{remark}
    In the half-space case, smoothness is assumed up to the boundary.
\end{remark}

\begin{proof}
Without loss of generality, up to applying a linear transformation, $L_0 = \Delta$ and $d = e_n = (0,...,0,1)$.

We first deal with the full space case. Observe that
\begin{align*}
    u(x,t) &= u(x,0) + \partial_t u(x,c(x,t))t, \quad |c(x,t)| \leq t \\
    &= u(x,0) + \partial_t u(0,0) t + [\partial_t u(x,c(x,t))-\partial_t u(0,0)]t\\
    &= u(x,0) + \partial_t u(0,0) t + [\Delta u(x,c(x,t))-\Delta u(0,0)]t\\
    &= u(x,0) + \mathcal{O}(|t| + |x|^{\alpha}|t| + |t|^{1+\alpha/2}) \\
    &= u(0,0) + \nabla_x u(0,0) \cdot x + \mathcal{O}(|t| + |x|^{\alpha}|t| + |t|^{1+\alpha/2} + |x|^2 + |x|^{2+\alpha} )
\end{align*}
if $D^2_xu$ is globally Hölder continuous and
$$
 u(x,t) = u(0,0) + \nabla_x u(0,0) \cdot x + \mathcal{O}(|t| + |x|^2)
$$
if it is globally bounded.

    Apply the derivative estimates with a cylinder $C(0,0;r) \subset \mathbb R^n \times (-\infty,0]$, i.e.
    $$
    \sup_{C(0,0;r/2)}|D_x^k D_t^l u| \leq \frac{C_{k,l,n}}{r^{k+2l+n+2}}||u||_{L^1(C(0,0;r))} \lesssim \frac{C_{k,l,n}}{r^{k+2l+n+2}} \cdot r^{n+4+\alpha}
    $$
    for any $r > 1$ and $k,l \in \mathbb N_{\geq 0}$. Letting $r \to +\infty$ shows that any derivative $D_x^k D_t^l u$, $k+2l > 2$, is identically zero on any $\mathbb R^n \times (-\infty,0]$. Hence,
    $$
    u(x,t) = u(x,0) + \partial_t u(x,0) t
    $$
    because $\partial_{tt} u = 0$. Then $D_x \partial_t u(x,0) = 0$ implies that $\partial_t u(x,0) = \partial_t u(0,0)$ is constant. Finally, $D^2_x u(x,0) = D^2_x (u(x,t) - \partial_t u(0,0) t) = 0$ implies that $u(x,0)$ is a polynomial of degree at most $2$.

In the half-space case, it suffices to extend $u$ by reflection. Let 
    $$
\tilde{u}(\tilde{x},x_n,t) = \begin{cases} u(\tilde{x},x_n,t), \quad &x_n > 0 \\
- u(\tilde{x},- x_n,t), \quad &x_n < 0 \\
\end{cases}
$$
Then $\tilde{u} \in C^{2,1}(\mathbb R^n \times (-\infty,0])$ by Proposition \ref{prop:extension by reflection}.  Moreover, $\tilde{u}$ solves the free heat equation on $\mathbb R^n$ in the classical sense with $||D^2_x \tilde{u}||_{L^{\infty}_{x,t}(\mathbb R^n \times (-\infty,0])} < +\infty$ or $[D^2_x\tilde{u}]_{\alpha,\alpha/2; \mathbb R^n \times (-\infty,0]} < +\infty$, which concludes the proof.
\end{proof}

\begin{lemma}[Interpolation Inequalities] \label{interpolation lemma}
    Let $\alpha \in (0,1)$ and $u \in C^{2+\alpha,1+\alpha/2}(\Omega \times I;\mathbb R)$, where $\Omega \subset \mathbb R^{n}$ is an open set which satisfies the interior ball condition (Definition \ref{def:smooth boundary}) with some ball of radius $2\rho$. Let $I$ be an open interval and assume $\rho^2 \leq |I|$. There exists $C = C(n,\alpha)$ for which
    \begin{align}
        [u]_{\alpha,\alpha/2} &\leq C \left[ \rho \left( [D^2_xu]_{\alpha} + [\partial_t u]_{\alpha/2} \right) + \rho^{-\frac{\alpha}{2}} ||u||_{L^{\infty}_{x,t}} \right] \label{Holder norm of u, interpolation}\\
        ||\nabla_x u||_{L^{\infty}_{x,t}} &\leq C \left[  \rho \left( [D^2_xu]_{\alpha} + [\partial_t u]_{\alpha/2} \right) + \rho^{-\frac{1}{1+\alpha}} ||u||_{L^{\infty}_{x,t}} \right] \label{Sup norm of nabla_x u, interpolation}  \\ 
        [\nabla_x u]_{\alpha,\alpha/2} &\leq C \left[ \rho \left( [D^2_xu]_{\alpha} + [\partial_t u]_{\alpha/2} \right) +  \rho^{-(1+\alpha)} ||u||_{L^{\infty}_{x,t}} \right] \label{Holder norm of nabla_x u, interpolation}  \\
        ||D^2_x u||_{L^{\infty}_{x,t}} &\leq C \left[  \rho \left( [D^2_xu]_{\alpha} + [\partial_t u]_{\alpha/2} \right) + \rho^{-\frac{2}{\alpha}} ||u||_{L^{\infty}_{x,t}} \right] \label{Sup norm of D^2_x u, interpolation} \\
        ||\partial_t u||_{L^{\infty}_{x,t}} &\leq C \left[  \rho \left( [D^2_xu]_{\alpha} + [\partial_t u]_{\alpha/2} \right) + \rho^{-\frac{2}{\alpha}} ||u||_{L^{\infty}_{x,t}} \right], \label{Sup norm of partial_t u, interpolation} 
    \end{align}
     where the (semi-)norms are taken on the whole $\Omega \times I$. In particular, the inequalities still hold if one replaces $\rho$ by a smaller $\rho' \leq \rho$ and they hold for any $\rho > 0$ if $I$ is unbounded together with either $\Omega = \mathbb R^n$ or $\Omega = \mathbb R^n \cap \{x_n > 0\}$.
\end{lemma}

\begin{proof}
 By replacing $\rho$ with $\rho - \varepsilon$, $\varepsilon \to 0$, one can assume that the balls one gets from the Definition \ref{def:smooth boundary} always have their closure inside $\Omega$ and $\rho^2 < |I|$. And after parabolic rescaling, one can assume without loss of generality that $\rho = 1$.

We start with (\ref{Sup norm of D^2_x u, interpolation}). Fix $i,j \in \{1,...n\}$ and $t \in I$. Assume $x \in \overline{B_{2}(x)} \subset \Omega$ (i.e. $x$ is at distance $> 2$ from the boundary). Let $\delta_k f = f(x+ e_k)-f(x)$. It follows from the mean-value theorem that
    $$
    \delta_i \delta_j u(x,t) = \delta_i \partial_j u(y,t) = \partial_i \partial_j u(z,t)
    $$
    for some $y \in \{x + s \rho e_j : s \in [0,1]\} \subset \overline{B_{1}(x)}$ and $z \in  \{y + s \rho e_i : s \in [0,1]\} \subset \overline{B_{2}(x)}$, i.e. for any such $x$, there is some $z \in \overline{B_{2}(x)}$ for which
    $$
    |\partial_i \partial_j u(z,t)| \leq 4  ||u||_{L^{\infty}_{x,t}}.    
    $$
   Then
 \begin{align}
     |\partial_i \partial_j u(x,t)| &\leq |\partial_i \partial_j u(x,t) - \partial_i \partial_j u(z,t)| + |\partial_i \partial_j u(z,t)| \notag \\
     &\leq |x-z|^{\alpha}[D^2_xu]_{\alpha,\alpha/2} + 4 ||u||_{L^{\infty}_{x,t}} \notag \\
     &\leq 4^{\alpha}[D^2_xu]_{\alpha,\alpha/2} + 4  ||u||_{L^{\infty}_{x,t}}. \label{eq: estimate away from boundary, D^2_x u infty case}
 \end{align}
    If $x$ is closer to the boundary, then $x \in \overline{B_{2}(z)} \subset \Omega$ for some $z$ further away from the boundary satisfying the above estimate (\ref{eq: estimate away from boundary, D^2_x u infty case}). In that case,
    \begin{align*}
     |\partial_i \partial_j u(x,t)| &\leq |\partial_i \partial_j u(x,t) - \partial_i \partial_j u(z,t)| + |\partial_i \partial_j u(z,t)| \\
     &\leq |x-z|^{\alpha}[D^2_xu]_{\alpha,\alpha/2} + 4^{\alpha}[D^2_xu]_{\alpha,\alpha/2} + 4 ||u||_{L^{\infty}_{x,t}} \\
     &\leq 2 \cdot 4^{\alpha}[D^2_xu]_{\alpha,\alpha/2} + 4 ||u||_{L^{\infty}_{x,t}}. 
 \end{align*}

 As for (\ref{Sup norm of partial_t u, interpolation}), note that $I$ satisfies the interior ball condition with a closed ball of radius $1$ since $1 < |I|$.  Fix $x \in \Omega$.  Let $\delta_t f = f(t + 1)-f(t)$. If $t \in \overline{B_{1}(t)} \subset I$, it follows from the mean-value theorem that
    $$
    \delta_t u(x,t) = \partial_t u(x,s)
    $$
    for some $s \in \{t + \lambda: \lambda \in [0,1]\} \subset \overline{B_{1}(t)}$, i.e. for any such $t$, there is some $s \in \overline{B_{1}(t)}$ for which
    $$
    |\partial_t u(x,t)| \leq |\partial_t u(x,t)- \partial_t u(x,s)| + |\partial_t u(x,s)| \leq 2^{\alpha/2} [\partial_t u]_{\alpha,\alpha/2} + 2  ||u||_{L^{\infty}_{x,t}}.
    $$
    If $t$ is closer to the boundary, then $x \in \overline{B_{1}(s)} \subset \Omega$ for some $s$ further away from the boundary satisfying the above estimate. In that case,
    \begin{align*}
     |\partial_t u(x,t)| &\leq 2 \cdot  2^{\alpha/2}[\partial_t u]_{\alpha,\alpha/2} + 2  ||u||_{L^{\infty}_{x,t}}
 \end{align*}
 as for $|\partial_i \partial_j u|$. One proceeds similarly to get (\ref{Sup norm of nabla_x u, interpolation}). Hence, (\ref{Holder norm of u, interpolation}) follows from the fact that $u(x,t)$ is Lipschitz ($\nabla_{x,t} u$ is bounded).

 Finally, it remains to prove (\ref{Holder norm of nabla_x u, interpolation}). The part
 $$
|\nabla_x u(x,t) - \nabla_x u(y,t)|
 $$
 is easily estimated as $D^2_x u \in L^{\infty}_{x,t}(\Omega \times I)$. The interesting part is the difference 
 $$
|\nabla_x u(x,t) - \nabla_x u(x,s)|
 $$
 when $d = |t-s|^{1/2} < 1$ (otherwise, bound it with $||\nabla_x u||_{L^{\infty}_{x,t}} \cdot |t-s|^{\frac{\alpha}{2}}$). In that case, $x \in \overline{B_{2}(y)} \subset \Omega$ for some $y \in \Omega$ and one can find $z \in \{x + t(y-x) : t \in [0,1]\}$ with $d  = |x-z|^{1/2}$ and $\overline{B_{d}(z)} \subset \Omega$. Then
 \begin{align*}
     |\nabla_x u(x,t) - \nabla_x u(x,s)| &\leq |\nabla_x u(x,t) - \nabla_x u(z,t) | \\
     &\phantom{\leq} + |\nabla_x u(z,t) - \nabla_x u(z,s)  | \\
     &\phantom{\leq} + |\nabla_x u(z,s) - \nabla_x u(x,s)| \\
     &\leq 2|x-z| \cdot ||D^2_x u||_{L^{\infty}_{x,t}} + |\nabla_x u(z,t) - \nabla_x u(z,s)  | \\
     &= 2|t-s|^{1/2} ||D^2_x u||_{L^{\infty}_{x,t}} + |\nabla_x u(z,t) - \nabla_x u(z,s)  | \\
     &\leq 2|t-s|^{\alpha/2} ||D^2_x u||_{L^{\infty}_{x,t}} + |\nabla_x u(z,t) - \nabla_x u(z,s)  |.
     \end{align*}
 As for the term involving $z$, fix $i \in \{1,...,n\}$. By construction of $z$, the segment $[z, z + de_i] \subset \Omega$ stays in $\Omega$. Then there exists $\lambda_1,\lambda_2 \in [0,1]$ for which
 \begin{align*}
     |\partial_{x_i} u(z,t) - \partial_{x_i}  u(z,s)| &\leq |\partial_{x_i}  u(z,t) - d^{-1}\left(u(z,t) - u(z+de_i,t) \right) | \\
     &\phantom{\leq} + |d^{-1}\left(u(z,t) - u(z+de_i,t) \right) - d^{-1}\left(u(z,s) - u(z+de_i,s) \right)  | \\
     &\phantom{\leq} + |d^{-1}\left(u(z,s) - u(z+de_i,s) \right)  - \partial_{x_i}  u(x,s)| \\
     &\leq |\partial_{x_i}  u(z,t) - \partial_{x_i} u(z + \lambda_1 de_i,s) | \\
     &\phantom{\leq} + d^{-1} |u(z,t)-u(z,s)| + d^{-1}|u(z+de_i,t)-u(z+de_i,s)| \\
     &\phantom{\leq} + |\partial_{x_i} u(z + \lambda_2 de_i,s)  - \partial_{x_i}  u(x,s)| \\
     &\leq 2d||D^2_x u||_{L^{\infty}_{x,t}} + 2d^{-1}|t-s| \cdot ||\partial_t u||_{L^{\infty}_{x,t}} \\
     &\leq 2|t-s|^{\alpha/2} \cdot ||D^2_x u||_{L^{\infty}_{x,t}} + 2|t-s|^{\alpha/2} \cdot||\partial_t u||_{L^{\infty}_{x,t}},
 \end{align*}
 which concludes the proof.

 \end{proof}
 
\begin{theorem}[Global Schauder Estimates for the Heat Equation]\label{thm:global schauder estimates}
Assume that $u \in C^{2+\alpha,1+\alpha/2}_{x,t}(\mathbb R^n \times (-\infty,0])$, $\alpha \in (0,1)$, is a solution to the linear heat equation
\begin{equation}
    u_t - L_0 u = f(x,t), \quad (x,t) \in \mathbb R^n \times (-\infty,0],
\end{equation}
where $f \in C^{\alpha,\alpha/2}_{x,t}(\mathbb R^n \times (-\infty,0])$ and $L_0$ is as in (\ref{eq: L_0 constant coefficient parabolic}).

There exists $C = C(n, \alpha, L_0)$ for which
$$
[D^2_x u]_{\alpha,\alpha/2;\mathbb R^n \times (-\infty,0]} + [\partial_t u]_{\alpha,\alpha/2;\mathbb R^n \times (-\infty,0]} \leq C[f]_{\alpha,\alpha/2;\mathbb R^n \times (-\infty,0]}.
$$
The result still holds if one replaces $\mathbb R^n$ by the half-space $\mathbb R^n \cap \{\langle d,x \rangle \geq 0\}$, $|d| = 1$, and one assumes a Dirichlet boundary condition $u = 0$ for all $(x,t) \in \mathbb R^n \cap \{\langle d,x \rangle = 0\} \times (-\infty,0]$.
\end{theorem}

\begin{remark}
    By $C^{2+\alpha,1+\alpha/2}_{x,t}(\mathbb R^n \times (-\infty,0])$, it is assumed a priori that $D^2_xu$ is globally bounded and Hölder-continuous. One should note that the backwards heat equation is ill-posed and there is no guarantee, given $f$, that such a solution exists.
    
    If the estimate from the theorem holds for $u$, it will not hold, in general, for $u + v$ where $v$ is a solution to the free heat equation. For example,  $v$ can be chosen to be any time-independent harmonic function which is not a polynomial.

    In the half-space case, Hölder-regularity is assumed up to the boundary.
\end{remark}

\begin{proof} 
We follow the argument by contradiction of Simon \cite{Simon1997}. Without loss of generality, up to applying a linear transformation, one can assume that $L_0 = \Delta$ and $d = e_n$.  The constant $C(n,\alpha,L_0)$ depends only on $\Lambda$ (which gives lower and upper bounds for the eigenvalues of the matric $(a^{i,j})$ which defines $L_0$).

   It suffices to prove $[D^2_x u]_{\alpha,\alpha/2}\leq C[f]_{\alpha,\alpha/2}$ and use $\partial_t u = \Delta u + f$ to bound $[\partial_t u]_{\alpha,\alpha/2}$. As for the inequality of norms, we can use the interpolation Lemma \ref{interpolation lemma} with $\rho = 1$, as well as the relationships $\partial_t u = \Delta u + f$.

   Without loss of generality, one can assume that $[f]_{\alpha, \alpha/2} \neq 0$, i.e. $f$ is not constant. If $f$ is constant, let $b \in \mathbb R^n \setminus \{0\}$ be perpendicular to $a$ and $c \in \mathbb R^n$ be perpendicular to $(b_1^2, b_2^2,...,b_n^2)$, $|c| = 1$. Then
   $$
   v(x,t) = \frac{-f}{2\sum_{i=1}^n c_i}\sum_{i=1}^n c_i x_i^2
   $$
   is a particular solution for which $[D^2_xv]_{\alpha, \alpha/2} = 0 = [f]_{\alpha, \alpha/2}$ and $v = 0$ on $\langle a,x\rangle = 0$. If $u$ is as in the theorem, then both $u$ and $v$ satisfy the conditions from Liouville's theorem \ref{Liouville theorem}. Hence, $u-v$ differ by a small-degree polynomial for which $D^2_x(u-v)$ is constant, which means $[D^2u]_{\alpha, \alpha/2} = [D^2v]_{\alpha, \alpha/2} = 0$.

    Assume the result does not hold. Then there exists a sequence of non-constant forcing terms  $f_k \in C^{\alpha,\alpha/2}_{x,t}(\mathbb R^n \times (-\infty,0])$ and solutions $u_k \in C^{2+\alpha,1+\alpha/2}_{x,t}(\mathbb R^n \times (-\infty,0])$ for which
    $$
    +\infty > [D^2_x u_k]_{\alpha,\alpha/2} > k[f_k]_{\alpha,\alpha/2} > 0.
    $$
    Normalizing $u_k \mapsto u_k/[f_k]_{\alpha,\alpha/2}$, we can assume without loss of generality that
    $$
    M_k := [D^2_x u_k]_{\alpha,\alpha/2} > k > 0, \quad 0 < [f_k]_{\alpha,\alpha/2} \leq 1.
    $$
    By definition of the Hölder semi-norm, for each $k \in \mathbb N$, there are pairs $p_k = (x_k, t_k), q_k = (y_k,s_k) \in \mathbb R^n \times (-\infty,0]$ satisfying
    $$
    \frac{|D^2_xu_k(p_k) - D^2_xu_k(q_k)|}{|p_k-q_k|} \geq \frac{1}{2}M_k,
    $$
    where $\lambda_k := |p_k - q_k| = |x_k - y_k| + |t_k-s_k|^{\frac{1}{2}}$ denotes the parabolic metric. The pidgeonhole principe implies that there are indices $i_k, j_k \in \{1,...,n\}$ for which
    \begin{equation}
    \frac{|D_{x_{i_k}}D_{x_{j_k}}u_k(p_k) - D_{x_{i_k}}D_{x_{j_k}}u_k(q_k)|}{|p_k-q_k|} \geq \frac{c_n}{2}M_k \label{Holder norm absurd hypothesis}
    \end{equation}
    for some $c_n > 0$ and passing to a subsequence of $k \in \mathbb N$, we can assume that $i_k = i, j_k = j$ for all $k \in \mathbb N$.

    We translate and rescale $u_k$ as follows by considering $\hat{u}_k(x,t) : \mathbb R^n \times (-\infty,0] \to \mathbb R$ defined by
 \begin{equation}
         \hat{u}_k(x,t) = M_k^{-1} \lambda_k^{-2-\alpha} u(x_k + \lambda_kx, t_k + \lambda_k t) \label{eq:hat u_k blow up simon}
 \end{equation}
    for which one has $[D^2_x \hat{u}_k]_{\alpha,\alpha/2} \leq M^{-1}_k[D^2_x u_k]_{\alpha,\alpha/2} = 1$. Next, we remove a small-order Taylor polynomial from $\hat{u}_k$ and define
\begin{equation}
        v_k(x,t) = \hat{u}_k(x,t) - \hat{u}_k(0,0) - \nabla_x \hat{u}_k(0,0) \cdot x - \frac{1}{2}x^T H_x\hat{u}_k(0,0)x - \partial_t \hat{u}_k(0,0)t, \label{eq:v_k blow up simon}
\end{equation}
   which solves
    \begin{equation}
\partial_t v_k - \Delta v_k = M_k^{-1}\lambda_k^{-\alpha}[f(x_k + \lambda_k x,t_k + \lambda_k t)-f(x_k,t_k)], \quad (x,t) \in \mathbb R^n \times (-\infty,0] \label{v_k heat equation}
    \end{equation}
    and satisfies
    \begin{equation}
        v_k(0,0) = \partial_tv_k(0,0) = D_x v_k(0,0) = D^2_x v_k(0,0) = 0. \label{v_k initial data}
    \end{equation}
    Moreover,
     \begin{equation}
     [D^2_x v_k]_{\alpha,\alpha/2} = [D^2_x \hat{u}_k]_{\alpha,\alpha/2} \leq 1 \label{D2v_k Holder}
     \end{equation}
     and
    \begin{equation}
    [\partial_t v_k]_{\alpha,\alpha/2} \leq [D^2_x v_k]_{\alpha,\alpha/2} +  M_k^{-1}\lambda_k^{-\alpha}[f(x_k + \lambda_k x,t_k + \lambda_k t)]_{\alpha,\alpha/2} \leq 1 + M_k^{-1}[f_k]_{\alpha,\alpha/2} \leq 2. \label{Dtv_k Holder}
    \end{equation}
    Combining (\ref{v_k initial data}) together with (\ref{D2v_k Holder}) and (\ref{Dtv_k Holder}), one finds that the families $(v_k)_k$, $(\partial_t v_k)$, $(D_x v_k)$, $(D^2_x v_k)$ are uniformly bounded and equicontinuous on any fixed cylinder $Q_R = B_R(0) \times (-R^2,0]$. By Arzelà-Ascoli, there exists $v_{\infty} \in C^{2+\alpha, 1+\alpha/2}( \mathbb R^n \times (-\infty,0] )$ for which $\lim \limits_{k \to +\infty} v_k = v_{\infty}$ in $C^{2,1}(Q_R)$-norm for any $R > 0$.
    
    Finally, we let $k \to +\infty$ in (\ref{v_k heat equation}). Since
    $$
M_k^{-1}\lambda_k^{-\alpha}|f(x_k + \lambda_k x,t_k + \lambda_k t)-f(x_k,t_k)| \leq M_k^{-1} [f_k]_{\alpha,\alpha/2} \leq M_k^{-1} \to 0,
    $$
    it follows that
    $$
    \partial_t v_{\infty} - \Delta v_{\infty} = 0, \quad (x,t) \in \mathbb R^n \times (-\infty,0]
    $$
    with globally Hölder-continuous $D^2_x v_{\infty}$. Hence, Liouville's theorem applies and 
    $$
v_{\infty}(x,t) = v_{\infty}(0,0) + \nabla_x v_{\infty}(0,0) \cdot x + \frac{1}{2}x^T H_x v_{\infty}(0,0)x + \partial_t v_{\infty}(0,0)t = 0
    $$
    is identically zero. This is a contradiction because 
    $$
    |D_{x_i}D_{x_j}v_{k}(0) - D_{x_i}D_{x_j}v_k(\lambda_k^{-1}(q_k-p_k))| = M_k^{-1}\lambda_k^{-\alpha}|D_{x_i}D_{x_j}u_{k}(p_k) - D_{x_i}D_{x_j}u_k(q_k)| \geq \frac{c_n}{2} > 0 \quad \forall k
    $$
    using (\ref{Holder norm absurd hypothesis}).
    
    Since $\lambda_k^{-1}(q_k-p_k)$ is a bounded (in parabolic norm) sequence, $\lambda_k^{-1}|q_k-p_k| = 1$, we can take the limit $k \to +\infty$ along a subsequence to find out that
    $$
|D_{x_i}D_{x_j}v_{\infty}(0) - D_{x_i}D_{x_j}v_{\infty}(p_{\infty})| \geq \frac{c_n}{2} > 0
    $$
    for some $p_{\infty} \in \mathbb R^n \times (-\infty,0]$, i.e. $D_{x_i}D_{x_j}v_{\infty}$ cannot be identically zero.

    In the half-space setting, we proceed in a similar fashion. We find sequences $(f_k)_k, (u_k)_k, (M_k)_k$, $(\lambda_k)_k$, $(p_k)_k = ((x_k,t_k))_k$, $(q_k)_k = ((y_k,s_k))_k$ satisfying (\ref{Holder norm absurd hypothesis}), where $x_k, y_k \in \mathbb R^n \cap \{x_n > 0\}$ and $-\infty < s_k \leq t_k \leq 0$. 
    
    If
    $$
     \sup_{k} \frac{\langle e_n, x_k \rangle}{\lambda_k} = \sup_{k} \frac{(x_k)_n}{\lambda_k} = +\infty
    $$
    we pass to a subsequence (still denoted with index $k$) for which $\lim \limits_{k \to +\infty} \langle a, x_k \rangle \lambda_k^{-1} = +\infty$. In that situation, we proceed as in the full space case. Define $\hat{u}_k$ and $v_k$ as in (\ref{eq:hat u_k blow up simon}) and (\ref{eq:v_k blow up simon}). Note that $\hat{u}_k$ and $v_k$ are only defined for $(x,t) \in \mathbb R^n \cap \{\langle a,x \rangle > -\langle a, x_k\rangle \lambda_k^{-1}\} \times (-\infty,0]$, but in the limit, it does not really matter so the rest of the argument works the same.

    If
    $$
    \sup_{k} \frac{(x_k)_n}{\lambda_k} < +\infty
    $$ 
    we instead let $\hat{u}_k(x,t) : \mathbb R^n 
\cap \{x_n > 0\} \times (-\infty,0] \to \mathbb R$ defined by
 \begin{equation*}
         \hat{u}_k(x,t) = M_k^{-1} \lambda_k^{-2-\alpha} u((\tilde{x}_k,0) + \lambda_kx, t_k + \lambda_k t), \quad x_k = (\tilde{x}_k, (x_k)_n)
 \end{equation*}
and  for which one has 
    $$[D^2_x \hat{u}_k]_{\alpha,\alpha/2} \leq M^{-1}_k[D^2_x u_k]_{\alpha,\alpha/2} = 1.$$
    Next, we remove a small-order Taylor polynomial from $\hat{u}_k$ and define
\begin{equation*}
        v_k(x,t) = \hat{u}_k(x,t) - \hat{u}_k(0,0) - \nabla_x \hat{u}_k(0,0) \cdot x - \frac{1}{2}x^T H_x\hat{u}_k(0,0)x - \partial_t \hat{u}_k(0,0)t
\end{equation*}
   which solves
    \begin{equation}
\partial_t v_k - \Delta v_k = M_k^{-1}\lambda_k^{-\alpha}[f((\tilde{x_k},0) + \lambda_k x,t_k + \lambda_k t)-f((\tilde{x_k},0),t_k)], \quad (x,t) \in \mathbb R^n \cap \{x_n > 0\} \times (-\infty,0] \label{v_k heat equation, half-space}
    \end{equation}
    and satisfies
    \begin{equation*}
        v_k(0,0) = \partial_tv_k(0,0) = D_x v_k(0,0) = D^2_x v_k(0,0) = 0.
    \end{equation*}
    Moreover,
     \begin{equation*}
     [D^2_x v_k]_{\alpha,\alpha/2} = [D^2_x \hat{u}_k]_{\alpha,\alpha/2} \leq 1 
     \end{equation*}
     and
    \begin{equation*}
    [\partial_t v_k]_{\alpha,\alpha/2} \leq [D^2_x v_k]_{\alpha,\alpha/2} +  M_k^{-1}\lambda_k^{-\alpha}[f((\tilde{x_k},0) + \lambda_k x,t_k + \lambda_k t)]_{\alpha,\alpha/2} \leq 1 + M_k^{-1}[f_k]_{\alpha,\alpha/2} \leq 2.
    \end{equation*}
    Combining those, one finds that the families $(v_k)_k, (\partial_t v_k), (D_x v_k), (D^2_x v_k)$ are uniformly bounded and equicontinuous on any fixed compact subset $K$ of $\mathbb R^n \cap \{x_n \geq 0\} \times (-\infty,0]$. By Arzelà-Ascoli, there exists $v_{\infty} \in C^{2+\alpha, 1+\alpha/2}( \mathbb R^n \times (-\infty,0] )$ for which $\lim \limits_{k \to +\infty} v_k = v_{\infty}$ in $C^{2,1}(K)$-norm for any compact subset $K$.
    
    Finally, we let $k \to +\infty$ in (\ref{v_k heat equation, half-space}). Since
    $$
M_k^{-1}\lambda_k^{-\alpha}|f((\tilde{x_k},0) + \lambda_k x,t_k + \lambda_k t)-f((\tilde{x_k},0),t_k)| \leq M_k^{-1} [f_k]_{\alpha,\alpha/2} \leq M_k^{-1} \to 0,
    $$
    it follows that
    $$
    \partial_t v_{\infty} - \Delta v_{\infty} = 0, \quad (x,t) \in \mathbb R^n \cap \{x_n > 0\} \times (-\infty,0]
    $$
    with globally Hölder-continuous $D^2_x v_{\infty}$. Moreover, differentiating the boundary condition $u(\tilde{x},0,t) = 0$, we find that $u = \partial_t u = \partial_{x_i} u = \partial_{x_i}\partial_{x_n} u = 0$ for $i \in \{1,...,n-1\}$ when $x_n = 0$. It follows that $v_k(\tilde{x},0,t) = 0$ as well for any $(\tilde{x},t) \in \mathbb R^{n-1} \times (-\infty,0]$ and this is preserved in the limit $k \to +\infty$.
    
    Hence, Liouville's theorem applies and 
    $$
v_{\infty}(x,t) = v_{\infty}(0,0) + \nabla_x v_{\infty}(0,0) \cdot x + \frac{1}{2}x^T H_x v_{\infty}(0,0)x + \partial_t v_{\infty}(0,0)t = 0
    $$
    is identically zero. This is a contradiction because 
\begin{align*}
        |D_{x_i}D_{x_j}v_{k}([0, (x_k)_n \lambda_k^{-1}], 0) - D_{x_i}D_{x_j}v_k([\lambda_k^{-1}(\tilde{y}_k - \tilde{x}_k),(y_k)_n \lambda_k^{-1}], \lambda_k^{-1}(s_k - t_k))| \\
        = M_k^{-1}\lambda_k^{-\alpha}|D_{x_i}D_{x_j}u_{k}(p_k) - D_{x_i}D_{x_j}u_k(q_k)| \geq \frac{c_n}{2} > 0 \quad \forall k
\end{align*}
    using (\ref{Holder norm absurd hypothesis}).
    
    Since $\lambda_k^{-1}(q_k-p_k)$ and $(x_k)_n \lambda_k^{-1}$ are bounded sequences as well, it follows that $(y_k)_n \lambda_k^{-1}$, $\lambda_k^{-1}(\tilde{y}_k - \tilde{x}_k$, $\lambda_k^{-1}(s_k-t_k)$ are bounded sequences as well. We can take the limit $k \to +\infty$ along a subsequence to find out that
    $$
|D_{x_i}D_{x_j}v_{\infty}(p_{\infty}^1 ) - D_{x_i}D_{x_j}v_{\infty}(p_{\infty}^2)| \geq \frac{c_n}{2} > 0
    $$
    for some $p_{\infty}^i \in \mathbb R^n \cap \{x_n \geq 0\} \times (-\infty,0]$, i.e. $D_{x_i}D_{x_j}v_{\infty}$ cannot be identically zero.
    \end{proof}

\begin{corollary}[Interior Schauder Estimates for the Heat Equation with a priori estimates] \label{interior schauder estimates, ver.1}
Assume that $u \in C^{2+\alpha,1+\alpha/2}(C(x,t;r))$ is a classical solution to the linear parabolic equation
\begin{equation}
    u_t - L_0 u = f(x,t), \quad (x,t) \in C(x,t;r), 
\end{equation}
where $f \in C^{\alpha,\alpha/2}_{x,t}(C(x,t;r))$ and $L_0$ is as in (\ref{eq: L_0 constant coefficient parabolic}).

There exists $C = C(n, \alpha, L_0)$ (independent of $x,t,r$) for which
\begin{align*}
    ||u||_{2+\alpha,1+\alpha/2;C(x,t;r/2)} &\leq C\left(1+\frac{1}{r^{2+\alpha}}\right)\left( ||f||_{\alpha,\alpha/2;C(x,t;r)} + ||u||_{L^{\infty}_{x,t}(C(x,t;r))} \right).
\end{align*}
The result still holds if one replaces the cylinder $C(x,t;r)$ by a half-cylinder
$$
    C^+(x,t;r) := \{(y,s) \in \mathbb R^n \times \mathbb R: |x-y| < r, \langle d, y-x \rangle \geq 0, t-r^2 < s < t\}, \quad |d| = 1,
    $$
and one assumes a boundary condition $u = 0$ on $C(x,t;r) \cap \{(y,t):  \langle d, y-x \rangle = 0\}$.
\end{corollary}

\begin{remark}
     It is a priori assumed that $D^2_x u$ is globally Hölder-continuous on $C(x,t;r)$. This assumption will be removed later.
\end{remark}

\begin{proof}
We use the same localization argument as in Theorem \ref{thm:Interior L^p Estimates for the Heat Equation} or Krylov's book (\cite[Theorem 8.11.1]{krylov}).

    After translation and parabolic scaling, one can assume without loss of generality that $x = t = 0$ and $r = 1$. It suffices to prove the inequality for $[D^2_x u]_{\alpha,\alpha/2}$ and use $\partial_t u = \Delta u + f$ to bound $[\partial_t u]_{\alpha,\alpha/2}$.

    We cover $C_{1} = C(0,0;1)$ with an increasing sequence of cylinders
    $$
    C_{1} = \bigcup_{l=1}^{\infty} P_l, \quad P_l = C(0,0;R_l), \quad R_l = \sum_{i=1}^l 2^{-i} \in [1/2,1]
    $$
    and consider smooth cut-off functions $\chi_l(x,t) \in C^{\infty}(\mathbb R^n \times \mathbb R;[0,1])$ for which $\chi_l = 1$ on $P_l$,
    $$
    \supp(\chi_l) \subset \{(y,s) \in \mathbb R^{n+1} : |y| \leq R_{l+1}, -R_{l+1}^2 \leq s \leq R_{l+1}^2\} = -P_{l+1} \cup P_{l+1}
    $$
    and
    $$
    ||\chi_l||_{2+\alpha,1+\alpha/2} \lesssim 2^{3l} =: \rho^{-l}.
    $$
    Then $v_l = u\chi_l \in  C_c^{2+\alpha, 1+\alpha/2}(\mathbb R^n \times (-\infty,0])$ solves
    $$
        \partial_t v_l - L_0 v_l = f\chi_l +u \partial_t \chi_l + u L_0 \chi_l + \sum_{i,j =1}^n a^{i,j}\left(\partial_{x_i} u \partial_{x_j} \chi_l +  \partial_{x_j} u \partial_{x_i} \chi_l \right), \quad (x,t) \in \mathbb R^n \times (-\infty,0].
    $$
    Global Schauder estimates show that
\begin{align}
   ||u\chi_l||_{2+\alpha,1+\alpha/2;P_{l+1}} &= ||u\chi_l||_{2+\alpha,1+\alpha/2;\mathbb R^n \times (-\infty,0]} \notag \\
   &\leq C [f\chi_l +u \partial_t \chi_l + u \Delta \chi_l +  \sum_{i,j} \partial_{x_i} u \partial_{x_j} \chi_l]_{\alpha,\alpha/2;\mathbb R^n \times (-\infty,0]} + C||u \chi_l||_{L^{\infty}_{x,t}(\mathbb R^n \times (-\infty,0])} \notag \\
    &= [f\chi_l +u \partial_t \chi_l + u L_0 \chi_l + \sum_{i,j} \partial_{x_i} u \partial_{x_j} \chi_l]_{\alpha,\alpha/2;P_{l+1}} + C||u \chi_l||_{L^{\infty}_{x,t}(P_{l+1})} \notag \\
    &\leq \rho^{-l} C\left(  ||f||_{\alpha,\alpha/2;P_{l+1}} + ||u||_{\alpha,\alpha/2;P_{l+1}} + ||\nabla_x u||_{\alpha,\alpha/2;P_{l+1}} \right) \label{schauder estimates proof, upper bound on D^2 u chi}
\end{align}
    thanks to the product rule
    $$
    [fg]_{\alpha,\alpha/2;U} \leq [f]_{\alpha,\alpha/2;U} ||g||_{L^{\infty}_{x,t}(U)} + ||f||_{L^{\infty}_{x,t}(U)}[g]_{\alpha,\alpha/2;U}.
    $$
    We introduce back $\chi_{l+1}$
    \begin{align*}
    ||u\chi_l||_{2+\alpha,1+\alpha/2;P_{l+1}} &\leq \rho^{-l}C\left(  ||f||_{\alpha,\alpha/2;P_{l+1}} + ||u||_{\alpha,\alpha/2;P_{l+1}} + ||\nabla_x u||_{\alpha,\alpha/2;P_{l+1}} \right) \\
    &= \rho^{-l}C\left(  ||f||_{\alpha,\alpha/2;P_{l+1}} + ||u \chi_{l+1}||_{\alpha,\alpha/2;P_{l+1}} + ||\nabla_x (u\chi_{l+1})||_{\alpha,\alpha/2;P_{l+1}} \right) \\
    &\leq \rho^{-l}C\left(  ||f||_{\alpha,\alpha/2;C_1} + ||u \chi_{l+1}||_{\alpha,\alpha/2;\mathbb R^n \times (-\infty,0]} + ||\nabla_x (u \chi_{l+1})||_{\alpha,\alpha/2;\mathbb R^n \times (-\infty,0]} \right) 
    \end{align*}
    and from now on, $C = C(n,\alpha) > 1$ is fixed and $\tilde{C} = \tilde{C}(n,\alpha) > 1$ denotes the constant from the interpolation Lemma \ref{interpolation lemma}. Apply the lemma with a radius $\tilde{\rho} = (2C\tilde{C})^{-1}\rho^3 \rho^l$. Then
   \begin{align*}
||u\chi_l||_{2+\alpha,1+\alpha/2;P_{l+1}} &\leq \rho^{-l} 2C\tilde{C}\left( ||f||_{\alpha,\alpha/2;C_1}  + \varepsilon [u\chi_{l+1}]_{2+\alpha,1+\alpha/2;\mathbb R^n \times (-\infty,0]} + \varepsilon^{-2} ||u \chi_{l+1}||_{L^{\infty}_{x,t}(\mathbb R^n \times (-\infty,0])} \right) \\
    &\leq \rho^{-l} 2C \tilde{C} \left( ||f||_{\alpha,\alpha/2;C_1}  + \varepsilon [u\chi_{l+1}]_{2+\alpha,1+\alpha/2;P_{l+2}} + \varepsilon^{-2} ||u||_{L^{\infty}_{x,t}(C_1)} \right) \\
    &\leq \rho^{-3l} \underbrace{(2C \tilde{C})^3}_{=: \hat{C}} \left( ||f||_{\alpha,\alpha/2;C_1}  + ||u ||_{L^{\infty}_{x,t}(C_1)} \right) + \frac{\rho^3}{2} [u\chi_{l+1}]_{2+\alpha,1+\alpha/2;P_{l+2}} 
  \end{align*}
  Starting with $l = 1$ and iterating, we obtain
\begin{align}
||u\chi_1||_{2+\alpha,1+\alpha/2;P_{2}} &\leq  \rho^{-3}\hat{C} \left( ||f||_{\alpha,\alpha/2;C_1}  + ||u ||_{L^{\infty}_{x,t}(C_1)} \right) + \frac{\rho^3}{2} [u\chi_{2}]_{2+\alpha,1+\alpha/2;P_{3}} \notag \\
&\leq \left( \rho^{-3}\hat{C} + \frac{\rho^3}{2} \cdot \rho^{-6}\hat{C} \right) \left( ||f||_{\alpha,\alpha/2;C_1}  + ||u ||_{L^{\infty}_{x,t}(C_1)} \right) + \left(\frac{\rho^3}{2}\right)^2 [u\chi_{3}]_{2+\alpha,1+\alpha/2;P_{4}} \notag \\
&\leq \rho^{-3} \hat{C} \left( 1 + \frac{1}{2} + ... + \frac{1}{2^{k-2}} \right) \left( ||f||_{\alpha,\alpha/2;C_1}  + ||u||_{L^{\infty}_{x,t}(C_1)} \right) \notag \\
&\phantom{\leq} + \left(\frac{\rho^3}{2}\right)^{k-1} [u\chi_{k}]_{2+\alpha,1+\alpha/2;P_{k+1}} \label{schauder estimates proof, final inequality}
  \end{align}
  for any $k > 1$. Since $\rho < 1$ and 
  $$
  ||u\chi_{k}||_{2+\alpha,1+\alpha/2;P_{k+1}} \lesssim \rho^{-k} \left(  ||f||_{\alpha,\alpha/2;C_1} + ||u||_{1+\alpha,\alpha/2;C_1} \right) 
  $$
  using (\ref{schauder estimates proof, upper bound on D^2 u chi}), we can let $k \to +\infty$ in (\ref{schauder estimates proof, final inequality}). Observing that $u = u \chi_1$ on $C(0,0;1/2) \subset P_2$, this concludes the proof.

  In the half-cylinder case, up to applying a rotation, one can assume that $d = e_n$. Then the proof is exactly the same. We use the same cut-off $\chi_l$ supported on $-P_{l+1} \cup P_l$, but use the Global Schauder Estimates on the half-space to deduce the analog to (\ref{schauder estimates proof, upper bound on D^2 u chi}), which is
  \begin{align*}
   ||u\chi_l||_{2+\alpha,1+\alpha/2;P^+_{l+1}} &\leq \rho^{-l} C\left(  ||f||_{\alpha,\alpha/2;P^+_{l+1}} + ||u||_{\alpha,\alpha/2;P^+_{l+1}} + ||\nabla_x u||_{\alpha,\alpha/2;P^+_{l+1}} \right)
\end{align*}
where $u \chi_l \in C^{2+\alpha,1+\alpha/2}_c(\mathbb R^n \cap \{x_n \geq 0\} \times (-\infty,0])$ and $P^+_{l+1} = C^+(0,0;R_l)$. Then we introduce back $\chi_{l+1}$, use interpolation and iterate as above.  The constant $C(n,\alpha,L_0)$ depends only on $\Lambda$ and the entries of the matrix $(a_{i,j})$ which defines $L_0$.
\end{proof}

\begin{corollary}[Interior Schauder Estimates for the Heat Equation with smooth forcing term] \label{cor:interior schauder estimates with smooth forcing}
Assume that $u \in L^{\infty}(C(x,t;r))$ is a distributional solution to the linear parabolic equation
\begin{equation}
    u_t - L_0 u = f(x,t), \quad (x,t) \in C(x,t;r),
\end{equation}
where $f \in C^{\alpha,\alpha/2}_{x,t}(C(x,t;r)) \cap C^{\infty}_{x,t}(C(x,t;r))$ and $L_0$ is as in (\ref{eq: L_0 constant coefficient parabolic}).

There exists $C = C(n, \alpha, L_0)$ (independent of $x,t,r$) for which
\begin{align*}
    ||u||_{2+\alpha,1+\alpha/2;C(x,t;r/2)} &\leq C\left(1+\frac{1}{r^{2+\alpha}}\right)\left( ||f||_{\alpha,\alpha/2;C(x,t;r)} + ||u||_{L^{\infty}_{x,t}(C(x,t;r))} \right).
\end{align*}
\end{corollary}

\begin{proof}
    By hypoellipticity, $u \in C^{\infty}(C(x,t;r))$. Take $r_n = r - 1/n$ and start at $n \geq n_0$ large enough so that $r > 1/n$. Take also a decreasing sequence $0 < \varepsilon_n < 2r/n - 1/n^2$ with $\varepsilon_n \to 0$. One obtains a sequence of increasing cylinders $(C(x,t-\varepsilon_n;r_n))_{n \geq n_0}$ with the inclusion
    $$
    \overline{C(x,t-\varepsilon_n;r_n)} \subset C(x,t;r)
    $$
    It then suffices to apply Corollary \ref{interior schauder estimates, ver.1} on these cylinders and let $n \to +\infty$.
    
\end{proof}

\begin{corollary}[Global Schauder Estimates for compactly supported solutions]\label{cor:global schauder estimates for compactly supported solutions}
    Assume that $u \in C^{0}_{b}(\mathbb R^n \times (-\infty,0])$, $\alpha \in (0,1)$, is a distributional solution to the linear heat equation
\begin{align*}
    u_t - L_0 u &= f(x,t), \quad (x,t) \in \mathbb R^n \times (-\infty,0),
\end{align*}
where $L_0$ is as in (\ref{eq: L_0 constant coefficient parabolic}).

Further assume that $u$ is compactly supported in the sense that there is some $x_1 \in \mathbb R^n$ and $r > 0$ for which $\overline{\supp(u(\cdot,t))} \subset B_r(x_1)$ for all $t \leq 0$, $u(\cdot,t) = 0$ for $t \leq -r^2$ and $f \in C^{\alpha,\alpha/2}_{x,t}(C(x,0;r))$. 

There exists $C = C(n, \alpha, L_0)$ for which
$$
    ||u||_{2+\alpha,1+\alpha/2;\mathbb R^n \times (-\infty,0]} \leq C\left(1+\frac{1}{r^{2+\alpha}}\right)\left( ||f||_{\alpha,\alpha/2;\mathbb R^n \times (-\infty,0]} + ||u||_{L^{\infty}_{x,t}(\mathbb R^n \times (-\infty,0])} \right).
$$
\end{corollary}

\begin{proof}
After translation, one can assume $x_1 = 0$. After applying a linear transformation (and up to taking a larger $r$), one can assume that $L_0 = \Delta$. After parabolic rescaling, one can assume without loss of generality that $r = 1$.

On $\mathbb R^n \times [-4,0)$, $u \in C^0_b(\mathbb R^n \times [-4,0))$ is the unique solution for the problem 
    \begin{align*}
     \partial_t u - \Delta u &= f(x,t) \in C^0_b(\mathbb R^n \times [-4,0)), \quad (x,t) \in \mathbb R^n \times (-4,0), \\
u(x,-4) &= 0, \quad x \in \mathbb R^n,
 \end{align*}
 which is given by convolution with the Gaussian heat kernel (up to a time translation), as shown in Lemma \ref{lemma:uniqueness heat equation}.
 
 Note that given the compact support assumption of $u(\cdot,t)$, the same holds for $f$. Setting $f(x,t) = 0$ for $x \in B_r(x_1)^c$, $f(x,t) = f(x,-t)$ for $t > 0$, $f(x,t) = 0$ for $t \leq -1$ yields a $C^{\alpha,\alpha/2}(\mathbb R^{n+1})$ extension of $f$.  Approximate $f$ by a sequence of $C^{\infty}_c(\mathbb R^{n+1})$-functions $(f_n)_{n \geq 0}$ with $||f_n-f||_{\alpha,\alpha/2;\mathbb R^n \times (-\infty,0)} \to 0$. Let $(u_n) \in C^0_b(\mathbb R^n \times [-4,0))$ be the corresponding sequence of solutions for
    \begin{align*}
     \partial_t u_n - \Delta u_n &= f_n(x,t), \quad (x,t) \in \mathbb R^n \times (-4,0) \\
u_n(x,-4) &= 0, \quad x \in \mathbb R^n
 \end{align*}
 given by convolution with the heat kernel. Thanks to (\ref{Gaussian heat kernel convolution, young inequality}), one has 
   $$
   ||u_n-u||_{L^{\infty}(\mathbb R^n \times (-4,0))} \to 0.
   $$
    By Corollary \ref{cor:interior schauder estimates with smooth forcing},
    \begin{align*}
            ||u_n-u_m||_{2+\alpha,1+\alpha/2;C(0,0;1)} &\leq C\left(1+\frac{1}{r^{2+\alpha}}\right)\left( ||f_n-f_m||_{\alpha,\alpha/2;C(0,0;2)} + ||u_n-u_m||_{L^{\infty}_{x,t}} \right),
    \end{align*}
    meaning that $(u_n)$ is a Cauchy sequence in the Banach space $C^{2+\alpha,1+\alpha/2}(C(0,0;1))$, where $C(0,0;1) \subset C(0,0;2) \subset \mathbb R^n \times (-4,0)$. Hence, $u \in C^{2+\alpha,1+\alpha/2}(C(0,0;1))$ and $u_n \to u$ in $C^{2+\alpha,1+\alpha/2}(C(0,0;1))$. Once more, by Corollary \ref{cor:interior schauder estimates with smooth forcing},
    \begin{align*}
            ||u_n||_{2+\alpha,1+\alpha/2;C(0,0;1)} &\leq C\left(1+\frac{1}{r^{2+\alpha}}\right)\left( ||f_n||_{\alpha,\alpha/2;C(0,0;2)} + ||u_n||_{L^{\infty}_{x,t}(C(0,0;2))} \right) \\
            &\leq 2 C\left(1+\frac{1}{r^{2+\alpha}}\right)\left( ||f||_{\alpha,\alpha/2;C(0,0;2)} + ||u||_{L^{\infty}_{x,t}(C(0,0;2))} \right) \\
            &= 2C\left(1+\frac{1}{r^{2+\alpha}}\right)\left( ||f||_{\alpha,\alpha/2;C(0,0;1)} + ||u||_{L^{\infty}_{x,t}(C(0,0;1))} \right)
    \end{align*}
    for all $n$ large enough. Take the limit $n \to +\infty$ to conclude.
\end{proof}

\begin{theorem}[Interior Schauder Estimates for the Heat Equation] \label{interior schauder estimates}
Assume that $u \in C^{1+\alpha,\alpha/2}(C(x,t;r))$ is a solution to the linear heat equation
\begin{equation}
    u_t - L_0 u = f(x,t), \quad (x,t) \in C(x,t;r), 
\end{equation}
where $f \in C^{\alpha,\alpha/2}_{x,t}(C(x,t;r))$ and $L_0$ is as in (\ref{eq: L_0 constant coefficient parabolic}).

There exists $C = C(n, \alpha, L_0)$ (independent of $x,t,r$) for which
\begin{align*}
    ||u||_{2+\alpha,1+\alpha/2;C(x,t;r/2)} &\leq C\left(1+\frac{1}{r^{2+\alpha}}\right)\left( ||f||_{\alpha,\alpha/2;C(x,t;r)} + ||u||_{L^{\infty}_{x,t}(C(x,t;r))} \right).
\end{align*}
\end{theorem}

\begin{remark}\label{weakening of interior Schauder estimates}
    If $L_0 = \Delta$, the a priori regularity on $u, \nabla_x u$ can be weakened. Thanks to Interior $L^p$-estimates (Theorem \ref{thm:Interior L^p Estimates for the Heat Equation}) and the Parabolic Sobolev Embedding (Corollary \ref{cor: local parabolic sobolev embedding}), it suffices that $u \in L^p(C(x,t;r))$, $\nabla_x u \in L^{p}(C(x,t;r))$ for $p$ large enough (depending on $\alpha$) to get $C^{1+\alpha,\alpha/2}$-regularity on a smaller cylinder $C(x,t;r/4)$. Then one can apply the interior Schauder estimates to get $C^{2+\alpha,1+\alpha/2}$-regularity on $C(x,t;r/8)$.
\end{remark}

\begin{proof}
    The proof is the same as the proof of Corollary \ref{interior schauder estimates, ver.1}. The only difference is that we use the Global Schauder Estimates for compactly supported solutions (Corollary \ref{cor:global schauder estimates for compactly supported solutions}) on $v_l = u\chi_l$ instead of Theorem \ref{thm:global schauder estimates}.
\end{proof}

The method of 'freezing' the coefficients allows to generalize Corollary \ref{interior schauder estimates, ver.1} to more general parabolic operators.

\begin{theorem}[Boundary Schauder Estimates for parabolic equations on half-cylinders] \label{boundary schauder estimates, ver.1}
Let $\mu \in (0,1), \rho  > 0,x_0 \in \mathbb R^n, t_0 \in \mathbb R$. Consider 
$$
    C^+(x_0,t_0;\rho) := \{(y,s) \in \mathbb R^n \times \mathbb R: |x_0-y| < \mu, y_n > (x_0)_n , t-\rho^2 < s < t\}
    $$
a half-cylinder and a parabolic operator
    $$
    L = c(x,t) + \sum_{i=1}^n b^i(x,t) D_{x_i} + \sum_{i,j=1}^n a^{i,j}(x,t) D_{x_i}D_{x_j}
    $$
    defined on $\overline{C^+(x_0,t_0;\rho) }$ and satisfying
    \begin{align*}
        ||a,b,c||_{\mu,\mu/2;\overline{C^+(x_0,t_0;\rho) }} &\leq \Lambda, \quad \Lambda > 0, \\
        \Lambda^{-1}|\xi|^2 \leq \sum_{i,j=1}^n a^{i,j}(x,t)\xi_i\xi_j &\leq \Lambda |\xi|^2, \quad (x,t) \in C^+(x_0,t_0;\rho) , \xi \in \mathbb R^n.
    \end{align*}

Assume that $u \in C^{2+\mu,1+\mu/2}(C^+(x_0,t_0;\rho))$ is a classical solution to the linear parabolic equation
\begin{align*}
    u_t - L  u &= f(x,t), \quad (x,t) \in C^+(x_0,t_0;\rho),  \\
    u &= u_0(x,t), \quad (x,t) \in C^+(x_0,t_0;\rho) \cap \{x_n = (x_0)_n \},
\end{align*}
where $f \in C^{\mu,\mu/2}_{x,t}(C^+(x_0,t_0;\rho) )$, $u_0 \in C^{2+\mu,1+\mu/2}_{x,t}(C^+(x_0,t_0;\rho) )$.

There exists $C = C(n, \mu, L)$ (independent of $x_0,t_0,\rho$) and $\kappa = \kappa(n,p,\mu,L)$ for which
\begin{align}
||D^2_xu||_{\mu,\mu/2;C^{+}_{\rho/2}} + ||\partial_t u||_{\mu,\mu/2;C^{+}_{\rho/2}} &\leq C \bigg( ||f||_{\mu,\mu/2;C^{+}_{\rho}} + ||\partial_t u_0||_{\mu,\mu/2;C^{+}_{\rho}} \notag \\
&+\sum_{i,j \neq n} ||D_{x_i}D_{x_j} u_0||_{\mu,\mu/2;C^{+}_{\rho}}  + \rho^{-2-\mu}||u||_{L^{\infty}_{x,t}(C^{+}_{\rho}) } \bigg) \label{boundary schauder, norm}
\end{align}
whenever $\rho \leq \kappa$.
\end{theorem}

\begin{proof}
Without loss of generality, one can assume $u_0 = 0$ by replacing $u$ by $u + u_0$ and $f$ by $f + \partial_t u_0 - L u_0$. One can also assume that $b^i = c = 0$. Indeed, extend the coefficients $c$, $b^i$, $a^{i,j}$, $a^{n,n}$, $i,j \neq n$, by odd reflection, the coefficients $b^n, a^{i,n}, a^{n,i}$, $i \neq n$, by even reflection and $u$ by odd reflection to the whole cylinder $C(x_0,t_0;\rho)$ (Proposition \ref{prop:extension by reflection}), and use interpolation to obtain
\begin{align}
    ||D^2_xu||_{L^{\infty}(C^+)} + ||\partial_t u||_{L^{\infty}(C^+)} 
    &\lesssim \varepsilon [D^2_x u]_{\mu,\mu/2;C^+} + \varepsilon [\partial_t u]_{\mu,\mu/2;C^+}  + ||u||_{L_{x,t}^{\infty}(C^+)} \label{interpolation and cut-off, lower order terms}
\end{align}
for $\varepsilon > 0$ arbitrarily small. We can proceed similarly to estimate terms $c(x,t) u$ and $b^i(x,t) \partial_{x_i} u$ coming from $Lu$.

First, we observe that (\ref{boundary schauder, norm}) is true if $L$ has constant coefficients by Corollary \ref{interior schauder estimates, ver.1}. For general parabolic operators $L$, assume without loss of generality that $x = t = 0$, $\rho = 1$. Fix $(x_0,t_0) \in C^+_{1/2}$. Write
    $$
L_0 = \sum_{i,j=1}^n a^{i,j}(x_0,t_0) D_{x_i}D_{x_j}.
$$
    One estimates 
\begin{align*}
        ||u||_{2+\mu,1+\mu/2;C^+_{1/2}} 
        &\lesssim C\left( ||(\partial_t - L_0)u||_{\mu,\mu/2;C^+_{1}} + ||u||_{L^{\infty}_{x,t}(C^+_{1})} \right) 
\end{align*}
    using the constant-coefficient case for $L = L_0$. Then
\begin{align*}
    ||(\partial_t - L_0)u||_{\mu,\mu/2;C^+_{1}} &\leq ||(\partial_t - L)u||_{\mu,\mu/2;C^+_{1}} \\
    &+ \bigg|\bigg| \sum_{i,j=1}^n \left[ a^{i,j}(x_0,t_0) - a^{i,j}(x,t) \right] D_{x_i}D_{x_j}u \bigg|\bigg|_{\mu,\mu/2;C^+_{1}}  \\
    &\lesssim  ||(\partial_t - L)u||_{\mu,\mu/2;C^+_{1}}  \\
    &+ \sum_{i,j=1}^n ||a^{i,j}(x_0,t_0) - a^{i,j}(x,t)||_{L_{x,t}^{\infty}(C^+_{1})} \cdot ||D^2_x u||_{\mu,\mu/2;C^+_{1}} \\ 
    &+ \sum_{i,j=1}^n ||a^{i,j}(x_0,t_0) - a^{i,j}(x,t)||_{\mu,\mu/2;C^+_{1}} \cdot ||D^2_x u||_{L_{x,t}^{\infty}(C^+_{1})} \\
    &\lesssim  ||(\partial_t - L)u||_{\mu,\mu/2;C^+_{1}} + \Lambda  \kappa^{\alpha} ||D^2_x u||_{\mu,\mu/2;C^+_{1}}  +  \Lambda ||D^2_x u||_{L^{\infty}_{x,t}(C^+_{1})}
\end{align*}
and we estimate $||D^2_x u||_{L^{\infty}(C^+_{1})}$ as in (\ref{interpolation and cut-off, lower order terms}). If $\kappa$ is small enough, the theorem is proved.
\end{proof}

 Finally, one obtains a boundary Schauder estimate for bounded domains with smooth boundary by flattening and covering the boundary.

\begin{theorem}[Boundary Schauder Estimates for parabolic equations on smooth domains] \label{boundary schauder smooth boundary}
    Assume that $u \in C^{2+\mu,1+\mu/2}_{x,t}(\Omega \times (0,T))$ is a classical solution to the linear parabolic equation
\begin{align*}
    u_t - L u &= f(x,t), \quad (x,t) \in \Omega \times (0,T), \\
    u &= u_0(x,t), \quad (x,t) \in \partial \Omega \times (0,T),
\end{align*}
where $f \in C^{\mu,\mu/2}_{x,t}(\Omega \times (0,T))$, $u_0 \in C^{2+\mu,1+\mu/2}_{x,t}(\Omega \times (0,T))$, $\Omega$ is open, bounded with $C^{2+\mu}$-boundary and $L$ is a parabolic operator satisfying the hypotheses from Theorem \ref{boundary schauder estimates, ver.1} on $\overline{\Omega} \times [0,T]$ for some constant $\Lambda > 0$.

For any interval $I = (\delta,T)$, $\delta > 0$, there exists $C =  C(n,\mu,\Lambda, \Omega, I)$ and $\varepsilon = \varepsilon(n,\mu,\Lambda, \Omega, I)$ such that the estimate
    \begin{align*}
         ||(|D^2_xu| + |\partial_t u|)||_{\mu,\mu/2;\Omega_{\varepsilon}  \times I} 
         &\leq C \bigg( ||f||_{\mu,\mu/2;\Omega \times (0,T)} + ||\partial_t u_0||_{\mu,\mu/2;\Omega \times (0,T)}\\
         &+ \sum_{i,j \neq n} ||D_{y_i}D_{y_j} (u_0)||_{\mu,\mu/2; \Omega \times (0,T)}  + ||u||_{L^{\infty}_{x,t}(\Omega \times (0,T)) } \bigg)
     \end{align*}
     holds on the $\varepsilon$-neighborhood
     $$
    \Omega_{\varepsilon} = \{x \in \Omega: \dist(x,\partial \Omega) < \varepsilon\}
    $$
    of the boundary. 

    For the heat operator $L = \Delta$, if $\Omega$ has $C^{\infty}$-boundary, then the assumption $u \in C^{2+\mu,1+\mu/2}_{x,t}(\Omega \times (0,T))$ can we weakened to $u \in C^0([0,T);C^1(\overline{\Omega}))$. 
\end{theorem}

\begin{remark}
    To obtain smoothness at time $t = 0$, one needs to assume some additional compatibility condition between the initial and boundary data $u_0$ and the forcing term $f$. See for example \cite[Chapter 7.1.3, Theorem 7]{evans10}.
\end{remark}

\begin{proof}
    The theorem is a direct consequence of Theorem \ref{boundary schauder estimates, ver.1} by locally straightening the boundary. 

    For the heat operator $L = \Delta$, replacing $u$ with $u+u_0$ and $f$ by $f - (\partial_t u_0 - \Delta u_0)$, one can assume that $u_0 = 0$. 
    
    Letting $\chi(t) \in C^{\infty}_c(\mathbb R;[0,1])$ be a cut-off which is $1$ on $I = (\delta,2T)$, $0$ on $(-\infty,\delta/2)$, it follows that $v = u \cdot \chi(t)$ solves
\begin{align*}
    v_t - L v &= f \cdot \chi(t) + u \cdot \chi'(t), \quad (x,t) \in \Omega \times (0,T), \\
    v &= 0, \quad (x,t) \in \partial \Omega \times (0,T) \cup \Omega \times \{0\},
\end{align*}
where $u \cdot \chi'(t) \in C^0([0,T) \times \overline{\Omega})$. Applying Corollary \ref{cor: local parabolic sobolev embedding} and Theorem \ref{thm:boundary parabolic sobolev embedding} for $0 < T' < T$, it follows that $v \in C^{1+\mu,\mu/2}(\overline{\Omega} \times [0,T])$ with 
$$
||u||_{C^{1+\mu,\mu/2}(\overline{\Omega} \times [\delta,T])} 
 \leq ||v||_{C^{1+\mu,\mu/2}(\overline{\Omega} \times [0,T])} \leq C(n,\mu,\Omega,T,\delta) \left( ||f||_{L^{\infty}(\overline{\Omega} \times [0,T])} + ||u||_{L^{\infty}} \right)
$$
for any $\delta \in (0,T)$. Hence, $F = f \cdot \chi(t) + u \cdot \chi'(t) \in C^{\mu,\mu/2}_{x,t}(\overline{\Omega} \times [0,T])$ with
$$
||F||_{C^{1+\mu,\mu/2}(\overline{\Omega} \times [0,T])} 
 \leq  C(n,\mu,\Omega,T,\delta) \left( ||f||_{C^{\mu,\mu/2}_{x,t}(\overline{\Omega} \times [0,T])} + ||u||_{L^{\infty}(\overline{\Omega} \times [0,T])} \right)
$$
as well. Let $F_n \in C^{\infty}(\overline{\Omega} \times [0,T])$, $F_n(x,t) = 0$ for $t \leq \delta/4$ converge to $F$ in $C^{\mu,\mu/2}_{x,t}(\overline{\Omega} \times [0,T])$. Such an approximation can be constructed by first extending $F_n$ as a global $C^{\mu,\mu/2}(\mathbb R^{n+1})$ function
$$
F_n(x,t) = \inf\{F_n(y,s) + [F_n]_{C^{\mu,\mu/2}_{x,t}(\overline{\Omega} \times [0,T])} \cdot |(x,t) - (y,s)|^{\mu}: (y,s) \in \overline{\Omega} \times [0,T]\},
$$
and then using a convolution with an approximate identity and a smooth cut-off to remove small times. Here, $|(x,t) - (y,s)| = |x-y| + |s-t|^{1/2}$ denotes the parabolic distance.

Let $v_n$ solve $(\partial_t  - \Delta)v_n= F_n$ on $\Omega \times (0,T)$, $v_n = 0$ on $\partial \Omega \times (0,T) \cup \Omega \times \{0\}$. Then $v_n \to v$ uniformly (see e.g. the proof Theorem \ref{thm:boundary parabolic sobolev embedding}) and $v_n \in C^{\infty}(\overline{\Omega} \times [0,T])$ (\cite[Chapter 7.1.3, Theorem 7]{evans10}). It follows from the boundary case with a priori estimates and the interior case (Theorem \ref{interior schauder estimates}) that each $v_n$, $v_n-v_m$ is $C^{2+\mu,1+\mu/2}$ with norms
\begin{align}
    ||v_n||_{C^{2+\mu,1+\mu/2}(\overline{\Omega} \times [\delta,T])} &\leq  
  C(n,\mu,\Omega,T,\delta) \left( ||F_n||_{C^{\mu,\mu/2}_{x,t}(\overline{\Omega} \times [0,T])} + ||v_n||_{L^{\infty}(\overline{\Omega} \times [0,T])} \right) \notag \\
  &\leq C(n,\mu,\Omega,T,\delta) \left( ||F||_{C^{\mu,\mu/2}_{x,t}(\overline{\Omega} \times [0,T])} + ||v||_{L^{\infty}(\overline{\Omega} \times [0,T])} \right) \notag \\
  &\leq C(n,\mu,\Omega,T,\delta) \left( ||f||_{C^{\mu,\mu/2}_{x,t}(\overline{\Omega} \times [0,T])} + ||u||_{L^{\infty}(\overline{\Omega} \times [0,T])} \right) \label{eq:boundary schauder estimate last inequality}
\end{align}
for all $n$ large enough, and similarly for $v_n-v_m$. In the limit, we obtain that $v \in C^{2+\mu,1+\mu/2}(\overline{\Omega} \times [\delta,T])$ with
$$
||u||_{C^{2+\mu,1+\mu/2}(\overline{\Omega} \times [\delta,T])} = ||v||_{C^{2+\mu,1+\mu/2}(\overline{\Omega} \times [\delta,T])}
$$
satisfying (\ref{eq:boundary schauder estimate last inequality}).
\end{proof}

\chapter{On the rigidity of blow-up for the Harmonic Map Heat Flow from $B^2$ to $S^2$}
~\newline~

The content of this chapter is taken from my article \cite{samuelian2025blowuptrees}, to appear soon in \textit{Calculus of Variations and Partial Differential Equations}. For the sake of completeness, the proof of the Comparison Principle has been added.

Shortly after defending this thesis, I submitted a preprint proving nonexistence of finite-time blow-up for the $k$-equivariant harmonic map heat flow from $B^2$ to $S^2$ when $k \geq 2$, see \cite{Samuelian2026HMHF}.

\section{Introduction}

We are interested in the behaviour of finite-time and finite-energy solutions to the harmonic map heat flow (\ref{heat map flow}) from the 2D-ball to the 2-sphere
\begin{align}
    v_t &= \Delta v + |\nabla v|^2 v, \quad (x,t) \in B^2 \times  (0,+\infty) \label{heat map flow} \\
    v(x,t) &= v_0(x,t), \quad (x,t) \in \partial B^2 \times [0,+\infty) \cup B^2 \times \{0\}, \notag
\end{align}
where $v: B^2 \rightarrow S^2 \subset \mathbb R^3$, $\nabla v = (Jv)^T$ is the transposed Jacobian matrix, the energy is defined as
\begin{equation}\label{energy of v}
    E(v(\cdot,t)) = E(v(t)) = \frac{1}{2}\int_{B^2} |\nabla v(\cdot,t)|^2 dx 
\end{equation}
and the boundary data is the restriction of a function $v_0 \in C^{\infty}(\overline{B^2} \times [0,+\infty); S^2)$ which is smooth and $k$-equivariant, meaning it has the form:
\begin{equation}
    v_0(re^{i\theta},t) = \left( e^{ik\theta} \sin h_0(r,t), \cos h_0(r,t) \right) \in \mathbb C \times \mathbb R \simeq \mathbb R^{3},  \label{k-equivariance}
\end{equation}
with $h_0(r,t) \in C^{0}([0,1] \times [0,+\infty),\mathbb R)$ being the \textit{inclination coordinate} of $v_0$ (which can be chosen, without loss generality, smooth and such that $h(r,t) = o(r^k)$ for fixed time, see Lemma \ref{k-equivariance, smoothness of h}). It follows that $v(x,t)$ is also $k$-equivariant (Proposition \ref{prop:equivalent fomrulation heat map flow}) with an inclination coordinate $h(r,t)$ solving the following one-dimensional nonlinear heat equation \begin{align}
h_t &= h_{rr} + \frac{h_r}{r} - k^2 \frac{\sin(2h)}{2r^2}, \quad (r,t) \in (0,1) \times (0,+\infty), \label{heat map formulation for h(r,t)} \\
h(r,t) &= h_0(r,t), \quad (r,t) \in \{0,1\} \times [0,+\infty) \cup (0,1) \times \{0\}. \notag
\end{align}

The study of the Harmonic Map Heat Flow goes back to the foundational work of by Eells and Sampson (\cite{Eells1964}). The equation was first introduced as a mean of finding harmonic maps which are homotopic to some given smooth map $v_0 \in C^{\infty}(M,N)$ inbetween two compact manifolds without boundary, $N$ having non-positive curvature (the analog work with boundaries is due to Hamilton, \cite{Hamilton1975}). It turns out that finding such harmonic maps is related to the limiting behaviour of solutions to the Harmonic Map Heat Flow with initial condition $v_0$ near global time. Struwe (\cite{struwe2008variational}) was able to remove the non-negative curvature assumption on $N$ by assuming that $N$ is a 2D-surface instead (with partial results in higher dimensions, \cite{Struwe1988OnTE}). Moreover, he described the formation of both finite-time and global-time singularities in terms of concentration of energy (Lemma \ref{lemma:no concentration of energy away origin, general case} and Theorem \ref{thm:smoothness around all but finitely many points}) and he was able to show that one can extract an harmonic map from $S^2$ to $N$ near such singularities (the analog work with boundaries and a boundary data $v_0(x,t) = v_0(x)$ is due to Chang, \cite{ChangStruweAnalog}). The description of the blow-up event was further improved by Qing (\cite{Qing1995OnSO}), who showed that near a singularity, the solution decomposes asymptotically as a sum of rescaled harmonic maps and a radiation term (Theorem \ref{thm: full bubble decomposition}). Such a rescaled harmonic map is called a \textit{bubble} of energy because the sum of the energy of all the bubbles accounts for all the energy loss which occurs at the blow-up event (see (\ref{pythagorean decomposition of energy})). Further exposition on the topic of harmonic maps can be found in the textbook \cite{analysisofharmonicmaps}. 

Finite time blow-up was shown to exist in the $1$-equivariant setting by Chang, Ding and Ye (\cite{Chang1992FinitetimeBO}) and infinite-time blow-up was exhibited by Gustafson, Nakanishi and Tsai for $k  = 2$ (\cite{stability}) and Chang, Ding for $k = 1$ (\cite{ChangDing2011}). However, the precise form of the decomposition (How many bubbles ? Is the decomposition unique along all possible time sequences ? Is the decomposition continuous in time ?) remained unknown for a long time. In the $1$-equivariant setting on $B^2$, Van der Hout (\cite{VANDERHOUT2003}) was able to prove that only single-bubbling can occur. In the $\mathbb R^2$ case instead of $B^2$, it has been recently proved by Jendrej and Lawrie (\cite{jendrej}) that the bubbling happens continuously in time (in the $k$-equivariant setting) using dispersive PDEs techniques. In the general non-symmetrical case, one only has a partial result, due to Jendrej, Lawrie and Schlag (\cite{Jendrej_Lawrie_Schlag_2025}): along any sequence of times, one can find a subsequence along which bubbling happens. However, they only treat the case $v_0(x,t) = v_0(x)$ and their methods rely on the monotonicity of the energy, which is not available in our setting. Even more recently, it has been shown by Kim and Merle (\cite{kim2025classificationglobaldynamicsenergycritical}) that there is no blow-up at all when $k \geq 3$, still in the $\mathbb R^2$ case with $v_0(x,t) = v_0(x)$. This was already hinted by a previous stability result from Gustafson, Nakanishi and Tsai (\cite{stability}) for $k \geq 3$. It is possible that no blow-up occurs for the ball when $k \geq 3$, but in any case, our rigidity result covers the unknown $k = 2$ case.

The aim of this paper is to provide a self-contained proof that if $v(x,t)$ blows up in finite time, then there can be only one bubble in the multi-bubble decomposition (\ref{multi bubble decomposition along a sequence of times}), following the intersection-comparison argument from Van der Hout (\cite{VANDERHOUT2003}) and Matano (\cite{matano}). This rigidity result could be used as a starting point to show nonexistence of finite time blow-up for $k \geq 2$. 

Our main contribution to the theory is the extension of the single bubble result from \cite{VANDERHOUT2003} for $k$-equivariant boundary data $v_0(x,t)$ from $k = 1$ to $k \geq 1$. The novelty relies on the proof of Lemma \ref{lemma:geometry of h_1 < h < h_2}, which is stated in \cite{VANDERHOUT2003} without proof and is actually non-trivial. This lemma studies the topology of level sets $h_1 < h < h_2$ for solutions $h, h_i$ of (\ref{heat map formulation for h(r,t)}) and shows that these level sets can be connected to the parabolic boundary of the square. The main difficulty follows from the fact that the boundary of a connected component $\Gamma$ of such a level set is not necessarily connected and both $h = h_1$ or $h = h_2$ could happen. To circumvent this problem, we rely on the notion of accessible boundary points from \cite{newman1964elements}, which are dense in the boundary. Even though the boundary is not connected, most boundary points can be reached from inside the domain. Hence, they can be connected together through some path inside the domain and this is the key element used to prove Lemma \ref{lemma:geometry of h_1 < h < h_2}.

We also extend the global existence result from Chang and Ding (\cite{ChangDing2011}) from $k = 1$ with a time-independent boundary data $v_0(x,t) = v_0(x)$ to $k \geq 1$ with a more general time-dependent boundary data $v_0(x,t)$ (Proposition \ref{prop:criterion for global solution}). We also show the existence of a solution which blows-up at infinity for any $k \geq 1$ (Proposition \ref{prop:blow up infinity global solution}), as a generalization of the result of \cite{ChangDing2011} from $k = 1$ to $k \geq 1$. The rigidity theorem from Van der Hout for $k = 1$ (\cite{VANDERHOUT2003}) relies on these results from Chang and Ding (\cite{ChangDing2011}), but they were not shown to be true in the new context of a time-dependent boundary data. Hence, we fill the gap and show that they remain true in our setting.

In addition, we strive to give a complete and modern presentation of classical $L^p$, parabolic Sobolev and Schauder estimates for which the main reference is the book from Ladyzhenskaia, Solonnikov, and Ural’tseva (\cite{ladyzhenskaia1968linear}) where some proofs (mainly Theorem \ref{thm:parabolic sobolev embedding}) are omitted. The standard references (e.g. \cite{krylov} for Schauder estimates) on that matter also assume some a priori regularity (e.g. a priori $C^{2+\alpha}$ regularity for Schauder estimates) of the solution, which we weaken so that, as expected, we obtain the higher regularity for ``free'' when dealing with parabolic equations. We also prove that, under the finite energy assumption, the concentration of energy argument from Struwe (\cite{struwe2008variational}), which is essentially local but was done on a manifold without boundary, works in the interior of the ball (which has a boundary) without imposing any boundary condition and with the absence of monotonicity for the energy. This is the content of Theorem \ref{thm:smoothness around all but finitely many points} which is to contrast with the analog work from Chang (\cite{ChangStruweAnalog}), where he assumes a time-independent boundary data $v_0(x,t) = v_0(x)$. We then prove that under the $k$-equivariant boundary data assumptions (which implies finiteness of the energy), there is no singularity at the boundary and interior singularities can only occur at the origin (Theorem \ref{thm:smoothness around all but finitely many points}, Proposition \ref{prop:smoothness at the boundary}). The work from Qing (\cite{Qing1995OnSO}), which is also local, finally allows getting a decomposition in terms of multi-bubbles along appropriate sequences (Theorem \ref{thm: full bubble decomposition}). No reference for the bubble-tree decomposition in the context of a time-dependent boundary data $v_0(x,t)$ could be found, so we have decided to include it in this work as it requires slightly adapting Struwe's work.

The generalization of the arguments from Van der Hout (\cite{VANDERHOUT2003}) and Chang, Ding (\cite{ChangDing2011}) holds because the Comparison Principle (Theorem \ref{comparison principle}), the intersection-number argument (Theorem \ref{lemma:intersection argument with chains mod 4}) and the different characterizations of blow-up (e.g. Proposition \ref{prop:criterion for global solution}) all generalize from $k = 1$ to $k \geq 1$. This is mostly due to the fact that the $k^2$ factor in the equation (\ref{heat map formulation for h(r,t)}) for $h(r,t)$ does not change the sign of the nonlinearity, gives a better decay $h(r,t) = r^k \tilde{h}(r,t)$ near $r = 0$ and that the family of stationary solutions given by the arctangent function still exists when $k \geq 1$.

The main theorems are the following:

\begin{theorem}[Nonexistence of multi-bubbles]\label{thm: main theorem, nonexistence bubble trees}
   Let $v(x,t)$ solve (\ref{heat map flow}) with smooth and $k$-equivariant boundary data, $k \geq 1$. Assume that $v$ blows-up at time $T < +\infty$. Consider a sequence $T_n \to T^-$, as well as finitely many smooth, non-constant, $k$-equivariant harmonic maps $\{\omega_1, ..., \omega_p\}$ from $\mathbb R^2 \cup \{\infty\} \simeq S^2$ to $S^2$ and positive sequences $(R_n^{(i)})_{n \geq 0}$, $i \in \{1,...,p\}$, all converging to zero, satisfying for all $i \neq j$,
    \begin{align}
        \frac{R_n^{(i)}}{R_n^{(j)}} + \frac{R_n^{(j)}}{R_n^{(i)}}  \to +\infty, \quad n \to +\infty,
    \end{align}
    and
 \begin{align}
     v(R_n^{(i)}x,T_n)  \to \omega_i(x), \quad n \to +\infty,
 \end{align}
    pointwisely almost everywhere on $\mathbb R^2$. Then $p = 1$.
    
    In particular, there is only one bubble (i.e., harmonic map) in the multi-bubble decomposition from Theorem \ref{thm: full bubble decomposition} and Corollary \ref{rmk:recovering the bubbles}.
\end{theorem}

\begin{proof}
    The full multi-bubble decomposition from Theorem \ref{thm: full bubble decomposition} implies the weaker formulation used in the above statement in the $k$-equivariant setting (see Corollary \ref{rmk:recovering the bubbles}).  The proof of Theorem \ref{thm: main theorem, nonexistence bubble trees} can be found in Corollary \ref{cor:end of proof}.
\end{proof}

\begin{theorem}[Global existence criterion]
    Let $v(x,t)$ solve (\ref{heat map flow}) on $[0,T)$ with smooth, $k$-equivariant boundary data, $k \geq 1$. Assume that the inclination coordinate $h_0$ of $v_0$ satisfies $h_0(0,t) = 0$, $|h_0(r,t)| \leq \pi$ on $\{0,1\} \times [0,T) \cup [0,1] \times \{0\}$. 
    
    If $T < +\infty$, then $v(x,t)$ and $h(r,t)$ cannot blow-up at time $T$. In particular, if $|h_0(r,t)| \leq \pi$ on $\{0,1\} \times [0,+\infty) \cup [0,1] \times \{0\}$, then $v(x,t)$ and $h(r,t)$ are global. 
\end{theorem}

\begin{proof}
    See Proposition \ref{prop:criterion for global solution}.
\end{proof}

\begin{remark}
    Chang, Ding and Ye (\cite{Chang1992FinitetimeBO}) have shown in the 1-equivariant setting that if $h_0(r,t) = h_0(r)$, $h_0(0) = 0$ and $|h_0(1)| > \pi$, then the solution $h(r,t)$ blows-up in finite-time.
\end{remark}

\begin{theorem}[Existence of $k$-equivariant solutions blowing-up at infinity]
    Let $k \geq 1$. Assume that $\psi = r^k\tilde{\psi}(r)$, $\tilde{\psi} \in C^{\infty}([0,1])$ with $\partial^{(2n+1)}_r \tilde{\psi}(0) = 0$ for all $n \in \mathbb N_{\geq 0}$, $0 \leq \psi \leq \pi$, $\tau(\psi)(1) = 0$ and $\tau(\psi)(r) \geq 0$ for all $r \in [0,1]$, where
    $$
\tau(\psi) := \psi_{rr} + \frac{1}{r}\psi_r -k^2 \frac{\sin(2\psi)}{2r^2}.
    $$
   Further assume there exists $r^* \in (0,1]$ for which $\psi(r^*) = \pi$. For example, one can choose
$$
\psi(r) = 4 \arctan \left( r^k \right).
$$
   Then the solution $h(r,t) \in C^{\infty}([0,1] \times (0,T))$ of (\ref{heat map formulation for h(r,t)}) with boundary data $h_0(r,t) = \psi(r)$ on $[0,1] \times \{0\} \cup \{0,1\} \times [0,+\infty)$ is global, has finite energy (\ref{Energy of v in terms of h}), satisfies $0 < h < \pi$ and $h_t > 0$ on $(0,1) \times (0,T)$, as well as $h_r > 0$ on $(0,1] \times (0,T)$. Moreover, $h(r,t)$ and the corresponding $k$-equivariant solution
   $$
   v(re^{i\theta},t) = \left( e^{ik \theta} \sin h(r,t), \cos h(r,t) \right)
   $$
   of (\ref{heat map flow})  blow-up at infinity in the sense that
   $$
\limsup \limits_{t \to +\infty}||h_r(\cdot,t)||_{C^1([0,r_0])} = \limsup \limits_{t \to +\infty}||\nabla v(\cdot,t)||_{C^1(B_{r_0}(0))} = +\infty \quad \forall r_0 \in (0,1].
$$
\end{theorem}

\begin{proof}
    See Theorem \ref{thm:existence of global solution} and Proposition \ref{prop:blow up infinity global solution}.
\end{proof}

\section{Outline of the proof}

This paper aims to be self-contained. That is why the third section goes back to the problem of existence and regularity for nonlinear parabolic problems. The aim of this section is, through a fixed point argument and a Schauder iteration, to deduce the regularity that one has for 'free' on a solution to the harmonic map heat flow. It is proved in Corollary \ref{cor: existence of solution for heat map flow} that a solution of (\ref{heat map flow}) has regularity
$$
    v(t,x) \in C^0([0,T_{\max}),C^1(\overline{B^2})) \cap C^{\infty}((0,T_{\max}) \times \overline{B^2}) 
$$
and that $\partial_t v, \partial_t \nabla_x v, \Delta v$ have additional Hölder-regularity at $t = 0$ (Proposition \ref{prop:holder continuity partial_t v}) even though we do not assume any kind of $m$-th order compatibility conditions (such as \cite[Chapter 7.1.4, Theorem 7]{evans10}). This is not trivial for nonlinear equations and those results only rely on the smoothness of the boundary data but not on its symmetry. 

In the fourth section, we add the assumptions of finite energy and symmetry of the boundary data. Under those assumptions, one obtains additional smoothness on $\overline{B^2} \setminus \{0\} \times \{T\}$ (Proposition \ref{prop:smoothness at the boundary}), thanks to the concentration of energy argument from Struwe (\cite{struwe2008variational}), as well as existence and rigidity for the limits $\lim \limits_{x \to 0} v(x,T), \lim \limits_{r \to 0^+} h(r,T)$ (Proposition \ref{value of v and h at origin}). Using the so-called ``barrier argument'', we also refine the description of blow-up. If $h(0,t) = 0$ for all $t \in [0,T)$, we prove that it cannot be that the boundary data $h_0$ is bounded by $\pi$ (Proposition \ref{prop:criterion for global solution}, which generalizes the result from \cite{ChangDing2011} to $k \geq 2$ and time-dependent boundary data) and one has 
$$
    \limsup_{\substack{r \to 0^+ \\ t \to T^-}}|h(r,t)| = \pi.
$$
As a bonus, we also generalize the other result from \cite{ChangDing2011} and show that there exists global solutions which blow-up at infinity for any $k \geq 1$.

In the fifth section, we recall the maximum principle and a comparison principle for the harmonic map heat flow problem. 

In the sixth section, following the lap-number argument from \cite{VANDERHOUT2003} and \cite{matano}, we are concerned on the level sets $h_1 < h_2$, as well as $h_1 < h < h_2$ for solutions $h, h_i$ of (\ref{heat map formulation for h(r,t)}). We show that those levels sets are connected to the parabolic boundary (Lemma \ref{geometry of h_1 < h_2} and Lemma \ref{lemma:geometry of h_1 < h < h_2}). Moreover, at a fixed time, there cannot be too many successive intersections between $h$ and $h_i$ on $[0,1]$ and this intersection number cannot increase over time (Lemma \ref{lemma:intersection argument with chains mod 4}). These topological results play a key role in the single bubble result, as one can obtain barriers from them restricting the development of bubbles.

Armed with this a priori knowledge on a solution which blows-up, we are prepared to prove the single bubble result in the seventh section. This knowledge, together with the comparison principle, imply (up to replacing $h$ by $-h$) that 
    $$
    \limsup_{\substack{r \to 0^+ \\ t \to T^-}} h(r,t) = \pi, \quad \lim \limits_{r \to 0^+} h(r,T) \in \{0,\pi\}.
    $$
    The difficult part is proving that $\lim \limits_{r \to 0^+} h(r,T) = \pi$, which follows by studying chains satisfying the property $P(t,h,h_1 = \pi - 2\arctan(\alpha^kr^k),h_2= \pi,M)$ from Lemma \ref{lemma:intersection argument with chains mod 4} when $t$ and $\alpha$ are large enough. A similar argument shows that the scale $(R_n)_{n \in \mathbb N}$ of any bubble lies on the left of the quantity 
  $$
 r^+(t) := \sup\{r \in [0,1]: h(s,t) \leq \pi \quad \forall s \in [0,r]\}
  $$
  for all $t$ large enough. After proving $\lim \limits_{r \to 0^+} h(r,T) = \pi$, one can assume without loss of generality that $h(r,T) > \pi/2$ on $(0,1]$ and $h(1,t) > \pi/2$ for $t \in [t_0,T]$. Another use of these chains and the Maximum Principle allows to exclude multi-bubbling.

Finally, the appendices build on top of each other and provide the aforementioned modern presentation of classical $L^p$, parabolic Sobolev and Schauder estimates, as well as the concentration of energy argument from Struwe.

\section{Local Existence and regularity} 
We first consider general nonlinear heat equations and prove a local existence and uniqueness result, before deducing existence of solutions for (\ref{heat map flow}). The goal of this section is to obtain a maximum of regularity for free on the solutions of the harmonic map heat flow, which is not trivial because the equation is nonlinear and it is not expected that the boundary data and the nonlinear term satisfy the usual compatibility conditions (e.g. \cite[Chapter 7.1.4, Theorem 7]{evans10}) required to get infinite differentiability.

We start with a classical fixed point argument to deduce the existence of a solution to nonlinear heat equations (Theorem \ref{thm:taylor local existence}). Combined with a Schauder iteration (Theorem \ref{thm:resgularty solution nonlinear pb}), one will deduce smoothness of the solution away from time zero. The fixed point argument is standard but we refine it to give a global existence result and a short-time estimate on the norm of the solution. This short-time estimate is used to prove that if $v$ is the solution to the harmonic map heat flow, then  $\partial_t v$, $\nabla_x \partial_t v$ and $\Delta v$ have some additional Hölder-regularity at time $t = 0$ (Proposition \ref{prop:holder continuity partial_t v}). This regularity is actually not needed for the bubble-tree result but only to show existence of a solution which blows-up at infinity (Proposition \ref{prop:blow up infinity global solution}).

Then, we turn to the particular structure of solutions under the $k$-equivariant symmetry assumption and show that the harmonic map heat flow problem reformulates as a one-dimensional nonlinear heat equation (Proposition \ref{prop:equivalent fomrulation heat map flow}).

\begin{notation}
    By $|D_x^kD_t^l u|$, we consider the Euclidean norm of the vector $(\partial_{x_1}^{k_1} ... \partial_{x_n}^{k_n}\partial_t^l u)_{k_1,...,k_n}$ for any indices $k_1 + ... + k_n = k$.

    For $x,t \in \mathbb R^n \times \mathbb R$ and $r > 0$, $C(x,t;r)$ denotes the open cylinder
    $$
    C(x,t;r) := \{(y,s) \in \mathbb R^n \times \mathbb R: |x-y| < r, t-r^2 < s < t\}.
    $$

    For $\alpha \in (0,1)$ and $u(x,t)$, $[u]_{\alpha,\alpha/2;U}$ denotes the Hölder semi-norm
\begin{align}
        \sup_{(x,t) \neq (y,s) \in U} \frac{|u(x,t)-u(y,s)|}{\left(|x-y| + |s-t|^{1/2}\right)^{\alpha}} &\sim \sup_{(x,t) \neq (y,t) \in U} \frac{|u(x,t)-u(y,t)|}{|x-y|^{\alpha}} \notag \\
    &+ \sup_{(y,t) \neq (y,s) \in U} \frac{|u(y,t)-u(y,s)|}{|s-t|^{\alpha/2}} \label{equivalence Holder norm}
\end{align}
and for $k,l \in \mathbb N_{\geq 0}$, $||u||_{k+\alpha,l+\alpha/2;U}$ denotes the norm
\begin{align}
    ||u||_{k+\alpha,l+\alpha/2;U} &= ||u||_{L^{\infty}_{x,t}(U)}+\sum_{i=1}^k ||D^i_x u||_{L^{\infty}_{x,t}(U)} + \sum_{j=1}^l  ||\partial^j_t u||_{L^{\infty}_{x,t}(U)} \notag \\
&+ [D^k_x u]_{\alpha,\alpha/2;U} + [\partial_t^l u]_{\alpha,\alpha/2;U}. \label{eq:holder norm}
\end{align}
Let $U \subset \mathbb R^n$ be any set (not necessarily open). When we write $u \in C^{k+\alpha,l+\alpha/2}(U)$ for $\alpha \in (0,1)$, we assume the norm (\ref{eq:holder norm}) to be finite. When we write $u \in C^{k,l}(U)$, we do not assume boundedness of $u$ and its derivatives, only continuity on $U$ (up to points $U \cap \partial U$). When we write $C_{loc}^{k+\alpha,l+\alpha/2}(U)$, it is to precise that we only assume (Hölder-)continuity on compact subsets $K \subset U$. Finally, $C^{k+\alpha,l+\alpha/2}_0(U)$ is the space of  $C^{k+\alpha,l+\alpha/2}(U)$ functions which vanishes at the boundary $\partial U$.
\end{notation}

\begin{theorem}[Local Existence of Solutions for nonlinear heat equation]\label{thm:taylor local existence}
    Consider the following general nonlinear heat equation 
\begin{align}
     V_t(x,t) &= \Delta V(x,t) + F(x,t,V(x,t),\nabla_x V(x,t)) , \quad (x,t) \in \Omega \times  (0,+\infty), \label{general nonlinear heat equation} \\
    V(x,t) &= 0, \quad (x,t) \in \partial \Omega \times [0,+\infty), \notag \\
    V(x,0) &= V_0(x), \quad x \in \Omega, \notag
\end{align}
where $\Omega \subset \mathbb R^n$ is a bounded, open, connected set with smooth boundary, $V: \overline{\Omega} \rightarrow \mathbb R^m$, the nonlinearity $F \in C^{0}(\overline{\Omega} \times \mathbb R_{\geq 0} \times \mathbb R^m \times \mathbb R^{n \times m};\mathbb R^m)$ is continuous, locally Lipschitz with respect to $(P,Q)$, i.e., 
    $$
    ||F(x,t,P_1,Q_1) - F(x,t,P_2,Q_2) ||_{L^{\infty}(\Omega)} \leq ||F||_{Lip,K} \cdot \left(|P_1-P_2| + |Q_1-Q_2|\right) 
    $$
    on any compact sets $K \subset \overline{\Omega} \times \mathbb R \times \mathbb R^m \times \mathbb R^{n \times m}$,
    and the initial data is $$
V_0 \in C^1_0(\overline{\Omega}) = \{f \in C^1(\overline{\Omega};\mathbb R^m): f_{|\partial \Omega} = 0\}.
$$
Then there exists a unique local solution
\begin{equation*}
    V(t,x) \in C^0([0,T_{\max}),C^1_0(\overline{\Omega}))
\end{equation*}
solving the fixed point problem
\begin{align}
    V(t) = e^{t\Delta} V_0 + \int_0^t e^{(t-s)\Delta} F(\cdot,s,V(\cdot,s),\nabla_x V(\cdot,s)) ds. \label{integral problem}
\end{align}
In fact, for all $0 < \alpha < 1$, the solution always has additional regularity
$$
V(x,t) \in C^{1+\alpha,\alpha/2}_{loc}(\Omega \times [0,T_{\max}) )
$$
(meaning it is Hölder-continuous on $K \times [0,T']$ for $T' < T$ and compact subsets $K \subset \Omega$) and if $V_0 = 0$, then one has regularity up to the boundary
     $$
     V(x,t) \in C^{1+\alpha,\alpha/2}_{loc}(\overline{\Omega} \times [0,T_{\max}) ).
     $$
      Moreover, if $T_{max} < +\infty$, then  
     $$
     \limsup \limits_{t \to T_{\max}}||V(\cdot,t)||_{C^1(\overline{\Omega})} = +\infty.
     $$
    Finally, if $F$ is globally bounded and $F(x,t,P,Q)$ is globally Lipschitz with respect to $(P,Q)$, i.e., 
    $$
    ||F(x,t,P_1,Q_1) - F(x,t,P_2,Q_2) ||_{L^{\infty}(\Omega)} \leq ||F||_{Lip} \cdot \left(|P_1-P_2| + |Q_1-Q_2|\right) 
    $$
    for some $||F||_{Lip}  < +\infty$, then the solution is global. Moreover, one has the short time estimate 
    $$
     \sup_{0 < t < T_0}||V(\cdot,t)||_{C^1(\overline{\Omega})} \leq 2 \sup \limits_{0 < t < 1}||e^{t \Delta}V_0(\cdot)-V_0(\cdot)||_{C^1(\overline{\Omega})} + 4C_{\frac{1}{2}}||V_0(\cdot)||_{C^1(\overline{\Omega})} + 1,
     $$
    where $2C_{\frac{1}{2}}T_0^{\frac{1}{2}}\max\{||F||_{\infty} + 1 , ||F||_{Lip} + 1\} = 1/2$ and $C_{\frac{1}{2}} = C_{\frac{1}{2}}(\Omega)$ is a constant depending on $\Omega$.
\end{theorem}

\begin{proof}
One can show (\cite[Chapter 8.1]{taylor2023partialvolI} and \cite[Chapter 13.7, (7.52), Proposition 7.4]{taylor2023partial}) that $e^{t\Delta}$ maps $L^{\infty}(\Omega)$ to $C^1_0(\overline{\Omega})$ for $t > 0$, is strongly continuous from $C^0_0(\overline{\Omega})$ to itself and from $C^1_0(\overline{\Omega})$ to itself. Moreover,
\begin{align*}
||e^{t \Delta} f||_{L^{\infty}(\Omega)} &\leq C_{0} ||f||_{L^{\infty}(\Omega)} \quad \forall t > 0, \\
||e^{t \Delta} f||_{C^1(\overline{\Omega})} &\leq C_{\frac{1}{2}} t^{-\frac{1}{2}} ||f||_{L^{\infty}(\Omega)} \quad \forall t > 0.
\end{align*}

Given any $\alpha > 0$, one can find $T_0 > 0$ small enough for which
$$
||e^{t\Delta} V_0(\cdot) - V_0(\cdot) ||_{C^1(\overline{\Omega})} \leq \alpha/2 \quad \forall t \in [0,T_0] 
$$
by continuity of the semigroup and then $K_1, K_2 > 0$ such that
\begin{align*}
    ||F(\cdot,t,V,\nabla V)||_{L^{\infty}(\Omega)} &\leq K_1 \quad \forall t \in [0,T_0], \\
    ||F(\cdot,t,V_1,\nabla V_1) - F(\cdot,t,V_2,\nabla V_2) ||_{L^{\infty}(\Omega)} &\leq K_2 ||V_1(\cdot,t) - V_2(\cdot,t)||_{C^1(\overline{\Omega})} \quad \forall t \in [0,T_0],
\end{align*}
for all
$$
V, V_1, V_2 \in Z(\alpha,T_0) = \{V \in C^0([0,T_0];C^1_0(\overline{\Omega})): V(x,0) = V_0(x), ||V(\cdot,t) - V_0(\cdot)||_{C^1(\overline{\Omega})} \leq \alpha\}
$$
by boundedness of $F$ and the local Lipschitz property around the closure of $\Omega \times [0,T_0] \times B_{\alpha + ||V_0||_{C^1(\overline{\Omega})}}(0) \times B_{\alpha + ||V_0||_{C^1(\overline{\Omega})}}(0)$.

Let $0 < T_1 \leq T_0$ be taken so that $2C_{\frac{1}{2}}T_1^{\frac{1}{2}}K_1 \leq \alpha/2$ and $2C_{\frac{1}{2}}T_1^{\frac{1}{2}}K_2 < 1$. Then the integral operator
$$
V \mapsto e^{t\Delta} V_0 + \int_0^t e^{(t-s)\Delta} F(\cdot,s,V(\cdot,s),\nabla_x V(\cdot,s)) ds
$$
is a contraction from $Z(\alpha,T_1)$ into itself. The Banach Fixed Point Theorem shows existence of a unique local solution $V \in Z(\alpha,T_1)$ to the corresponding integral problem (\ref{integral problem}). 

Let $T_{\max}$ be maximal with $V \in C^0([0,T_{\max}),C^1_0(\overline{\Omega}))$. We observe that $V$ has additional regularity
\begin{equation*}
    V(x,t) \in C^{1+\alpha,\alpha/2}_{loc}(\Omega \times [0,T_{\max}))
\end{equation*}
for all $0 < \alpha < 1$.

Indeed, fix any compact set $K \subset \Omega$ and any $0 < T' < T_{\max}$. We have that $V(s,y)$, $\nabla_x V(s,y)$ are bounded on $[0,T'] \times \overline{\Omega}$ thanks to our initial regularity given by the fixed point theorem. 

By continuity of $F$, $\partial_t V - \Delta V = F(y,s,V(y,s),\nabla_x V(y,s))$ is bounded as well on $[0,T'] \times \overline{\Omega}$. The parabolic Sobolev embedding (Corollary \ref{cor: local parabolic sobolev embedding}) implies that $V \in C^{1+\alpha,\alpha/2}(K \times [0,T'])$ for all $0 < \alpha < 1$. If $V_0 = 0$, one can apply Theorem \ref{thm:boundary parabolic sobolev embedding} instead to get the estimate up to the boundary.

If $T_{max} < +\infty$ is the maximal time of existence of the solution, then one must have 
$$
\limsup \limits_{t \to T_{\max}}||V(\cdot,t)||_{C^1(\overline{\Omega})} = +\infty.
$$
Otherwise, let $\alpha = \sup \limits_{0 < t < T_{\max}}||V(\cdot,t)||_{C^1(\overline{\Omega})} + 2 < +\infty$ and
$$
\tilde{K}_0 = ||F||_{L^{\infty}(\overline{\Omega} \times [0,2T_{\max}]  \times \overline{B_{\alpha}(0)} \times  \overline{B_{\alpha}(0)})}.
$$
Let $T_{\max}/2 < T^* < T_{\max}$. Observe that
\begin{align*}
||e^{t\Delta} V(\cdot,T^*)||_{C^1(\overline{\Omega})} &= ||e^{(T^*+t)\Delta} V_0(\cdot)
 + \int_0^{T^*} e^{(T^*+t-s)\Delta}F(\cdot,s,V(\cdot,s),\nabla_x V(\cdot,s))ds||_{C^1(\overline{\Omega})}  \\
 &\leq C_{\frac{1}{2}}\left(T^* + t\right)^{-\frac{1}{2}} \alpha  + \int_0^{T^*}C_{\frac{1}{2}}(T^* +t-s)^{-\frac{1}{2}} \tilde{K}_0 ds \\
 &\leq C_{\frac{1}{2}}\left(T^*\right)^{-\frac{1}{2}}\alpha + 2 \tilde{K}_0 C_{\frac{1}{2}} \left[\left(T^*+t\right)^{\frac{1}{2}} - t^{\frac{1}{2}} \right] \\
 &\leq C_{\frac{1}{2}}\left(T^*\right)^{-\frac{1}{2}}\alpha + 2 \tilde{K}_0 C_{\frac{1}{2}} \left(T^*\right)^{\frac{1}{2}} \\
 &\leq 2C_{\frac{1}{2}}\left(T_{\max}\right)^{-\frac{1}{2}}\alpha + 2 \tilde{K}_0 C_{\frac{1}{2}} \left(T_{\max}\right)^{\frac{1}{2}}, \quad t \geq 0, \\
 ||e^{t\Delta} V(\cdot,T^*) - V(\cdot,T^*) ||_{C^1(\overline{\Omega})} &\leq 2C_{\frac{1}{2}}\left(T_{\max}\right)^{-\frac{1}{2}}\alpha + 2 \tilde{K}_0 C_{\frac{1}{2}} \left(T_{\max}\right)^{\frac{1}{2}} + \alpha =: 2\beta, \quad t \geq 0.
\end{align*}
Then let 
    \begin{align*}
\tilde{K}_1 &= ||F||_{L^{\infty}(\overline{\Omega} \times [0,2T_{\max}]  \times \overline{B_{\alpha+\beta}(0)} \times  \overline{B_{\alpha+\beta}(0)})}, \\
\tilde{K}_2 &= ||F||_{Lip,\overline{\Omega} \times [0,2T_{\max}]  \times \overline{B_{\alpha+\beta}(0)} \times  \overline{B_{\alpha+\beta}(0)}},
    \end{align*}
    which are independent of $T^*$ and observe that for all $T_{\max}/2 < T^* < T_{\max}$ and $0 < t < T_{\max}$,
\begin{align*}
    ||F(\cdot,t+T^*,W,\nabla W)||_{L^{\infty}(\Omega)} &\leq \tilde{K}_1, \\
    ||F(\cdot,t+T^*,W_1,\nabla W_1) -F(\cdot,t+T^*,W_2,\nabla W_2) ||_{L^{\infty}(\Omega)} &\leq \tilde{K}_2 ||W_1(\cdot,t) - W_2(\cdot,t)||_{C^1(\overline{\Omega})},
\end{align*}
    on
\begin{align*}
    \tilde{Z}&(\beta,T^*,T_{\max}) \\
   & = \{W \in C^0([0,T_{\max}];C^1_0(\overline{\Omega})): W(x,0) = V(x,T^*), ||W(\cdot,t) - V(\cdot,T^*)||_{C^1(\overline{\Omega})} \leq \beta \}.
\end{align*}
Fix $0 < T_2 \leq T_{\max}$ so that $2C_{\frac{1}{2}}T_2^{\frac{1}{2}}\max\{\tilde{K}_1, \tilde{K}_2\} = 1/2$ and pick $T^*$ close enough to $T_{\max}$ so that $T^* + T_2 > T_{\max}$. Then solving the heat equation 
\begin{align*}
     W_t(x,t) &= \Delta W(x,t) + F(x,t+T^*,W(x,t),\nabla_x W(x,t)) , \quad (x,t) \in \Omega \times  (0,+\infty), \\
    W(x,t) &= 0, \quad (x,t) \in \partial \Omega \times [0,+\infty), \notag \\
    W(x,0) &= W(x,T^*), \quad x \in \Omega, \notag
\end{align*}
on $\tilde{Z}(\beta,T^*,T_2)$ allows extending $V$ past time $T_{\max}$, which contradicts its maximality.

Finally, if $F$ is globally bounded and globally Lipschitz, one can set
    \begin{align*}
        K_1 = ||F||_{L^{\infty}}+1, \quad K_2 = ||F||_{Lip}+1.
    \end{align*}
    Then for $2C_{\frac{1}{2}}T_0^{\frac{1}{2}}\max\{K_1, K_2\} = 1/2$ and any initial data $V_0 \in C^1_0$, set 
    $\alpha = 2\sup \limits_{0 < t < T_0}||e^{t \Delta}V_0(\cdot)-V_0(\cdot)||_{C^1(\overline{\Omega})} + 1$. The Banach Fixed Point then applies on $Z(\alpha, T_0)$. In particular, we have a minimal time of existence $T_0$ which is independent of the initial data. Hence, we can use $V(\cdot,T_0)$ as our new initial data to obtain a solution up to $2T_0$ and iterate this procedure. 
\end{proof}

\begin{remark}
    We remark that blow-up of the $C_x^1$-norm at infinity can occur if the solution is global (e.g. the global solution $\overline{h}$ constructed in Proposition \ref{prop:blow up infinity global solution}). 
\end{remark}

\begin{theorem}[Improved Regularity]\label{thm:resgularty solution nonlinear pb}
        Consider the solution $V(t,x)$ of (\ref{general nonlinear heat equation}). If the nonlinearity $F$ is smooth, then $V$ has regularity 
     $$
     V(t,x) \in C^{\infty}((0,T_{\max}) \times \Omega).
     $$
    Moreover, if $V_0 = 0$, then $V \in C^{\infty}((0,T_{\max}) \times \overline{\Omega})$ as well.
\end{theorem}

\begin{proof}
  Fix any cylinder $C(x,t;r) \subset \overline{C(x,t;r)} \subset (0,T_{\max}) \times \Omega$. As $V$ and $\nabla_x V$ are in $C^{\alpha,\alpha/2}_{x,t}(C(x,t;r))$ and $F$ is smooth, $F(y,s,V(y,s),\nabla_x V(y,s))$ is $C^{\alpha,\alpha/2}$-Hölder continuous as well on $C(x,t;r)$. The Schauder estimates (Theorem \ref{interior schauder estimates}) imply that $V(x,t) \in C^{2+\alpha,1+\alpha/2}(C(x,t;r/2))$. As $C(x,t;r)$ was arbitrary, $V(x,t) \in C^{2+\alpha,1+\alpha/2}_{loc}(\Omega \times (0,T_{\max}))$.
 In case $V_0 = 0$, we know that $V$ and $\nabla_x V$ are actually in $C^{\alpha,\alpha/2}_{loc}(\overline{\Omega} \times  [0,T_{\max}))$. The boundary Schauder estimates (Theorem \ref{boundary schauder smooth boundary}) show that $V(x,t) \in C^{2+\alpha,1+\alpha/2}_{x,t}(\Omega_{\varepsilon} \times   [T_1,T_2])$ for any $0 < T_1 < T_2 < T_{\max}$ and some $\varepsilon$-neighborhood $\Omega_{\varepsilon, T_1, T_2}$ of the boundary, hence  $V(x,t) \in C^{2+\alpha,1+\alpha/2}_{loc}(\overline{\Omega} \times  (0,T_{\max}))$.

    Finally, we do the so-called Schauder iteration to deduce higher regularity. Observe that $w = \partial_{x_k} V \in  C^{1+\alpha,\alpha/2}_{x,t}$ solves
    $$
    \partial_t w = \Delta w + \partial_{x_k} \left[ F(x,t,V,\nabla_x V) \right],
    $$
    where the forcing term is Hölder-continuous (either on a cylinder $C(x,t;r)$ or on $\overline{\Omega} \times [T_1,T_2]$ if $V_0 = 0$) by smoothness of $F$ and Hölder-continuity of $V, \nabla_x V$. Schauder estimates imply that $w, \nabla_x V \in  C^{2+\alpha,1+\alpha/2}_{loc}(\Omega \times (0,T_{\max}))$ and $w, \nabla_x V \in  C^{2+\alpha,1+\alpha/2}_{loc}(\overline{\Omega} \times (0,T_{\max}))$ if $V_0 = 0$.
    
    Next, $\tilde{w} = \partial_t V \in  C^{1+\alpha,\alpha/2}$ solves
    $$
    \partial_t \tilde{w} = \Delta \tilde{w} + \partial_{t} \left[ F(x,t,V,\nabla_x V) \right],
    $$
    and we can apply Schauder estimates again to deduce that $\partial_t V \in  C^{2+\alpha,1+\alpha/2}$. We carry on this procedure to obtain smoothness of $V$.
\end{proof}

\begin{corollary}[Existence of solutions for the Harmonic Map Heat flow]\label{cor: existence of solution for heat map flow}
For the harmonic map heat flow with smooth (not necessarily $k$-equivariant) boundary data, one can rewrite (\ref{heat map flow}) as a Dirichlet problem via $v(x,t) = v_0(x,t) + \tilde{v}(x,t)$, where $\tilde{v}(x,t)$ must solve
\begin{align}
     \tilde{v}_t &= \Delta \tilde{v} + F(x,t,\tilde{v},\nabla \tilde{v}) , \quad (x,t) \in B^2 \times  (0,+\infty), \label{heat map flow, dirichlet bound condition} \\
    \tilde{v} &= 0, \quad (x,t) \in \partial B^2 \times [0,+\infty) \cup B^2 \times \{0\}, \notag
\end{align}
and $F$ is a smooth nonlinearity
\begin{equation}
    F(x,t,\tilde{v},\nabla \tilde{v}) = - \partial_t v_0 + \Delta v_0 + |\nabla v_0|^2 (v_0 + \tilde{v}) + |\nabla \tilde{v} |^2 (v_0 + \tilde{v}) + 2 \nabla \tilde{v} \cdot \nabla v_0 (v_0 + \tilde{v}), \label{nonlinearity heat map flow}
\end{equation}
where $\nabla \tilde{v} \cdot \nabla v_0$ designates the dot product of the two (transposed) Jacobian matrices. We get a unique local solution of (\ref{heat map flow}) of regularity 
\begin{equation}
    v(t,x) \in C^0([0,T_{\max}),C^1(\overline{B^2})) \cap C^{\infty}((0,T_{\max}) \times \overline{B^2}), \label{Taylor existence of solution}
\end{equation}
 and if $T_{max} < +\infty$ is the maximal time of the solution, then one must have 
 $$
 \lim \limits_{t \to T_{\max}}||v(\cdot,t)||_{C^1(\overline{B^2})} = +\infty.
 $$
 One also has $v(x,t) \in C^{1+\alpha,\alpha/2}_{loc}(\overline{B^2} \times [0,T_{\max}))$ for any $0 < \alpha < 1$.
\end{corollary}

 \begin{proposition}[Hölder-Continuity of $\partial_t v$, $\Delta v$ at time zero]\label{prop:holder continuity partial_t v}
    Let
$$
     v(t,x) \in C^0([0,T_{\max}),C^1(\overline{B^2})) \cap C^{\infty}((0,T_{\max}) \times \overline{B^2}) \cap C_{loc}^{\alpha/2,1+\alpha}([0,T_{\max}) \times \overline{B^2})
$$
solve (\ref{heat map flow}) with smooth boundary data. Then $\partial_t v(x,t) \in C_{loc}^{1+\alpha,\alpha/2}(\overline{B^2} \times [0,T_{\max}))$, as well as $\Delta v(x,t) \in C_{loc}^{\alpha,\alpha/2}(\overline{B^2} \times [0,T_{\max}))$.
\end{proposition}

\begin{proof}
We take inspiration from \cite[Chapter V, Theorem 5.1]{ladyzhenskaia1968linear} by studying difference quotients. 

Write $T_{\max} = T$, $v = v_0 + \tilde{v}$ as in (\ref{heat map flow, dirichlet bound condition}). The nonlinearity can be rewritten as
$$
F(x,t,\tilde{v},\nabla \tilde{v}) = A(x,t) + B(x,t) \tilde{v} + \sum_{i,j} C_{i,j}(x,t) \partial_{x_i} \tilde{v}_j(x,t) \tilde{v} + \sum_{i,j}D_{i,j}(x,t) \left( \partial_{x_i} \tilde{v}_j(x,t) \right)^2 \tilde{v},
$$
where the coefficients are $C^{\infty}(\overline{B^2} \times [0,+\infty))$. Let $0 < h < T/4$ be fixed and $0 < t < T/2$. Consider the difference quotient $\tilde{v}^h(x,t) = h^{-1}[\tilde{v}(x,t+h)-\tilde{v}(x,t)]$. The goal is to bound the norm of this difference uniformly in $h$ in some parabolic Hölder space.

Note that 
$$A^h(x,t) = h^{-1} \left( A(x,t+h) - A(x,t) \right) = \int_0^1 (\partial_s A)(x,t+sh)ds \in C^{0}(\overline{B^2}\times [0,T/2])$$ with norm independent of $0 < h < T/4$. Similarly, one observes that
\begin{align*}
    &h^{-1} \left[ D_{i,j}(x,t+h) \left( \partial_{x_i} \tilde{v}_j(x,t+h) \right)^2 \tilde{v}(x,t+h) -  D_{i,j}(x,t) \left( \partial_{x_i} \tilde{v}_j(x,t) \right)^2 \tilde{v}(x,t)  \right] \\
    &=  h^{-1} \left[ D_{i,j}(x,t+h) \left( \partial_{x_i} \tilde{v}_j(x,t+h) \right)^2 \tilde{v}(x,t+h) - D_{i,j}(x,t+h) \left( \partial_{x_i} \tilde{v}_j(x,t+h) \right)^2 \tilde{v}(x,t) \right. \\
    &\phantom{=}  \left. +  D_{i,j}(x,t+h) \left( \partial_{x_i} \tilde{v}_j(x,t+h) \right)^2 \tilde{v}(x,t)  - D_{i,j}(x,t) \left( \partial_{x_i} \tilde{v}_j(x,t+h) \right)^2 \tilde{v}(x,t) \right. \\
    &\phantom{=} \left. + D_{i,j}(x,t)\left( \partial_{x_i} \tilde{v}_j(x,t+h) \right)^2 \tilde{v}(x,t) - D_{i,j}(x,t) \left( \partial_{x_i} \tilde{v}_j(x,t) \right)^2 \tilde{v}(x,t)  \right] \\
    &= D_{i,j}(x,t+h) \left( \partial_{x_i} \tilde{v}_j(x,t+h) \right)^2 \tilde{v}^h(x,t) + D_{i,j}^h(x,t) \left( \partial_{x_i} \tilde{v}_j(x,t+h) \right)^2 \tilde{v}(x,t) \\
    &\phantom{=}+ D_{i,j}(x,t)\partial_{x_i} \tilde{v}_j^h(x,t+h)  \left( \partial_{x_i} \tilde{v}_j(x,t+h) +\partial_{x_i} \tilde{v}_j(x,t) \right) \tilde{v}(x,t) \\
    &= D_{i,j,0}(x,t) + D_{i,j,1}(x,t) \tilde{v}^h(x,t) + D_{i,j,2}(x,t)  \partial_{x_i} \tilde{v}_j^h(x,t),
\end{align*}
where the coefficients are $C^0(\overline{B^2} \times [0,T/2])$ with norm independent of $0 < h < T/4$ given the regularity $v \in C^0([0,T), C^1(\overline{B^2}))$. We proceed similarly for the difference involving the coefficients $B$ and $C$. In other words, $\tilde{v}^h(x,t)$ is the unique local solution to a linear parabolic equation
\begin{align*}
    \partial_t \tilde{v}^h &=   \Delta  \tilde{v}^h + G_{h,0}(x,t) \tilde{v}^h(x,t) + \sum_{i,j} G_{h,i,j}(x,t) \partial_{x_i} \tilde{v}^h_j(x,t), \quad x \in B^2, 0 < t < T/2, \\
    \tilde{v}^h(x,t) &= 0, \quad (x,t) \in B^2 \times \{0\} \cup \partial B^2 \times [0,T/2],
\end{align*}
where $G_{h,0}(t,x), G_{h,i,j}(t,x) \in C^0([0,T/2] \times \overline{B^2})$ are coefficients (which we can extend as constant for $t \geq T/2$) with norm independent of $0 < h < T/4$.

It follows from Theorem \ref{thm:taylor local existence} that there is a small time $0 <T_0 < T/2$, independent of $h$, for which $v^h \in C^0([0,T_0];C^1(\overline{B^2}))$ with
$$
||\tilde{v}^h(\cdot,t)||_{C^1(\overline{B^2})} \leq 1, \quad \forall 0 < t < T_0,
$$
and the parabolic Sobolev embedding (Theorem \ref{thm:boundary parabolic sobolev embedding}) then implies
$$
||\tilde{v}^h||_{C^{1+2\gamma,\gamma}(\overline{B^2} \times [0,T_0])} \leq C(v_0,\gamma)
$$
for any $0 < \gamma < 1/2$, for some constant $C(v_0,\gamma)$ which is independent of $h$. Take any sequence $h_k \to 0^+$. One has $\lim \limits_{k \to +\infty}\tilde{v}^{h_k}(x,t) \to \partial_t \tilde{v}(x,t)$ on $\overline{B^2} \times (0,T_0]$ using our a priori smoothness away from $t = 0$. Moreover, by Arzela-Ascoli, the convergence is uniform on $\overline{B^2} \times [0,T_0]$, meaning that the limit $\partial_t \tilde{v}(x,t)$ extends continuously to $t = 0$. Similarly, $\partial_t \nabla_x \tilde{v}$ extends continuously at $t = 0$. Moreover, the limits $\partial_t v$ and $\partial_t \nabla_x v$ must be $C^{2\gamma, \gamma}$ Hölder-continuous as well. 

As $v \in C_{loc}^{1+\alpha,\alpha/2}([0,T_{\max}) \times \overline{B^2})$, it follows from the equation that $\Delta v \in  C_{loc}^{\alpha,\alpha/2}([0,T_{\max}) \times \overline{B^2})$.
\end{proof}

It turns out that if $v_0(x,t)$ is $k$-equivariant, then so is $v(x,t)$. First, we need a lemma about the structure of $k$-equivariant functions.

\begin{lemma}[Smoothness of inclination coordinate] \label{k-equivariance, smoothness of h}
     Let $I$ be any kind of interval and $R > 0$. If $v_0 \in C^{\infty}(\overline{B_R(0)} \times I, S^2)$ is smooth and $k$-equivariant in the sense that
     \begin{align}
         v_0(re^{i\theta},t) = \left( e^{ik\theta} \sin h_0(r,t), \cos h_0(r,t) \right) \in S^2 \label{k equivariant form of v_0}
     \end{align}
     for some function $h_0: [0,R] \times I \rightarrow \mathbb R$. Then one can always find some $h_0$ for which (\ref{k equivariant form of v_0}) holds and $h_0(r,t) = r^k \tilde{h}_0(r,t)$, where $\tilde{h}_0 \in C^{\infty}([0,R] \times I,\mathbb R)$ and $\partial_r^{(2n+1)}\tilde{h}_0(0,t) = 0$ for all $n \in \mathbb N_{\geq 0}$.

     Conversely, if $h_0$ satisfies the above properties, then $v_0$ defined as in (\ref{k equivariant form of v_0}) is smooth, $v_0 \in C^{\infty}(\overline{B_R(0)} \times I, S^2)$.
\end{lemma}

\begin{remark}
    In the direction where $v_0$ is smooth and $k$-equivariant, no smoothness is a priori assumed from $h_0(r,t)$. The lemma shows that $h_0$ can be replaced by a smooth representative.
    
    The representative is not unique. Letting $\theta = 0$, one observes that two continuous inclination coordinates for $v_0$ must differ by a multiple of $\pi$. Hence, we will always choose the one with $h_0(0,0) = 0$.
\end{remark}

\begin{proof}

First, assume that $v_0$ is smooth and $k$-equivariant.

Using a rescaling of the space variable, one can assume that $B_R(0) = B_1(0) = B^2$. Write $v_0 = (v_{0,1}, v_{0,2}, v_{0,3}) \in S^2$.

    It suffices to verify that $F(r,t) = \sin h_0(r,t)$ is smooth, $F(r,t) = r^k \tilde{F}(r,t)$ with $\tilde{F}$ smooth, as well as $\partial_r^{(2n+1)}\tilde{F}_0(0,t) = 0$ for all $n \in \mathbb N$.  

    If this is proved, then it is a classical result from topology that one can lift the smooth map $\tilde{v}_0(r,t) = (\sin h_0(r,t), \cos h_0(r,t)) = (v_{0,1}(re^{i0},t), v_{0,3}(re^{i0},t)) \in S^1$ to a continuous map $h_0(r,t): [0,R] \times I \rightarrow \mathbb R$ satisfying $\tilde{v}_0 = \pi \circ h_0$, where $\pi(\theta) = (\sin \theta, \cos \theta)$. Then $h_0$ is smooth around any points where $F(r,t) = \sin h_0(r,t) \neq \pm 1$ by the inverse function theorem, and the expansion of $h_0$ at $(0,t)$ follows from the ones of $F$ and $\arcsin(\cdot)$. At points where $F(r,t) = \pm 1$, one has $\cos h_0(r,t) = 0$ and the smoothness of $h_0(r,t)$ is obtained from the smoothness of $G(r,t) = \cos h_0(r,t) = v_{0,3}$.

    Since 
    $$
\underbrace{(v_{0,1}(x,y,t),v_{0,2}(x,y,t))}_{=: \phi(x,y,t)} = F(r,t)e^{ik\theta} = F\left(\sqrt{x^2+y^2},t\right) \frac{(x,y)}{\sqrt{x^2+y^2}}
    $$
    is smooth on $\overline{B^2} \times I$, setting $x = r$, $y = 0$ proves the smoothness of $F(r,t)$ on $[0,R] \times I$. Moreover, considering $\phi(1/n,0,t)$ and $\phi(0,1/n,t)$ with $n \to +\infty$ shows that $F(0,t) = 0$. In the following, we fix $t \in I$ and omit it from the notation $\phi(x,y,t)$ and $F(r,t)$. 
    
    Now, we prove that $F(r) = \mathcal{O}(r^k)$ near $r = 0$. Indeed,
    $$
    \Delta_{(x,y)}^j \phi = \left( \partial_{rr} + \frac{1}{r}\partial_r - \frac{1}{r^2}\partial_{\theta\theta} \right)^j \phi = e^{ik\theta} \left( \partial_{rr} + \frac{1}{r}\partial_r - \frac{k^2}{r^2}\cdot  \right)^j F 
    $$
    is smooth on $\mathbb R^2$ for all $j \in \mathbb N$. In particular, the radial part 
    $$
    \left( \partial_{rr} + \frac{1}{r}\partial_r - \frac{k^2}{r^2}\cdot  \right)^j F(r)
    $$
    is smooth on $[0,1]$ and zero at $r = 0$.
    
    On a sum $\sum_{n=n_0}^{N}a_nr^{n-i}$, the operator $\left( \partial_{rr} + \frac{1}{r}\partial_r - \frac{k^2}{r^2}\cdot  \right)^j$  acts as
    $$
\left( \partial_{rr} + \frac{1}{r}\partial_r - \frac{k^2}{r^2}\cdot  \right)^j \sum_{n=n_0}^{N}a_nr^{n-i} = \sum_{n=n_0}^{N}\prod_{l=0}^{j-1}\left[(n-i-2l)^2 - k^2\right] a_n r^{n-i-2j}.
    $$    
    Write $F(r)$ with a one-sided Taylor expansion at $r = 0$, $F(r) = \sum_{n=1}^Na_nr^n + \mathcal{O}(r^{N+1})$ where $N > 2k$ is fixed and the remainder is given by the integral formula (from which it is seen that $\partial_{r}^i \mathcal{O}(r^{N+1}) = \mathcal{O}(r^{N+1-i})$ for $i \leq N$). 

    Assume $a_0, ..., a_{2j} = 0$ for $j \geq 0$ and $2j < k-1$. We prove that $a_{2j+1} = 0$ and $a_{2j+2} = 0$ as well if $2j+2 < k$. In other words, we prove by induction that $a_0 = a_1 = ... = a_{k-1} = 0$. Since
\begin{align*}
    \left( \partial_{rr} + \frac{1}{r}\partial_r - \frac{k^2}{r^2}\cdot  \right)^{j+1} F(r) &= \sum_{n=1}^{N} \prod_{l=0}^{j}\left[(n-2l)^2 - k^2\right] a_n r^{n-2(j+1)} + \mathcal{O}(r^{k+1}) \\
    &= \sum_{n=2j+1}^{N} \prod_{l=0}^{j}\left[(n-2l)^2 - k^2\right] a_n r^{n-2(j+1)} + \mathcal{O}(r^{k+1})
\end{align*}
is smooth on $[0,1]$, the first term 
$$
 \prod_{l=0}^{j}\left[(2j+1-2l)^2 - k^2\right] a_{2j+1} r^{-1}
$$
is smooth as well on $[0,1]$. Since $1 \leq 2j+1-2l < k$, the only possibility is that $a_{2j+1} = 0$. Moreover, if $2j+2 < k$, then 
\begin{align*}
    \left( \partial_{rr} + \frac{1}{r}\partial_r - \frac{k^2}{r^2}\cdot  \right)^{j+1} F(r) 
\end{align*}
is zero at $r = 0$. Hence,
$$
 \prod_{l=0}^{j}\left[(2j+2-2l)^2 - k^2\right] a_{2j+2}r^0 = 0,
$$
meaning that $a_{2j+2} = 0$ as well.

We have proved that $F(r) = r^k \tilde{F}(r)$, where $\tilde{F}(r)$ is smooth. We show that $\tilde{F}(r)$ must be even. Fix $J > 2k$ and $N > 2J$. Write $F(r) = r^k \sum_{n=0}^N a_nr^n + \mathcal{O}(r^{N+1})$. We prove that $a_n = 0$ for any odd $n \leq J$. Indeed,
\begin{align*}
    \left( \partial_{rr} + \frac{1}{r}\partial_r - \frac{k^2}{r^2}\cdot  \right)^{J+1} F(r) &= \sum_{n=0}^{N} \prod_{l=0}^{J}\left[(n+k-2l)^2 - k^2\right] a_n r^{n+k-2(J+1)} + \mathcal{O}(r^{J+1})
\end{align*}
is smooth on $[0,1]$. If $n \leq J$ is odd, then $k+n-2(J+1) \leq -1$,
$$
\prod_{l=0}^{J}\left[(n+k-2l)^2 - k^2\right] a_n r^{n+k-2(J+1)}
$$
is smooth on $[0,1]$ and $\prod_{l=0}^{J}\left[(n+k-2l)^2 - k^2\right] \neq 0$, meaning that $a_n = 0$.

Conversely, assume that $h_0 = r^k \tilde{h}_0$ is smooth with $\partial_{r}^{(2n+1)}\tilde{h}_0(0,t) = 0$ for $n \geq 0$. Let 
$$
v_0(z,t) = \left( \frac{z^k}{|z|^k} \sin h_0(|z|,t), \cos h_0(|z|,t) \right), \quad z = (x,y) \in \overline{B_R(0)} \subset \mathbb C,
$$
which is clearly smooth away from $z = 0$. It suffices to verify the $\mathbb R^2$-smoothness of the radial functions
$$
z \mapsto \frac{\sin h_0(|z|,t)}{|z|^k}, \quad z \mapsto \cos h_0(|z|,t)
$$
near $z = 0$, which is the case if and only if the odd (one-sided) $r$-derivatives of
$$
\frac{\sin h_0(r,t)}{r^k} = \frac{\sin \left( r^k\tilde{h}_0(r,t) \right)}{r^k} = \frac{\sin \left( r^k\tilde{h}_0(r,t) \right)}{r^k \tilde{h}_0(r,t)}  \tilde{h}_0(r,t) , \quad \cos h_0(r,t) = \cos \left( r^k\tilde{h}_0(r,t) \right)
$$
vanish at zero. Extend $\tilde{h}_0$ for $r < 0$ by even reflection $\tilde{h}_0(-r,t) = \tilde{h}(r,t)$ for $r < 0$. Then $\tilde{h}_0 \in C^{\infty}([-R,R] \times I,\mathbb R)$ is even with respect to $r$ and it also follows that 
$$
\frac{\sin \left( r^k\tilde{h}_0(r,t) \right)}{r^k \tilde{h}_0(r,t)}  \tilde{h}_0(r,t) , \quad  \cos \left( r^k\tilde{h}_0(r,t) \right)
$$
extend as $C^{\infty}([-R,R] \times I,\mathbb R)$ and even functions, which concludes the proof.
\end{proof}

\begin{proposition}[Solutions of harmonic map heat flow with $k$-equivariant boundary data] \label{prop:equivalent fomrulation heat map flow}
Let $h_0(r,t) = r^k \tilde{h}_0(r,t)$, where $\tilde{h}_0 \in C^{\infty}([0,1] \times [0,+\infty),\mathbb R)$ and $\partial_r^{(2n+1)}\tilde{h}_0(0,t) = 0$ for all $n \in \mathbb N_{\geq 0}$. There exists a unique solution
$$
h(r,t) = r^k \tilde{h} + h_0, \quad \tilde{h}(t,r) \in C^0([0,T),C^1([0,1])) \cap C^{\infty}((0,T) \times [0,1]),
$$
which solves 
\begin{align}
h_t &= h_{rr} + \frac{h_r}{r} - k^2 \frac{\sin(2h)}{2r^2}, \quad (r,t) \in (0,1) \times (0,+\infty), \tag{\ref{heat map formulation for h(r,t)}} \\
h(r,t) &= h_0(r,t), \quad (r,t) \in \{0,1\} \times [0,+\infty) \cup (0,1) \times \{0\}. \notag
\end{align}
Moreover, 
$$
     v(t,x) \in C^0([0,T),C^1(\overline{B^2})) \cap C^{\infty}((0,T) \times \overline{B^2}) 
    $$
defined as 
\begin{equation} \label{k-equivariant form of v(x,t)}
    v(re^{i\theta},t) = \left( e^{ik\theta} \sin h(r,t), \cos h(r,t) \right)
\end{equation}
solves (\ref{heat map flow}) with smooth and $k$-equivariant data
\begin{equation*}
    v_0(re^{i\theta},t) = \left( e^{ik\theta} \sin h_0(r,t), \cos h_0(r,t) \right).
\end{equation*}
    Conversely, solutions 
    $$
     v(t,x) \in C^0([0,T),C^1(\overline{B^2})) \cap C^{\infty}((0,T) \times \overline{B^2}) 
    $$
    to equation (\ref{heat map flow}) with smooth and $k$-equivariant boundary data $v_0$ are exactly of the form:
\begin{equation*} 
    v(re^{i\theta},t) = \left( e^{ik\theta} \sin h(r,t), \cos h(r,t) \right),
\end{equation*}
where 
$$
h(t,r) \in C^0([0,T),C^1([0,1])) \cap C^{\infty}((0,T) \times [0,1])
$$
solves (\ref{heat map formulation for h(r,t)}) with boundary data $h_0$, where $h_0$ is the inclination coordinate of the boundary data $v_0$, chosen of the form (see Lemma \ref{k-equivariance, smoothness of h}) $h_0(r,t) = r^k \tilde{h}_0(r,t)$, where $\tilde{h}_0 \in C^{\infty}([0,1] \times [0,+\infty),\mathbb R)$ and $\partial_r^{(2n+1)}\tilde{h}_0(0,t) = 0$ for all $n \in \mathbb N_{\geq 0}$.

In particular, $v$ and $h$ have the same maximal time of existence $T_{\max} \leq +\infty$.

Finally, if $T \leq +\infty$ and $m \in \mathbb N$ is such that $||h_0||_{L^{\infty}_{r,t}(\Sigma_T)} \leq m \pi$ on the parabolic boundary $\Sigma_T = \{0,1\} \times [0,T] \cup [0,1] \times \{0\}$, then $|h(r,t)| \leq m\pi$ on $[0,1] \times (0,T)$ as well.
\end{proposition}

\begin{remark}\label{remark:uniqueness of inclination coordinate}
    It is possible to solve equation (\ref{prop:equivalent fomrulation heat map flow}) even if the boundary data $h_0(r,t)$ is not of the form $h_0 = r^k \tilde{h}_0$.

     Indeed, if $h$ solves (\ref{prop:equivalent fomrulation heat map flow}) (without caring about the boundary data), then so does $h + m \pi, m \in \mathbb Z$. Hence, the equation can be solved with boundary data $h_0 = r^k \tilde{h}_0 + m\pi$ as well. This is expected as continuous inclination coordinates are unique up to a multiple of $\pi$. Hence, the multiple of $\pi$ does not change $v$ or $v_0$.
     
     However, the constant solutions $\pi/2 + m \pi$, $m \in \mathbb Z$, also solve (\ref{prop:equivalent fomrulation heat map flow}) with a boundary data which is not of the previous form $h_0 = r^k \tilde{h}_0 + m\pi$. This is consistent as $\pi/2$ corresponds to a solution
     $$
    v(x,t) = (e^{ik\theta},0) = \left( \frac{x}{|x|}, 0 \right)
     $$
     of the harmonic map heat flow for which $v(x,0) = v_0(x,0)$ is not smooth at the origin.
\end{remark}

\begin{proof}
Looking for a solution of the form $h = r^k \tilde{h} + h_0$, $\tilde{h}(r,t)$ must solve
\begin{align*}
    \tilde{h}_t &= \tilde{h}_{rr} + \frac{2k+1}{r}\tilde{h}_r - G(r^2,\tilde{h}) + H_0, \quad (r,t) \in (0,1) \times [0,+\infty), \\
    \tilde{h}(r,t) &= 0, \quad (r,t) \in \{0,1\} \times [0,+\infty) \cup (0,1) \times \{0\},
\end{align*}
where, writing $h_0 = r^k \tilde{h}_0$ as in Lemma \ref{k-equivariance, smoothness of h},
\begin{align*}
    G(r^2,\tilde{h}) &= r^{-k} \left( \frac{\sin(2h) - 2h - \sin(2h_0) + 2h_0}{2r^2} \right) = r^{-k-2} \sum_{n=1}^{\infty} \frac{(2r^k\tilde{h} + 2h_0)^{2n+1}-(2h_0)^{2n+1}}{2(2n+1)!} \\
    &= r^{-k-2} \sum_{n=1}^{\infty} \frac{2r^k\tilde{h} \sum_{j=0}^{2n} (2r^k \tilde{h})^j (2h_0)^{2n-j}}{2(2n+1)!} = \sum_{n=1}^{\infty} \frac{4^n \tilde{h} (r^2)^{kn -1} \sum_{j=0}^{2n}\tilde{h}^j \tilde{h}_0^{2n-j}}{(2n+1)!}
\end{align*}
is a smooth nonlinearity and 
$$
H_0(r,t) =\left[  -h_{0,t} + h_{0,rr} + \frac{h_{0,r}}{r} - k^2 \frac{\sin(2h_0)}{2r^2} + \frac{h_0}{r^2}\right] r^{-k} 
$$
is smooth as well given the regularity of $h_0$ obtained in Lemma \ref{k-equivariance, smoothness of h}. The PDE for $\tilde{h}(r,t)$ is a radial heat equation in dimension $2k+2$ with a smooth nonlinearity. It follows from Theorem \ref{thm:taylor local existence} that there exists a unique smooth solution
$$
\tilde{h}(t,x) \in C^0([0,T),C^1(\overline{B^2})) \cap C^{\infty}((0,T) \times \overline{B^2})
$$ for (\ref{heat map formulation for h(r,t)}). Moreover, if $T < +\infty$, then 
$$\lim \limits_{t \to T}||\tilde{h}(\cdot,t)||_{C^1(\overline{B^2})} = +\infty.$$

Such a solution must be radial. This is the case because $h$ is obtained as the limit of a Banach Fixed Point Iteration for the problem (\ref{integral problem}). Initializing the iteration with the radial initial data, one must check that at each iteration, the approximate solution is still radial. This is the case since the nonlinearity is radial and the heat semi-group maps radial functions to radial functions (one can solve the Dirichlet problem for the heat equation on $B^2$ in radial coordinates using separation of variables and Sturm-Liouville theory). Hence,
$$
\tilde{h}(t,r) \in C^0([0,T),C^1([0,1])) \cap C^{\infty}((0,T) \times [0,1])
$$
is radial and we get the same regularity and blow-up criterion for $h$.

If $|h_0| \leq m\pi$ on the parabolic boundary of $(0,1) \times (0,T)$ by smoothness of the boundary data, then the Comparison Principle (Theorem \ref{comparison principle}) shows that $|h(r,t)| \leq m \pi$ is bounded on $[0,1] \times (0,T')$ for any $0 < T' < T$ since $m\pi$ is always a finite energy solution.

If $h$ solves (\ref{heat map formulation for h(r,t)}) on $[0,T)$, letting
$$
v(re^{i\theta},t) = \left( e^{ik\theta} \sin h(r,t), \cos h(r,t) \right),
$$
one verifies directly that $v$ is smooth (Proposition \ref{k-equivariance, smoothness of h}) and solves the harmonic map heat flow problem (\ref{heat map flow}) on $[0,T)$. 

In the converse direction, one starts with a solution $v$ of (\ref{heat map flow}) on $[0,T)$. Solving (\ref{heat map formulation for h(r,t)}) with the inclination coordinate $h_0$ of the data $v_0$ yields an inclination coordinate $h$ on $[0,1] \times [0,T_{h,\max})$ which we use to form another solution $\tilde{v}$ of (\ref{heat map flow}) defined up to time $T_{h,\max}$. If $T_{h,\max} < T \leq +\infty$, then $\tilde{v}$ must blow-up at $T_{h,\max}$. Indeed, as seen above, the Comparison Principle (Theorem \ref{comparison principle}) shows that $h(r,t)$ is bounded on $[0,1] \times (0,T_{h,\max})$ by a large multiple $m \pi$. Hence, only the derivative $h_r$ blows-up at time $T_{h,\max}$. Yet,
\begin{equation} 
    |\nabla \tilde{v}|^2 = h_r^2 + \frac{k^2}{r^2} \sin(h)^2, \label{nabla v in terms of h}
\end{equation}
meaning that $\nabla \tilde{v}$ must also blow-up at time $T_{h,\max}$ as well. By uniqueness of the problem (\ref{heat map flow}), $\tilde{v} = v$ up to $T_{h,\max} < T$, meaning that $v$ blows-up before time $T$, a contradiction. Hence, $T \leq T_{h,\max}$, $h$ exists up to time $T$ and $v = \tilde{v}$ on $[0,T)$ by uniqueness.
\end{proof}

\begin{proposition}[Continuity of $h_t$ and $h_{rr}$]\label{prop:continuity of h_t}
Let
$$
     v(t,x) \in C^0([0,T),C^1(\overline{B^2})) \cap C^{\infty}((0,T) \times \overline{B^2}) 
$$
solve (\ref{heat map flow}) with smooth and $k$-equivariant boundary data. Let 
$$
h(t,r) \in C^0([0,T),C^1([0,1])) \cap C^{\infty}((0,T) \times [0,1])
$$
be the corresponding inclination coordinate (\ref{heat map formulation for h(r,t)}). Then $h_t(t,r) \in C^0([0,T) \times [0,1])$ and $h_{rr} \in C^0([0,T) \times (0,1])$.
\end{proposition}

\begin{proof}
    Continuity of $h$ and $h_r$ on $[0,T) \times [0,1]$ is already known. Observe that $v_t(t,re^{i 0}) \in C^0([0,T) \times [0,1])$ by Proposition \ref{prop:holder continuity partial_t v} and
    $$
    v_t(t,re^{i 0}) = \left( \cos(h) h_t,\cos(h) h_t, - \sin(h) h_t \right).
    $$
    As $\cos(h), \sin(h) \in C^0([0,T) \times [0,1])$, continuity of $h_t$ follows. Similarly, one has the regularity $\Delta v(t,re^{i 0}) \in C^0([0,T) \times [0,1])$ and
    $$
   \Delta v(t,re^{i 0}) = \begin{pmatrix}
       \cos(h) h_{rr} - \sin(h) h_r^2 + \cos(h) h_r r^{-1} - k^2 \sin(h)r^{-2} \\
       \cos(h) h_{rr} - \sin(h) h_r^2 + \cos(h) h_r r^{-1} - k^2 \sin(h)r^{-2} \\
       - \sin(h) h_{rr} - \cos(h) h_r^2 - \sin(h) h_r r^{-1}  
   \end{pmatrix},
    $$
    from which continuity away from $r = 0$ follows.
\end{proof}

We finish with an estimate of the blow-up rate of $h$.

\begin{proposition}[Blow-Up of $k$-equivariant solutions] \label{blow up point is at zero k equivariant}
Let $T < +\infty$ and assume
$$
     v(t,x) \in C^0([0,T),C^1(\overline{B^2})) \cap C^{\infty}((0,T) \times \overline{B^2}) 
$$
solve (\ref{heat map flow}) with smooth and $k$-equivariant boundary. Let 
$$
h(t,r) \in C^0([0,T),C^1([0,1])) \cap C^{\infty}((0,T) \times [0,1])
$$
be the corresponding inclination coordinate (\ref{heat map formulation for h(r,t)}). For all $n \in \{0,1,2\}$, there exists $C = C(k^2,n,h) > 0$ for which
\begin{align*}
    |\partial_r^{(n)} h(\lambda,t)| \leq ||h_{r}||_{L^{\infty}_{r,t}([\lambda,1] \times [T/2,T])} &\leq \frac{\tilde{C}(k^2,n,h)}{\lambda^n} \quad \forall \lambda \in (0,1], t \in [T/2,T].
\end{align*}
In particular, if $v$ blows-up at time $T$, then the blow-up of $h_r$ happens exactly at $r = 0$, i.e., 
$$
\limsup_{\substack{r \to 0^+ \\ t \to T^-}} |h_r(r,t)| = +\infty.
$$
\end{proposition}

\begin{proof}
The inclination coordinates solves a 2D radial heat equation with smooth and bounded forcing term on $[1/2,1] \times [T/2,T]$. The estimate follows from the boundary Schauder estimates (Theorem \ref{boundary schauder smooth boundary}) on that region. Hence, it suffices to prove the statement for all $\lambda$ small enough. Let $\lambda \in (0,1)$ and 
$$
h_{\lambda}(r,t) = h( \lambda r,\lambda^2 t) \in C^{\infty}( (0,\lambda^{-1}) \times (0,\lambda^{-2}T) ).
$$
It follows from the Comparison Principle (Proposition \ref{prop:equivalent fomrulation heat map flow}) that $h$, hence $h_{\lambda}$, is bounded by a large multiple $m\pi$.

Moreover, $h_{\lambda}$ solves the radial 2D nonlinear heat equation
$$
h_t = h_{rr} + \frac{h_r}{r} - k^2 \frac{\sin(2h)}{2r^2}, \quad (r,t) \in  (1/2,\lambda^{-1}) \times (0,\lambda^{-2}T)
$$
with a bounded forcing term. For any open cylinder 
$$
C(x,t;1)  = B_1(x) \times (t-1,t) \subset \{(x,t) \in \mathbb R^2 \times \mathbb R_{\geq 0}:   (|x|,t) \in (3/2,\lambda^{-1}-1) \times (0,\lambda^{-2}T)\},
$$ 
the Schauder estimates (Theorem \ref{interior schauder estimates}, Remark \ref{weakening of interior Schauder estimates}) show the existence of $\gamma \in (0,1)$ and a constant $C = C(n)$ for which 
\begin{align*}
    ||h_{\lambda}||_{C^{2+2\gamma,1+\gamma}_{x,t}(C(x,t;1/8) )} \leq C \left( ||h_{\lambda}||_{L^p(C(x,t;1))} + \bigg| \bigg|  k^2 \frac{\sin(2h_{\lambda})}{2r^2} \bigg| \bigg|_{L^p(C(x,t;1))} \right) \leq C(n,m\pi,k^2)
\end{align*}
on the cylinder $C(x,t;1/8) = B_{1/8}(x) \times (t-1/64,t)$. We deduce that
 \begin{align*}
          &||(|\partial_t h_{\lambda}| + |\partial_r h_{\lambda}| + |\partial_{rr} h_{\lambda}| )||_{L^{\infty}_{r,t}((5/2,\lambda^{-1}-2) \times (1,\lambda^{-2}T))} \\
          &\lesssim ||(|\partial_t h_{\lambda}| + |\nabla_x h_{\lambda}| + |D^2_xh_{\lambda}| )||_{L^{\infty}_{x,t}(\{5/2 \leq |x| \leq \lambda^{-1}-2\} \times (1,\lambda^{-2}T))} \\
          &\leq C(n, m\pi, k^2),
     \end{align*}
which proves the desired bounds on $h_r$, $h_{rr}$ on $(r,t) \in (5\lambda/2,1-2\lambda) \times (\lambda, T)$ for all $\lambda$ small enough, which is enough to conclude the proof.
\end{proof}

\section{Finite energy solutions and further blow-up characterizations}

In this section, we restrict to solutions of (\ref{heat map flow}) on $[0,T)$ with finite energy and $k$-equivariant symmetry.

Under these assumptions, one obtains additional regularity of solutions at $t = T$ (Proposition \ref{prop:smoothness at the boundary}) as well as existence and rigidity for the value of the limit $\lim \limits_{r \to 0^+} h(r,T)$ (Proposition \ref{value of v and h at origin}). As seen as a corollary of the fixed point argument in Theorem \ref{thm:taylor local existence}, the usual blow-up criterion is the blow-up of the $C^1$-norm of the solution. In this section, we refine the description of blow-up using the so-called ``barrier argument''. Indeed, a blow-up event implies that the limit 
$$
\lim \limits_{\substack{r \to 0^+ \\ t \to T^{-}}} h(r,t)
$$
does not exist, as one can extract at least one non-constant harmonic map along sequences $r_n \to 0^+, t_n \to T^{-}$, as shown in the blow-up description from Struwe (Theorem \ref{thm:smoothness around all but finitely many points}). The barrier argument consists of finding an upper bound $|h(r,t)| \leq g(r,t)$ which is continuous at $(r,t) = (0,T)$ with $g(0,T) = 0$, which leads to a contradiction. This majorant is usually a stationary solution $2\arctan(\alpha^kr^k)$ and it suffices to show the upper bound on the parabolic boundary of $[0,1] \times [0,T]$ thanks to the Comparison Principle (Theorem \ref{comparison principle}).

Under the aforementioned assumptions, if $h(0,t) = 0$ for all $t \in [0,T)$, we prove that it cannot be that the boundary data $h_0$ is bounded by $\pi$ (Proposition \ref{prop:criterion for global solution}, which generalizes the result from \cite{ChangDing2011} to $k \geq 2$ and time-dependent boundary data) and that 
$$
    \limsup_{\substack{r \to 0^+ \\ t \to T^-}}|h(r,t)| = \pi.
$$
Moreover, we also generalize the other result from \cite{ChangDing2011} and show that there are global solutions which blow-up at infinity for any $k \geq 2$.
 
 \begin{definition}[Finite energy solutions]
One says that $v$ has finite energy on $(0,T)$ if 
\begin{align*}
   \sup_{0 < t < T} E(v(\cdot,t)) < + \infty.
\end{align*}
From (\ref{nabla v in terms of h}), it follows that
\begin{align}
    E(v(\cdot,t)) :=  E(v(t)) := \frac{1}{2}\int_{B^2} |\nabla v|^2 dx = \pi \int_0^1 \left( h_r^2 + \frac{k^2}{r^2} \sin(h)^2 \right) r dr = E(h(\cdot,t)). \label{Energy of v in terms of h}
\end{align}
\end{definition}

\begin{remark}[Equivariant solutions have finite energy]\label{k-equivariance implies finite energy}
    If $v_0(x,t) = v_0(x)$ is a time-independent boundary data and $v_0$ has finite energy, then the energy of the solution $v$ is non-increasing in time. Hence, $v$ has finite energy.

    In the time-dependent case, we lose monotonicity. However, if $v_0(x,t)$ is $k$-equivariant, then $v(x,t)$ has finite energy as well even though it might not be monotone in time. Indeed, if $v(x,t)$ is defined on $[0,T)$, then
\begin{align*}
    \frac{d}{dt}E(h(\cdot,t)) &= \pi \int_0^1 \left( 2 h_r h_{rt} + 2 \frac{k^2}{r^2} \sin(h) \cos(h) h_t \right) rdr  \\
     &= 2\pi [h_rh_t r]_{r=0}^{r=1} - 2\pi \int_0^1 h_t \left(- h_{rr}r -  h_r  +  \frac{k^2}{r} \sin(h) \cos(h) \right) dr  \\
     &= 2\pi h_r(1,t)h_t(1,t) - 2\pi \int_0^1 h_t(r,t)^2 rdr \\
     &= 2\pi h_r(1,t)\partial_t h_0(1,t) - 2\pi \int_0^1 h_t(r,t)^2 rdr \\
     &\leq 2\pi h_r(1,t)\partial_t h_0(1,t), \quad t \in [0,T),
\end{align*}
 and we know from Schauder estimates (see Proposition \ref{blow up point is at zero k equivariant}) that $h_r(1,t)$ is bounded on $[0,T)$ under the $k$-equivariance symmetry. Hence,
 $$
\sup_{t \in [0,T)} E(h(\cdot,t)) \leq E(h_0) + C(k,h,T)
 $$
 after integrating with respect to time.
\end{remark}

 First, observe that the symmetry assumption implies some further smoothness at $t = T$. In particular, from now on, it makes sense to consider and write $v(T,x)$ and $h(T,r)$ away from zero.

 {
\begin{proposition}[Smoothness at $t = T$] \label{prop:smoothness at the boundary, earlier statement}
Let
$$
     v(t,x) \in C^0([0,T),C^1(\overline{B^2})) \cap C^{\infty}((0,T) \times \overline{B^2}) 
$$
solve (\ref{heat map flow}) with smooth and $k$-equivariant boundary data. Let 
$$
h(t,r) \in C^0([0,T),C^1([0,1])) \cap C^{\infty}((0,T) \times [0,1])
$$
be the corresponding inclination coordinate (\ref{heat map formulation for h(r,t)}). Assume that $v$ blows-up at time $T < +\infty$.

Then 
$$
v(t,x) \in C^{\infty}((0,T] \times \overline{B^2} \setminus \{0\} ), \quad h(t,r) \in  C^{\infty}((0,T] \times (0,1] ).
$$
\end{proposition}

\begin{proof}
    See Proposition \ref{prop:smoothness at the boundary} in Appendix D. This is a consequence of the energy concentration argument from Struwe (\cite{struwe2008variational}).
\end{proof}

Even though we have an estimate of the blow-up rate, it turns out that the limit $\lim \limits_{r \to 0^+} h(r,T)$ exists and is a multiple of $\pi$. 

\begin{proposition}[Value at the origin]\label{value of v and h at origin}
    Let $v(x,t)$ solve (\ref{heat map flow}) with smooth and $k$-equivariant boundary data. Assume that $v$ blows-up at time $T < +\infty$. Then for all $t \in [0,T]$, one has
    $$
    \lim_{r \to 0^+} h(r,t) = k(t) \pi, \quad k(t) \in \mathbb Z.
    $$
    In particular, $k(t) = c$ is constant on $[0,T)$.
\end{proposition}

\begin{remark}
    It was already known from Proposition \ref{prop:equivalent fomrulation heat map flow} and the uniqueness of a continuous inclination coordinate (up to a multiple of $\pi$) that $h(0,t) = m \pi$ is constant on $[0,T)$. The interesting part is that the limit $\lim \limits_{r \to 0^+} h(r,T)$ also exists and is a multiple of $\pi$. Moreover, there might be a jump at time $T$, i.e.,  $k(T) \neq c$.
\end{remark}

\begin{proof}
    We note that $h(r,t)$ is continuous on $[0,1] \times [0,T) \cap (0,1] \times [0,T]$ (Proposition \ref{prop:equivalent fomrulation heat map flow} and Proposition \ref{prop:smoothness at the boundary, earlier statement} for the additional regularity at $t = T$) and, for fixed $t \in [0,T]$, belongs to $L^2([0,1],rdr)$ since the energy is finite. The existence of the limit and its value follow from the following two lemmas. The fact that $k(t)$ is constant on $[0,T)$ follows from the continuity of $h$ on $[0,1] \times [0,T)$.
\end{proof}

\begin{lemma}\label{value at r = 0 is a multiple of pi}
Assume that $\theta(r) \in L^2([0,R],rdr)$ with a weak derivative $\theta'(r) \in L^2([0,R],rdr)$. Moreover, assume there exists $k \in \mathbb Z$ for which 
$$
\frac{\theta-k\pi}{r} \in L^2([0,R],rdr)
$$
as well. Then $\theta \in C^0([0,R])$ and $\lim \limits_{r \to 0}\theta(r) = k \pi$.
\end{lemma}

\begin{proof}
    For any compact interval $I \subset (0,R]$, one has $\theta, \theta' \in L^2(I)$. Hence $\theta \in W^{1,2}(I) \subset C^{0}(I)$.

    Then $(\theta-k\pi)^2 \in W^{1,2}(I)$ as well and $\partial_r \left[ (\theta-k\pi)^2 \right] = 2(\theta-k\pi) \theta'$ using the product rule.

    Since
    $$
2(\theta-k\pi) \theta' = \frac{\theta-k\pi}{r}r^{1/2} \cdot \theta' r^{1/2} \in L^1([0,R]),
    $$
    one actually has $(\theta-k\pi)^2 \in W^{1,1}([0,R]) \subset C^0([0,R])$ and $(\theta-k\pi)^2$ is continuous at $r = 0$.

    Finally, since
    $$
    \int_0^R \frac{(\theta(r)-k\pi)^2}{r} dr < +\infty, 
    $$
    one must have $\lim \limits_{r \to 0} (\theta - k\pi)^2 = 0$, which concludes the proof.
\end{proof}

\begin{lemma}\label{value at r = 0 is a multiple of pi, ver.2}
    Assume that $\theta(r) \in L^2([0,R],rdr)$ with a weak derivative $\theta'(r) \in L^2([0,R],rdr)$. Then 
    $$
    \exists k \in \mathbb Z: \frac{\theta-k\pi}{r} \in L^2([0,R],rdr) \iff \frac{\sin \theta}{r} \in  L^2([0,R],rdr).
    $$
\end{lemma}

\begin{proof}
    The direction $\Rightarrow$ follows from the inequality $|\sin(\theta)| = |\sin(\theta-k\pi)| \leq |\theta-k\pi|$.
    
    We prove $\Leftarrow$. For any compact interval $I \subset (0,R]$, $\theta \in W^{1,2}(I) \subset C^0(I)$. Applying the chain rule, $\cos \theta \in W^{1,2}(I)$ with 
    $$
    \partial_r \cos(\theta) = -\sin(\theta) \theta' = -\frac{\sin \theta}{r}r^{1/2} \cdot \theta' r^{1/2} \in L^1([0,R]),
    $$
    meaning that $\cos \theta \in W^{1,1}([0,R]) \subset C^0([0,R])$. Hence, $\sin(\theta)^2 = 1 - \cos(\theta)^2 \in C^0([0,R])$ and $\lim \limits_{r \to 0} \sin(\theta(r))^2 = 0$ because
    $$
    \int_0^R \frac{\sin(\theta)^2}{r} dr < +\infty. 
    $$
    Note that $\sin(\theta) \in W^{1,2}(I) \subset C^0(I)$ as well. Hence, we have proved that $\sin(\theta)$ is continuous at zero with $\lim \limits_{r \to 0} \sin(\theta(r)) = 0$. One checks, using an argument by contradiction and the Intermediate Value property of $\theta$ away from $r = 0$, that this implies $\lim \limits_{r \to 0} \theta(r) = k \pi$ for some $k \in \mathbb Z$.

    Finally, $|\sin(\theta(r))| = |\sin(\theta(r)-k\pi)| \gtrsim |\theta(r) - k\pi|$ near $r = 0$ which implies the desired integrability.
\end{proof}

Next, we present a theorem which, in a more general setting, states that solutions to nonlinear parabolic problems with analytic nonlinearity are real-analytic with respect to the space variable. This result is attributed to Friedman (\cite[Theorem 2]{friedman}), but looking at the paper, the proof is omitted. In any case, we only need a weaker statement for our purpose, taking inspiration from Friedman's techniques for the elliptic case.

\begin{theorem}[Space-analyticity]\label{thm:space analyticity}
Let $v(t,x)$ solve (\ref{heat map flow}) on $(0,T)$, $T  < +\infty$, with smooth, $k$-equivariant boundary data. Let
$$
h(t,r) \in C^0([0,T),C^1([0,1])) \cap C^{\infty}((0,T] \times (0,1]) 
$$
be the corresponding inclination coordinate (\ref{heat map formulation for h(r,t)}) chosen with $h(t,0) = h_0(t,0) = 0$.

Then $r \mapsto h(r,t)$ is real-analytic at $r_0 \in (0,1]$ for any fixed $t \in (0,T]$.
\end{theorem}

\begin{proof}
After parabolic rescaling, one can assume that $r_0 = 1$.

    Let $H = h - h(t,1)$. Then that $H(t,x) \in C^0([0,T);C^1(\overline{B^2}))$ (the regularity on $\overline{B^2}$ can be found in Proposition \ref{prop:equivalent fomrulation heat map flow} and its proof) satisfies a 2D nonlinear heat equation
    \begin{align*}
        H_t &= \Delta H  - k^2\frac{\sin(2H + 2h(t,1))}{2|x|^2} - h_t(t,1), \quad x \in \overline{B^2} \setminus \{0\} \times (0,T), \\
        H &= 0, \quad (x,t) \in \overline{B^2} \times \{0\} \cup \partial{B^2} \times (0,T).
    \end{align*}
   Fix a $\varepsilon$-neighborhood $\Omega_{\varepsilon}$ of the boundary $\partial B^2$ on which Theorem \ref{thm:boundary parabolic sobolev embedding} applies and let $G = \overline{\Omega_{\varepsilon}}$. 

    Suppose by induction that one has
    \begin{align*}
        ||H||_{L^{\infty}(G \times (0,T))} + ||\nabla H||_{L^{\infty}(D \times (0,T))} \leq H_0, \\
        ||D^q_x H||_{L^{\infty}(G \times (0,T))} \leq H_0 H_1^{q}q!, \quad 2 \leq q \leq n,
    \end{align*}
    for $n \geq 2$. Letting 
    $$
X(t,x,h) = - k^2\frac{\sin(2H + 2h(t,1))}{2r^2} - h_t(t,1),
    $$
    one has
  \begin{align*}
      |\partial^n_{x_1} X| &\leq k^2 \sum_{i=0}^n \binom{n}{i} (n-i)!C(G)^{n-i} \sum_{\sum_{j=1}^i jm_j = i} \frac{i!}{m_1! \cdot ... \cdot m_i!} \prod_{j=1}^i \left( \frac{2 |\partial^{j}_{x_1} H|}{j!} \right)^{m_j}  \\
     &\leq n!  k^2  \sum_{i=0}^n  C(G)^{n-i}\sum_{\sum_{j=1}^i jm_j = i} \frac{1}{m_1! \cdot ... \cdot m_i!} \prod_{j=1}^i \left( \frac{2 H_0 H_1^{j} j!}{j!} \right)^{m_j}  \\ 
     &\leq n!    k^2 H_1^n  \sum_{i=0}^n  C(G)^{n-i} \sum_{\sum_{j=1}^i jm_j = i} \prod_{j=1}^i \frac{(2H_0)^{m_j}}{m_j!}
  \end{align*}
 on $G \times (0,T)$ for $n \geq 2$ using the product rule and Faa di bruno's formula. The constant $C(G)$ essentially comes from the $r^{-2-n+i}$ factor when differentiating $X$. When $i = 0$, we use the convention that
 $$
\sum_{\sum_{j=1}^i jm_j = i} \frac{i!}{m_1! \cdot ... \cdot m_i!} \prod_{j=1}^i \left( \frac{2 |\partial^{j}_{x_1} H|}{j!} \right)^{m_j}  = 1.
 $$
 
 There exists $B_0 = B_0(H_0) > 0$  for which
 $$
\frac{(2H_0)^{j}}{j!} \leq B_0 \quad \forall j \in \mathbb N_{\geq 0}.
 $$
 Hence,
   \begin{align*}
      |\partial^n_{x_1} X| &\leq n!  H_1^n k^2 \sum_{i=0}^n C(G)^{n-i} B_0^i \sum_{\sum_{j=1}^i jm_j = i}1  \leq n!  k^2 H_1^n \sum_{i=0}^n   C(G)^{n-i}B_0^i \binom{2i-1}{i-1} \\
      &\leq n!   k^2 H_1^n  B_1 \sum_{i=0}^n   C(G)^{n-i} B_0^i 4^{i} \leq  n! k^2 H_1^n  B_1 C(G)^n \frac{ C(G)^{-(n+1)} B_0^{n+1} 4^{n+1} - 1}{4C(G)^{-1} B_0-1},
   \end{align*}
  where the binomial coefficient is treated as $1$ when $i = 0$ and $B_1$ is a universal constant obtained, for example, via Stirling's formula. Taking $B_0 = B_0(H_0,G)$ larger if needed, one can assume that $4C(G)^{-1}B_0 \gg 1$ and
  $$
|\partial^n_{x_1} X| \leq 2k^2 B_1 4^n H_1^n B_0^{n} n!
  $$
  on $G \times (0,T)$ for $n \geq 2$.
  
Then Theorem \ref{thm:boundary parabolic sobolev embedding} yields that
    $$
    ||\partial^n_{x_1} h||_{L^{\infty}(G \times (0,T))} \leq 2C(G)k^2 B_1 4^n H_1^n B_0^{n} n!.
    $$
    Setting
    $$
    H_0 =||h||_{L^{\infty}(G \times (0,T))} + ||\nabla h||_{L^{\infty}(G \times (0,T))}  +  ||D^2_x h||_{L^{\infty}(G \times (0,T))} + 1,
    $$
    as well as $H_1 =8 k^2C(G) B_1 B_0(H_0,G)+ 1$ as constants initializing our induction, one deduces
    $$
    ||\partial^n_{x_1} h||_{L^{\infty}(K \times (0,T))} \leq H_0H_1^{n+1}(n+1)!
    $$
    for all $n \geq 2$ and similarly for $\partial^n_{x_2}h$. As $h$ is smooth on $G \times (0,T]$, we obtain that $h(x,t)$ is real-analytic with respect to $x \in G$ for any fixed $t \in (0,T]$. In particular, $h(r,t)$ is real-analytic at $r = 1$ for any fixed $t \in (0,T]$.
\end{proof}

We come now to the so-called barrier argument. First, we show that one can construct a barrier in some subset of the parabolic boundary. 

\begin{lemma}\label{lemma:semi-barrier with global solution}
     Let $v(x,t)$ solve (\ref{heat map flow}) on $[0,T)$, $0 < T \leq +\infty$, with smooth, $k$-equivariant boundary data. Assume that the inclination coordinate $h$ is chosen so that $h(0,t) = 0$. Then for all $t_0 \in (0,T)$, there exists $\alpha_0 > 0$ and $r_0 \in (0,1)$ for which
      \begin{equation}
       \theta_{\alpha}(r) \geq |h(r,t_0)|, \quad \alpha \geq \alpha_0, \quad r \in [0,r_0], \label{eq:comparison between theta_alpha and any solution}
   \end{equation}
   where $\theta_{\alpha}(r) = 2\arctan\left( (\alpha r)^k \right)$.
\end{lemma}

\begin{proof}
    Fix an arbitrary $t_0 \in (0,T)$. As $v(x,t)$ is smooth and $k$-equivariant on $\overline{B^2} \times \{t_0\}$, one has
    $$
    h(0,t_0) = \partial_r h(0,t_0) = ... = \partial_r^{(k-1)} h(0,t_0) = 0
    $$
    by Lemma \ref{k-equivariance, smoothness of h}. Similarly, one computes
    $$
    \theta_{\alpha} (0) = \partial_r \theta_{\alpha} (0) = ... = \partial_r^{(k-1)} \theta_{\alpha} (0) = 0,
    $$
    while $\partial_r^{(k)} \theta_{\alpha}(0) = 2k!\alpha^k$. Fix $\alpha_0 \gg 1$ for which $\partial_r^{(k)} \theta_{\alpha_0}(0) > |\partial_r^{(k)} h(0,t_0)|$. By continuity, this inequality remains valid on $[0,r_0]$ for some $r_0 > 0$. In particular, for $r \in [0,r_0]$,
\begin{align*}
        \theta_{\alpha_0}(r) &= \int_0^r \int_0^{s_{k-1}} ... \int_0^{s_1} \partial_r^{(k)} \theta_{\alpha_0}(s_0) ds_0 ds_1 ... ds_{k-1} \\
        &\geq \int_0^r \int_0^{s_{k-1}} ... \int_0^{s_1} \partial_r^{(k)}h(s_0,t_0) ds_0 ds_1 ... ds_{k-1} = h(r,t_0),
\end{align*}
    and similarly, $\theta_{\alpha_0}(r) \geq -h(r,t_0)$.
    
    Moreover, if $\alpha \geq \alpha_0$, then
    $$
\theta_{\alpha}(r) \geq  \theta_{\alpha_0}(r) \geq |h(r,t_0)|, \quad r \in [0,r_0],
    $$
    by monotonicity of the $\arctan$ function. 
\end{proof}

\begin{proposition}[Further Property of Blow-up]\label{prop:further properties of blow-up of h}
    Let $v(x,t)$ solve (\ref{heat map flow}) with smooth and $k$-equivariant boundary data. Assume that $v$ blows-up at time $T < +\infty$. Then one has
    $$
    \limsup_{\substack{r \to 0^+ \\ t \to T^-}}|h(r,t)-h(0,t)| \geq \pi.
    $$
\end{proposition}

\begin{proof}
     Replacing $h$ and $h_0$ by $h - m\pi$ and $h_0 - m\pi$, where $h(0,t) = m\pi$ for $t \in [0,T)$ (Proposition \ref{value of v and h at origin}), one can assume that $h(0,t) = 0$ for $t \in [0,T)$ and one needs to prove 
     $$
    \limsup_{\substack{r \to 0^+ \\ t \to T^-}}|h(r,t)| \geq \pi.
    $$
    By the Comparison Principle (Proposition \ref{prop:equivalent fomrulation heat map flow}, Theorem \ref{comparison principle}), $h$ is bounded by a large multiple of $\pi$. Hence, the limsup exists. 

    Assume for a contradiction that 
    $$
 \limsup_{\substack{r \to 0^+ \\ t \to T^-}}|h(r,t)| = M < \pi.
    $$
    We use a barrier argument using a stationary solution to obtain a contradiction. There exists $0 < \delta < \pi - M$ and a neighborhood $(0,r_0) \times (t_0,T)$ with $r_0 < 1$, $t_0 > 0$, on which $|h(r,t)| < M + \delta$. Consider the stationary solution $\theta_{\alpha}(r) = 2 \arctan(\alpha^k r^k)$ of (\ref{heat map formulation for h(r,t)}). Up to choosing a smaller $r_0$, there is some $\alpha_0 > 0$ for which
    $$
\theta_{\alpha}(r) \geq  \theta_{\alpha_0}(r) \geq |h(r,t_0)|, \quad \alpha \geq \alpha_0, r \in [0,r_0],
    $$
    by Lemma \ref{lemma:semi-barrier with global solution}.
    
    Fix $\alpha \geq \alpha_0$ so that $\theta_{\alpha}(r_0) > M$. One has found a solution $\theta_{\alpha}$ satisfying $h(0,t) = 0 = \theta_{\alpha}(0)$ on $\{0\} \times (t_0,T)$, $|h(r_0,t)| \leq M < \theta_{\alpha}(r_0)$ on $\{r_0\} \times (t_0,T)$ and $|h(r,t_0)| \leq \theta_{\alpha}(r)$ on $[0,r_0] \times \{t_0\}$.
    
    The Comparison Principle (Theorem \ref{comparison principle}) applied with $h, \theta_{\alpha}$ and $-h, \theta_{\alpha}$ on $[0,r_0] \times [t_0,T']$ for any $T' < T$ shows that $|h(r,t)| \leq \theta_{\alpha}(r)$ on $[0,r_0] \times [t_0,T)$. In particular, $\lim \limits_{r \to 0^+} |h(r,t)| = 0$ uniformly with respect to $t \in [t_0,T)$. 
     As
    $$
      v(re^{i\theta},t) = \left( e^{ik\theta} \sin h(r,t), \cos h(r,t) \right),
    $$
    one deduces that $\lim \limits_{x \to 0} v(x,t) = (0,0,1)$ uniformly with respect to $t \in [t_0,T)$ as well.
    
    If $v(x,t)$ blows-up at time $T$, blows-up happens at $(0,T)$ (Proposition \ref{prop:smoothness at the boundary, earlier statement}). Then $v_m(x) = v(x_m + R_mx, T_m + R_m^2s)$ converges strongly in $H^{1,2}_{loc}(\mathbb R^2, S^2)$ to a non-constant harmonic map along some appropriate sequences $x_m \to 0, (x_m) \subset B^2, R_m \to 0^+, T_m \to T^-$ and some fixed $s \leq 0$ (Theorem \ref{thm:smoothness around all but finitely many points}). Up to taking a subsequence, the convergence holds pointwise almost everywhere for $x \in \mathbb R^2$. But for fixed $x \in \mathbb R^2$, $x_m + R_mx \to 0$ and since $T_m + R_m^2s \to T^-$, the uniform convergence $\lim \limits_{x \to 0} v(x,t) = (0,0,1)$ implies that $\lim \limits_{m \to +\infty} v(x_m + R_mx, T_m + R_m^2s) = (0,0,1)$ is a constant function of $x$, which is a contradiction.
\end{proof}

\begin{proposition}[Criterion for existence of global solution]\label{prop:criterion for global solution}
    Let $v(x,t)$ solve (\ref{heat map flow}) on $[0,T)$ with smooth, $k$-equivariant boundary data.  Assume that the initial inclination coordinate satisfies $h_0(0,t) = 0$, $|h_0(r,t)| \leq \pi$ on $\{0,1\} \times [0,T) \cup [0,1] \times \{0\}$. If $T < +\infty$, then $v(x,t)$ and $h(r,t)$ cannot blow-up at time $T$. In particular, if  $|h_0(r,t)| \leq \pi$ on $\{0,1\} \times [0,+\infty) \cup [0,1] \times \{0\}$, then $v(x,t)$ and $h(r,t)$ are global. 
\end{proposition}

\begin{proof}
This is a generalization of the global solution constructed in \cite{ChangDing2011}'s paper for $k > 1$. The proof can be made shorter thanks to the Comparison Principle (Theorem \ref{comparison principle}).

    Assume for a contradiction that $T < +\infty$ and the solution blows-up at time $T$. We use a barrier argument, as in Proposition \ref{prop:further properties of blow-up of h}.
    
    It follows from the Comparison Principle (Theorem \ref{comparison principle}) that $|h(r,t)| \leq \pi$ on $[0,1] \times [0,T)$. Write
    $$
    (h-\pi)_t = (h-\pi)_{rr} +\frac{(h-\pi)_r}{r} -k^2 (h-\pi) \frac{\sin(2h)-\sin(2\pi)}{2r^2(h-\pi)}
    $$
    Assume for a contradiction that $h(r_0,t_0) = \pi$ at some point $(r_0,t_0) \in (0,1) \times (0,T]$ ($t_0 = T$ is allowed as $h$ is smooth on $(0,1] \times (0,T]$ by Proposition \ref{prop:smoothness at the boundary, earlier statement}). On any fixed region $E_{r_1,t_1} = (r_1,1) \times (t_1,t_0]$, $0 < r_1 < r_0$, $0 < t_1 < t_0 \leq T$, the coefficient in front of $h-\pi$ in the above PDE is bounded (in particular, from above). The Maximum Principle (Theorem \ref{maximum principle}) applied on $H = h-\pi$ implies that $h(r,t) = \pi$ is constant on $[0,1] \times [0,t_0]$, which is a contradiction to $h(0,t) = 0$. Similarly, if $h(r_0,t_0) = -\pi$, applying the Maximum Principle shows that $h = -\pi$ is constant on $[0,1] \times [0,t_0]$.  Hence, $|h(r,t)| < \pi$ on $(0,1) \times (0,T]$.

Fix an arbitrary $t_0 \in (0,T)$. It follows from Lemma \ref{lemma:semi-barrier with global solution} that
    $$
\theta_{\alpha}(r) \geq  \theta_{\alpha_0}(r) \geq |h(r,t_0)|, \quad \alpha \geq \alpha_0, r \in [0,r_0]
    $$
    for some $\alpha_0 > 0, r_0 \in (0,1)$. Moreover,
    $$
    \sup_{t \in [t_0,T]} |h(r_0,t)| = M < \pi
    $$
    as a consequence of the Comparison Principle. Fixing $\alpha \geq \alpha_0$ so that $\theta_{\alpha}(r_0) > M$, one can conclude as in the proof of Proposition \ref{prop:further properties of blow-up of h}.
\end{proof}

\begin{corollary}[An even more precise blow-up description, \cite{vanderhout_scale}]\label{cor:comparison of h with itself and exact value of limsup}
 Let $v(x,t)$ solve (\ref{heat map flow}) with smooth and $k$-equivariant boundary data on $[0,T)$ with $T < +\infty$ (which may or may not be the maximal existence time). Assume that $h$, $h_0$ are chosen so that $h(0,t) = 0$ on $[0,T)$. There exists $\tau_0 = \tau_0(v,T)$ such that for all $0 < \tau \leq \tau_0$, $t \in [\tau,T)$, one has
    $$
    h(r,t-\tau) - \pi \leq h(r,t) \leq h(r,t-\tau) + \pi, \quad r \in [0,1].
    $$
In particular, if $T = T_{\max} < +\infty$, for any sequence $(r_n,t_n) \to (0^+,T^-)$, one has
$$
-\pi \leq \liminf_{n \to +\infty} h(r_n,t_n) \leq \limsup_{n\to +\infty} h(r_n,t_n) \leq \pi,
$$
hence
$$
    \limsup_{\substack{r \to 0^+ \\ t \to T^-}}|h(r,t)| = \pi.
$$
\end{corollary}

\begin{proof}
    As $h(r,t) \in C^{0}([0,1] \times [0,T'])$ is uniformly continuous for any $0 < T' < T$, there is $\tau_0 > 0$ small enough for which $h(r,\tau) \leq h(r,0) + \pi$ for all $r \in [0,1]$, $0 < \tau \leq \tau_0$. Similarly, $h(1,t) = h_0(1,t) \leq h_0(1,t-\tau) + \pi = h(1,t-\tau) + \pi$ for all $\tau \leq \tau_0 \leq t \leq T$ by uniform continuity of the boundary data. Finally, $0 = h(0,t) \leq h(0,t-\tau) + \pi = \pi$ for all $\tau \leq \tau_0 \leq t \leq T$. For $0 < \tau \leq \tau_0$ fixed, the Comparison Principle (Theorem \ref{comparison principle}) implies that 
    $$
    h(r,t-\tau) - \pi \leq h(r,t) \leq h(r,t-\tau) + \pi, \quad r \in [0,1], t \in [\tau,T),
    $$
    where we used a similar argument with $-\pi$ instead of $\pi$. The limsup value follows from Proposition \ref{prop:further properties of blow-up of h}.
\end{proof}

\begin{theorem}[Existence of a global solution, \cite{ChangDing2011}]\label{thm:existence of global solution}
    Assume that $\psi$ is non-constant, $\psi = r^k\tilde{\psi}$, $\tilde{\psi} \in C^{\infty}([0,1])$ with $\partial^{(2n+1)}_r \tilde{\psi}(0) = 0$ for all $n \in \mathbb N_{\geq 0}$, $0 \leq \psi \leq \pi$, and
    $$
\tau(\psi) := \psi_{rr} + \frac{1}{r}\psi_r -k^2 \frac{\sin(2\psi)}{2r^2} \geq 0, \quad \tau(\psi)(1) = 0.
    $$
   Then the solution $h(r,t) \in C^{\infty}([0,1] \times (0,T))$ of (\ref{heat map formulation for h(r,t)}) with boundary data $h_0(r,t) = \psi(r)$ on $[0,1] \times \{0\} \cup \{0,1\} \times [0,+\infty)$ is global, i.e.,  $T = +\infty$,  has finite energy (\ref{Energy of v in terms of h}), satisfies $0 < h < \pi$ and $h_t > 0$ on $(0,1) \times (0,T)$, as well as $h_r > 0$ on $(0,1] \times (0,T)$.
\end{theorem}

\begin{remark}
    One can choose $\psi(r) = 4 \arctan \left( (\alpha r)^k \right)$, $\alpha \in (0,1]$, which satisfies
    $$
\tau(\psi) = -k^2 \frac{\sin(2\psi)}{2r^2}  + 2 k^2 \frac{\sin(\psi)}{2r^2} = k^2r^{-2} \sin(\psi) (1-\cos(\psi)) \geq 0,
    $$
    and all the other required properties. If one wants $\psi(1) = \pi$ (as in Proposition \ref{prop:blow up infinity global solution}) to have a blow-up at infinity, then one needs to take $\alpha = 1$.
\end{remark}

\begin{proof}
    This is a generalization of \cite{ChangDing2011}'s result for $k > 1$, but the proof is essentially the same because the sign of the nonlinearity does not change and we still have access to a family of stationary solutions given by the arctan function. Let $\psi$ be as in the theorem.
    
    Let $h(r,t) \in C^{\infty}([0,1] \times (0,T))$ solve (\ref{heat map formulation for h(r,t)}) with boundary data $h_0(r,t) = \psi(r)$ on $[0,1] \times \{0\} \cup \{0,1\} \times [0,+\infty)$ with maximal time of existence $T =  +\infty$ by Proposition \ref{prop:criterion for global solution}. 
    
    Then $h, h_r, h_t$ are well-defined and continuous on $[0,1] \times [0,T)$, $h_{rr}$ is well-defined and continuous on $(0,1] \times [0,T)$ (Proposition \ref{prop:continuity of h_t}). Differentiating $h(0,t) = \psi(0)$, $h(1,t) = \psi(1)$ with respect to $t$ yields $h_t(0,t) = h_t(1,t) = 0$ on $[0,T)$ and (\ref{heat map formulation for h(r,t)}) at $t = 0$, $r > 0$ yields $h_t(r,0) = \tau(\psi)(r) \geq 0$ on $(0,1]$, $h_t(1,0) = \tau(\psi)(1) = 0$. In the following, we prove that $h_t > 0$ on $(0,1) \times (0,T)$ and $h_r > 0$ on $(0,1] \times (0,T)$.
    
First, $u = h_t$ solves a 2D radial heat equation
    $$
    u_t = \Delta u - k^2 \frac{2\cos(2h)}{2|x|^2} u, \quad (r,t) \in  (0,1) \times (0,T).
    $$
    For any fixed $0 \leq t < T$, $\cos(2h) \approx 1$ near $r = 0$ and away from $r = 0$, the coefficient of $u$ is bounded. Hence, it is bounded from above. The Maximum Principle (Theorem \ref{maximum principle}) applied on $E_{t_1} = (B^2 \setminus B_{\varepsilon}(0)) \times (0,t_1]$, $t_1 < T$, with $v = ue^{-ct}$ for $c$ sufficiently large shows $v$ cannot achieve a minimum $m \leq 0$ inside $E_{t_1}$. Hence, either the minimum is achieved at the parabolic boundary, on which case it must be $m = 0$ or $m = \inf_{0 < t < t_1} h_t(\varepsilon,t)$, either there is a positive minimum $m > 0$ achieved inside $E_{t_1}$. The latter is not possible as points $h_t = 0$ at the parabolic boundary $\{r = 1\}$ would be a smaller minimum. Taking $\varepsilon \to 0$, the minimum on $B^2 \times (0,t_1]$ is  $m = 0$. In other words, $h_t > 0$ on $(0,1) \times (0,T)$.

   Writing
    $$
    -k^2 \frac{\sin(2h)}{2r^2} =  -k^2 \frac{\sin(2h)}{2hr^2} h,
    $$
    one has $\sin(2h)/2h \approx 1$ near $r = 0$ and away from zero, the coefficient of $h$ is bounded. As for $h_t$, the Maximum Principle implies $0 < h < M \leq \pi$ on $(0,1) \times (0,T)$, where $0$ and $0 < M \leq \pi$ are respectively the minimum and maximum of $h$ achieved on the parabolic boundary.
    
Next, write $h(r,t)=r^k \tilde{h}(r,t)$. Then $\tilde{h}$ is non-negative, extends smoothly across $r=0$ (Proposition \ref{prop:equivalent fomrulation heat map flow}) and satisfies
$$\tilde{h}_t = \tilde{h}_{rr} + \frac{2k+1}{r}\tilde{h}_r + \frac{k^2}{r^2} \left( 1-\frac{\sin(2h)}{2h} \right)\tilde{h}, \quad (r,t) \in  (0,1) \times (0,T).
$$
Since $h(r,t)=O(r^k)$ near $r = 0$ for fixed $0 < t < T$ and
$$
1-\frac{\sin(2h)}{2h}=O(h^2),
$$
the zeroth-order coefficient
$$
\frac{k^2}{r^2}
\left(
1-\frac{\sin(2h)}{2h}
\right)
$$
in front of $\tilde{h}$ extends as a bounded function near $r=0$ (it is bounded from above and below, by contrast with $h$ and $h_t$). Regarding $\tilde{h}$ as a radial function on the ball $B^{(2k+2)}_1(0)$, a minimum $\tilde{h}(0,t) = 0$ cannot be achieved at the origin by the Maximum Principle (Theorem \ref{maximum principle}). In other words, $\tilde{h}(0,t) > 0$ for any $t \in (0,T)$. Moreover, for fixed $t \in (0,T)$, smoothness and radiality give
$$
\tilde{h}(r,t)=\tilde{h}(0,t)+O(r^2),
\qquad
\tilde{h}_r(r,t)=O(r),
$$
near $r = 0$. Therefore,
\begin{align*}
h_r(r,t)
&=
k r^{k-1}\tilde{h}(r,t)+r^k \tilde{h}_r(r,t) \\
&=
k \tilde{h}(0,t)r^{k-1}+O(r^{k+1}).
\end{align*}
This implies $h_r(r,t) > 0$ on $(0,\delta_t]$ for some $\delta_t > 0$.


     Actually, $h_r(r,t) > 0$ on $(0,1] \times (0,T)$. Otherwise, there is $t \in (0,T)$ and $r^* \in (0,1]$ for which $h_r(r,t) > 0$ on $(0,r^*)$, $h_r(r^*,t) = 0$. Multiplying (\ref{heat map formulation for h(r,t)}) by $r^2h_r$, one finds
     \begin{align}
         r^2h_rh_t = \frac{1}{2}\partial_r \left[r^2 h_r^2 + k^2 \frac{\cos(2h)}{2} \right], \quad (r,t) \in (0,1) \times (0,T), \notag \\
         2\int_0^r s^2 h_r(s,t)h_t(s,t)ds = r^2 h_r^2 + k^2 \frac{\cos(2h)-1}{2} = r^2 h_r^2 - k^2 \sin(h)^2, \quad (0,1) \times (0,T). \label{eq:sacks-uhlenbeck identity}
     \end{align}
     At $(r^*,t)$, (\ref{eq:sacks-uhlenbeck identity}) has strictly positive left-hand side with a non-positive right-hand side, which is a contradiction.

As $h$ does not depend on time at the boundary $r = 1$, one checks that
    $$
\frac{d}{dt} E(h(\cdot,t)) = -\int_0^1 |h_t|^2 r dr, \quad t > 0,
    $$
    which means that the energy is decreasing over time. Hence, $h$ is a finite-energy solution on $[0,1] \times [0,T)$.
    

\end{proof}

\begin{proposition}[Blow-Up at infinity of global solution, \cite{ChangDing2011}]\label{prop:blow up infinity global solution}
Let $\overline{h}$ be the global solution constructed in Theorem \ref{thm:existence of global solution} for the $k$-equivariant problem, $k \geq 1$. Assume there is $r^* \in (0,1]$ for which $\psi(r^*) = \pi$. Then $\overline{h}$ blows-up at $T = +\infty$ in the sense that for any $m \in \mathbb N$, there exists sequences $T_n \to +\infty$, $R_n \to 0^+$, $l = l(m)$, $\alpha > 0$ for which 
    $$
\lim \limits_{n \to +\infty } ||\overline{h}_r(\cdot,T_n)||_{C^0([0,r_0])} = +\infty \quad \forall r_0 \in (0,1],
$$
and $\overline{h}_n(r) = \overline{h}(R_nr, R_n^2t + T_n): [0,R_n^{-1}2^{-l}] \times [-1,0] \to [0,\pi]$ converges in $C^{m}_{loc}(\mathbb R \times [-1,0])$ to the stationary map
$$
\overline{H}_{\infty}(r) = 2 \arctan \left( (\alpha r)^k \right).
$$
Moreover, for any $r \in (0,1]$,
$$
\lim \limits_{t \to +\infty} \overline{h}(r,t) = \pi.
$$
\end{proposition}

\begin{proof}
     As $\overline{h}$ does not depend on time at the boundary $r = 1$, one checks that
    $$
\frac{d}{dt} E(\overline{h}(\cdot,t)) = -\int_0^1 |\overline{h}_t|^2 r dr, \quad t > 0.
    $$
    If $\overline{v}$ is the $k$-equivariant solution of the heat-map flow corresponding to $\overline{h}$ (Proposition \ref{prop:equivalent fomrulation heat map flow}), one deduces the energy inequality
  \begin{align}
        E(\overline{h}(\cdot,t)) = E(\overline{v}(\cdot,t)) &\leq E(\overline{v}(\cdot,0)) < +\infty, \label{eq: energy inequality, bounded energy} \\
        \int_0^{+\infty} \int_{0}^1 |\partial_t \overline{h}|^2r dr dt = \int_0^{+\infty} \int_{B^2} |\partial_t \overline{v}|^2 dx dt &\leq E(\overline{v}(\cdot,0)) < +\infty, \label{eq: energy inequality, bounded partial_t v L^2 norm}
  \end{align}
    (see also \cite[Chapter 3, Lemma 5.8]{struwe2008variational}). The energy inequality implies
    \begin{equation}
        \lim \limits_{n \to +\infty} \int_{0}^1 |\partial_t \overline{h}(\cdot,T_n)|^2 rdr = \lim \limits_{n \to +\infty} \int_{B^2} |\partial_t \overline{v}(\cdot,T_n)|^2 dx = 0  \label{eq: decay of partial_t overline v along T_n}
    \end{equation}
    along some sequence $T_n \to +\infty$. It follows that along any such sequence, there is a blow-up of $\overline{h}_r$ at $r = 0$ in the sense that
    $$
\lim \limits_{n \to +\infty } ||\overline{h}_r(\cdot,T_n)||_{C^0([0,r_0])} = +\infty \quad \forall r_0 \in (0,1].
$$
    Otherwise, $\overline{h}_r$ is bounded on $[0,r_0]$ along the sequence $T_n \to +\infty$. Along a subsequence, $\overline{v}(x,T_n)$ must then converge weakly in $H^{1,2}(\overline{B}^2, S^2)$ to a harmonic map $\overline{v}_{\infty}$ which must be smooth on $B^2$ (\cite[Theorem 3.6]{sacks-uhlenbeck}). In fact, the convergence is uniform on $\overline{B^2}$ up to taking another subsequence. Away from zero, this follows from the fact that $h(r,T_n)$ converges weakly in $H^{1,2}([r_0,1])\hookrightarrow C^{0,1/4}([r_0,1]) \subset \subset C^0([r_0,1])$, hence uniformly, to a $H^{1,2}([r_0,1])$-map. Moreover, on $[0,r_0]$, $h(r,T_n)$ converges uniformly along a subsequence to a $C^0([0,r_0])$-map by Arzela-Ascoli as $h_r$ is bounded. Finally, the harmonic map satisfies $\overline{v}_{\infty}(\partial B_{r^*}(0)) = \{(0,0,-1)\}$, meaning that $\overline{v}_{\infty}$ must be constant (\cite[Theorem 3.2]{lemaire}), which contradicts $\overline{v}_{\infty}(0) = (0,0,1)$. 

    Let
    $$
    \theta_T = \max_{t \in [0,T]} ||\partial_r \overline{h}(\cdot,t)||_{C^0([0,1])}\to  +\infty, \quad T \to +\infty.
    $$
    One can find $T_n \to +\infty$, $r_n \in [0,1]$, such that
    $$
     R_n^{-1} := \theta_{T_n} = ||\partial_r \overline{h}(\cdot,T_n)||_{C^0([0,1])} = |\partial_r \overline{h}(r_n,T_n)| \to +\infty.
    $$
    Consider the sequence $\overline{V}_n(x,t) = \overline{v}(R_nx, R_n^2t + T_n)$, $t \in [-1,0]$, where we take $n$ large enough so that $T_n - R_n^2 > T_{n-1}$, $T_0 - R_0^2 > 0$. It follows from (\ref{eq: energy inequality, bounded partial_t v L^2 norm}) that
    \begin{equation}
        \lim \limits_{n \to +\infty} \int_{-1}^{0} \int_{B_{R_n^{-1}}(0)} |\partial_t \overline{V}_n|^2 dxdt = \lim \limits_{n \to +\infty} \int_{T_n-R_n^2}^{T_n} \int_{B^2} |\partial_t \overline{v}|^2 dxdt = 0. \label{eq: convergence of partial_t overline v^2 to zero in L^2}
    \end{equation}
   Recall now from (\ref{nabla v in terms of h}) that
\begin{equation*} 
    |\partial_r \overline{h}|^2 \leq |\nabla_x \overline{v}|^2 \leq  |\partial_r \overline{h}|^2 + \frac{k^2}{r^2}\sin(\overline{h})^2, \quad (r = |x|,t) \in (0,1] \times [0,+\infty),
\end{equation*}
where 
$$
\sin(\overline{h}(r,t)) = r \cos(\overline{h}(c_{r,t},t)) \partial_r \overline{h}(c_{r,t},t) \leq r||\partial_r \overline{h}(\cdot,t)||_{C^0([0,1])}.
$$
 Hence,
\begin{align}
       ||\nabla_x V_n||_{C^0\left(B_{R_n^{-1}}(0) \times [-1,0] \right)} &\lesssim_k 1, \label{eq: uniform bounded C1 norm form v_n} \\
    ||(\partial_t - \Delta) \overline{V}_n ||_{C^0\left(B_{R_n^{-1}}(0) \times [-1,0] \right)} &= ||\left( |\nabla_x \overline{V}_n|^2 \overline{V}_n \right) ||_{C^0\left(B_{R_n^{-1}}(0) \times [-1,0] \right)} \lesssim_k 1,  \notag
\end{align}
as $|\tilde{V}_n| = 1$. It follows from Parabolic Sobolev Embedding (Corollary \ref{cor: local parabolic sobolev embedding}) and a Schauder iteration as in the proof of Theorem \ref{thm:resgularty solution nonlinear pb} that
$$
\overline{V}_n, \nabla_x \overline{V}_n, ..., D^{m}_x \overline{V}_n \in C^{2+2\gamma,1+\gamma} \left(B_{R_n^{-1} 2^{-l}}(0) \times [-1,0] \right)
$$
for some $l = l(m) \in \mathbb N$ with norm uniformly bounded with respect to $n$ (for $m$ fixed).

 By Arzela-Ascoli, up to taking a subsequence, the sequence of rescaled solutions $\overline{W}_n = \overline{v}(R_n r, R_n^2t + T_n): [0,R_n^{-1}2^{-l}] \times [-1,0] \rightarrow S^2$ converges in $C^{k+2,k+1}(K \times [-1,0])$ for any compact subset $K \subset \mathbb R^2$ to some harmonic map $\overline{V}_{\infty} \in C^{k+2}(\mathbb R^2)$ which must be stationary by (\ref{eq: convergence of partial_t overline v^2 to zero in L^2}) and smooth (\cite[Theorem 3.6]{sacks-uhlenbeck}, Theorem 3.6). 

Similarly, the sequence of rescaled solutions $\overline{H}_n = \overline{h}(R_n r, R_n^2t + T_n): [0,R_n^{-1}2^{-l}] \times [-1,0] \rightarrow [0,\pi]$ converges uniformly on all compact sets $K \times [-1,0] \subset [0,+\infty) \times [-1,0]$ to a map $\overline{H}_{\infty} \in C^0([0,+\infty))$ solving
\begin{align}
     0 = H_{rr} + \frac{H_r}{r} - k^2 \frac{\sin(2H)}{2r^2}, \quad r \in (0,+\infty), \label{time independent ode for H}
\end{align}
with $H_{\infty}(0) = 0$. But Lemma \ref{k-equivariance, smoothness of h} shows that 
$$
\overline{V}_{\infty}(re^{i\theta}) = \left(e^{ik \theta} \sin H(r), \cos H(r) \right) 
$$
for some lifting $H \in C^{\infty}([0,+\infty))$. Hence, $\overline{H}_{\infty}(r)$ and $H(r)$ must differ by a fixed multiple of $2\pi$ and $\overline{H}_{\infty}$ is smooth as well. It follows that $\overline{H}_{\infty} \in C^{\infty}([0,+\infty);[0,\pi])$ is a smooth solution for the ODE (\ref{time independent ode for H}) with initial condition $H_{\infty}(0) = 0$. 

Those solutions $V_{\infty}$, $H_{\infty}$ must be non-constant. Indeed, one has
\begin{equation*}
    R_n^{-1} = |\partial_r \overline{h}(r_n,T_n)| \leq |\nabla_x \overline{v}(r_n,T_n)|
\end{equation*}
as a consequence from (\ref{nabla v in terms of h}) and
$$
R_n^{-1} = |\partial_r \overline{h}(r_n,T_n)| \lesssim r_n^{-1}
$$
by Proposition \ref{blow up point is at zero k equivariant}. Along a subsequence, $r_nR_n^{-1} \to r^* \in \mathbb R$ converges. Hence, it holds that
$$
  |\nabla_x \overline{W}_n(R_n^{-1}r_n,T_n)| \geq |R_n \nabla_x \overline{v}(r_n,T_n)| \geq 1 \quad \forall n \in \mathbb N,
$$
meaning that $|\nabla_x \overline{V}_{\infty}(r^*)| \geq 1$ is not identically zero. Hence,  $\overline{H}_{\infty}$ is non-constant as well and exactly of the form:
 $$
\overline{H}_{\infty}(r) = 2\arctan \left( (\alpha r)^k \right), \quad \alpha > 0,
 $$
 (see (\ref{all possible form for limiting H}) for a short argument based on the comparison principle). Finally, the convergence the convergence $\overline{H}_n \to \overline{H}_{\infty}$ holds in $C^{k+2,k+1}(K \times [-1,0])$ (and not only in $C^0(K \times [-1,0])$) for any compact subset $K \subset [0,+\infty)$ by locally expressing $\overline{H}_n$ in terms of $\overline{V}_n$ and the arccos, arcsin functions (e.g. see the proof of Theorem \ref{thm:smoothness around all but finitely many points}).

Finally, the limit $\overline{h}(r,+\infty)$ exists as $\overline{h}$ is increasing for fixed $r \in (0,1]$ and bounded by $0 \leq \overline{h} \leq  \pi$.  If $R \in (0,+\infty)$ is also fixed, there is $m$ large enough so that $r > R_m R$, meaning that
    $$
    0 \leq \overline{h} (R_m R, T_m + R_m^2s) \leq \overline{h}(r, T_m + R_m^2s) \leq \pi
    $$
    as $\overline{h}$ is increasing for fixed time. Taking the limit $m \to +\infty$ yields
    $$
    0 \leq 2 \arctan \left( (\alpha R)^k \right) \leq \lim \limits_{m \to +\infty} h(r, T_m + R_m^2s) \leq \pi,
    $$
    and letting $R \to +\infty$ shows that $\overline{h}(r,+\infty) = \pi$.
\end{proof}

\section{Comparison Principle, Maximum Principle}
We state two fundamental results for the study of parabolic PDEs, which are mostly independent from the previous sections and which will be used extensively in the following sections. The Comparison Principle (Theorem \ref{comparison principle}) is stated only for our specific nonlinear problem, but the proof is standard and can be applied to a wide class of nonlinear heat equations.

\begin{theorem}[Maximum Principle] \label{maximum principle}
    Let $E \subset \mathbb R^n \times \mathbb R$ be a non-empty, open and connected set. Let $E_{t_1} = \{(x,t) \in E: t \leq t_1\}$. Consider a parabolic operator
    $$
    L = \sum_{i=1}^n b^i(x,t) D_{x_i} + \sum_{i,j=1}^n a^{i,j}(x,t) D_{x_i}D_{x_j} - \partial_t,
    $$
    where $a_{i,j},b_i \in L^{\infty}(E_{t_1})$ and
    \begin{align*}
    \sum_{i,j=1}^n a^{i,j}(x,t)\xi_i\xi_j &\geq \mu |\xi|^2, \quad (x,t) \in E_{t_1}, \xi \in \mathbb R^n
    \end{align*}
    for some $\mu > 0$.

    Assume that $u \in C^{2,1}_{x,t}(E_{t_1})$ (meaning that $u$ is continuously differentiable with respect to $t$, twice-continuously differentiable with respect to $x$ and the derivatives extend continuously at $t = t_1$) and $L[u] \geq 0$. Assume that $M = \sup_{E_{t_1}} u(x,t) < +\infty$. 
    
    If $u(x_1,t_1) = M$ at some point $(x_1,t_1) \in E_{t_1}$, then $u = M$ at every point $(x,t) \in E_{t_1}$ which can be attained from $(x_1,t_1)$ by a finite sequence of horizontal line segments 
    $$H(p_0,p_1,T) = \{(x,T) \in E_{t_1}: u(p_0,T) = M, x = p_0 + h(p_1-p_0), 0 \leq h \leq 1\}$$ 
    and vertical line segments going backwards in time
    $$V(x_0,T_0,T_1) = \{(x_0,t) \in E_{t_1}: u(x_0,T_1) = M, T_0 \leq t \leq T_1\}$$ 
    that stay completely in $E_{t_1}$ (i.e., they are connected segments in the topological sense).

    The theorem remains true if $(L+c(x,t))[u] \geq 0$ with $c(x,t) \leq 0$ measurable and $M \geq 0$, as well as if $(L+c(x,t))[u] \geq 0$ with any $c(x,t)$ bounded from above and $M = 0$.
\end{theorem}

\begin{proof}
    See \cite[Chapter 3, Theorem 5, Theorem 7]{protter2012maximum} and the remark that follows it.
\end{proof}

\begin{theorem}[Comparison Principle for sub- and super-solutions, \cite{vanderhoutcomparison}]\label{comparison principle}
    Let $0 < r_0,T < +\infty$, $0 \leq t_0 < +\infty$, $R_T = (0,r_0) \times (t_0,t_0+T)$ and denote by $\Sigma_T = \partial R_T \setminus [0,r_0] \times \{t_0+T\}$ its parabolic boundary. 
    
    Assume that $\theta, \psi \in C^0(\overline{R_T}) \cap W^{1,2}_{loc}(R_T)$ (meaning that $\theta_r, \psi_r, \theta_t, \psi_t$ are defined in a weak sense and square-integrable on compact subsets of $R_T$) are respectively subsolutions and supersolutions of (\ref{heat map formulation for h(r,t)}), i.e., 
    $$
    \int_{R_T} r \theta_t \phi drdt \leq -\int_{R_T} r \left( \phi_r \theta_r + k \frac{\sin 2\theta}{2r^2} \phi \right) drdt
    $$
    and
    $$
    \int_{R_T} r \psi_t \phi drdt \geq -\int_{R_T} r \left( \phi_r \psi_r + k \frac{\sin 2\psi}{2r^2} \phi \right) drdt
    $$
    for any $\phi \in W^{1,2}(R_T)$ with compact support in $R_T$. 
    
    Suppose that either $\sup \limits_{t} E(\theta(\cdot,t)) < +\infty$ or $\sup \limits_{t} E(\psi(\cdot,t)) < +\infty$. 
    
    If $\theta \leq \psi$ on $\Sigma_T$, then $\theta \leq \psi$ on $\overline{R_T}$.
\end{theorem}

\begin{remark}
     This is a standard argument, which would also work when replacing 
    $$
    N(h,r) = k \frac{\sin 2 h}{2r^2}
    $$ 
    by any nonlinearity $N(h,r)$ for which one can prove 
    $$
    -(N(\theta(r),r) - N(\psi(r),r)) \cdot \chi_{\{\theta(r)-\psi(r) \geq 0\}}(r) \lesssim|\theta(r)-\psi(r)|, \quad r \in (0,r_0)
    $$
    given the assumptions on $\theta, \psi$. For example, take $N$ continuous on $\mathbb R \times [0,r_0]$ and Lipschitz with respect to $h$.
\end{remark}

\begin{proof}
   We follow the argument in the Appendix of \cite{vanderhoutcomparison}. After parabolic rescaling and shifting time to zero, one can assume that $R_T = (0,1) \times (0,T)$. Without loss of generality, assume that $\theta$ has finite energy (else, compare $-\theta$ and $-\psi$). Since $\theta$ is continuous on $\overline{R_T}$, Lemmas \ref{value at r = 0 is a multiple of pi} and \ref{value at r = 0 is a multiple of pi, ver.2} imply $\theta(0,t) = k \pi$ on $[0,T]$ for some $k \in \mathbb Z$. Replacing $\theta$ and $\psi$ by $\theta - k\pi$ and $\psi - k\pi$, one can assume that $\theta(0,t) = 0$.

    Before starting the proof, we introduce the positive part of a function : $u^+(x) = u(x)$ if $u(x) > 0$ and $u^+(x) = 0$ if $u(x) \leq 0$. The weak derivative $(u^+)'$ of $u^+$ is given by $(u^+)'(x) = u'(x)$ if $u(x) > 0$ and $(u^+)'(x) = 0$ if $u(x) \leq 0$. Hence, one observes that
    \begin{equation}
        u^+ \cdot u = (u^+)^2, \quad (u^+)' \cdot u = (u^+)' \cdot u^+ = (u') \cdot u^+
    \end{equation}
    and
    \begin{equation}
        0 \leq (u^+)' \cdot u' = (u^+)' \cdot (u^+)'  = \begin{cases}
            (u')^2 &\text{ if } u > 0 \\
            0 &\text{ if } u \leq 0
        \end{cases} 
    \end{equation}
    The idea is to prove 
\begin{align}
            \int_0^1 r(1-r) \left[ (\theta-\psi)^+ \right]^2(r,\tau) dr \lesssim - \int_0^{\tau}\int_0^1 r(1-r)  (\theta-\psi)^+ (N(\theta,r)-N(\psi,r))drdt \label{gronwall pre-step, comparison proof}
\end{align}
for any $0 < \tau < T$. Intuitively, this amounts to choosing $\phi = (1-r) (\theta-\psi)^+ \chi_{[0,\tau]}$ as a test function and observing that
\begin{align*}
    \int_{R_T} r(1-r)  (\theta - \psi)_t (\theta-\psi)^+ \chi_{[0,\tau]} drdt &=  \int_{R_T} r (1-r) \left[(\theta - \psi)^{+} \right]_t (\theta-\psi)^+ \chi_{[0,\tau]} drdt  \\
    &= \frac{1}{2}\int_0^1 r (1-r) \left[ (\theta-\psi)^+ \right]^2(r,\tau) dr,
\end{align*}
    while
    $$
-\int_{R_T} r(1-r)  (\theta - \psi)_r \left[ (\theta-\psi)^+ \right]_r \chi_{[0,\tau]} drdt \leq 0.
    $$
    If inequality (\ref{gronwall pre-step, comparison proof}) is proved, we can apply Grönwall's Lemma as follows and deduce that $(\theta-\psi)^+ = 0$.

    Since $\theta, \psi$ are uniformly continuous on $R_T$, $\theta - \psi \leq 0$ and $\theta + \psi \geq 0$ on $\{0\} \times [0,T]$, there is $\delta > 0$ for which
    $$
    (\theta-\psi)^+ \leq \frac{\pi}{2}, \quad (\theta + \psi) \chi_{\{\theta-\psi \geq 0\}} \leq \frac{\pi}{2} \quad (r,t) \in [0,\delta] \times [0,T].
    $$
    Then
    \begin{align*}
         \frac{\sin 2\theta - \sin 2\psi}{2r^2} \chi_{\{\theta-\psi \geq 0\}} = \frac{\cos(\theta+\phi)\sin(\theta-\phi)}{r^2} \chi_{\{\theta-\psi \geq 0\}} \geq 0 \quad (r,t) \in [0,\delta] \times [0,T] 
    \end{align*}
    and
    $$
 \left| \frac{\cos(\theta+\phi)\sin(\theta-\phi)}{r^2} \chi_{\{\theta-\psi \geq 0\}} \right|  \leq \frac{(\theta-\phi)^+}{\delta^2} \quad (r,t) \in [\delta,1] \times [0,T], 
    $$
    meaning that
\begin{align*}
    - \int_0^{\tau}\int_0^1 r (1-r) (\theta-\psi)^+ \frac{\sin 2\theta - \sin 2\psi}{2r^2}drdt &= - \int_0^{\tau}\int_0^{\delta} ... - \int_0^{\tau}\int_{\delta}^1 ... \leq - \int_0^{\tau}\int_{\delta}^1 ... \\
    &\leq \frac{1}{\delta^2} \int_0^{\tau} \int_{\delta}^1 r(1-r)\left[(\theta-\psi)^+\right]^2 dr dt \\
    &\leq \frac{1}{\delta^2} \int_0^{\tau} \int_{0}^1 r(1-r)\left[(\theta-\psi)^+\right]^2  dr dt
\end{align*}
which concludes the proof.

The issue is that we cannot directly use the test function $\phi = (1-r)(\theta-\psi)^+ \chi_{[0,\tau]}$. Instead, let $\phi_n^{a}(s) \in C^{\infty}_c(\mathbb R)$ be a smooth cut-off function satisfying
    \begin{enumerate}
        \item $\phi_n^a(s) = 1, \quad s \in [2/n, a-2/n]$ 
        \item $\phi_n^a(s) \in [0,1], \quad s \in [1/n,2/n] \cup [a-2/n,a-1/n]$ 
        \item $\phi_n^a(s) = 0, \quad s \in [1/n,a-1/n]^c$ 
        \item $|\partial_s \phi_n^a(s)| \leq 2n, \quad s \in \mathbb R$
    \end{enumerate}
    and consider the test function $\varphi_{m,n} = (1-r)(\theta-\psi)^+(r,t) \phi^1_n(r)^2\phi^{\tau}_m(t)$. Then
 \begin{align*}
      \int_0^1 r (1-r)\left[ (\theta-\psi)^+ \right]^2(r,\tau) dr &= \int_{R_T} r (1-r)(\theta - \psi)_t (\theta-\psi)^+ \chi_{[0,\tau]} drdt  \\
      &= \lim_{n \to +\infty} \lim_{m \to +\infty} \int_{R_T} r (\theta - \psi)_t \varphi_{m,n} drdt \\
      &\leq - \lim_{n \to +\infty} \lim_{m \to +\infty} \int_{R_T} r (\theta - \psi)_r \partial_r \varphi_{m,n} drdt \\
      &\phantom{\leq}- k \int_0^{\tau}\int_0^1 r (1-r) (\theta-\psi)^+ \frac{\sin 2\theta - \sin 2\psi}{2r^2}drdt
 \end{align*}
 and it suffices to prove that the additional first term with $r$-derivatives is non-positive. It is made of three terms
\begin{align*}
I_1 &= -\int_{R_T} r (1-r) (\theta - \psi)_r \left[ (\theta-\psi)^+ \right]_r \phi^1_n(r)^2 \phi^{\tau}_m(t) drdt  \\
   I_2 &=  -\int_{R_T} r (1-r)(\theta - \psi)_r (\theta-\psi)^+  \phi^1_n(r) \partial_{r}\phi^1_n(r) \phi^{\tau}_m(t) drdt\\
    I_3 &= \int_{R_T} r (\theta - \psi)_r (\theta-\psi)^+  \phi^1_n(r)^2 \phi^{\tau}_m(t) drdt,
\end{align*}
where
\begin{align*}
       2|I_2| &= 2\left| \int_{R_T} r  (1-r)\left[ (\theta - \psi)^+ \right]_r (\theta-\psi)^+  \phi^1_n(r) \partial_{r}\phi^1_n(r) \phi^{\tau}_m(t) drdt \right| \\
    &\leq \int_{R_T} r  (1-r)\left[ (\theta - \psi)^+ \right]_r^2 \phi^{1}_n(r)^2 \phi^{\tau}_m(t)drdt  \\
    &+  \int_{R_T} r (1-r) \left[ (\theta - \psi)^+ \right]^2 \left[ \partial_r\phi^{1}_n(r)\right]^2 \phi^{\tau}_m(t)drdt \\
    &\leq -I_1 + I_2'.
\end{align*}
Hence, the first term $I_1$ cancels out because of $I_2$. As for $I_2'$, 
\begin{align*}
     |I_2'| &\leq 2 \int_0^T \int_{\substack{0 \leq r \leq 2/n \\ 1-2/n \leq r \leq 1}} r (1-r) \left[ (\theta - \psi)^+ \right]^2 n^2 dr dt \\
     &\leq 2 \int_0^T \left[ (\theta - \psi)^+ \right]^2(r(n,t),t)dt + 2 \int_0^T \left[ (\theta - \psi)^+ \right]^2(\tilde{r}(n,t),t)dt,
\end{align*}
where $r(n,t) \in [0,2/n]$ and $\tilde{r}(n,t) \in [1-2/n,1]$ by the Mean Value Theorem. Since $(\theta-\psi)^+$ is uniformly continuous on $R_T$ and vanishes at the boundary $\{r = 0\} \cup \{r = 1\}$, letting $n \to +\infty$ shows that $I_2' = 0$. 

After an integration by parts for $I_3$, the boundary term vanishes and the other integral is treated as $I_2'$.
\end{proof}

\section{On the geometry of level sets $h_1 < h_2$}

In this section, following the lap-number argument from \cite{VANDERHOUT2003} and \cite{matano}, we are concerned on the level sets $h_1 < h_2$, as well as $h_1 < h < h_2$. We show that under some mild assumptions, those levels sets are connected to the parabolic boundary (Lemma \ref{geometry of h_1 < h_2} and Lemma \ref{lemma:geometry of h_1 < h < h_2}). Moreover, at a fixed time, there cannot be too many successive intersections between $h$ and $h_i$ on $[0,1]$ and this intersection number cannot increase over time (Lemma \ref{lemma:intersection argument with chains mod 4}. These results play a key role in the bubble-tree result.

Even though the results are formulated in greater generality, we will mostly be concerned with the cases where $h_i = \pi/2, \pi$ or $h_i$ is a stationary solution given in terms of the arctan function.

We note that the proof of Lemma \ref{lemma:geometry of h_1 < h < h_2} is novel. The original paper from \cite{VANDERHOUT2003} omits this proof but one cannot repeat the argument from Lemma \ref{geometry of h_1 < h_2}. Indeed, repeating the argument by contradiction from Lemma \ref{geometry of h_1 < h_2} shows that the maximum of $h-h_1$, resp. $h_2-h$, is achieved on $\overline{\Gamma} \setminus \Gamma$. However, one cannot conclude as $h = h_1$ or $h = h_2$ on $\overline{\Gamma} \setminus \Gamma$ and both can happen at different points because this set is not necessarily connected. We propose a more subtle argument based on the notion of accessible boundary points from \cite{newman1964elements}. These points, which are dense in the boundary, allow to circumvent the aforementioned problem. Even though the boundary is not connected, most boundary points can be reached from inside the domain and are connected together.

\begin{lemma}[\cite{VANDERHOUT2003}, Lemma 2.1] \label{geometry of h_1 < h_2}
Let $k \in \mathbb N$ and let $h_i(r,t) \in C^0([0,1] \times [0,T)) \cap C^{2,1}((0,1) \times (0,T))$ (meaning that $u$ is continuously differentiable with respect to $t$, twice-continuously differentiable with respect to $r$ and the derivatives are continuous at the included boundary) , $i = 1,2$, be two classical solutions for the nonlinear problem
\begin{align*}
    h_{i,t} &= h_{i,rr}+\frac{h_{i,r}}{r} - k^2 \cdot \frac{\sin(2h_i)}{2r^2}, \quad 0 < r < 1, \quad 0 < t < T,
\end{align*}
with a common time of existence $T \leq +\infty$. Let $0 \leq t_1 < t_2 < T$ and $0 \leq R_1 < R_2 \leq 1$. Let $A$ be a (path)-connected component of $\{(r,t) \in [R_1,R_2] \times [t_1,t_2]: h_1 > h_2\}$. Then
$$
A \cap \left[ \{r = R_1\} \cup \{r = R_2\} \cup \{t = t_1\} \right] \neq \emptyset,
$$
where, if $R_1 = 0$, one needs to further assume either $h_1(0,t) \neq h_2(0,t)$ for all $t \in [0,T)$ or $h_1(0,t) = h_2(0,t) = m\pi \in \pi \mathbb Z$ for all $t \in [0,T)$.
\end{lemma}

\begin{notation}
    In the following, we write $\overline{S}$, $\partial_{\mathbb R^2} S$, $\mathrm{int}_{\mathbb R^2}(S)$ for the closure, boundary and interior of $S \subset [R_1,R_2] \times [t_1,t_2]$ with respect to $\mathbb R^2$ topology. One should note that the closure $\overline{S}$ with respect the $\mathbb R^2$ topology or the closed square topology is the same, but the boundary might differ. 
\end{notation}

\begin{remark}
    In our Euclidean setting, connectedness, path-connectedness and arc-connectedness (i.e., existence of injective paths between any two points) are all equivalent.
\end{remark}

\begin{proof}
We follow the proof in \cite[Lemma 2.1]{VANDERHOUT2003}. Assume for a contradiction that this is not the case.

     First, observe that $A$ is open in the closed square topology. Moreover, $h_1 = h_2$ on $\overline{A} \setminus A$ (i.e., the boundary with respect to the closed square topology). If we had $h_1 > h_2$ at a point $(r,t) \in \overline{A} \setminus A$, then there would be an open (in the closed square topology) path-connected neighborhood $U \subset [R_1,R_2] \times [t_1,t_2]$ of $(r,t)$ on which $h_1 > h_2$. Since $(r,t)$ belongs to the boundary of $A$ (in the closed square topology), such a neighborhood must intersect $A$ by definition, meaning we can find a path from $A$ to $(r,t) \notin A$ on which $h_1 > h_2$. This contradicts the maximality of the connected component $A$.

    We consider $H = (h_1-h_2)$. Then
    $$
     H_{t} = H_{rr}+\frac{H_{r}}{r} - k^2 \cdot H \cdot \frac{\sin(2h_1)-\sin(2h_2)}{2r^2(h_1-h_2)}, \quad (r,t) \in A,
    $$
    which is a radial heat equation in 2D. If $R_1 > 0$, the coefficients in front of $H$ is bounded. If $R_1 = 0$ and $h_1(0,t) \neq h_2(0,t)$, then $\mathrm{dist}(\overline{A},\{r = 0\}) > 0$. Otherwise, there would be $x \in \left( \overline{A} \setminus A \right)\cap \{r = 0\}$ which is not possible because $h_1 = h_2$ on $\overline{A} \setminus A$, while $h_1 \neq h_2$ at $r = 0$. In that case, the coefficient in front of $H$ is bounded as well. If $R_1 = 0$ and $\mathrm{dist}(\overline{A},\{r = 0\}) > 0$, no matter $h_1$ and $h_2$, the situation is the same. Finally, if $R_1 = 0$ and $h_1(0,t) = h_2(0,t) = m\pi$ for all $t \in [0,T)$ while $\mathrm{dist}(\overline{A},\{r = 0\}) = 0$, then replacing $h_i$ by $h_i - m\pi$, we can assume that $m = 0$. Near $r = 0$, given $t_1, t_2$, there is $r_0$ for which $|2h_i(r,t)| \leq \pi/2$ for $(r,t) \in [0,r_0] \times [t_1,t_2]$ by uniform continuity. On $\overline{A} \cap  [0,r_0] \times [t_1,t_2]$, one has
    $$
    \frac{\sin(2h_1)-\sin(2h_2)}{2r^2(h_1-h_2)}  \geq 0,
    $$
    as $h_1 \geq h_2$ on $\overline{A}$ and $\sin(x)$ is increasing on $[-\pi/2,\pi/2]$. Away from $r_0$, the coefficient is bounded. 

    All in all, the coefficient in front of $H$ is bounded from above. Hence, one can apply the maximum principle (Theorem \ref{maximum principle}) on $e^{-ct}H$ for $c \geq 0$ large enough so that the coefficient in front of $e^{-ct}H$ is bounded from above by zero.
    
    Let $M = \sup_A e^{-ct}H$ be achieved at a point $(r_0,t_0) \in \overline{A}$. As $A$ is non-empty, $M > 0$. We prove that the maximum is also achieved on $\overline{A} \setminus A$, meaning that $M \leq 0$ on $A$: a contradiction.

    If $M$ is achieved on $\overline{A} \setminus A$, we are done. Moreover, $M$ cannot be achieved on the empty set
    $$
    A \cap \left[ \{r = R_1\} \cup \{r = R_2\} \cup \{t = t_1\} \right] = \emptyset.
    $$
    Hence, we can assume that $t_1 < t_0 \leq t_2$ and $R_1 < r_0 < R_2$. we will prove that a maximum cannot occur at such point. Since $A$ is open with respect to the closed square topology, one can find a lower half-ball
    $$
    \{(r-r_0)^2 + (t-t_0)^2 < \varepsilon, t \leq t_2\} \subset A \cap \{t_1 < t \leq t_2, R_1 < r < R_2\}
    $$
    for $\varepsilon > 0$ small enough. On this half-ball, the maximum principle (Theorem \ref{maximum principle}) applies and $e^{-ct}H = M$ is constant. Hence, we can assume that $t_1 < t_0 < t_2$ without loss of generality. Now, the same theorem implies that $e^{-ct}H = M$ for any point of $A$ which can be connected to $(r_0, t_0)$ using a horizontal segment which stays inside $A$. Take such a point $(r,t_0)$ with $r \in [R_1,R_2]$ minimal. This point must be in $\overline{A}$ (because it is the limit of a sequence in $A$) and it is connected to the point $(r_0,t_0)$ with an horizontal segment staying completely in $A$ except for the endpoint $(r,t_0)$. Since $A$ is a connected component and $e^{-ct}H = M > 0$ at $(r,t_0)$, one must have $(r,t_0) \in A$. And it cannot be that $r > R_1$ since $A$ is open (in the square topology) but $r$ is minimal. As $A \cap \{r = R_1\} = \emptyset$, we reach a contradiction.
\end{proof}

\begin{lemma}
    Let $A \subset [0,R] \times [t_1,t_2]$ be open (with respect to the subspace topology). Then 
    $$\mathrm{int}_{\mathbb R^2}(A) = A \setminus \partial_{\mathbb R^2} ([0,R] \times [t_1,t_2]), \quad A \subset \overline{\mathrm{int}_{\mathbb R^2}(A)}, \quad \overline{A} = \overline{\mathrm{int}_{\mathbb R^2}(A)}, \quad \partial_{\mathbb R^2} \mathrm{int}_{\mathbb R^2}(A) = \partial_{\mathbb R^2} A
    $$
    Moreover, if $A$ is connected, then $\mathrm{int}_{\mathbb R^2}(A)$ is connected as well.
\end{lemma}

\begin{proof}
     This is an elementary topological argument.
\end{proof}

\begin{definition}[Accessible points]
Let $U \subset \mathbb R^2$ be an open and connected domain. A boundary point $x \in \partial U$ is called \textit{accessible} if there exists a simple curve $\gamma: [0,1] \rightarrow \overline{U}$ for which $\gamma(0) = x$ and $\gamma((0,1]) \subset U$.
\end{definition}

\begin{lemma}[\cite{newman1964elements}, Chapter 6.4]
    The subset of accessible points $A \subset \partial U$ is dense in $\partial U$.
\end{lemma}

\begin{proof}
    Let $x \in \partial U$ and $\varepsilon > 0$ be fixed. By definition, $B(x,\varepsilon) \cap U \neq \emptyset$ and there is some $y \neq x$ in $B(x,\varepsilon) \cap U$. Consider the segment $\{y + s(x-y): s \in [0,1]\} \subset B(x,\varepsilon)$ going from $y$ to $x$. The first point $y + s(x-y)$ which intersects the boundary $\partial U$ is accessible.
\end{proof}

\begin{lemma}\label{lemma:geometry of h_1 < h < h_2}
     Let $0 \leq t_1 < t_2 < T$ and $0 < R \leq 1$. Let $(h, h_1)$ and $(h,h_2)$ be two pairs of solutions for which Lemma \ref{geometry of h_1 < h_2} applies on each pair on the square $[0,R] \times [t_1,t_2]$. If $\Gamma$ is a non-empty connected component of $\{(r,t) \in [0,R] \times [t_1,t_2]: h_1 < h < h_2\}$ and if
    $$
    \max_{(x,t) \in \overline{\Gamma}} h_1 < C < \min_{(x,t) \in \overline{\Gamma}} h_2
    $$
    for some $C \in \mathbb R$, then 
    $$
    \Gamma \cap \left[ \{r = 0\} \cup \{r = R\} \cup \{t = t_1\} \right] \neq \emptyset.
    $$
    No matter $\Gamma$, one can always choose $h$ with $h(0,t) = 0$ for all $t \in [0,T)$ satisfying the regularity from Lemma \ref{geometry of h_1 < h_2}, as well as $h_1 = \pi/2$ or $h_1 = \chi_{\alpha}(r) = \pi - 2 \arctan(\alpha^k r^k)$, $\alpha > 0$, together with $h_2 = \pi$. 
\end{lemma}

\begin{proof}
This corresponds to \cite[Corollary 2.2]{VANDERHOUT2003}, in which the proof is omitted. We believe that the previous proof cannot be repeated and we need a more subtle topological argument in this case using the notion of accessible points.

First, observe that either $h = h_1$ or $h = h_2$ on $\overline{\Gamma} \setminus \Gamma$ (as in Lemma \ref{geometry of h_1 < h_2}).

In the case $h_1 = \chi_{\alpha}$ and $h_2 = \pi$, observe that $\mathrm{dist}(\overline{\Gamma}, \{r = 0\}) > 0$ (because $h \neq \chi_{\alpha}$ and $h \neq \pi$ on $\{r = 0\}$), meaning that 
$$
\chi_{\alpha}(r) \leq \max_{(r,t) \in \overline{\Gamma}} \chi_{\alpha}(r) < C < \min_{(r,t) \in \overline{\Gamma}} \pi = \pi, \quad (r,t) \in \overline{\Gamma}
$$
on $\overline{\Gamma}$ for some constant $C(\Gamma)$. The same holds with $h_1 = \pi/2$.

Now, suppose for a contradiction that the statement of the lemma is not true. We observe that if $h = h_1$ (resp. $h = h_2$) on the whole $\overline{\Gamma} \setminus \Gamma$, then $\Gamma$ is actually a connected component of $\{h >  h_1\}$ (resp. $\{h < h_2 \}$) and is connected to the desired boundary of the square. Hence, such a case does not happen.


Our previous observation allows us to assume that there is a point $(r_0,s_0)$ in $ \partial_{\mathbb R^2}\Gamma$ with $h = h_1$ and another point $(r_1,s_1)$ with $h = h_2$. Slightly perturbating these points, we can find accessible boundary points of $\mathrm{int}_{\mathbb R^2}(\Gamma)$ satisfying $h(r_0,s_0) < C < h(r_1,s_1)$. These two points can be connected through a continuous path (it does not matter whether it is injective) lying entirely in $\mathrm{int}_{\mathbb R^2}(\Gamma) \subset \{t < t_2\}$ except for those endpoints. By continuity, there exists a point $(r,t)$ on that curve for which $h(r,t) = C$ and $t < t_2$. Now, take any point $(r^*,t^*) \in \overline{\mathrm{int}_{\mathbb R^2}(\Gamma)}$ which achieves the infimum
$$
\inf\{t_1 \leq t < t_2: (r,t) \in \mathrm{int}_{\mathbb R^2}(\Gamma), h(r,t) = C\}.
$$
Then $(r^*,t^*) \notin \partial_{\mathbb R^2}  \mathrm{int}_{\mathbb R^2}(\Gamma) = \partial_{\mathbb R^2} \Gamma$. Otherwise, slightly perturbating $(r^*,t^*)$ leads to an accessible boundary point $(r,t) \in \partial_{\mathbb R^2} \Gamma$ with $h_1 < h(r,t) \approx C < h_2$, hence in $\partial_{\mathbb R^2}  \Gamma\cap \Gamma$. This means that either $r = R$ or $t = t_1$ (which contradicts our absurd hypothesis), either $0 < r < R$ and $t > t_1$, which is not possible by maximality of the (path-)connected component $\Gamma$.

Hence, $(r^*,t^*) \in  \mathrm{int}_{\mathbb R^2}(\Gamma)$. In particular, $t_2 > t^* > t_1$ and any point 
$$
 \mathrm{int}_{\mathbb R^2}(\Gamma) \cap \{t_1 \leq t < t^*\}
$$
satisfies $h(r,t) \neq C$. Take the connected component of the intersection between the open, connected domains 
$$
 \mathrm{int}_{\mathbb R^2}(\Gamma) \cap (0,R) \times (t_1,t^*),
$$
which contains $(r^*,t^*-\varepsilon)$ for any $\varepsilon > 0$ small enough.

Such a connected component is part of a larger connected component $(r^*,t^*) \in \tilde{\Gamma} \subset \Gamma$ of $\{(r,t) \in [0,R] \times [t_1,t^*]: h_1< h < h_2\}$. If 
$$
\tilde{\Gamma} \cap \left[ \{r = 0\} \cup \{r = R\} \cup \{t = t_1\} \right] \neq \emptyset,
$$
then $(r^*,t^*)$ can be connected to the desired boundary: a contradiction. Otherwise, reiterating the above argument with $\tilde{\Gamma}$ instead of $\Gamma$, we find a point $(r,t)$ in $\mathrm{int}_{\mathbb R^2}(\tilde{\Gamma}) \subset \mathrm{int}_{\mathbb R^2}(\Gamma)$ with $h(r,t) = C$ and $t_1 < t < t^*$, which contradicts the minimality of $t^*$.
\end{proof}

\begin{lemma}[\cite{VANDERHOUT2003}, Lemma 2.3]\label{lemma:intersection argument with chains mod 4}
     Let $h,h_1,h_2$ be defined on $[0,T)$ and satisfy Lemma \ref{lemma:geometry of h_1 < h < h_2} with $R = 1, t_1 = 0$, any $0 < t_2 < T$ and any $\Gamma$. Further assume that $h$ has finite energy (\ref{Energy of v in terms of h}), $h(0,t) < h_1(0,t)$, $h(1,t) > h_1(1,t)$ on $[0,T)$ and $h_1(r,t) \leq h_2(r,t)$ for all $r \in (0,1]$ and $t \in [0,T)$.

     There exists $M(h) \in \mathbb N_{\geq 1}$ such that for all $t \in [0,T)$, any sequence $0 \leq r_1 < r_2 < ... < r_N \leq 1$ with the property
     $$
     P(t,N)
     = \begin{cases}
         N = 1 \mod 4; \\
         h(r_{4j+1},t) < h_1(r_{4j+1},t); \quad h_1(r_{4j+2},t) < h(r_{4j+2},t) < h_2(r_{4j+2},t); \\
         h_2(r_{4j+3},t) < h(r_{4j+3},t); \quad  h_1(r_{4j},t) < h(r_{4j},t) < h_2(r_{4j},t);
     \end{cases}
    $$
    one has $N \leq M$. Moreover, if $M(t,h,h_{1},h_2)$ denotes the maximal length of such sequence, then $N$ is non-decreasing with respect to $t$. 
    
    If $h_1(r,t) = h_1(r,t,\alpha)$ depends on a parameter $\alpha > 0$ for which $h_1(r,t,\alpha_1) \geq h_1(r,t,\alpha_2)$ whenever $\alpha_1 \leq \alpha_2$, if the hypotheses of the Lemma for $h,h_1,h_2$ are satisfied for any fixed $\alpha > 0$ with a strict inequality $h_1(r,t,\alpha) < h_2(r,t)$ when $r \in (0,1]$, $t \in [0,T)$, then $M(t,h,h_1(\cdot,\alpha),h_2)$ is non-decreasing with respect to $\alpha$ as well. 
    
    Finally, a sequence satisfying $P(t,1)$ always exists.
\end{lemma}

\begin{proof}
We follow the proof in \cite[Lemma 2.3]{VANDERHOUT2003}. As $h(0,t) < h_1(0,t)$, a trivial sequence satisfying $P(t,1)$ always exists. As $h(1,t) > h_1(1,t)$, there is also 
    $$
    0 < r_1^*(t) = \sup \{ r \in [0,1]: h(s,t) < h_1(s,t) \ \forall s \in [0,r] \} < 1
    $$
for which $h(r_1^*(t),t) = h_1(r_1^*(t),t)$, $h(r,t) < h_1(r,t)$ on $[0,r_1^*(t))$. Hence, any sequence $r_1,...,r_N$ with $N > 1$ satisfying $P(t,N)$ must have $r_2 > r_1^*(t)$, as well as $r_N < 1$.

For $t = 0$, there exists $C = C(h) > 0$ for which
$$
 E(r_{j},h(\cdot,0)) = \pi \int_0^{r_{j}} \left( h_{r}^2 + \frac{k^2}{r^2} \sin(h)^2 \right) r dr \geq jC \quad \forall j \geq 2,
$$
if $r_1, ..., r_N$ satisfies $P(t,N)$ with $N \geq 2$. Otherwise, for all $C > 0$, there would be some $j \geq 2$ with
$$
 E(r_1^*(0),h(\cdot,0)) \leq E(r_{j},h(\cdot,0)) = \pi \int_0^{r_{j}} \left( h_{r}^2 + \frac{k^2}{r^2} \sin(h)^2 \right)  r dr \leq jC.
$$
Letting $C \to 0$, we find $E(r_1^*(0),h(\cdot,0)) = 0$ which implies that $h(r,0)$ is constant on $[0,r_1^*(0)]$. By hypothesis, $h(0,0) < h_1(0,0)$ but $h(r_1^*(0),0) = h_1(r_1^*(0),0)$, which is a contradiction. In particular, 
$$
N \leq C(h)^{-1} \sup_{t \in [0,T]}E(h(\cdot,t)) = M(h)
$$
for any such sequence.

Finally, we show that if $r_1,...,r_N$ satisfies $P(t_2,N)$ for $t_2 \in (0,T)$, then one can always find a sequence satisfying $P(t_1,N)$ with $0 \leq t_1 < t_2$, which proves that the maximal length is non-decreasing and bounded by the one for time $t = 0$. As $h(0,t) < h_1(0,t)$  and $h, h_1$ are uniformly continuous on $[0,1] \times [t_1,t_2]$, we may always replace $r_1$ in the chain by a sufficiently small value $r_1 \in (0,r_2)$ which satisfies $h(r_1,t) < h_1(r_1,t)$ for all $t \in [t_1,t_2]$.

We define first an ordering on the boundary $\{r = 0\} \cup \{r = 1\} \cup \{t = t_2\}$ of $[0,1] \times [t_2,t_1]$. We say that $(r_1,s_1) < (r_2,s_2)$ whenever $(r_2,s_2)$ 'lies on the right' of $(r_1,s_1)$ along the parabolic boundary. More precisely, $(r_1,s_1) < (r_2,s_2)$ when $r_1 < r_2$, $s_1 = s_2$ or $r_1 = r_2 = 0$, $s_1 > s_2$ or $r_1 = r_2 = 1$, $s_1 < s_2$.

By Lemma \ref{lemma:geometry of h_1 < h < h_2} applied on the square $[0,1] \times [t_1,t_2]$, for each $(r_i,t_2)$, there is an injective path $\Sigma_i$ connecting $(r_i,t_2)$ to $\xi_i \in \{r = 0\} \cup \{r = 1\} \cup \{t = t_1\}$ while preserving the strict inequality between $h$ and $h_1,h_2$ associated to $(r_i,t_2)$. Given our hypotheses on $r_1$, $\Sigma_1$ can always be chosen as a vertical path from $(r_1,t_2)$ to $(r_1,t_1)$. The important point being that $\Sigma_1$ does not intersect $\{r = 0\}$ and $\{t = t_2\}$ (except for the initial point $(r_1,t_2)$).

\begin{figure}
    \centering
    \includegraphics{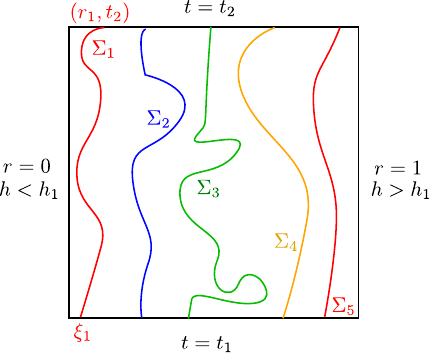}
    \caption{\label{fig:figure1} Sequence of five non-intersecting paths going from $(r_i,t_2)$ to $\xi_i \in \{t = t_1\}$}
    \label{fig:sigma_paths}
\end{figure}

Without loss of generality, one can assume that the others paths $\Sigma_i$, $i \neq 1$, do not intersect $\{t = t_2\}$ except for the initial point $(r_i,t_2)$. Indeed, the strict inequality between $h$ and $h_1,h_2$ is still satisfied at $(r_i,t)$ for all $t$ sufficiently close to $t_2$. We may connect $(r_i,t_2)$ to $(r_i,t)$ through a vertical line and then $(r_i,t)$ to the parabolic boundary through Lemma $\ref{lemma:geometry of h_1 < h < h_2}$ applied on $[0,1] \times [t_1,t]$. By choosing a subpath if needed, one can assume that $\xi_i$ is the first intersection of $\Sigma_i$ with the parabolic boundary of $[0,1] \times [t_1,t]$. In particular, $\Sigma_i \subset (0,1) \times (t_1,t_2)$ except for the endpoints $(r_i,t_2)$ and $\xi_i$. These properties are also satisfied by $\Sigma_1$.

Observe next that $(r_{4j + 1},t_2)$ can only be connected to $\{r = 0\} \cup \{t = t_1\}$ while $(r_{4j+2},t_2)$, $(r_{4j+3},t_2)$ and $(r_{4j},t_2)$ can only be connected to $\{t = t_1\} \cup \{r = 1\}$. We also know that $(r_1,t_2)$ is connected to $\xi_1 \in \{t = t_1\}$.

Clearly, two consecutive paths $\Sigma_{i}, \Sigma_{i+1}$ associated to $(r_i,t_2)$, $(r_{i+1},t_2)$ must be non-intersecting. If $C_i$ is the interior of the square enclosed by $\Sigma_i$ and the left-oriented path going from $(r_i,t_2)$ to $\xi_i$ through $\partial([0,1] \times [t_1,t_2])$, it follows from the Jordan Curve Theorem that either $\Sigma_{i+1} \subset C_i$ or $\Sigma_{i+1} \subset (\overline{C_i})^c$, which are two disjoint, open path-connected sets. As $(r_{i+1},t_2)$ lies on the right of $(r_i,t_2)$, it can only be the latter. In particular, $\xi_{i+1}$ must be in $(\overline{C_i})^c$ as well, meaning it lies on the right of $\xi_i$. By induction, $\xi_1 < \xi_2 < ... < \xi_N$. As $\xi_1 \in \{t = t_1\}$ and $\xi_N \in \{r = 0\} \cup \{t = t_1\}$ lies on the right of $\xi_1$, we must have $\xi_N \in \{t = t_1\}$. It then follows from the ordering that $\xi_i \in \{t = t_1\}$ for all $i = 1, ..., N$ and $\xi_1 < \xi_2 < ... < \xi_N$ satisfies $P(t_1,N)$ by construction.

If $h_1$ depends on a parameter $\alpha> 0$, then the maximal number $N(t,h,h_1(\cdot,\alpha),h_2)$ is well-defined for fixed $\alpha > 0$. Let $r_1,...,r_N$ satisfy $P(t,N,\alpha_2)$ (i.e., $P(t,N)$ with $h$, $h_1 = h_1(\cdot,\alpha_2)$, $h_2$) and we find a $s_1, ..., s_N$ sequence satisfying $P(t,N,\alpha_2)$ whenever $0 < \alpha_1 < \alpha_2$. Without loss of generality, $r_i > 0$ for all $i$.

As our initialization, we set $s_{4j+1} = r_{4j+1}$, $s_{4j+3} = r_{4j+3}$ for all $j$. Since
$$
h(s_{4j+1},t) < h_1(s_{4j+1},t, \alpha_2) \leq  h_1(s_{4j+1},t, \alpha_1), \ 
 h_1(s_{4j+3},t, \alpha_1) < h_2(s_{4j+3},t) < h(s_{4j+3},t),
$$
there is $s^* \in (s_{4j+1},s_{4j+3})$ maximal with $h(s^*,t) = h_1(s^*,t,\alpha_1)$, $h(s,t) > h_1(s^*,t,\alpha_1)$ on $(s^*,s_{4j+3})$. As $h_2(s^*,t) > h_1(s^*,t,\alpha_1) = h(s^*,t)$, one can set $s_{4j+2} \in (s^*,s_{4j+3})$ close enough to $s^*$ so that $h_2(s_{4j+2},t) > h_1(s_{4j+2},t,\alpha_1)$, $h_2(s_{4j+2},t) >  h(s_{4j+2},t)$ and $h(s_{4j+2},t) > h_1(s_{4j+2},t,\alpha_1)$.

Similarly, since
$$
h_1(s_{4j+3},t, \alpha_1) < h_2(s_{4j+3},t) < h(s_{4j+3},t), \ h(s_{4j+5},t) < h_1(s_{4j+5},t, \alpha_2) \leq  h_1(s_{4j+5},t, \alpha_1),
$$
there is $s^* \in (s_{4j+3},s_{4j+5})$ minimal with $h(s^*,t) = h_1(s^*,t,\alpha_1)$, $h(s,t) > h_1(s^*,t,\alpha_1)$ on $(s_{4j+3}, s^*)$. As $h_2(s^*,t) > h_1(s^*,t,\alpha_1) = h(s^*,t)$, one can set $s_{4j} \in (s_{4j+3},s^*)$ close enough to $s^*$ so that $h_2(s_{4j},t) > h_1(s_{4j},t,\alpha_1)$, $h_2(s_{4j},t) >  h(s_{4j},t)$ and $h(s_{4j},t) > h_1(s_{4j},t,\alpha_1)$. This finishes the proof.
\end{proof}




\section{Nonexistence of multi-bubbles}
    Let $v(x,t)$ solve (\ref{heat map flow}) with smooth and $k$-equivariant boundary data, $k \geq 1$. Assume that $v$ blows-up at time $T < +\infty$. In particular, $v$ has finite energy (Remark \ref{k-equivariance implies finite energy}). Let $h$ be the corresponding inclination coordinate chosen so that $h(0,t) = 0$. We give a further description of the blow-up event using the intersection-comparison argument from \cite{VANDERHOUT2003} and prove that there can only be a single bubble in the asymptotic decomposition (\ref{multi bubble decomposition along a sequence of times}). Recall from Proposition \ref{prop:smoothness at the boundary, earlier statement} that $v(x,T)$ and $h(r,T)$ are well-defined except at the origin.

    The starting point is Proposition \ref{prop:criterion for global solution} and Corollary \ref{cor:comparison of h with itself and exact value of limsup}. From those two results, one deduces (up to replacing $h$ by $-h$) that 
    $$
    \limsup_{\substack{r \to 0^+ \\ t \to T^-}} h(r,t) = \pi, \quad \lim \limits_{r \to 0^+} h(r,T) \in \{0,\pi\},
    $$
    and that there are points with $h(r,t) > \pi$ for each $t$.
    
    The difficult part is actually showing that $\lim \limits_{r \to 0^+} h(r,T) = \pi$. Together with Corollary \ref{cor:comparison of h with itself and exact value of limsup}, one deduces Corollary \ref{cor: liminf h(r,t) = 0}, i.e.,
    $$
    \liminf_{\substack{r \to 0^+ \\ t \to T^-}}  h(r,t) = 0.
    $$
    The idea is that if $h(0^+,T) = 0$ (hence $h(r,T) \leq \pi/2$ w.l.o.g.) and if $t_0,\alpha_0$ are large enough so that the maximal length $M(t,h,h_1 = \chi_{\alpha} = \pi - 2\arctan(\alpha^k r^k),h_2 = \pi)$ of a chain $r_1(\alpha,t) < ... < r_M(\alpha,t)$ satisfying $P(t,h,h_{1} = \chi_{\alpha},h_2 = \pi,M)$ from Lemma \ref{lemma:intersection argument with chains mod 4} is minimal for $t \geq t_0, \alpha \geq \alpha_0$, then $M \geq 5$. Intuitively, this is because if we had $M = 1$, then there would be a curve $h(r(t),t) = \chi_{\alpha}(r(t))$ splitting the open square into two regions: $h \leq \chi_{\alpha}$ on the left of the curve and $h \geq \chi_{\alpha}$ on the right. One must have $\lim \limits_{t \to T^-} r(t) = 0$ as otherwise, Proposition \ref{prop:criterion for global solution} applied on a sub-rectangle of the region $h \leq \chi_{\alpha}$ would imply that $h$ is global. Thus $\{t = T\}$ must be on the right of the curve, contradicting $h(0^+,T) = 0$.
    
    Hence, one can let $\alpha^* > \alpha_0$ large enough such that $0 < \chi_{\alpha^*}(r_{M-2}(\alpha_0,t_0)) < h(r_{M-2}(\alpha_0,t_0),t_0) - \pi$ and let
   \begin{align*}
         R_1 &= \inf\{r \in (0,r_{M-2}(\alpha_0,t_0)]: h(s,t_0) - \pi > \chi_{\alpha^*}(s) \quad \forall s \in (r,r_{M-2}(\alpha_0,t_0)]\}, \\
  R_2  &= \sup\{r \in (r_{M-2}(\alpha_0,t_0),1]: h(s,t_0)  - \pi > \chi_{\alpha^*}(s) \quad \forall s \in [r_{M-2}(\alpha_0,t_0),r)\}.
   \end{align*}
   We will prove that $r_{M-2}(\alpha^*,t) \to 0$ as $t \to T^-$. This is because there must be a curve $h(r(t),t) = \chi_{\alpha^*}(r(t))$ with $r_M(\alpha^*,t) \geq r(t) \geq r_{M-2}(\alpha^*,t)$ satisfying $h(r(t)+\varepsilon,t) < \chi_{\alpha^*}(r(t) + \varepsilon)$ for $\varepsilon \leq \varepsilon(t) \ll 1$, as well as $r(t) \to 0$ as $t \to T^-$ since $h(r,T) \leq \pi/2$ on $(0,1]$. 
   
   Fix $t^*$ large enough for which $r_{M-2}(\alpha^*,t^*) < r(t^*) < R_1$. One can also redefine $r_M(\alpha^*,t^*) = r(t^*) + \varepsilon < R_1$, and then $r_{M-1}(\alpha^*,t^*) \in (r_{M-2}(\alpha^*,t^*), r_M(\alpha^*,t^*))$. Next, one constructs paths going from time $t^*$ to $t_0$ as in Lemma \ref{lemma:geometry of h_1 < h < h_2} so that one obtains a new chain $s_i(\alpha^*,t_0)$ at time $t_0$ satisfying $P(t_0,\alpha^*,M)$, while still having the other chain $r_i(\alpha^*,t_0)$ at our disposal. We will observe that $s_M(\alpha^*,t_0)$ cannot exist: $s_M \notin (R_1,R_2)$ by definition, $s_M$ cannot be larger than $R_2$ because of the inequality deduced in Corollary \ref{cor:comparison of h with itself and exact value of limsup} and $s_M$ cannot be smaller than $R_1$ as this would induce chains of length $M+4$ satisfying $P(t_0,\alpha_0,M+4)$ or $P(t_0,\alpha^*,M+4)$, which is a contradiction.

  All of this reduces the possible harmonic maps in Theorem \ref{thm: main theorem, nonexistence bubble trees} to two possibilities: harmonic maps $\omega_i$ having inclination coordinate $2\arctan( \alpha^k r^k)$ or $\pi - 2 \arctan( \alpha^k r^k)$. As one can now assume without loss of generality that $h(r,T) > \pi/2$ on $(0,1]$ and $h(1,t) > \pi/2$ for $t \in [t_0,T]$, these chains allow to find a barrier $h(r,t) \geq \chi_{\alpha_0}(r)$ when $t \geq t_0, 1 \geq r \geq r^{+}(t)$, which implies that the scale $(R_n)_{n \in \mathbb N}$ of any bubble must lie on the left of the quantity 
  $$
 r^+(t) := \sup\{r \in [0,1]: h(s,t) \leq \pi \quad \forall s \in [0,r]\}.
  $$
  For fixed $t \in [t_0,T)$, the Maximum Principle also shows that $r \mapsto h(r,t)$ must first go from $0$ to $\pi/2$ and then from $\pi/2$ to $\pi$ without going below $\pi/2$ again. This prevents $h(r,t)$ from developing two harmonic maps or a harmonic map $\pi - 2 \arctan( \alpha^k r^k)$ on the left of $r^+(t)$, which concludes the proof.

We now prove the bubble-tree result with all the necessary details.

    \begin{lemma}\label{lemma:existence of value larger than pi on parabolic boundary}
        For any $(r_0,t_0) \in (0,1) \times (0,T)$, there is $(r_1,t_1) \in \{r_0\} \times [t_0,T] \cup (0,r_0] \times \{t_0\}$ for which
        $$
        |h(r_1,t_1)| > \pi.
        $$
    \end{lemma}

    \begin{proof}
   Fix any $(r_0,t_0) \in (0,1) \times (0,T)$. After parabolic rescaling, one can assume that $\{0,r_0\} \times [t_0,T] \cup [0,r_0] \times \{t_0\}$ is the new parabolic boundary. The new boundary data is given by any extension (e.g. given by Whitney's extension theorem combined with a smooth cut-off) of 
    $$
    \tilde{h}_0(r,t) = \begin{cases}
        h(r,t) &\text{ if } (t,r) \in   \{r_0\} \times  [t_0,T],  \\
        h(r,t) &\text{ if } (t,r) \in [0,r_0] \times  \{t_0\},  \\
        h_0(r,t) = 0 &\text{ if } (t,r) \in   \{0\} \times [t_0,T],
    \end{cases}
    $$
    which is smooth on the connected set where it is defined (Proposition \ref{prop:smoothness at the boundary, earlier statement}).

   For any such $\tilde{h}_0$, as $\tilde{h}_0(0,t) = 0$, there must be a point with $|\tilde{h}_0(r_1,t_1)| = |h(r_1,t_1)| > \pi$ by Proposition \ref{prop:criterion for global solution}, $(r_1,t_1) \in \{r_0\} \times [t_0,T] \cup (0,r_0] \times \{t_0\}$. Otherwise, $h$ would be global.
   \end{proof}

As a Corollary, replacing $h$ by $-h$ and $v = (v_1,v_2,v_3)$ by $(-v_1,-v_2,v_3)$ if needed, one can assume that there is $(r_n,t_n) \to (0,T)$, $r_n \in (0,1)$, $0 < t_n \leq T$ for which 
\begin{equation}
    h(r_n,t_n) > \pi  \quad \forall n. \label{eq: sequence of h larger than pi}
\end{equation}
Slightly perturbating $t_n$, we may even assume $t_n < T$. In particular,
    \begin{equation}
        \limsup_{\substack{r \to 0^+ \\ t \to T^-}}h(r,t) = \pi \label{eq: limsup value of h}
    \end{equation}
    by Corollary \ref{cor:comparison of h with itself and exact value of limsup}. 

   \begin{lemma}\label{lemma:limit of h(r,T) at zero}
   For all $r_0 \in (0,1)$, there exists $r \in (0,r_0)$ with $h(r_0,T) > 0$. In particular, $\lim \limits_{r \to 0^+} h(r,T) \in \{0,\pi\}$.
   \end{lemma}

   \begin{proof}
   The limit at $0$ exists and is a multiple of $\pi$ by Proposition \ref{value of v and h at origin}.

   As     $$
\limsup_{\substack{r \to 0^+ \\ t \to T^-}}h(r,t) = \pi,
    $$
    the limit must be $\leq \pi$. If the first part of the lemma is proved, then it follows that the limit is either $0$ or $\pi$.

   If the lemma was wrong, there would be $r_0 \in (0,1)$ with $h(r,T) \leq 0$ for all $r \in (0,r_0]$.  By Corollary \ref{cor:comparison of h with itself and exact value of limsup}, there exists $\tau_0 = \tau_0(v,T)$ such that for all $0 < \tau \leq \tau_0$, one has
    $$
    h(r,T-\tau) - \pi \leq h(r,T), \quad r \in (0,1].
    $$    
    Let $(r_n,t_n) \to (0,T)$ with $h(r_n,t_n) > \pi$, $r_n \in (0,1)$, $t_n \in (0,T)$. Write $t_n = T - \tau_n$ and assume that $0 < \tau_n \leq \tau_0$ and $0 < r_n < r_0$ by taking $n$ large enough. Then
    $$
    \pi < h(r_n, T-\tau_n) \leq h(r_n,T) + \pi \leq \pi,
    $$
    which is a contradiction. 
   \end{proof}

   \begin{theorem}\label{thm:limit r to 0 at time T is pi}
     The limit of $h(\cdot,T)$ at zero is exactly  $\lim \limits_{r \to 0^+} h(r,T) = \pi$.
   \end{theorem}

   \begin{proof}
       Assume not. Choose $r_0$ small enough so that $|h(r,T)| < \pi/4$ on $(0,r_0]$, $h(r_0,T) > 0$ by Lemma \ref{lemma:limit of h(r,T) at zero} and $t_0$ large enough so that $0 < h(r_0,t) < \pi/4$ on $[t_0,T]$. After parabolic rescaling and translating the time, one can assume that $(r_0,t_0) = (1,0)$ for simplicity.

    Let $\alpha_0$ and $t_0$ be large enough so that the maximal number of intersections $M(t,h,h_1 = \chi_{\alpha},h_2 = \pi)$ from Lemma \ref{lemma:intersection argument with chains mod 4} is minimal when $t \geq t_0$, $\alpha \geq \alpha_0$, $\chi_{\alpha} = \pi - 2\arctan \left( (\alpha r)^k \right)$.

    By taking $\alpha_0$ even larger, one can assume that $\chi_{\alpha}(1) < h(1,t)$ for all $t \in [t_0,T]$, $\alpha \geq \alpha_0$.

    As there exists $(r_n,t_n)_{n \in \mathbb N} \to (0,T)$, $r_n \in (0,1)$, $t_n \in (0,T)$,  with $h(r_n,t_n) > \pi$, for any $t \in [t_0,T)$, one can fix $t_n$ large enough so that $t_n > t$ and $\varepsilon$ small enough so that $t_n + \varepsilon < T$. Applying Lemma \ref{geometry of h_1 < h_2}, one must have a point $h(r(t),t) > \pi$ on the parabolic boundary of $[0,1] \times [t,t_n + \varepsilon]$. 

    Then $\lim \limits_{t \to T^-} r(t) = 0$. Indeed, if $(t_n) \to T^{-}$ is any sequence and $(t_{n_m})$ is any subsequence, then there is a further subsequence (still denoted with index $n_m$) for which $r(t_{n_m}) \to r^* \in [0,1]$. If $r^* > 0$, then $h(r(t_{n_m}),t_{n_m}) \to h(r^*,T) \in [-\pi/4,\pi/4]$ by continuity on $(0,1] \times [0,T]$, which is a contradiction. Hence, $r^* = 0$.

    In particular, $M(t_0,h,h_1 = \chi_{\alpha_0},h_2 = \pi) \geq 5$. Assume not and let  $t \geq t_0, \alpha \geq \alpha_0$. As $h(r_1,t) = 0 < \chi_{\alpha}(r_1)$ near $r_1 = 0$, $\chi_{\alpha}(r_2) < h(r_2,t) < \pi$ at a point $r_2 < r^-(t)$ close enough to 
    $$
    r^-(t) = \sup\{r \in [0,1]: h(s,t) < \pi \quad \forall s \in [0,r]\} \in (0, r(t)),
    $$
    as well as $h(r_3,t) > \pi$ at $r_3 = r(t)$ and $\chi_{\alpha}(r_4) < h(r_4,t) < \pi$ at a point $r_4 > r^+(t)$ close enough to
        $$
    r^+(t) = \sup\{r \in [r(t),1]: h(s,t) > \pi \quad \forall s \in [r(t),r]\} \in (r(t),1),
    $$
    then one would have $h(r,t) \geq \chi_{\alpha}(r)$ for all $r \geq r(t) > 0$, $t \in [t_0,T)$. If $r \in (0,1]$ is fixed, then there is $t_1(r) \in [t_0,T)$ large enough for which $r \geq r(t) > 0$ for all $[t_1,T)$. Hence, $h(r,t) \geq \chi_{\alpha}(r)$ on $[t_1(r),T)$. Taking the limit $t \to T^{-}$, we find out that $h(r,T) \geq \chi_{\alpha}(r)$ for all $r \in (0,1]$, which contradicts $\lim \limits_{r \to 0} h(r,T) = 0$.

    Hence, $M(t_0,h,h_1 = \chi_{\alpha_0},h_2 = \pi)  \geq 5$. For $\alpha \geq \alpha_0$, $t \in [t_0,T)$, fix chains $r_1(\alpha,t) < ... < r_M(\alpha,t)$ satisfying $P(t,\alpha,M)$ with $M$ minimal. Consider $r_{M-2}(\alpha, t)$ for which $h(r_{M-2}(\alpha,t),t) > \pi$. Arguing as with $r(t)$ above,
    $$\lim \limits_{t \to T^{-}} r_{M-2}(\alpha,t) = 0.$$
    Let 
    $$
    z(\alpha,t) = \sup\{r \in (r_{M-2}(\alpha,t),1]: h(s,t) \geq \chi_{\alpha}(s) \quad \forall s \in [r_{M-2}(\alpha,t),r]\}.
    $$
    Note that $z(\alpha,t) < 1$ since
    $$
    r_{M-2}(\alpha,t) < z(\alpha,t) < r_{M}(\alpha,t) \leq 1.
    $$
    In particular, $h(r,t) < \chi_{\alpha}(r)$ for $(r,t)$ with all $r > z(\alpha,t)$ close enough to $z(\alpha,t)$.

    Then one must have $\lim \limits_{t \to T^-} z(\alpha,t) = 0$ for all $\alpha \geq \alpha_0$. Otherwise, there would be $t_n \to T^{-}$, $t_n > t_0$ and $z^* \in (0,1]$ for which $z(\alpha,t_n) \geq z^*$ for all $n$ so that
    $$
    h(s,t_n) \geq \chi_{\alpha}(s) \quad \forall s \in [r_{M-2}(\alpha,t_n),z^*] \quad \forall n \in \mathbb N.
    $$
    As $\lim \limits_{n \to +\infty} r_{M-2}(\alpha,t_n) = 0$, for $s \in (0,z^*]$ fixed, one deduces $s \in [r_{M-2}(\alpha,t_n), z^*]$ for all $n \geq n_0$ large enough so that
    $$h(s,t_n) \geq \chi_{\alpha}(s) \quad \forall n \geq n_0,$$
    hence $h(s,T) \geq \chi_{\alpha}(s)$ for all $s \in (0,z^*]$ by passing at the limit, which contradicts $\lim \limits_{r \to 0} h(r,T) = 0$ as before.

     It follows from Corollary \ref{cor:comparison of h with itself and exact value of limsup} that there exists $\tau_0 = \tau_0(v,T)$ such that for all $0 < \tau \leq \tau_0$, $t \in [0,T-\tau]$, one has
    \begin{equation}
        h(r,t) - \pi \leq h(r,t+\tau), \quad r \in (0,1]. \label{eq: contradiction at time t^*}
    \end{equation}
    Fix $t_0$ large enough so that $0 < T-t_0 < \tau_0$. Fix $\alpha^* > \alpha_0$ for which $0 < \chi_{\alpha^*}(r_{M-2}(\alpha_0,t_0)) < h(r_{M-2}(\alpha_0,t_0),t_0) - \pi$.

  Let 
    $$
  R_1 = \inf\{r \in (0,r_{M-2}(\alpha_0,t_0)]: h(s,t_0) - \pi > \chi_{\alpha^*}(s) \quad \forall s \in (r,r_{M-2}(\alpha_0,t_0)]\}
  $$
  and
  $$
  R_2  = \sup\{r \in (r_{M-2}(\alpha_0,t_0),1]: h(s,t_0)  - \pi > \chi_{\alpha^*}(s) \quad \forall s \in [r_{M-2}(\alpha_0,t_0),r)\}.
  $$
  One has $R_1 > 0$ so there is $t^* \in (t_0,T)$ large enough for which $0 < z(\alpha^*,t^*)< R_1$. In particular, $z(\alpha^*,t^*) < z(\alpha_0,t_0)$. By definition of $z(\alpha^*,t^*)$, there is $\varepsilon > 0$ for which  $0 < z(\alpha^*,t^*)+ \varepsilon < R_1$ and
  $$
  h(z(\alpha^*,t^*) + \varepsilon,t^*) < \chi_{\alpha^*}(z(\alpha^*,t^*) + \varepsilon).
  $$
  As 
  $$r_M(\alpha^*,t^*) > z(\alpha^*,t^*)  > r_{M-2}(\alpha^*,t^*),$$
  one can redefine $r_M(\alpha^*,t^*)  = z(\alpha^*,t^*) + \varepsilon < R_1$ and find another point $r_{M-1}(\alpha^*,t^*)$ inbetween $(r_{M-2}(\alpha^*,t^*), z(\alpha^*,t^*))$ without altering the length of the chain. As $h(0,t) < \chi_{\alpha^*}(0)$ on $[t_0,t^*]$, there is $r_1 > 0$ for which $h(r,t) < \chi_{\alpha^*}(r)$ on $[t_0,t^*] \times [0,r_1]$. Then we can replace $r_1(\alpha^*,t)$ by the same value $r_1$ for all $t \in [t_0,t^*]$.
  
  Arguing as in the proof of Theorem \ref{lemma:intersection argument with chains mod 4}, there are injective paths  $\Sigma_1, ..., \Sigma_M$ for which $\Sigma_i$ connects $(r_i(\alpha^*,t^*), t^*)$ to $(s_i,t_0)$ and $s_1 < ... < s_M$ satisfies $P(t_0,\alpha^*,M)$ as well.

  Then one must have either $s_M \leq R_1$ or $s_M \geq R_2$. In the former case, $s_M \notin [r_{M-2}(\alpha^*,t_0),z(\alpha^*,t_0)]$ as $h \geq \chi_{\alpha^*}$ on that interval at time $t_0$. One cannot have $z(\alpha^*,t_0) < s_M \leq R_1 \leq r_{M-2}(\alpha_0,t_0)$ either as $h < \chi_{\alpha^*} \leq \chi_{\alpha_0}$ at $(s_M,t_0)$, $h > \pi$ at $s_{M+2} = R_1$, $\chi_{\alpha_0} < h < \pi$ at some $s_{M+1}$ inbetween $(s_M, s_{M+1})$ and at $s_{M+3} = r_{M-1}(\alpha_0,t_0)$, as well as $h < \chi_{\alpha_0}$ at $s_{M+4} = r_M(\alpha_0,t_0)$, leading to a chain satisfying $P(t_0,\alpha_0,M+4)$, contradicting the maximality of $M$. 
  
  Hence, $s_M < r_{M-2}(\alpha^*,t_0) < z(\alpha^*,t_0)$. Then $h < \chi_{\alpha^*}$ at $(s_M,t_0)$, $h > \pi$ at $(s_{M+2},t_0) := (r_{M-2}(\alpha^*,t_0),t_0)$, $\chi_{\alpha*} < h < \pi$ for a point $s_{M+1}$ inbetween $(s_M,r_{M-2}(\alpha^*,t_0))$ and a point $s_{M+3} \in (r_{M-2}(\alpha^*,t_0),  z(\alpha^*,t_0))$, and $h < \chi_{\alpha^*}$ at $(z_{M+4},t_0) := (z(\alpha^*,t_0)+ \varepsilon',t_0)$, $\varepsilon' > 0$ small enough. Hence, we would have found a chain satisfying $P(t_0,\alpha^*,M+4)$, contradicting the maximality of $M$ once again.

  Thus, $s_M \geq R_2$. The path $\Sigma_M$ from $(r_M(\alpha^*,t^*),t^*)$ to $(s_M, t_0)$ must path through a point $(\tilde{r},\tilde{t}) \in (R_1,R_2) \times [t_0,t^*]$. 
 
 Finally, note that $t_0 \in [0,T-\tau]$ for any $0 < \tau \leq T - t_0 < \tau_0$. Then $\tilde{t} = t_0 + \tau$ for some $0 < \tau \leq T-t_0$. Evaluating (\ref{eq: contradiction at time t^*}) at $r = \tilde{r}$, $t = t_0$ and $t+\tau = \tilde{t}$ yields
 $$
 h(\tilde{r},t_0) - \pi \leq h(\tilde{r},\tilde{t}) < \chi_{\alpha^*}(\tilde{r}).
 $$
 As $\tilde{r} \in (R_1,R_2)$, this contradicts the definition of $R_1, R_2$.
    \end{proof}

\begin{corollary}\label{cor: liminf h(r,t) = 0}
The limit inferior at the singular point is exactly 
    $
\liminf \limits_{\substack{r \to 0^+ \\ t \to T^-}}h(r,t) = 0.
    $
\end{corollary}

\begin{proof}
    Assume not. As $h(0,t) = 0$, one must have
    $$
m:= \liminf_{\substack{r \to 0^+ \\ t \to T^-}}h(r,t) < 0,
    $$
    and there exists $(r_n,t_n) \to (0,T)$ with  $\lim \limits_{n \to +\infty} h(r_n,t_n) = m < 0$. Moreover, $r_n > 0$ (as $h(0,t) = 0$) and $t_n < T$ (as $\lim \limits_{r \to 0^+} h(r,T) = \pi$) for all $n$ large enough.
    
    By Corollary \ref{cor:comparison of h with itself and exact value of limsup}, there exists $\tau_0 = \tau_0(v,T)$ such that for all $0 < \tau \leq \tau_0$, one has
    $$
    h(r,T) \leq h(r,T-\tau) + \pi, \quad r \in (0,1].
    $$    
    Write $t_n = T - \tau_n$ and assume that $0 < \tau_n \leq \tau_0$ and $0 < r_n < r_0$ by taking $n$ large enough. Then
    $$
     h(r_n, T) \leq h(r_n,T-\tau_n) + \pi.
    $$
    This is a contradiction since the left-hand side converges to $\pi$, while the right-hand side converges to $m + \pi < \pi$.
\end{proof}

    Next, choose $r_0$ small enough and $t_0$ large enough so that 
    \begin{equation}
            5\pi/4 \geq h(r,t) \geq -\pi/2, \quad (r,t) \in (0,r_0] \times [t_0,T). \label{eq: bounds on h}
    \end{equation}
   Choose $r_0$ smaller if needed so that 
    \begin{equation}
            h(r,T) > \pi/2 \quad r \in (0,r_0] \label{eq: lower bound of h at t = T}
    \end{equation}
    as well and $t_0$ larger so that 
\begin{equation}
        h(r_0,t) > \pi/2 \quad t \in [t_0,T]. \label{eq: lower bound of h at r = 1}
\end{equation}
    
    After parabolic rescaling and translating the time, one can assume that $(r_0,t_0) = (1,0)$ for simplicity. 
    
    Let $t_0 > 0$ and $\alpha_0 > 0$ be large enough so that $\chi_{\alpha_0}(1) < \pi/2$ and the length $\tilde{M}$ of a chain satisfying the property $P(t,h,h_1 = \chi_{\alpha},h_2 = \pi,\tilde{M})$ from Lemma \ref{lemma:intersection argument with chains mod 4} is minimal when $t \geq t_0$, $\alpha \geq \alpha_0$.
    
    Let 
    \begin{align*}
        r^+(t) &= \sup\{r \in [0,1]: h(s,t) \leq \pi \quad \forall s \in [0,r]\}, \\
        r^-(t) &= \sup\{r \in [0,1]: h(s,t) \leq \pi/2 \quad \forall s \in [0,r]\} < 1.
    \end{align*}

\begin{lemma} \label{lemma: lower bound for h on the right of r^+}
One has
         $$
         h(r,t) \geq \chi_{\alpha_0}(1), \quad \forall t \geq t_0, \ r \geq r^+(t).
         $$
\end{lemma}

\begin{proof}
     Assume there is $t^* \geq t_0$ and $r^* \geq r^+(t^*)$ with 
     $$
h(r,t)< \chi_{\alpha_0}(1) = \inf_{r \in [0,1]} \chi_{\alpha_0}(r) < \pi/2.
     $$
    In particular, $r^+(t^*) \leq r^* < 1$ by (\ref{eq: lower bound of h at r = 1}), where the problem has been rescaled so that $r_0 = 1$.
    
    Then at time $t^*$, $h < \chi_{\alpha_0}$ near $r_1 = 0$, $h > \pi > \chi_{\alpha_0}$ near $r_3 \in (r^+(t^*),r^+(t^*) + \delta)$ with $\delta \ll 1$, $r^+(t^*) + \delta < r^*$, $\chi_{\alpha_0} < h < \pi$ for a point $r_2 \in (r_1,r_2)$, $h < \chi_{\alpha_0}$ for a point $r_5 = r^*$ and $\chi_{\alpha_0} < h < \pi$ for some $r_4 \in (r_3,r_5)$. Hence, $\tilde{M}(t,h,h_1 = \chi_{\alpha},h_2 = \pi) \geq 5$ for $t \geq t_0$ and $\alpha \geq \alpha_0$.

    For $t \geq t_0,\alpha \geq \alpha_0$, consider chains $r_1(t,\alpha) < ... < r_{\tilde{M}}(t,\alpha)$ of length $\tilde{M} \geq 5$ satisfying $P(t,h,h_1 = \chi_{\alpha},h_2 = \pi,\tilde{M})$ from Lemma \ref{lemma:intersection argument with chains mod 4}. In particular, $r^+(t) < r_3(t,\alpha) \leq 1$ for all $t\geq t_0$.

    Let also
    $$
    r^+(t) < z(t,\alpha) = \sup\{r \in [r^+(t),1]: h(s,t) \geq \chi_{\alpha}(s) \ \forall s \in [r^+(t),r]\} < r_{5}(t,\alpha) \leq 1.
    $$
    We show that $z(t,\alpha) > r_3(t,\alpha)$. Assume not. Then, there is $\rho_3 \in (r^+(t),r_3(t,\alpha))$ with $h(\rho_3,t) < \chi_{\alpha}$. There also exists $\rho_4 \in (\rho_3, r_3(t,\alpha))$ with $\chi_{\alpha}(\rho_4) < h(\rho_4,t) < \pi$ since $h(r_3(t,\alpha),t) > \pi$. Moreover, there exists $\rho_1 \in (r^+(t),\rho_3)$ close to $r^+(t)$ with $h(\rho_1,t) > \pi$ and $\rho_2 \in (\rho_1,\rho_3)$ with $\chi_{\alpha} < h  < \pi$. Hence, $r_1 < r_2 < \rho_1 < \rho_2 < \rho_3 < \rho_4 < r_3 < ... < r_{\tilde{M}}$ is chain of length $\tilde{M} + 4$, which contradicts the maximality of the chain.
    
    Then $\lim \limits_{t \to T^{-}}z(t,\alpha) = 0$. Otherwise, there is $(t_n) \to T^{-}$ with $z(t_n,\alpha) \to z^*(\alpha) > 0$. As $z^*(t_n,\alpha) < 1$, one must have $h(z(t_n,\alpha),t_n) = \chi_{\alpha}(t_n)$ for all $n$. Then 
    $$
    \chi_{\alpha}(z(t_n,\alpha)) = h(z(t_n,\alpha),t_n) \to h(z^*(\alpha),T) = \chi_{\alpha}(z^*(\alpha))
    $$
    with $h(z^*(\alpha),T) \geq \pi/2$. Hence, $z^*(\alpha) \leq \alpha^{-1}$ meaning that $\lim \limits_{\alpha \to +\infty} z^*(\alpha) = 0$. Yet, observe that $z(t_n,\alpha)$ and $z^*(\alpha)$ are increasing with respect to $\alpha$, which would imply $z^*(\alpha) = 0$ is constant with respect to $\alpha$, a contradiction.

It follows from Corollary \ref{cor:comparison of h with itself and exact value of limsup} that there is $\tau_0 > 0$ such that for all $0 < \tau \leq \tau_0$, one has
    $$
    h(r,T-\tau) - \pi \leq h(r,T), \quad r \in (0,1].
    $$
     Choose $t_0$ initially larger so that $0 < T-t_0 < \tau_0$. Fix $\alpha^* > \alpha_0$ for which $0 < \chi_{\alpha^*}(r_{3}(\alpha_0,t_0)) < h(r_{3}(\alpha_0,t_0),t_0) - \pi$.

  Letting
    \begin{align*}
          R_1 &= \inf\{r \in (0,r_{3}(\alpha_0,t_0)]: h(s,t_0) - \pi > \chi_{\alpha^*}(s) \quad \forall s \in (r,r_{3}(\alpha_0,t_0)]\}, \\
  R_2  &= \sup\{r \in (r_{3}(\alpha_0,t_0),1]: h(s,t_0)  - \pi > \chi_{\alpha^*}(s) \quad \forall s \in [r_{3}(\alpha_0,t_0),r)\},
    \end{align*}
    and arguing as in the end of the proof of Theorem \ref{thm:limit r to 0 at time T is pi}, one deduces a contradiction.
\end{proof}

\begin{corollary}\label{cor: bubbles lie on the left of r^+}
   Assume that there exists $R_n \to 0^+, T_n \to T^-$ for which $h(R_nr,T_n)$ converges to 
     $$
    h(R_nr,T_n) \to h_{\infty}(r) = m_i \pi + \varepsilon_i \cdot 2 \arctan \left( (\alpha_i r)^k \right), \quad m_i \in \mathbb Z, \varepsilon_i \in \{-1,1\}, \alpha_i > 0,
    $$
    uniformly on any compact subsets of $(0,+\infty)$. Then there exists $R_0 \in (0,+\infty)$ and $n_0 \in \mathbb N$ for which $R_0R_n \leq r^+(T_n)$ for all $n \geq n_0$.
\end{corollary}

\begin{proof}
    As 
    $$
\liminf_{\substack{r \to 0^+ \\ t \to T^-}}h(r,t) = 0, \quad \limsup_{\substack{r \to 0^+ \\ t \to T^-}}h(r,t) = \pi,
    $$
    the only possibilities for the limit $h_{\infty}$ are
    $$
 \pi -  2 \arctan\left( (\alpha r)^k \right), \quad 2 \arctan \left( (\alpha r)^k \right)
    $$
    In both cases, as $h_{\infty}((0,+\infty)) = (0,\pi)$, there must be $R_0 \in (0,+\infty)$ for which
    $$
    h(R_nR_0,T_n) \to h_{\infty}(R_0) < \chi_{\alpha_0}(1),
    $$
    where $\alpha_0$ is as in Lemma \ref{lemma: lower bound for h on the right of r^+}. For all $n$ large enough, one must have $T_n \geq t_0$ and $h(R_nR_0,T_n) < \chi_{\alpha_0}(1)$, hence $R_0R_n \leq r^+(T_n)$ by Lemma \ref{lemma: lower bound for h on the right of r^+}.
\end{proof}


     
\begin{lemma}\label{lemma:continuity of r^+}
The maps $t \mapsto r^-(t)$, $t \mapsto r^+(t)$ are continuous on $[t_0,T)$.
\end{lemma}

\begin{proof}
The proof is exactly the same for both functions, so we only do the one for $r^{+}(t)$.

    Fix $t_1 \in [t_0,T)$ and we show continuity at $t_1$. Given $\varepsilon > 0$, there must be $r_1 \in [r^+(t_1)-\varepsilon, r^+(t_1)]$ which $h(r_1,t_1) < \pi$ as otherwise, $h(r,t_1) = \pi$ would be constant on $(0,1]$ since $r \mapsto h(r,t_1)$ is real-analytic on $(0,1]$ (Theorem \ref{thm:space analyticity}), which contradicts $h(0,t_1) = 0$. By continuity of $h$, there is $\delta > 0$ small enough for which $h(r_1,t) < \pi$ on $[t_1, t_1 + \delta]$. The Comparison Principle (Theorem \ref{comparison principle}) implies $h(r,t) \leq \pi$ on $[0,r_1] \times [t_1,t_1+\delta]$. As in the proof of Proposition \ref{prop:criterion for global solution}, applying the Maximum Principle (Theorem \ref{maximum principle}), one finds that $h < \pi$ on $(0,r_1) \times (t_1,t_1+\delta)$. Hence, $r^+(t) \geq r_1 \geq r^+(t_1) - \varepsilon$ for all $t \in [t_1,t_1+\delta]$.

If $r^+(t_1) = 1$, we are done as $r^+(t) \leq 1 = r^+(t_1)$ for all $t \in [t_1,t_1+\delta]$. Otherwise, $r^+(t_1) < 1$ and there is $r_2 \in (r^+(t_1), r^+(r_1) + \varepsilon]$ with $h(r_2,t_1) > \pi$. Then there is $\delta > 0$ small enough for which $h(r_2,t) > \pi$ on $[t_1, t_1 + \delta]$ meaning that $r^+(t) \leq r_2 \leq r^+(t_1) + \varepsilon$ on $[t_1, t_1 + \delta]$. This concludes the proof of continuity.
\end{proof}


\begin{lemma}\label{cor:only one bubble in the left of r^+}
    For $t \in [t_{0},T)$, $h < \pi/2$ on $[0,r^{-}(t))$ and $\pi/2 < h < \pi$ on $(r^{-}(t),r^+(t))$. 
\end{lemma}

\begin{proof}
We start by proving that $h(r,t) < \pi/2$ for all $r \in (0, r^{-}(t))$, $t \in [t_0,T)$. Indeed, fix $t_1 \in [t_0,T)$, $0 < r_1 < \min\{r^{-}(t) : t \in [t_0,t_1]\}$ and denote by $E_{r_1,t_1}$ the open, connected, bounded domain enclosed by 
$$
\{r = r_1, t_0 < t < t_1\} \cup \{t \in \{t_0,t_1\}, r_1 \leq r \leq r^{-}(t)\} \cup \{(t,r^{-}(t)) : t_0 < t < t_1\}.
$$
Then $h \leq \pi/2$ on this domain by definition. As in the proof of Proposition \ref{prop:criterion for global solution}, applying the Maximum Principle (Theorem \ref{maximum principle}) to $(h-\pi/2)$ shows that $h < \pi/2$ on $E_{r_1,t_1}$. As $t_1$ is arbitrary and $r_1$ is arbitrarily small, $h(r,t) < \pi/2$ on $(0,r^{-}(t))$ for any $t \in [t_0,T)$.

Similarly, one proves that $h < \pi$ on $(0,r^{+}(t))$ by looking at $h-\pi$ and the domain $E_{r_1,t_1}$ enclosed by 
$$
\{r = r_1, t_0 < t < t_1\} \cup \{t \in \{t_0,t_1\}, r_1 \leq r \leq r^{+}(t)\} \cup \{(t,r^{+}(t)) : t_0 < t < t_1\}.
$$
Finally, one proves that $h > \pi/2$ on $(r^{-}(t),r^{+}(t))$ by looking at $\pi/2-h$ and the domain $E_{t_1}$ enclosed by 
$$
 \{(t,r^{-}(t)) : t_0 < t < t_1\} \cup \{t \in \{t_0,t_1\}, r^{-}(t) \leq r \leq r^{+}(t)\} \cup \{(t,r^{+}(t)) : t_0 < t < t_1\}.
$$
which lies at a positive distance from the line $\{r = 0\}$ (since $h(0,t) = 0$).  One still needs to prove that $h(r,t_0) \geq \pi/2$, $r^{-}(t_0) \leq r \leq r^+(t_0)$, at least at one time $t_0$ close to $T$. This is another chain argument similar to the previous ones (see \cite{VANDERHOUT2003}, Lemma 2.12).
\end{proof}

\begin{corollary}\label{cor:end of proof}
Theorem \ref{thm: main theorem, nonexistence bubble trees} holds and $\omega_1$ has inclination coordinate of the form:
$$
\pm 2\arctan \left( (\alpha r)^k \right) + 2\pi \mathbb Z, \quad \alpha > 0.
$$
\end{corollary}

\begin{remark}
   The sign comes from the fact that, without loss of generality, one replaces $h$, $(v_1,v_2,v_3)$ by $-h$, $(-v_1,-v_2,v_3)$ if needed to get (\ref{eq: sequence of h larger than pi}). Moreover, inclination coordinates are unique up to a multiple of $2\pi$. Up until now, we have always considered inclination coordinates chosen so that $h(0,t) = 0$.
\end{remark}
\begin{proof}
     Suppose for a contradiction that there are at least two harmonic maps (i.e. $p \geq 2$ in the statement of Theorem \ref{thm: main theorem, nonexistence bubble trees}) and order their scales so that $R_n^{(1)} \ll R_n^{(2)} \to 0^+$. Thanks to the Sobolev embedding $H^{1,2}_{loc}((0,+\infty)) \hookrightarrow C^{0,1/4}_{loc}((0,+\infty))$, $h(R_n^{(i)}r,T_n)$ converges uniformly on any compact set to a limit $H_{\infty}^{(i)} \in C^{0}((0,+\infty))$. One must have
\begin{equation*}
    \omega_i(re^{i\theta}) = \left(e^{ik \theta} \sin H_{\infty}^{(i)}(r), \cos H_{\infty}^{(i)}(r) \right)  
\end{equation*}
 by pointwise convergence $v(R_n^{(i)}x,T_n) \to \omega_i(x)$ given by the statement of Theorem \ref{thm: main theorem, nonexistence bubble trees}. But Lemma \ref{k-equivariance, smoothness of h} shows that 
\begin{equation*}
    \omega_i(re^{i\theta}) = \left(e^{ik \theta} \sin H^{(i)}(r), \cos H^{(i)}(r) \right) 
\end{equation*}
for some lifting $H^{(i)} \in C^{\infty}([0,+\infty))$. Hence, $H_{\infty}^{(i)}(r)$ and $H^{(i)}(r)$ must differ by a fixed multiple of $2\pi$ and $H_{\infty}$ is smooth as well. Repeating the ODE argument that lies beneath equation (\ref{eq: H vs H_infty, 1}) in the proof of Theorem \ref{thm:smoothness around all but finitely many points} shows that $H_{\infty}^{(i)}(r)$ is of the form (\ref{inclination coordinate of harmonic map}). In other words, there are some parameters for which
    $$
    h(R_n^{(i)}r,T_n) \to H_{\infty}^{(i)}(r) = m_i \pi + \varepsilon_i \cdot 2 \arctan \left( (\alpha_i r)^k \right), \quad m_i \in \mathbb Z, \varepsilon_i \in \{-1,1\}, \alpha_i > 0,
    $$
    uniformly on compact subsets of $(0,+\infty)$. Replacing $R_n^{(i)}$ by $R_n^{(i)}\alpha_i^{-1}$, one can assume that $\alpha_i = 1$.
    
    As 
    $$
\liminf_{\substack{r \to 0^+ \\ t \to T^-}}h(r,t) = 0, \quad \limsup_{\substack{r \to 0^+ \\ t \to T^-}}h(r,t) = \pi,
    $$
    the only possibilities for the limits $H_{\infty}^{(i)}$ are
    $$
 \pi -  2 \arctan (r^k), \quad 2 \arctan (r^k),
    $$
    where it might a priori be the case that we get the same stationary solution for both scales $R_n^{(1)}$ and $R_n^{(2)}$. 



    
   Fix $\pi/2 > \varepsilon > 0$ small enough so that $\varepsilon < \sup_{r \in [0,1]} \chi_{\alpha_0}(r)$, where $\alpha_0$ is as in Lemma \ref{lemma: lower bound for h on the right of r^+}. In case of a limiting map $2\arctan(r^k)$, for $R_n^{(i)}$, one chooses $0 < r_0 < r_1 < +\infty$ so that $2\arctan(r_0^k) < \varepsilon$ and $2\arctan(r_1^k) > \pi-\varepsilon$ and $N(\varepsilon)$ large enough so that $h(r_0R_n^{(i)}, T_n) < \varepsilon$ and $h(r_1R_n^{(i)}, T_n) > \pi - \varepsilon$ for all $n \geq N(\varepsilon)$. We proceed similarly for a limiting map $\pi - 2 \arctan(r^k)$. Hence, for all $n \geq N(\varepsilon)$, $h(r,T_n)$ achieves all the values inbetween $(\varepsilon, \pi-\varepsilon)$ twice, once inside an interval $[r_0R_n^{(1)},r_1R_n^{(1)}]$ and once inside the interval $[r_2R_n^{(2)},r_3R_n^{(2)}]$. As $r_i = r_i(\varepsilon)$ does not depend on $n$, for $n$ large enough, $r_1R_n^{(1)} < r_2 R_n^{(2)}$, i.e.,  the intervals $[r_0R_n^{(1)},r_1R_n^{(1)}]$, $[r_2R_n^{(2)},r_3R_n^{(2)}]$ are disjoint. By Lemma \ref{cor:only one bubble in the left of r^+}, they cannot both lie on the left of $r^{+}(T_n)$, as otherwise, they would both contain $r^{-}(T_n)$, the only possible point for which $h(r,T_n) = \pi/2$ on $(0,r^+(T_n))$. 
   
   Hence, there is at least one interval lying on the right of $r^+(T_n)$. By Lemma \ref{lemma: lower bound for h on the right of r^+}, $h(r,T_n) \geq \chi_{\alpha_0}(1) > \varepsilon$ on $(r^{+}(T_n),1]$, i.e., $h$ cannot achieve all the values $(\varepsilon, \chi_{\alpha_0}(1)) \subset (\varepsilon, \pi-\varepsilon)$ on that side. This is a contradiction.

    Finally, we cannot have a single harmonic map $\omega_1$ with inclination coordinate $\pi - 2 \arctan( r^k)$. Indeed, if $R_n^{(1)}$ and $T_n^{(1)}$ are the corresponding sequences, then 
 we would have $[r_0R_n^{(1)},r_1R_n^{(1)}] \subset [0,r^+(T_n)]$ at $n = N(\varepsilon)$, i.e., the interval must lie on the left of $r^+$ (otherwise it would contradict Lemma \ref{lemma: lower bound for h on the right of r^+} as above), as well as $h(0, T_n) = 0$, $h(r_0R_n^{(1)},T_n) > \pi - \varepsilon > \pi/2$, $h(r_1R_n^{(1)},T_n) = \varepsilon < \pi/2$, contradicting Lemma \ref{cor:only one bubble in the left of r^+}.

\end{proof}




    

\appendix
\chapter{Energy concentration and bubbles for the Harmonic Map Heat Flow} \label{section:struwe energy hmhf}

The content of this chapter is taken from the appendices of my article \cite{samuelian2025blowuptrees}, to appear soon in \textit{Calculus of Variations and Partial Differential Equations}. For the sake of completeness, proofs which were omitted in the original paper have been added.

\section{Introduction}
In this chapter, we follow the argument of Struwe (\cite[Chapter 3, Theorem 5.6]{struwe2008variational}), which holds for a wider class of evolution of harmonic maps, and proves that a finite-energy, smooth solution $v \in C^{\infty}(B^2 \times (0,T))$ of 
\begin{align*}
    v_t = \Delta v + |\nabla v|^2v, \quad (x,t) \in B^2  \times (0,T),
\end{align*}
blows-up in finite-time (in the sense that one cannot extend $v$ to a $C^{\infty}(B^2 \times (0,T])$ function) if and only if there is a concentration of energy at a finite number of points from $B^2$ (Theorem \ref{thm:smoothness around all but finitely many points}). In the $k$-equivariant setting, such a concentration can only happen at the origin (and $v$ is actually smooth at the boundary at the blow-up time $T$, Proposition \ref{prop:smoothness at the boundary}). Moreover, along some appropriate sequences converging to a point of concentration $(x,T)$, one can extract a non-constant harmonic map (called 'bubble') from $v$ and even grab all the bubbles at the same time (Theorem \ref{thm: full bubble decomposition}), yielding the so-called bubble-tree decomposition.

We note that Struwe's work concerns only manifolds without boundary, which is not the case here. As his arguments are essentially local, we show that one can still deduce a similar result in the interior of the ball, without assuming any kind of regularity at the boundary (the analog result with a boundary condition which does not depend on time has been proved in \cite{ChangStruweAnalog}) and only assuming finiteness of the energy. The main difference in our setting is that the energy $E(v(\cdot,t))$ is not necessarily decreasing in time, so Struwe's argument needs to be adapted slightly. 

We will still deduce that under $k$-equivariance symmetry and smooth boundary data, no singularity can occur at the boundary and the solution is actually smooth at the boundary at the time of blow-up (Proposition \ref{prop:smoothness at the boundary}). 

For $U \subset  B^2$, let
    $$
 E(v(t);U) = \frac{1}{2}\int_{U} |\nabla v|^2 dx
    $$
be the local energy of $v$ on $U$.

\section{The energy estimates of Struwe} 
In the following Lemmas, unless specified otherwise, the solution $v(x,t)$ need not to be equivariant or defined at the boundary $\partial B_R(x_0)$ or at time $t \in \{0,T\}$. 

\begin{lemma}[\cite{struwe2008variational}, Chapter 3, Lemma 5.7] \label{Sobolev, estimate L4 norm in terms of H1 norm}
    Let $u \in H^1(B_R(x_0),\mathbb R)$, $B_R(x_0) \subset \mathbb R^2$, and $\phi(x) \in C^{\infty}_c(B_{R}(x_0))$, $0 \leq \phi \leq 1$, be any cut-off function with $|\nabla \phi|^2 \leq C/R^2$. Then
    $$
    \int_{\mathbb R^2} |u|^4 \phi^2 dx \leq 8 \left( \int_{B_{R}(x_0)} |u|^2 dx \right) \left( \int_{B_{R}(x_0)} |\nabla u|^2 \phi^2 dx + C^2R^{-2} \int_{B_{R}(x_0)} |u|^2 dx\right).
    $$
\end{lemma}
\begingroup
\allowdisplaybreaks
\begin{lemma}\label{Sobolev comparison, int (Delta u)^2 > int (D^2_x u)^2}
     Let $u \in H^2(B_R(x_0),\mathbb R)$, $B_R(x_0) \subset \mathbb R^2$, and $\phi(x) \in C^{\infty}_c(B_{R}(x_0))$, $0 \leq \phi \leq 1$, be any cut-off function with $|\nabla \phi|^2 \leq C/R^2$. Then
    $$
\int_{B_R(x_0) } \phi^2  (\Delta u)^2 dx \geq \frac{1}{2} \int_{B_R(x_0) } \phi^2  |D^2_x u|^2 dx - 12C^2R^{-2} \int_{B_R(x_0) } |\nabla u|^2 dx.
    $$
\end{lemma}

\begin{proof}
    Without loss of generality, assume that $u$ is smooth. Then
        \begin{align*}
        \int_{B_R(x_0) } \phi^2  (\Delta u)^2 dx &=  \int_{B_R(x_0) } \phi^2  \sum_{i,j} \partial_{x_i}^2 u  \partial_{x_j}^2 u dx \\
        &= -\int_{B_R(x_0) }   \sum_{i,j} \phi^2 \partial_{x_i} u  \partial_{x_i} \partial_{x_j}^2 u + 2\phi \partial_{x_i} \phi  \partial_{x_i} u \partial_{x_j}^2 u dx \\
            &= \int_{B_R(x_0) }   \sum_{i,j} \phi^2 \left( \partial_{x_i} \partial_{x_j} u  \right)^2 + 2 \phi \partial_{x_j} \phi \partial_{x_i} u  \partial_{x_i}\partial_{x_j} u - 2\phi \partial_{x_i} \phi  \partial_{x_i} u \partial_{x_j}^2 u dx \\
        &=  \int_{B_R(x_0) }   \phi^2 |D^2_x u|^2 - 2 \phi (\Delta u) \nabla \phi \cdot \nabla u + \sum_{i,j} 2\phi \partial_{x_j} \phi \partial_{x_i} u  \partial_{x_i}\partial_{x_j} u dx.
    \end{align*}
    Applying Young Inequality for products,
        \begin{align*}
        \int_{B_R(x_0) } \phi^2  (\Delta u)^2 dx &\geq \int_{B_R(x_0) }   \phi^2 |D^2_x u|^2 - \delta \phi^2 (\Delta u)^2 -\delta^{-1} |\nabla \phi|^2 |\nabla u|^2 dx \\
        &+ \int_{B_R(x_0) }   -\sum_{i,j}  \tilde{\delta} \phi^2 (\partial_{x_i}\partial_{x_j} u)^2 -\sum_{i,j}  \tilde{\delta}^{-1}  (\partial_{x_j} \phi)^2 (\partial_{x_i} u)^2 dx \\
        &\geq \int_{B_R(x_0) }   \phi^2 |D^2_x u|^2 - \delta \phi^2 (\Delta u)^2 - \delta^{-1} |\nabla \phi|^2 |\nabla u|^2 dx \\
        &+ \int_{B_R(x_0) } - \tilde{\delta} \phi^2 |D^2_xu|^2 - \tilde{\delta}^{-1} |\nabla \phi|^2 |\nabla u|^2 dx.
    \end{align*}
    Choose $\delta = 1/2$ and $\tilde{\delta} = 1/4$ to obtain the desired inequality.
 
\end{proof}
\endgroup
\begin{lemma}
    Let $x_0 \in B^2$ and $0 < R < 1-|x_0|$. Assume that $v(x,t) \in C^2(B_R(x_0) \times (0,T), S^2)$ solves 
    \begin{align}
    v_t = \Delta v + |\nabla v|^2v, \quad (x,t) \in B_R(x_0)  \times (0,T). \label{local heat map flow}
    \end{align}
    Let $\phi(x) \in C^{\infty}_c(B_{R}(x_0) )$ be any smooth function.
    
    For any $s \geq 2$ and $t \in (0,T)$,
    \begin{align}
     \int_{B_R(x_0) }  \frac{1}{s}  \frac{d}{dt} \left( |\nabla v|^s \phi^2 \right) dx =   -\int_{B_R(x_0) } &\phi^2 |\nabla v|^{s-2} |v_t|^2 + 2\phi |\nabla v|^{s-2} \Delta v^T \nabla v^T \nabla \phi \label{integration by parts for partial_t |nabla v|^s phi^2} \\
     &+ \phi^2 \frac{s-2}{4} |\nabla v|^{s-4} |\nabla(|\nabla v|^2)|^2dx, \notag
    \end{align}
    where $(\nabla v)^T$ denotes the Jacobian matrix.
    
    In particular,
     \begin{align}
     \int_{B_R(x_0) } \frac{1}{s}\frac{d}{dt} \left( |\nabla v|^s \phi^2 \right) &+ \phi^2 |\nabla v|^{s-2} |\Delta v|^2 + \phi^2 \frac{s-2}{4} |\nabla v|^{s-4} |\nabla(|\nabla v|^2)|^2 dx \label{integration by parts for partial_t |nabla v|^s phi^2 inequality form} \\
     &\leq  - \int_{B_R(x_0)}  2\phi |\nabla v|^{s-2} \Delta v^T \nabla v^T \nabla \phi dx +  \int_{B_R(x_0) } \phi^2 |\nabla v|^s |\Delta v| dx. \notag
    \end{align}
\end{lemma}

\begin{proof}
Let $\cdot$ denote the dot product in $\mathbb R^{2 \times 2} \sim \mathbb R^4$. Observe that 
\begin{align*}
    \partial_t |\nabla v|^s = \partial_t (\nabla v \cdot \nabla v)^{s/2} = s (\nabla v \cdot \nabla v)^{s/2-1} \nabla v \cdot \nabla v_t,
\end{align*}
so integrating by parts leads to
    \begin{align}
     \int_{B_R(x_0) } \frac{1}{s}  \frac{d}{dt} \left( |\nabla v|^s \phi^2 \right) dx =   -\int_{B_R(x_0) } &\phi^2 |\nabla v|^{s-2} v_t^T \Delta v + 2\phi |\nabla v|^{s-2} v_t^T \nabla v^T \nabla \phi \label{proof: integration by parts for partial_t |nabla v|^s phi^2}  \\
     &+ \phi^2 \frac{s-2}{2} |\nabla v|^{s-4} v_t^T \nabla v^T \nabla(|\nabla v|^2)dx. \notag
    \end{align}
Since the graph of $v = (v_1,v_2)$ is a level-set $v_1(x,t)^2 + v_2(x,t)^2 = 1$, $v(x,t)$ is orthogonal to $\partial_{x_i} v(x,t)$, the columns of $\nabla v^T$, as well as $\partial_t v(x,t)$. Hence,
\begin{align}
    v_t^T \nabla v^T &= \Delta v^T \nabla v^T + |\nabla v|^2 v^T \nabla v^T=  \Delta v^T \nabla v^T  + 0 = \frac{1}{2}\nabla (|\nabla v|^2)^T \label{v_t nabla v^T identity} \\
    |v_t|^2 &= \Delta v^T v_t + |\nabla v|^2 v^T v_t =  \Delta v^T v_t  + 0 \label{|v_t| identity} \\ 
    v_t^T \Delta v &= |\Delta v|^2 + |\nabla v|^2 v^T \Delta v \label{|v_t| second identity}
\end{align}
after transposing and respectively multiplying (\ref{local heat map flow}) by $\nabla v^T$, $v_t$ and $\Delta v$. Using equalities (\ref{v_t nabla v^T identity}) and (\ref{|v_t| identity}) in (\ref{proof: integration by parts for partial_t |nabla v|^s phi^2}) proves identity (\ref{integration by parts for partial_t |nabla v|^s phi^2}). 

As for the inequality, one just needs to prove that
$$
\int_{B_R(x_0) } \phi^2 |\nabla v|^{s-2} \left( |\Delta v|^2 - |v_t|^2\right) dx \leq \int_{B_R(x_0) } \phi^2 |\nabla v|^{s}  |\Delta v| dx.
$$
Yet, (\ref{|v_t| identity})  and  (\ref{|v_t| second identity}) show that 
\begin{align*}
    |\Delta v|^2 - |v_t|^2 = -|\nabla v|^2 v^T \Delta v \leq |\nabla v|^2 |v| |\Delta v|  =  |\nabla v|^2 |\Delta v| 
\end{align*}
since $|v| = 1$. This concludes the proof.
\end{proof}

\begin{lemma}[\cite{struwe2008variational}, Chapter 3, Lemma 5.9] \label{lemma:estimate energy at a later time}
    Let $x_0 \in B^2$ and $0 < R < 1-|x_0|$. Assume that $v(x,t) \in C^2(B_R(x_0) \times (0,T), S^2)$ solves (\ref{local heat map flow})  with
    $$
    \sup_{t_1 < t < t_2} E(v(t);B_R(x_0) ) < +\infty
    $$
    for some $0 < t_1 < t_2 < T$. Then
    $$
    E(v(t_2);B_{R/2}(x_0) ) \leq E(v(t_1);B_{R}(x_0)) + \frac{32(t_2-t_1)}{R^2}\sup_{t_1 < t < t_2} E(v(t);B_{R}(x_0)),
    $$
    as well as
    \begin{align}
             \int_{t_1}^{t_2} \int_{B_{R/2}(x_0) } \frac{1}{2}  |v_t|^2dxdt  &+ E(v(t_2);B_{R/2}(x_0) )  \notag \\
             &\leq E(v(t_1);B_{R}(x_0)) + \frac{128(t_2-t_1)}{R^2}\sup_{t_1 < t < t_2} E(v(t);B_{R}(x_0)). \label{estimate on v_t L^2 norm}
    \end{align}
\end{lemma}

\begin{proof}
Let $\phi(x) \in C^{\infty}_c(B_{R}(x_0) )$, $0 \leq \phi \leq 1$, be a smooth cut-off function which is $1$ on $B_{R/2}(x_0) $ and $|\nabla \phi| \leq 4/R$ on $B_R(x_0) $.
    
    For a fixed time $0 < t_1 \leq t \leq t_2 < T$, equality (\ref{proof: integration by parts for partial_t |nabla v|^s phi^2}) with $s = 2$ implies
    \begin{align}
   \int_{B_R(x_0) }  |v_t|^2\phi^2 +  \frac{1}{2}\frac{d}{dt} \left( |\nabla v|^2 \phi^2 \right) dx &= -  \int_{B_R(x_0) } 2 \phi (\nabla \phi)^T \nabla v v_t dx  dt \notag \\
    &\leq \int_{B_R(x_0) } \delta |v_t|^2\phi^2dx + \int_{B_R(x_0) }  \delta^{-1} |\nabla v|^2|\nabla \phi|^2 dx. \label{inequality involving v_t}
    \end{align}
    
    First, we make the choice $\delta = 1$ in (\ref{inequality involving v_t}). Then
    \begin{align*}
        E(v(t_2);B_{R/2}(x_0) ) - E(v(t_1);B_{R}(x_0) ) 
        &\leq \frac{1}{2} \int_{B_{R}(x_0) }  |\nabla v(x,t) \phi(x)|^2 dx \bigg\rvert^{t = t_2}_{t = t_1} 
        \\
        &= \frac{1}{2}  \int_{t_1}^{t_2} \frac{d}{dt}\int_{B_R(x_0) } \left( |\nabla v|^2 \phi^2 \right) dx dt  \\
        &= \int_{t_1}^{t_2} \int_{B_R(x_0) } \frac{1}{2}\frac{d}{dt} \left( |\nabla v|^2 \phi^2 \right) dx dt \\
        &\leq \int_{t_1}^{t_2}  \int_{B_R(x_0) }  |\nabla v|^2|\nabla \phi|^2 dxdt \\
        &\leq \frac{16}{R^2} \int_{t_1}^{t_2}  \int_{B_R(x_0) }  |\nabla v|^2 dxdt,
    \end{align*}
where the interchange 
    $$
\int_{B_R(x_0) } \frac{1}{2}\frac{d}{dt} \left( |\nabla v|^2 \phi^2 \right) dx = \frac{d}{dt}\int_{B_R(x_0) } \frac{1}{2} |\nabla v|^2 \phi^2 dx
    $$
    is justified for any fixed $0 < t < T$ when $u \in C^2(B_R(x_0)  \times (0,T),S^2)$.

    If we make the choice $\delta = 1/2$ in (\ref{inequality involving v_t}), we find that 
     \begin{align*}
   \int_{B_{R/2}(x_0) }  \frac{1}{2} |v_t|^2 +  \frac{1}{2} \frac{d}{dt} \left( |\nabla v|^2 \phi^2 \right) dx &\leq 4 \int_{B_R(x_0) }  |\nabla v|^2|\nabla \phi|^2 dx
    \end{align*}
    and
    \begin{align*}
    \int_{t_1}^{t_2} \int_{B_{R/2}(x_0) }  \frac{1}{2}|v_t|^2dxdt &+  E(v(t_2);B_{R/2}(x_0) ) - E(v(t_1) ;B_{R/2}(x_0) ) \\
   &\leq 4 \int_{t_1}^{t_2}   \int_{B_R(x_0) }  |\nabla v|^2|\nabla \phi|^2 dx dt
   \leq  \frac{64}{R^2} \int_{t_1}^{t_2}  \int_{B_R(x_0) }  |\nabla v|^2 dxdt 
    \end{align*}
    by repeating the same argument.
\end{proof}

\begin{lemma}[\cite{struwe2008variational}, Chapter 3, Lemma 5.10] \label{lemma: estimate of nabla^2 v given small energy}
    There exists $0 < \varepsilon_1 < 1$ and $C > 0$ (independent of $x_0$, $R$, $t_1$, $t_2$ and $T$) such that if $v \in C^2(B_R(x_0) \times (0,T), S^2)$, $0 < R < 1-|x_0|$, solves (\ref{local heat map flow})  with
    $$
    \sup_{t_1 < t < t_2} E(v(t);B_R(x_0) ) < +\infty
    $$
    for some $0 \leq t_1 < t_2 \leq T$, then
        \begin{align*}
        \int_{t_1}^{t_2} \int_{B_{R/2}(x_0)} |D^2_x v|^{2}dxdt &\leq C \sup_{ t_1 < t < t_2} E(v(t);B_R(x_0) )  \left( 1 + \frac{(t_2-t_1)}{R^2} \right).
         \end{align*}
\end{lemma}

\begin{proof}
Let $\phi(x) \in C^{\infty}_c(B_{R}(x_0) )$, $0 \leq \phi \leq 1$, be a smooth cut-off function which is $1$ on $B_{R/2}(x_0) $ and $|\nabla \phi| \leq 4/R$ on $B_R(x_0) $.

    For a fixed time $0 \leq t_1 < t < t_2 \leq T$, inequality (\ref{integration by parts for partial_t |nabla v|^s phi^2 inequality form}) with $s = 2$ implies
    \begin{align*}
     I = \int_{B_R(x_0) } \frac{1}{2}\frac{d}{dt} \left( |\nabla v|^2 \phi^2 \right) + \phi^2  |\Delta v|^2 dx &\leq    \int_{B_R(x_0) } \phi^2 |\nabla v|^2 |\Delta v| dx - \int_{B_R(x_0) } 2 \phi (\nabla \phi)^T \nabla v \Delta v dx \\
     &\leq  \frac{1}{2}\int_{B_R(x_0) } \phi^2 \left( \delta_1 |\Delta v|^2  + \delta_1^{-1} |\nabla v|^4  \right) dx \\
     &+ \int_{B_R(x_0) }  \delta_2 \phi^2 |\Delta v|^2 + \delta_2^{-1} |\nabla \phi^2| | \nabla v|^2 dx 
    \end{align*}
    for any $\delta_1, \delta_2 > 0$. Choose $\delta_1 = 1/8$, $\delta_2 = 1/16$ so that 
     \begin{align*}
     I \leq \frac{1}{8} \int_{B_R(x_0) }  \phi^2 |\Delta v|^2 dx + C \int_{B_R(x_0) } \phi^2 |\nabla v|^4 dx + \frac{C}{R^2} \int_{B_R(x_0) } | \nabla v|^2 dx. 
    \end{align*}
    Apply Lemma \ref{Sobolev, estimate L4 norm in terms of H1 norm} to each $|\partial_{x_i} v_j|^4$ to get
     \begin{align*}
     I \leq \frac{1}{8} \int_{B_R(x_0) }  \phi^2 |D^2_x v|^2 dx  + C \varepsilon_1 \left( \int_{B_R(x_0) } | D^2_x v|^2 \phi^2 dx + R^{-2}  \int_{B_R(x_0) } |\nabla u|^2 dx  \right) + \frac{C}{R^2} \int_{B_R(x_0) } | \nabla v|^2 dx. 
    \end{align*}
    Choose $0 < \varepsilon_1 < 1$ small enough so that $C\varepsilon_1 \leq 1/8$ to obtain 
    $$
    I \leq \frac{1}{4} \int_{B_R(x_0) }  \phi^2 |D^2_x v|^2 dx + \frac{C}{R^2} \int_{B_R(x_0) } | \nabla v|^2 dx. 
    $$
    Using Lemma \ref{Sobolev comparison, int (Delta u)^2 > int (D^2_x u)^2}, we also obtain
    $$
 \int_{B_R(x_0) } \frac{1}{2}\frac{d}{dt} \left( |\nabla v|^2 \phi^2 \right) + \frac{1}{2} \int_{B_R(x_0) } \phi^2  |D^2_x u|^2 dx - \frac{12C^2}{R^2} \int_{B_R(x_0) } |\nabla u|^2 dx \leq I
    $$
    and integrating with respect to time yields
    \begin{align*}
        \int_{t_1}^{t_2}   \int_{B_R(x_0) } \frac{1}{4} \phi^2  |D^2_x v|^2 dx dt &\leq \frac{C}{R^2}\int_{t_1}^{t_2}   \int_{B_R(x_0) } |\nabla v|^2dxdt  - \frac{1}{2} \int_{B_R(x_0)}\left(  |\nabla v(t_2)|^2 - |\nabla v(t_1)|^2 \right)\phi^2 dx  \\
        &\leq \left( \frac{1}{2} + \frac{C (t_2-t_1)}{R^2}  \right) \sup_{t_1 < t < t_2} E(v(t);B_R(x_0) ),
    \end{align*}
    which finishes the proof.
\end{proof}

\begin{lemma}[\cite{struwe2008variational}, Chapter 3, Lemma 5.11] \label{lemma: struwe, integrability of |nabla v|^2s}
 Let $v \in C^2(B_R(x_0) \times (0,T), S^2)$, $0 < R < 1-|x_0|$, solve (\ref{local heat map flow}) with
    $$
E_0 := \sup_{0 < t < T}E(v(t), B_R(x_0)) < +\infty, \quad \nabla^2 v \in L^2(B_R(x_0) \times (0,T)).
    $$
    For all $2 \leq s \leq 2^k$, $k \in \mathbb N_{\geq 1}$, there exists $C(k) > 0$ (independent of $x_0$, $R$, $T$) and $0 < \tau = \tau(v,k,R) < T$ such that for any $0 < t_1 < t_2 \leq T$ with $t_2 - t_1 \leq \tau$,
    \begin{align*}
        \int_{t_1}^{t_2}\int_{B_{R/2^k}(x_0)} |\nabla v|^{2s} dx dt &\leq C(k) R^{-4k} \tau^2 \sup_{0 < t < T}E(v(t), B_R(x_0)) + C(k) \sup_{2 \leq  s' \leq s}\left(\int_{B_R(x_0)} |\nabla v(t_1)|^{s'}  dx  \right)^2. 
\end{align*}
More precisely, there is some absolute constant $\tilde{C}(k) > 0$ such that $0 < \tau < T$ can be taken as any constant satisfying 
\begin{equation}
    8^2 E_0 \left( \int_{t_1}^{t_2} \int_{B_{R}(x_0)} |\nabla^2v|^2 dx dt + 16^2 \tau R^{-2} E_0 \right) \leq \frac{1}{2\tilde{C}(s) + 1}, \quad \forall 0 < t_1 < t_2 < T, \ t_2 - t_1 \leq \tau. \label{eq: condition which defines tau}
\end{equation}
\end{lemma}

\begin{proof}
One assumes without loss of generality that $t_2 < T$. If the Lemma is proved in that case, then one can apply it with $t_1 = T-\tau$, $t_1 < t_2 < T$ and then pass to the limit $t_2 \to T$.

    Let $\phi(x) \in C^{\infty}_c(B_{R}(x_0) )$, $0 \leq \phi \leq 1$, be a smooth cut-off function which is $1$ on $B_{R/2}(x_0) $ and $|\nabla \phi| \leq 4/R$ on $B_R(x_0) $.

    For a fixed time $0 < t_1 < t < t_2 < T$ and $s \geq 2$, inequality (\ref{integration by parts for partial_t |nabla v|^s phi^2 inequality form}), implies
    \begin{align*}
     \int_{B_R(x_0) } \frac{1}{s}\frac{d}{dt} \left( |\nabla v|^s \phi^2 \right) &+ \phi^2 |\nabla v|^{s-2} |\Delta v|^2 + \phi^2 \frac{s-2}{4} |\nabla v|^{s-4} |\nabla(|\nabla v|^2)|^2 dx  \\
     &\leq  - \int_{B_R(x_0)}  2\phi |\nabla v|^{s-2} \Delta v^T \nabla v^T \nabla \phi dx +  \int_{B_R(x_0) } \phi^2 |\nabla v|^s |\Delta v| dx \notag \\
     &\leq \int_{B_R(x_0)}  |\nabla v|^{s-2} \left( \frac{1}{4} \phi^2 |\Delta v|^2 + 4 |\nabla v|^2 |\nabla \phi|^2 \right) dx \notag\\
     &+ \int_{B_R(x_0)} \phi^2 |\nabla v|^{s-2} \left( \frac{1}{4}|\Delta v|^2+4|\nabla v|^4 \right) dx, \notag
    \end{align*}
    where we observe that 
    $$
    \nabla \left( |\nabla v|^{\frac{s}{2}} \right) = \frac{s}{4} |\nabla v|^{\frac{s-4}{2}} \nabla \left( |\nabla v|^2 \right).
    $$
    Integrating $I$ with respect to time on $[t_1,t_2] \subset (0,T)$, one gets
    \begin{align}
        \sup_{t_1 \leq t \leq t_2} \frac{1}{s}\int_{B_R(x_0)} |\nabla v(t)|^s \phi^2 dx &+ C(s) \int_{t_1}^{t_2}  \int_{B_R(x_0) } |\nabla \left( |\nabla v|^{\frac{s}{2}} \right)|^2 \phi^2 dxdt \leq \int_{t_1}^{t_2} \int_{B_R(x_0)}  16R^{-2} |\nabla v|^{s} dxdt  \label{Struwe continuation proof, inequality |nabla v|^s}\\
        &+4\int_{t_1}^{t_2} \int_{B_R(x_0)} \phi^2 |\nabla v|^{s+2} dxdt + \frac{1}{s}\int_{B_R(x_0)} |\nabla v(t_1)|^s \phi^2  dx.  \notag 
    \end{align}
   The following application of Hölder's inequality will also be useful,
\begin{equation}
        \int_{t_1}^{t_2}\int_{B_R(x_0)} \phi^2 |\nabla v|^{s+2}  dx  \leq \left( \int_{t_1}^{t_2}\int_{B_R(x_0)} \phi^2 |\nabla v|^{2s}  dx dt\right)^{\frac{1}{2}} \left(\int_{t_1}^{t_2} \int_{B_R(x_0)}  \phi^2 |\nabla v|^{4} dxdt \right)^{\frac{1}{2}}.  \label{Struwe continuation proof, inequality |nabla v|^s+2}
\end{equation}
    Finally, we use Lemma \ref{Sobolev, estimate L4 norm in terms of H1 norm} with $u = |\nabla v|^{\frac{s}{2}}$, 
\begin{align}
        \int_{t_1}^{t_2}\int_{\mathbb R^2} \phi^2 |\nabla v|^{2s} dx dt &\leq 8 \left( \sup_{t_1 < t < t_2} \int_{B_{R}(x_0)} |\nabla v|^s dx \right) \left( \int_{t_1}^{t_2}\int_{B_{R}(x_0)} |\nabla (|\nabla v|^{\frac{s}{2}})|^2 \phi^2 dx dt \right.  \label{Struwe continuation proof, inequality |nabla v|^2s} \\
        &\left. + 4R^{-2} \int_{t_1}^{t_2}\int_{B_{R}(x_0)} |\nabla v|^s dx dt\right) \notag
\end{align}
and we plug (\ref{Struwe continuation proof, inequality |nabla v|^s}) into (\ref{Struwe continuation proof, inequality |nabla v|^2s}) to deduce
\begin{align*}
        \int_{t_1}^{t_2}\int_{\mathbb R^2} \phi^2 |\nabla v|^{2s} dx dt &\leq C(s) \left( \int_{t_1}^{t_2} \int_{B_R(x_0)}  R^{-2} |\nabla v|^{s} dxdt \right)^2 + C(s) \left(\int_{t_1}^{t_2} \int_{B_R(x_0)} \phi^2 |\nabla v|^{s+2} dxdt \right)^2 \\
        &+ C(s) \left(\int_{B_R(x_0)} |\nabla v(t_1)|^s \phi^2  dx  \right)^2 
\end{align*}
     and then (\ref{Struwe continuation proof, inequality |nabla v|^s+2}),
\begin{align}
        \int_{t_1}^{t_2}\int_{\mathbb R^2} \phi^2 |\nabla v|^{2s} dx dt &\leq C(s) \left( \int_{t_1}^{t_2} \int_{B_R(x_0)}  R^{-2} |\nabla v|^{s} dxdt \right)^2 + C(s) \left(\int_{B_R(x_0)} |\nabla v(t_1)|^s \phi^2  dx  \right)^2 \label{Struwe continuation proof, better inequality |nabla|^2s}  \\
        & + C(s) \left( \int_{t_1}^{t_2}\int_{B_R(x_0)} \phi^2 |\nabla v|^{2s}  dx dt\right)\left(\int_{t_1}^{t_2} \int_{B_R(x_0)}  \phi^2 |\nabla v|^{4} dxdt  \right). \notag
\end{align}
Observe that Lemma \ref{Sobolev, estimate L4 norm in terms of H1 norm} with $u = \partial_{x_i}v_j$, $i,j \in \{1,2\}$, leads to  
\begin{align}
   \int_{\mathbb R^2} |\nabla v|^4 \phi^2 dx &\leq 8 \sum_{i,j} \int_{\mathbb R^2} |\partial_{x_i}v_j|^4 \phi^2 dx \notag \\
   &\leq 8^2 \sum_{i,j} \left( \int_{B_{R}(x_0)} |\partial_{x_i}v_j|^2 dx \right) \left( \int_{B_{R}(x_0)} |\nabla \partial_{x_i}v_j|^2 \phi^2 dx + 16^2 R^{-2} \int_{B_{R}(x_0)} |\partial_{x_i}v_j|^2 dx\right) \notag \\
   &\leq 8^2 \left( \sup_{t_1 < t < t_2} \int_{B_{R}(x_0)} |\nabla v|^2 dx \right) \left( \int_{B_{R}(x_0)} |\nabla^2v|^2 \phi^2 dx + 16^2 R^{-2} \int_{B_{R}(x_0)} |\nabla v|^2 dx\right). \label{l4 norm of gradient} 
\end{align}
Under the hypotheses that $v$ has finite energy and $\nabla^2 v \in L^2(B_R(x_0) \times (0,T))$, for any $\delta > 0$, there is $\tau = \tau(v,R,\delta) > 0$ small enough such that if we integrate on $0 < t_1 \leq t \leq t_2 < T$ with $t_2 - t_1 \leq \tau$, then
$$
\int_{t_1}^{t_2} \int_{\mathbb R^2} |\nabla v|^4 \phi^2 dx \leq \delta.
$$
Hence, (\ref{Struwe continuation proof, better inequality |nabla|^2s}) is further estimated as
\begin{align*}
        \int_{t_1}^{t_2}\int_{\mathbb R^2} \phi^2 |\nabla v|^{2s} dx dt &\leq C(s) \left( \int_{t_1}^{t_2} \int_{B_R(x_0)}  R^{-2} |\nabla v|^{s} dxdt \right)^2 + C(s) \left(\int_{B_R(x_0)} |\nabla v(t_1)|^s \phi^2  dx  \right)^2  \\
        & + \frac{1}{2} \left( \int_{t_1}^{t_2}\int_{B_R(x_0)} \phi^2 |\nabla v|^{2s}  dx dt\right)
\end{align*}
if $\delta = 1/(2C(s)+1)$ and $t_2-t_1 \leq \tau(v,R,s)$, leading to
\begin{align*}
        \int_{t_1}^{t_2}\int_{B_{R/2}(x_0)} |\nabla v|^{2s} dx dt &\leq 2C(s) \left( \int_{t_1}^{t_2} \int_{B_R(x_0)}  R^{-2} |\nabla v|^{s} dxdt \right)^2 + 2C(s) \left(\int_{B_R(x_0)} |\nabla v(t_1)|^s  dx  \right)^2.
\end{align*}
Starting from $s = 2$, we iterate the formula to deduce the Lemma result when $s$ is a power of two. For $k \in \mathbb N_{\geq 1}$, the inequalities obtained for $s = 2, ..., 2^{k}$ all hold on the smallest ball $B_{R/2^k}(x_0)$ with the largest constant $C(k)$ and the smallest $\tau(k)$. On that region, we can use interpolation to get the inequality for the exponents $2^{k-1} \leq s \leq 2^k$.
\end{proof}

\begin{corollary} \label{cor:extension of v to time T}
     Let $v \in C^2(B_R(x_0) \times (0,T), S^2)$, $0 < R < 1-|x_0|$, solve
    $$
    v_t = \Delta v + |\nabla v|^2v, \quad (x,t) \in B_R(x_0)  \times (0,T),
    $$
    with
    $$
E_0 := \sup_{0 < t < T}E(v(t), B_R(x_0)) < +\infty, \quad \nabla^2 v \in L^2(B_R(x_0) \times (0,T)).
    $$
    Then $v$ extends to $v \in C^{\infty}(B_{R}(x_0) \times (0,T], S^2)$.
\end{corollary}

\begin{proof}
We show that for any $x \in B_R(x_0)$ and $t \in (0,T]$, for any $n,m \in \mathbb N_{\geq 0}$, there is an open neighborhood $U_{x}$ of $x$ and a neighbourhood $(t-\tau, t]$ of $t$ around which $v \in C^{n,m}_{x,t}(U_x \times (t-\tau,t])$. Without loss of generality, we may assume that $x = x_0$ and $t = T$ (the argument can be repeated on any ball $B_{R'}(x) \subset B_R(x_0)$ and for any other $t$).

    Choose any $8 < p < +\infty$. By Lemma \ref{lemma: struwe, integrability of |nabla v|^2s}, there is some $0 < r < R$ and some $0 < \tau < T$ small enough for which $|\nabla v| \in L^2(B_r(x_0) \times (T-\tau,T)) \cap L^p(B_r(x_0) \times (T-\tau,T))$ (the integrability of $|\nabla v(\cdot,t_1)|^s$, $0 < t_1 < T$, in a interior region $B_r(x_0)$ follows from continuity of $\nabla v$). One can always take $r$ and $\tau$ smaller if needed so that $\tau = r^2$ and $B_r(x_0) \times (T-\tau,T)$ corresponds to an open cylinder $C(x_0,T;r)$. Since $|v| = 1$ is bounded, one has $|\nabla v|^2v \in L^{p/2}_{x,t}(C(x,T;r))$.
    
    The interior $L^p$-estimates (Theorem \ref{thm:Interior L^p Estimates for the Heat Equation}) applied to each coordinate of $v$ imply that $\nabla^2 v \in L^{p/2}_{x,t}(C(x,T;r/2))$. Then $v_t = \Delta v + |\nabla v|^2 v \in L^{p/2}_{x,t}(C(x,T;r/2))$ as well and the parabolic Sobolev embedding (Corollary \ref{cor: local parabolic sobolev embedding}) implies that $v \in C^{1+\alpha,\alpha/2}(C(x,T;r/4))$ for some $\alpha > 0$. As $|\nabla v|^2 v$ is a product of Hölder-continuous functions on $C(x,T;r/4)$, it is Hölder-continuous as well. The Schauder estimates (Theorem \ref{interior schauder estimates}) finally imply that $\nabla^2 v, \partial_t v \in  C^{\alpha,\alpha/2}(C(x,T;r/8))$, i.e. $v \in C^{2+\alpha,1+\alpha/2}(C(x,T;r/8))$. In particular, those derivatives can be extended by continuity to the boundary $t = T$.

    Finally, we do the so-called Schauder iteration to deduce higher regularity. Observe that $w = \partial_{x_k} v \in  C^{1+\alpha,\alpha/2}(C(x,T;r/8))$ solves
    $$
    \partial_t w = \Delta w + 2 \sum_{i=1}^2 v \partial_{x_i}v \partial_{x_i} \partial_{x_k} v + |\nabla v|^2 \partial_{x_k} v = \Delta w + f(x,t),
    $$
    where the forcing term $f(x,t)$ is a product of Hölder-continuous function on $C(x,T;r/8)$, hence is Hölder-continuous as well. Schauder estimates imply that $\nabla^2 w$, $\partial_t w \in  C^{\alpha,\alpha/2}(C(x,T;r/16))$, hence, $\nabla_x v \in C^{2+\alpha,1+\alpha/2}(C(x,T;r/16))$. Next, $\tilde{w} = \partial_t v \in  C^{1+\alpha,\alpha/2}(C(x,T;r/8))$ solves
    $$
    \partial_t \tilde{w} = \Delta \tilde{w} + 2 \sum_{i=1}^2 v \partial_{x_i}v \partial_{x_i} \partial_{t} v + |\nabla v|^2 \partial_{t} v = \Delta \tilde{w} + g(x,t)
    $$
    and we can apply Schauder estimates again to deduce that $\partial_t v \in  C^{2+\alpha,1+\alpha/2}(C(x,T;r/32))$.
    
    We carry on this procedure to extend any derivative at the point $(x_0,T)$.
\end{proof}

\begin{lemma} \label{lemma: struwe analog estimate}
    Let $\phi(x) \in C^{\infty}_c(B_{R}(x_0) )$, $0 \leq \phi \leq 1$, be a smooth cut-off function which is $1$ on $B_{R/2}(x_0) $ and $|\nabla \phi| \leq 4/R$ on $B_R(x_0)$. There exists $0 < \varepsilon_1 < 1$ and $C > 0$ (independent of $x_0$, $R$, $t_1$, $t_2$, $T$ and $\phi$) such that if $v \in C^2(B_R(x_0) \times (0,T), S^2)$, $0 < R < 1-|x_0|$, solves (\ref{local heat map flow}) with
    $$
    \sup_{t_1 < t < t_2} E(v(t);B_R(x_0) ) \leq \varepsilon_1
    $$
    for some $0 < t_1 < t_2 < T$, then for all $t \in [t_1,t_2]$, one has
    \begin{align*}
        \int_{B_{R/2}(x_0)} |D^2_x v(x,t)|^{2}dx &\leq  \frac{C}{R^2} \int_{B_R(x_0) } | \nabla v(x,t)|^2 dx + C \int_{B_R(x_0) }|\partial_t v(x,t)|^2 \phi^2 dx
    \end{align*}
    and
    \begin{align*}
        \int_{t_1}^{t_2}\int_{B_{R/4}(x_0)} |D^2_x v(x,t)|^{2}dxdt &\leq C\sup_{t_1 < t < t_2} E(v(t);B_R(x_0) )\left(  1+\frac{(t_2-t_1)}{R^2}  \right).
    \end{align*}
\end{lemma}

\begin{proof}
The proof is similar to \cite[Chapter 3, Lemma 5.10]{struwe2008variational}. For a fixed time $0 < t_1 < t < t_2 < T$, $\partial_t v = \Delta v + |\nabla v|^2 v$ and $|v| = 1$ implies
    \begin{align*}
     I = \int_{B_R(x_0) } \phi^2  |\Delta v|^2 dx &\leq   C \int_{B_R(x_0) } |\nabla v|^4 \phi^2  dx + C \int_{B_R(x_0) }|\partial_t v|^2 \phi^2 dx.
    \end{align*}
    Apply Lemma 5.7 from \cite{struwe2008variational} to get
     \begin{align*}
     I \leq C \varepsilon_1 \left( \int_{B_R(x_0) } | D^2_x v|^2 \phi^2 dx + R^{-2}  \int_{B_R(x_0) } |\nabla u|^2 dx  \right) +  C \int_{B_R(x_0) }|\partial_t v|^2 \phi^2 dx.
    \end{align*}
    Choose $0 < \varepsilon_1 < 1$ small enough so that $C\varepsilon_1 \leq 1/4$ to obtain 
    $$
    I \leq \frac{1}{4} \int_{B_R(x_0) }  \phi^2 |D^2_x v|^2 dx + \frac{C}{R^2} \int_{B_R(x_0) } | \nabla v|^2 dx + C \int_{B_R(x_0) }|\partial_t v|^2 \phi^2dx.
    $$
    A simple integration by parts and Young's inequality for products shows that if $u \in H^2(B_R(x_0),\mathbb R)$, $B_R(x_0) \subset \mathbb R^2$, and $\phi(x) \in C^{\infty}_c(B_{R}(x_0))$, $0 \leq \phi \leq 1$, is any cut-off function with $|\nabla \phi|^2 \leq C/R^2$, then
    $$
\int_{B_R(x_0) } \phi^2  (\Delta u)^2 dx \geq \frac{1}{2} \int_{B_R(x_0) } \phi^2  |D^2_x u|^2 dx - 12C^2R^{-2} \int_{B_R(x_0) } |\nabla u|^2 dx
    $$
    Hence, we also obtain
    $$
 \frac{1}{2} \int_{B_R(x_0) } \phi^2  |D^2_x u|^2 dx - \frac{12C^2}{R^2} \int_{B_R(x_0) } |\nabla u|^2 dx \leq I,
    $$
    which finishes the proof for fixed time. Integrating on $[t_1,t_2]$ and using (\ref{estimate on v_t L^2 norm}) shows the other inequality.
\end{proof}

\section{Energy concentration and bubbles}

\begin{lemma}\label{lemma:no concentration of energy away origin, k equivariant}
    Let $v(x,t) \in C^{2,1}_{x,t}(B^2 \times (0,T))$ be a finite-energy and $k$-equivariant solution of (\ref{heat map flow}) on $[0,T)$.
    For any $x_0 \in \overline{B^2} \setminus \{0\}$, for any $\varepsilon > 0$, there is $\delta > 0$ for which $0 < R < \delta$ implies
    $$
    \limsup_{t \to T} E(v(t);B_R(x_0) \cap B^2) \leq \varepsilon.
    $$
    In other words, around any $x_0 \neq 0$, there is no concentration of energy as $t \to T$.
\end{lemma}

\begin{proof}
    Let $E_0 = \sup_{t \in [0,T]} E(v(t);B^2) < +\infty$ be the maximal energy of $v(x,t)$.

    Assume for a contradiction that there is $x_0 \neq 0$ for which the lemma does not hold, i.e.,  there is $\varepsilon > 0$ and a decreasing sequence $(R_n)_{n > 0} \subset (0,1]$, $R_n \to 0$, for which
    $$
\limsup_{n \to +\infty} \limsup_{t \to T} E(v(t);B_{R_n}(x_0) \cap B^2) > \varepsilon.
    $$
    As $|\nabla v|^2 = h_r^2 + \frac{k^2}{r^2}\sin(h)^2$ is radial, $E(v(t);B_R(x_0) \cap B^2) = E(v(t);B_R(x_1) \cap B^2)$ whenever $|x_0| = |x_1|$.

    Choose $M$ so that $M\varepsilon > E_0$. Choose $N$ large enough so that $R_N$ is small enough and one can choose $M$ distinct points $y_1,...,y_M$ on the circle $|y| = |x_0|$ satisfying $B_{R_N}(y_i) \cap B_{R_N}(y_j) = \emptyset$ for all $i \neq j$.

    Then
    \begin{align*}
    E_0 = \sup_{t \in [0,T]} E(v(t);B^2) \geq E(v(t);B^2) &\geq  \sum_{k = 1}^M E(v(t);B_{R_N}(y_k) \cap B^2) \\
    &=  M \cdot E(v(t);B_{R_N}(y_1) \cap B^2) \\
    &\geq M \cdot E(v(t);B_{R_n}(y_1) \cap B^2), \quad t \in [0,T],
\end{align*}
holds for all $n \geq N$, since $(R_n)_{n \geq 0}$ is decreasing. Applying $\limsup \limits_{n \to +\infty} \limsup\limits_{t \to T}$ shows that $E_0 > M\varepsilon > E_0$, a contradiction.
\end{proof}

\begin{lemma}\label{lemma:no concentration of energy away origin, general case}
    Let $v(x,t) \in C^{2,1}_{x,t}(B^2 \times (0,T))$ be a finite-energy solution to
    \begin{align*}
    v_t = \Delta v + |\nabla v|^2v, \quad (x,t) \in B^2  \times (0,T).
    \end{align*}
   For all $\varepsilon > 0$, for all but finitely many $x_0 \in B^2$ (this number might depend on $\varepsilon$), there is $\delta > 0$ for which $0 < R < \delta$ implies
    $$
    \limsup_{t \to T} E(v(t);B_R(x_0)) \leq \varepsilon.
    $$
    In other words, there is no concentration of energy as $t \to T$ around all but finitely many points in $B^2$. If $v(x,t)$ is $k$-equivariant, the only point where concentration of energy could happen is the origin.
\end{lemma}

\begin{proof}
    The proof is as in Lemma \ref{lemma:no concentration of energy away origin, k equivariant}. If there were too many points $x_0 \in B^2$ where the energy concentrates with a value $> \varepsilon$, then the sum of their local contributions to the energy would be larger than
    $$
    \sup_{0 < t < T} E(v(t);B^2),
    $$
    which is a contradiction.
\end{proof}

\begin{theorem}[\cite{struwe2008variational}, Chapter 3, Theorem 5.6]\label{thm:smoothness around all but finitely many points}
    Let $v(x,t) \in C^2(B^2 \times (0,T), S^2)$ be a finite-energy solution to
    \begin{align*}
    v_t = \Delta v + |\nabla v|^2v, \quad (x,t) \in B^2  \times (0,T),
    \end{align*}
    with $T < +\infty$. There exists finitely many points $x_1, ..., x_m \in B^2$ such that $$v(x,t) \in C^{\infty}(B^2 \times (0,T)) \cap C^{\infty}(B^2 \setminus \{x_1,...,x_m\} \times (0,T]),$$
    and if $v(x,t)$ is $k$-equivariant, then $\{x_1,...,x_m\} = \emptyset$ or $\{0\}$.

    Moreover, around any singular points $(x_0,T)$, $x_0 \in \{x_1,..,x_m\}$, for any sequence $r_n \to 0^+$, there are sequences $R_n \to 0^+$, $0 < R_n < r_n$, $T_n \to T^-$, $x_n \to x_0$ such that for almost every $- \varepsilon_1/(64E_0) < s  < 0$, up to passing to a subsequence (which might depend on $s$), $v_n(x) = v(x_n + R_nx, T_n + R_n^2s)$  converges weakly in $H^{2,2}_{loc}(\mathbb R^2,S^2)$, strongly in $H^{1,2}_{loc}(\mathbb R^2 ,S^2)$ and uniformly on any compact set to some smooth non-constant harmonic map $v_{\infty}(x)$, i.e.,  $v_{\infty}(x)$ solves $-\Delta v_{\infty} = |\nabla v_{\infty}|^2 v_{\infty}$.

    Finally, if $v$ solves (\ref{heat map flow}) with smooth and $k$-equivariant initial data and if $(0,T)$ is a singular point, then for almost every $- \varepsilon_1/(64E_0) < s  < 0$, up to passing to another subsequence, $V_n(x) = v(R_nx, T_n + R_n^2s)$ also converges weakly in $H^{2,2}_{loc}(\mathbb R^2,S^2)$, strongly in $H^{1,2}_{loc}(\mathbb R^2 ,S^2)$ and uniformly on any compact set to $V_{\infty}(x) = v_{\infty}(x+x^*)$ for some $x^* \in \mathbb R^2$. Moreover, if $h(r,t)$, $H(r)$ are respectively the inclination coordinates (\ref{heat map formulation for h(r,t)}) of $v$, $V_{\infty}$, then  $h(R_n r, T_n + R_n^2s)$ converges uniformly to $H(r)$ on all compact subsets of $[0,+\infty)$ and $H(r)$ must be of the form:
\begin{align}
    H(r) = m\pi \pm 2\arctan \left( (\alpha r)^k \right), \alpha > 0, \label{inclination coordinate of harmonic map}
\end{align}
where $h(0,t) = m\pi$ on $[0,T)$.
\end{theorem}

\begin{remark}
    In particular, $\lim \limits_{\substack{x \to x_0 \\ t \to T^-}} v(x,t)$ does not exist for any singular point $(x_0,T)$.
\end{remark}

\begin{proof}
    First, $u \in C^{\infty}(B^2 \times (0,T))$ because $u \in C^2(B_R(x_0) \times (t_0,t_1))$ satisfies the hypotheses of Corollary \ref{cor:extension of v to time T} (after applying a space and time translation) for any fixed $\overline{B_R(x_0)} \subset B^2$ and $0 < t_0 < t_1 < T$. 
    
    Around any point $x_0 \in B^2 \setminus \{x_1,...,x_m\}$ where there is no concentration of energy, one can take $0 < R < 1 - |x_0|$ small enough so that
    $$
    \sup_{0 < t < T} E(v(t);B_R(x_0)) < \varepsilon_1,
    $$
    where $\varepsilon_1$ is as in Lemma \ref{lemma: struwe analog estimate}.
    
    Then this lemma yields $\nabla^2 v \in L^2(B_{R/2}(x_0) \times (0,T))$ and we can apply Corollary \ref{cor:extension of v to time T} on $B_{R/2}(x_0) \times (0,T)$ to deduce that $v(x.t)$ can be extended to $B_{R/2}(x_0) \times (0,T]$. Hence, $v(x,t) \in C^{\infty}(B^2 \setminus \{x_1,...,x_m\} \times (0,T])$, where $\{x_1, ..., x_m\} = \emptyset$ or $\{0\}$ is $v$ has smooth and $k$-equivariant initial data (Lemma \ref{lemma:no concentration of energy away origin, k equivariant}).

Now fix a singular point $(x_s,T) \in B^2 \times \{T\}$. By concentration, for all decreasing sequences $r_n \to 0$,
\begin{align}
    \lim_{n \to +\infty} \limsup_{t \to T} E(v(t);B_{r_n}(x_s)) > \varepsilon_1, \label{singularity, concentration}
\end{align}
where the limit exists because $n \mapsto \limsup \limits_{t \to T} E(v(t);B_{r_n}(x_s))$ is bounded by $E_0$ and decreasing, hence convergent. 

Let $0 < \delta < (1-|x_s|)/8$ be small enough so that there is no other singular point in $\overline{B_{2\delta}(x_s)}$. One can prove that there are sequences $0 < R_n < r_n$, $T_n \to T^-$, $x_n \to x_s$ and an index $N \in \mathbb N$ such that one has $T_n \geq T/2, x_n \in B_{\delta}(x_s)$ and
\begin{align}
    E(v(T_n);B_{R_n}(x_n)) = \sup_{x \in B_{2\delta}(x_s)} \sup_{T/2 < t < T_n} E(v(t);B_{R_n}(x)) = \varepsilon_1  \label{energy equation leading to harmonic map extraction}
\end{align}
for all $n \geq N$. 

Let $\tau(r)$ be defined via 
 $$
\tau(r) = t_0 r^2, \quad t_0 =  \frac{\varepsilon_1}{64 E_0},
$$
and set $\tau_n = \tau(R_n)$.

Given our sequences and starting at a larger index for which $T_n - \tau(R_n) \geq T/2$, one has 
$$
E(v(T_n);B_{R_n}(x_n)) = \sup_{x \in B_{2\delta}(x_s)} \sup_{T_n - \tau(R_n) < t < T_n} E(v(t);B_{R_n}(x)) = \varepsilon_1 
$$
for all $n$.

We now define
$$
v_n(x,t) = v(x_n + R_nx, T_n + R_n^2 t): B_{\delta/R_n}(0) \times [-t_0, 0] \to S^2,
$$
which satisfies 
$$
\sup_{-t_0 < t < 0} E(v_n(t); B_{\delta/R_n}(0)) \leq E_0 = \sup_{0 < t < T} E(v(t); B^2).
$$
First, observe that
\begin{align*}
   || \partial_t v_n||_{L^2([-t_0,0] \times B_{\delta/R_n}(0))}^2 &= \int_{T_n-\tau_n}^{T_n} \int_{B_{\delta}(x_n)} |\partial_t v(y, t) |^2 dydt \\
   &\leq  \int_{T_n-\tau_n}^{T_n} \int_{B_{2\delta}(x_s)} |\partial_t v(y, t) |^2 dydt.
\end{align*}
Yet, by \cite[Lemma 5.9, Chapter 3,]{struwe2008variational} (the estimate for $|v_t|$ does not appear in the statement but in the proof of the lemma),
$$
 \int_{0}^{T} \int_{B_{2\delta}(x_s)} |\partial_t v(y, t) |^2 dydt \leq C(\delta, T, E_0) < +\infty,
$$
meaning that 
\begin{align}
    \lim_{n \to +\infty} || \partial_t v_n||_{L^2([-t_0,0] \times B_{\delta/R_n}(0))}^2 \leq \lim_{n \to +\infty} \int_{T_n-\tau_n}^{T_n} \int_{B_{2\delta}(x_s)} |\partial_t v(y, t) |^2 dydt = 0 \label{l^2 norm of partial_t v_n goes to zero}
\end{align}
since the measure of the set $[T_n - \tau_n,T_n] \times B_{2\delta}(x_s)$ goes to zero. 

In particular, $t \mapsto || \partial_t v_n(\cdot,t)||_{L^2(B_{\delta/R_n}(0))}^2$ is a sequence converging to zero in $L^1([-t_0,0])$. Up to passing to a subsequence (which we still denote with the index $n$), for all $t \in [-t_0,0] \setminus N$, $|N| = 0$, one has pointwise convergence 
\begin{align}
     \lim_{n \to +\infty} ||\partial_t v_n(\cdot, t)||_{L^2(B_{\delta/R_n}(0))}^2 = 0. \label{convergence of partial_t v_n(,t)}
\end{align}

Furthermore, if $x \in B_{\delta R_n^{-1}-1}(0)$ is fixed, then $B_1(x) \subset B_{\delta/R_n}(0)$, $x_n + R_n y \in B_{2\delta}(x_s)$ for all $y \in B_{\delta R_n^{-1}-1}(0)$, and, for all $t \in [-t_0,0]$,
\begin{align}
    E(v_n(t);B_1(x)) &= \int_{x_n + B_{R_n}(R_n x)} |\nabla v(y, T_n + R_n^2 t) |^2dy \notag \\
    &\leq \sup_{T_n - \tau_n < t < T_n}E(v(t);B_{R_n}(x_n+R_nx))  \notag \\
    &\leq \sup_{y \in B_{2\delta}(x_s)}  \sup_{T_n - \tau_n < t < T_n}E(v(t);B_{R_n}(y)) = \varepsilon_1. \label{norm less than espilon on ball of radius 1}
\end{align}
Hence, for all $x \in B_{\delta R_n^{-1}-1}(0)$, for all $t \in (-t_0,0) \setminus N$,
\begin{align}
     \int_{B_{1/2}(x)} |D^2_x v_n|^{2}dy \leq C\int_{B_1(x_0) } | \nabla v_n|^2 dx + C \int_{B_1(x_0) }|\partial_t v_n|^2 dx \quad \forall n \geq 0 \label{local hessian bound}
\end{align}
by Lemma \ref{lemma: struwe analog estimate}.

Using a covering argument (in the same spirit of those of Vitali and Besikovitch), one can extend this estimate to $B_{\delta R_n^{-1}-1}(0)$, i.e., 
\begin{align}
     \int_{B_{\delta R_n^{-1}-1}(0)} |D^2_x v_n(x,t)|^{2}dx  &\leq C_2 \left(   C \int_{B_{\delta /R_n}(0) } | \nabla v_n|^2 dx + C \int_{B_{\delta /R_n}(0) } |\partial_t v_n|^2 dx   \right) \notag \\
    &\leq C \left( E_0 + \int_{B_{\delta /R_n}(0) } |\partial_t v_n(x,t)|^2 dx   \right) \leq C(t) < +\infty \label{unif bounded hessian}
    \end{align}
for all $t \in (-t_0,0)  \setminus N$ and $n \geq 0$, where $C_2$ is a dimensional constant. If such a $t$ is fixed, the upper bound is uniform with respect to $n$ thanks to (\ref{convergence of partial_t v_n(,t)}).

Finally, fix any $s \in (-t_0,0)  \setminus N$ and let
$$
\tilde{v}_n(x) = v(x_n + R_nx, T_n + R_n^2 s): B_{\delta/R_n}(0) \to S^2.
$$
By construction, on any compact set $K \subset \mathbb R^2$, $|\tilde{v}_n| = 1$, $\partial_t v_n(\cdot,s) \to 0$ in $L^2(K)$, $D^2_x \tilde{v}_n$ is uniformly bounded in $L^2(K)$ as $n \to +\infty$, meaning that $\tilde{v}_n \to \tilde{v} \in H^{2,2}_{loc}(\mathbb R^2,S^2)$ weakly along some subsequence (which we still denote with the index $n$) and $\tilde{v}_n \to \tilde{v} \in H^{1,2}_{loc}(\mathbb R^2,S^2)$ strongly by Rellich–Kondrachov Theorem. The convergence is also uniform on any compact set because of the embeddings $H^{2,2}(K) \hookrightarrow C^{0,\alpha}(K) \subset \subset C^0(K)$, $K \subset \mathbb R^2$. Since
$$
\partial_t v_n = \Delta v_n + |\nabla v_n|^2 v_n, \quad (x,t) \in B_{\delta/R_n}(0) \times [-t_0, 0],
$$
fixing $t = s$ and passing to the limit (in the sense of distributions), we find that $\tilde{v}$ is harmonic. Finally, $\tilde{v}$ is smooth by \cite[Theorem 3.6]{sacks-uhlenbeck}, and non-constant because 
\begin{align}
    E(\tilde{v}; B_2(0)) &= \lim_{n \to +\infty} E(\tilde{v}_n; B_2(0)) = \lim_{n \to +\infty} E(v_n(s), B_2(0)) \notag \\
    &\geq \lim_{n \to +\infty} \left( E(v_n(0); B_1(0)) - 32 s E_0  \right) \notag \\
    &\geq \lim_{n \to +\infty} \left( E(v(T_n); B_{R_n}(x_n)) - 32 t_0 E_0  \right) \notag \\
    &\geq \varepsilon_1 - \varepsilon_1/2 > 0 \label{non constant limit}
\end{align} 
using (\ref{estimate on v_t L^2 norm}) (on the translated time-interval $[-t_0,0]$ instead of $[0,T]$). This concludes the proof in the general setting.

If $v$ solves (\ref{heat map flow}) with smooth and $k$-equivariant initial data, then the only possible singular point in $B^2$ is $(x_s,T) = (0,T)$. In that case, we show that $x_nR_n^{-1}$ is a bounded sequence and converges to some $x^*$ up to taking a subsequence. Then if $\tilde{v}_n(x) = v(x_n + R_nx, T_n + R_n^2 s) \to v_{\infty}(x)$ uniformly on compact sets as above, 
$$\tilde{V}_n(x) =  v(R_nx, T_n + R_n^2 s) = v_n(x-R_n^{-1}x_n) \to v_{\infty}(x-x^*) = V_{\infty}(x)$$
uniformly on compact sets as well. Similarly, $\tilde{V}_n(x) \to V_{\infty}(x)$ converges weakly in $H^{2,2}_{loc}(\mathbb R^2)$ by continuity of the translation operator in $L^2$ and the limit $V_{\infty}$ is smooth and harmonic.

Thanks to radial symmetry, one has
\begin{equation}
    E(v(T_n); B_{R_n}(Ox_n)) = E(v(T_n); B_{R_n}(x_n)) > \varepsilon_1 \label{eq:radial symmetry, eneryg rotation matrices}
\end{equation}
for all $n \geq 0$ and all rotation matrices $O$.

One can fit $m \sim |x_n|R_n^{-1}$ disjoint balls $B_{R_n}(O_1x_n)$, ..., $B_{R_n}(O_mx_n)$ in the annulus
$$
C_n = \{x \in \mathbb R^2: |x_n| - R_n < |x| < |x_n| + R_n\},
$$
meaning that
$$
E_0 \geq E(v(T_n); C_n) > \varepsilon_1 |x_n|R_n^{-1}.
$$
Hence, $x_nR_n^{-1}$ is bounded and, up to passing to a subsequence, converges in $\mathbb R^2$ to a limit $x^*$. This finishes the proof of the first part of the statement in the equivariant setting.


Finally, let $\tilde{H}_n(r) = h(R_n r, T_n + R_n^2s)$ be the inclination coordinate (\ref{heat map formulation for h(r,t)}) of $\tilde{V}_n$, where $h$ denotes the smooth inclination coordinate of $v$, $h(0,t) = m \pi$ for all $t \in [0,T)$ for some $m \in \mathbb Z$ (Proposition \ref{value of v and h at origin}, which does not require Proposition \ref{prop:smoothness at the boundary} when omitting $t = T$). Observe that $\tilde{H}_n$ is uniformly bounded in $L^{\infty}(K)$ (Proposition \ref{prop:equivalent fomrulation heat map flow}) and $|\partial_r \tilde{H}_n|$ is uniformly bounded in $L^2(K)$ for any compact set $K \subset (0,+\infty)$ because the solution $v$ has finite energy (\ref{Energy of v in terms of h}). Up to passing to a subsequence, $\tilde{H}_n$ converges weakly in $H^{1,2}_{loc}((0,+\infty)) \hookrightarrow C^{0,1/4}_{loc}((0,+\infty)) \subset \subset C^0_{loc}((0,+\infty))$ and uniformly on any compact set to a limit $H_{\infty} \in C^{0}((0,+\infty))$. One must have
\begin{equation}
    V_{\infty}(re^{i\theta}) = \left(e^{ik \theta} \sin H_{\infty}(r), \cos H_{\infty}(r) \right)  \label{eq: H vs H_infty, 1}
\end{equation}
 by pointwise convergence $\tilde{V}_n \to V_{\infty}$, $\tilde{H}_n \to H_{\infty}$. But Lemma \ref{k-equivariance, smoothness of h} shows that 
\begin{equation}
    V_{\infty}(re^{i\theta}) = \left(e^{ik \theta} \sin H(r), \cos H(r) \right) \label{eq: H vs H_infty, 2}
\end{equation}
for some lifting $H \in C^{\infty}([0,+\infty))$. Hence, $H_{\infty}(r)$ and $H(r)$ must differ by a fixed multiple of $2\pi$ and $H_{\infty}$ is smooth as well. It follows that $H_{\infty} \in C^{\infty}([0,+\infty))$ is a smooth non-constant solution for the ODE
\begin{align}
     0 = H_{rr} + \frac{H_r}{r} - k^2 \frac{\sin(2H)}{2r^2}, \quad r \in (0,+\infty), \label{time independent ode for H v2}
\end{align}
 with initial condition $H_{\infty}(0) = \tilde{m}\pi$, where $\tilde{m}$ depends on the $k$-equivariant harmonic map $V_{\infty}$ (it is shown later that $\tilde{m} = m = h(0,t)$).  All the solutions to this problem are of the form:
\begin{equation} \label{all possible form for limiting H}
    \tilde{m}\pi \pm 2\arctan \left( (\alpha r)^k \right), \alpha > 0.
\end{equation}
 Indeed, fix $r_0 > 0$ with $H_{\infty}(r_0) \in ((\tilde{m}-1)\pi,(\tilde{m}+1)\pi) \setminus \{\tilde{m}\pi\}$ (this exists as $H_{\infty}$ is non-constant). Let $\varepsilon = 1$ if $H_{\infty}(r_0) > \tilde{m}\pi$ and $\varepsilon = -1$ otherwise. By selecting $\alpha$ appropriately, one can find a solution
 $$
\chi_{\alpha}(r) = \tilde{m}\pi + \varepsilon \cdot 2\arctan \left( (\alpha r)^k \right)
 $$
 of (\ref{time independent ode for H v2}), $\chi_{\alpha}(0) =\tilde{m}\pi$, $\chi_{\alpha}(r_0) = H_{\infty}(r_0)$. The Comparison Principle (Theorem \ref{comparison principle}) applied to $\theta = H_{\infty}$ and $\psi = \chi_{\alpha}$ (and then reversing the roles of $\theta, \psi$) shows that $H_{\infty} = \chi_{\alpha}$ on $[0,r_0]$. As $H_{\infty}(r_0/2) = \chi_{\alpha}(r_0/2)$ and $\partial_r H_{\infty}(r_0/2) = \partial_r \chi_{\alpha}(r_0/2)$, standard existence and uniqueness theory for regular initial-value problems shows that $H_{\infty} = \chi_{\alpha}$ on all $[0,+\infty)$.

 Finally, observe that $\tilde{H}_n(r)$ converges uniformly near $r = 0$ as well (up to taking a final further subsequence). First note that $\tilde{H}_n(0) = h(0,T_n + R_n^2s) = m\pi$ and $\tilde{V}_n(0) = (0,0,1) \in S^2$ for all $n$. As $\tilde{V}_n \to V_{\infty}$ uniformly around $x = 0$, the first coordinate belongs to $\tilde{V}_{n,1}(x) \in [-\pi/4,\pi/4]$ for all $n \geq n_0$ and $x \in B_{K}(0)$ for some $K,n_0 > 0$. On $[0,K]$, one must have
 $$
 \tilde{H}_n(r) = m_n \pi + \varepsilon_n \arcsin \left( \tilde{V}_{n,1}(re^{i0}) \right).
$$
As $\tilde{V}_{n,1}(0) = 0$, $\tilde{H}_n(r) = m \pi$, one must have $m_n = m$, which does not depend on $n$. Up to taking a subsequence, one can assume that $\varepsilon_n \in \{-1,1\}$ is constant as well. Hence, $\tilde{H}_n(r)$ converges uniformly on $[0,K]$.
 \end{proof}

As a Corollary from Struwe's result, one deduces:

\begin{proposition}[Smoothness at $t = T$] \label{prop:smoothness at the boundary}
Let
$$
     v(t,x) \in C^0([0,T),C^1_0(\overline{B^2})) \cap C^{\infty}((0,T) \times \overline{B^2}) 
$$
solve (\ref{heat map flow}) with smooth and $k$-equivariant initial data. Let 
$$
h(t,r) \in C^0([0,T),C^1([0,1])) \cap C^{\infty}((0,T) \times [0,1])
$$
be the corresponding inclination coordinate (\ref{heat map formulation for h(r,t)}). Assume that $v$ blows-up at time $T < +\infty$.

Then 
$$
v(t,x) \in C^{\infty}((0,T] \times \overline{B^2} \setminus \{0\} ), \quad h(t,r) \in  C^{\infty}((0,T] \times (0,1] ).
$$

\end{proposition}

\begin{proof}
Recall that $v$ has finite energy (Remark \ref{k-equivariance implies finite energy}).

    If $T < +\infty$, the energy concentration argument from Struwe (Theorem \ref{thm:smoothness around all but finitely many points}) shows that $v(x,t) \in C^{\infty}(B^2 \setminus \{0\} \times (0,T])$. Hence, writing $v(x,t)$ in its symmetric form (\ref{k-equivariant form of v(x,t)}), the corresponding 
    $$
    h(t,r) \in C^0([0,T),C^1([0,1])) \cap C^{\infty}((0,T) \times [0,1])
    $$
    has additional regularity $h(t,r) \in C^{\infty}((0,T] \times (0,1))$ at time $T$ (same argument as in the beginning of the proof of Lemma \ref{k-equivariance, smoothness of h}). We prove that $h$ is smooth at $r =1, t = T$ as well.

    The function defined as
$$
    \tilde{h}_0(t,r) = \begin{cases}
        h(t,r) &\text{ if } (t,r) \in  [T/2,T] \times  \{1/2\} \\
        h(t,r) &\text{ if } (t,r) \in \{T/2\} \times [1/2,1] \\
        h_0(t,r) &\text{ if } (t,r) \in [T/2,2T] \times \{1\}
    \end{cases}
 $$
is smooth on the connected, closed set where it is defined and can be extended to a smooth function $\tilde{h}_0(t,r) \in C^{\infty}_c(\mathbb R^2)$ (e.g. by Whitney's Extension Theorem and multiplication by a cut-off). Then $h$ locally solves
\begin{align}
   h_t &= h_{rr} + \frac{h_r}{r} - k^2 \frac{\sin(2h)}{2r^2}, \quad (r,t) \in (1/2,1) \times (T/2,+\infty),  \notag \\
h(r,t) &= \tilde{h}_0(r,t), \quad (r,t) \in \{1/2,1\} \times [T/2,+\infty) \cup (1/2,1) \times \{T/2\} ,  \label{equation extending h}
\end{align}
on $[T/2,T)$. 

We rewrite (\ref{extension of h using a tilde h_0}) as a Dirichlet problem $\tilde{h} = h + \tilde{h}_0$,
\begin{align}
   \tilde{h}_t &= \tilde{h}_{rr} + \frac{\tilde{h}_r}{r} - k^2 \frac{\sin(2\tilde{h}-2\tilde{h}_0)}{2r^2} + H_0, \quad (r,t) \in (1/2,1) \times (T/2,+\infty),  \notag \\
\tilde{h}(r,t) &= 0, \quad (r,t) \in \{1/2,1\} \times [T/2,+\infty) \cup (1/2,1) \times \{T/2\}, \label{extension of h using a tilde h_0}
\end{align}
where $H_0 = \tilde{h}_{0,t} - \tilde{h}_{0,rr} - \frac{\tilde{h}_{0,r}}{r}$ and shift the time to $0$. This corresponds to a 2D nonlinear heat equation on an annulus $C$ with a smooth, globally bounded and globally Lipschitz (w.r.t. to $P,Q$) nonlinearity
\begin{align*}
    (x,t,P,Q) \mapsto F(x,t,P,Q) &= - k^2 \frac{\sin(2P-2\tilde{h}_0)}{2|x|^2} + H_0(x,t) \\
    &= - k^2 \frac{\sin(2P-2\tilde{h}_0)}{2|x|^2} \chi_{(1/4,2)}(|x|) + H_0(x,t),
\end{align*}
where $\chi_{(1/4,2)} \in C^{\infty}_c(\mathbb R;[0,1])$ is a smooth cut-off which is $1$ on $(1/2,1)$ and zero on $(1/4,2)^c$. We can then apply Theorem \ref{thm:taylor local existence} and Theorem \ref{thm:resgularty solution nonlinear pb} to deduce that $\tilde{h} \in C^{\infty}((T/2,+\infty),\overline{C})$ must be global. In particular, $h$ extends as a function $h \in C^{\infty}((T/2,T] \times [0,1])$.
\end{proof}

An argument similar to Theorem \ref{thm:smoothness around all but finitely many points} allows to extract all the bubbles at the same time, leading to the following well-known bubble decomposition of blowing-up solutions along an appropriate time sequence. One should note that the decomposition holds only in the energy space and not uniformly. The harmonic maps appearing in the decomposition are called \textit{bubbles} of energy, as they account for all the energy that $v(t)$ loses at the blow-up time (see (\ref{pythagorean decomposition of energy})).

\begin{theorem}[Bubble decomposition along a sequence of times]\label{thm: full bubble decomposition}
    Let $v(x,t) \in C^{2}(B^2 \times (0,T), S^2)$ be a finite-energy solution to
    \begin{align*}
    v_t = \Delta v + |\nabla v|^2v, \quad (x,t) \in B^2  \times (0,T)
    \end{align*}
    with $T < +\infty$. 

    Around any singular points $(x_s,T)$, $x_s \in B^2$, there exists a sequence $T_n \to T^-$, there exists finitely many smooth non-constant harmonic maps $\{\omega_1, ..., \omega_p\}$ from $\mathbb R^2 \cup \{\infty\} \simeq S^2$ to $S^2$, there exists sequences $(x_n^{(i)})_{n \geq 0}$, $i \in \{1,...,p\}$, all converging to $x_s$, there exists positive sequences $(R_n^{(i)})_{n \geq 0}$, $i \in \{1,...,p\}$, all converging to zero, satisfying for all $i \neq j$,
    \begin{align}
        \frac{R_n^{(i)}}{R_n^{(j)}} + \frac{R_n^{(j)}}{R_n^{(i)}} + \frac{|x_n^{(i)}-x_n^{(j)}|^2}{R_n^{(i)}R_n^{(j)}} \to +\infty, \quad n \to +\infty \label{different rate of convergence of scaling of bubbles}
    \end{align}
    and
 \begin{align}
     v(x,T_n) - \sum_{i=1}^p \left[ \omega_i \left(  \frac{x-x_n^{(i)}}{R_n^{(i)}}\right) - \omega_i(\infty) \right] \to v(x,T), \quad n \to +\infty \label{multi bubble decomposition along a sequence of times}
 \end{align}
    strongly in $H^{1,2}(B_{\delta}(x_s))$ for all balls $B_{\delta}(x_s)$ such that $\overline{B_{2\delta}(x_s)} \subset B^2$ contains no other singular points from $\{x_1,...,x_m\}$. Moreover,
\begin{equation}
        \lim_{t \to T^-} E(v(t);B_{\delta}(x_s)) = E(v(T);B_{\delta}(x_s)) + \sum_{i=1}^p E(\omega_i;\mathbb R^2). \label{pythagorean decomposition of energy}
\end{equation}
\end{theorem}

\begin{proof}
    See \cite[Theorem 1 and Remark 6]{Qing1995OnSO}. The argument is essentially local around the singularity and does not depend on the boundary.
\end{proof}

\begin{corollary}[Recovering the bubbles and decomposition in the $k$-equivariant setting]\label{rmk:recovering the bubbles}
    Let $v$ be as in Theorem \ref{thm: full bubble decomposition}. Order the indices $\{1,...,p\}$ so that $R_n^{(i)}$ has faster or comparable decay to $R_n^{(j)}$, $j > i$, meaning that $R_n^{(i)}/R_n^{(j)}$ is bounded. For $j \in \{1,...,p\}$, up to taking a subsequence, there is $M = M(j)$ large enough and some constant $c_j \in \mathbb R^3$ for which
    \begin{align*}
    v_n^{(j)}(x) &= \left(v(y,T_n) - \sum_{i=1}^{j-1} \left[ \omega_i \left(  \frac{y-x_n^{(i)}}{R_n^{(i)}}\right) - \omega_i(\infty) \right] \right)_{| y = x_n^{(j)} + R_n^{(j)}x} \\
    &=: V_n^{(j-1)}(x_n^{(j)} + R_n^{(j)}x)
\end{align*}
converges strongly to $\omega_{j}(x) + c_{j}$ in $H^{1,2}_{loc}(\mathbb R^2,B_{M}(0))$. When $j = 1$, the convergence and the limit $\omega_1(x) + c_1$ are actually in $H^{1,2}_{loc}(\mathbb R^2,S^2)$.

In the $k$-equivariant setting, keeping the same sequences $T_n$ and $R_n$, (\ref{multi bubble decomposition along a sequence of times}) can be rewritten as:
$$
v(x,T_n) - \sum_{i=1}^p \left[ \tilde{\omega}_i \left(  \frac{x}{R_n^{(i)}}\right) - \tilde{\omega}_i(\infty) \right] \to v(x,T), \quad n \to +\infty.
$$
where $\tilde{\omega}_i(x) = \omega_i(x+x_i^*)$, $x_i^* \in \mathbb R^2$, is a smooth, non-constant, harmonic map from $\mathbb R^2 \cup \{\infty\} \simeq S^2$ to $S^2$ for which $|\nabla \tilde{\omega}_i|$ is radial. Moreover, for $i \in \{1,...,p\}$, there exists $\tilde{c}_i \in \mathbb R^3$ for which
$$
v(R_n^{(i)}x,T_n) \to \tilde{\omega}_i(x) + \tilde{c}_i
$$
in $H^{1,2}_{loc}(\mathbb R^2,S^2)$, where $\tilde{\omega}_i(x) + \tilde{c}_i$ is a $k$-equivariant harmonic map from $\mathbb R^2 \cup \{\infty\} \simeq S^2$ to $S^2$ with inclination coordinate of the form (\ref{inclination coordinate of harmonic map}).
\end{corollary}

\begin{remark}
    On the level of the inclination coordinate $h$ of $v$, we have uniform convergence of $h(R_n^{(i)}r,T_n)$ to the inclination coordinate of $\tilde{\omega}_i(x) + \tilde{c}_i$ on any compact subset $K \subset (0,+\infty)$ thanks to the one-dimensional Sobolev embedding  $H^{1,2}_{loc}((0,+\infty)) \hookrightarrow C^{0,1/4}_{loc}((0,+\infty))$.
\end{remark}

\begin{proof}
Up to taking a subsequence, $R_n^{(i)}/R_n^{(j)}$ converges to a finite value in $\mathbb R_{\geq 0}$. Up to passing to further subsequences, we will also assume that for all $i \neq j$,
$$
\frac{x_n^{(i)}-x_n^{(j)}}{R_n^{(i)}}
$$
either converges in $\mathbb R^2$ or is unbounded.

Let
\begin{align*}
    A_n^{(i)} = \int_{B_{KR_n^{(1)}}(x_n^{(1)})}\left| \left(R_n^{(i)}\right)^{-1} \nabla \omega_i \left(\frac{x-x_n^{(i)}}{R_n^{(i)}} \right) \right|^2 dx = \int_{B_{KR_n^{(1)}/R_n^{(i)}}\left( \frac{x_n^{(1)}-x_n^{(i)}}{R_n^{(i)}}\right)}\left| 
    \nabla \omega_i \left(x \right) \right|^2 dx.
\end{align*}
If $i = 1$, one has $A_n^{(1)} = \int_{B_K(0)} |\nabla \omega_1(x)|^2 dx > 0$. Assume now that $i \neq 1$. If $R_n^{(1)}/R_n^{(i)} \to 0$, then $A_n^{(i)} \to 0$ because the domain of integration is shrinking to a zero measure set. Else, $R_n^{(1)}/R_n^{(i)} \to c > 0$ and one has
$$
\left( \frac{|x_n^{(1)}-x_n^{(i)}|}{R_n^{(i)}} \right)^2 \gtrsim \frac{|x_n^{(1)}-x_n^{(i)}|^2}{R_n^{(1)}R_n^{(i)}} \to +\infty,
$$
i.e., the domain $B_{KR_n^{(1)}/R_n^{(i)}}\left( \frac{x_n^{(1)}-x_n^{(i)}}{R_n^{(i)}}\right)$ is drifting to infinity, while keeping a measure which is uniformly bounded with respect to $n$. In other words, the characteristic function of the domain $B_{KR_n^{(1)}/R_n^{(i)}}\left( \frac{x_n^{(1)}-x_n^{(i)}}{R_n^{(i)}}\right)$ converges, in a pointwise sense, to zero as $n \to +\infty$ and it follows from Dominated Convergence Theorem that $A_n^{(i)} \to 0$ as well.

Similarly, the mixed term
\begin{align*}
    B_n^{(i,j)} = \int_{B_{KR_n^{(1)}}(x_n^{(1)})} \left| \left(R_n^{(i)}\right)^{-1} \nabla \omega_i \left(\frac{x-x_n^{(i)}}{R_n^{(i)}} \right) \right| \cdot \left| \left(R_n^{(j)}\right)^{-1} \nabla \omega_j \left(\frac{x-x_n^{(j)}}{R_n^{(j)}} \right) \right| dx
\end{align*}
converges to zero as $n \to +\infty$ by Cauchy-Schwarz inequality and the convergence of $A_n^{(i)}$.

One also has 
\begin{align*}
    C_n = \int_{B_{KR_n^{(1)}}(x_n^{(1)})}\left| \nabla v(x,T) \right|^2 dx \to 0, \quad n \to +\infty,
\end{align*}
as the domain is shrinking to a zero measure set and
\begin{align*}
    D_n^{(i)} = \int_{KB_{R_n^{(1)}}(x_n^{(1)})}\left| \nabla v(x,T) \right| \cdot  \left| \left(R_n^{(i)}\right)^{-1} \nabla \omega_i \left(\frac{x-x_n^{(i)}}{R_n^{(i)}} \right) \right|  dx \to 0, \quad n \to +\infty,
\end{align*}
for any $i \in \{1,...,p\}$ by Cauchy-Schwarz inequality.

It follows from the decomposition that 
$$
\int_{B_{KR_n^{(1)}}(x_n^{(1)})} \left| \nabla v(x,T_n) - \left(R_n^{(1)}\right)^{-1} \nabla \omega_1 \left( \frac{x-x_n^{(1)}}{R_n^{(1)}} \right) \right|^2 dx \to 0,
$$
i.e., if $v_n^{(1)}(x) = v(x_n^{(1)} + R_n^{(1)}x, T_n)$, then $\nabla v_n^{(1)}(x) \to \nabla \omega_1$ in $L^2(B_K(0))$. By Poincaré's inequality, up to taking a subsequence, $v_n^{(1)}$ converges strongly in $H^{1,2}_{loc}$ to $\omega_1(x) + c_1 \in H^{1,2}_{loc}$ for some constant $c_1$. As $|v| = 1$, convergence is in $H^{1,2}_{loc}(\mathbb R^2,S^2)$.

 Similarly, one proves that 
\begin{align*}
    v_n^{(j+1)}(x) &= \left(v(y,T_n) - \sum_{i=1}^j \left[ \omega_i \left(  \frac{y-x_n^{(i)}}{R_n^{(i)}}\right) - \omega_i(\infty) \right] \right)_{| y = x_n^{(j+1)} + R_n^{(j+1)}x} \\
    &=: V_n^{(j)}(x_n^{(j+1)} + R_n^{(j+1)}x)
\end{align*}
converges strongly to $\omega_{j+1} + c_{j+1}$ in $H^{1,2}_{loc}(\mathbb R^2,B_{M}(0))$ for some fixed $M = M(j)$ large enough (the maps $v$ and $\omega_i$ are bounded, but we can no longer ensure $|v_n^{(j+1)}| = 1$).


In the $k$-equivariant setting, one can proceed as in the end of the proof of Theorem \ref{thm:smoothness around all but finitely many points} (see the energy argument below equation (\ref{eq:radial symmetry, eneryg rotation matrices})) to show that $x_n^{(1)}/R_n^{(1)}$ is bounded thanks to radial symmetry and, up to taking a subsequence, converges. Hence, one can always replace $(x_n^{(1)})$ by $(0)$ in the decomposition and $\omega_1(x)$ by a translate $\omega_1(x + x_1^*)$ without changing the conclusion of Theorem \ref{thm: full bubble decomposition}. As $|\nabla v|$ is radial, so is the case for the limiting harmonic map, i.e., $x \mapsto |\nabla \omega_1(x+x_1^*)|$ is radial. Assume then by induction that $(x_n^{(i)}) = (0)$ and that $|\nabla \omega_i|$ is radial for $i \leq j$. In particular, $|\nabla V_n^{(j)}|$ is radial.

Then 
$$
E(V_n^{(j)};B_{KR_n^{(j+1)}}(x_n^{(j+1)})) = E(v_n^{(j+1)};B_{K}(0)) \to E(\omega_{j+1};B_K(0))  > 0, \quad n \to +\infty,
$$
for some fixed $K = K(j)$ large enough and 
$$
E(V_n^{(j)};B_{KR_n^{(j+1)}}(x_n^{(j+1)})) \leq E(V_n^{(j)};B_{1}(0)) \leq E_0 + \sum_{i=1}^j E(\omega_i;\mathbb R^2),
$$
where $E_0$ is our usual upper bound on the energy of $v(x,t)$ on $B_1(0) \times (0,T)$. As in the end of the proof of Theorem \ref{thm:smoothness around all but finitely many points}, radial symmetry shows that $x_n^{(j+1)}/R_n^{(j+1)}$ is bounded and, up to taking a subsequence, converges. Then one can always replace $(x_n^{(j+1)})$ by $(0)$ in the decomposition and $\omega_{j+1}(x)$ by a translate $\omega_{j+1}(x + x_{j+1}^*)$. As $|\nabla V_n^{(j)}|$ is radial, the translate of the harmonic map is such that $x \mapsto |\nabla \omega_{j+1}(x+x_{j+1}^*)|$ is radial, which finishes the induction proof.

We have now proved that (\ref{multi bubble decomposition along a sequence of times}) can be rewritten as:
$$
v(x,T_n) - \sum_{i=1}^p \left[ \tilde{\omega}_i \left(  \frac{x}{R_n^{(i)}}\right) - \tilde{\omega}_i(\infty) \right] \to v(x,T), \quad n \to +\infty.
$$
in $H^{1,2}_{loc}(B_{\delta}(0))$, where $\tilde{\omega}_i(x) = \omega_i(x+x_i^*)$ is a non-constant harmonic map for which $|\nabla \tilde{\omega}_i|$ is radial. It is then easier to recover each $\tilde{\omega}_i$. Simply observe that 
$$
v(R_n^{(i)}x,T_n) \to \tilde{\omega}_i(x)
$$
in $\dot{H}^{1,2}_{loc}(\mathbb R^2, B_{M}(0))$. By Poincaré's inequality, up to taking subsequences, there exists $c_i$ so that $v(R_n^{(i)}x,T_n) \to \tilde{\omega}_i(x) + c_i$ strongly in $H^{1,2}_{loc}(\mathbb R^2,S^2)$. As the image of $v(R_n^{(i)}x,T_n)$ is in $S^2$, it follows that $\tilde{\omega}_i + c_i$ is also a non-constant harmonic map from $\mathbb R^2 \cup \{\infty\}$ to $S^2$.

Repeating the ODE argument (see (\ref{time independent ode for H v2})) at the end of Theorem \ref{thm:smoothness around all but finitely many points} shows that $\tilde{\omega}_i(x) + c_i$ is $k$-equivariant with inclination coordinate of the form (\ref{inclination coordinate of harmonic map}).
\end{proof}

\backmatter
\cleardoublepage
\phantomsection
\nocite{*}
\bibliographystyle{alpha}
\bibliography{tail/References.bib}

\addcontentsline{toc}{chapter}{Bibliography}

\end{document}